\documentclass[reqno,a4paper]{amsart} 
\usepackage[english]{babel}

\usepackage{graphicx,amsmath,amsfonts,amssymb,amsthm, mathrsfs, tikz,tikz-cd, amscd,amsbsy,amstext,mathtools, enumerate, enumitem,stmaryrd, tensor,a4wide, fancyhdr,ifthen, empheq, verbatim, latexsym, epigraph,tabularray, etoolbox, patchcmd}

\mathtoolsset{showonlyrefs,showmanualtags}

\usepackage{mathpazo}

\usepackage{nameref,zref-xr} 
\zxrsetup{toltxlabel} 
\zexternaldocument*[foundation-]{Nilpotent_COHA_v2} 
\zexternaldocument*[torsion-pairs-]{Torsion_pairs_final} 

\usetikzlibrary{positioning,shapes,shadows,arrows}

\tikzstyle{roundbox} = [rectangle, draw, text centered, rounded corners, minimum height=2em]
\tikzstyle{connector} = [draw, -latex']

\usepackage{color}
\definecolor{MyDarkBlue}{rgb}{0.15,0.25,0.45}

\usepackage[utf8x]{inputenc}

\usepackage[toc,page]{appendix}
\usepackage[mathscr]{euscript} 
\allowdisplaybreaks

\newif\ifpersonal
\personaltrue 
\ifpersonal

\newcommand{\todo}[1]{\textcolor{red}{(Todo: #1)}}
\else
\newcommand\ast{\personal}[1]{\ignorespaces}
\newcommand\ast{\todo}[1]{\ignorespaces}
\fi

\newcommand{\calA}{\mathcal{A}}

\newcommand{\calC}{\mathcal{C}}

\newcommand{\calE}{\mathcal{E}}
\newcommand{\calF}{\mathcal{F}}
\newcommand{\calG}{\mathcal{G}}
\newcommand{\calH}{\mathcal{H}}
\newcommand{\calI}{\mathcal{I}}

\newcommand{\calK}{\mathcal{K}}
\newcommand{\calL}{\mathcal{L}}
\newcommand{\calM}{\mathcal{M}}
\newcommand{\calN}{\mathcal{N}}
\newcommand{\calO}{\mathcal{O}}
\newcommand{\calP}{\mathcal{P}}
\newcommand{\calQ}{\mathcal{Q}}

\newcommand{\calS}{\mathcal{S}}
\newcommand{\calT}{\mathcal{T}}
\newcommand{\calU}{\mathcal{U}}
\newcommand{\calV}{\mathcal{V}}
\newcommand{\calW}{\mathcal{W}}

\newcommand{\A}{\mathbb{A}}
\newcommand{\C}{\mathbb{C}}
\newcommand{\E}{\mathbb{E}}
\newcommand{\F}{\mathbb{F}}
\newcommand{\G}{\mathbb{G}}
\newcommand{\HH}{\mathbb{H}}
\newcommand{\I}{\mathbb{I}}
\newcommand{\J}{\mathbb{J}}
\newcommand{\K}{\mathbb{K}}
\newcommand{\LL}{\mathbb{L}}
\newcommand{\N}{\mathbb{N}}
\newcommand{\PP}{\mathbb{P}}
\newcommand{\Q}{\mathbb{Q}}
\newcommand{\R}{\mathbb{R}}
\newcommand{\bS}{\mathbb{S}}

\newcommand{\V}{\mathbb{V}}

\newcommand{\Y}{\mathbb{Y}}
\newcommand{\Z}{\mathbb{Z}}

\newcommand{\scrC}{\mathscr{C}}
\newcommand{\scrD}{\mathscr{D}}
\newcommand{\scrE}{\mathscr{E}}
\newcommand{\scrF}{\mathscr{F}}

\newcommand{\scrO}{\mathscr{O}}

\newcommand{\scrS}{\mathscr{S}}
\newcommand{\scrT}{\mathscr{T}}
\newcommand{\scrU}{\mathscr{U}}
\newcommand{\scrV}{\mathscr{V}}

\newcommand{\sfB}{\mathsf{B}}

\newcommand{\sfE}{\mathsf{E}}
\newcommand{\sfG}{\mathsf{G}}
\newcommand{\sfH}{\mathsf{H}}

\newcommand{\sfK}{\mathsf{K}}

\newcommand{\sfP}{\mathsf{P}}

\newcommand{\sfT}{\mathsf{T}}
\newcommand{\sfU}{\mathsf{U}}

\newcommand{\sff}{\mathsf{f}}

\newcommand{\frakM}{\mathfrak{M}}

\newcommand{\frakQ}{\mathfrak{Q}}
\newcommand{\frakS}{\mathfrak{S}}
\newcommand{\frakU}{\mathfrak{U}}
\newcommand{\frakX}{\mathfrak{X}}

\newcommand{\frakY}{\mathfrak{Y}}

\newcommand{\frakZ}{\mathfrak{Z}}

\newcommand{\frake}{\mathfrak{e}}
\newcommand{\frakg}{\mathfrak{g}}
\newcommand{\frakh}{\mathfrak{h}}
\newcommand{\frakm}{\mathfrak{m}}
\newcommand{\frakn}{\mathfrak{n}}
\newcommand{\frakp}{\mathfrak{p}}
\newcommand{\frakr}{\mathfrak{r}}
\newcommand{\fraks}{\mathfrak{s}}

\newcommand{\bfC}{\mathbf{C}}
\newcommand{\bfD}{\mathbf{D}}

\newcommand{\bfH}{\mathbf{H}}

\newcommand{\bfS}{\mathbf{S}}
\newcommand{\bfT}{\mathbf{T}}

\newcommand{\bfY}{\mathbf{Y}}
\newcommand{\bfX}{\mathbf{X}}
\newcommand{\bfW}{\mathbf{W}}

\newcommand{\bfc}{\mathbf{c}}
\newcommand{\bfd}{\mathbf{d}}
\newcommand{\bfe}{\mathbf{e}}

\newcommand{\bfh}{\mathbf{h}}

\newcommand{\bfn}{\mathbf{n}}

\newcommand{\bfv}{\mathbf{v}}
\newcommand{\bfw}{\mathbf{w}}

\newcommand{\bfy}{\mathbf{y}}
\newcommand{\bfz}{\mathbf{z}}

\newcommand{\op}{^{\mathsf{op}}}

\newcommand{\Cat}{\mathsf{Cat}}

\newcommand{\CAlg}{\mathsf{CAlg}}
\newcommand{\dAff}{\mathsf{dAff}}

\newcommand{\PrLomega}{\mathsf{Pr}^{\mathsf{L},\omega}}

\newcommand{\BettiD}{\tensor*[^{\mathsf{top}}]{\bfD}{}}

\newcommand{\PreSt}{\mathsf{PreSt}}
\newcommand{\dSt}{\mathsf{dSt}}

\newcommand{\Corr}{\mathsf{Corr}}

\newcommand{\Ind}{\mathsf{Ind}}
\newcommand{\Pro}{\mathsf{Pro}}

\newcommand{\Fun}{\mathsf{Fun}}

\newcommand{\Sym}{\mathsf{Sym}}

\newcommand{\Red}[1]{\tensor*[^{\mathsf{red}}]{#1}{}}

\newcommand{\catMod}{\mathsf{Mod}}
\newcommand{\catmod}{\mathsf{mod}}

\newcommand{\catP}{\mathsf{P}}
\newcommand{\catPps}{\mathsf{P}_{\mathsf{ps}}}

\newcommand{\Mod}{\textrm{-} \mathsf{Mod}}

\newcommand{\catPnil}{\mathsf{P}_C}

\newcommand{\modPii}[1]{\mathsf{mod}_{#1}(\Pi_\qv)}

\newcommand{\modPi}{\catmod(\Pi_\qv)}
\newcommand{\ModPi}{\catMod(\Pi_\qv)}
\newcommand{\nilpPi}{\mathsf{nilp}(\Pi_\qv)}

\newcommand{\catCoh}{\mathsf{Coh}}

\newcommand{\catCohps}{\mathsf{Coh}_{\mathsf{ps}}}

\newcommand{\catQCoh}{\mathsf{QCoh}}

\newcommand{\catDb}{\mathsf{D}^\mathsf{b}}
\newcommand{\catDbps}{\mathsf{D}^\mathsf{b}_\mathsf{ps}}

\newcommand{\ps}{\mathsf{ps}}

\newcommand{\elldLambdaqv}{\tensor*[^\ell]{\mathbf{\Lambda}}{_\qv}}
\newcommand{\iLambda}{\tensor*[^{(i)}]{\Lambda}{}}
\newcommand{\Lambdai}{\Lambda^{(i)}}

\newcommand{\stackBun}{\mathfrak{Bun}}
\newcommand{\stackCoh}{\mathfrak{Coh}}

\newcommand{\stackCohpsnil}{\mathfrak{Coh}_{\mathsf{ps}}^{\mathsf{nil}}}

\newcommand{\trunc}[1]{\mathsf{t}_0(#1)}

\newcommand{\dstackCoh}{\mathbf{Coh}}

\newcommand{\dstackCohps}{\mathbf{Coh}_{\mathsf{ps}}}
\newcommand{\dstackCohpsnil}{\mathbf{Coh}_{\mathsf{ps}}^{\mathsf{nil}}}
\newcommand{\dstackRep}{\mathbf{Rep}}
\newcommand{\stackRep}{\mathfrak{Rep}}

\newcommand{\dstackCohext}{\mathbf{Coh}^{\mathsf{ext}}}

\newcommand{\dstackRepext}{\mathbf{Rep}^{\mathsf{ext}}}

\newcommand{\dLambda}{\mathbf{\Lambda}}
\newcommand{\dLambdaext}{\mathbf{\Lambda}^{\mathsf{ext}}}
\newcommand{\idLambda}{\tensor*[^{(i)}]{\mathbf{\Lambda}}{}}
\newcommand{\dLambdai}{\mathbf{\Lambda}^{(i)}}

\newcommand{\sfIrr}{\mathsf{Irr}}

\newcommand{\stackPic}{\mathfrak{Pic}}
\newcommand{\dstackPic}{\mathbf{Pic}}

\newcommand{\Hilb}{\mathsf{Hilb}}

\newcommand{\sfFun}{\mathsf{Fun}}

\newcommand{\ii}[1]{\tensor*[^i]{#1}{}}

\newcommand{\iii}[1]{\tensor*[^{(i)}]{#1}{}}

\newcommand{\Tmax}{T_{\mathsf{max}}}
\newcommand{\Ttilde}{\widetilde{T}}

\newcommand{\qv}{\mathcal{Q}}
\newcommand{\qvfin}{\mathcal{Q}_\sff}
\newcommand{\doubleqv}{\overline{\qv}}
\newcommand{\doubleOmega}{\overline{\Omega}}
\newcommand{\Kr}{\mathsf{Kr}}

\newcommand{\Lalpha}{\check{\alpha}}
\newcommand{\Lomega}{\check{\omega}}
\newcommand{\Ltheta}{\check{\theta}}
\newcommand{\Lthetafin}{\check{\theta}_\sff}
\newcommand{\Lvarphi}{\check{\varphi}}
\newcommand{\Llambda}{\check{\lambda}}
\newcommand{\Lmu}{\check{\mu}}
\newcommand{\Lrho}{\check{\rho}}

\newcommand{\rootsetfin}{\Delta_\sff}

\newcommand{\rootlattice}{\bfY}
\newcommand{\corootlattice}{\check{\bfY}}
\newcommand{\rootlatticefin}{\rootlattice_\sff}
\newcommand{\corootlatticefin}{\check{\rootlattice_\sff}}

\newcommand{\weightlattice}{\bfX}
\newcommand{\coweightlattice}{\check{\weightlattice}}

\newcommand{\coweightlatticefin}{\coweightlattice_\sff}

\newcommand{\frakgfin}{\frakg_\sff}
\newcommand{\frakhfin}{\frakh_\sff}
\newcommand{\fraknfin}{\frakn_\sff}

\newcommand{\sfex}{\mathsf{ex}}
\newcommand{\sfaf}{\mathsf{af}}

\newcommand{\almostsame}{\stackrel{\bullet}{=}}

\newcommand{\ix}{\tensor[^{(i)}]{x}{}}

\newcommand{\rsv}{X}
\newcommand{\ssv}{X_{\mathsf{con}}}

\newcommand{\frakgl}{\mathfrak{gl}}
\newcommand{\fraksl}{\mathfrak{sl}}
\newcommand{\frakgell}{\mathfrak{g}_{\mathsf{ell}}}
\newcommand{\fraknell}{\mathfrak{n}_{\mathsf{ell}}}

\newcommand{\fraknellzero}{\mathfrak{n}_{\mathsf{ell}, 0}}

\newcommand{\eY}{\mathbb{Y}^{\mathfrak{e}}}
\newcommand{\zeroeY}{\mathbb{Y}^{\mathfrak{e},\, 0}}

\newcommand{\DiYqv}{\mathbf{D}_i\Y_\qv}

\newcommand{\gr}{\mathsf{gr}}
\newcommand{\iY}{\tensor*[^{(i)}]{\mathbb{Y}}{}}
\newcommand{\iYzero}{\iY_0}
\newcommand{\Yi}{\mathbb{Y}^{(i)}}
\newcommand{\Yizero}{\Yi_0}

\newcommand{\Yzero}{\Y_0}
\newcommand{\eYzero}{\Y^{\mathfrak{e}}_0}

\newcommand{\YMO}{\Y_{\mathsf{MO}}}

\newcommand{\ad}{\mathsf{ad}}

\newcommand{\uce}{\mathfrak{uce}}

\newcommand{\iU}{\tensor[^{(i)}]{\sfU}{}}
\newcommand{\Ui}{{\sfU}^{(i)}}

\newcommand{\HBM}{\mathsf{H}^{\mathsf{BM}}}

\newcommand{\Hbullet}{\mathsf{H}^\bullet}
\newcommand{\HBMbullet}{\mathsf{H}_\bullet^{\mathsf{BM}}}
\newcommand{\HBMbulletT}{\mathsf{H}_\bullet^T}
\newcommand{\HBMbulletA}{\mathsf{H}_\bullet^A}

\newcommand{\HBMDGamma}{\mathsf{H}^{\bfD,\Gamma}}

\newcommand{\HBMGamma}{\mathsf{H}^{\mathsf{BM},\Gamma}}

\newcommand{\coha}{\mathbf{HA}}

\newcommand{\shuffle}{\mathbf{Sh}}

\newcommand{\cohaqv}{\coha_\qv}
\newcommand{\cohaqvd}{\coha_{\bfd}}

\newcommand{\Dicohaqv}{\mathbf{D}_i\cohaqv^{\Ttilde}}

\newcommand{\cohazero}{\coha_0}

\newcommand{\icohaqv}{\tensor*[^{(i)}]{\coha}{^T_\qv}}
\newcommand{\cohai}{\coha^{T, (i)}}
\newcommand{\cohaiqv}{\coha^{T, (i)}_\qv}
\newcommand{\icohad}{\tensor*[^{(i)}]{\coha}{^T_{\bfd}}}
\newcommand{\cohaid}{\coha^{T, (i)}_{\bfd}}

\newcommand{\icohatildeqv}{\tensor*[^{(i)}]{\coha}{^{\Ttilde}_\qv}}

\newcommand{\cohatildeiqv}{\coha^{\Ttilde, (i)}_\qv}
\newcommand{\icohatilded}{\tensor*[^{(i)}]{\coha}{^{\Ttilde}_{\bfd}}}
\newcommand{\cohatildeid}{\coha^{\Ttilde, (i)}_{\bfd}}

\newcommand{\icohatildezero}{\tensor*[^{(i)}]{\coha}{^{\Ttilde}_0}}
\newcommand{\cohatildeizero}{\coha^{\Ttilde, (i)}_0}

\newcommand{\wcohaqv}{\tensor*[^{(w)}]{\coha}{^T_\qv}}
\newcommand{\cohaqvw}{\cohaqv^{T, (w)}}

\newcommand{\cohavw}{\tensor*[^{(v)}]{\coha}{^{T, (w)}_\qv}}

\newcommand{\ires}{\tensor*[^{(i)}]{\mathsf{res}}{}}
\newcommand{\resi}{\mathsf{res}^{(i)}}
\newcommand{\iresd}{\tensor*[^{(i)}]{\mathsf{res}}{_\bfd}}
\newcommand{\resid}{\mathsf{res}^{(i)}_\bfd}

\newcommand{\Exp}{\mathsf{Exp}}

\newcommand{\bfLambda}{\mathbf{\Lambda}}

\newcommand{\id}{{\mathsf{id}}}
\newcommand{\Spec}{\mathsf{Spec}}

\newcommand{\ch}{\mathsf{ch}}

\newcommand{\bfDelta}{\mathbf{\Delta}}

\newcommand{\Spc}{\mathsf{Spc}}

\newcommand{\ev}{\mathsf{ev}}

\newcommand{\Hom}{\mathsf{Hom}}
\newcommand{\calHom}{\mathcal{H}\mathsf{om}}
\newcommand{\Ext}{\mathsf{Ext}}

\newcommand{\Aut}{\mathsf{Aut}}
\newcommand{\End}{\mathsf{End}}
\newcommand{\Tor}{\mathsf{Tor}}

\newcommand{\Pic}{\mathsf{Pic}}
\newcommand{\GL}{\mathsf{GL}}
\newcommand{\SL}{\mathsf{SL}}

\DeclareMathOperator*{\fcolim}{``colim''}
\DeclareMathOperator*{\flim}{``lim''}

\newcommand{\bfjmathhat}{\hat{\boldsymbol{\jmath}}}

\newcommand{\red}{\mathsf{red}}

\newcommand{\pr}{\mathsf{pr}}

\makeatletter
\newcommand{\colim@}{%
	\vtop{\m@th\ialign{##\cr
			\hfil$\operator@font colim$\hfil\cr
			\noalign{\nointerlineskip\kern1.5\ex@}\cr
			\noalign{\nointerlineskip\kern-\ex@}\cr}}%
}
\newcommand{\colim}{%
	\mathop{\mathpalette\colim@{\textstyle}}\nmlimits@
}
\makeatother

\usepackage{thmtools}

\declaretheoremstyle[
spaceabove=1em, spacebelow=1em,
headfont=\bfseries, notefont=\normalfont, bodyfont=\itshape,
headpunct={.}, notebraces={(}{)}, postheadspace={ }
]{basic-theorem}

\declaretheoremstyle[
spaceabove=1em, spacebelow=1em,
headfont=\bfseries, notefont=\normalfont, bodyfont=\normalfont,
headpunct={.}, notebraces={(}{)}, postheadspace={ }, qed=$\oslash$
]{basic-definition}

\declaretheoremstyle[
spaceabove=1em, spacebelow=1em,
headfont=\itshape, notefont=\normalfont, bodyfont=\normalfont,
headpunct={.}, notebraces={(}{)}, postheadspace={ }, qed=$\triangle$
]{basic-remark}

\patchcmd{\part}{\normalfont}{\large\scshape}{}{}

\theoremstyle{basic-theorem}

\newtheorem{theorem}{Theorem}[section]
\newtheorem{corollary}[theorem]{Corollary}
\newtheorem{lemma}[theorem]{Lemma}
\newtheorem{proposition}[theorem]{Proposition}
\newtheorem{conjecture}[theorem]{Conjecture}

\declaretheorem[style=basic-definition, numberlike=theorem]{definition}

\declaretheorem[style=basic-definition, numberlike=theorem]{notation}
\declaretheorem[style=basic-definition, numberlike=theorem]{warning}

\declaretheorem[style=basic-definition, numberlike=theorem]{construction}
\declaretheorem[style=basic-remark, numberlike=theorem]{remark}

\newtheorem{assumption}{Assumption}



\numberwithin{equation}{section}



\declaretheorem[style=basic-theorem,name=Theorem]{theoremintroduction}

\declaretheorem[style=basic-remark, numbered=no, name=Remark]{remarkunnum}
\declaretheorem[style=basic-theorem, numbered=no, name=Corollary]{corollaryunnum}

\usepackage{hyperref}

\hypersetup{anchorcolor=black, linktocpage=true,
	colorlinks=true,
	citecolor=MyDarkBlue,
	linkcolor=MyDarkBlue,
	urlcolor=black,
	pdftitle={Cohomological Hall algebras of one-dimensional sheaves on surfaces and Yangians}, 
	breaklinks=true,
	plainpages=true
}
\usepackage[hyperpageref]{backref}

\title[Cohomological Hall algebras and Yangians]{Cohomological Hall algebras of one-dimensional sheaves on surfaces and Yangians}

\author[D.-E.~Diaconescu]{Duiliu-Emanuel Diaconescu}
\address[Duiliu-Emanuel Diaconescu]{New High Energy Theory Center - Serrin Building, Rutgers, The State University Of New Jersey, 126 Frelinghuysen Rd., Piscataway, NJ 08854-8019, USA}
\curraddr{}
\email{\href{mailto:duiliu@physics.rutgers.edu}{duiliu@physics.rutgers.edu}}

\author[M.~Porta]{Mauro Porta}
\address[Mauro Porta]{Institut de recherche mathématique avancée (IRMA), Université de Strasbourg, France and Institut Universitaire de France (IUF)}
\curraddr{}
\email{\href{mailto:porta@math.unistra.fr}{porta@math.unistra.fr}}

\author[F.~Sala]{Francesco Sala}
\address[Francesco Sala]{Università di Pisa, Dipartimento di Matematica, Largo Bruno Pontecorvo 5, 56127 Pisa (PI), Italy}
\address{Kavli IPMU (WPI), UTIAS, The University of Tokyo, Kashiwa, Chiba 277-8583, Japan}
\curraddr{}
\email{\href{mailto:francesco.sala@unipi.it}{francesco.sala@unipi.it}}

\author[O.~Schiffmann]{Olivier Schiffmann}
\address[Olivier Schiffmann]{Laboratoire de Math\'ematiques d'Orsay, Universit\'e de Paris-Sud Paris-Saclay, B\^at. 425, 91405 Orsay Cedex, France, UMR8628 (CNRS), and Simion Stoilow Institute of Mathematics, Bucharest, Romania}
\email{\href{mailto:olivier.schiffmann@universite-paris-saclay.fr}{olivier.schiffmann@universite-paris-saclay.fr}}

\author[E.~Vasserot]{Eric Vasserot}
\address[Eric Vasserot]{Université de Paris, 75013 Paris, France, UMR7586 (CNRS) and Institut Universitaire de France (IUF)}
\email{\href{mailto:eric.vasserot@imj-prg.fr}{eric.vasserot@imj-prg.fr}}

\thanks{The work of the first-named author is partially supported by NSF grant DMS-1802410. The work of the third-named author is partially supported by JSPS KAKENHI Grant Numbers JP21K03197 and JP26K06721. The work of the fourth-named author is partially supported by the PNRR Grant `\textit{Cohomological Hall algebras and smooth surfaces and applications}', CF~44/14.11.2022. The work of the fifth-named author is partially supported by ANR-18-CE40-0024 ``Categorification in topology and representation theory - Catore''. This work was also partially supported by the ``National Group for Algebraic and Geometric Structures, and their Applications'' (GNSAGA – INDAM)} 



\begin{document}

\begin{abstract}
	
	This paper provides the first algebraic characterization of an algebra of cohomological Hecke operators associated with modifications of coherent sheaves on a smooth surface $X$ along a fixed proper curve $Z \subset X$ (possibly singular and reducible), establishing a direct connection with Yangians. It is based on the theory of equivariant nilpotent cohomological Hall algebras $\coha^T_{X,Z}$, developed by the same authors in \cite{DPSSV-1}.
	
	\medskip
	
	More precisely, let $X$ be a resolution of a Kleinian singularity (for example, $X\coloneqq T^\ast\PP^1$) and let $Z$ be the exceptional divisor.
	One of the main results of this paper is an explicit isomorphism
		\begin{align}
		\coha^T_{X,Z} \simeq \Y^+_\infty(\frakg)\ ,
	\end{align}
	where $\Y^+_\infty(\frakg)$ is a completed, \textit{nonstandard}, positive half of the affine Yangian $\Y(\frakg)$ of the corresponding affine ADE Lie algebra $\frakg$. Under this identification, the natural action of $\Pic(X)$ on $\coha^T_{X,Z}$ corresponds to the action of the extended affine braid group on $\Y(\frakg)$. Furthermore, the generators of $\coha^T_{X,Z}$--given by fundamental classes of substacks of zero-dimensional sheaves and of pushforwards of line bundles on $Z$--are expressed explicitly in terms of Yangian generators.
	
	\medskip
	
	Our main tools, which may be of independent interest, are
	\begin{enumerate}[label=(\roman*)]\itemsep0.2cm
		\item A `continuity' theorem describing the behavior of cohomological Hall algebras of objects in the heart of $t$-structures $\tau_n$ when the sequence $(\tau_n)_n$ converges, in an appropriate sense, to a fixed $t$-structure $\tau_\infty$.
		\item The definition of a multi-parameter Yangian $\Y_\qv$ for an arbitrary quiver $\qv$, given by generators and relations.
		\item A theorem relating the algebraic action of the braid group $B_\qv$ on the Yangian $\Y_\qv$ to the action of $B_\qv$ on the equivariant 2-dimensional cohomological Hall algebra $\coha^T_\qv$ of $\qv$, where the latter can be described in terms of derived reflection functors of the bounded derived category of modules over the preprojective algebra of $\qv$.
	\end{enumerate}
\end{abstract}

\maketitle
\thispagestyle{empty}

\tableofcontents

\section*{Introduction}

Let $X$ be a smooth quasiprojective variety. A \textit{punctual modification} of a coherent sheaf $\calE$ on $X$ is a subsheaf $\calF \subset \calE$ for which $\calE/\calF$ is zero-dimensional. Hecke operators associated with punctual modifications have traditionally played a central role in studying of moduli stacks or moduli spaces of (semi)stable vector bundles on smooth projective curves. There are several instances in which one can consider algebras of such Hecke operators. The oldest and most studied context involves functions on the moduli stack of vector bundles on a smooth projective curve defined over a finite field, as in the function field version of the original Langlands program, and their constructible sheaf analogs, as in the geometric Langlands program. In these cases, the resulting algebra is fully understood: it decomposes as a product of spherical affine Hecke algebras, which are commutative polynomial algebras. 

The study of \textit{cohomological} Hecke operators associated with punctual modifications of torsion-free sheaves on smooth \textit{surfaces} has played a central role in geometric representation theory and algebraic geometry. This line of research originates in the pioneering works of Nakajima \cite{Nakajima1997, Nakajima_Lectures} and Grojnowski \cite{Grojnowski} in the 1990s, which investigated the cohomology of Hilbert schemes of points on surfaces.
The structure of the resulting algebra of operators has been intensively studied by many authors over the last three decades, culminating in a complete understanding in \cite{MMSV}.
We will revisit the cornerstones of this theory below.

The goal of this paper, together with its companions \cite{DPSSV-1, DPSSV-2} is to go \textit{beyond} the punctual case, developing a framework to study modifications along curves and providing the first systematic structural results in this case.
More precisely, let $X$ be a smooth quasiprojective complex surface; a \textit{curve modification} of a coherent sheaf $\calE$ on $X$ is a subsheaf $\calF \subset \calE$ with $\calE/\calF$ set-theoretically supported on a proper (possibly singular and reducible) curve.
In the context of this work, we fix once and for all the curve $Z$ along which modifications happen.
A prototypical example arises when $X\coloneqq \mathsf{Tot}(\calL)$ is the total space of a line bundle $\calL$ over a smooth projective curve $C$, with modifications along the zero section $Z\coloneqq C$. 

We begin by briefly reviewing the known results concerning Hecke operators associated with punctual modifications and the motivations related to those of curve modifications, before presenting our main results.

\subsection*{Punctual modifications of coherent sheaves on smooth surfaces} 

Let $X$ be a smooth quasiprojective complex surface. The first occurrence of cohomological Hecke operators associated with punctual modifications on surfaces can be traced back to the foundational works of Nakajima \cite{Nakajima1997, Nakajima_Lectures} and Grojnowski \cite{Grojnowski}. Using such operators, they constructed an action of the Heisenberg-Clifford algebra $\mathsf{Heis}_X$, modeled on the cohomology $\sfH^\bullet(X)$ of $X$, on the direct sum $\V_X$ of the Borel-Moore homologies of Hilbert schemes $\Hilb^\ell(X)$ of $\ell$ points on $X$. This yields a canonical isomorphism between $\V_X$ and the Fock space representation of $\mathsf{Heis}_X$, thereby creating an extremely powerful dictionary between the representation theory of Heisenberg algebras and the geometry of Hilbert schemes of points on smooth surfaces. This dictionary was later further enriched by the work of many authors. Notably, Lehn~\cite{Lehn_Virasoro} provided a geometric construction of the action of the Virasoro algebra on $\V_X$. 

By the work of Baranosky~\cite{Baranovsky2000} (see also \cite{Licata_Savage}), this construction admits natural extensions to other moduli spaces, such as the \textit{instanton moduli space} $M_{\PP^2}(r)$, parametrizing rank $r$ \textit{framed} torsion-free sheaves on $\PP^2$ (see \cite[Chapter~2]{Nakajima_Lectures}), and the moduli spaces of Gieseker-stable torsion-free sheaves with fixed rank and first Chern class on a smooth projective complex surface $X$, which are higher-rank versions of $\Hilb^\ell(\C^2)$ and $\Hilb^\ell(X)$, respectively. 

As noted in \cite{Baranovsky2000}, higher-rank moduli spaces no longer yield irreducible or cyclic $\mathsf{Heis}_X$-modules. For the instanton moduli space $M_{\PP^2}(r)$, a solution emerged from the Alday-Gaiotto-Tachikawa correspondence \cite{AGT}, which relates two-dimensional conformal field theories to supersymmetric gauge theories on $\R^4$. This correspondence suggests replacing $\mathsf{Heis}_{\C^2}$ with the affine W-algebra $\calW(\frakgl(r))$, endowing the equivariant cohomology of $M_{\PP^2}(r)$ with the structure of a Verma module of $\calW(\frakgl(r))$, as proved in \cite{SV_Cherednik, MO_Yangian, BFN-AGT, RSYZ_vertex, DSYZ}. A key aspect of \cite{SV_Cherednik}'s construction involves the (equivariant) \textit{cohomological Hall algebra} (COHA) of zero-dimensional sheaves on $\C^2$. This COHA serves as the biggest algebra of cohomogical Hecke operators associated with punctual modifications which acts on the (equivariant) cohomology of instanton moduli spaces. 

The COHA $\coha_0(X)$ of zero-dimensional sheaves on a smooth surface $X$ has been introduced in \cite{Zhao_Hall, COHA_surface}, generalizing \cite{SV_Cherednik} when $X=\C^2$, and in \cite{Minets_Higgs, Sala_Schiffmann} when the surface is the cotangent of a curve. This algebra naturally acts on the homology of moduli spaces of (semi)stable sheaves via punctual modifications (see, e.g. \cite{COHA_surface, MMSV}), recovering the `positive half' of the algebras of cohomological Hecke operators associated with punctual modifications introduced in \cite{Nakajima1997, Nakajima_Lectures,Grojnowski, Baranovsky2000} and discussed above.

The algebraic structure of $\coha_0(X)$ has been determined explicitly for $\A^2$ in \cite{SV_Cherednik,DavisonCOHAA2}, and for an arbitrary cohomologically pure surface $X$ in \cite{MMSV}\footnote{In \textit{loc.cit.}, the cohomological Hall algebra is denoted $\bfH_0(X)$.}. Concretely, $\coha_0(X)$ is isomorphic to the positive half of the \textit{deformed $\sfH^\bullet(X)$-colored $W_{1+\infty}$-algebra} defined in \cite{MMSV} (see \cite{NegutIHES} and references therein for the $K$-theoretical version). Equivalently, it realizes a $\sfH^\bullet(X)$-colored version of the \textit{Yangian of} $\widehat{\frakgl}_1$, which is a deformation of the enveloping algebra of the Lie superalgebra $\widehat{\frakgl}_1[z]\otimes \sfH^\bullet(X)$. The zero-degree modes subalgebra $\widehat{\frakgl}_1\otimes \sfH^\bullet(X)$ recovers the Heisenberg-Clifford algebra, while degree $\leq 1$ modes yield Lehn's Virasoro algebra. 

The explicit algebraic description of $\coha_0(X)$ has played a crucial role in several important applications. If $X=\C^2$,  it was essential for proving the Alday-Gaiotto-Tachikawa conjectures for pure supersymmetric gauge theories on $\R^4$ in \cite{SV_Cherednik}. If $X$ is the cotangent bundle of a smooth projective complex curve, this explicit algebraic description proved fundamental in establishing the proof of the $P=W$ conjecture in \cite{HMMS}. 

\subsection*{Curve modifications of coherent sheaves on smooth surfaces}

Similarly to the case of punctual modifications, interest in Hecke operators associated with \textit{curve} modifications also dates back to the 1990s. For example, in their work on a `surface analog' of the geometric Langlands program, the fifth-named author, together with Ginzburg and Kapranov \cite{GKV_Langlands}, introduced an algebra of Hecke operators associated with modifications of coherent sheaves on a smooth surface $X$ along an ADE configuration of $\PP^1$'s, acting on the space of functions of certain moduli stacks of coherent sheaves on $X$. For an elliptic surface and a special fiber of affine type \textit{ADE}, this algebra was extended to a two-parameter version of the quantum toroidal algebra of the same ADE quiver. Around the same time, motivated by gauge-theoretic considerations, Nakajima \cite{Nakajima1996} speculated about the existence of an algebra of cohomological Hecke operators associated with curve modifications. The potential role of such operators in gauge theories on four-manifolds remains an active area of investigation in theoretical physics, see, e.g., \cite{Feigin_Gukov}.

The construction of algebras of cohomological Hecke operators for curve modifications of \textit{perverse} coherent sheaves was initiated in Nakajima's work on ALE spaces \cite{Nakajima1994}, preceding his investigations of punctual modifications on Hilbert schemes. In this setting, where perverse coherent sheaves correspond to specific quiver representations, the algebras of these Hecke operators can be explicitly described via two-dimensional \textit{cohomological Hall algebras associated with quivers}. These were originally introduced by the fourth and fifth-named authors in \cite{SV12, SV_Elliptic, SV_Cherednik} (see also \cite{YZ_preproj_COHA}), with a different terminology. Recently, a complete algebraic characterization of these COHAs in terms of Maulik-Okounkov Yangians \cite{MO_Yangian} was given independently in \cite{BD_Okounkov}, building on \cite{DM_DT}, and \cite{SV_Yangians=COHA}, extending \cite{SV_generators, SV_Yangians}. Furthermore, the study of modification of perverse coherent sheaves has been extended to the three-dimensional setting \cite{RSYZ, BR_Perverse}, relating them to Kontsevich-Soibelman's three-dimensional COHAs \cite{KS_Hall} and conjecturally to Li-Yamazaki \textit{quiver Yangians} \cite{Li_Yamazaki}. Finally, quiver varieties have played a pivotal role in realizing geometric actions of Lorentzian Kac-Moody algebras (resp.\ toroidal algebras) on the homology of moduli spaces of rank one sheaves on K3 surfaces (resp.\ elliptic surfaces) in \cite{DeHority_K3, DeHority_Elliptic}, making use of the birational geometry of these moduli spaces. 

The original approach outlined in \cite{Nakajima1996} uses Hecke modifications associated with \textit{parabolic} bundles. This framework was later formalized in \cite{Nakajima_Lie, Yoshioka_Lie}, yielding actions of Lie algebras, associated with exceptional objects, on the cohomology of moduli spaces of Gieseker-stable sheaves on K3 surfaces. This yields the first instance where Hecke operators associated with curve modifications of coherent sheaves, as opposed to perverse coherent sheaves, are constructed. As the preceding examples show, these Lie algebras constitute only a small part of the full algebra of cohomological Hecke operators one seeks to construct. In this context, the theory of two-dimensional COHAs offers a promising framework to realize the full algebra of cohomological Hecke operators associated with curve modifications.

The existing theory of COHAs for coherent sheaves on smooth surfaces \cite{COHA_surface}, which generalizes earlier work on both the zero-dimensional case \cite{Zhao_Hall} and properly supported coherent sheaves on cotangent bundles of smooth curves \cite{Minets_Higgs, Sala_Schiffmann}, needs to be modified for the effective study of these Hecke operators, especially for modifications along a \textit{fixed} curve. This adaptation is crucial for establishing a direct connection to Yangians, analogous to the established correspondence between (possibly nilpotent) quiver COHAs and Maulik-Okounkov Yangians \cite{SV_Yangians, SV_Yangians=COHA}. Such a connection to Yangian-type algebras is a first important step towards studying and exploiting the action of these algebras on moduli spaces of interest.

\subsection*{Our work}

A central objective of this work is to develop a unified framework achieving two fundamental aims: first, the construction of the largest algebra of cohomological Hecke operators associated with curve modifications of coherent sheaves along a fixed proper curve $Z$, which may be singular and reducible; and second, the derivation of an explicit connection to Yangian-type quantum groups. This program is initiated and in part achieved through the sequence of works \cite{DPSSV-1, DPSSV-2}, and the present paper. 

In \cite{DPSSV-1}, we developed the general theory of the COHA $\coha_{X,Z}$ of coherent sheaves on a smooth quasiprojective complex surface $X$ with set-theoretic support on a proper curve $Z\subset X$. One may view this algebra as a \textit{nilpotent} version of the COHA introduced in \cite{COHA_surface}. The companion paper \cite{DPSSV-2} establishes a PBW-type theorem relating the COHA $\coha_{X,Z}$ for a reducible curve $Z$ with the COHAs $\coha_{X,Z_i}$ of the irreducible components $Z_i$ of $Z$. This allows in principle to reconstruct $\coha_{X,Z}$ from the individual pieces $\coha_{X,Z_i}$, and from local commutation relations between intersecting irreducible components $Z_i, Z_j$.

\medskip

The present paper provides the first computation of $\coha_{X,Z}$ and consequently of the algebra of cohomological Hecke operators of curve modifications of coherent sheaves, in the important case of Kleinian singularities. Using the results of \cite{DPSSV-2}, this paves the way to the study of the COHA $\coha_{X,Z}$ whenever $Z$ is a union of $(-2)$ rational curves which intersect transversally.

More precisely, fix a finite subgroup $G \subset \SL(2,\C)$, let $\pi\colon\rsv \to \ssv$ be the minimal resolution of the Kleinian singularity $\ssv \coloneqq \C^2/G$, and let $Z=C\coloneqq \pi^{-1}(0)\subset \rsv$ be the exceptional divisor. We also allow an arbitrary diagonal torus $T \subset \GL(2,\C)$ which commutes with the action of $G$. Let $\qv=(I,\Omega)$ denote the McKay quiver of $G$, and let $\frakgfin$ and $\frakg\coloneqq\widehat{\frakg}_\sff$ be the corresponding simple Lie algebra and its affinization. The affine Yangian $\Y_\qv$ is a deformation of the enveloping algebra $\sfU(\frakgell)$ of the universal central extension
\begin{align}
	\frakgell\coloneqq \frakgfin[s^{\pm 1},t] \oplus K \quad \text{with } K\coloneqq \bigoplus_{\ell\in \N}\Q c_\ell \oplus \bigoplus_{\genfrac{}{}{0pt}{}{\ell\in \N, \; \ell\geqslant 1}{k\in \Z, \; k\neq 0}} \Q c_{k, \ell}
\end{align}
of the double loop algebra $\frakgfin[s^{\pm 1},t]$. Consider the following nonstandard positive half of $\frakgell$
\begin{align}
	\frakgell^+\coloneqq \frakn_+[s^{\pm 1},t] \oplus s^{-1}\frakh[s^{-1},t] \oplus K_-\quad\text{where } K_-\coloneqq \bigoplus _{k<0} \Q c_{k, \ell}\ .
\end{align}

The first main result of the paper is the following.
\begin{theoremintroduction}[Theorems~\ref{thm:coha-surface-as-limit2} \& \ref{thm:coha-surface-as-limit3}]\label{thmintro:Kleinian_Sing} 
	The following hold.
	\begin{enumerate}[label=(\roman*)]\itemsep0.2cm
		\item Assume that $T=\{\id\}$. There is an algebra isomorphism
		\begin{align}
			\begin{tikzcd}[ampersand replacement=\&]
				\Theta_{\rsv,C}\colon\coha_{\rsv,C} \ar{r}{\sim}\& \widehat{\sfU}(\frakgell^+)
			\end{tikzcd}
		\end{align}
		where $\widehat{\sfU}(\frakgell^+)$ is an explicit completion of $\sfU(\frakgell^+)$.
		
		\item There is an algebra isomorphism
		\begin{align}
			\begin{tikzcd}[ampersand replacement=\&]
				\Theta_{\rsv,C} \colon\coha^T_{\rsv,C} \ar{r}{\sim} \& \Y^+_\infty
			\end{tikzcd}
		\end{align}
		where $\Y^+_\infty$ is a filtered deformation of $\widehat{\sfU}(\frakgell^+)$ defined over $\sfH^\bullet_T$. 
	\end{enumerate}	
\end{theoremintroduction}
The algebra $\Y^+_\infty$ is defined as a limit of quotients of the standard negative half $\Y^-_\qv \subset \Y_\qv$ in which the transition morphisms are given by truncated braid group operators.

\begin{remarkunnum}
	The isomorphisms $\Theta_{\rsv,C}$ intertwine the action of $\coweightlatticefin$ on $\coha^T_{\rsv,C}$ and $\widehat{\sfU}(\frakgell^+)$ or $\Y^+_\infty$ given by tensor product by line bundles, under the isomorphism $\Pic(\rsv)\simeq \coweightlatticefin$, and braid group action, respectively. Furthermore, the isomorphisms $\Theta_{\rsv, C}$ also intertwine the action of multiplication by tautological Chern classes on $\dstackCohpsnil(\widehat{\rsv}_C)$ and the action of explicit derivations on $\widehat{\sfU}(\frakgell^+)$ and $\Y^+_\infty$. 
\end{remarkunnum}

Furthermore, we can compute explicitly the image under $\Theta_{\rsv, C}$ of the fundamental classes of some natural irreducible components of $\stackCohpsnil(\widehat{\rsv}_C)$. More precisely, let $C_i$ be an irreducible component of $C_\red$ and let 
\begin{align}
	\gamma_i \colon \fraksl_2 \longrightarrow \frakgfin \quad \text{and} \quad \gamma_{i,\mathsf{ell}}\colon \frakgl_{2,\mathsf{ell}} \longrightarrow \frakgell 
\end{align}
be the canonical embeddings of the corresponding root subalgebras. Let $\frakZ_{i,n}$ be the irreducible component of $\stackCohpsnil(\widehat{\rsv}_C)$ parametrizing sheaves scheme-theoretically supported on $C_i$, of rank $1$ and degree $n-1$ on $C_i$, and let $\frakY_{i,d}$ be the irreducible component parametrizing zero-dimensional sheaves of length $d$ with punctual support at a point of $C_i$.

The second main result of the paper is the following.
\begin{theoremintroduction}[Theorem~\ref{thm:class_irred_component}, Corollary~\ref{cor:functoriality-Theta}]\label{thmintro:generators}
	\hfill
	\begin{enumerate}[label=(\roman*)]\itemsep0.2cm
		\item \label{item:intro-i} Assume that $G=\Z/ 2\Z$, $\rsv=T^\ast \PP^1$, and $C=\PP^1$. Then, for any $n \in \Z$ and $d \in \N$ 
		\begin{align}
			\sum_{n \in \Z} (-1)^n\Theta_{T^\ast \PP^1, \PP^1}([\frakZ_{n}])u^{-n}&= \left( \sum_{n \in\Z} x^+s^{-n}u^{-n}\right) \cdot \exp\left(\sum_{k \geq 1}hs^{-k}\frac{u^{-k}}{k}\right)\ ,\\
			\Theta_{T^\ast \PP^1, \PP^1}([\frakY_{d}])&=(-1)^{d-1}hs^{-d}\ ,
		\end{align}
		where $x^+,f,x^-$ is the standard basis of $\fraksl_2$. 
		
		\item \label{item:intro-ii} For an arbitrary $G$ and any irreducible component $C_i$ of $C$, we have a commutative square
		\begin{align}
			\begin{tikzcd}[ampersand replacement=\&]
				\coha_{T^\ast \PP^1,\PP^1} \arrow[swap]{d}{\Theta_{T^\ast \PP^1,\PP^1}} \arrow{r}{\sim}\& \coha_{\rsv,C_i} \arrow{r}{} \& \coha_{\rsv,C} \arrow{d}{\Theta_{\rsv,C}}\\
				\widehat{\sfU}(\frakgl_{2,\mathsf{ell}}^+) \arrow{rr}{\gamma_{i,\mathsf{ell}}} \& \&\widehat{\sfU}(\frakgell^+)
			\end{tikzcd}
		\end{align}
		of algebra morphisms, where the top left arrow is induced by the isomorphism of formal neighborhoods $\widehat{T^\ast \PP^1}_{\PP^1} \simeq \widehat{\rsv}_{C_i}$, while the top right arrow by the functoriality of COHAs with respect to closed embeddings.
	\end{enumerate}	
\end{theoremintroduction}
Combining the claims~\ref{item:intro-i} and \ref{item:intro-ii} yields an explicit description of $\Theta_{\rsv,C}([\frakZ_{i,n}])$ and $\Theta_{\rsv,C}([\frakY_{i,d}])$ for any $i, n,d$.  

\begin{corollaryunnum}[Corollary~\ref{cor:generation_HAXC}]
	$\coha_{\rsv,C}^A$ is topologically generated, as an algebra over $\bS_{\rsv,C}$\footnote{The \textit{ring of universal tautological classes} $\bS_{\rsv,C}$ in Formula~\ref{eq:universal-tautological}.}, by the fundamental classes $[\frakY_{i,d}]$, $[\frakZ_{i,n}]$ for $i \in I_\sff, d \in \N$ and $n\in \Z$. 
\end{corollaryunnum}

\subsubsection*{Strategy of proof} 

The construction of the isomorphisms $\Theta_{\rsv, C}$ is obtained by a `limit' process, using the \textit{derived McKay equivalence} \cite{Kleinian_derived,VdB_Flops} $\catDbps(\catCoh(\rsv))\simeq \catDbps(\ModPi)$ together with the weak action of $B_\sfex$ on these triangulated categories. More precisely, we show that the derived McKay correspondence restricts to an equivalence
\begin{align}
	\sfP_C(\rsv/\ssv) \simeq \nilpPi
\end{align}
where $\sfP_C(\rsv/\ssv)$ is the heart of \textit{perverse coherent} sheaves on $\rsv$ with set-theoretic support in $C$. Composing this equivalence with the autoequivalence $\bfS_{\textrm{-}\Lthetafin}$ given by the tensor product by line bundles $\calL_{\textrm{-}\Lthetafin}$, with $\Lthetafin\in \coweightlatticefin\simeq \Pic(\rsv)$, and letting $\Lthetafin$ tend to infinity in the strictly dominant chamber yields a sequence of hearts 
\begin{align}
	\bfS_{\textrm{-}\Lthetafin} (\nilpPi)\ ,  
\end{align}
in $\catDb(\catCoh_C(\rsv))$, each equivalent to $\nilpPi$, which converges to $\catCoh(\widehat{\rsv}_C)$.\footnote{A generalization of this result when $\Lthetafin$ is only dominant is presented in \cite{SSS_2025}.} By the main result of Part~\ref{part:COHA-stability-condition}, Theorem~\ref{thmintro:variation} below, we thus get an algebra isomorphism between $\coha_{\rsv,C}^A$ and the limit cohomological Hall algebra of the sequence of hearts $\bfS_{\textrm{-}\Lthetafin} (\nilpPi)$. Note that the COHA of each of these hearts is isomorphic to the negative half $\Y_\qv^-$ of the Yangian of $\frakg_\qv$ via the morphism $\Phi$ introduced in \S\ref{sec:relation-Yangian-COHA}. To understand this limit COHA, we identify the operation of tensoring by a line bundle $\calL_{\textrm{-}\Lthetafin}$ with the action $T_{\textrm{-}\Lthetafin}$ of the lattice element of the affine braid group $B_\sfex$. This, combined with the main result of Part~\ref{part:Yangians}, Theorem~\ref{thmintro:Si-equals-Ti} below, allows us to identify the limit COHA with an adequate limit of $\Y^+_\infty$ quotients of affine Yangians, with transition morphisms given by the action of the algebraic braid group operators.

The precise computation of the image under the isomorphisms $\Theta_{\rsv, C}$ of fundamental classes $[\frakZ_{i,n}]$ and $[\frakY_{i,d}]$ uses the functoriality of nilpotent cohomological Hall algebras, which allows us to eventually reduce the problem to some explicit computations in the case of type $A_1$, i.e., the case of the surface $T^\ast\PP^1$. On the way, we perform several geometric computations which are interesting in their own right, and valid in the general context of a symplectic surface $S$ and a possibily reduced curve $C$ whose irreducible components are smooth and intersect transversally.

\medskip

Let us now describe two auxiliary results of Parts~\ref{part:COHA-stability-condition} and \ref{part:Yangians} which could be of independent interest.

\subsubsection*{Part~\ref{part:COHA-stability-condition}. Variation of $t$-structures and $2$-dimensional Cohomological Hall algebras}

In Part~\ref{part:COHA-stability-condition}, we address the following question: \textit{how are cohomological Hall algebras related when constructed from hearts of different $t$-structures on the same triangulated category?}

The answer to the above question is well-known in the setting of the usual Hall algebras of functions on moduli stacks of objects in a finitary hereditary category: Cramer \cite{Cramer_Hall} proved that even though the Hall algebras of derived equivalent hereditary categories need not be isomorphic, their Drinfeld doubles are. In other words, they may be viewed as two different positive halves of the same quantum group (see also \cite[Lecture~5]{Schiffmann-lectures-Hall-algebras}). Although we believe that something along these lines should also hold for cohomological Hall algebras\footnote{In fact, Theorem~\ref{thmintro:Kleinian_Sing} provides a first nontrivial example of such a principle for cohomological Hall algebras. It also illustrates the need to consider adequate completions.}, the situation is much more subtle. For instance, there is no general construction of a comultiplication, hence no recipe to define a natural \textit{Drinfeld double}.

One can ask a simpler question. Given a sequence of $t$-structures $\tau_n$ on a fixed triangulated category $\scrC$, assume that $\tau_n$ converges in a suitable sense to a $t$-structure $\tau_\infty$. Assume also that one may define cohomological Hall algebras $\coha_{\tau_n}$, for any $n\in \N$, and $\coha_{\tau_\infty}$. Note that the construction of cohomological Hall algebras associated with hearts of $t$-structures on triangulated categories is given in \cite{DPS_Torsion-pairs} and we refer to \textit{loc.cit.} for details. In what sense does $\coha_{\tau_n}$ `converge' to $\coha_{\tau_\infty}$? In Part~\ref{part:COHA-stability-condition}, we develop a framework to answer this question for variations of $t$-structures arising from slicings on a triangulated category $\scrC$. Under suitable assumptions on the moduli stacks of $\tau_n$-objects we define a \textit{limiting cohomological Hall algebra} associated to the sequence $\tau_n$ and we prove a `stabilization' theorem relating such an algebra to $\coha_{\tau_\infty}$.

Let $\calP$ be a slicing on $\scrC$ in the sense of Bridgeland (cf.\ Definition~\ref{def:slicing}). Let $(a_k)$ be a decreasing sequence of real numbers converging to $a_\infty\in (0, 1]$. We define for each $k \in \N \cup \{\infty\}$ a $t$-structure 
\begin{align}
	\tau_k \coloneqq \big( \calP(> a_k - 1 ), \calP(\leqslant a_k) \big) \ , 
\end{align} 
for a pair of integers $k \geqslant \ell$ we write $I_{k,\ell}\coloneqq (a_\ell-1, a_k]$ and we consider the derived stack
\begin{align}
	\dstackCohps(\scrC, I_{\ell,k}) \coloneqq \dstackCohps(\scrC, \tau_k) \cap \dstackCohps(\scrC, \tau_\ell)  
\end{align}
as well as the induced simplicial derived stack
\begin{align}
	\calS_\bullet \dstackCohps\big(\scrC, I_{\ell,k}\big) \ .
\end{align}
When $\ell=k$ (including $\infty$) we simply write $\dstackCohps(\scrC, I_k)$ and $\calS_\bullet\dstackCohps(\scrC, I_k)$. Here, the stack $\dstackCohps(\scrC, \tau_n)$ is the derived moduli stack of pseudo-perfect objects of $\scrC$ which are $\tau_n$-flat for $n\in \N\cup\{\infty\}$. 

\begin{theoremintroduction}[Theorem~\ref{thm:limiting_vs_limit}, Corollary~\ref{cor:limiting_COHA_vs_COHA_limit}]\label{thmintro:variation}
	Let $\scrC,\calP$ and $(a_k)_k$ be as above. Under Assumption~\ref{assumption:limiting_2_Segal_stack}, we have:
	\begin{enumerate}[label=(\roman*)]\itemsep0.2cm
		\item For any $k \leq \ell$, $\calS_\bullet \dstackCohps\big(\scrC, I_{\ell,k}\big)$ satisfies the $2$-Segal condition.
		
		\item The ind-pro derived stack
		\begin{align}
			\calS_\bullet \dstackCohps(\scrC, \tau_\infty^+) \coloneqq \fcolim_{\ell} \flim_{k\geqslant \ell} \calS_\bullet \dstackCohps(\scrC, I_{\ell,k}) 
		\end{align}
		satisfies the $2$-Segal condition.
		
		\item \label{item-introduction:iii} Under also Assumption~\ref{assumption:limiting_CoHA_I}, we get 
		\begin{align}
			\coha_{\tau_\infty^+}\coloneqq \HBM_\bullet\big( \dstackCohps\big(\scrC, \tau_\infty^+ \big)) \coloneqq \flim_{\ell} \colim_{k \geqslant \ell} \HBM_\bullet\big( \dstackCohps\big(\scrC, I_{\ell,k}\big) \big)
		\end{align}
		carries a canonical graded associative algebra structure. We call it the \textit{limiting cohomological Hall algebra} associated to $\tau_n$.
		
		\item \label{item-introduction:iv} Under Assumptions~\ref{assumption:limiting_2_Segal_stack}, \ref{assumption:limiting_CoHA_I}, and \ref{assumption:quasi-compact_interval}, there is an equivalence of ind-pro derived stacks
		\begin{align}
			\calS_\bullet\dstackCohps\big(\scrC, \tau_\infty\big) \longrightarrow \calS_\bullet\dstackCohps\big(\scrC, \tau_\infty^+\big) 
		\end{align}
		which induces an isomorphism of cohomological Hall algebras
		\begin{align}
			\coha_{\tau_\infty} \longrightarrow \coha_{\tau_\infty^+}\ . 
		\end{align}
	\end{enumerate}	
\end{theoremintroduction}
Under the presence of a compatible torus action, the above theorem holds equivariantly as well.

Point~\ref{item-introduction:iv} establishes that the limiting COHA $\coha_{\tau_\infty^+}$ coincides with the COHA $\coha_{\tau_\infty} $ of the limit $t$-structure $\tau_\infty$ under specific conditions: Assumption~\ref{assumption:limiting_2_Segal_stack} ensures the openness of the derived stacks $\dstackCohps(\scrC, \tau_\ell)$ in the derived stack of pseudo-perfect complexes; Assumption~\ref{assumption:limiting_CoHA_I} guarantees the derived lci nature of the map $q$, which is automatic for our choice of $\scrC$, and the properness of $p$. These assumptions are necessary in order to obtain a limiting COHA. Assumption~\ref{assumption:quasi-compact_interval} is the most relevant one here and essentially requires quasi-compactness for the derived stacks $\dstackCohps(\scrC, I_{\ell, k})$ parametrizing objects of a fixed class, whose HN factors belong to the interval $I_{\ell, k}$.

\subsubsection*{Part~\ref{part:Yangians}. COHAs of quivers, reflection functors, Yangians, and braid group actions}

Let $\qv=(I,\Omega)$ be a quiver without edge loops. Let $\Pi_\qv$ denote the associated \textit{preprojective algebra} and let $\Lambda_\qv$ be the stack of nilpotent finite-dimensional $\Pi_\qv$-modules. There is a torus $T$ acting on $\Lambda_\qv$ by rescaling the edges of $\qv$. We consider the twisted $T$-equivariant nilpotent cohomological Hall algebra 
\begin{align}
	\coha_\qv^T\coloneqq \HBMbulletT(\Lambda_\qv)\ , 
\end{align}
which is a $(\Z \times \N I)$-graded associative $\Hbullet_T$-algebra, where $\Hbullet_T$ denotes the equivariant cohomology of a point. Following \cite{NSS_KHA} in the K-theoretical context, we define a multi-parameter Yangian $\Y_\qv$ by generators and relations and prove in Theorem~\ref{thm:Phi} that there exists a surjective algebra morphism 
\begin{align}
	\Phi\colon \Y^-_\qv \longrightarrow \coha^T_\qv	
\end{align}
from the negative half $\Y^-_\qv\subset\Y_\qv$ to $\coha^T_\qv$ which we conjecture to be an isomorphism. This conjecture is proved in the case of affine quivers in Part~\ref{part:ADE-case}. An analogous result for the full COHA, rather than the nilpotent one, is established in \cite{Jindal_COHA} for affine type A quivers. In the K-theoretical context, this conjecture is addressed in \cite{VV_KHA}. 

One of the key features of the bounded derived category $\catDb(\modPi)$ of the abelian category $\modPi$ of finite-dimensional representations of $\Pi_\qv$ is that it carries an action of the braid group $B_\qv$ associated with $\qv$ by derived autoequivalences, generated by the derived \textit{reflection functors} $\R S_i$ for each $i \in I$. This action has important consequences for the representation theory of the associated Kac-Moody Lie algebra $\frakg$ and for the associated quantum groups, see e.g. \cite{BK12}. 

The main result of Part~\ref{part:Yangians} relates this action of $B_\qv$ on $\catDb(\modPi)$ to the standard algebraic action of $B_\qv$ on the Yangian $\Y_\qv$ given by the usual triple exponential
\begin{align}
	T_i\coloneqq\exp\big(\ad\big(x_i^+\big)\big)\circ\exp\big(-\ad\big(x_i^-\big)\big)\circ\exp\big(\ad\big(x_i^+\big)\big)\ .
\end{align}
Just as the functors $\R S_i$ do not preserve the heart of the standard $t$-structure of $\catDb(\modPi)$, the action of $B_\qv$ does not preserve $\Y^-_\qv$. Let $\Lambdai_\qv$ and $\iLambda_\qv$ be the open substacks of finite-dimensional $\Pi_\qv$-representations which are injective and surjective at $i$, respectively. Then, the functor $\R S_i$ induces an equivalence
\begin{align}
	\begin{tikzcd}[ampersand replacement=\&]
		S_{i,\ast}\colon \iLambda_\qv \ar{r}{\sim} \& \Lambdai_\qv
	\end{tikzcd}\ .
\end{align}
Let $(\alpha_i)_{i \in I}$ denote the canonical basis of $\N I$.
\begin{theoremintroduction}[Corollary~\ref{cor:Phi-i},Proposition~\ref{prop:overlineTi}, and Theorem~\ref{thm:compatibility-braid-reflection-functors}]\label{thmintro:Si-equals-Ti} 
	The following holds:
	\begin{enumerate}[label=(\roman*)]\itemsep0.2cm
		\item  The restriction map induces isomorphisms
		\begin{align}
			\begin{tikzcd}[ampersand replacement=\&, row sep=tiny]
				\HBMbulletT(\iLambda_\qv) \ar{r}{\sim} \& \icohaqv \coloneqq \coha^T_\qv \Big/ \Big(\coha^T_\qv\Big)_{\alpha_i} \cdot \coha_\qv^T\\
				\HBMbulletT(\Lambdai_\qv) \ar{r}{\sim}\& \cohai_{\qv} \coloneqq \coha^T_\qv \Big/ \coha^T_\qv \cdot \Big(\coha^T_\qv\Big)_{\alpha_i}
			\end{tikzcd} \ ,
		\end{align}
		where $\Big(\coha^T_\qv\Big)_{\alpha_i}$ is the graded piece of degree $\alpha_i$. 
		
		\item \label{thmintro:Si-equals-Ti-II} The action of the braid operator $T_i \in \Aut(\Y_\qv)$ induces an isomorphism
		\begin{align}
			\begin{tikzcd}[ampersand replacement=\&]
				\overline{T}_i\colon \iY_\qv \coloneqq \Y^-_\qv \Big/ \Y_i \cdot \Y^-_\qv \ar{r}{\sim} \& \Yi_\qv \coloneqq \Y_\qv^-\Big/\Y^-_\qv \cdot \Y_i
			\end{tikzcd}\ .
		\end{align}
		
		\item \label{item-introduction:Si-equals-Ti} The following diagram commutes, up to an explicit sign 
		\begin{align}
			\begin{tikzcd}[ampersand replacement=\&]
				\iY_\qv \arrow{r}{\overline{T}_i}\arrow[swap,twoheadrightarrow]{d}{\Phi} \& \Yi_\qv \arrow[twoheadrightarrow]{d}{\Phi} \\
				\icohaqv \arrow{r}{S_{i,\, \ast}} \& \coha^{T,(i)}_\qv
			\end{tikzcd}\ .
		\end{align}
	\end{enumerate}
\end{theoremintroduction}
The sign involved in the commutativity of the above diagram comes from the twisting of the multiplication in $\coha^T_\qv$.

\subsection*{Further directions}

To finish this introduction, we very briefly gather some of the natural directions of research suggested by the present work and to which we hope to return in the future.

\subsubsection*{A description by generators and relations of $\coha^T_{X,C}$} 

As proved in \cite[Theorem~5.1]{Sala_Schiffmann} for $X=T^\ast \PP^1$ and conjectured for an arbitrary minimal resolution $X$ of a Kleinian singularity, the fundamental classes $[\frakZ_{i,n}]$ and $[\frakY_{i,d}]$ are expected to generate topologically $\coha^T_{X,C}$. Then, Theorem~\ref{thmintro:generators} provides an approach for determining the relations among these generators by using their realization in terms of Yangian generators and the relations between Yangian generators. We plan to compute the relations between these geometric generators in the future.

\subsubsection*{Representations of $\coha^T_{X,Z}$ and a realization of the double} 

The formalism developed in \cite{DPS_Torsion-pairs} for constructing representations of cohomological Hall algebras should extend naturally to the nilpotent COHA case. This generalization would enable the construction of representations of the algebra $\coha^T_{X,Z}$ by considering appropriate moduli stacks of semistable sheaves or complexes of sheaves, including their framed versions when $X$ is an open surface. Such geometric representations of $\coha^T_{X,Z}$ would be essential for defining the right notion of `double' of $\coha^T_{X,Z}$by studying the commutator between \textit{positive} and \textit{negative} operators. In the case of Kleinian resolution of singularities, we anticipate that this double algebra recovers a completion of the Maulik-Okounkov Yangian \cite{MO_Yangian} associated with a \textit{non-dominant} coweight.

This approach to the representation theory of $\coha^T_{X,Z}$ should recover and generalize known constructions of Lie algebra actions on the homology of moduli spaces of Gieseker-stable sheaves on K3 or elliptic surfaces, as studied in \cite{Nakajima_Lie, Yoshioka_Lie, DeHority_K3, DeHority_Elliptic}.

We further hope and anticipate that the action of $\coha^T_{X,Z}$ could help precisely determine the \textit{Ext}-operator on the cohomology of smooth moduli spaces of stable sheaves on smooth surfaces. This application would be particularly valuable when combined with Theorems~\ref{thmintro:Kleinian_Sing} and \ref{thmintro:generators} in the context of Alday-Gaiotto-Tachikawa conjectures for ALE spaces, see \cite{Negut_W-algebras,Negut_AGT_relations_S}.

\subsubsection*{Triple loop Yangians} 

When $X=T^\ast \PP^1$ and $Z=\PP^1$, the algebra $\coha_{X,Z}^T$ is isomorphic to the quantum group $\Y^+_\infty(\widehat{\frakgl}_2)$, which may be thought of as a double loop version of the Drinfeld positive half of the Yangian $\Y(\widehat{\frakgl}_2)$. The case of a Kleinian resolution of type $A_2$ gives rise to two copies of $\Y^+_\infty(\widehat{\frakgl}_2)$ which satisfy commutation relations specified by their respective embeddings in $\Y^+_\infty(\widehat{\frakgl}_3)$.

For $\qv$ a quiver without edge-loops or multiple edges, it is therefore tempting to define a double loop Yangian associated to $\qv$ generated by subalgebras $\Y^+_\infty(\widehat{\frakgl}_2)_i$ isomorphic to $\Y^+_\infty(\widehat{\frakgl}_2)$ for each vertex of $i \in \qv$ and by imposing that $\Y^+_\infty(\widehat{\frakgl}_2)_i$ and $\Y^+_\infty(\widehat{\frakgl}_2)_j$ commute if $i,j$ are not adjacent and satisfy the above type $A_2$ relations if $i,j$ are adjacent in $\qv$. When $\qv$ is taken to be itself an affine quiver, the above quantum group could be taken to be a definition of a positive half of a \textit{triple loop} Yangian. It would be interesting to compare this with the quiver Yangians defined by Li and Yamazaki in the context of the algebra of BPS states for non-compact CY3 in \cite{Li_Yamazaki}.

\subsubsection*{Limiting COHAs} 

It is worth noting that the \textit{limiting COHA} (cf.\ Theorem~\ref{thmintro:variation}) has potential applications extending beyond the scope of this paper. In particular, the theorem becomes valuable in situations where directly defining a two-dimensional COHA associated with the heart of a $t$-structure proves difficult. 

Consider the case where $X$ is a minimal resolution of the singular surface $C^2/G$, for $G \subset \GL(2,\C)$ a finite subgroup. The \textit{derived $\GL(2,\C)$ McKay correspondence} yields an equivalence $\catQCoh(X) \simeq A\Mod$, where $A$ is the endomorphism algebra from Van~der~Bergh's construction \cite{VdB_Flops}. When $G \subset \SL(2, \C)$, this algebra $A$ coincides with the preprojective algebra of the associated McKay quiver. According to \cite[Corollary~1.1]{Wemyss_McKay}, the algebra $A$ has global dimension three when $G \not\subset \SL(2, \C)$, despite being derived equivalent to a surface. In such cases, we anticipate constructing a \textit{limiting} COHA for $A$ in the sense of Theorem~\ref{thmintro:variation}--\ref{item-introduction:iii} and characterizing it.

\subsection*{Structure of the paper}

Part~\ref{part:COHA-stability-condition} consists of four sections. \S\ref{sec:slicing} recall the notions of \textit{slicing} and \textit{Bridgeland pre-stability condition} of a triangulated category. In \S\ref{sec:flat-objects} we introduce the derived moduli stacks of objects of interest and discuss their 2-Segal space structure, which is used in \S\ref{sec:limiting-COHA} to introduce the limiting COHA $\coha_{\tau_\infty^+}$. Furthermore, in \S\ref{sec:stabilization} we compare the limiting COHA with $\coha_{\tau_\infty}$.

Part~\ref{part:Yangians} consists of seven sections and one appendix. In \S\ref{sec:quiver}, we recall the Lie theory associated with a quiver, in particular the notions of Weyl and braid groups. \S\ref{sec:nilpotent-quiver-COHA} concerns the construction of the two-dimensional (nilpotent) COHA of a quiver. In \S\ref{sec:reflection-functors} we recall the theory of (derived) reflection functors which will be used later in the part. In \S\ref{sec:Yangians} we introduce the multi-parameter Yangian associated with an arbitrary quiver, while in \S\ref{sec:braid-group-action-Yangian} we define the action of the braid group of a quiver on the multi-parameter Yangian associated with the same quiver. \S\ref{sec:relation-Yangian-COHA} establishes the relation between the Yangian and the nilpotent COHA associated with a quiver. Finally, in \S\ref{sec:compatibility-braid} we prove the main result of this Part, which is the compatibility between the braid group actions on the Yangian and the nilpotent quiver COHA. Appendix~\ref{sec:sign-twisting} concerns the sign twist appearing in the compatibility result.

Part~\ref{part:ADE-case} is formed by eight sections and three Appendices. In \S\ref{sec:affine-quivers} we recall the structure of the root and weight lattices in the affine quiver case, moreover we recall the notions of affine and extended affine braid groups. \S\ref{sec:Yangians-affine-quivers} concerns the introduction of a two-parameter Yangian associated with an affine quiver: we discuss its classical limit and we compare with the multi-parameter Yangian introduced in \S\ref{sec:Yangians}. In \S\ref{sec:quotients-affine-coha} we introduce quotients of the Yangian and of the affine quiver COHA which will play a crucial role in the description of $\coha^T_{\rsv, C}$ in terms of affine Yangians. \S\ref{sec:perverse-coherent-sheaves} provides a summary of the \textit{derived McKay equivalence} following \cite{VdB_Flops} and a description of the category of (nilpotent) perverse coherent sheaves. In \S\ref{sec:braid-group-derived-categories}, we define a weak action of the extended braid group on the bounded derived category $\catDb(\modPi)$, while in \S\ref{sec:stab-braiding} we discuss the relation between such an action and a slicing induced by a King stability condition on $\catDb(\modPi)$. Our first main result, the characterization of $\coha^T_{\rsv, C}$ in terms of affine Yangians, is given in \S\ref{sec:limit-affine-Yangian}, while in \S\ref{sec:explicit-computations} we provide a description of the (conjectural) geometric generators of $\coha^T_{\rsv, C}$, which are fundamental classes of certain stacks, in terms of Yangian generators. Finally, Appendices~\ref{sec:characterization-classical-limit} and \ref{sec:hopfact} contain proofs of statements given in this part, while Appendix~\ref{sec:purity-formality} discusses some purity results for quotient stacks.

\subsection*{Acknowledgments}

We would like to thank the following people for numerous conversations regarding the subject of this paper: Marian Aprodu, Ben Davison, Drago\v{s} Fr\u{a}til\u{a}, Nicolas Guay, Lucien Hennecart, Benjamin Hennion, David Hernandez, Ryosuke Kodera, Emanuele Macrì, Anton Mellit, Hiraku Nakajima, Andrei Neguţ, Simon Riche, Marco Robalo, Yan Soibelman, Michael Wemyss, and Curtis Wendlandt.

The third-named author would like to thank Andrei Neguţ for very helpful discussions, which took place under the MIT-UNIPI Project (XI call). Moreover, he acknowledges the MIUR Excellence Department Project awarded to the Department of Mathematics, University of Pisa, CUP I57G22000700001. Finally, he is a member of GNSAGA of INDAM.

A part of the paper was finalized during a research visit of the second-named and fourth-named authors at the Department of Mathematics of the University of Pisa under GNSAGA – INDAM research visits program and the 2023 Visiting Fellow program of the University of Pisa, respectively. They thank the University of Pisa for the great working conditions provided during their visit. Another part of the paper was finalized during research visits of the third-named author at Kavli IPMU, the University of Tokyo, in 2023, 2024, and 2025: he thanks Kavli IPMU  for the great working conditions provided during his visits. Finally, the fourth-named author would like to thank the Simion Stoilow Institute of Mathematics for the great working conditions provided during the preparation of this paper.

\numberwithin{section}{part}

\newpage
\part{Variation of $t$-structures and $2$-dimensional Cohomological Hall algebras}\label{part:COHA-stability-condition} 

In this part we introduce the construction of the \textit{limiting COHA}. Given a slicing $\calP$ on a triangulated category $\scrD_0$ and a decreasing sequence of real numbers $\{a_k\}_{k \in \N}$ with $a_0 = 1$ and limit $a_\infty \in (0,1)$, we associate to the abelian category $\calP((a_\infty-1, a_\infty])$ an associative algebra, defined as a suitable \textit{limit} of the cohomological Hall algebras attached to the abelian categories $\calP((a_k-1, a_k])$ for $k\in \N$.

\begin{remarkunnum}
	The construction of the \textit{limiting COHA} can equivalently be formulated in terms of an increasing sequence of real numbers $\{a_k\}_{k \in \N}$ with $a_0 = 0$ and limit $a_\infty \in (0,1)$, together with the abelian categories $\calP([a_k-1, a_k))$ for $k \in \N \cup {\infty}$.
\end{remarkunnum}

\subsubsection*{Notation} 

In this part, we follow the notation introduced in \cite[\S1.6]{Porta_Sala_Hall}, with the difference that we denote with $\Spc$ instead of $\calS$ the $\infty$-category of spaces. In particular, we use the \textit{implicitly derived convention}: given a morphism of derived stacks $f \colon X \to Y$, we let 
\begin{align}
	f^\ast \colon \catQCoh(Y) \to \catQCoh(X)
\end{align}
be the \textit{derived} pullback functor, and we let 
\begin{align}
	f_\ast \colon \catQCoh(X) \to \catQCoh(Y)
\end{align}
be the \textit{derived} pushforward. All fiber products, Hom sheaves and spaces, and tensor products will be understood in the derived sense, unless otherwise stated.

\section{Slicings}\label{sec:slicing}

Let $\scrC$ be a triangulated category.
\begin{definition}[{\cite[Definition~3.3]{B07}}]\label{def:slicing}
	A \textit{slicing $\calP$ of $\scrC$} consists of full additive subcategories $\calP(\phi) \subset \scrC$ for each $\phi\in \R$, such that
	\begin{itemize}\itemsep0.2cm
		\item for all $\phi\in \R$, we have $\calP(\phi + 1) = \calP(\phi)[1]$,
		\item  if $\phi_1 > \phi_2$ and $E_j\in \calP(\phi_j)$ for $j=1,2$, then $\Hom_\scrC(E_1, E_2) = 0$,
		\item (HN filtrations) for every nonzero object $E \in \scrC$ there exists a finite sequence of morphisms
		\begin{align}\label{eq:HN-filtration-slicing}
			\begin{tikzcd}[ampersand replacement=\&]
				0\eqqcolon E_{s+1} \ar{r}{f_s} \& E_s \ar{r} \& \cdots \ar{r}{f_1} \& E_1\coloneqq E
			\end{tikzcd}
		\end{align}
		such that the cofiber of $f_i$ is in $\calP(\phi_i)$ for some sequence $\phi_s > \phi_{s-1} >\cdots > \phi_1$ of real numbers. \qedhere
	\end{itemize}
\end{definition}

\begin{remark}
	Note that the sequence~\eqref{eq:HN-filtration-slicing} is unique up to isomorphism. 
\end{remark}

We write $\phi^+(E) \coloneqq \phi_1$ and $\phi^-(E) \coloneqq \phi_s$. For an interval $I\subset \R$, we write
\begin{align}
	\calP(I)\coloneqq \big\{E\in\scrC\, \vert\,\phi^+(E), \phi^-(E)\in I\big\}= \langle \calP(\phi)\, \vert\, \phi\in I\rangle \subset \scrC\ .
\end{align}

\begin{notation}
	For any $\phi\in \R$, we set
	\begin{align}
		\calP(>\phi)&\coloneqq \calP((\phi, +\infty)) \quad\text{and}\quad \calP(\geqslant \phi)\coloneqq \calP([\phi, +\infty))\ , \\
		\calP(<\phi)&\coloneqq \calP((-\infty, \phi)) \quad\text{and}\quad \calP(\leqslant \phi)\coloneqq \calP((-\infty, \phi])\ . \tag*{\qedhere} 
	\end{align}
\end{notation}

\begin{remark}\label{rem:slicing-t-structures}
	For any $\phi\in\R$, one has pairs of orthogonal subcategories $(\calP(>\phi)), \calP((\leqslant \phi))$ and $(\calP(\geqslant \phi)), \calP(<\phi))$. Note that the subcategories $\calP(>\phi)$ and $\calP(\geq \phi)$ are closed under left shifts and thus define $t$-structures on $\scrC$, where $\calP((\phi, \phi + 1])$ and $\calP([\phi, \phi + 1))$ are the corresponding hearts, respectively.
	
	In the following, we shall denote by $\tau_\phi$ the $t$-structure with heart $\calP((\phi, \phi + 1])$.
\end{remark}

\begin{definition}\label{def:pre-stability-condition}
	Let $\scrC$ be a triangulated category. Fix a finite rank free abelian group $\Lambda$ and a group homomorphism 
	\begin{align}
		v\colon K_0(\scrC) \longrightarrow \Lambda \ .
	\end{align}
	A \textit{pre-stability condition}\footnote{In \cite{B07}, this is the notion of a \textit{stability condition}.} on $\scrC$ with respect to $\Lambda$ is a pair $\sigma = (\Lambda, v, \calP, Z)$ where $\calP$ is a slicing of $\scrC$ and $Z\colon \Lambda\to \C$ is a group homomorphism, called a \textit{central charge}, that satisfy the following condition: for all $0 \neq  E \in  \calP(\phi)$, we have $Z(v(E)) \in \R_{>0}\cdot e^{\imath\pi\phi}$.
	
	We will often abuse notation and write $Z(E)$ for $Z(v(E))$. The nonzero objects of $\calP(\phi)$ are called \textit{$\sigma$-semistable of phase $\phi$}.
\end{definition}

\section{2-Segal derived stack of flat objects}\label{sec:flat-objects}

In this section, we shall introduce a suitable derived stack of flat objects and its simplicial version, which would inherit the structure of a 2-Segal derived stack.

\subsection{Recollection of $2$-Segal stacks}

We refer the reader to \cite[\S4.1]{Porta_Sala_Hall} and \cite[\S~I.11]{DPS_Torsion-pairs} for background material on the notion of $2$-Segal object introduced by Dyckerhoff-Kapranov in \cite{Dyckerhoff_Kapranov_Higher_Segal}.

Let $k$ be an algebraically closed field of characteristic zero. We start by reviewing the \textit{Wadhausen construction}. Given an integer $n$ we let $[n]$ denote the linearly ordered poset $\{0 < 1 < \cdots < n\}$ and we set
\begin{align}
	\sfT_n \coloneqq \Fun([1],[n]) \ . 
\end{align}
The collection of the various $\sfT_n$ determines a functor
\begin{align}
	\sfT_\bullet \colon \bfDelta \longrightarrow \Cat_\infty \ . 
\end{align}
Given a stable $\infty$-category $\calC$, we set
\begin{align}
	\calS_n \subseteq \Fun(\sfT_n, \calC) 
\end{align}
be the full subcategory spanned by those functors $F \colon \sfT_n \to \calC$ satisfying the following two conditions:
\begin{enumerate}\itemsep=0.2cm
	\item $F(i,i) \simeq 0$;
	\item for every $0 \leqslant i < j \leqslant n-1$, the square
	\begin{align}
		\begin{tikzcd}[ampersand replacement=\&]
			F(i,j) \arrow{r} \arrow{d} \& F(i+1, j) \arrow{d} \\
			F(i,j+1) \arrow{r} \& F(i+1, j+1)
		\end{tikzcd}
	\end{align}
	is a pullback in $\calC$.
\end{enumerate}
The $\infty$-categories $\calS_n \calC$ depends simplicially on $n$, i.e., they assemble into a simplicial object
\begin{align}
	\calS_\bullet \calC \colon \bfDelta\op \longrightarrow \Cat_\infty \ , 
\end{align}
known as the \textit{Waldhausen construction} of $\calC$. One of the main observations of \cite{Dyckerhoff_Kapranov_Higher_Segal} is that this satisfies the $2$-Segal condition. The functoriality of this construction allow to replace $\calC$ by a functor
\begin{align}
	\scrC \colon \dAff_k\op \longrightarrow \Cat_\infty^{\mathsf{st}} \ , 
\end{align}
giving rise to a $2$-Segal object in derived prestacks
\begin{align}
	\calS_\bullet \scrC \colon \bfDelta\op \longrightarrow \PreSt_k \ , 
\end{align}
which we refer to as the Waldhausen construction of $\scrC$. Besides, this construction is obviously functorial in $\scrC$.

\begin{notation}
	Given a stable $\infty$-category $\calC$, we refer to $F \in \calS_n \calC$ as an \textit{$n$-flag of objects in $\calC$}.
	Given $0 \leqslant i \leqslant j \leqslant n$, we set
	\begin{align}\label{eq:ev}
		\ev_{i,j}(F) \coloneqq F(i,j) \ , 
	\end{align}
	We also denote by
	\begin{align}
		\partial_i(F) \in \calS_{n-1} \calC 
	\end{align}
	the $(n-1)$-flag obtained restricting $F$ along the morphism $[n-1] \to [n]$ in $\bfDelta$ that misses $i$.
	In particular, when $n = 2$ we have the overlap of notation
	\begin{align}\label{eq:partial}
		\partial_0 = \ev_{1,2} \ , \qquad \partial_1 = \ev_{0,2} \ , \qquad \partial_2 = \ev_{0,1} \ . 
	\end{align}
\end{notation}

\begin{remark}
	Let $(\Lambda,+)$ be an abelian group. The group structure can be encoded into the structure of a simplicial set satisfying the ($1$- and hence) $2$-Segal condition, that we denote
	\begin{align}
		\Lambda^\bullet \colon \bfDelta\op \longrightarrow \mathbf{Set} \ . 
	\end{align}
	We write $\pi \colon \bfLambda\op \to \bfDelta\op$ for the associated cocartesian fibration. 
	
	We can modify the above construction by replacing $\bfDelta\op$ with $\bfLambda\op$ and accordingly define a natural $\Lambda$-graded variant of the $2$-Segal condition (cf.\ \cite[\S\ref*{torsion-pairs-sec:Lambda_graded}]{DPS_Torsion-pairs}).
\end{remark}

Let $\scrD \in \PrLomega_k$ be a compactly generated $k$-linear $\infty$-category of finite type. Given a derived affine scheme $S \in \dAff_k$, we write
\begin{align}
	\scrD_S \coloneqq \catQCoh(S) \otimes_k \scrD \ , 
\end{align}
and we write $\scrD_S^{\ps}$ for the full subcategory of $\scrD_S$ spanned by $S$-pseudo-perfect objects. The functor
\begin{align}\label{eq:stack-categorical}
	\calM^{\mathsf{cat}}_\scrD \colon \dAff_k\op \longrightarrow \Cat_\infty^{\mathsf{st}} 
\end{align}
that sends $S$ to $\scrD_S^\ps$ is a étale hypersheaf, and Toën-Vaquié \cite{Toen_Vaquie_Moduli} \textit{moduli of objects} $\calM_\scrD$ is by definition given by the rule
\begin{align}
	\calM_\scrD(S) \coloneqq \big( \calM^{\mathsf{cat}}_\scrD(S) \big)^\simeq \ , 
\end{align}
where for an $\infty$-category $\scrC$, we denote by $\scrC^\simeq$ its maximal $\infty$-groupoid. Applying the Waldhausen construction to \eqref{eq:stack-categorical} and passing to the maximal $\infty$-groupoid at the end, we obtain a $2$-Segal object
\begin{align}
	\calS_\bullet \calM_{\scrD} \colon \bfDelta \longrightarrow \PreSt_\C \ .
\end{align}
When $n = 1$, there is a canonical identification
\begin{align}
	\calS_1 \calM_{\scrD}\simeq \calM_{\scrD} \ . 
\end{align}

\begin{notation}
	Given a full categorical substack $\calM_0^{\mathsf{cat}} \subseteq \calM^{\mathsf{cat}}_\scrD$, we write
	\begin{align}
		\scrD_0 \coloneqq \calM_0(\Spec(k)) \subseteq \scrD^\ps \ . 
	\end{align}
	We also write
	\begin{align}
		\calM_0 \coloneqq \big( \calM_0^{\mathsf{cat}} \big)^\simeq 
	\end{align}
	for the associated derived stack.
\end{notation}

\begin{definition}
	A full categorical substack $\calM_0^{\mathsf{cat}}$ of $\calM_\scrD^{\mathsf{cat}}$ is said to be \textit{admissible} if:
	\begin{enumerate}\itemsep=0.2cm
		\item the associated derived stack $\calM_0$ is an admissible indgeometric derived stack in the sense of \cite[Definition~\ref*{foundation-def:admissible_indgeometric_stack}]{DPSSV-1}\footnote{Roughly speaking, $\calM_0$ is a derived stack, which is essentially characterized by the property that their underlying reduced stack is geometric.},
		
		\item the map $\calM_0 \to \calM_\scrD$ is formally étale. \qedhere
	\end{enumerate}
\end{definition}

Following \cite{DPS_Torsion-pairs}, we define a simplicial object $\calS_\bullet \calM_0$ by declaring that for every $[n] \in \bfDelta$ the square
\begin{align}
	\begin{tikzcd}[ampersand replacement=\&]
		\calS_n \calM_0 \arrow{r} \arrow{d} \& \calS_n \calM_\scrD \arrow{d}{(\ev_{i,j})} \\
		\displaystyle\prod_{0 \leqslant i < j \leqslant n} \calM_0 \arrow{r} \& \displaystyle\prod_{0 \leqslant i < j \leqslant n} \calM_\scrD
	\end{tikzcd}
\end{align}
is a pullback, where the maps $\ev_{i,j}$ are introduced in Formula~\eqref{eq:ev}. 

\begin{definition}
	A full categorical substack $\calM_0^{\mathsf{cat}}$ of $\calM_{\scrD}^{\mathsf{cat}}$ is said to be \textit{closed under extensions} if the square
	\begin{align}
		\begin{tikzcd}[ampersand replacement=\&]
			\calS_2 \calM_0 \arrow{r} \arrow{d}{\partial_0 \times \partial_2} \& \calS_2 \calM_\scrD \arrow{d}{\partial_0 \times \partial_2}  \arrow{d} \\
			\calM_0 \times \calM_0 \arrow{r} \& \calM_\scrD \times \calM_\scrD
		\end{tikzcd}
	\end{align}
	is a pullback, where $\partial_i$ is defined in Formula~\eqref{eq:partial} for $i=0, 1, 2$.
\end{definition}

Now, we introduce a $\Lambda$-graded $2$-Segal object structure on $\calS_\bullet \calM_0$.

Let $\calM_0^{\mathsf{cat}}$ be a categorical substack of $\calM_{\scrD}^{\mathsf{cat}}$ closed under extensions. We set $\scrD_0 \coloneqq \calM_0^{\mathsf{cat}}(\Spec(k))$ and we fix a morphism of abelian groups
\begin{align}
	v \colon K_0(\scrD_0) \longrightarrow \Lambda \ . 
\end{align}
Given $\bfv \in \Lambda$, we let $\calM_0(\bfv)$ be the substack of $\calM_0$ defined by the following property: given $S \in \dAff_k$, an $S$-point $x \colon S \to \calM_0$ belongs to $\calM_0(\bfv)$ if and only if for every $k$-point $y \colon \Spec(k) \to S$ the image in $\Lambda$ of the composite $x \circ y \colon \Spec(k) \to \calM_0$ via $v$ coincides with $\bfv$.

\medskip

More generally, given a positive integer $n \geqslant 0$ and
\begin{align}
	\underline{\bfv} = (\bfv_{0,1}, \bfv_{1,2}, \ldots, \bfv_{n-1,n}) \in \Lambda^n \ , 
\end{align}
we define $\calS_n \calM_0(\underline{\bfv})$ as the fiber product
\begin{align}
	\begin{tikzcd}[ampersand replacement=\&]
		\calS_n \calM_0(\underline{\bfv}) \arrow{r} \arrow{d} \& \calS_n \calM_0 \arrow{d}{(\ev_{0,1}, \ev_{1,2}, \ldots, \ev_{n-1,n})} \\
		\displaystyle\prod_{i = 0}^{n-1} \calM_0(\bfv_{i,i+1}) \arrow{r} \& \displaystyle \prod_{i = 0}^{n-1} \calM_0 
	\end{tikzcd} \ .
\end{align}
Concretely, $\calS_n \calM_0(\underline{\bfv})$ parametrizes families of flags $\{M_{i,j}\}_{0 \leqslant i<j \leqslant n}$ such that
\begin{align}
	v(M_{i,j}) = \sum_{i \leqslant k < j} \bfv_{k,k+1} \in \Lambda \ . 
\end{align}
Since $v$ is additive, it is straightforward to check that the assignment 
\begin{align}
	([n],\underline{\bfv}) \longmapsto \calS_n \calM_0(\underline{\bfv})
\end{align}
forms a $\Lambda$-graded $2$-Segal object with values in $\dSt_k$. In particular, \cite[Theorem~I.1.1]{DPS_Torsion-pairs} shows that $\calS_\bullet \calM_0$ determines an $\E_1$-monoid in $\Corr^{\times_\Lambda}(\Fun(\Lambda, \dSt_k))$, the $(\infty, 2)$-category of \textit{correspondences} with values in $\Fun(\Lambda, \dSt_k)$ with monoidal structure induced by $\Lambda$ (for details, see \cite{DPS_Torsion-pairs}). In what follows, we tacitly consider $\calS_\bullet \calM_0$ as a $\Lambda$-graded $2$-Segal derived stack.

\subsection{Derived stack of flat objects}

Let $\calM_0^{\mathsf{cat}}$ be a categorical substack of $\calM_\scrD^{\mathsf{cat}}$ and let $\tau$ be a $t$-structure on $\scrD_0 \coloneqq \calM_0^{\mathsf{cat}}(\Spec(k))$. 
\begin{construction}\label{construction:derived-stack-D_0}
	We let $\dstackCoh(\scrD_0,\tau)$ be the full substack of $\calM_0$ defined as follows: given $S \in \dAff_k$, an $S$-point $x \colon S \to \calM_0$ belongs to $\dstackCoh(\scrD_0,\tau)$ if and only if for every $k$-point $y \colon \Spec(k) \to S$ the composite $x \circ y \colon \Spec(k) \to \calM_0$ is $\tau$-flat. 	
	Given a morphism of abelian groups
	\begin{align}
		v \colon K_0(\scrD_0) \longrightarrow \Lambda \ , 
	\end{align}
	a positive integer $n \geqslant 0$ and $\underline{\bfv} \in \Lambda^n$, we write
	\begin{align}
		\calS_n \dstackCoh(\scrD_0,\tau;\underline{\bfv}) \coloneqq \calS_n \dstackCoh(\scrD_0,\tau) \cap \calS_n \calM_0(\underline{\bfv}) \ . \tag*{\qedhere} 
	\end{align}
\end{construction}
This gives rise to a $\Lambda$-graded $2$-Segal derived stack $\calS_\bullet \dstackCoh(\scrD_0,\tau)$.

\begin{definition}\label{def:open}
	Let $\calM_0^{\mathsf{cat}}$ be a full categorical substack of $\calM_\scrD^{\mathsf{cat}}$.
	\begin{enumerate}\itemsep=0.2cm
		\item We say that a $t$-structure $\tau$ on $\scrD_0$ is \textit{open} if the canonical map
		\begin{align}
			\dstackCohps(\scrD_0,\tau) \longrightarrow \calM_0 
		\end{align}
		is representable by open Zariski immersions.
		
		\item We say that a pair $(\Lambda, v \colon K_0(\scrD_0) \to \Lambda)$ is \textit{open} if for every $\bfv \in \Lambda$, the canonical map
		\begin{align}
			\calM_0(\bfv) \longrightarrow \calM_0 
		\end{align}
		is representable by open Zariski immersions.
		
		\item We say that a triple $(\tau,\Lambda,v)$ is \textit{open} if both $\tau$ and $(\Lambda,v)$ are open, and moreover the derived stack $\dstackCohps(\scrD_0,\tau;\bfv)$ is both open and closed inside $\dstackCohps(\scrD_0,\tau)$. \qedhere
	\end{enumerate}
\end{definition}

\section{The limiting $2$-Segal stack}\label{sec:limiting-COHA}

Fix a full categorical substack $\calM_0^{\mathsf{cat}}$ of $\calM_\scrD^{\mathsf{cat}}$. We also fix a slicing $\calP$ on $\scrD_0$ in the sense of Definition~\ref{def:slicing} and an open pair $(\Lambda, v)$ in the sense of Definition~\ref{def:open}. Finally, we introduce a decreasing sequence of real numbers $\{a_k\}_{k \in \N}$ starting at $a_0 = 1$ and converging at $a_\infty \in (0,1)$. Following Remark~\ref{rem:slicing-t-structures}, for $k \in \N \cup \{\infty\}$, we denote by $\tau_k$ the $t$-structure on $\scrD_0$ defined by
\begin{align}
	\tau_k \coloneqq \big( \calP(> a_k - 1 ), \calP(\leqslant a_k) \big) \ , 
\end{align}
and we make the following assumption:
\begin{assumption}\label{assumption:limiting_2_Segal_stack}
	The triples $\{(\tau_k,\Lambda,v)\}_{k \in \N}$ are open.
\end{assumption}

Under this assumption, we have
\begin{align}
	\dstackCohps(\scrD_0,\tau) = \coprod_{\bfv \in \Lambda} \dstackCohps(\scrD_0,\tau;\bfv) \ . 
\end{align}
For this reason, we commit an abuse of notation and consider $\dstackCohps(\scrD_0,\tau)$ as a $\Lambda$-graded $2$-Segal derived stack.

\begin{notation}
	Given two positive integers $\ell \leqslant k \leqslant \infty$, we denote by $I_{\ell,k}$ the interval $(a_\ell-1, a_k]$. Our assumptions on the sequence $\{a_k\}$ guarantee that $I_{\ell,k}$ has length less than one. When $\ell = k$, we simply write $I_k \coloneqq (a_k-1,a_k]$ instead of $I_{k,k}$.
\end{notation}

\begin{construction}\label{construction:limiting_COHA-1}
	Let $k \geqslant \ell$ be two positive integers, we write
	\begin{align}\label{eq:derived-stack-sequence}
		\dstackCohps(\scrD_0, I_{\ell,k}) \coloneqq \dstackCohps(\scrD_0,\tau_k) \cap \dstackCohps(\scrD_0,\tau_\ell) \ . 
	\end{align}
	Notice that Assumption~\ref{assumption:limiting_2_Segal_stack} automatically implies that this is an open substack of $\calM_0$, and even of $\dstackCohps(\scrD_0,\tau_0)$. We let
	\begin{align}
		\calS_\bullet \dstackCohps\big(\scrD_0, I_{\ell,k}\big) 
	\end{align}
	be the induced simplicial derived stack, which we further consider over $\calS_\bullet \bfLambda$. Notice that for $[n] \in \bfDelta$ and $\underline{\bfv} \in \Lambda^n$ and $k' \geqslant k$, there are canonical open immersions
	\begin{align}
		\calS_n \dstackCohps\big( \scrD_0, I_{\ell,k'}; \underline{\bfv} \big) \longrightarrow \calS_n \dstackCohps\big( \scrD_0, I_{\ell,k}; \underline{\bfv} \big) \ . 
	\end{align}
	We set
	\begin{align}
		\calS_n \dstackCohps\big( \scrD_0, I_{\ell,\infty}^+; \underline{\bfv} \big) \coloneqq \flim_{k \geqslant \ell} \calS_n \dstackCohps\big( \scrD_0, I_{\ell,k} ; \underline{\bfv} \big) \in \Pro(\dSt_k) \ . 
	\end{align}
	This allows to consider $\calS_\bullet \dstackCohps\big( \scrD_0, I_{\ell,\infty}^+\big)$ as a $\Lambda$-graded simplicial object.
\end{construction}

\begin{lemma}\label{lem:limiting_2_Segal_I}
	In the setting of the above construction, for any $\ell \in \N$, the $\Lambda$-graded simplicial pro-derived stack $\calS_\bullet \dstackCohps\big( \scrD_0, I_{\ell,\infty}^+ )$ satisfies the $\Lambda$-graded $2$-Segal condition.
\end{lemma}

\begin{proof}
	To begin with, observe that each $\calS_\bullet \dstackCohps\big( \scrD_0, I_{\ell,k} \big) \in \Fun(\bfDelta\op, \dSt_k)$ satisfies the $\Lambda$-graded $2$-Segal condition. At this point, the conclusion follows from the fact the $2$-Segal condition is closed under limits and that whenever $\calE$ is an $\infty$-category with finite limits, the inclusion
	\begin{align}
		\calE \longrightarrow \Pro(\calE) 
	\end{align}
	preserves them.
	\end{proof}

\begin{construction}\label{construction:limiting_COHA-2}
	For $\ell' \geqslant \ell$, observe that there is a canonical map
	\begin{align}
		\calS_\bullet \dstackCohps\big(\scrD_0, I_{\ell,\infty}^+ \big) \longrightarrow \calS_\bullet \dstackCohps\big(\scrD_0, I_{\ell',\infty}^+ \big) \ . 
	\end{align}
	We set
	\begin{align}
		\calS_\bullet  \dstackCohps(\scrD_0, \tau_\infty^+) \coloneqq \fcolim_{\ell} \calS_\bullet \dstackCohps(\scrD_0, I_{\ell,\infty}^+) \in \Fun\big( \bfDelta\op, \Ind(\Pro(\dSt_k)) \big) \ . 
	\end{align}
	As before, we can perform this construction at the $\Lambda$-graded level, and we can therefore consider $\calS_\bullet  \dstackCohps(\scrD_0, \tau_\infty^+)$ as a $\Lambda$-graded simplicial object in $\Ind(\Pro(\dSt_k))$.
\end{construction}

\begin{lemma}\label{lem:limiting_2_Segal_II}
	The $\Lambda$-graded simplicial ind-pro-derived stack $\calS_\bullet \dstackCohps\big( \scrD_0, \tau_\infty^+\big)$ satisfies the $\Lambda$-graded $2$-Segal condition.
\end{lemma}

\begin{proof}
	We know from Lemma~\ref{lem:limiting_2_Segal_I} that each $\calS_\bullet \dstackCohps\big( \scrD_0, I_{\ell,\infty}^+\big)$ satisfies the $\Lambda$-graded $2$-Segal condition. At this point, the conclusion follows from the fact that the $2$-Segal condition is closed under limits and that whenever $\calE$ is an $\infty$-category with finite limits, the inclusion
	\begin{align}
		\calE \longrightarrow \Ind(\calE) 
	\end{align}
	preserves them.
	\end{proof}

\begin{definition}
	We refer to $\calS_\bullet \dstackCohps\big( \scrD_0, \tau_\infty^+ \big)$ as the \textit{limiting $\Lambda$-graded $2$-Segal derived stack} associated to the sequence $\{a_k\}_k$.
\end{definition}

In order to attach a cohomological Hall algebra to $\dstackCohps\big( \scrD_0, \tau_\infty^+ \big)$, we need some extra assumptions.

\begin{notation}
	For positive integers $\ell_1 \leqslant \ell_2 \leqslant k$ and $\underline{\bfv} = (\bfv_{0,1}, \bfv_{1,2})$, we write
	\begin{align}
		\calS_2 \dstackCohps\big( \scrD_0, I_{\ell_2,k} \times I_{\ell_1,k} \times I_{\ell_2,k}; \underline{\bfv} \big) 
	\end{align}
	for the open substack of $\calS_2\dstackCohps(\scrD_0,\tau_k;\underline{\bfv})$ parametrizing extensions $F_{0,1} \to F_{0,2} \to F_{1,2}$ such that $F_{0,1}$ and $F_{1,2}$ are $\tau_{\ell_2}$ and $\tau_k$-flat, and $F_{0,2}$ is $\tau_{\ell_1}$ and $\tau_k$-flat.
\end{notation}

\begin{assumption}\label{assumption:limiting_CoHA_I}
	\hfill
	\begin{enumerate}\itemsep=0.2cm
		\item \label{assumption:limiting_CoHA_I-1} For any $k\in N$, the map 
		\begin{align}
			\partial_0 \times \partial_2 \colon \calS_2 \dstackCohps\big( \scrD_0, \tau_k \big) \longrightarrow \dstackCohps\big( \scrD_0, \tau_k \big) \times \dstackCohps\big( \scrD_0, \tau_k \big)
		\end{align}
		is derived lci.
		
		\item \label{assumption:limiting_CoHA_I-2} For any $\ell_1\in \N$ and for every pair $\underline{\bfv} = (\bfv_{0,1},\bfv_{1,2}) \in \Lambda^2$, there exists an integer $N \geqslant \ell_1$ such that for every $k' \geqslant k \geqslant \ell_2 \geqslant N$, the square
		\begin{align}
			\begin{tikzcd}[ampersand replacement=\&]
				\calS_2 \dstackCohps\big( \scrD_0, I_{\ell_2,k'} \times I_{\ell_1,k'} \times I_{\ell_2,k'} ; \underline{\bfv} \big) \arrow{r}{\partial_1} \arrow{d} \& \dstackCohps\big( \scrD_0, I_{\ell_1,k'}; \bfv_{0,1} + \bfv_{1,2}\big) \arrow{d} \\
				\calS_2 \dstackCohps\big( \scrD_0, I_{\ell_2,k} \times I_{\ell_1,k} \times I_{\ell_2,k} ; \underline{\bfv} \big) \arrow{r}{\partial_1} \& \dstackCohps\big( \scrD_0, I_{\ell_1,k}; \bfv_{0,1} + \bfv_{1,2}\big)
			\end{tikzcd} 
		\end{align}
		is a pullback and the horizontal arrows are representable by proper algebraic spaces.
			\end{enumerate}
\end{assumption}

\begin{remark}
	Note that if for  any $k \in \N$, the smooth Serre functor\footnote{See e.g. \cite{DPS_Torsion-pairs} for its definition.} $\mathsf S^!_\scrD[2]$ of $\scrD$ restricts to $\scrD_0$ and it is $\tau_k$-exact for $k\in \N$, then Assumption~\ref{assumption:limiting_CoHA_I}--\eqref{assumption:limiting_CoHA_I-1} holds.
\end{remark}

\begin{proposition}\label{prop:limiting_CoHA}
	Let $\bfD$ be a motivic formalism, and let $\calA \in \CAlg(\bfD(\Spec(k)))$ and let $\Gamma \subseteq \Pic(\bfD^\ast(\Spec(k)))$ be an abelian subgroup. Assume that $\calA$ is oriented and that $\Gamma$ is closed under Thom twists. Under Assumptions~\ref{assumption:limiting_2_Segal_stack} and \ref{assumption:limiting_CoHA_I},
	\begin{align}
		\HBMDGamma_0\big( \dstackCohps\big( \scrD_0, \tau_\infty^+ \big); \calA \big) \coloneqq \bigoplus_{\bfv \in \Lambda} \flim_{\ell} \colim_{k \geqslant \ell} \HBMDGamma_0\big( \dstackCohps\big( \scrD_0, I_{\ell,k}; \bfv \big) \big) \in \Fun(\Lambda, \Pro(\catMod_k^\heartsuit))
	\end{align}
	carries a canonical $\Lambda$-graded algebra structure.
\end{proposition}

\begin{proof}
	In virtue of \cite[Theorem~\ref*{foundation-thm:BM_Lambda_graded_functoriality}]{DPSSV-1} and of Lemma~\ref{lem:limiting_2_Segal_II}, it suffices to argue that Assumption~\ref{assumption:limiting_CoHA_I} ensures that for every pair $\underline{\bfv} = (\bfv_{0,1}, \bfv_{1,2}) \in \Lambda^2$, writing $\bfv \coloneqq \bfv_{0,1} + \bfv_{1,2}$, the correspondence
	\begin{align}
		\begin{tikzcd}[column sep=small, ampersand replacement=\&]
			\&[-25pt] \calS_2 \dstackCohps\big( \scrD_0, \tau_\infty^+; \underline{\bfv} \big) \arrow{dl}[swap]{\partial_0 \times \partial_2} \arrow{dr}{\partial_1} \\
			\dstackCohps\big( \scrD_0, \tau_\infty^+; \bfv_{1,2} \big) \times \dstackCohps\big( \scrD_0, \tau_\infty^+; \bfv_{0,1} \big) \& \& \dstackCohps\big(\scrD_0, \tau_\infty^+; \bfv \big)
		\end{tikzcd}
	\end{align}
	For the left diagonal morphism, it suffices to observe that for fixed positive integers $k \geqslant \ell$, Assumption~\ref{assumption:limiting_CoHA_I}-(1) guarantees that the canonical map
	\begin{align}
		\partial_0 \times \partial_2 \colon \calS_2 \dstackCohps\big( \scrD_0, I_{\ell,k}; \underline{\bfv} \big) \longrightarrow \dstackCohps\big( \scrD_0, I_{\ell,k}; \bfv_{1,2} \big) \times \dstackCohps\big( \scrD_0, I_{\ell,k}; \bfv_{0,1} \big) 
	\end{align}
	is derived lci, and Assumption~\ref{assumption:limiting_2_Segal_stack} guarantees that the transition maps as $k$ and $\ell$ vary are open immersions.
	We can therefore apply $\HBMDGamma_0(-;\calA)$ and pass to the colimit in $k$ and the limit in $\ell$ to obtain a well defined morphism
	\begin{multline}
		(\partial_0 \times \partial_2)^! \colon \bigoplus_{\bfv_{0,1} + \bfv_{1,2} = \bfv} \HBMDGamma_0\big( \dstackCohps\big( \scrD_0, I_{\infty}^+;\bfv_{1,2} \big); \calA \big) \otimes \HBMDGamma_0\big( \dstackCohps\big( \scrD_0, I_{\infty}^+; \bfv_{0,1} \big); \calA \big) \\
		 \longrightarrow \bigoplus_{\bfv_{0,1} + \bfv_{1,2} = \bfv} \HBMDGamma_0\big( \calS_2 \dstackCohps\big( \scrD_0, I_{\infty}^+; \underline{\bfv} \big); \calA \big) \ .
	\end{multline}
	As for $\partial_1$, we first observe that a simple cofinality argument yields
	\begin{align}
		\fcolim_{\ell} \calS_2\dstackCohps\big( \scrD_0, I_{\ell,\infty}^+; \underline{\bfv} \big) \simeq \fcolim_{\ell_1} \fcolim_{\ell_2 \geqslant \ell_1} \calS_2\dstackCohps\big( \scrD_0, I_{\ell_2,\infty}^+ \times I_{\ell_1,\infty}^+ \times I_{\ell_2,\infty}^+; \underline{\bfv} \big) \ . 
	\end{align}
	Fix $\ell_1$. Assumption~\ref{assumption:limiting_CoHA_I}--\eqref{assumption:limiting_CoHA_I-2} guarantees that there exists an integer $N$ such that for every $k \geqslant \ell_2 \geqslant N$ the map
	\begin{align}
		\calS_2 \dstackCohps\big( \scrD_0, I_{\ell_2,k} \times I_{\ell_1,k} \times I_{\ell_2,k}; \underline{\bfv} \big) \longrightarrow \dstackCohps\big( \scrD_0, I_{\ell_1,k}; \bfv \big) 
	\end{align}
	is proper, and therefore it induces a well defined morphism
	\begin{align}
		\partial_{1,\ast} \colon \HBMGamma_0\big( \calS_2 \dstackCohps\big( \scrD_0, I_{\ell_2,k} \times I_{\ell_1,k} \times I_{\ell_2,k}; \underline{\bfv} \big); \calA \big) \longrightarrow \HBMGamma_0\big( \dstackCohps\big( \scrD_0, I_{\ell_1,k}; \bfv \big) ; \calA \big) \ . 
	\end{align}
	Since the square obtained varying $k$ is a pullback by assumption, we can pass to the colimit in $k$ to obtain a well defined morphism
	\begin{align}
		\partial_{1,\ast} \colon \HBMGamma_0\big( \calS_2 \dstackCohps\big( \scrD_0, I_{\ell_2,\infty}^+ \times I_{\ell_1,\infty}^+ \times I_{\ell_2,k}^+; \underline{\bfv} \big); \calA \big) \longrightarrow \HBMGamma_0\big( \dstackCohps\big( \scrD_0, I_{\ell_1,\infty}^+; \bfv \big); \calA \big) \ . 
	\end{align}
	By composition with the natural open restriction map 
	\begin{align}
		\HBMGamma_0\big( \calS_2 \dstackCohps\big( \scrD_0, I_{\infty}^+ \times I_{\ell_1,\infty}^+ \times \tau_\infty^+; \underline{\bfv} \big) ; \calA \big) \longrightarrow \HBMGamma_0\big( \calS_2 \dstackCohps\big( \scrD_0, I_{\ell_2,\infty}^+ \times I_{\ell_1,\infty}^+ \times I_{\ell_2,\infty}^+; \underline{\bfv} \big) ; \calA \big) 
	\end{align}
	we also obtain a well defined morphism
	\begin{align}
		 \partial_{1,\ast} \colon \HBMGamma_0\big( \calS_2 \dstackCohps\big( \scrD_0, I_{\infty}^+ \times I_{\ell_1,\infty}^+ \times \tau_\infty^+; \underline{\bfv} \big) ; \calA \big) \longrightarrow \HBMGamma_0\big( \dstackCohps\big( \scrD_0, I_{\ell_1,\infty}^+ \big); \calA\big) \ . 
	\end{align}
	At last, for $\ell_1' \geqslant \ell_1$, the square
	\begin{align}
		\begin{tikzcd}[ampersand replacement=\&]
			\calS_2 \dstackCohps\big( \scrD_0, I_{\infty}^+ \times I_{\ell_1,\infty}^+ \times \tau_\infty^+; \underline{\bfv} \big) \arrow{r} \arrow{d} \& \calS_2 \dstackCohps\big( \scrD_0, I_{\infty}^+ \times I_{\ell_1,\infty}^+ \times \tau_\infty^+; \underline{\bfv} \big) \arrow{d} \\
			\dstackCohps\big( \scrD_0, I_{\ell_1',\infty}^+; \underline{\bfv} \big) \arrow{r} \& \dstackCohps\big( \scrD_0, I_{\ell_1',\infty}^+; \underline{\bfv} \big)
		\end{tikzcd}
	\end{align}
	is obviously a pullback. This allows to pass to the (formal) limit in $\ell_1$, and to define
	\begin{align}
		\partial_{1,\ast} \colon \bigoplus_{\bfv_{0,1} + \bfv_{1,2} = \bfv} \HBMGamma_0\big(\calS_2 \dstackCohps\big(\scrD_0,\tau_\infty^+;\underline{\bfv}\big);\calA \big) \longrightarrow \HBMGamma_0\big(\dstackCohps(\scrD_0,\tau_\infty^+;\bfv\big); \calA \big) \ . \tag*{\qedhere} 
	\end{align}
\end{proof}

\begin{remark}\label{rem:limiting_CoHA}
	When $\bfD^\ast \coloneqq \BettiD$ is the topological formalism of \cite[\S\ref*{torsion-pairs-ex:topological_formalism}]{DPS_Torsion-pairs}, $\calA \coloneqq \Q$, and $\Gamma \coloneqq \Z\langle 1/2 \rangle$ (see Remark~\ref*{torsion-pairs-rem:explicit_bigrading} of \textit{loc.cit.}), we simply write
	\begin{align}
		\coha_{\scrD_0, \tau_\infty^+} 
	\end{align}
	for the cohomological Hall algebra $\HBMDGamma_0\big( \dstackCohps\big( \scrD_0, \tau_\infty^+ \big); \calA \big)$. Moreover, instead of seeing it as a pro-ind-($\Lambda$-graded)-vector space, we see it as a $\Lambda$-graded vector space endowed with the quasi-compact topology\footnote{In the sense of \cite[Remark~\ref*{foundation-rem:quasi-compact_topology}]{DPSSV-1}.}. We shall call it the \textit{limiting cohomological Hall algebra}.
	
	If by varying of $\ell, k$ the derived stacks $\dstackCohps(\scrD_0, I_{\ell,k})$ admit an action of a torus $T$ such that Constructions~\ref{construction:limiting_COHA-1} and \ref{construction:limiting_COHA-2} are $T$-equivariant, we denote by 
	\begin{align}
		\coha^T_{\scrD_0, \tau_\infty^+}
	\end{align}
	the corresponding $T$-equivariant limiting cohomological Hall algebra.
\end{remark}

\section{Stabilization and comparison}\label{sec:stabilization}

Observe that for fixed $k \geqslant \ell$ there is a canonical map
\begin{align}
	\calS_\bullet \dstackCohps\big(\scrD_0, I_{\ell,\infty}\big) \longrightarrow \calS_\bullet \dstackCohps\big( \scrD_0, I_{\ell,k} \big) \ . 
\end{align}
This map propagates first to a map
\begin{align}
	\calS_\bullet \dstackCohps\big( \scrD_0, I_{\ell,\infty}\big) \longrightarrow \calS_\bullet \dstackCohps\big( \scrD_0, I_{\ell,\infty}^+\big) 
\end{align}
and then to a map
\begin{align}\label{eq:limiting_vs_limit}
	\phi \colon \fcolim_\ell \calS_\bullet \dstackCohps\big( \scrD_0, I_{\ell,\infty}\big) \longrightarrow \calS_\bullet \dstackCohps\big( \scrD_0, \tau_\infty^+\big) 
\end{align}
of $\Lambda$-graded simplicial objects in $\Ind(\Pro(\dSt_k))$. In practice, we are interested in knowing that $\phi$ is an equivalence. This can be ensured under some more stringent conditions.

\begin{construction}\label{construction:HN}
	Fix a slicing $\calP$ and an open pair $(\Lambda,v)$ on $\scrD_0$. Moreover, fix a group homomorphism $Z\colon \Lambda\to \C$ such that $(\Lambda,v,\calP,Z)$ is a pre-stability condition on $\scrD_0$ (cf.\ Definition~\ref{def:pre-stability-condition}).
	
	Fix now $\ell\in \N$. We say that an object $E\in \calP(I_\ell)$ is of class $\widetilde{\bfv} = (\bfv_1, \bfv_2, \ldots, \bfv_s)\in \Lambda^s$ if for the HN factor $E_i/E_{i+1}$ of the Harder-Narasimhan filtration \eqref{eq:HN-filtration-slicing} of $E$ we have $v(E_i/E_{i+1})=\bfv_i$ for $i=1, \ldots, s$. In such a case, we set $\phi_i=\phi(E_i/E_{i+1})\eqqcolon \phi(\bfv_i)$ for $i=1, \ldots, s$.
	
	Given $\bfv \in \Lambda$, we denote by $\mathsf{HN}(\bfv)$ the set of tuples $\widetilde{\bfv} = (\bfv_1, \bfv_2, \ldots, \bfv_s)$ of variable length satisfying
	\begin{align}
		\bfv_1 + \cdots + \bfv_s = \bfv 
	\end{align}
		and
	\begin{align}
		\phi(\bfv_s) > \phi(\bfv_{s-1}) > \cdots > \phi(\bfv_1) \ . 
	\end{align}
	Following \cite[\S3]{Bayer_Short_proof}, we write $\mathsf{HNP}(\widetilde{\bfv})$ for the convex hull of the set
	\begin{align}
		\{0, Z(\bfv_1), Z(\bfv_2), \ldots, Z(\bfv_s), Z(\bfv)\} \subset \C \ . 
	\end{align}
	Given two elements $\widetilde{\bfv}, \widetilde{\bfw} \in \mathsf{HN}(\bfv)$, we write
	\begin{align}
		\widetilde{\bfv} \preceq \widetilde{\bfw} \iff \mathsf{HNP}(\widetilde{\bfv}) \subset \mathsf{HNP}(\widetilde{\bfw}) \ . 
	\end{align}
	This makes $\mathsf{HN}(\bfv)$ into a countable poset.
	 
	We let $\dstackCohps^{\mathsf{HN}}(\scrD_0,\tau_\ell;\widetilde{\bfv})$ be the full substack of $\dstackCoh(\scrD_0,\tau)$ defined as follows: given $S \in \dAff_k$, an $S$-point $x \colon S \to \dstackCoh(\scrD_0,\tau)$ belongs to $\dstackCohps^{\mathsf{HN}}(\scrD_0,\tau_\ell;\widetilde{\bfv})$ if and only if for every $k$-point $y \colon \Spec(k) \to S$ the object $E\in \scrD_0$ associated to the composite $x \circ y \colon \Spec(k) \to \calM_0$ admits an Harder-Narasimhan filtration of the form \eqref{eq:HN-filtration-slicing} such that $v(E_i/E_{i+1})=\bfv_i$ for $i=1, \ldots, s$. We shall call  $\dstackCohps^{\mathsf{HN}}(\scrD_0,\tau_\ell;\widetilde{\bfv})$ an \textit{Harder-Narasimhan stratum} of $\dstackCoh(\scrD_0,\tau)$.
\end{construction}

\begin{assumption}\label{assumption:quasi-compact_interval}
	\hfill
	\begin{enumerate}\itemsep=0.2cm
		\item In addition to the slicing $\calP$ and the open pair $(\Lambda,v)$ on $\scrD_0$, we assume the existence of a central charge $Z\colon \Lambda\to \C$ in such a way that $(\Lambda,v,\calP,Z)$ is a pre-stability condition on $\scrD_0$.
		
		\item For every $\ell \in \N$, every $\bfv \in \Lambda$ and every $\widetilde{\bfv} \in \mathsf{HN}(\bfv)$, the corresponding Harder-Narasimhan stratum
		\begin{align}
			\dstackCohps^{\mathsf{HN}}(\scrD_0,\tau_\ell;\widetilde{\bfv}) 
		\end{align}
		is a locally closed substack of $\dstackCohps(\scrD_0,\tau_\ell)$.
		
		\item For every Mukai vector $\bfv \in \Lambda$ and every pair of positive integers $k \geqslant \ell$, the derived stack $\dstackCohps(\scrD_0, I_{\ell,k};\bfv)$ is quasi-compact.
	\end{enumerate}
\end{assumption}

Under this assumption, the Harder-Narasimhan strata form a stratification
\begin{align}
	\dstackCohps(\scrD_0, \tau_\ell;\bfv) = \bigsqcup_{\widetilde{\bfv} \in \mathsf{HN}(\bfv)} \dstackCohps^{\mathsf{HN}}(\scrD_0,\tau_\ell;\widetilde{\bfv}) \ . 
\end{align}
of $\dstackCohps(\scrD_0, \tau_\ell;\bfv)$.

\begin{theorem}\label{thm:limiting_vs_limit}
	Under Assumptions~\ref{assumption:limiting_2_Segal_stack} and \ref{assumption:quasi-compact_interval}, for every $\ell \in \N$, $[n] \in \bfDelta$ and $\bfv \in \Lambda^n$, there exists $N \in \N$ such that for every $k \geqslant N$, the canonical comparison map
	\begin{align}
		\calS_n \dstackCohps\big( \scrD_0, I_{\ell,\infty} ; \bfv \big) \longrightarrow \calS_n \dstackCohps\big( \scrD_0, I_{\ell,k} ; \bfv \big) 
	\end{align}
	is an equivalence.
	Moreover, $N$ does not depend on the simplicial variable $[n]$.
	In particular:
	\begin{enumerate}\itemsep=0.2cm
		\item \label{item:limiting_vs_limit-1} The canonical comparison map
		\begin{align}
			\dstackCohps\big( \scrD_0, I_{\ell,\infty}\big) \longrightarrow \dstackCohps\big(\scrD_0, I_{\ell,\infty}^+\big) 
		\end{align}
		is an equivalence of $\Lambda$-graded pro-derived stacks.
		
		\item \label{item:limiting_vs_limit-2} The map \eqref{eq:limiting_vs_limit} is an equivalence of $\Lambda$-graded simplicial objects.
	\end{enumerate}
\end{theorem}

\begin{proof}
	It is clearly enough to prove the first statement, as items \eqref{item:limiting_vs_limit-1} and \eqref{item:limiting_vs_limit-2} are obvious consequences.
	Fix $\ell \in \N$, $[n] \in \bfDelta$ and
	\begin{align}
		\underline{\bfv} = (\bfv_{0,1}, \bfv_{1,2}, \ldots, \bfv_{n-1,n}) \in \Lambda^n \ . 
	\end{align}
	We have to prove that the induced map
	\begin{align}
		\calS_n \dstackCohps\big( \scrD_0, I_{\ell,\infty}; \underline{\bfv}\big) \longrightarrow \calS_n \dstackCohps\big(\scrD_0, I_{\ell,\infty}^+;\underline{\bfv}\big) 
	\end{align}
	is an equivalence in $\Pro(\dSt_k)$.
	Thanks to the $\Lambda$-graded $2$-Segal property, it is enough to deal with the cases $n = 1$ and $n = 2$.
	Since $\calM_0^{\mathsf{cat}}$ is closed under extensions by assumption and $v \colon K_0(\scrD_0) \to \Lambda$ is additive, we see that it is in fact enough to deal with the case $n = 1$.
	
	\medskip
	
	Fix two positive integers $k \geqslant \ell$ and $\bfv \in \Lambda$. As shown in Construction~\ref{construction:HN}, there is a stratification of $\dstackCohps(\scrD_0, \tau_\ell;\bfv)$ in locally closed substacks parametrized by the Harder-Narasimhan types $\mathsf{HN}(\bfv)$:
	\begin{align}
		\dstackCohps(\scrD_0, \tau_\ell;\bfv) = \bigsqcup_{\widetilde{\bfv} \in \mathsf{HN}(\bfv)} \dstackCohps^{\mathsf{HN}}(\scrD_0,\tau_\ell;\widetilde{\bfv}) \ . 
	\end{align}
	Since $\dstackCohps(\scrD_0,I_{\ell,k};\bfv)$ is quasi-compact for every $\ell < k \leqslant \infty$, it must encounter only finitely many strata of the above stratification. In particular, the set $J(k)$ of phases $\theta$ that satisfy $\theta > a_\infty$ and that occur for families in $\dstackCohps(\scrD_0, I_{\ell,k};\bfv)$ is finite.
	Define
	\begin{align}
		\theta_{\mathsf{min}}(k) \coloneqq \begin{cases}
			\mathsf{min}( J(k) ) & \text{if } J(k) \neq \emptyset \\
			a_\ell & \text{if } J(k) = \emptyset \ .
		\end{cases} 
	\end{align}
	Since $\ell$ is fixed, for $k' \geqslant k \geqslant \ell$ one has $J(k') \subset J(k)$, and therefore one has
	\begin{align}
		\theta_{\mathsf{min}}(k') \geqslant \theta_{\mathsf{min}}(k) \ . 
	\end{align}
	On the other hand, one also as $a_\infty < \theta_{\mathsf{min}}(k)$, and furthermore for $k > \ell$ one has
	\begin{align}
		\theta_{\mathsf{min}}(k) \leqslant a_k \quad \text{if and only if}\quad  J(k) \neq \emptyset \ . 
	\end{align}
	In particular, if $N$ is such that $a_N < \theta_{\mathsf{min}}(\ell)$, it follows that for every $k \geqslant N$ one has $J(k) = \emptyset$.
	This immediately implies that for $k \geqslant N$ the natural inclusion
	\begin{align}
		\dstackCohps\big( \scrD_0; I_{\ell,\infty} ; \bfv \big) \subset \dstackCohps\big( \scrD_0; I_{\ell,k}; \bfv \big) 
	\end{align}
	is an equality, and the conclusion follows.
\end{proof}

\begin{corollary}\label{cor:assumption-holds}
	Under Assumptions~\ref{assumption:limiting_2_Segal_stack} and \ref{assumption:quasi-compact_interval}, for any $\ell_1\in \N$ and for every pair $\underline{\bfv} = (\bfv_{0,1},\bfv_{1,2}) \in \Lambda^2$, there exists an integer $N \geqslant \ell_1$ such that for every $k' \geqslant k \geqslant \ell_2 \geqslant N$, the square
	\begin{align}
		\begin{tikzcd}[ampersand replacement=\&]
			\calS_2 \dstackCohps\big( \scrD_0, I_{\ell_2,k'} \times I_{\ell_1,k'} \times I_{\ell_2,k'} ; \underline{\bfv} \big) \arrow{r}{\partial_1} \arrow{d} \& \dstackCohps\big( \scrD_0, I_{\ell_1,k'}; \bfv_{0,1} + \bfv_{1,2}\big) \arrow{d} \\
			\calS_2 \dstackCohps\big( \scrD_0, I_{\ell_2,k} \times I_{\ell_1,k} \times I_{\ell_2,k} ; \underline{\bfv} \big) \arrow{r}{\partial_1} \& \dstackCohps\big( \scrD_0, I_{\ell_1,k}; \bfv_{0,1} + \bfv_{1,2}\big)
		\end{tikzcd} 
	\end{align}
	is a pullback.
\end{corollary}

\begin{proof}
	It follows from Theorem~\ref{thm:limiting_vs_limit} that we can choose $N$ such that the vertical maps become equivalences.
\end{proof}

\begin{corollary}\label{cor:limiting_COHA_vs_COHA_limit}
	Under Assumptions~\ref{assumption:limiting_2_Segal_stack}, \ref{assumption:limiting_CoHA_I}, and Assumption~\ref{assumption:quasi-compact_interval}, and the additional assumptions
	\begin{enumerate}\itemsep=0.2cm
		\item the $t$-structure $\tau_\infty\coloneqq \tau_{a_\infty}$ is open; 
		
		\item \label{item:2} the map
		\begin{align}
			\partial_0 \times \partial_2 \colon \calS_2 \dstackCohps\big( \scrD_0, \tau_\infty \big) \longrightarrow \dstackCohps\big( \scrD_0, \tau_\infty \big) \times \dstackCohps\big( \scrD_0, \tau_\infty \big)
		\end{align}
		is derived lci;
		
		\item the map
		\begin{align}
			\partial_1 \colon \calS_2 \dstackCohps\big( \scrD_0, \tau_\infty \big) \longrightarrow \dstackCohps\big(\scrD_0, \tau_\infty\big) 
		\end{align}
		of $\Lambda$-graded derived stacks is representable by proper algebraic spaces;
	\end{enumerate}
	the cohomological Hall algebra $\coha_{\scrD_0, \tau_\infty}$ associated to $\dstackCohps(\scrD_0, \tau_\infty)$ is well defined and it is isomorphic to $\coha_{\scrD_0, \tau_\infty^+}$. A similar statement holds in the $T$-equivariant setting.
\end{corollary}

\begin{remark}
	Note that if the smooth Serre functor $\mathsf S_\scrD^![2]$ restricted to $\scrD_0$ is $\tau_\infty$-exact, then the condition~\eqref{item:2} above holds. 
\end{remark}
%
%

\newpage
\part{COHAs of quivers, reflection functors, Yangians, and braid group actions}\label{part:Yangians}

\section{Lie theory associated to quivers}\label{sec:quiver}

We fix a quiver $\qv$ with vertex set $I$ and edge set $\Omega$. We assume that both $I$ and $\Omega$ are finite, and that $\Omega$ does not contain any edge-loop. 

\subsection{(Co)root and (co)weight lattices}

We denote by $\rootlattice_\qv\coloneqq \Z I$ the free abelian group on $I$ and we denote its canonical basis by $\Delta_\qv\coloneqq \{\alpha_i\, \vert\, i\in I\}$. We call the $\alpha_i$'s the \textit{simple roots}, $\Delta_\qv$ the set of simple roots, and $\rootlattice_\qv$ the \textit{root lattice} of $\qv$. Let $\{\Lomega_i\}_{i\in I}$ be the dual basis in $\coweightlattice_\qv\coloneqq (\Z I)^\vee$. Thus, $\Lomega_i$ is the \textit{fundamental coweight} defined by the linear form $\bfd\in \rootlattice_\qv\mapsto d_i\in \Z$ for any $i\in I$. Set 
\begin{align}\label{eq:Lrho}
	\Lrho\coloneqq\sum_{i\in I}\Lomega_i\ .
\end{align}
We consider the canonical pairing
\begin{align}\label{eq:canonical-pairing}
	(-, -)\colon \coweightlattice_\qv\times \rootlattice_\qv\longrightarrow \Z\ .
\end{align} 

For any $i\in I$, we define the simple coroot $\Lalpha_i\in \coweightlattice_\qv$ by setting $(\Lalpha_i,\alpha_i)\coloneqq 2$ and $(\Lalpha_i,\alpha_j)$ is the negative of the number of edges between the vertices $i$ and $j$ in the quiver $\qv$ for $i\neq j$. The matrix $A=(a_{i,j})_{i,j\in I}$, with $a_{i,j}\coloneqq (\Lalpha_i,\alpha_j)$, is a symmetric generalized Cartan matrix. It gives rise to a derived Kac-Moody algebra $\frakg_\qv$, which is $\Z I$-graded. 
\begin{notation}
	Let $\frakn_\qv$ be the negative nilpotent Lie subalgebra of $\frakg_\qv$, graded by the negative roots.
\end{notation}

\subsection{Weyl and braid groups}

The \textit{Weyl group} $W_\qv$ of $\qv$ is the subgroup of $\mathsf{GL}(\Z I)$ generated by the simple reflections
\begin{align}
	s_i\colon \alpha_j \to \alpha_j - (\Lalpha_j, \alpha_i) \alpha_i\ .
\end{align}
This is a Coxeter group, subject to $s_i^2=1$ for any $i\in I$, and the \textit{braid relations}
\begin{align}\label{eq:braid-relations}
	\begin{cases}
		s_is_j=s_js_i  &\text{if } a_{i, j}=0 \ ,\\
		s_is_js_i=s_js_is_j &\text{if } a_{i, j}=-1\ .
	\end{cases}
\end{align}
The length $\ell(w)$ of $w \in W_\qv$ is the minimal $\ell$ such that $w = s_{i_1}\cdots s_{i_\ell}$. Such an expression is called a \textit{reduced expression}. The Weyl group $W_\qv$ acts on $\C I=\rootlattice_\qv\otimes_\Z \C$ and $(\C I)^\vee=\coweightlattice_\qv\otimes_\Z \C$ by
\begin{align}
	s_i(\bfd)=\bfd-(\Lalpha_i,\bfd)\alpha_i\quad\text{and}\quad s_i(\lambda)=\lambda-(\lambda,\alpha_i)\Lalpha_i \ ,
\end{align}
with $\bfd\in\C I$ and $\lambda\in (\C I)^\vee$.

The \textit{Braid group} $B_\qv$ of $\qv$ is the group with generators $T_i$ with $i\in I$ and the braid relations \eqref{eq:braid-relations} with $T_i$ in place of $s_i$. If $w = s_{i_1}\cdots s_{i_\ell}\in W_\qv$ is a reduced expression, we set $T_w\coloneqq T_{i_1}\cdots T_{i_\ell}$. These elements are well-defined, i.e., they do not depend on the chosen reduced expression, and satisfy $T_{uv}= T_u T_v$ when $\ell(uv)=\ell(u)+\ell(v)$.

\section{The nilpotent cohomological Hall algebra of a quiver}\label{sec:nilpotent-quiver-COHA}

In this section, we introduce the 2-dimensional cohomological Hall algebra associated to nilpotent finite-dimensional representations of the preprojective algebra of a fixed quiver $\qv$.

\subsection{Moduli stack of (nilpotent) finite-dimensional representations of the preprojective algebra of a quiver}

In this section, we shall introduce the moduli stacks of finite-dimensional representations of the preprojective algebra of a quiver and their nilpotent counterparts.

\subsubsection{Quiver representations}

Let $\C \qv$ be the \textit{path algebra} of $\qv$. As customary, we call a \textit{representation of $\qv$} a (right) $\C\qv$-module.

We denote by $\catMod(\qv)$ be the category of representations of $\qv$ and by $\catmod(\qv)$ the subcategory of $\catMod(\qv)$ consisting of finite-dimensional representations. For each dimension vector $\bfd\in\N I$, let $\catmod(\qv)_\bfd$ be the subcategory of representations of dimension $\bfd$. Assigning to a finite-dimensional representation of $\qv$ its dimension vector yields a map from the complexified Grothendieck group $\sfK_0(\catmod(\qv))\otimes\C$ to $\C I$. The Euler form on $\C I$ is given by
\begin{align}\label{eq:euler-form-quiver}
	\langle \bfd_1, \bfd_2\rangle \coloneqq \sum_{i\in I}d_{1,i}\, d_{2,i}  - \sum_{e\colon i\to j\in\Omega} d_{1,i}\, d_{2,j}=\dim \Hom(M_1,M_2)-\dim \Ext^1(M_1,M_2)
\end{align}
where $M_1, M_2$ are two representations of $\qv$ of dimension $\bfd_1, \bfd_2$, respectively. We introduce the symmetric bilinear form
\begin{align}\label{eq:symmetric-euler-form-quiver}
	(\bfd_1, \bfd_2)\coloneqq \langle \bfd_1,\bfd_2\rangle+\langle \bfd_2, \bfd_1\rangle
\end{align}
for $\bfd_1, \bfd_2\in \N I$.

\subsubsection{Preprojective algebras}\label{subsubsec:preprojective-algebra}

Recall that  $\doubleqv=(I, \overline{\Omega})$ is the \textit{double} quiver of $\qv$. The \textit{preprojective algebra} $\Pi_\qv$ is the quotient of the path algebra $\C \doubleqv$ by the relations 
\begin{align}
	\sum_{e\in\Omega}(e^\ast e -e e^\ast)=0\ .
\end{align}
Let $\ModPi$ be the abelian category of representations of $\Pi_\qv$ and let $\modPi$ be the subcategory of $\ModPi$ consisting of finite dimensional representations. There is an obvious forgetful functor $\modPi \to \catmod(\qv)$. For each dimension vector $\bfd\in\N I$, let $\modPi_\bfd$ be the subcategory of representations of $\Pi_\qv$ of dimension $\bfd$. Since $\qv$ has no edge-loops, there is for any $i\in I$ a unique $\Pi_\qv$-module $\sigma_i$ with dimension $\alpha_i$. Furthermore, it is simple.

We shall need the following crucial property of $\modPi$ (cf.\ \cite[Lemma~2.18]{SY13} or \cite[Proposition~3.1]{SV_generators}):
\begin{theorem}\label{thm:PiQis2CY}
	Assume that $\qv$ is not a finite type Dynkin quiver. Then the algebra $\Pi_\qv$ is Calabi-Yau of homological dimension two, i.e., there are functorial isomorphisms
	\begin{align}\label{eq:Pi_Qis2CY}
		\Ext_{\Pi_\qv}^i(M,N)^\vee \simeq \Ext_{\Pi_\qv}^{2-i}(N,M)
	\end{align}
	for any two finite-dimensional $\Pi_\qv$-modules $M,N$. Moreover, we have
	\begin{align}\label{eq:eulerformPi_Q}
		\sum_{i=0}^2 (-1)^i\dim \Ext_{\Pi_\qv}^i(M,N)=(\bfd_M,\bfd_N)\ ,
	\end{align}
	where $\bfd_M, \bfd_N$ are the dimension vectors of $M$ and $N$.
\end{theorem}

\begin{remark}
	The condition on the quiver is not an issue for us: if $\qv$ is of finite type then we may embed it in a larger quiver $\qv'$ (which is not of Dynkin type) and extend trivially any representation of $\qv$ to $\qv'$.
\end{remark}

\subsubsection{Nilpotent representations}

Let $\calI$ be the two-sided ideal of $\Pi_\qv$ which is generated by all the arrows of $\doubleqv$.
\begin{definition}[{\cite[Definition~2.14]{SY13}}]\label{def:nilpotent-representation}  
	A representation $M$ of $\Pi_\qv$ is \textit{nilpotent} if there exists $\ell \geq 1$ such that $M \calI^\ell=0$.  
\end{definition}

\begin{remark}\label{rem:nilpotency}
	Let $M$ be a finite-dimensional representation of $\Pi_\qv$. It is nilpotent if and only if $M$ has finite length and its composition factors consist of $\sigma_0, \ldots, \sigma_e$.
\end{remark}
Nilpotent representations form Serre subcategories of $\ModPi$ and $\modPi$. In particular, they are stable by extensions. We denote by $\nilpPi$ the category of nilpotent finite-dimensional representations of $\Pi_\qv$.

\subsubsection{Stability conditions}\label{subsec:stabquiver}

Fix $\Ltheta\in \coweightlattice_\qv$. In what follows, we shall call $\Ltheta$ also a \textit{stability condition}. 
\begin{definition}
	The \textit{($\Ltheta$-)degree} of a dimension vector $\bfd=(d_i)_{i\in I}\in\N I\smallsetminus \{0\}\subset \rootlattice_\qv$ is defined as
	\begin{align}
		\deg_{\Ltheta}(\bfd)\coloneqq (\Ltheta, \bfd) \ .
	\end{align}
	The \textit{($\Ltheta$-)slope} of $\bfd$ is
	\begin{align}
		\mu_{\Ltheta}(\bfd)\coloneqq \frac{\deg_{\Ltheta}(\bfd)}{(\Lrho, \bfd)}\ ,
	\end{align}
	where $\Lrho$ is introduced in Formula ~\eqref{eq:Lrho}.
	
	Let $M$ be a finite-dimensional representation of $\Pi_\qv$. Its \textit{degree} (resp.\ \textit{slope}) is the degree (resp.\ slope) of its dimension vector.
\end{definition}

\begin{definition}
	A representation $M$ of $\Pi_\qv$ of dimension $\bfd\neq 0$ is \textit{$\Ltheta$-semistable} if for any subrepresentation $N$ of $M$ of dimension $\bfd'$ we have
	\begin{align}
		\mu_{\Ltheta}(\bfd')\leq \mu_{\Ltheta}(\bfd)\ . 
	\end{align}
	A nonzero representation $M$ is called \textit{$\Ltheta$-stable} if the strict inequality holds for any nonzero proper subrepresentation $N \subset M$.
\end{definition}

A finite-dimensional $\Pi_\qv$-module $M$ has a unique filtration 
\begin{align}\label{eq:HN-filtration-representations}
	0\eqqcolon M_{s+1}\subset M_s\subset \cdots\subset M_1=M
\end{align}
with subsequent quotients being $\Ltheta$-semistable of strictly decreasing $\Ltheta$-slopes $\ell_s>\cdots>\ell_1$. This filtration is called the \textit{Harder-Narasimhan (HN) filtration} of $M$ and the successive quotients are the \textit{HN factors} of $M$. For each representation $M$, we define the rational numbers $\mu_{\Ltheta\textrm{-}\mathsf{min}}(M)$ and $\mu_{\Ltheta\textrm{-}\mathsf{max}}(M)$ to be
the minimal and maximal $\Ltheta$-slopes of the HN-factors of $M$: we have
\begin{align}
	\mu_{\Ltheta\textrm{-}\mathsf{min}}(M)\coloneqq \ell_1 =\mu_{\Ltheta}(M_1/M_2)\quad \text{and}\quad \mu_{\Ltheta\textrm{-}\mathsf{max}}(M)\coloneqq \ell_s=\mu_{\Ltheta}(M_s)\ .
\end{align}

\subsubsection{Moduli stacks of representations}\label{subsubsec:moduli-stacks}

We denote by $\stackRep(\Pi_\qv)$ the classical moduli stack parame\-tri\-zing finite-dimensional representations of $\Pi_\qv$. It splits as a disjoint union
\begin{align}
	\stackRep(\Pi_\qv)=\bigsqcup_\bfd\stackRep_\bfd(\Pi_\qv)
\end{align}
into closed and open connected components, according to the dimension vector $\bfd \in \rootlattice_\qv$. Each $\stackRep_\bfd(\Pi_\qv)$ is a finite-type classical geometric stack. 

The stack $\stackRep_\bfd(\Pi_\qv)$ may be realized as a quotient stack as follows. The group $\G_\bfd\coloneqq \prod_i \GL(d_i)$ acts by conjugation on the space 
\begin{align}
	\sfE_\bfd\coloneqq \bigoplus_{e\colon i \to j \in \overline{\Omega}} \Hom(\C^{d_i}, \C^{d_j})
\end{align}
of $\bfd$-dimensional representations of $\doubleqv$. Consider the $\G_\bfd$-equivariant map
\begin{align}
	\mu_\bfd \colon \sfE_\bfd \longrightarrow \bigoplus_i \frakgl(d_i)\ , \quad \big(x_e\big)_{e \in \overline{\Omega}} \longmapsto \sum_{e \in\Omega} (x_e x_{e^\ast}-x_{e^\ast}x_e)\ .
\end{align}
Then  
\begin{align}
	\stackRep_\bfd(\Pi_\qv) \simeq \mu_\bfd^{-1}(0)/\G_\bfd\ .
\end{align}
For any $i \in I$, there is a \textit{tautological vector bundle} $\calV_i$ on $\stackRep_\bfd(\Pi_\qv)$ of rank $d_i$, which is pulled back from $\sfB \GL(d_i)$, and for any arrow $e \colon i \to j \in \overline{\Omega}$ there is a \textit{tautological map} $x_e\colon \calV_i \to \calV_j$. The maps $(x_e)_{e \in \overline{\Omega}}$ satisfy the preprojective relations. We will refer to the data $(\calV_i, x_e)_{i \in I, e \in \overline{\Omega}}$ as the \textit{tautological $\Pi_\qv$-module} on $\stackRep_\bfd(\Pi_\qv)$, and sometimes abusively denote it $\calV_\bfd$. The collection of tautological $\Pi_\qv$-modules on $\stackRep_\bfd(\Pi_\qv)$ gives rise to a locally free sheaf $\calV$ on $\stackRep(\Pi_\qv)$.

We let $\Lambda_\qv$ stand for the closed substack of $\stackRep(\Pi_\qv)$ parametrizing nilpotent representations. We have 
\begin{align}
	\Lambda_\qv=\bigsqcup_{\bfd \in \N I} \Lambda_{\bfd}\ .
\end{align}
This is a finite type classical geometric stack. Furthermore, it is pure of dimension
\begin{align}
	\dim \Lambda_{\bfd}=-\langle \bfd, \bfd \rangle\ .
\end{align}

There exists a derived enhancement $\dstackRep_\bfd(\Pi_\qv)$ of $\stackRep_\bfd(\Pi_\qv)$ which is of (virtual) dimension equal to $-(\bfd,\bfd)$ (cf.\ \cite[\S2.1.4]{VV_KHA} or \cite[\S~I.2]{DPS_Torsion-pairs}). In particular, thanks to Theorem~\ref{thm:PiQis2CY} we deduce\footnote{Recall that if $\qv$ is of finite Dynkin type then we embed it into a larger quiver.} that $\dstackRep(\Pi_\qv)$ is a derived lci derived geometric stack, of finite type over $\C$.

Let $\Ltheta \in \coweightlattice$ be a stability condition. The subfunctor of $\Ltheta$-semistable $\Pi_\qv$-represen\-ta\-tions of dimension $\bfd\in \rootlattice_\qv$ forms an open substack $\stackRep_\bfd^{\Ltheta\textrm{-}\mathsf{ss}}(\Pi_\qv)$ of $\stackRep_\bfd(\Pi_\qv)$. The former admits a canonical enhancement, so there is also a derived open substack  $\dstackRep_\bfd^{\Ltheta\textrm{-}\mathsf{ss}}(\Pi_\qv)$ of $\dstackRep_\bfd(\Pi_\qv)$. Note that the construction of a canonical derived enhancement of an open embedding of a geometric classical stack into a geometric derived stack follows from \cite[Proposition~2.1]{STV}. We shall use this result also in the rest of this Part without explicitly mentioning it.

\begin{definition}
	The \textit{derived} moduli stack $\dLambda_\qv$ of nilpotent finite-dimensional representations of $\Pi_\qv$ is the formal completion of $\Lambda_\qv$ inside $\dstackRep(\Pi_\qv)$. 
\end{definition}
We denote by $\dLambda_\bfd^{\Ltheta\textrm{-}\mathsf{ss}}$ the derived moduli stack of $\Ltheta$-semistable nilpotent representations of $\Pi_\qv$ of dimension $\bfd$.

\begin{remark}
	One can show that $\stackRep_\bfd(\Pi_\qv)$ is equipped with a symplectic structure and that $\Lambda_\bfd$ is a generically Lagrangian substack (cf.\ \cite{BSV_Nilpotent}). We will not use this. 
\end{remark}

\subsection{Definition of the preprojective Cohomological Hall algebra}\label{sec:defCOHAnilp}

We shall now define, following \cite{SV_generators} (see also \cite{SV_Cherednik}) the \textit{cohomological Hall algebra} (COHA for short) associated to the category $\modPi$. Since we will mostly be interested in the nilpotent case in this paper, we focus on that. Rather than the Hamiltonian formalism used in \textit{loc. cit.}, we will use here the formalism developed in \cite{Porta_Sala_Hall, DPS_Torsion-pairs, DPSSV-1}. The two approaches coincide thanks to \cite[Proposition~6.1.5]{COHA_surface} and \cite[\S5]{Porta_Sala_Hall}. Note that \cite{COHA_surface} considers the COHA of the one-loop quiver, while, in \cite{Porta_Sala_Hall}, the authors compared their construction of the K-theoretical Hall algebra via Derived Algebraic Geometry with the construction in \cite{COHA_surface}. Both arguments used in \textit{loc.cit.} can be readapted to the present setting.

\subsubsection{Equivariant Borel-Moore homology of the stack of nilpotent $\Pi_\qv$-representations}

Consider the action of the torus 
\begin{align}
	\Tmax\coloneqq \Big\{(\gamma_e)_{e\in \overline{\Omega}}\in (\C^\ast)^{\doubleOmega}\, \Big\vert\, \gamma_e\gamma_{e^\ast}=\gamma_f\gamma_{f^\ast} \quad\forall \, e, f\in \Omega \Big\} \simeq (\G_m)^\Omega\times\G_m 
\end{align}
on the stack $\Lambda_\qv$, which is given by the following formula
\begin{align}\label{eq:torus-action-stack}
	(q_e)_e\cdot x_h\coloneqq q_h x_h  
\end{align}
for $h \in \overline{\Omega}$. Fix a subtorus $T \subseteq \Tmax$ and let $\Hbullet_T$ be the cohomology ring of the \textit{classifying stack} $\sfB T\coloneqq \mathsf{pt}/T$. Set
\begin{align}
	\cohaqv^T\coloneqq \HBMbulletT(\dLambda_\qv)\simeq \HBMbulletT(\Lambda_\qv) \quad \text{and} \quad \cohaqvd^T\coloneqq \HBMbulletT(\dLambda_\bfd)\simeq \HBMbulletT( \Lambda_\bfd )
\end{align}
for $\bfd\in \N I$. The obvious $\N I$-grading on $\cohaqv^T$ is called the \textit{horizontal grading}. The space $\cohaqv^T$ carries a shifted cohomological grading, called the \textit{vertical grading}, such that
\begin{align}
	\coha_{\bfd, k}^T\coloneqq \sfH^T_{k-(\bfd, \bfd)}( \dLambda_\bfd )
\end{align}
for $k\in \Z$. The top homological degree is in zero vertical degree, i.e., we have
\begin{align}
	\coha_0^T\coloneqq \sfH^T_{\mathsf{top}}( \dLambda_\qv )=\bigoplus_\bfd \sfH^T_{-(\bfd,\bfd)}(\dLambda_\bfd)\ .
\end{align}

For a future use, let us quote the following important structural results. 
\begin{theorem}[{\cite[Theorem~A]{SV_generators}}, \cite{Davison_Purity}]\label{thm:purityandKac}
	For any $\bfd \in \N I$, the following holds:
	\begin{enumerate}\itemsep0.2cm
		\item \label{lem:free} $\cohaqvd^T$ is free as a $\Hbullet_T$-module.
		\item $\cohaqvd^T$ is concentrated in even degrees and carries a pure mixed Hodge structure (in particular, it is equivariantly formal).
		\item for any $i \geq 0$, the equivariant cycle map $\sfH^{\mathsf{mot}, T}_i(\dLambda_\bfd; 0) \to \sfH^T_{2i}(\dLambda_\bfd)$ is surjective.
		\item \label{item:graded-dim-coha-Kac-pols} The graded dimension of $\coha_\qv^T$ is given by the following generating function
		\begin{align}
			\sum_{\bfd,k} \dim \coha_{\bfd,2k}^T\,z^\bfd q^k= (1-q^{-1})^{-\dim A}\Exp\Big( \frac{1}{1-q^{-1}}\sum_{\bfd\in \N I} A_{\bfd}(q^{-1}) z^\bfd \Big) \ ,
		\end{align}
		where $A_{\qv,\bfd}(t)$ is Kac's polynomial counting absolutely indecomposable representations of $\qv$ over finite fields.  
	\end{enumerate}
\end{theorem}

\begin{remark}
	It is in fact possible to prove that the cycle map $\sfH^{\mathsf{mot}, T}_i(\dLambda_\bfd; 0) \to \sfH^T_{2i}(\dLambda_\bfd)$ is an isomorphism.\end{remark}

\subsubsection{Hall multiplication}

Consider the convolution diagram
\begin{align}\label{eq:convolution-def-COHA-quiver}
	\begin{tikzcd}[ampersand replacement=\&]
		\dstackRep(\Pi_\qv) \times \dstackRep(\Pi_\qv) \& \dstackRepext(\Pi_\qv) \ar{r}{p} \ar[swap]{l}{q} \& \dstackRep(\Pi_\qv)
	\end{tikzcd}\ ,
\end{align}
where $\dstackRepext(\Pi_\qv)$ is the derived stack parametrizing short exact sequence $0 \to M \to R \to N \to 0$ of finite-dimensional $\Pi_\qv$-modules, and where the maps $q$ and $p$ associate to such a sequence the pair of representations $(N,M)$ and the representation $R$ respectively. We use the notation $\dstackRepext(\Pi_\qv)\coloneqq \calS_2\dstackRep(\Pi_\qv)$ and $q=\partial_0\times \partial_2$, $p=\partial_1$, where the maps $\partial_i$ are defined in Formula~\eqref{eq:partial}, because it is more standard.

Note that $p$ is a proper representable map, while $q$ is a derived lci map, see \cite[\S2.2]{VV_KHA} or \cite[\S~I.4]{DPS_Torsion-pairs}. 

As nilpotent representations form a Serre subcategory, the diagram~\eqref{eq:convolution-def-COHA-quiver} restricts by base change to a diagram involving $\dLambda_\qv$ in place of $\dstackRep(\Pi_\qv)$:
\begin{align}\label{eq:definition-COHA-product-nil}
	\begin{tikzcd}[ampersand replacement=\&]
		\dLambda_\qv \times \dLambda_\qv \& \dLambdaext_\qv \ar{r}{p} \ar[swap]{l}{q} \& \dLambda_\qv
	\end{tikzcd}\ .
\end{align}

By arguing as in \cite[\S\ref*{foundation-sec:COHA}]{DPSSV-1} (and using \cite[\S~II.5.2]{DPS_Torsion-pairs} instead of \cite[\S4]{Porta_Sala_Hall}), we obtain a well-defined associative algebra structure on both
\begin{align}
	\HBMbulletT(\dstackRep(\Pi_\qv)) \quad \text{and} \quad \coha_\qv^T\coloneqq \HBMbulletT(\Lambda_\qv)=\bigoplus_\bfd \HBMbulletT(\Lambda_\bfd) \ ,
\end{align}
by the map $p_\ast q^!$. 

We summarize the main properties of $\coha^T_\qv$ as follows.

\begin{theorem}[{\cite[\S5]{SV_generators}}] 
	The multiplication $p_\ast q^!$ endows $\cohaqv^T$ with the structure of an associative $\N I$-graded algebra. Moreover, the multiplication is of degree zero with respect to the vertical grading, i.e., it factors as $m \colon \coha_{\bfd,\, k}^T \otimes \coha_{\bfd',\, k'}^T \to \coha_{\bfd + \bfd',\, k+k'}^T$. In particular, $\coha_0^T$ is a $\N I$-graded subalgebra of $\cohaqv^T$.
\end{theorem}

\subsubsection{Orientations and sign twists}\label{sec:sign twists}

Let $\qv'=(I, \Omega')$ be another orientation of $\qv$ and let $\Gamma \subset \Omega$ be the set of arrows in $\Omega$ which are reversed in $\Omega'$. There is a natural isomorphism
\begin{align}
	\begin{tikzcd}[ampersand replacement=\&, row sep=tiny]
		u_{\qv,\qv'}\colon\stackRep(\Pi_\qv) \ar{r}{\sim} \& \stackRep(\Pi_{\qv'}) \\
		 (x_e)_{e \in \overline{\Omega}}  \ar[mapsto]{r} \&(\epsilon_{e}x_e)_{e \in \overline{\Omega}}
	\end{tikzcd}\ ,
\end{align}
where
\begin{align}
	\epsilon_{e}\coloneqq 
	\begin{cases} 1 & \text{if } e \in \overline{\Omega} \smallsetminus \Gamma\ ,\\ 
		-1 & \text{if } e \in \Gamma\ ,  
	\end{cases}
\end{align}
which extends to the derived enhancements $\dstackRep(\Pi_\qv) \simeq\dstackRep(\Pi_{\qv'})$, and restricts to isomorphisms $\Lambda_\qv \simeq \Lambda_{\qv'}$ and $\dLambda_\qv \simeq \dLambda_{\qv'}$. It is easy to check that this isomorphism is compatible with the diagram \eqref{eq:definition-COHA-product-nil}. It follows that $u_{\qv,\qv'}$ induces an isomorphism of algebras 
\begin{align}
	\begin{tikzcd}[ampersand replacement=\&, row sep=tiny]
		\cohaqv^T \ar{r}{\sim} \& \coha_{\qv'}^T
	\end{tikzcd}\ ,
\end{align}
which we will denote with the same notation $u_{\qv,\qv'}$.

Note that the composition $u_{\qv',\qv}\circ u_{\qv,\qv'}$ is not the identity map at the level of stacks as it reverses the sign of $x_e$ for $e \in \Gamma \cup \Gamma^{\ast}$, but the corresponding isomorphism at the level of COHAs is the identity. In other words, $u_{\qv,\qv'}$ are \textit{canonical} identifications between the COHAs $\cohaqv^T$ for all orientations of the underlying graph of $\qv$. 

In order to relate COHAs with Yangians in \S\ref{sec:relation-Yangian-COHA}, it will however be necessary to twist the multiplication on $\cohaqv^T$ by appropriate signs. 
\begin{definition}
	A \textit{twist} is a bilinear form 
	\begin{align}
		\Theta \colon \Z I \times \Z I \longrightarrow \Z/2\Z
	\end{align} 
	satisfying
	\begin{align}\label{eq:sign-twist-form}
		\Theta(\bfd,\bfe)+\Theta(\bfe,\bfd)=(\bfd,\bfe)\ .
	\end{align}
\end{definition}

We set
\begin{align}\label{eq:def_coha_product-quiver}
	m_{\bfd_1,\bfd_2}^{\Theta}\coloneqq(-1)^{\Theta(\bfd_1,\bfd_2)} p_\ast q^! \colon \coha_{\bfd_1}^T \otimes \coha_{\bfd_2}^T \longrightarrow \coha_{\bfd_1 + \bfd_2}^T\ .
\end{align}
\begin{notation}
	$\cohaqv^{T,\Theta}$ denotes the vector space $\HBMbulletT(\dLambda_\qv)$ endowed with the twisted multiplication $\oplus_{\bfd_1, \bfd_2}m_{\bfd_1,\bfd_2}^{\Theta}$. By abuse of notation, from now on, the untwisted COHA will be denoted by $\cohaqv^{T,\mathsf{untw}}$, while we shall use the notation $\cohaqv^T$ for the COHA whose multiplication is twisted by the Euler form $\Theta=\langle -,- \rangle_\qv$, which is introduced in Formula~\ref{eq:euler-form-quiver}.
\end{notation}

As shown in Proposition~\ref{prop:twist-sign-antisym}, the isomorphism class of $\coha^{T,\Theta}_\qv$ is independent of the choice of $\Theta$. Furthermore, the isomorphism restricts to a sign twist on each weight space $\coha_{\bfd}^T$ and can be chosen to be the identity for $\bfd \in \{\alpha_i\;\vert\;i \in I\}$. By \cite[Theorem~B]{SV_generators}, $\coha_\qv^T$ is generated as an algebra by $\coha^T_{\bfd}$ for $\bfd \in \{\alpha_i\;\vert\;i \in I\}$ and hence any two choices of twistings as above yield \textit{canonically} isomorphic algebras. We obtain the following result.

\begin{proposition}\label{prop:sign-twists-I} 
	There are compatible graded algebra isomorphisms 
	\begin{align}
		\begin{tikzcd}[ampersand replacement=\&, row sep=tiny]
			\gamma_{\qv,\qv'}\colon \cohaqv^T\ar{r}{\sim} \& \coha_{\qv'}^T
		\end{tikzcd}
	\end{align}
	between the COHAs of any two orientations $\qv, \qv'$ of the same underlying graph. These are uniquely characterized by the condition that 
	\begin{align} 
		\gamma_{\qv,\qv'}\vert_{\coha^T_{\alpha_i}}= u_{\qv,\qv'}\vert_{\coha^T_{\alpha_i}}\ ,
	\end{align} 
	for any $i \in I$.
\end{proposition}

\begin{remark}
	Proposition~\ref{prop:sign-twists-I} extends to arbitrary pairs $(\Theta,\Theta')$ of bilinear forms satisfying Equation~\eqref{eq:sign-twist-form}.
\end{remark}

\section{Derived reflection functors and cohomological Hall algebras of quivers}\label{sec:reflection-functors}

In this section, we recall the theory of (derived) reflection functors and describe the relations between it and the theory of nilpotent quiver COHAs.

For simplicity, we shall assume here that the quiver $\qv$ is not a finite ADE quiver\footnote{Again, in the case that $\qv$ is a Dynkin quiver, we may embedd it into a larger quiver.}. 

\subsection{(Derived) reflection functors for preprojective algebras}\label{subsec:reflection-functors}

In this section, we shall define (derived) reflection functors for representations of $\Pi_\qv$. We follow the treatment given in \cite{SY13}. Note that the results we cite from \cite{SY13} have been already proven in \cite{BIRS09} for finite-dimensional modules over the \textit{completed} preprojective algebra $\Lambda_\qv$, i.e., the completion of $\Pi_\qv$ with respect to the augmentation ideal. In \cite{SY13}, the authors reproved some of the results of \cite{BIRS09} to show that they also hold for $\Pi_\qv$. Derived reflection functors were also treated in \cite{BK12,BKT14}.

Fix a vertex $i\in I$. Let $e_i$ be the primitive idempotent of $\Pi_\qv$ attached to the vertex $i\in I$. We define a two-sided ideal $I_i$ of $\Pi_\qv$ by 
\begin{align}
	I_i\coloneqq \Pi_\qv(1-e_i)\Pi_\qv\ .
\end{align}
As $\qv$ has no edge loops, $I_i$ is an ideal of codimension $1$ and there is an exact sequence
\begin{align}
	\begin{tikzcd}[ampersand replacement=\&]
		0 \ar{r} \& I_i \ar{r} \& \Pi_\qv \ar{r} \& \sigma_i \ar{r} \& 0
	\end{tikzcd} \ ,
\end{align}
where $\sigma_i$ is the unique simple $\Pi_\qv$-module of dimension $\alpha_i$.
\begin{theorem}[{\cite[Theorem~2.20, Lemma~2.21, Proposition~2.25, and Theorem~2.26]{SY13}}]\label{thm:reflection-functors}
	Let $w$ be an element of $W_\qv$ with a reduced expression $w=s_{i_1}s_{i_2}\cdots s_{i_r}$. 
	\begin{itemize}\itemsep0.2cm
		\item The multiplication in $\Pi_\qv$ gives rise to an isomorphism of $(\Pi_\qv,\Pi_\qv)$-bimodules 
		\begin{align}
			I_{i_1}\otimes_{\Pi_\qv}^\LL I_{i_2}\otimes_{\Pi_\qv}^\LL\cdots\otimes_{\Pi_\qv}^\LL I_{i_r} \simeq I_{i_1}\otimes_{\Pi_\qv}I_{i_2}\otimes_{\Pi_\qv}\cdots\otimes_{\Pi_\qv}I_{i_r} \simeq I_{i_1}\cdots I_{i_r}\ .
		\end{align}
		
		\item The product $I_{i_1}\cdots I_{i_r}$ depends only on $w$; we can thus denote it by $I_w$. It has finite codimension in $\Pi_\qv$.
		
		\item $I_w$ is a tilting $(\Pi_\qv,\Pi_\qv)$-bimodule of projective dimension at most one and $\End_{\Pi_\qv}(I_w)\simeq \Pi_\qv$.
	\end{itemize}
\end{theorem}
For any $w\in W_\qv$, we define the following subcategories of $\modPi$:\begin{align}
	\scrT_w&\coloneqq \{M\in \modPi\,\vert \,I_w\otimes_{\Pi_\qv}M=0\}\ ,\\[2pt]	
	\scrF_w&\coloneqq \{M\in \modPi\,\vert \,\Tor^1_{\Pi_\qv}(I_w,M)=0\} \ .
\end{align}
and
\begin{align}
	\scrT^{\, w}&\coloneqq\big\{M\in \modPi\,\vert \,\Ext^1_{\Pi_\qv}(I_w,M)=0\}\ ,\\[2pt]	
	\scrF^{\, w}&\coloneqq\big\{M\in \modPi\,\vert \,\Hom_{\Pi_\qv}(I_w,M)=0\}\ .
\end{align}

\begin{remark}\label{rem:reflection-functors}
	Since $I_w$ is of projective dimension one, $\Ext^2_{\Pi_\qv}(I_w,M)=0$ and $\Tor^2_{\Pi_\qv}(I_w,M)=0$ for any module $M$.
\end{remark}

We define the following endofunctors of $\ModPi$:
\begin{align}\label{eq:SS'}
	S_w\coloneqq \Hom_{\Pi_\qv}(I_w,-)\quad\text{and}\quad S'_w\coloneqq I_w\otimes_{\Pi_\qv}- \ .
\end{align}

\begin{proposition}[{\cite[Proposition~2.7]{SY13} and \cite[Theorem~5.4]{BKT14}}]\label{prop:torsion-pairs}
	\hfill
	\begin{itemize}\itemsep0.2cm
		\item The pair $(\scrT^{\, w}, \scrF^{\, w})$ is a torsion pair of $\modPi$. For any finite-dimensional representation $M$ of $\Pi_\qv$, the evaluation map
		\begin{align}
			I_w\otimes_{\Pi_\qv} \Hom_{\Pi_\qv}(I_w, M)\longrightarrow M 
		\end{align} 
		is injective and its image is the torsion object of $M$ with respect to the torsion pair $(\scrT^{\, w}, \scrF^{\, w})$. 
		
		\item The pair $(\scrT_w, \scrF_w)$ is a torsion pair of $\modPi$. For any finite-dimensional representation $M$ of $\Pi_\qv$, the coevaluation map
		\begin{align}
			M\longrightarrow \Hom_{\Pi_\qv}(I_w, I_w\otimes_{\Pi_\qv} M)
		\end{align}
		is surjective and its kernel is the torsion object of $M$ with respect to the torsion pair $(\scrT_w, \scrF_w)$.
		
		\item There are mutually inverse equivalences
		\begin{align}\label{eq:reflection-functors-equivalences-A}
			\begin{tikzcd}[ampersand replacement=\&]
				\scrT^{\, w} \arrow[shift left = 5pt]{r}{S_w}  \& \scrF_w \arrow[shift right = -5pt]{l}{S'_w}
			\end{tikzcd}\ .
		\end{align} 
	\end{itemize}
\end{proposition}

\begin{remark}\label{rem:openness}
	Note that $M\in \scrT^{\, w}$ if and only if $I_w\otimes_{\Lambda_\qv} \Hom_{\Lambda_\qv}(I_w, M)\to M$ is surjective, while $M\in \scrF_w$ if and only if $M\to \Hom_{\Lambda_\qv}(I_w, I_w\otimes_{\Lambda_\qv} M)$ is injective. Thus, $\scrT^{\, w}$ and $\scrF_w$ are open, in the sense of \cite[Definition~\ref*{torsion-pairs-def:open_torsion_pair}]{DPS_Torsion-pairs}, thanks to the semicontinuity theorem.
\end{remark}

\begin{proposition}[{\cite[Proposition~5.7]{BKT14}}]
	Let $u,v,w\in W_\qv$ be such that $\ell(uvw)=\ell(u)+\ell(v)+\ell(w)$. Then, one has mutually inverse equivalences
	\begin{align}\label{eq:reflection-functors-equivalences-B}
		\begin{tikzcd}[ampersand replacement=\&]
			\scrF_u\cap \scrT^{\, vw} \arrow[shift left = 5pt]{r}{S_v}  \& \scrF_{uv}\cap \scrT^{\, w}\arrow[shift right = -5pt]{l}{S'_v}
		\end{tikzcd}\ .
	\end{align} 
\end{proposition}

The ideal $I_w$ gives rise to two quasi-inverse auto-equivalences $\R S_w\coloneqq \R\Hom_{\Pi_\qv}(I_w,-)$ and $\LL S_w\coloneqq I_w\otimes^\LL_{\Pi_\qv}- $, see \cite[Theorem~2.3 and Lemma~2.22]{SY13},
\begin{align}\label{eq:derived-reflection-functors}
	\begin{tikzcd}[ampersand replacement=\&]
		\catDb(\ModPi) \arrow[shift left = 5pt]{r}{\R S_w}  \& \catDb(\ModPi)\arrow[shift right = -5pt]{l}{\LL S_w}\ ,\\
		\catDbps(\ModPi) \arrow[shift left = 5pt]{r}{\R S_w}  \& \catDbps(\ModPi)\arrow[shift right = -5pt]{l}{\LL S_w} \ , \\
		\catDb_{\ps,\,\mathsf{nil}}(\ModPi) \arrow[shift left = 5pt]{r}{\R S_w}  \& \catDb_{\ps,\,\mathsf{nil}}(\ModPi)\arrow[shift right = -5pt]{l}{\LL S_w} \ .
	\end{tikzcd}	
\end{align} 

\begin{remark}
	We have $\catDbps(\ModPi)\simeq \catDb(\modPi)$ and $\catDb_{\ps,\,\mathsf{nil}}(\ModPi)\simeq \catDb(\nilpPi)$. The latter is proved in \cite[Theorem~1.4]{Lewis_stability}.
\end{remark}

\begin{proposition}[{\cite[Proposition~2.27]{SY13}}]\label{prop:derived-reflection-functors-factorization} 
	For $w$ an element of $W_\qv$ with a reduced expression $w=s_{i_r}\cdots s_{i_2} s_{i_1}$, we have
	\begin{align}
		\R S_w \simeq \R S_{i_r} \circ \cdots \circ \R S_{i_2}\circ \R S_{i_1} \quad \text{and}\quad 
		\LL S_w \simeq \LL S_{i_r}\circ \cdots \circ \LL S_{i_2}\circ \LL S_{i_1}\ . 
	\end{align}
\end{proposition}

\begin{remark}
	As a consequence of the previous proposition, it is straightforward to deduce that the assignments $T_i \mapsto \R S_i$ and $T_i^{-1} \mapsto \LL S_i$, for $i \in I$ give rise to a (weak) action of the braid group $B_\qv$ on $\catDb(\ModPi)$. At the K-theoretical level, the corresponding action of $B_\qv$ on the dimension vector $\Z I$ is the standard action of $W_\qv$.
\end{remark}

The following result is a direct consequence of Theorem~\ref{thm:reflection-functors} and Remark~\ref{rem:reflection-functors}, which will be used later on.
\begin{proposition}\label{prop:derived-underived-reflection-functors}
	Let $u,v,w\in W_\qv$ be such that $\ell(uvw)=\ell(u)+\ell(v)+\ell(w)$. Then the derived functor $\R S_v$ restricted to $\scrF_u\cap \scrT^{\, vw}$ coincides with the reflection functor $S_v$, while the derived functor $\LL S_v$ restricted to $\scrF_{uv}\cap \scrT^{\, w}$ coincides with the reflection functor $S'_v$.
\end{proposition}

\subsubsection{An equivalent description of $S_i$}\label{subsec:equivalent-description-S-i}

In this section, we provide an equivalent description of the functors $S_i$ and $S_i'$ for $i\in I$ at the level of the categories of finite-dimensional representations of $\Pi_\qv$ following \cite[\S2.2]{BK12}. First, it is convenient to view a $\Pi_\qv$-module as a pair $(V, x)$ consisting of an $I$-graded vector space
\begin{align}
	V=\bigoplus_{i\in I}V_i
\end{align}
with linear maps $x_e\colon V_i\longrightarrow V_j$ and $x_{e^\ast}\colon V_j\longrightarrow V_i$ for each arrow $e\colon i\to j$ in $\Omega$ satisfying the preprojective relation. We will denote the pair $(V, x)$ simply by $x$. 

For each $i\in I$, we associate to a representation $(V,  x)$ of $\Pi_\qv$ the following diagram 
\begin{align}\label{eq:prepro}
	\begin{tikzcd}[ampersand replacement=\&]
		\widetilde V_i \ar{r}{\ix} \& V_i \ar{r}{x^{(i)}} \& \widetilde  V_i
	\end{tikzcd} \ ,
\end{align}
where
\begin{align}
	\widetilde V_i\coloneqq \bigoplus_{\genfrac{}{}{0pt}{}{e\in \overline{\Omega}}{t(e)=i}} V_{s(e)} \ ,
\end{align}
and
\begin{align}
	\ix\coloneqq \bigoplus_{\genfrac{}{}{0pt}{}{e\in \overline{\Omega}}{t(e)=i}}\varepsilon(e) x_e\quad \text{and}\quad x^{(i)}\coloneqq \bigoplus_{\genfrac{}{}{0pt}{}{e\in \overline{\Omega}}{t(e)= i}} x_{e^\ast}\ .
\end{align}
Here, $\varepsilon(e)=1$ if $e\in\Omega$ and $-1$ else. The preprojective relation at the vertex $i$ is $\ix\circ x^{(i)}=0$. 

We now define two new $\Pi_\qv$-modules $S_i(x)$ and $S'_i(x)$ as follows. The underlying vector spaces for $S_i(x)$ at the vertices $j \neq i$, as well as the maps $V_j \to V_k$ for $j,k \neq i$, coincide with those of $x$. We replace $V_i$ by $\ker\big(\ix\big)$, and define the maps adjacent to the vertex $i$ by replacing Formula~\eqref{eq:prepro} with the canonical diagram
\begin{align}
	\begin{tikzcd}[ampersand replacement=\&]
		\widetilde V_i	\ar{r}{x^{(i)}\circ \ix} \& \ker\big(\ix\big) \ar{r}\& \widetilde V_i  
	\end{tikzcd}\ .
\end{align}
We define $S'_i(x)$ in the same way, replacing $V_i$ by $\mathsf{Coker}\big(x^{(i)}\big)$ and using now the canonical diagram  
\begin{align}
	\begin{tikzcd}[ampersand replacement=\&]
			\widetilde V_i \ar{r}\& \mathsf{Coker}\big(x^{(i)}\big) \ar{r}{x^{(i)}\circ \ix} \& \widetilde V_i
	\end{tikzcd}\ .
\end{align}
It is not difficult to check that the endofunctors 
\begin{align}
	S_i\ , \ S'_i\colon \modPi\longrightarrow \modPi
\end{align}
introduced in \eqref{eq:SS'} coincide with those which associate to $(V, x)$ the representations $S_i(x)$ and $S_i'(x)$, respectively. In particular, $	\scrT^{\, s_i}$ and $\scrF_{s_i}$ are equivalent to the full subcategories of $\modPi$ consisting of the modules such that the maps $\ix$ and $x^{(i)}$ are surjective and injective, respectively. 

\subsection{(Derived) reflection functors and COHAs}

For any $v, w\in W_\qv$, let $\stackRep(\Pi_\qv)_w^v$ be the open substack of $\stackRep(\Pi_\qv)$ parameterizing finite-dimensional representations of $\Pi_\qv$, which belong to $\scrT^{\, v}\cap \scrF_w$. Let $\dstackRep(\Pi_\qv)_w^v$ be its canonical derived enhancement. Define
\begin{align}\label{eq:def-of-Lambda-uv}
	\Lambda_w^v\coloneqq \dLambda_\qv\times_{\stackRep(\Pi_\qv)} \stackRep(\Pi_\qv)_w^v\quad\text{and}\quad
	\dLambda_w^v\coloneqq \dLambda_\qv\times_{\dstackRep(\Pi_\qv)} \dstackRep(\Pi_\qv)_w^v\ .
\end{align}
In addition, set
\begin{align}
	\iLambda_\qv\coloneqq \Lambda^{s_i}\quad &\text{and}\quad \Lambdai_\qv \coloneqq \Lambda_{s_i}\ ,\\[2pt]
	\idLambda_\qv\coloneqq \dLambda^{s_i}\quad &\text{and}\quad \dLambdai_\qv \coloneqq \dLambda_{s_i}
\end{align}
for any $i\in I$. The derived equivalences \eqref{eq:reflection-functors-equivalences-B} induces mutually inverse auto-equivalences of $\dstackRep(\Pi_\qv)$. Thanks to Proposition~\ref{prop:derived-underived-reflection-functors}, the equivalences  \eqref{eq:derived-reflection-functors} induces mutually inverse equivalences
\begin{align}
	\begin{tikzcd}[ampersand replacement=\&]
		\dstackRep(\Pi_\qv)_u^{vw} \arrow[shift left = 5pt]{r}{\R S_v}  \& \dstackRep(\Pi_\qv)_{uv}^w\arrow[shift right = -5pt]{l}{\LL S_v}
	\end{tikzcd}\ ,
\end{align} 
hence also equivalences
\begin{align}\label{eq:equivalences-stacks}
	\begin{tikzcd}[ampersand replacement=\&]
		\dLambda_u^{vw} \arrow[shift left = 5pt]{r}{\R S_v}  \& \dLambda_{uv}^w\arrow[shift right = -5pt]{l}{\LL S_v}
	\end{tikzcd}\ .
\end{align} 
Note that the isomorphism of classical stacks, which are open in their respective full stack, lifts automatically to the derived level by using the construction of the canonical derived enhancement of \cite[Proposition~2.1]{STV}.

For $w, v\in W_\qv$, set
\begin{align}
	\wcohaqv\coloneqq \HBMbulletT( \dLambda^w )\ , \quad \cohaqvw\coloneqq \HBMbulletT( \dLambda_w )\ , \quad\text{and}\quad
	\cohavw\coloneqq \HBMbulletT( \dLambda^v_w )\ .
\end{align}
The stack equivalences \eqref{eq:equivalences-stacks} yield mutually inverse isomorphisms
\begin{align}\label{eq:isomorphism-S-w}
	\begin{tikzcd}[ampersand replacement=\&]
		\wcohaqv \arrow[shift left = 5pt]{r}{(S_w)_\ast}  \& \cohaqvw \arrow[shift right = -5pt]{l}{(S'_w)_\ast}
	\end{tikzcd}	
\end{align}
and
\begin{align}\label{eq:isomorphism-S-uvw}
	\begin{tikzcd}[ampersand replacement=\&]
		\tensor*[^{(vw)}]{\coha}{^{T, (u)}_\qv} \arrow[shift left = 5pt]{r}{(S_v)_\ast}  \& \tensor*[^{(w)}]{\coha}{^{T, (uv)}_\qv} \arrow[shift right = -5pt]{l}{(S'_v)_\ast}
	\end{tikzcd}\ .	
\end{align}
When $w=s_i$ with $i\in I$, we set
\begin{align}
	\icohaqv\coloneqq \wcohaqv\quad\text{and} \quad \cohaiqv\coloneqq \cohaqvw\ , 
\end{align}
and $S_{i,\, \ast }\coloneqq (S_{s_i})_\ast$ and $S'_{i,\, \ast }\coloneqq (S'_{s_i})_\ast$.

\begin{proposition}\label{prop:reflection-stacks}
	\hfill
	\begin{enumerate}\itemsep0.2cm
		\item \label{item:reflection-stacks-1} There are right and left $\cohaqv^T$-modules structures on $\icohaqv$ and $\cohaiqv$, respectively, such that both restriction maps $\ires\colon \cohaqv^T\to\icohaqv$ and $\resi\colon \cohaqv^T\to\cohaiqv$ are $\cohaqv^T$-module homomorphisms.
		\item  \label{item:reflection-stacks-2} The mixed Hodge structures on both $\icohaqv$ and $\cohaiqv$ are pure and the maps $\ires$ and $\resi$ are surjective.
		\item  \label{item:reflection-stacks-3} $\ker\big(\resid\colon \cohaqvd^T\to\cohaid\big)=\coha_{\bfd -\alpha_i}^T\ast \coha_{\alpha_i}^T$ and $\ker\big(\iresd\colon \cohaqvd^T\to\icohad\big)=\coha_{\alpha_i}^T\ast \coha_{\bfd-\alpha_i}^T$.
		\item \label{item:reflection-stacks-4} The graded dimensions of $\icohaqv$ and $\cohaiqv$ are equal and given by
		\begin{align}
			\dim \icohaqv = \Exp\Big( \frac{1}{1-q^{-1}}\sum_{\bfd\in \N I\smallsetminus \N \alpha_i} A_{\bfd}(q^{-1}) z^\bfd \Big) = \dim \cohaiqv\ ,
		\end{align}
		where $z^\bfd\coloneqq \prod_{i\in I} z_i^{d_i}$ and $A_{\bfd}$ is the Kac's polynomial of the quiver $\qv$ with dimension vector $\bfd$.
	\end{enumerate}
\end{proposition}

\begin{proof}
	Since $\scrT^{s_i}$ is the torsion part of a torsion pair, it is closed under quotients, while since $\scrF_{s_i}$ is the torsion-free part of a torsion pair, it is closed under subobjects. This implies that upon the open base change, the convolution diagram \eqref{eq:definition-COHA-product-nil} gives rise to diagrams
	\begin{align}
		\begin{tikzcd}[ampersand replacement=\&]
			\idLambda_\qv \times \dLambda_\qv \& \dLambdaext_\qv\times_{\dLambda_\qv} \idLambda_\qv \ar{r} \ar{l} \& \idLambda_\qv  
		\end{tikzcd}
	\end{align}
	and
	\begin{align}
		\begin{tikzcd}[ampersand replacement=\&]
			\dLambda_\qv \times \dLambdai_\qv \& \dLambdaext_\qv \times_{\dLambda_\qv} \dLambdai_\qv \ar{r} \ar{l} \& \dLambdai_\qv 
		\end{tikzcd}\ .
	\end{align}
	From this, \eqref{item:reflection-stacks-1} is a consequence of standard open base change properties.
	
	\medskip
	
	We now prove \eqref{item:reflection-stacks-2}. Fix a dimension vector $\bfd\in \N I$. We shall concentrate on $\cohaid$, since the other case is similar. By \cite[Lemma~2.23-(4)]{SY13}, a finite-dimensional $\Pi_\qv$ module belongs to $\scrT^{\, s_i}$ if and only if it does not contain any submodule supported at the vertex $i$. Fix a stability condition $\Ltheta\in \coweightlattice$ such that $\Ltheta_i=1$ and $\Ltheta_j< 1$ for each $j\neq i$. Then, we also get that a finite-dimensional $\Pi_\qv$ module $M$ belongs to $\scrT^{s_i}$ if and only $M$ does not have any submodule of $\Ltheta$-slope equal to $1$. Equivalently, the maximal $\Ltheta$-slope of $M$ is strictly less than $1$. Thus, the classical stack $\Lambdai_\bfd$ is a union of the (finitely many) Harder-Narasimhan strata with $\mu_{\Ltheta\textrm{-}\mathsf{max}}<1$. Likewise, the closed complement $\Lambda_\bfd \smallsetminus \Lambdai_\bfd$ is the union of all Harder-Narasimhan strata with $\mu_{\Ltheta\textrm{-}\mathsf{max}}=1$. 
	
	Since the category of finite-dimensional representations of $\Pi_\qv$ is 2-Calabi-Yau, each Harder-Narasimhan stratum is a vector bundle stack (in the sense of \cite[Definition~\ref*{torsion-pairs-def:vector-bundle-stack}]{DPS_Torsion-pairs}) over a product of stacks of $\Ltheta$-semistable nilpotent representations of $\Pi_\qv$ of a fixed dimension vector. Now, By \cite[Theorem~B]{Davison_Integrality} the mixed Hodge structure of the $T$-equivariant Borel-Moore homology of $\Lambda_\bfd^{\Ltheta\textrm{-}\mathsf{ss}}$ is a direct summand of $\cohaqvd$ and the latter is cohomologically pure by \cite[Theorem~A]{SV_generators} or \cite[Theorem~A]{Davison_Integrality}. This proves that both $\Lambdai_\bfd$ and $\Lambda_\bfd \smallsetminus \Lambdai_\bfd$ have a pure mixed Hodge structure, from which we deduce \eqref{item:reflection-stacks-2}. 
	
	\medskip
	
	We next prove the first statement in \eqref{item:reflection-stacks-3} (the proof of the second one is similar). The kernel of the restriction $\resid\colon \cohaqvd^T\to\cohaid$ is the image of the pushforward map $\HBMbulletT(\Lambda_\bfd\smallsetminus \Lambdai_\bfd)\to \cohaqvd^T$. 
	
	Because of \cite[Lemma~2.23-(4)]{SY13}, the map
	\begin{align}
		\trunc{p_{\bfd-\alpha_i,\, \alpha_i}}\colon \trunc{\dLambdaext_{\alpha_i,\bfd -\alpha_i}} \longrightarrow \Lambda_\bfd
	\end{align}
	factors via the inclusion of the stack $\Lambda_\bfd\smallsetminus \Lambda_{\bfd}^{(i)}$ into $\Lambda_\bfd$. Thus, the multiplication $m_{\bfd-\alpha_i,\, \alpha_i}=(p_{\bfd-\alpha_i,\, \alpha_i})_\ast\circ q_{\bfd-\alpha_i,\, \alpha_i}^!$ factors as
	\begin{align}
		\begin{tikzcd}[ampersand replacement=\&]
			\coha_{\bfd-\alpha_i}^T\otimes\coha_{\alpha_i}^T\ar{r}{\nu}\& \HBMbulletT( \Lambda_{\bfd}\smallsetminus \Lambda_{\bfd}^{(i)} )\ar[r]	\& \cohaqvd^T
		\end{tikzcd}\ .
	\end{align}
	Thus we are reduced to proving that the map $\nu$ is surjective. The following argument is similar to that used in \cite[\S5.8]{SV_generators}.
	
	Fix $i\in I$ and $\ell\in\N$. For any dimension vector $\bfd$, define the derived stack $\dLambda_{(i, \ell),\, \bfd}$ by the pullback square
	\begin{align}\label{eq:pull-back-i-l}
		\begin{tikzcd}[ampersand replacement = \&]
			\dLambda_{(i, \ell),\, \bfd}\arrow{r}{j} \arrow[swap]{d}{q^{(i)}} \& \dLambdaext_{\bfd-\ell\alpha_i,\ell\alpha_i} \arrow{d}{q} \\
			\dLambdai_{\bfd-\ell\alpha_i} \times \dLambda_{\ell\alpha_i}\arrow{r} \& \dLambda_{\bfd-\ell\alpha_i}\times  \dLambda_{\ell\alpha_i}
		\end{tikzcd} \ .
	\end{align}
	For any short exact sequence of nilpotent representations of the form $0\to N\to E\to M\to 0$ such that $N$ is of the form $\sigma_i^{\oplus \ell}$ for some $\ell\in \N$ and $M$ belongs to $\scrF_{s_i}$, $N$ is the torsion part of $E$ and $M$ is the torsion-free part of $E$ with respect to the torsion pair $(\scrT_{s_i}, \scrF_{s_i})$ of $\modPi$. In particular, such a short exact sequence is canonical, and $N$ is the unique submodule of $E$ of dimension vector $\ell\alpha_i$. Now, $\Hom_{\Pi_\qv}(\sigma_i, M)=0$ by the defining properties of the torsion pair $(\scrT_{s_i}, \scrF_{s_i})$. The latter and the 2-Calabi-Yau property of the category of finite-dimensional representations of $\Pi_\qv$ imply that the map $q^{(i)}$ is a vector bundle stack morphism. Since $\dLambda_{\bfd-\ell\alpha_i}^{(i)}\times \dLambda_{\ell\alpha_i}$ is cohomologically pure, $\dLambda_{(i, \ell),\,\bfd}$ is cohomologically pure as well.
	
	In addition, the analysis above yields that the composition of the top horizontal map $j$ in \eqref{eq:pull-back-i-l} with the map $p\colon \dLambdaext_{\bfd-\ell\alpha_i,\ell\alpha_i} \to \Lambda_{\bfd}$ is the inclusion of a locally closed geometric substacks. In particular, $\dLambda_{(i, \ell),\, \bfd}$ parametrizes $\bfd$-dimensional nilpotent representations $M$ of $\Pi_\qv$ such that $\dim \Hom(\sigma_i, M)=\ell$. Let $\Lambda_{(i, \ell),\, \bfd}$ denote the truncation of $\dLambda_{(i, \ell),\, \bfd}$. Therefore,
	\begin{align}
		\Lambda_\bfd\smallsetminus \Lambda_{\bfd}^{(i)}=\bigsqcup_{\ell>0}\, \Lambda_{(i, \ell),\,\bfd}\ .
	\end{align}
	We also define
	\begin{align}
		\Lambda_{(i, >\ell),\, \bfd}\coloneqq \bigsqcup_{\ell'>\ell}\, \Lambda_{(i, \ell'),\, \bfd}\quad\text{and}\quad
		\Lambda_{(i, \geqslant\ell),\, \bfd}\coloneqq \bigsqcup_{\ell'\geqslant\ell}\, \Lambda_{(i, \ell'),\, \bfd}\ ,\\[2pt]
		\dLambda_{(i, >\ell),\, \bfd}\coloneqq \bigsqcup_{\ell'>\ell}\, \dLambda_{(i, \ell'),\, \bfd}\quad\text{and}\quad
		\dLambda_{(i, \geqslant\ell),\, \bfd}\coloneqq \bigsqcup_{\ell'\geqslant\ell}\, \dLambda_{(i, \ell'),\, \bfd}\ .
	\end{align}
	Set
	\begin{align}
		\cohaqvd^{T, (i, \ell)}\coloneqq \HBMbulletT( \dLambda_{(i, \ell),\, \bfd} )\ , \quad \cohaqvd^{T, (i, >\ell)}\coloneqq \HBMbulletT( \dLambda_{(i, >\ell),\,\bfd} )\ ,  \quad \cohaqvd^{T, (i, \geqslant\ell)}\coloneqq \HBMbulletT( \dLambda_{(i, \geqslant\ell),\,\bfd} )\ .
	\end{align}
	By cohomological purity, we have the short exact sequence
	\begin{align}\label{eq:short-exact-sequence}
		0\longrightarrow \cohaqvd^{T, (i, >\ell)} \longrightarrow \cohaqvd^{T, (i, \geqslant\ell)} \longrightarrow \cohaqvd^{T, (i, \ell)} \longrightarrow 0 \ .
	\end{align}
	We now prove by descending induction on $\ell$ that $\cohaqvd^{T, (i, >\ell)}$ is contained in the image of the map $\nu$. 
		Assume that  $\cohaqvd^{T, (i,>\ell)}\subset \mathsf{Im}(\nu)$ for some $\ell>0$. By Formula~\eqref{eq:short-exact-sequence}, to prove that $\cohaqvd^{T, (i, \geqslant \ell)}\subset \mathsf{Im}(\nu)$ we consider the following diagram
	\begin{align}
		\begin{tikzcd}[column sep=small, ampersand replacement = \&]
			\dLambda_{\bfd-\ell\alpha_i}^{(i)} \times \dLambda_{\ell\alpha_i}\arrow{d}{u} \& \dLambda_{(i, \ell),\, \bfd}\arrow{d}{j} \arrow[swap]{l}{q^{(i)}} \arrow{r}{p^{(i)}} \& \dLambda_{(i, \ell),\, \bfd}\arrow{d}{u'}\\
			\dLambda_{\bfd-\ell\alpha_i} \times \dLambda_{\ell\alpha_i} \&  \dLambdaext_{\bfd-\ell\alpha_i,\ell\alpha_i} \arrow[swap]{l}{q}  \arrow{r}{p} \& \dLambda_{(i, \geqslant \ell),\,\bfd}
		\end{tikzcd} \ ,
	\end{align}
	where both squares are cartesian. Note that by construction $q^{(i)}$ is an equivalence. Then, the composition
	\begin{align}
		(p^{(i)})_\ast\circ (q^{(i)})^!\colon \cohai_{\bfd-\ell\alpha_i}\otimes \coha_{\ell\alpha_i}^T\longrightarrow \cohaqvd^{T, (i, \ell)}
	\end{align}
	is an isomorphism. Since the map
	\begin{align}
		u^\ast \colon \coha_{\bfd-\ell\alpha_i}^T\otimes \coha_{\ell \alpha_i}^T\longrightarrow \cohai_{\bfd-\ell\alpha_i}\otimes \coha_{\ell\alpha_i}^T 
	\end{align}
	is surjective by purity, we deduce that the map	$(u')^\ast\circ p_\ast\circ q^!$ is also surjective.
	
	\medskip
	
	Finally, we prove \eqref{item:reflection-stacks-4}. By \cite[Corollary~2.2]{SV_generators}, since $\idLambda_\bfd$ and $\dLambdai_\bfd$ are pure by \eqref{item:reflection-stacks-2} and have polynomial count (cf.\ \cite[\S5]{BSV_Nilpotent}), the computation of the graded dimensions of $\cohaid$ and $\icohad$ reduces to a finite field points count: this is performed in \textit{loc.cit} (see in particular, Corollary 5.2 of \textit{loc.cit.}).
\end{proof}

\section{Yangians}\label{sec:Yangians}

In this section, we introduce an extended version of a Yangian associated to $\qv$, as an algebra given by some explicit set of generators and relations.

\medskip

Recall the torus
\begin{align}\label{eq:Tmax}
	\Tmax=\Big\{(\gamma_e)_{e\in \overline{\Omega}}\in (\C^\ast)^{\doubleOmega}\, \Big\vert\, \gamma_e\gamma_{e^\ast}=\gamma_f\gamma_{f^\ast} \quad\forall \, e, f\in \Omega \Big\} \ .
\end{align}
Let $\varepsilon_e$ and $\hbar$ be the first Chern classes of the tautological characters $\gamma_e$ and $\gamma_e\gamma_{e^\ast}$ of the group $\Tmax$. Then
\begin{align}
	\Hbullet_{\Tmax}\coloneqq H^\ast_{\Tmax}(\mathsf{pt})=\Q\Big[\varepsilon_e\, \Big\vert\, e\in \doubleOmega\Big]\Big/\Big\{\varepsilon_e+\varepsilon_{e^\ast}=\hbar \quad \forall\, e\in \Omega\Big\}\ .
\end{align}

\medskip

In what follows, we use the standard notation for (anti) Lie brackets: set $\{a, b\}\coloneqq ab+ba$ and $[a,b]\coloneqq ab-ba$. 

\subsection{Definition of the Yangian}\label{subsec:def-Yangians}

For any $i, j\in I$, with $i\neq j$, set
\begin{align}\label{eq:zeta}
	\zeta_{i,j}(t)\coloneqq \prod_{\genfrac{}{}{0pt}{}{e\in \doubleOmega}{e\colon i\to j}} (t+\varepsilon_e)\ .
\end{align}

The following ``functional equation'' holds:
\begin{align}\label{eq:functional-equation}
	\zeta_{i,j}(t)=(-1)^{a_{i, j}} \zeta_{j,i}(-t-\hbar)\ .
\end{align}

\begin{definition}\label{def:Yangian}
	The \textit{extended (multi-parameter) Yangian $\eY_\qv$ of $\qv$} is the unital associative $\Hbullet_{\Tmax}$-algebra generated by $x_{i, \ell}^\pm, h_{i, \ell}, \kappa_{i, \ell}$, with $i \in I$ and $\ell \in \N$, subject to the following relations:
	\begin{itemize}\itemsep0.2cm
		\item for any $i\in I$ and $\ell\in \N$
		\begin{align}\label{eq:Yangian-Lie-algebra-1} 
			\kappa_{i, \ell}&\in Z\big(\eY_\qv\big) \ ,
		\end{align}
		
		\item for any $i, j\in I$ and $r,s\in \N$
		\begin{align}\label{eq:Yangian-Lie-algebra-2}
			\Big[h_{i,r}, h_{j,s}\Big]  & = 0\ ,\\[4pt]  \label{eq:Yangian-Lie-algebra-3}
			\Big[x_{i,r}^{+}, x_{j,s}^{-}\Big] &= \delta_{i,j} h_{i, r+s} \ ,
		\end{align}
		
		\item for $i\in I$
		\begin{align}\label{eq:Yangian-4}
			(u-w\mp\hbar)h_i(u)x_i^\pm(w)&\almostsame (u-w\pm\hbar)x_i^\pm(w)h_i(u)\ ,\\[4pt]  \label{eq:Yangian-5}
			(u-w\mp\hbar)x_i^\pm(u)x_i^\pm(w)&\almostsame (u-w\pm\hbar)x_i^\pm(w)x_i^\pm(u)\ ,
		\end{align}
		
		\item for $i, j\in I$, with $i\neq j$
		\begin{align}\label{eq:Yangian-6}
			\zeta_{i,j}(u-w) h_i(u)x_j^+(w)&\almostsame  \zeta_{i, j}(u-w-\hbar)x_j^+(w)h_i(u)\ ,\\[4pt]\label{eq:Yangian-7}
			\zeta_{i,j}(u-w-\hbar) h_i(u)x_j^-(w)&\almostsame \zeta_{i, j}(u-w)x_j^-(w)h_i(u) \ ,\\[4pt]  \label{eq:Yangian-8}
			\zeta_{i,j}(u-w) x_i^+(u)x_j^+(w)&\almostsame \zeta_{i, j}(u-w-\hbar)x_j^+(w)x_i^+(u) \ ,\\[4pt]\label{eq:Yangian-9}
			\zeta_{i,j}(u-w-\hbar) x_i^-(u)x_j^-(w)&\almostsame  \zeta_{i,j}(u-w)x_j^-(w)x_i^-(u)\ ,
		\end{align}	
		
		\item Serre relations: 
		\begin{align}\label{eq:Yangian-Serre}
			\sum_{\sigma \in \frakS_m}\Big[x_i^{\pm}(u_{\sigma(1)}), \Big[x_i^{\pm}(u_{\sigma(2)}), \Big[\cdots, \Big[x_i^{\pm}(u_{\sigma(m)}), x_j^{\pm}(w)\Big]\cdots\Big]\Big]\Big] = 0 
		\end{align}
		for $i, j\in I$, with $i \neq j$, where $m\coloneqq 1 - a_{i,j}$ and $\frakS_m$ denotes the $m$-th symmetric group.
		
		\item for any $i, j\in I$, with $i\neq j$, and any $e\colon i\to j\in \doubleOmega$
		\begin{multline}\label{eq:Yangian-cubic}
			\frac{\zeta_{j,i}(w-u)\zeta_{j,i}(w-v)}{(v-w-\hbar+\varepsilon_e)}x_i^-(u)x_i^-(v)x_j^-(w)\\
			+\frac{(u-v-\hbar)\zeta_{i,j}(u-w)\zeta_{j,i}(w-v)}{(v-w-\hbar+\varepsilon_e)(u-w+\varepsilon_e)}x_i^-(v)x_j^-(w)x_i^-(u)\\
			+\shoveright{\frac{\zeta_{i,j}(v-w)\zeta_{i,j}(u-w)}{(u-w+\varepsilon_e)}x_j^-(w)x_i^-(u)x_i^-(v)\almostsame 0 \ ,}\\[4pt]
			\shoveleft{\frac{\zeta_{j,i}(w-u)\zeta_{j,i}(w-v)}{(v-w-\hbar+\varepsilon_e)}x_j^+(w)x_i^+(v)x_i^+(u)}\\
			+\frac{(u-v-\hbar)\zeta_{i,j}(u-w)\zeta_{j,i}(w-v)}{(v-w-\hbar+\varepsilon_e)(u-w+\varepsilon_e)}x_i^+(u)x_j^+(w)x_i^+(v)\\
			+\frac{\zeta_{i,j}(v-w)\zeta_{i,j}(u-w)}{(u-w+\varepsilon_e)}x_i^+(v)x_i^+(u)x_j^+(w)\almostsame 0 \ ,
		\end{multline}
		
	\end{itemize}
	Here, 
	\begin{align}
		h_i(u)\coloneqq 1+ \hbar\sum_{\ell\in \N} h_{i, \ell} u^{-\ell-1} \quad\text{and}\quad x_i^\pm(u) \coloneqq \sum_{\ell\in \N} x_{i, \ell}^\pm u^{-\ell-1} \in \eY_\qv\llbracket u^{-1}\rrbracket \ ,
	\end{align}	
	and  for $A(u_1, \ldots, u_s), B(u_1, \ldots, u_s)\in \eY_\qv\llbracket u_1^{\pm 1}, \ldots, u_s^{\pm 1}\rrbracket$ the equality $A(u_1, \ldots, u_s)\almostsame B(u_1, \ldots, u_s)$ means an equality of the coefficient of each monomial involving strictly negative powers of the variables associated to \textit{associated to} $x_{i, \ell}^\pm$ for $i\in I$ and $\ell\in \N$.
	
	The \textit{(multi-parameter) Yangian $\Y_\qv$ of $\qv$} is the unital associative algebra over $\Hbullet_{\Tmax}$, which is the quotient of $\eY_\qv$ by the two-sided ideal generated by the $\kappa_{i, \ell}$'s, with $i \in I$ and $\ell \in \N$.
\end{definition}

\begin{remark}
	The first definition of the Yangian of a quiver without edge-loops is as an algebra over $\Q[\hbar]$ (e.g. \cite{Drinfeld_Yangians_1,Drinfeld_Yangians_2,Varagnolo_Yangian,Guay07, GRW19}). In our terminology, it would be the one-parameter Yangian. Recently, also the two-parameters Yangian has been introduced, when the quiver is an affine type A quiver, as an algebra over $\Q[\varepsilon_1, \varepsilon_2]$, see \cite{Kodera_Fock19, BT-Yangians, Kodera21}. Since we are interested in the relation between Yangians and two-dimensional cohomological Hall algebras of quivers, for us it is more natural to consider the largest torus $\Tmax$ which acts on representations of the double quiver and from this point of view define the Yangian as an algebra over $\Hbullet_{\Tmax}$. We shall consider the two-parameters Yangian when the quiver is an affine ADE quiver in \S\ref{subsec:Yangians-affine}.
	
	The inclusion of the additional set of central elements $\{\kappa_{i,\ell}\}$ in the definition of the extended Yangian is motivated by an analogous construction for the Maulik-Okounkov Yangian (see e.g. \cite[\S5.3.4 and \S6.1.1]{MO_Yangian}). From their perspective, our \textit{extended} Yangian corresponds to the Yangian and our Yangian would be closer to their \textit{core} Yangian.
\end{remark}

\begin{remark}\label{rem:specialized-relations} 
	The set of equations above recovers the usual relations of Yangians in particular cases, as we shall show below.
	
	First, note that \eqref{eq:Yangian-Lie-algebra-3} is equivalent to
	\begin{align}\label{eq:Yangian-Lie-algebra-3-alternative-formulation}
		\hbar\Big[x_i^+(u), x_j^-(w)\Big]=\delta_{i, j}\frac{h_i(w)-h_i(u)}{u-w}\ .
	\end{align}
	
	Let $r,s\in \N$. By taking the $u^{-r-1}w^{-s-1}$ coefficient of \eqref{eq:Yangian-4} and \eqref{eq:Yangian-5} for $i\in I$, we get
	\begin{align}\label{eq:Yangian-4-alternative-formulation}
		\Big[h_{i, r+1}, x_{i, s}^{\pm}\Big] - \Big[h_{i, r}, x_{i, s+1}^{\pm}\Big]& = \pm \hbar \Big\{h_{i, r}, x_{i, s}^{\pm}\Big\}\ ,\\[4pt] \label{eq:Yangian-5-alternative-formulation}
		\Big[x_{i, r+1}^{\pm}, x_{i, s}^{\pm}\Big] - \Big[x_{i, r}^{\pm}, x_{i, s+1}^{\pm}\Big] &= \pm \hbar \Big\{x_{i, r}^{\pm}, x_{i, s}^{\pm}\Big\} \ ,
	\end{align}
	respectively.
	
	Let $i, j\in I$, with $i\neq j$. Let $r,s\in \N$. If $a_{i,j}=0$, the $u^{-r-1}w^{-s-1}$ coefficients of \eqref{eq:Yangian-6} -- \eqref{eq:Yangian-9} are	
	\begin{align}
		\Big[h_{i, r+1}, x_{j, s}^{\pm}\Big] - \Big[h_{i, r}, x_{j, s+1}^{\pm}\Big]& = 0\ ,\\[4pt] 
		\Big[x_{i, r+1}^{\pm}, x_{j, s}^{\pm}\Big] - \Big[x_{i, r}^{\pm}, x_{j, s+1}^{\pm}\Big] &= 0 \ .
	\end{align}
	
	If there exists a unique edge $e$ between $i$ and $j$ in $\Omega$, the $u^{-r-1}w^{-s-1}$ coefficients of \eqref{eq:Yangian-6} -- \eqref{eq:Yangian-9} are
	\begin{align}
		\Big[h_{i, r+1}, x_{j, s}^{\pm}\Big] - \Big[h_{i, r}, x_{j, s+1}^{\pm}\Big]& = \mp \frac{\hbar}{2} \Big\{h_{i, r}, x_{j, s}^{\pm}\Big\}-m_{i,j}\frac{\varepsilon_e-\varepsilon_{e^\ast}}{2}\Big[h_{i,r},x_{j,s}^{\pm}\Big]\ ,\\[4pt] 
		\Big[x_{i, r+1}^{\pm}, x_{j, s}^{\pm}\Big] - \Big[x_{i, r}^{\pm}, x_{j, s+1}^{\pm}\Big] &= \mp \frac{\hbar}{2} \Big\{x_{i, r}^{\pm}, x_{j, s}^{\pm}\Big\}-m_{i,j}\frac{\varepsilon_e-\varepsilon_{e^\ast}}{2}\Big[x_{i,r}^{\pm},x_{j,s}^{\pm}\Big] \ ,
	\end{align}
	where
	\begin{align}
		m_{i,j}\coloneqq \begin{cases}
			1 & \text{if } i\to j\in \Omega\ , \\
			-1 & \text{if } j\to i \in \Omega \ .
		\end{cases}
	\end{align}
	
	If there are two edges $e_1, e_2$ from $i$ to $j$, the $u^{-r-1}w^{-s-1}$ coefficients of \eqref{eq:Yangian-6} -- \eqref{eq:Yangian-9} are
	\begin{multline}
		\Big[h_{i,r+2}, x_{j,s}^{\pm}\Big] - 2\Big[h_{i,r+1}, x_{j,s+1}^{\pm}\Big] + \Big[h_{i,r}, x_{j,s+2}^{\pm}\Big] =\\[4pt]
		\shoveright{\mp(\varepsilon_{e_1}+\varepsilon_{e_2})h_{i, r+1} x_{j,s}^\pm \mp(\varepsilon_{e_1^\ast}+\varepsilon_{e_2^\ast})x_{j,s}^\pm h_{i, r+1} \pm (\varepsilon_{e_1}+\varepsilon_{e_2}) h_{i,r}x_{j, s+1}^\pm \pm (\varepsilon_{e_1^\ast}+\varepsilon_{e_2^\ast}) x_{j, s+1}^\pm h_{i,r}}\\[4pt]
		\shoveright{\mp \varepsilon_{e_1}\varepsilon_{e_2} h_{i, r} x_{j,s}^\pm \pm \varepsilon_{e_1^\ast}\varepsilon_{e_2^\ast} x_{j,s}^\pm h_{i, r}\ ,}\\[6pt]
		\shoveleft{\Big[x_{i,r+2}^\pm, x_{j,s}^{\pm}\Big] - 2\Big[x_{i,r+1}^\pm, x_{j,s+1}^{\pm}\Big] + \Big[x_{i,r}^\pm, x_{j,s+2}^{\pm}\Big] =}\\[4pt]
		\shoveright{\mp(\varepsilon_{e_1}+\varepsilon_{e_2})x_{i, r+1}^\pm x_{j,s}^\pm \mp(\varepsilon_{e_1^\ast}+\varepsilon_{e_2^\ast})x_{j,s}^\pm x_{i, r+1}^\pm \pm (\varepsilon_{e_1}+\varepsilon_{e_2}) x_{i,r}^\pm x_{j, s+1}^\pm \pm (\varepsilon_{e_1^\ast}+\varepsilon_{e_2^\ast}) x_{j, s+1}^\pm x_{i,r}^\pm}\\[4pt]
		\mp \varepsilon_{e_1}\varepsilon_{e_2} x_{i, r}^\pm x_{j,s}^\pm \pm \varepsilon_{e_1^\ast}\varepsilon_{e_2^\ast} x_{j,s}^\pm x_{i, r}^\pm \ .
	\end{multline}
	These relations simplify to the standard ones appearing in the affine Yangian of $\widehat{\fraksl}_2$ under the specialization $\varepsilon_{e_1}+\varepsilon_{e_2}=\hbar$ (cf.\ \cite[Definition~5.1]{Kodera_Fock19} and \cite[\S1.2]{BT-Yangians}).
	
	Finally, note that \eqref{eq:Yangian-6} and \eqref{eq:Yangian-7} yield a relation involving only the generators $h_{i,0}, \ldots, h_{i,k}$ and $x^{\pm}_{j,r}, \ldots, x^{\pm}_{j,r+k}$ for any $r$ and any $k<\vert a_{i,j}\vert$. For instance, when $k=0$, we get
	\begin{align}\label{eq:Yangian-6-alternative-formulation}
		\Big[h_{i,0}, x_{j,r}^\pm\Big]=\pm a_{i,j}x_{j,r}^\pm\ ,
	\end{align}
	while when $k=1$, we have
	\begin{align}
		\Big[h_{i,1},x_{j,r}^\pm\Big]=\mp 2\begin{pmatrix} \vert a_{i,j}\vert\\ 2\end{pmatrix} x_{j,r+1}^\pm \mp \frac{\vert a_{i,j}\vert}{2}\hbar\Big\{h_{i,0},x_{j,r}^\pm\Big\}\pm\sum_{\genfrac{}{}{0pt}{}{e\in \doubleOmega}{e\colon i\to j}} \Big(\varepsilon_e-\frac{\hbar}{2}\Big)x_{j,r}^\pm \ .
	\end{align} 
	Again, in the case $\vert a_{i,j}\vert=2$ and $\varepsilon_{e_1}+\varepsilon_{e_2}=\hbar$ we recover the usual relation
	\begin{align}\label{eq:missing-relation}
		\Big[h_{i,1},x_{j,r}^\pm\Big]=\mp 2x_{j,r+1}^\pm\mp\hbar\Big\{h_{i,0},x_{j,r}^\pm\Big\}
	\end{align}
	appearing in \cite[\S1.2]{BT-Yangians}.
\end{remark}

\begin{remark} 
	The cubic relations \eqref{eq:Yangian-cubic} are additive analogues of similar relations appearing in \cite[\S3]{NSS_KHA}. In \textit{loc.cit.}, it was shown that together with the usual quadratic relations, these provide a presentation of the (localized) K-theoretical Hall algebra of an arbitrary quiver $\qv$, for a \textit{full} set of equivariant parameters. Although we do expect a similar result to hold in the cohomological Hall algebra context, we also imposed the usual Serre relation. As an illustration of this, we prove below that the Serre relation~\eqref{eq:Yangian-Serre} for $i,j$ two vertices joined by a single edge $e\colon i \to j$ is a formal consequence of the other relations. In order to unburden the notation, we set
	\begin{align}
		U\coloneqq x_i^-(u)\ , \quad V\coloneqq x_i^-(v)\ , \quad W\coloneqq x_j^-(w)\ .
	\end{align}
	The quadratic relations read
	\begin{align}\label{eq:cubicserre1}
		(u-v+\hbar)UV+(v-u+\hbar)VU&\almostsame 0\ ,\\[4pt]\label{eq:cubicserre2}  
		(u-w-\varepsilon_{e^\ast})UW-(u-w+\varepsilon_{e})WU&\almostsame 0\ ,\\[4pt]\label{eq:cubicserre3} 
		(v-w-\varepsilon_{e^\ast})VW-(v-w+\varepsilon_{e})WV&\almostsame 0\ ,
	\end{align}
	while the cubic relation is
	\begin{align}\label{eq:cubicserre4}
		(w-u+\varepsilon_{e^\ast})UVW-(u-v-\hbar)VWU + (v-w+\varepsilon_e)WUV\almostsame 0\ .
	\end{align}
	Let us denote by $R_{UV}, R_{UW}, R_{VW}$ and $R_{UVW}$ the left-hand-side of \eqref{eq:cubicserre1}, \eqref{eq:cubicserre2}, \eqref{eq:cubicserre3} and \eqref{eq:cubicserre4} respectively, and let us denote by $R_{VUW}$ the cubic relation with the roles of $U$ and $V$ exchanged. Then a direct computation shows that
	\begin{multline}
		R_{UVW}+R_{VUW}+ 2\{R_{UV},W\} + [U,R_{VW}]+[V,R_{UW}]\\
		\almostsame\hbar \big((UV+VU)W -2(UWV+VWU) + W(UV +VU)\big)
	\end{multline} 
	which is the usual Serre relation \eqref{eq:Yangian-Serre} with $a_{i,j}=-1$.
\end{remark}

When $\qv$ is clear from the context, we will sometimes simply denote $\Y_\qv, \eY_\qv, \ldots$ by $\Y, \eY,\ldots$.

\subsection{Gradings and truncations of the Yangians}

In this section, we shall introduce a grading and a truncation of the Yangian, which will play a relevant role later on.

\subsubsection{Horizontal and vertical gradings}\label{sec:yangiangradingsandtruncation}

Let $\Y_\qv^\pm$ and $\zeroeY_\qv$ be the subalgebras of $\eY_\qv$ generated by $\{x_{i,\ell}^\pm\,\vert\, i \in I,\, \ell \in \N\}$ and $\{h_{i,\ell}, \kappa_{i,\ell}\,\vert\, i \in I,\, \ell \in \N\}$. Let $\Y_\qv^{\mathfrak{e},\,\geqslant 0}$ be the subalgebra generated by $\zeroeY_\qv$ and $\Y_\qv^+$ and let $\Y_\qv^{\frake,\,\leqslant 0}$ be the subalgebra generated by $\zeroeY_\qv$ and $\Y_\qv^-$. Similarly, we introduce $\Y_\qv^{\geqslant 0}$, $\Y_\qv^0$,  and $\Y_\qv^{\leqslant 0}$.

The algebras $\eY_\qv$, $\Y_\qv^\pm$, and $\zeroeY_\qv$ are $\Z\times\Z I$-graded with
\begin{align}\label{degree_grading_yangian}
	\deg\big(x^\pm_{i,\ell}\big)\coloneqq (-2\ell,\pm \alpha_i)\ , \ \deg\big(h_{i, \ell}\big) \coloneqq (-2\ell,0)\ ,\ \deg\big(\kappa_{i, \ell}\big) \coloneqq (-2\ell,0)\ ,
\end{align}
and $\deg\big(\epsilon_e\big)\coloneqq (-2,0)$. 

Let $\eY_\qv[\ell,\bfd]$, $\Y_\qv^\pm[\ell,\bfd]$, and $\zeroeY_\qv[\ell,\bfd]$ (resp.\ $\eY_{\bfd}$, $\Y_{\bfd}^\pm$, and $\zeroeY_{\bfd}$) be the graded pieces of degree $(-2\ell, \bfd)$ (resp.\  of degree $\bfd$) of $\eY_\qv$, $\Y_\qv^\pm$, and $\zeroeY_\qv$, respectively. The $\Z$-grading is called the \textit{vertical grading}, the $\Z I$-grading the \textit{horizontal one}.

It is also fruitful to consider the filtration on the same algebras obtained by putting the generators in the same degree as in \eqref{degree_grading_yangian} but the scalars $\epsilon_e$ in degree $0$. We shall call this the \textit{standard} filtration of $\eY_\qv, \Y_\qv, \ldots$. 

Set $x^\pm_i\coloneqq x^\pm_{i,0}$, $h_i\coloneqq h_{i,0}$, and $\kappa_i\coloneqq \kappa_{i,0}$ with $i\in I$. The elements $x^\pm_i, h_i, \kappa_i$ with $i \in I$ generate the subalgebra $\eYzero$ of $\eY_\qv$ consisting of all elements of vertical degree zero. Similarly, we define $\Yzero$. Set $\Yzero^\pm\coloneqq \Y_\qv^\pm\cap \eYzero=\Y_\qv^\pm\cap \Yzero$.

\begin{lemma} 
	The algebra $\Y_0$ is isomorphic to the enveloping algebra $\sfU(\frakg_\qv)$ of the Kac-Moody Lie algebra associated to $\qv$. 
\end{lemma}
\begin{proof}
	Since all the relations defining $\Y_\qv$ are homogeneous, $\Y_0$ is isomorphic to the algebra generated by $\{x_i^{\pm}, h_i\}_{i \in I}$ modulo the degree zero part of the relations. It suffices to observe that these reduce to the usual Serre presentation of $\frakg_\qv$.
\end{proof}

\subsubsection{Augmentation ideal and truncations of the Yangian}\label{subsubsec:truncations-Yangian}

For $i\in I$, let $\Y_i$ be the subalgebra generated by the elements $x_{i,\ell}^\pm$ and $h_{i,\ell}$ with $\ell \in \N$. The subalgebras $\Y_i^\pm$ and $\Y_i^0$, etc., are defined as above. 

Let $\I$ be the \textit{augmentation ideal} in $\Y_\qv$. We abbreviate
\begin{align}
	\I^\pm\coloneqq \I\cap\Y_\qv^\pm\ ,\quad\I^{\geqslant 0}\coloneqq \I\cap\Y_\qv^{\geqslant 0}\ ,\quad\I^{\leqslant 0}\coloneqq \I\cap\Y_\qv^{\leqslant 0}\ ,\quad \I_i^\pm\coloneqq \Y_i\cap\I^\pm\ .
\end{align}

For a future use, we now introduce some truncations of the algebra $\Y_\qv^-$. Fix $i\in I$ and define
\begin{align}\label{eq:iYi}
	\iY_\qv\coloneqq \Y_\qv^-/\,\I_i^-\Y_\qv^-\qquad\text{and}\qquad \Yi_\qv\coloneqq \Y_\qv^-/\,\Y_\qv^-\I_i^-\ . 
\end{align}
There are right and left $\Y_\qv^-$-modules structures on $\iY_\qv$ and $\Yi_\qv$, respectively,  such that the canonical maps
\begin{align}\label{eq:projections}
	\iii{\pi}\colon\Y_\qv^-\longrightarrow \iY_\qv\quad \text{and} \quad \pi^{(i)}\colon \Y_\qv^-\longrightarrow\Yi_\qv
\end{align} 
are $\Y_\qv^-$-module homomorphisms. Let $\iYzero$ and $\Yizero$ be the subspaces of vertical degree zero.

\subsection{Triangular decomposition}

In this section, we prove the triangular decomposition of $\eY_\qv$:
\begin{theorem}\label{thm:triangular-decomposition}
	The multiplication yields the following isomorphisms
	\begin{align}\label{eq:triangular}
		\begin{tikzcd}[ampersand replacement=\&]
			\Y_\qv^+ \otimes \zeroeY_\qv \otimes \Y_\qv^- \ar{r}{\sim}\& \eY_\qv
		\end{tikzcd}\ ,
	\end{align}	
	and
	\begin{align}\label{eq:half-triangular} 
		\begin{tikzcd}[ampersand replacement=\&]
			\Y_\qv^+\otimes\zeroeY_\qv  \ar{r}{\sim}\& \Y_\qv^{\frake,\geqslant 0} \& \text{and} \& \zeroeY_\qv\otimes \Y_\qv^-  \ar{r}{\sim} \& \Y_\qv^{\frake,\leqslant 0}
		\end{tikzcd}\ .
	\end{align}
\end{theorem}

Theorem~\ref{thm:triangular-decomposition} yields the existence of a canonical projection
\begin{align}
	\pr\colon\eY_\qv\longrightarrow \Y_\qv^-=\eY_\qv/\I^{\geqslant 0}\eY_\qv\ .
\end{align}

Furthermore, note that there exists a version of Theorem~\ref{thm:triangular-decomposition} for the Yangian $\Y_\qv$, which yields the canonical projection
\begin{align}\label{eq:projection}
	\pr\colon\Y_\qv\longrightarrow \Y_\qv^-=\Y_\qv/\I^{\geqslant 0}\Y_\qv\ .
\end{align}

Our method of proof
 Theorem~\ref{thm:triangular-decomposition} will follow the standard approach in the case of quantum loop algebras (see e.g. Hernandez \cite[\S3]{Hernandez-quantum}). First, let us recall the following useful lemma.

Let $A$ be an algebra with a triangular decomposition $(A^+, H, A^-)$. Let $\tau^+$ (resp.\ $\tau^-$) be a two-sided ideal of $A^+$ (resp.\ $A^-$). Let $B \coloneqq A/(A\cdot(\tau^+ +\tau^-)\cdot A)$ and denote by $B^\pm$ the image of $A^\pm$ in $B$.

\begin{lemma}[{\cite[Lemma~4]{Hernandez-quantum}}]\label{lem:Humphreys}
	If $\tau^-\cdot A \subseteq A\cdot \tau^-$ and $A\cdot \tau^+ \subseteq \tau^+\cdot A$, then $(B^+, H, B^-)$ is a triangular decomposition of $B$ and the algebra $B^\pm$ is isomorphic to $A^\pm /\tau^\pm$.
\end{lemma}

\begin{definition}
	Let $\widetilde{\widetilde{\eY_\qv}}$ be the unital associative $\Hbullet_{\Tmax}$-algebra generated by $x_{i, \ell}^\pm, h_{i, \ell}, \kappa_{i, \ell}$, with $i \in I$ and $\ell \in \N$, subject to the relations \eqref{eq:Yangian-Lie-algebra-1}, \eqref{eq:Yangian-Lie-algebra-2}, \eqref{eq:Yangian-Lie-algebra-3}, \eqref{eq:Yangian-4}, \eqref{eq:Yangian-6}, and \eqref{eq:Yangian-7}.
	
	Let $\widetilde{\eY_\qv}$ be the quotient of $\widetilde{\widetilde{\eY_\qv}}$ by the relations \eqref{eq:Yangian-5}, \eqref{eq:Yangian-8} and \eqref{eq:Yangian-9}.
\end{definition}

It is easy to see that the triple $\Big(\Big(\widetilde{\widetilde{\eY_\qv}}\Big)^+, \zeroeY_\qv, \Big(\widetilde{\widetilde{\eY_\qv}}\Big)^-\Big)$ gives rise to a triangular decomposition of $\widetilde{\widetilde{\eY_\qv}}$, where $\Big(\widetilde{\widetilde{\eY_\qv}}\Big)^\pm$ is generated by the $x_{i,\ell}^\pm$ without relations (see the end of Remark~\ref{rem:specialized-relations}).

Let $\widetilde{\tau}^\pm$ be the two-sided ideal of $\Big(\widetilde{\widetilde{\eY_\qv}}\Big)^\pm$ generated by
\begin{align}
	\Big[(u-w-\hbar)x_i^+(u)x_i^+(w)- (u-w+\hbar)x_i^+(w)x_i^+(u)\Big]_{<0}& \quad \text{for } i\in I\ ,\\[4pt]
	\Big[\zeta_{i,j}(u-w) x_i^+(u)x_j^+(w)- \zeta_{i, j}(u-w-\hbar)x_j^+(w)x_i^+(u)\Big]_{<0}& \quad \text{for } i, j\in I\ , \ i\neq j\ ,
\end{align}
and
\begin{align}
	\Big[(u-w+\hbar)x_i^-(u)x_i^-(w) -(u-w-\hbar)x_i^-(w)x_i^-(u)\Big]_{<0}& \quad \text{for } i\in I\ ,\\[4pt]
	\Big[\zeta_{i,j}(u-w-\hbar) x_i^-(u)x_j^-(w)- \zeta_{i,j}(u-w)x_j^-(w)x_i^-(u)\Big]_{<0}&  \quad \text{for } i, j\in I\ , \ i\neq j \ ,
\end{align}
respectively, where $[-]_{< 0}$ denotes truncation in negative degree of the power series.
\begin{lemma}\label{lem:Hernandez}
	The following inclusions hold:
	\begin{align}
		\widetilde{\tau}^-\cdot \Big(\widetilde{\widetilde{\eY_\qv}}\Big)^-\subseteq \Big(\widetilde{\widetilde{\eY_\qv}}\Big)^- \cdot \widetilde{\tau}^-\quad\text{and}\quad \Big(\widetilde{\widetilde{\eY_\qv}}\Big)^+\cdot \widetilde{\tau}^+\subseteq \widetilde{\tau}^+\cdot \Big(\widetilde{\widetilde{\eY_\qv}}\Big)^+\ .
	\end{align}
	In particular, the triple
	\begin{align}
		\Big( \Big(\widetilde{\eY_\qv}\Big)^+, \zeroeY_\qv,  \Big(\widetilde{\eY_\qv}\Big)^-\Big)
	\end{align}
	is a triangular decomposition of $\widetilde{\eY_\qv}$, where $\Big(\widetilde{\eY_\qv}\Big)^\pm$ is the subalgebra generated by the $x_{i, \ell}^\pm$.
\end{lemma}

\begin{proof}
	First, a direct computation shows that
	\begin{align}
		\Big[(u-w+\hbar)x_i^+(u)x_i^+(w)- (u-w-\hbar)x_i^+(w)x_i^+(u), x_k^-(z)\Big]\almostsame 0 \ ,
	\end{align}
	by using the relations \eqref{eq:Yangian-4-alternative-formulation} and \eqref{eq:Yangian-5-alternative-formulation}. In addition, by using the relation \eqref{eq:Yangian-Lie-algebra-3-alternative-formulation}, one can show that 
	\begin{align}
		\Big[\zeta_{i,j}(u-w-\hbar) x_i^+(u)x_j^+(w)- \zeta_{i, j}(u-w)x_j^+(w)x_i^+(u), x_k^-(z)\Big]\almostsame 0\ , 
	\end{align}
	where the last equality follows from the relation \eqref{eq:Yangian-6} and the ``functional equation'' \eqref{eq:functional-equation}. A similar but tedious computation shows that
	\begin{align}
		\widetilde{\tau}^-\cdot \zeroeY_\qv\subseteq \zeroeY_\qv\cdot \widetilde{\tau}^-\ .
	\end{align}
	Using the triangular decomposition $\widetilde{\widetilde{\eY_\qv}}=\Big(\widetilde{\widetilde{\eY_\qv}}\Big)^+\cdot \zeroeY_\qv\cdot \Big(\widetilde{\widetilde{\eY_\qv}}\Big)^-$ we easily deduce that 
	\begin{align}
		\widetilde{\tau}^-\cdot \widetilde{\widetilde{\eY_\qv}}\subseteq \widetilde{\widetilde{\eY_\qv}} \cdot \widetilde{\tau}^-\ .
	\end{align}
	The case of $\widetilde{\tau}^-$ may be shown in the same way. We conclude using Lemma~\ref{lem:Humphreys}. 
\end{proof}

Now, set
\begin{align}
	\tau^\pm\coloneqq \ker\Big( \Big(\widetilde{\eY_\qv}\Big)^\pm \longrightarrow \Y^\pm_\qv \Big)\ .
\end{align}
\begin{proof}[Proof of Theorem~\ref{thm:triangular-decomposition}] The proof boils down to some direct but lengthy computations, which we sketch below.
	Using Lemma~\ref{lem:Humphreys} again, it is enough to show that
	\begin{align}
		\tau^-\cdot \widetilde{\eY_\qv}\subseteq \widetilde{\eY_\qv} \cdot \tau^-\quad\text{and}\quad \widetilde{\eY_\qv}\cdot \tau^+\subseteq \tau^+\cdot \widetilde{\eY_\qv}\ .
	\end{align}
	
	First, using relations \eqref{eq:Yangian-4}, \eqref{eq:Yangian-6}, and \eqref{eq:Yangian-7}, one shows that
	\begin{align}
		\tau^-\cdot \zeroeY_\qv \subset \zeroeY_\qv\cdot \tau^- \quad \text{and} \quad \zeroeY_\qv\cdot \tau^+ \subset \tau^+\cdot\zeroeY_\qv \ .
	\end{align}
	It follows that 
	\begin{align}
		\tau^-\cdot \Big(\eY_\qv\Big)^{\leqslant 0} \subset \Big(\eY_\qv\Big)^{\leqslant 0}\cdot \tau^- \quad \text{and} \quad \Big(\eY_\qv\Big)^{\leqslant 0} \cdot\tau^+ \subset \tau^-\cdot\Big(\eY_\qv\Big)^{\geqslant 0} \ ,
	\end{align}
	where $\Big(\widetilde{\eY_\qv}\Big)^{\geqslant 0}$ and $\Big(\widetilde{\eY_\qv}\Big)^{\leqslant 0}$ are the subalgebras of $\widetilde{\eY_\qv}$ generated by $\Big(\widetilde{\eY_\qv}\Big)^+$ and $\zeroeY_\qv$ and by $\Big(\widetilde{\eY_\qv}\Big)^-$ and $\zeroeY_\qv$, respectively. The isomorphisms in \eqref{eq:half-triangular} for $\Y_\qv^{\frake,\, \leqslant 0}$ and $\Y_\qv^{\frake,\, \geqslant 0}$ thus hold by Lemma~\ref{lem:Humphreys}. 
	
	Using the relations \eqref{eq:Yangian-Lie-algebra-3}, one see that
	\begin{align}
		[\tau^+, x_{i, \ell}^-] \in\Big(\widetilde{\eY_\qv}\Big)^{\geqslant 0}\quad \text{and}\quad
		[\tau^-, x_{i, \ell}^+] \in \Big(\widetilde{\eY_\qv}\Big)^{\leqslant 0} \ .
	\end{align}	
	By projecting to $\Y_\qv^{\geqslant 0}$ and $(\Y_\qv^{\leqslant 0}$, respectively, we get that
	\begin{align}
		[\tau^+, x_{i, \ell}^-] &\in \ker\Big( \Big(\widetilde{\eY_\qv}\Big)^{\geqslant 0} \longrightarrow \Y_\qv^{\frake,\,\geqslant 0} \Big) = \tau^+\cdot\zeroeY_\qv\cdot\ , \\[4pt]
		[\tau^-, x_{i, \ell}^+] &\in \ker\Big( \Big(\widetilde{\eY_\qv}\Big)^{\leqslant 0} \longrightarrow \Y_\qv^{\frake,\,\leqslant 0} \Big) = \zeroeY_\qv \cdot \tau^-\ .
	\end{align}
	Here the equalities on the right-hand-side follow from the isomorphisms in \eqref{eq:half-triangular}. Therefore, the triangular decomposition \eqref{eq:triangular} of $\eY_\qv$ holds thanks again to Lemma~\ref{lem:Humphreys}.
\end{proof}

\begin{remark}
	Note that a priori we could have defined 
	\begin{align}
		\widetilde{\widetilde{\tau}}^\pm\coloneqq \ker\Big( \Big(\widetilde{\widetilde{\eY_\qv}}\Big)^\pm \longrightarrow \Y^+_\qv \Big)\ ,
	\end{align}
	and proved that
	\begin{align}
		\Y^\pm_\qv\simeq \Big(\widetilde{\widetilde{\eY_\qv}}\Big)^\pm / \widetilde{\widetilde{\tau}}^\pm 
	\end{align}
	gives rise to a triangular decomposition of $\eY_\qv$ by using the same arguments as in the proof of Theorem~\ref{thm:triangular-decomposition}. By introducing the algebra $\widetilde{\eY_\qv}$ and proving the triangular decomposition for it, we deduce that $\Y^\pm_\qv$ can be seen as the algebra generated by the $x_{i,\ell}^\pm$, with $i\in I$ and $\ell\in \N$, subject to the relations \eqref{eq:Yangian-5}, \eqref{eq:Yangian-8}, \eqref{eq:Yangian-9}, respectively, the cubic and the Serre relations and `only' the extra relations obtained by applying $\ad(x_{i,\ell}^\pm)$ to the cubic and the Serre relations.
\end{remark}

\section{Braid group action on the Yangian}\label{sec:braid-group-action-Yangian}

In this section, we shall introduce the braid group action on the Yangian, its negative nilpotent half, and certain quotients of the latter.

\subsection{The braid group actions on $\Y_\qv$ and on $\Y_\qv^-$}\label{sec:braid-on-Yangian}

There is a $B_\qv$-action on $\Y_\qv$ by algebra automorphisms defined as
\begin{align}\label{eq:Ti}
	B_\qv&\longrightarrow \Aut(\Y_\qv) \ , \\
	T_i&\longmapsto \exp\big(\ad\big(x_i^+\big)\big)\circ\exp\big(-\ad\big(x_i^-\big)\big)\circ\exp\big(\ad\big(x_i^+\big)\big)
\end{align}
for $i\in I$. These operators are well-known in the context of Kac-Moody algebras. In the context of Yangians they were first considered in \cite{GNW18}, see Formula~(3.15) in \textit{loc.cit.}, for Yangians of affine ADE quivers and further studied in \cite{Kodera_Braid19}. This gives well-defined operators since the adjoint representation of the Yangian, restricted to each root $\fraksl_2$, is integrable.
\begin{lemma}\label{lem:braid-group-action-weight}
	Let $\bfd\in\Z I$.
	\begin{enumerate}\itemsep0.2cm
		\item \label{item:braid-group-action-weight-2} We have $T_i(\Y_\bfd) =\Y_{s_i(\bfd)}$ and $T_i(\Y^+_\bfd)\subset\sum_\bfe \Y^{\geqslant 0}_\bfe\Y_i^-$, where $\bfe$ runs over the set $s_i(\bfd)+\N\alpha_i$.
		
		\item \label{item:braid-group-action-weight-1} $T_i(\I_i^-) \subset \I_i^+\Y_i$, and $T_i(\Y_\qv^-)\subset\Y_i^{\geqslant 0}\Y_\qv^-$.
	\end{enumerate}
	
	Moreover, let $w=s_{i_\ell}\cdots s_{i_2}s_{i_1}$ be a reduced expression of $w\in W_\qv$ and $\beta_r\coloneqq s_{i_\ell}\cdots s_{i_{r+1}}(\alpha_{i_r})$ for $r=1,\ldots,\ell$. We have:
	\begin{enumerate}\setcounter{enumi}{2}
		\item \label{item:braid-group-action-weight-3} $T_w(\Y_\qv^-)\subset\sum_\bfe\Y^{\geqslant 0}_\bfe\Y_\qv^-$, where $\bfe$ runs over the set $\sum_{r=1}^{\ell}\N \beta_r$.
	\end{enumerate}
\end{lemma}

\begin{proof}
	The first claim is straightforward and well-known (see e.g. to \cite[Proposition~A.2]{Kodera_Braid19}). Thanks to Theorem~\ref{thm:triangular-decomposition} (the triangular decomposition), the second claim in \eqref{item:braid-group-action-weight-2} follows from the inclusion
	\begin{align}
		T_i\big(\Y^+_\bfd\big)\subset \big(\Y_\qv^{\geqslant 0}\Y_i^-\big)_{s_i(\bfd)}\ ,
	\end{align}
	where $(-)_{s_i(\bfd)}$ is the degree $s_i(\bfd)$ part of $(-)$. We next turn to \eqref{item:braid-group-action-weight-1}. Clearly, $T_i(\Y_i)=\Y_i$. Since $\I_i^-\subset\sum_{\ell\in \N}x_{i,\ell}^-\Y_i$, it is enough to check that $T_i(x^-_{i,\ell})\in\I_i^+$ for each integer $\ell$. Now, the relations \eqref{eq:Yangian-6} (see also \eqref{eq:Yangian-6-alternative-formulation}) and \eqref{eq:Yangian-7} (see also \eqref{eq:Yangian-4-alternative-formulation}) yield 
	\begin{align}
		x^-_{i, \ell}= -\frac{1}{2} \Big[h_{i, 0}, x^-_{i, \ell}\Big]= \frac{1}{2}\left(\Big[h_{i, 1}, x^-_{i, \ell-1}\Big]+\hbar \Big\{h_{i, 0}, x^-_{i, \ell-1}\Big\} \right)\ .
	\end{align}
	Thus, the claim follows by iteratively applying \cite[Proposition~A.2]{Kodera_Braid19}. The second statement in \eqref{item:braid-group-action-weight-1} follows from the above argument together with the equalities
	\begin{align}
		\Y_\qv^+\Y_\qv^0=\Y_\qv^0\Y_\qv^+\quad\text{and}\quad \Y_i^{\leqslant 0}\Y_\qv^+=\Y_\qv^+\Y_i^{\leqslant 0}\ .
	\end{align}
	Finally, \eqref{item:braid-group-action-weight-3} follows from \eqref{item:braid-group-action-weight-1} and \eqref{item:braid-group-action-weight-2} by induction.
\end{proof}

For any $w\in W_\qv$, we define the linear operator $\overline T_w$ given by the composition
\begin{align}\label{eq:overlineT}
	\overline T_w\colon
	\begin{tikzcd}[ampersand replacement=\&]
		\Y_\qv^-\arrow{r} \& \Y_\qv\arrow{r}{T_w} \&\Y_\qv\arrow{r}{\pr} \&\Y_\qv^-
	\end{tikzcd}\ ,
\end{align} 
where the last map is the projection \eqref{eq:projection}. We call $\overline T_w$ the \textit{truncated} braid group operator associated to $w$.

Let $B^+_\qv\subset B_\qv$ be the submonoid generated by the elements $T_w$ with $w\in W_\qv$. We call $B^+_\qv$ the \textit{positive braid monoid} of $\qv$. It is known that $B^+_\qv$ is isomorphic to the monoid generated by the $T_i$'s modulo the braid relations \eqref{eq:braid-relations}.
\begin{proposition}\label{prop:truncated-action}
	The assignment $T_w \mapsto \overline T_w$ for $w \in W_\qv$ gives rise to a representation of $B^+_\qv$ on $\Y_\qv^-$, i.e., to a morphism of groupoids $B^+_\qv \to \End(\Y_\qv^-)$.
\end{proposition}
\begin{proof}
	Let $w=s_{i_\ell}\cdots s_{i_2}s_{i_1}$ be a reduced expression and $\beta_r=s_{i_\ell}\cdots s_{i_{r+1}}(\alpha_{i_r})$ for $r=1,\ldots,\ell$. We shall argue by induction on $r$.
	
	Let $v\coloneqq s_{i_{\ell-1}}\cdots s_{i_2}s_{i_1}$ and $i=i_\ell$. We will prove that $\overline T_w=\overline T_i\overline T_v$. Fix an element $x\in\Y_\qv^-$.	By Lemma~\ref{lem:braid-group-action-weight}, we have $T_v(x)=y+z$, with $y=\overline T_v(x) \in\Y_\qv^-$ and
	\begin{align}
		z\in \sum_{\bfe}\I^{\geqslant 0}_\bfe\Y_\qv^-\quad\text{with}\quad 0\neq \bfe\in\sum_{r=1}^{\ell-1}\N s_{i_\ell}\cdots s_{i_{r+1}}(\alpha_{i_r})\ .
	\end{align}
		We have $T_w(x)=T_i(y)+T_i(z)$ and
	\begin{align}
		T_i(z)\in  \sum_{\bfe} T_i(\I^{\geqslant 0}_\bfe)T_i(\Y_\qv^-)\quad\text{with}\quad 0\neq \bfe\in\sum_{r=1}^{\ell-1}\N s_{i_\ell}\cdots s_{i_{r+1}}(\alpha_{i_r})\ .
	\end{align}
	Note that all weights $s_{i_\ell}\cdots s_{i_{r+1}}(\alpha_{i_r})$ in the formula above are positive roots, which are different from $\alpha_i$, hence $T_i(z)\in \I^{\geqslant 0}T_i(\Y_\qv^-)\subset \I^{\geqslant 0}\Y_\qv$. Thus $T_i(z)$ belongs to the kernel of the projection \eqref{eq:projection}. Therefore $\overline T_w(x)$ is the image of $T_i\overline T_v(x)$ by the projection \eqref{eq:projection}, proving the proposition.
\end{proof}

\subsection{Truncated braid group action}

\begin{proposition}\label{prop:overlineTi}
	There exists a unique map $\overline T_i\colon \iY_\qv \to \Yi_\qv$ such that the diagram
	\begin{align}\label{eq:commutative-diagram-A}
		\begin{tikzcd}[ampersand replacement=\&]
			\Y_\qv^-  \ar{r}{\overline T_i} \ar[swap]{d}{\iii{\pi}} \&\Y_\qv^- \ar{d}{\pi^{(i)}} \\
			\iY_\qv \ar{r}{\overline T_i} \& \Yi_\qv
		\end{tikzcd}
	\end{align}
	is commutative. Here, the vertical arrows are the canonical projections \eqref{eq:projections}.
\end{proposition}

\begin{proof}
	We have an inclusion $T_i(\I_i^-\Y_\qv^-)\subset\I_i^+\Y_\qv \subset \I^{\geqslant 0}\Y_\qv$ hence there is a factorization
	\begin{align}\label{eq:commutative-diagram-truncated-braid}
		\begin{tikzcd}[ampersand replacement=\&]
			\Y_\qv^-  \ar{r}{T_i} \ar[swap]{d}{\iii{\pi}} \ar{dr}{\overline T_i} \&\Y_\qv \ar{d}{\pr} \\
			\iY_\qv \ar{r}{} \& \Y^-_\qv
		\end{tikzcd}
	\end{align}
	We now compose with the projection $\pi^{(i)} \colon \Y_\qv^-\to \Yi_\qv$ to obtain the desired commutative diagram.
\end{proof}

\section{Relation between the Yangian and the nilpotent COHA of the quiver}\label{sec:relation-Yangian-COHA}

In this section, we construct a surjective algebra homomorphism from the negative half $\Y_\qv^-$ of the Yangian to the cohomological Hall algebra $\cohaqv^{\Ttilde}$ for $\Ttilde\coloneqq\Tmax$. We shall use the notation $\Ttilde$ and $\Tmax$ interchangeably.

\subsection{The algebra $\bS$}\label{sec:the-algebra-S}

Define the $\Z\times\N I$-graded commutative algebra
\begin{align}
	\bS \coloneqq \Hbullet_{\Ttilde}[ p_\ell(z_i)\,\vert\, i \in I, \ \ell \geq 0]\ ,
\end{align}
where 
\begin{align}
	\deg(\varepsilon_e)=(2,0)\ ,\quad\deg(p_\ell(z_i))=(2\ell,\alpha_i)\ .
\end{align}
For any $i\in I$, we set
\begin{align}\label{eq:defchzi}
	\ch(z_i) \coloneqq p_0(z_i)+\sum_{\ell\geq 1} \frac{p_\ell(z_i)}{\ell !}\ .
\end{align}
We endow the algebra $\bS$ with a cocommutative bialgebra structure by setting
\begin{align}
	\Delta(\ch(z_i))\coloneqq\ch(z_i) \otimes 1 + 1 \otimes \ch(z_i)
\end{align}
for each $i \in I$. 
\begin{remark}
	It can be helpful to think of $p_\ell(z_i)$ as the $\ell$-th power sum in an infinite set of variables $z_{i,1}, z_{i,2}, \ldots$ and of $p_0(z_i)$ as an extra formal variable. 
\end{remark}

For each dimension vector $\bfd\in\N I$, there is a canonical algebra homomorphism 
\begin{align}
	\iota \colon \bS \longrightarrow \Hbullet_{\Ttilde\times \G_\bfd}(\mathsf{pt})\ , \quad
	\ch(z_i) \longmapsto \ch(\calV_i)\ ,
\end{align}
where $\calV_i$ is the tautological vector bundle on $\sfB\G_\bfd$ pulled back from $\sfB\GL(d_i)$ introduced in \S\ref{subsubsec:moduli-stacks}. Observe that $\iota(p_0(z_i))=d_i$ for any $i\in I$. Composing the natural $\Hbullet_{\Ttilde\times \G_\bfd}$-action on $\cohaqvd^{\Ttilde}$ with the map $\iota$ and taking the sum over all dimensions $\bfd$, we get an action of $\bS$ on the $\Hbullet_{\Ttilde}$-algebra $\cohaqv^{\Ttilde}$, which we denote as
\begin{align}
	\cap\colon \bS\otimes\cohaqv^{\Ttilde}\longrightarrow\cohaqv^{\Ttilde} \ .
\end{align}

\begin{lemma}
	The action $\cap$ is Hopf, i.e., for any $c \in \bS$ and $a_1,a_2 \in \cohaqv^{\Ttilde}$ we have
	\begin{align}
		c \cap (a_1 \cdot a_2) = \sum_j (c'_j \cap a_1) \cdot (c''_j \cap a_2) \ ,
	\end{align}
	where $\Delta(c)=\sum_j c'_j \otimes c''_j$.
\end{lemma}
\begin{proof} 
	Fix $\bfd_1, \bfd_2 \in \N I$ and consider the convolution diagram \eqref{eq:definition-COHA-product-nil}
	\begin{align}
		\begin{tikzcd}[ampersand replacement=\&]
			\dLambda_{\bfd_1} \times \dLambda_{\bfd_2} \& \dLambda_{\bfd_1,\bfd_2}^{\mathsf{ext}} \ar{r}{p} \ar[swap]{l}{q} \& \dLambda_{\bfd_1+\bfd_2}\,.
		\end{tikzcd}
	\end{align}
	There is a short exact sequence of tautological bundles
	\begin{align}
		0 \longrightarrow \calV_{\bfd_2} \longrightarrow \calV_{\bfd_1 + \bfd_2} \longrightarrow \calV_{\bfd_1} \longrightarrow 0\ .
	\end{align}
	The lemma follows from the additivity of the Chern character and the projection formula, as in \cite[Proposition~4.2]{Sala_Schiffmann}.
\end{proof}

On the other hand, there is an action of $\bS$ on $\eY_\qv$
\begin{align}
	\cap\colon\bS\otimes\eY_\qv\longrightarrow\eY_\qv
\end{align}
such that for each vertex $i\in I$ and each positive integer $\ell$, the polynomial $p_\ell(z_i)$ acts as the derivation satisfying
\begin{align}\label{eq:S-action-Yangian}
	p_\ell(z_i) \cap x^\pm_{j,k} =\pm \delta_{i,j}\,x^\pm_{j,k+\ell}\ , \quad p_\ell(z_i)\cap h_{j,k} =0\ ,\quad  p_\ell(z_i)\cap \kappa_{j,k} =0\ ,
\end{align}
and $p_0(z_i)$ acts as $d_i\cdot \id$ on the $\bfd$-graded piece $\eY_\bfd$. This action preserves the subalgebras $\Y_\qv^\pm$ and $\zeroeY_\qv$. It is of degree zero with respect to the $\Z I$-grading.

\subsection{The homomorphism $\Phi$}\label{sec:Phi}

The goal of this section is to prove the following theorem.
\begin{theorem}\label{thm:Phi}
	There is a surjective $\Z\times\N I$-graded $\Hbullet_{\Ttilde}$-algebra homomorphism 
	\begin{align}\label{eq:Phi}
		\begin{split}
			\Phi\colon\Y_\qv^-&\longrightarrow \cohaqv^{\Ttilde}\ ,\\
			x^-_{i,\ell} &\longmapsto (z_{i,1})^\ell\cap[\Lambda_{\alpha_i}]
		\end{split}
	\end{align}
	for $i\in I$ and $\ell\in \N$. Moreover, we have:
	\begin{enumerate}\itemsep0.2cm
		\item \label{thm:Phi-1}
		$\Phi$ intertwines the vertical gradings on $\Y_\qv^-$ and  $\cohaqv^{\Ttilde}$.
		\item  \label{thm:Phi-2}
		$\Phi$ intertwines the $\bS$-actions on $\Y_\qv^-$ and $\cohaqv^{\Ttilde}$.
		\item  \label{thm:Phi-3}
		$\Phi$ restricts to an algebra isomorphism $\Y^-_0 \to\cohazero^{\Ttilde}$. 
	\end{enumerate}
\end{theorem}

Proposition~\ref{prop:sign-twists-I} yields the following.
\begin{corollary}\label{cor:Phi-sign-twists} 
	For any choice of twist $\Theta\colon \Z I \times \Z I \to \Z/2\Z$, there exists a surjective algebra morphism 
	\begin{align}
		\Phi_\Theta\colon\Y^-_\qv \longrightarrow \cohaqv^{\Ttilde,\Theta}
	\end{align}
	satisfying the same conditions as in Theorem~\ref{thm:Phi}. The isomorphisms $\gamma_{\Theta,\Theta'}$ are compatible with the morphisms $\Phi$, i.e., we have $\gamma_{\Theta,\Theta'}\Phi_{\Theta}=\Phi_{\Theta'}$.
\end{corollary}
When there is no risk of confusion, we will simply write $\Phi$ for $\Phi_\Theta$. Motivated by the main result of \cite{NSS_KHA}, we make the following conjecture:
\begin{conjecture}
	The map $\Phi$ induces an isomorphism $\Y_\qv^-\otimes \mathsf{Frac}\big(\Hbullet_{\Ttilde}\big) \stackrel{\sim}{\longrightarrow} \cohaqv^{\Ttilde} \otimes \mathsf{Frac}\big(\Hbullet_{\Ttilde}\big)$.
\end{conjecture}

\begin{warning}
	Note that with our conventions, the map $\Phi$ does not preserve weights, but rather sends the graded piece of $\Y^-_\qv$ of weight $-\alpha \in \rootlattice_\qv$ to the graded piece of $\cohaqv^{\Ttilde}$ of weight $\alpha$.
\end{warning}
To prove Theorem~\ref{thm:Phi}, we first recall the shuffle algebra presentation of the cohomological Hall algebra of a quiver.

\subsubsection{Shuffle algebra and cohomological Hall algebra of quiver}\label{subsubsec:shuffle-algebra}

For $\bfd\in \N I$, define
\begin{align}
	\shuffle_\bfd\coloneqq \Hbullet_{\Ttilde}\Big[z_{i,\ell}\,\Big\vert\, i \in I\ , \ 1 \leq \ell \leq d_i\Big]^{\frakS_\bfd}\ ,
\end{align}
where $\frakS_d$ is the permutation group of $d$ letters and $\frakS_\bfd\coloneqq \prod_i \frakS_{d_i}$  is permuting the set of variables $\{z_{i,1},\cdots, z_{i,d_i}\}$ for each vertex $i$. 

When considering elements of $\shuffle_\bfd$, we will indicate the range of variables in bracket: for instance, $P(z_{[1,\bfd]})$ is a function in the variables $(z_{i,\ell})$ for $i \in I$ and $1 \leq \ell \leq d_i$, while $P(z_{[\bfd+1,\bfd+\bfe]})$ is a function in variables $(z_{i,\ell})$ for $i \in I$ and $d_i+1 \leq \ell \leq d_i+e_i$.

Note that there is a canonical identification $\shuffle_\bfd=H^{\Ttilde}_\ast(\stackRep_\bfd(\doubleqv))$. We define a graded associative algebra structure on $\shuffle_\qv\coloneqq H_\ast^{\Ttilde}(\stackRep(\doubleqv))$ such that, for each $P\in \shuffle_\bfd$ and $Q \in \shuffle_\bfe$, we have 
\begin{multline}\label{eq:shuffle-definition}
	\big(P \star Q\big)\big(z_{[1,\bfd+\bfe]}\big)\coloneqq
	(-1)^{\langle \bfd,\bfe\rangle}\sum_\sigma \sigma \Big( \prod_{i \in I} \zeta_i\big(z_{[1,\bfd]},z_{[\bfd+1,\bfd+\bfe]}\big)
	\cdot \prod_{i,j \in I}\zeta_{i,j}\big(z_{[1,\bfd]},z_{[\bfd+1,\bfd+\bfe]}\big) \cdot\\
	\cdot P\big(z_{[1,\bfd]}\big)\cdot Q\big(z_{[\bfd+1,\bfd+\bfe]}\big)\Big)\ ,
\end{multline}
where $\sigma$ runs among all $I$-tuples of $(d_i,e_i)$-shuffles, and
\begin{align}
	\zeta_i\big({z}_{[1,\bfd]},{w}_{[1,\bfe]}\big)\coloneqq\prod_{\genfrac{}{}{0pt}{}{1 \leq k \leq d_i}{1 \leq \ell \leq e_i}}\frac{z_{i,k}-w_{i,\ell}- \hbar}{z_{i,k}-w_{i,\ell}} \quad\text{and}\quad \zeta_{i,j}\big({z}_{[1,\bfd]},{w}_{[1,\bfe]}\big) \coloneqq \prod_{\genfrac{}{}{0pt}{}{1 \leq k \leq d_i}{1 \leq \ell \leq e_j}}\zeta_{i,j}(z_{i,k}-w_{j,\ell})\ .
\end{align}
Here, $\zeta_{i,j}(t)$ has been introduced in Formula~\eqref{eq:zeta}, while $\langle -,-\rangle$ stands for the \textit{non-symmetrized} Euler pairing attached to $\qv$, introduced in Formula~\eqref{eq:euler-form-quiver}.

\begin{theorem}[{\cite[Theorem~A]{SV_generators}}]\label{thm:iota}
	The pushforward with respect to the closed embedding of classical geometric stacks $\Lambda_\qv \to \stackRep(\doubleqv)$ yields an injective $\Hbullet_{\Ttilde}$-algebra homomorphism $\iota\colon \cohaqv^{\Ttilde}\to \shuffle_\qv$ which is a $\bS$-module homomorphism\footnote{Note that the torus $\Ttilde$ fulfills the genericity condition which \cite[Theorem~A]{SV_generators} requires.}.  
\end{theorem}
Concretely, we have $\iota((z_{i,1}^\ell) \cap [\Lambda_{\alpha_i}])=z_{i,1}^\ell \in\shuffle_{\alpha_i}$ for any $i \in I$ and $\ell \geq 0$.

\subsubsection{Proof of Theorem~\ref{thm:Phi}}

We are ready to prove the main result of this section.
\begin{proof}
	The proof of relations \eqref{eq:Yangian-5}, \eqref{eq:Yangian-9}, and \eqref{eq:Yangian-cubic} are done by a direct computation using the shuffle formula \eqref{eq:shuffle-definition}. Since the first two (quadratic relations) are standard and follow directly from the presentation, we focus on the cubic relation. The images in the shuffle algebra of all the terms involved consist of series of the form 
	\begin{align}
		E(z_1,u_1) \cdots E(z_s,u_s)P(z_1, \ldots, z_s)=\big[E(z_1,u_1) \cdots E(z_s,u_s)P(u_1,\ldots, u_s)\big]_{<0}\ ,
	\end{align} 
	where $[-]_{<0}$ stands for the truncation to negative powers of all $u$ variables. Here $P (z_1, \ldots, z_s) \in \Q[z^{\pm 1}_1, \ldots, z^{\pm 1}_s]$ is any Laurent polynomial. Let us now fix $i \neq j \in I$. Using the above observation and the functional equation \eqref{eq:functional-equation} for the functions $\zeta_{i,j}(t)$, one sees that in the shuffle algebra,
	\begin{align}
		x_i^-(u)x_i^-(v)x_j^-(w)&=-\big[X(u,v,w) \frac{u-v-\hbar}{u-v}\zeta_{i,j}(u-w)\zeta_{i,j}(v-w)\big]_{<0}\ ,\\[4pt]
		x_i^-(v)x_j^-(w)x_i^-(u)&=(-1)^{a_{i,j}+1}\big[X(u,v,w) \frac{u-v+\hbar}{u-v}\zeta_{i,j}(v-w)\zeta_{ji}(w-u)\big]_{<0}\\
		&=-\big[X(u,v,w) \frac{u-v+\hbar}{u-v}\zeta_{i,j}(v-w)\zeta_{i,j}(u-w-\hbar)\big]_{<0}\ ,\\[4pt]
		x_j^-(w)x_i^-(u)x_i^-(v)&=-\big[X(u,v,w) \frac{u-v-\hbar}{u-v}\zeta_{ji}(w-u)\zeta_{ji}(w-v)\big]_{<0}\\
		&=-\big[X(u,v,w) \frac{u-v-\hbar}{u-v}\zeta_{i,j}(u-w-\hbar)\zeta_{i,j}(v-w-\hbar)\big]_{<0}\ ,
	\end{align}
	where
	\begin{align}
		X(u,v,w)=\big(E(z_{i,1},u)E(z_{i,2},v)+E(z_{i,2},u)E(z_{i,1},v)\big)E(z_j,w)\ .
	\end{align}
	The cubic relations follow.
	
	Let us now turn to the Serre relation. By a classical trick of Grojnowski and Nakajima (see \cite[\S10.4]{Quiver_var_finite_dim_quantum_affine}), it is enough to prove it for generators $x^-_{i,\bullet},x^-_{j,\bullet}$ all in vertical degree zero (indeed, by successively applying the derivations $p_\ell(z_i)$ to the degree zero Serre relations, one can generate the Serre relations in any multidegree, in the same way as symmetric polynomials are generated by power sum functions).
	
	In terms of the COHA, this amounts to checking a relation in the subalgebra $\cohazero^{\Ttilde}$. From the equivariant formality of $\cohaqvd^{\Ttilde}$ for any $\bfd \in \N I$ (cf.\ Theorem~\ref{thm:purityandKac}), it follows that the structure of this algebra is independent of the choice of the torus; we may therefore replace $\Ttilde$ by any smaller torus, in particular by the (standard) torus used in e.g. \cite{YZ_preproj_COHA}, where the Serre relation is checked\footnote{Note that the authors of \cite{YZ_preproj_COHA} use the opposite multiplication for the COHA.}. 
	Statements \eqref{thm:Phi-1}, \eqref{thm:Phi-2}, and \eqref{thm:Phi-3} are immediate from the definition of $\Phi$.
	
	The surjectivity of the map $\Phi$ follows from the generation theorem \cite[Theorem~B]{SV_generators} together with the obvious fact that the image of $\Phi$ contains $\coha_{\alpha_i}^{\Ttilde}$ for all $i \in I$. Finally, the restriction $\Phi_0$ of $\Phi$ to $\Y^-_0$ is an isomorphism since $\Y^-_0$ and $\cohazero^{\Ttilde}$ have the same graded dimension as $U(\frakn_\qv)$, where $\frakn_\qv$ is the negative nilpotent subalgebra of $\frakg_\qv$ (see e.g. \cite[Theorem~4.23]{Schiffmann-lectures-canonical-bases}).
\end{proof}

\begin{remark}
	To prove the degree zero Serre relations, one may alternatively use the characteristic cycle map and the compatibility with Lusztig's geometric Hall algebra (see \S\ref{sec:CCmap}).
\end{remark}

\subsection{Extension by the loop Cartan}\label{sec:extensionofPhi}

Motivated by the analogy with Maulik-Okounkov Yangians and the theory of Nakajima quiver varieties, we introduce an extension $\cohaqv^{\Ttilde, \leqslant 0}$ of $\cohaqv^{\Ttilde}$ by adding a 'Cartan' part and extend the morphism $\Phi$ to that setting. For this, let us consider the commutative polynomial algebra
\begin{align}
	\cohaqv^0\coloneqq\bS_z \otimes \bS_y\ ,
\end{align}
where
\begin{align}
	\bS_z\coloneqq \Hbullet_{\Ttilde}[p_\ell(z_i);\, i \in I, \ell \geq 0]\quad\text{and} \quad \bS_y\coloneqq \Hbullet_{\Ttilde}[p_\ell(y_i);\, i \in I, \ell \geq 0]\ .
\end{align}
Here we write $\bS_z$ for $\bS$ in order to avoid confusion.
We equip it with a cocommutative coproduct by setting 
\begin{align}
	\Delta(p_{\ell}(y_i))\coloneqq p_{\ell}(y_i) \otimes 1 + 1 \otimes p_{\ell}(y_i)\quad\text{and} \quad \Delta(p_{\ell}(z_i))\coloneqq p_{\ell}(z_i) \otimes 1 + 1 \otimes p_{\ell}(z_i) 
\end{align}
for any $i,\ell$. We also define an action 
\begin{align}
	\cap \colon \cohaqv^0 \otimes\cohaqv^{\Ttilde} \longrightarrow \cohaqv^{\Ttilde}
\end{align}
as follows: $\ch(z_i)$ acts by multiplication by $\ch(\calV_i)$ and $p_{\ell}(y_i)$ acts by zero. This allows us to form the semi-direct product
\begin{align}
	\cohaqv^{\Ttilde, \leqslant 0} \coloneqq \cohaqv^0 \ltimes \cohaqv^{\Ttilde}\ ,
\end{align}
which is the quotient of the free algebra generated by $\cohaqv^0$ and $\cohaqv^{\Ttilde}$ modulo the two-sided ideal of relations
\begin{align}
	x a=\sum_i x'_i \cap a x''_i\ ,
\end{align}
with $x \in \cohaqv^0$ and $a \in \cohaqv^{\Ttilde}$.

In order to extend $\Phi$ to $\cohaqv^{\Ttilde, \leqslant 0}$, we introduce symbols $\calW_i$ for $i \in I$ which we shall formally treat as vector bundles on $\Lambda_\qv$ satisfying 
\begin{align}
	\ch_\ell(\calW_i)\coloneqq \frac{p_{\ell}(y_i)}{\ell !} \in \cohaqv^{\Ttilde, \leqslant 0}\ .
\end{align}
Let us also set
\begin{align}
	\calF_i\coloneqq q^{-2} \calW_i -(1+q^{-2})\calV_i + \sum_{\genfrac{}{}{0pt}{}{e\in \Omega}{e\colon j\to i}}\gamma_e^{-1} \calV_j + \sum_{\genfrac{}{}{0pt}{}{e\in \Omega}{e\colon i\to j}}\gamma_{e^\ast}^{-1} \calV_j\ .
\end{align}

\begin{theorem}\label{thm:defkappa}
	The assignment
	\begin{align}
		\begin{tikzcd}[ampersand replacement=\&, row sep=small]
			h_i(z)\ar[mapsto]{r} \& 1+ \displaystyle \sum_{\bfv\in \N I}\frac{\lambda_{-1/z}\Big(\calF_i \Big)}{\lambda_{-1/z}\Big(q^2 \calF_i \Big)}\cap -\\
			\kappa_{i, \ell} \ar[mapsto]{r} \&p_{\ell}(y_i) 
		\end{tikzcd}\ ,
	\end{align}
	where $\lambda_{-1/z}(\calE)$ is the \textit{equivariant Chern polynomial} of $\calE$, extends the map $\Phi$ to a surjective algebra homomorphism
	$\Y_\qv^{\frake,\leqslant 0} \to \cohaqv^{\Ttilde, \leqslant 0}$, which restricts to an isomorphism of algebras $\Y_\qv^{\frake,0} \simeq \cohaqv^0$.
\end{theorem}
\begin{proof}
	The only nontrivial relations to check are  \eqref{eq:Yangian-5} and \eqref{eq:Yangian-8}. These follow from direct computations in all points similar to \cite[\S5]{Varagnolo_Yangian}.
\end{proof}
We will often identify $\Y_\qv^{\frake,0}$ with $\bS_z \otimes \bS_y$ without further mention.

\section{Compatibility between braid group action and derived reflection functors}\label{sec:compatibility-braid}

In this section, we relate the action of the truncated braid group on $\Y_\qv$  with the derived reflection functors. Note that Proposition~\ref{prop:reflection-stacks} has the following immediate consequence:
\begin{corollary}\label{cor:Phi-i}
	The surjective algebra homomorphism $\Phi\colon\Y_\qv^-\longrightarrow \cohaqv^{\Ttilde}$ factors to surjective homomorphisms $\iY_\qv \to \icohatildeqv$ and $\Yi_\qv  \to \cohatildeiqv$.
\end{corollary}

Recall from \S\ref{sec:the-algebra-S} the $\bS$-actions on $\cohaqv^{\Ttilde}$ and $\Y_\qv^-$. It is clear that these actions descend to $\bS$-actions on the quotients $\cohatildeiqv$ and $\icohatildeqv$, as well as $\iY_\qv$ and $ \Yi_\qv$ for each $i\in I$. Moreover the morphisms $\iY_\qv \to \icohatildeqv$ and $\Yi_\qv  \to \cohatildeiqv$ are $\bS$-equivariant.

\medskip

As usual, we say that $i \in I$ is a \textit{sink} of $\qv$ if there is no edge $i\to j\in \Omega$, while $i \in I$ is a \textit{source} of $\qv$ if there is no edge $j\to i\in \Omega$.

The main Theorem of this section relates the action of reflection functors on the COHA and the action of braid operators on the Yangian.

\begin{theorem}\label{thm:compatibility-braid-reflection-functors}
	Fix $i\in I$ and assume that $i$ is a sink of $\qv$. Then the diagram
	\begin{align}\label{eq:S_i=T_ione}
		\begin{tikzcd}[ampersand replacement=\&]
			\iY_\qv \arrow[swap]{d}{\overline{T}_i}\arrow[twoheadrightarrow]{r}{\Phi} \& \icohatildeqv \arrow{d}{S_{i,\, \ast}} \\
			 \Yi_\qv \arrow[twoheadrightarrow]{r}{\Phi} \& \coha^{\Ttilde, (i)}_{s_i\qv}
		\end{tikzcd}
	\end{align}
is commutative.
\end{theorem}

The Yangian $\Y_\qv$ is independent of the orientation of $\qv$, and the COHAs $\coha^T_\qv$ and $\coha^T_{s_i\qv}$ are canonically isomorphic. However, these isomorphisms are given by sign twists on each weight space (see \S\ref{sec:sign twists}), and are in fact \textit{not} compatible with the reflection functors. The assumption that $i$ is a source guarantees that there is no sign twist in the diagram. Note that there is a slight abuse of notations in diagram~\eqref{eq:S_i=T_ione}: the map $S_{i,\, \ast}$ is induced in Borel-Moore homology by the composition of equivalences
\begin{align}
	\begin{tikzcd}[ampersand replacement=\&, column sep=large]
		\idLambda_\qv \ar{r}{\R S_i} \& \dLambdai_\qv\ar{r}{u_{\qv,s_i\qv}}\& \dLambdai_{s_i\qv}
	\end{tikzcd}\ .
\end{align}

\medskip

The proof of Theorem~\ref{thm:compatibility-braid-reflection-functors} is divided into four steps:
\begin{itemize}\itemsep0.2cm
	\item In \S\ref{subsec:degree-zero-case} we prove the vertical degree zero version of Theorem~\ref{thm:compatibility-braid-reflection-functors}, i.e., Theorem~\ref{thm:compatibility-braid-reflection-functors-degree-zero}. The proof reduces to the compatibility of the braid group operator $\overline{T}_i$ with the operators at the level of classical Hall algebra of $\qv$ induced by the BGP functor and the use of the characteristic cycle map to lift such a relation at the level of the cohomological Hall algebra of $\qv$.
	
	\item In \S\ref{subsec:tautological-classes-braid} we prove the compatibility between the action of the Chern classes of tautological bundles and the action of the braid group, using the representations of the Yangians on the Borel-Moore homologies of Nakajima quiver varieties.
	
	\item In \S\ref{subsec:support-shuffle-algebras} we characterize the supports of $\icohatilded$ and $\cohatildeid$ as $\shuffle_\bfd$-modules. 
	
	\item \S\ref{subsec:partial-doubles} we show that the operators $\overline{T}_i$ descend to maps $\icohatildeqv \to \cohatildeiqv$. For this we introduce partial doubles of $\cohaqv^{\Ttilde}$ and $\Y_\qv^-$ at a fixed vertex $i\in I$.
\end{itemize}
All these steps are put together in \S\ref{subsec:proof-Theorem-compatibility} where we finally prove Theorem~\ref{thm:compatibility-braid-reflection-functors}.

\medskip

Theorem~\ref{thm:compatibility-braid-reflection-functors} admits a generalization to the case of $\coha^{T,\Theta}_\qv$ for an arbitrary twist $\Theta$. To any $\Llambda \in \coweightlattice_\qv$ one may associate an automorphism $\rho_{\Llambda} \in \Aut(\Y^-_\qv)$ given by 
\begin{align}
	\rho_{\Llambda} \colon u \longmapsto (-1)^{(\Llambda,\vert u\vert)}
\end{align}
and a surjective algebra morphism $\Phi^{\Llambda}_\Theta\colon \Y_\qv \to \cohaqv^{T,\Theta}$ given as 
\begin{align}
	\Phi^{\Llambda}_\Theta\coloneqq\Phi_\Theta \circ \rho_{\Llambda}\ .
\end{align}
There is a natural action of the Weyl group $W_\qv$ on the set of twisting forms $\Theta$, given as follows: if $w \in W$ and $\Theta \colon \Z I \times \Z I \to \Z_2$ is a twist, then $w\cdot \Theta$ is the twist defined by 
\begin{align}
	(\bfe,\bfd) \longmapsto \Theta(w^{-1}(\bfe),w^{-1}(\bfd))\ .
\end{align}

\begin{corollary}\label{cor:Ti=Si-for-Theta}
	Fix $i \in I$ and let $\Theta$ be an arbitrary twist. There exists an (explicit) $\Llambda \in \coweightlattice_\qv$ such that for any $\nu \in \coweightlattice_\qv$ the diagram
	\begin{align}\label{eq:S_i=T_ione theta}
		\begin{tikzcd}[ampersand replacement=\&, column sep=large, row sep=large]
			\iY_\qv \arrow[swap]{d}{\overline{T}_i}\arrow[twoheadrightarrow]{r}{\Phi^{\nu}_{\Theta}} \& \tensor*[^{(i)}]{\coha}{^{\Ttilde, \Theta}_\qv} \arrow{d}{S_{i,\, \ast}}\\ 
			\Yi_\qv \arrow[twoheadrightarrow]{r}{\Phi^{s_i(\nu)+\Llambda}_{s_i\cdot \Theta}} \& \coha^{\Ttilde, s_i\cdot \Theta,(i)}_{s_i\qv}
		\end{tikzcd}
	\end{align}
	is commutative, where $\coha^{\Ttilde, s_i\cdot \Theta,(i)}_{s_i\qv}\coloneqq \sfH_\bullet^{\Ttilde}(\dLambdai_{s_i\qv})$.
\end{corollary}

\begin{proof}
	Let us first assume that $\nu=0$. 
	
	Recall the algebra isomorphisms $\gamma_{\Theta',\Theta''}\colon \cohaqv^{\Ttilde,\Theta'} \simeq \cohaqv^{\Ttilde,\Theta''}$ for any two twists $\Theta'$ and $\Theta''$ (cf.\ Proposition~\ref{prop:sign-twists-I}). Fix an arbitrary twist $\Theta$ and denote by $\gamma_{\qv,\Theta}$ the algebra isomorphism corresponding to the Euler form $\langle -,-\rangle_\qv$ and $\Theta$. Similarly, for any vertex $i\in I$, we shall denote by $\gamma_{s_i \qv,s_i\cdot \Theta}$ the algebra isomorphism corresponding to the Euler form $\langle -,-\rangle_{s_i\qv}$ and $s_i\cdot\Theta$.
	
	 Let us fix the orientation of $\qv$ so that $i$ is a sink, and extend the commutative diagram \eqref{eq:S_i=T_ione} as follows:
	\begin{align}\label{eq:S_i=T_i-one-proof-theta}
		\begin{tikzcd}[ampersand replacement=\&, column sep=huge, row sep=large]
			\iY_\qv \arrow[swap]{d}{\overline{T}_i}\arrow[twoheadrightarrow]{r}{\Phi} \& \icohatildeqv \arrow{d}{S_{i,\, \ast}} \arrow{r}{\gamma_{\qv,\Theta}} \& \tensor*[^{(i)}]{\coha}{^{\Ttilde, \Theta}_\qv} \arrow{d}{S_{i,\, \ast}}\\
			\Yi_\qv \arrow[twoheadrightarrow]{r}{\Phi} \& \coha^{\Ttilde, (i)}_{s_i\qv} \arrow{r}{\gamma_{s_i\qv,s_i\cdot \Theta}} \& \coha^{\Ttilde, s_i\cdot \Theta,(i)}_{s_i\qv}
		\end{tikzcd}\ .
	\end{align}
	The composition 
	\begin{align}
		\begin{tikzcd}[ampersand replacement=\&,column sep=huge]
			\Y_\qv^ -\ar{r}{\Phi}\& \coha^{\Ttilde}_{s_i\qv} \ar{r}{\gamma_{s_i\qv,s_i\cdot \Theta}} \& \coha^{\Ttilde, s_i\cdot \Theta}_{s_i\qv}
		\end{tikzcd}
	\end{align}
	is an algebra morphism, which coincides with $\Phi_{s_i\cdot \Theta}$ up to a sign twist. It thus has the form $\Phi_{s_i\cdot \Theta}^{\Llambda}$ for a well-defined $\Llambda$, whose explicit form may be computed. 
	
	The case of $\nu \neq 0$ follows by a similar argument.
\end{proof}

\subsection{Degree zero case}\label{subsec:degree-zero-case} 

Let $\eYzero$ be the subalgebra $\eY_\qv$ consisting of all elements of vertical degree zero. Set $\Yzero^-\coloneqq \Y_\qv^-\cap \eYzero$. Let also consider the subalgebra $\cohazero^{\Ttilde}\coloneqq H^{\Ttilde}_\mathsf{top}( \dLambda_\qv )$ of $\cohaqv^{\Ttilde}$. The restriction of the map $\Phi$ to the subalgebra $\Yzero^-$ of $\Y_\qv^-$ factorizes to a map 
\begin{align}\label{eq:Phizero}
	\Phi_0\colon \Yzero^-\longrightarrow \cohazero^{\Ttilde}\ ,
\end{align}
as shown in Theorem~\ref{thm:Phi}-\eqref{thm:Phi-3}. Now, introduce
\begin{align}
	\iYzero \coloneqq \Yzero^-/\,\I_i^-\Yzero^- \qquad&\text{and}\qquad \Yizero \coloneqq \Yzero^-/\,\Yzero^-\I_i^- \ , \\[3pt] 
	\icohatildezero\coloneqq  H^{\Ttilde}_\mathsf{top}( \iLambda_\qv )  \qquad&\text{and}\qquad\cohatildeizero \coloneqq H^{\Ttilde}_\mathsf{top}( \Lambdai_\qv ) \ .
\end{align}

By Corollary~\ref{cor:Phi-i}, the map $\Phi_0$ yields a pair of surjective algebra homomorphisms
\begin{align}
	\Phi_0\colon \iYzero \longrightarrow \icohatildezero\quad\text{and}\quad \Phi_0\colon \Yizero \longrightarrow \cohatildeizero\ .
\end{align}
Both maps $\overline T_i$ and $S_{i,\, \ast}$ preserve the degree 0 subspaces. The first step to prove Theorem~\ref{thm:compatibility-braid-reflection-functors} is the following result.
\begin{theorem}\label{thm:compatibility-braid-reflection-functors-degree-zero}
	Fix $i\in I$ a sink of $\qv$. Then, the diagram 
	\begin{align}\label{eq:diag:T_i=S_i}
		\begin{tikzcd}[ampersand replacement=\&]
			\iYzero \arrow{r}{\overline T_i}\arrow[twoheadrightarrow]{d}{\Phi_0} \& \Yizero \arrow[twoheadrightarrow]{d}{\Phi_0} \\
			\icohatildezero \arrow{r}{S_{i,\, \ast}} \& \cohatildeizero
		\end{tikzcd}
	\end{align}
	is commutative.
\end{theorem}To prove  the theorem we proceed in the following way. We first consider the analogous statement for $\F_q$-Hall algebras $\bfH_{\qv,\varepsilon}$, which essentially follows from work of Lusztig and Sevenhant-Van den Bergh on braid group actions and BGP functors. We next use the Cebotarev density theorem and the trace map to lift this result at the level of categories of perverse sheaves $\sfK_0(\frakQ_\qv)$, first over $\F_q$ with $\overline{\Q_\ell}$-coefficients, then over $\C$ with $\C$-coefficients. Finally we apply the characteristic cycle map $\mathsf{CC}\colon \sfK_0(\frakQ_{\qv,\C}) \to \cohazero$ and use the known fact, see \cite{Hennecart_Geometric}, that it is an algebra homomorphism of algebras.
\subsubsection{The BGP functor} 

Fix a field $k$. For $i$ a sink of $\qv$ we let $\catmod(\qv)^i$ be full subcategory of $\catmod(\qv)$ whose objects are $(V, x)$ such that the map 
\begin{align}\label{eq:surjective-BGP}
	\bigoplus_{e\colon j\to i} x_e \colon \bigoplus_j V_j \to V_i
\end{align}
is surjective. Similarly, for $i$ a source of $\qv$ we let $\catmod(\qv)_i$ be full subcategory of $\catmod(\qv)$ whose objects are $(V, x)$ such that the map 
\begin{align}\label{eq:injective-BGP}
	\bigoplus_{e\colon i\to j} x_e\colon V_i\to \bigoplus_j V_j
\end{align}
is injective. Here, in the direct sum in the right-hand-side of the map \eqref{eq:surjective-BGP} and in the direct sum in the left-hand-side of the map \eqref{eq:injective-BGP}, $j$ runs over all vertices of $\qv$ that are adjacent to $i$.

For a fixed vertex $i$, let $s_i\qv$ be the quiver obtained by reversing all the edges which end at $i$. In particular, if $i$ is a sink (resp. source) of $\qv$, $i$ is a source (resp.\ sink) of $s_i\qv$. We shall now recall the construction of the \textit{BGP functor} following \cite{BGP}. Assume that $i$ is a sink of $\qv$. To a representation $(V, x)$  we associate the representation $(\widetilde{V}, \widetilde{x})$, where
\begin{align}
	\widetilde{V}_j\coloneqq \begin{cases}
		V_j & \text{if $j\neq i$}\ ,\\
		\displaystyle \ker\Big( \bigoplus_{e\colon j\to i} x_e \colon \bigoplus_j V_j \to V_i \Big) & \text{if $j=i$} \ ,
	\end{cases}
\end{align}
and the linear maps $\widetilde{x}$ are the obvious ones. Similarly, if we assume that $i$ a source of $\qv$, to a representation $(V, x)$ we associate the representation $(\widetilde{V}, \widetilde{x})$, where
\begin{align}
	\widetilde{V}_j\coloneqq \begin{cases}
		V_j & \text{if $j\neq i$}\ ,\\
		\displaystyle \mathsf{coker}\Big( \bigoplus_{e\colon i\to j}\colon V_i\to \bigoplus_j V_j \Big) & \text{if $j=i$} \ ,
	\end{cases}
\end{align}
and the linear maps $\widetilde{x}$ are the obvious ones. This construction gives rise to the functor
\begin{align}
	S_i\colon \catmod(\qv)\longrightarrow \catmod(s_i\qv)
\end{align}
if $i$ is a sink, and to the functor
\begin{align}
	S_i'\colon  \catmod(\qv)\longrightarrow \catmod(s_i\qv)
\end{align}
if $i$ is a source. In particular, if $i$ is a sink, the functor $S_i$ descends to an equivalence
\begin{align}
	S_i\colon \catmod(\qv)^i\longrightarrow \catmod(s_i\qv)_i\ ,
\end{align}
whose inverse is given by the restriction of $S'_i$. 

From now on, we let $i$ be a sink of $\qv$. Let 
\begin{align}\label{eq:embedding-f-sink}
	\iii{f}\colon \stackRep(\qv)^{(i)}\longrightarrow \stackRep(\qv)
\end{align}
be the embedding of the open substack consisting of finite-dimensional representations of $\qv$ for which \eqref{eq:surjective-BGP} is surjective; likewise, let 
\begin{align}\label{eq:embedding-f-source}
	f^{(i)}\colon \stackRep(s_i\qv)_{(i)} \longrightarrow \stackRep(s_i\qv)
\end{align}
be the embedding of the open substack consisting of representations of $s_i\qv$ for which \eqref{eq:injective-BGP} is injective. The BGP functor $S_i$ induces a morphism
\begin{align}\label{eq:Si-functor}
	S_i\colon \stackRep(\qv)\longrightarrow \stackRep(s_i\qv)\ .
\end{align}
It restricts to an equivalence
\begin{align}\label{eq:Si-equivalence}
	\begin{tikzcd}[ampersand replacement=\&]
		S_i\colon \stackRep(\qv)^{(i)}\arrow{r}{\sim}\& \stackRep(s_i\qv)_{(i)}
	\end{tikzcd} \ ,
\end{align}
whose inverse is given by $S'_i$.

\subsubsection{Classical Hall algebras and braid group action} 

Let $\upsilon$ be a formal variable and let $\sfU_\upsilon(\frakn_\qv)$ be Lusztig's integral form of the quantized enveloping algebra of $\frakn_\qv$ (see e.g. \cite[Appendix~4]{Schiffmann-lectures-Hall-algebras} for the definition). Now, let $q$ be a prime power and let us fix a square root $\varepsilon$ of $q$. Set
\begin{align}
	\sfU_\varepsilon(\frakn_\qv)\coloneqq \sfU_\upsilon(\frakn_\qv)\vert_{\upsilon=\varepsilon}\ .
\end{align}
We also set
\begin{align}
	\iU_\upsilon(\frakn_\qv)\coloneqq \sfU_\upsilon(\frakn_\qv)/f_iU_\upsilon(\frakn_\qv) \quad\text{and}\quad \Ui_\upsilon(\frakn_\qv)\coloneqq \sfU_\upsilon(\frakn_\qv)/\sfU_\upsilon(\frakn_\qv)f_i\ ,
\end{align}
where $f_i \in \sfU_{\varepsilon}(\frakn_\qv)$ is the $i$th Chevalley generator. We define $\iU_\varepsilon(\frakn_\qv), \Ui_\varepsilon(\frakn_\qv)$ in a similar fashion. 

Let $\stackRep_\qv(\F_q)$ be the set of isomorphism classes of $\F_q$-points of $\stackRep(\qv)$ and let $\sfFun_{\F_q}(\qv)$ be the space of all maps $\stackRep_\qv(\F_q)\to\C$. We endow it with a \textit{twisted} Hall multiplication, where the twist is defined as e.g. \cite[Lecture~1]{Schiffmann-lectures-Hall-algebras}. The \textit{composition subalgebra} $\bfH_{\qv,\varepsilon}\subset \sfFun_{\F_q}(\qv)$ is the subalgebra generated by the characteristic functions $1_i\colon \stackRep_{\alpha_i}(\F_q)\to\C$ for $i \in I$. The Ringel-Green Theorem provides an isomorphism 
\begin{align}\label{eq:ringel-green-theorem}
	\begin{tikzcd}[ampersand replacement=\&]
		\theta_\varepsilon\colon\bfH_{\qv,\varepsilon} \arrow{r}{\sim}\& \sfU_{\varepsilon}(\frakn_\qv)
	\end{tikzcd}
\end{align}
see e.g. \cite[Theorem~3.16]{Schiffmann-lectures-Hall-algebras}\footnote{Note that $\sfU_{\varepsilon}(\frakn_\qv)$ and $\sfU_{\varepsilon}(\frakn^+_\qv)$ are isomorphic as algebras, where $\frakn_\qv^+$ is the positive nilpotent Lie subalgebra of $\frakg_\qv$.}.

Recall that $i$ is a sink. The embeddings of stacks \eqref{eq:embedding-f-sink} and \eqref{eq:embedding-f-source} induces inclusions $\iii{f}$ and $f^{(i)}$ of the sets $\stackRep_\qv(\F_q)^{(i)}$ and $\stackRep_{s_i\qv}(\F_q)_{(i)}$ of isomorphism classes of $\F_q$-points of the stacks $\stackRep(\qv)^{(i)}$ and $\stackRep(s_i\qv)_{(i)}$ into $\stackRep_\qv(\F_q)$ and $\stackRep_{s_i\qv}(\F_q)$, respectively. This gives rise to the following quotients
\begin{align}
	\iii{\bfH}_{\qv,\varepsilon}\coloneqq (\iii{f})^\ast(\bfH_{\qv,\varepsilon})\qquad\text{and}\qquad \bfH^{(i)}_{s_i\qv,\varepsilon}\coloneqq (f^{(i)})^\ast(\bfH_{s_i\qv,\varepsilon})
\end{align}
of the vector spaces $\bfH_{\qv,\varepsilon}$ and $\bfH_{s_i\qv,\varepsilon}$, respectively. The map $\theta$ descends to linear isomorphisms
\begin{align}\label{eq:ithetai}
	\begin{tikzcd}[ampersand replacement=\&]
		\iii{\theta}_\varepsilon\colon 	\iii{\bfH}_{\qv,\varepsilon}\arrow{r}{\sim}\& \iU_{\varepsilon}(\frakn_\qv) 
	\end{tikzcd}   \qquad \text{and}\qquad
	\begin{tikzcd}[ampersand replacement=\&]
			\theta^{(i)}_\varepsilon\colon \bfH^{(i)}_{s_i\qv,\varepsilon} \arrow{r}{\sim}\& \Ui_{\varepsilon}(\frakn_\qv)
	\end{tikzcd} \ .
\end{align}
Indeed, the isomorphism $\theta_\varepsilon$ maps $1_i$ to $f_i$, hence it restricts to isomorphisms 
\begin{align}
	\begin{tikzcd}[ampersand replacement=\&]
		1_i\cdot \bfH_{\qv,\varepsilon} \arrow{r}{\sim}\& f_iU_{\varepsilon}(\frakn_\qv)
	\end{tikzcd}   \qquad \text{and}\qquad
	\begin{tikzcd}[ampersand replacement=\&]
		\bfH_{\qv,\varepsilon}\cdot1_i  \arrow{r}{\sim}\& \sfU_{\varepsilon}(\frakn_\qv)f_i 
	\end{tikzcd} \ .
\end{align}
Moreover, a function $g \in \sfFun_{\F_q}(\qv)$ is supported on the complement of $\stackRep_\qv(\F_q)^{(i)}$ (resp. of  $\stackRep_{s_i\qv}(\F_q)_{(i)}$) if and only if it is left (resp.\ right) divisible by $1_i$, and the same holds for $\bfH_{\qv,\varepsilon}$ by e.g. \cite[Theorem~14.3.2]{Lusztig_Quantum_groups}.

The functor $S_i$ introduced in Formula~\eqref{eq:Si-equivalence} induces a bijection 
\begin{align}
	\begin{tikzcd}[ampersand replacement=\&]
		\stackRep_\qv(\F_q)^{(i)} \ar{r}{\sim} \& \stackRep_{s_i\qv}(\F_q)_{(i)}
	\end{tikzcd}\ .
\end{align}
The direct image by $S_i$ yields a linear isomorphism 
\begin{align}
	\begin{tikzcd}[ampersand replacement=\&]
		S_{i,\ast}\colon \iii{\bfH}_{\qv, \varepsilon}\arrow{r}{\sim}\& \bfH^{(i)}_{s_i\qv,\varepsilon}
	\end{tikzcd} 
\end{align}
such that the following diagram of linear isomorphisms commutes
\begin{align}\label{eq:braid-BGP}
	\begin{tikzcd}[ampersand replacement=\&]
		\iii{\bfH}_{\qv,\varepsilon}  \ar{d}{S_{i,\ast}}\ar{r}{\iii{\theta}}\&
		\iU_\varepsilon(\frakn_\qv) \ar{d}{T'_{i,-1}} \\
		\bfH^{(i)}_{s_i\qv,\varepsilon}\ar{r}{\theta^{(i)}} \&
		\Ui_{\varepsilon}(\frakn_\qv)
	\end{tikzcd}\ .
\end{align}
This is shown in \cite[\S38, Lemma~38.1.3]{Lusztig_Quantum_groups}, see also \cite{VdB-S99} and \cite[Theorem~6.3-(a)]{XY-BGP} . Here, $T'_{i,-1}$ is Lusztig's \textit{braid group operator} considered in \cite[\S8]{Lusztig-Canonical-bases-Hall-algebras} and \cite[\S38]{Lusztig_Quantum_groups}. More precisely, it corresponds to Lusztig's braid group operator twisted by the Cartan involution swapping $e_i$ and $f_i$. The classical limit of $T'_{i,-1}$ is the standard operator $T_i\in B_\qv$. 
In the following paragraphs, we will successively prove variants of \eqref{eq:braid-BGP}: first at the level of categories of perverse sheaves over $\F_q$, then over $\C$, and finally at the level of $\coha_0^{\Ttilde}$. 

\subsubsection{Lusztig's categorification of the Hall algebra (${\F_q}$-case)} 

Let $k=\F_q$. Let $\catDb_m(\stackRep(\qv))$ be the bounded derived category of constructible mixed complexes (i.e., complexes of $\F_q$-Weil sheaves) on $\stackRep(\qv)$. We shall consider the usual (complexified) Grothendieck group $\sfK_0(\catDb_m(\stackRep(\qv)))$ and the (complexified) \textit{split} Grothendieck group $\sfK_\oplus(\catDb_m(\stackRep(\qv)))$. It has the structure of a $\C[\upsilon,\upsilon^{-1}]$-module given by $\upsilon^n \cdot [C]\coloneqq [C[n]]$ for any $n\in \Z$ (i.e., $\upsilon$ acting as the grading shift operator).

Let $\frakQ_{\qv,k}\subset \catDb_m(\stackRep(\qv))$ be the \textit{Lusztig category} introduced, e.g., in \cite[\S2 \& \S3]{Lusztig_Quivers}. Its complexified split Grothendieck group $\calK_{\qv,k}\coloneqq \sfK_\oplus(\frakQ_{\qv,k})$ is a free $\C[\upsilon,\upsilon^{-1}]$-module. By \cite[Theorem~10.17]{Lusztig_Quivers}, there is a $\C[\upsilon,\upsilon^{-1}]$-algebra structure on $\calK_{\qv,k}$ and a canonical isomorphism
\begin{align}
	\begin{tikzcd}[ampersand replacement=\&]
		\Psi_k\colon \calK_{\qv,k}\arrow{r}{\sim}\& \sfU_\upsilon(\frakn_\qv)
	\end{tikzcd}\ .
\end{align}

Let $\iii{\frakQ}_{\qv,k}$ and $\frakQ_{s_i\qv,k}^{(i)}$ be the essential images  of the categories $\frakQ_{\qv,k}$ and $\frakQ_{s_i\qv,k}$ by the pull-back functors induced by $\iii{f}$ and $f^{(i)}$, respectively. Set
\begin{align}
	\iii{\calK}_{\qv,k}\coloneqq \sfK_\oplus\big(\iii{\frakQ}_{\qv,k}\big)\quad \text{and} \quad \calK^{(i)}_{s_i\qv,k}\coloneqq \sfK_\oplus\big(\frakQ_{s_i\qv,k}^{(i)}\big)\ .
\end{align}
These are free $\C[\upsilon,\upsilon^{-1}]$-modules generated by the classes of shifts of simple perverse sheaves in $\frakQ_{\qv,k}$ whose support intersect $\iii{\stackRep}(\qv)$ and $\stackRep^{(i)}(\qv)$ respectively.

By an argument of Lusztig, see, e.g., \cite[Lemma~3.19]{Schiffmann-lectures-canonical-bases}, the kernel of the restriction maps $\calK_{\qv,k} \to \iii{\calK}_{\qv,k}$ and $\calK_{\qv,k} \to \calK^{(i)}_{\qv,k}$ are respectively the left and right ideals generated by $[\mathsf{IC}(\stackRep_{\alpha_i}(\qv))]=\Psi_k^{-1}(f_i)$. Thus $\Psi_k$ factors to isomorphisms
\begin{align}\label{eq:ipsik}
	\begin{tikzcd}[ampersand replacement=\&]
		\iii{\Psi}_k\colon \iii{\calK}_{\qv,k}  \arrow{r}{\sim}\& \iU_\upsilon(\frakn_\qv)
	\end{tikzcd}
	 \quad \text{and} \quad
	\begin{tikzcd}[ampersand replacement=\&]
		\Psi_k^{(i)}\colon \calK^{(i)}_{\qv,k} \arrow{r}{\sim}\& \Ui_\upsilon(\frakn_\qv)
	\end{tikzcd}\ .
\end{align}
The aim of this paragraph is to prove the following result:
\begin{proposition}\label{prop:sigma-equivalence}
	The direct image by $S_i$ yields an equivalence of graded additive categories $\iii{\frakQ}_{\qv,k}\simeq \frakQ_{s_i\qv,k}^{(i)}$. Taking the split Grothendieck groups we get a commutative square of linear isomorphisms
	\begin{align}\label{eq:sigma-equivalence}
		\begin{tikzcd}[ampersand replacement=\&]
			\iii{\calK}_{\qv,k} \ar[swap]{d}{S_{i,\ast}}\ar{r}{\iii{\Psi}_{k}} \& \iU_\upsilon(\frakn_\qv)  \ar{d}{T'_{i,-1}}\\
			\calK_{s_i\qv,k}^{(i)} \ar{r}{\Psi^{(i)}_{k}}\&  \Ui_\upsilon(\frakn_\qv)
		\end{tikzcd}\ .
	\end{align}
\end{proposition}

To prove the proposition, we need some preliminary result. First, we shall use the trace map to deduce this from \eqref{eq:braid-BGP}. By \cite{Lusztig-Canonical-bases-Hall-algebras}, see e.g. \cite[Theorem~3.25]{Schiffmann-lectures-canonical-bases}, the trace of the geometric Frobenius yields an algebra isomorphism 
\begin{align}\label{eq:Lusztig-trace-isomorphism}
	\begin{tikzcd}[ampersand replacement=\&]
		\tau_\varepsilon\colon \calK_{\qv,\varepsilon}\coloneqq\calK_{\qv,k} \vert_{\upsilon=\varepsilon}\arrow{r}{\sim}\&\bfH_{\qv,\varepsilon}
	\end{tikzcd}\ .
\end{align}
Note that the composition $\theta_\varepsilon\circ \tau_\varepsilon$ coincides with $\Psi_k\vert_{\upsilon=\varepsilon}$.

If $i$ is a sink, we likewise set
\begin{align}
	\iii{\calK}_{\qv,\varepsilon}\coloneqq \iii{\calK}_{\qv,k}\vert_{\upsilon=\varepsilon} \quad \text{and} \quad \calK^{(i)}_{s_i\qv,\varepsilon}\coloneqq \calK^{(i)}_{s_i\qv,k}\vert_{\upsilon=\varepsilon}\ .
\end{align}

\begin{lemma}
	The trace of the geometric Frobenius $\tau_\varepsilon$ yields linear isomorphisms
	\begin{align}\label{eq:itaui}
		\begin{tikzcd}[ampersand replacement=\&]
			\iii{\tau}_\varepsilon\colon \iii{\calK}_{\qv,\varepsilon} \arrow{r}{\sim}\& \iii{\bfH}_{\qv,\varepsilon}
		\end{tikzcd}
		\quad\text{and}\quad 
		\begin{tikzcd}[ampersand replacement=\&]
			\tau^{(i)}_\varepsilon\colon \calK_{s_i\qv,\varepsilon}^{(i)} \arrow{r}{\sim}\& \bfH^{(i)}_{s_i\qv,\varepsilon}
		\end{tikzcd} \  .
	\end{align} 
\end{lemma}
	
\begin{proof}
	Note that the following diagrams
	\begin{align}\label{eq:commutative-diagram}
		\begin{tikzcd}[ampersand replacement=\&]
			\calK_{\qv,\varepsilon} \ar{d}{\tau_\varepsilon}\ar{r}{\mathsf{res}}\&\iii{\calK}_{\qv,\varepsilon}\ar{d}{\iii{\tau}_\varepsilon}\\
			\bfH_{\qv,\varepsilon}\ar{r}{(\iii{f})^\ast} \&\iii{\bfH}_{\qv,\varepsilon}
		\end{tikzcd}
		\quad \text{and}\quad
		\begin{tikzcd}[ampersand replacement=\&]
			\calK_{s_i\qv,\varepsilon} \ar{d}{\tau_\varepsilon}\ar{r}{\mathsf{res}} \&\calK_{s_i\qv,\varepsilon}^{(i)}\ar{d}{\tau^{(i)}_\varepsilon}\\
			\bfH_{s_i\qv,\varepsilon} \ar{r}{(f^{(i)})^\ast}\&\bfH^{(i)}_{s_i\qv,\varepsilon}
		\end{tikzcd}
	\end{align}
	commute and the maps $(\iii{f})^\ast\circ\tau_\varepsilon$ and $(f^{(i)})^\ast\circ\tau_\varepsilon$ are surjective. Hence, $\iii{\tau}_\varepsilon$ and $\tau^{(i)}_\varepsilon$ are surjective as well. The injectivity follows from the fact that the domains and targets have the same graded dimension. Indeed, $\iii{\calK}_{\qv,\varepsilon}$ and $\iii{\bfH}_{\qv,\varepsilon}$ both have the same graded dimension as $\iU_{\varepsilon}(\frakn_\qv)$; the same argument works in the second case.
\end{proof}
	
\begin{proof}[Proof of Proposition~\ref{prop:sigma-equivalence}] 
	Let us consider the following commutative diagram of linear isomorphisms:
	\begin{align}\label{eq:braid-BGP-2}
		\begin{tikzcd}[ampersand replacement=\&]
			\iii{\calK}_{\qv,\varepsilon} \ar{r}{\iii{\tau}_\varepsilon}\& \iii{\bfH}_{\qv,\varepsilon}  \ar{d}{S_{i,\ast}}\ar{r}{\iii{\theta}_\varepsilon}\&
			\iU_{\varepsilon}(\frakn_\qv) \ar{d}{T'_{i,-1}} \\
			\calK_{s_i\qv,\varepsilon}^{(i)} \ar{r}{\tau^{(i)}_\varepsilon}\& \bfH^i_{s_i\qv,\varepsilon}\ar{r}{\theta^{(i)}_\varepsilon} \&
			\Ui_{\varepsilon}(\frakn_\qv)
		\end{tikzcd}\ .
	\end{align}
	The functor
	\begin{align}
		\begin{tikzcd}[ampersand replacement=\&]
			S_{i,\ast}\colon \catDb_m( \iii{\stackRep}(\qv))  \ar{r}{\sim} \& \catDb_m(\stackRep(s_i\qv))^{(i)}
		\end{tikzcd}
	\end{align}
	is an equivalence. Since both categories $\iii{\frakQ}_{\qv,k}$ and $\frakQ_{s_i\qv,k}^{(i)}$ consist of semisimple complexes and are stable by direct summands, it is enough to prove that $S_{i,\ast}$ restricts to an isomorphism $\iii{\calK}_{\qv,k}  \simeq \calK^{(i)}_{s_i\qv,k}$, which fits into a commutative diagram of the form \eqref{eq:sigma-equivalence}. 
	
	Now, the trace of Frobenius may not be injective on $\sfK_\oplus(\catDb_m(\ii{\stackRep}_{\qv,k}) )\vert_{\upsilon=\varepsilon}$, so we may not directly define 
	\begin{align}
		\begin{tikzcd}[ampersand replacement=\&]
			S_{i,\ast}\colon \sfK_\oplus(\iii{\frakQ}_{\qv,k})  \ar{r}{\sim} \& \sfK_\oplus(\frakQ_{s_i\qv} ^{(i)})
		\end{tikzcd}
	\end{align}
	such that the diagram \eqref{eq:sigma-equivalence} commutes by completing the the diagram \eqref{eq:braid-BGP-2}. To overcome this problem, we shall consider all extensions of the field $\F_q$ in $k$ at once and use the Cebotarev density theorem.
	
	For any $n\in \Z_{> 0}$, considering the set $\stackRep_\qv(\F_{q^n})$ instead of $\stackRep_\qv(\F_q)$, we define the algebra $\bfH_{\qv,\varepsilon^n}$ as before. Consider the subalgebra $\bfH_{\qv, \upsilon}$ of $\prod_{n>0} \bfH_{\qv,\varepsilon^n}$ generated by the families of characteristic functions of the simple modules. We also see $\sfU_\upsilon(\frakn_\qv)$ as the subalgebra of $\prod_{n>0} \sfU_{\varepsilon^n}(\frakn_\qv)$ via the diagonal embedding. Considering the maps $\theta_{\varepsilon^n}$ in Formula~\eqref{eq:ringel-green-theorem} and $\tau_{\varepsilon^n}$ in Formula~\eqref{eq:Lusztig-trace-isomorphism} for each extension $k/\F_{q^n}$, we recover the canonical $\C[\upsilon,\upsilon^{-1}]$-algebra isomorphism
	\begin{align}
		\Psi_{k}\colon \calK_{\qv,k}\longrightarrow \sfU_\upsilon(\frakn_\qv) \ .
	\end{align}
	Similarly, we define maps $\iii{\theta}$ and $\theta^{(i)}$ as in Formula~\eqref{eq:ithetai} and $\iii{\tau}$ and $\tau^{(i)}$ as in Formula~\eqref{eq:itaui} for $ \calK_{\qv,k}$ and $\bfH_{\qv,k}$. Thus, we get the following commutative diagram in which all maps are isomorphisms
	\begin{align}\label{eq:braid-commutativity-closure}
		\begin{tikzcd}[ampersand replacement=\&]
			\iii{\calK}_{\qv,k} \ar{r}{\iii{\tau}}\& \iii{\bfH}_{\qv,\upsilon}  \ar{d}{S_{i,\ast}}\ar{r}{\iii{\theta}}\&
			\iU_\upsilon(\frakn_\qv) \ar{d}{T'_{i,-1}} \\
			\calK_{s_i\qv,k}^{(i)} \ar{r}{\tau^{(i)}}\& \bfH^{(i)}_{s_i\qv,\upsilon}\ar{r}{\theta^{(i)}} \&
			\Ui_\upsilon(\frakn_\qv)
		\end{tikzcd}\ .
	\end{align}
	By the Cebotarev density theorem, the products of all traces of Frobenius for $n\in \N$ is injective on $\sfK_\oplus(\catDb_m(\ii{\stackRep}_{\qv,k}))$. Furthermore, the trace map commutes with the functor $S_{i,\ast}$ so the map $\iii{\tau}$ and $\tau^{(i)}$ are injective.	We thus deduce that the functor $S_{i,\ast}$ restricts to an equivalence $\iii{\frakQ}_{\qv,k}\simeq \frakQ_{s_i\qv,k}^{(i)}$ yielding the desired isomorphism 
	\begin{align}
		\begin{tikzcd}[ampersand replacement=\&]
			S_{i,\ast}\colon  \iii{\calK}_{\qv,k}\ar{r}{\sim} \& \calK^{(i)}_{s_i\qv,k}
		\end{tikzcd}\ .
	\end{align}
	Thus, the claim is proved.
\end{proof}

\subsubsection{Lusztig's categorification of the Hall algebra ($\C$-case)}\label{sec:LCH2}

Let $k=\C$. We denote by $\frakQ_{\qv, \C}$ Lusztig's category of semisimple complexes on $\stackRep(\qv)$. It is again an additive, graded subcategory of $\catDb_m(\stackRep_{\qv, \C})$. We shall denote by $\calK_{\qv,\C}$ its complexified Grothendieck group, a free $\C[\upsilon, \upsilon^{-1}]$-module. We have a $\C[\upsilon, \upsilon^{-1}]$-algebra isomorphism
\begin{align}
	\begin{tikzcd}[ampersand replacement=\&]
			\Psi_\C\colon \calK_{\qv, \C} \ar{r}{\sim}\& \sfU_\upsilon(\frakn_\qv)
	\end{tikzcd}\ .
\end{align} 

We define the categories $\iii{\frakQ}_{\qv,\C}, \frakQ^{(i)}_{s_i\qv,\C}$ as before, as well as their respective complexified Grothendieck groups $\iii{\calK}_{\qv,\C}, \calK_{s_i\qv,\C}^{(i)}$. Again, by e.g. \cite[Lemma~3.19]{Schiffmann-lectures-canonical-bases} (which applies also over $\C$), the map $\Psi_\C$ descends to algebra isomorphisms
\begin{align}
	\begin{tikzcd}[ampersand replacement=\&]
		\iii{\Psi}_\C\colon \iii{\calK}_{\qv,\C} \ar{r}{\sim} \& \iU_\upsilon(\frakn_\qv)
	\end{tikzcd}
	\quad \text{and} \quad  
	\begin{tikzcd}[ampersand replacement=\&]
		\Psi^{(i)}_\C\colon \calK^{(i)}_{s_i\qv,\C} \ar{r}{\sim}\& \Ui_\upsilon(\frakn_\qv)
	\end{tikzcd}\ .
\end{align}
Composing these isomorphisms with the isomorphisms appearing in Formula~\eqref{eq:ipsik}, we get commuting squares
\begin{align}
	\begin{tikzcd}[ampersand replacement=\&]
		\calK_{\qv,k} \ar[swap]{d}{}\ar{r}{{\Theta}} \& \calK_{\qv,\C}  \ar{d}{}\\
		\iii{\calK}_{\qv,k}\ar{r}{{\iii{\Theta}}} \& \iii{\calK}_{\qv,\C}
	\end{tikzcd} 
	\quad \text{and} \quad
	\begin{tikzcd}[ampersand replacement=\&]
		\calK_{s_i\qv,k} \ar[swap]{d}{}\ar{r}{{\Theta}} \& \calK_{s_i\qv,\C}  \ar{d}{}\\
		{\calK}^{(i)}_{s_i\qv,k}\ar{r}{\Theta^{(i)}} \& {\calK}^{(i)}_{s_i\qv,\C}
	\end{tikzcd}\ .
\end{align}

\begin{lemma}\label{lem:S_i-k-vs-C}
	The functor $S_{i,\ast}$ induces an equivalence of graded categories $\iii{\frakQ}_{\qv,\C} \stackrel{\sim}{\to} \frakQ^{(i)}_{s_i\qv,\C}$ whose class in the Grothendieck group fits in a commutative diagram
	\begin{align}\label{eq:S_i-k-vs-C}
		\begin{tikzcd}[ampersand replacement=\&]
			\iii{\calK}_{\qv,k} \ar[swap]{d}{\iii{\Theta}}\ar{r}{{S_{i,\ast}}} \& \calK^{(i)}_{s_i\qv,k}  \ar{d}{\Theta^{(i)}}\\
			\iii{\calK}_{\qv,\C}\ar{r}{S_{i,\ast}} \& \calK^{(i)}_{s_i\qv,\C}
		\end{tikzcd}\ .
	\end{align}
\end{lemma}

\begin{proof}
	Let $\stackRep_{\Z}(\qv) \to \Spec(\Z)$ be the stack of finite-dimensional representations of $\qv$ over $\Z$, and let $\catDb_c(\stackRep_\Z(\qv))$ be the derived category of constructible complexes of $\overline{\Q}_\ell$-sheaves on $\stackRep_\Z(\qv)$. One may define a graded additive subcategory $\frakQ_{\qv,\Z}$ of semisimple complexes whose restriction to each $\stackRep_{\F_q}(\qv)$ yields $\frakQ_{\qv,\overline{\F_q}}$ and whose restriction to the generic fiber yields, after tensoring with $\C$, the category $\frakQ_{\qv,\C}$. For each fixed dimension vector $\bfd \in \N I$, over a dense open subset $O_\bfd \subset \Spec(\Z)$, the category $\frakQ_{\qv,\Z}$ is additively generated by the collection of Lusztig sheaves $\LL_{\mathbf{i},\Z}$, where $\mathbf{i}$ runs over all $I$-compositions of $\bfd$, see, e.g., \cite[\S 3.5]{Schiffmann-lectures-canonical-bases}.
	
	The reflection functor $S_i$ is defined over $\Z$ and gives rise to an equivalence of open substacks
	\begin{align}
		\begin{tikzcd}[ampersand replacement=\&]
			S_{i,\ast} \colon \iii{\stackRep}_{\Z}(\qv) \ar{r}{\sim}\& \stackRep^{(i)}_{\Z}(s_i\qv)
		\end{tikzcd}\ .
	\end{align}
	By Proposition~\ref{prop:sigma-equivalence}, we have
	\begin{align}\label{eq:S_iLusztigsheaves}
		S_{i,\ast}\big((\LL_{\mathbf{i},k})\vert_{ \iii{\stackRep}_k(\qv)}\big) \in \frakQ^{(i)}_{s_i\qv,k}
	\end{align}
	for any $k=\overline{\F_q}$. The same thus holds over a suitable dense open subset of $\Spec(\Z)$ for any object of $\iii{\frakQ}_{\qv,\Z}$. This implies that \eqref{eq:S_iLusztigsheaves} holds also over the generic point of $\Spec(\Z)$, which gives the first statement of Lemma~\ref{lem:S_i-k-vs-C}. The commutativity of the diagram~\eqref{eq:S_i-k-vs-C} follows by a similar argument: $S_{i,\ast}\big((\LL_{\mathbf{i},k})\vert_{ \iii{\stackRep}_k(\qv)}\big)$ is a certain virtual linear combination of (restrictions of) Lusztig sheaves $\LL_{\mathbf{i'},k}$ with coefficients independent of $k$ by the commutative diagram~\eqref{eq:sigma-equivalence}. The same relation thus holds over $\C$.
\end{proof}

As an immediate corollary, we obtain:
\begin{proposition}\label{lem:sigma-equivalence-for-C}
	There is a commutative square of linear isomorphisms
	\begin{align}\label{eq:sigma-equivalence-for-C}
		\begin{tikzcd}[ampersand replacement=\&]
			\iii{\calK}_{\qv,\C} \ar[swap]{d}{S_{i,\ast}}\ar{r}{\iii{\Psi}_{\C}} \& \iU_\upsilon(\frakn_\qv)  \ar{d}{T'_{i,-1}}\\
			\calK_{s_i\qv,\C}^{(i)} \ar{r}{\Psi^{(i)}_{\C}}\&  \Ui_\upsilon(\frakn_\qv)
		\end{tikzcd}\ .
	\end{align}
\end{proposition}

\subsubsection{Characteristic cycle map}\label{sec:CCmap}

We shall now use the characteristic cycle map as a way to relate the Grothendieck group of $\frakQ_{\qv,\C}$ with the degree zero piece $\cohazero$ of the nilpotent quiver COHA. 

Note that $\sfK_0(\frakQ_{\qv,\C}) =(\calK_{\qv,\C})\vert_{\upsilon=-1}$. We shall need the following result due to Hennecart (for standard properties of characteristic cycle map, see \cite[\S2.1]{Hennecart_Geometric} and references therein). Recall that the twist $\calA^{\langle-,-\rangle}$ of an $\N I$-graded associative algebra $\calA$ is the algebra having the same underlying vector space but twisted multiplication
\begin{align}
	m^{\langle-,-\rangle}\Big\vert_{\calA_\bfd \otimes \calA_{\bfe}} \coloneqq (-1)^{\langle \bfd ,\bfe\rangle} m\vert_{\calA_\bfd \otimes \calA_{\bfe}}\ .
\end{align}
\begin{theorem}[{\cite[Theorem~1.2 and Corollary~9.3]{Hennecart_Geometric}}]\label{thm:cc-map}
	The characteristic cycle map 
	\begin{align}
		\mathsf{CC}\colon \sfK_0(\frakQ_{\qv, \C})^{\langle-,-\rangle} \longrightarrow \cohazero
	\end{align}
	is a $\C$-algebra isomorphism, which sends the generator $[\mathsf{IC}(\stackRep(\qv; \alpha_i))]$ to the generator $[\Lambda_\qv(\alpha_i)]$.
\end{theorem}

\begin{corollary}
	There is a sequence of algebra isomorphisms
	\begin{align}
		\begin{tikzcd}[ampersand replacement=\&]
			\mathsf{CC}'\colon (\calK_{\qv,\C})\Big\vert_{\upsilon=1} \arrow{r}{\sim} \& U(\frakn_\qv)  \arrow{r}{\sim} \& \sfU_{-1}(\frakn_\qv)^{\langle-,-\rangle} \arrow{r}{\sim} \&(\calK_{\qv,\C})^{\langle-,-\rangle}\Big\vert_{\upsilon=-1} \arrow{r}{\sim} \&\coha_0 
		\end{tikzcd}\ .
	\end{align}
\end{corollary}

\begin{proof} 
	The isomorphism $\sfU(\frakn_\qv) \simeq \sfU_{-1}(\frakn_\qv)^{\langle-,-\rangle}$, which is induced by the identity on the Chevalley generators, is proved in \cite[Proposition~8.7]{Hennecart_Geometric}. \end{proof}

In order to unburden the notation, we shall simply write $\calK_{\qv,\C,1}$ for $(\calK_{\qv,\C})\Big\vert_{\upsilon=1}$.

As $i$ is a sink of $\qv$, the classical cotangent stack of $\iii{\stackRep}(\qv)$ is equivalent to the classical truncation of $\iii{\stackRep}(\Pi_\qv)$; similary, as $i$ is a source of $s_i\qv$, the classical cotangent stack of $\stackRep(s_i\qv)^{(i)}$ is equivalent to the classical truncation of $\stackRep(\Pi_{s_i\qv})^{(i)}$. Moreover, the isomorphism $S_i : \iii{\stackRep}(\Pi_\qv) \to \stackRep(\Pi_\qv)^{(i)}$ is equal to the composition $\iii{\stackRep}(\Pi_\qv) \to \stackRep(\Pi_{s_i\qv})^{(i)} \to \stackRep(\Pi_\qv)^{(i)}$ where the first isomorphism is induced by the BGP reflection $S_i$.
The $\mathsf{CC}$ map being compatible with the restriction to open subsets, we have commutative diagrams
\begin{align}\label{eq:Yangian-CC2}
	\begin{tikzcd}[ampersand replacement=\&]
		\calK_{\qv,\C,1} \ar[d]  \ar{r}{\mathsf{CC}'}\&\cohazero \ar[d]\\
		\iii{\calK}_{\qv,\C,1} \ar{r}{\iii{\mathsf{CC}}'} \&   \iii{\coha}_0
	\end{tikzcd}
	\quad \text{and} \quad 
	\begin{tikzcd}[ampersand replacement=\&]
		\calK_{s_i\qv,\C,1} \ar[d]  \ar{r}{\mathsf{CC}'}\&\cohazero \ar[d]\\
		{\calK}^{(i)}_{s_i\qv,\C,1} \ar{r}{(\mathsf{CC}')^{(i)}} \&   {\coha}^{(i)}_0
	\end{tikzcd}\ .
\end{align}
Since the restriction map $\coha_0 \to \iii{\coha}_0$ is surjective by construction, the map $\iii{\mathsf{CC}}'$ is also surjective. As $\iii{\calK}_{\qv,\C,1}$ and $\iii{\coha}_0$ both have the same graded dimension as $\iii{U}(\frakn_\qv)$, the map $\iii{\mathsf{CC}}'$ is an isomorphism. The same applies to $(\mathsf{CC}')^{(i)}$. Finally, using the functoriality of $\mathsf{CC}$ with respect to stack isomorphisms, we obtain the following result.
\begin{corollary}\label{cor:sigma-equivalence-for-CC}
	The following diagram of isomorphisms is commutative:
	\begin{align}\label{eq:sigma-equivalence-for-CC}
		\begin{tikzcd}[ampersand replacement=\&]
			\iii{\calK}_{\qv,\C,1}\ar[swap]{d}{S_{i,\ast}}  \ar{r}{\iii{\mathsf{CC}}'}\&\iii{\coha}_0 \ar{d}{S_{i,\ast}}\\
			{\calK}^{(i)}_{s_i\qv,\C,1} \ar{r}{(\mathsf{CC}')^{(i)}} \& {\coha}^{(i)}_0\ .
		\end{tikzcd}
	\end{align}
	Here we have identified the preprojective algebras of $\qv$ and $s_i\qv$ via the map $u_{\qv,s_i\qv}$.
\end{corollary}

\subsubsection{Proof of Theorem~\ref{thm:compatibility-braid-reflection-functors-degree-zero}}

The claim simply follows from the composition of the specialization at $\upsilon=1$ of the commutative square \eqref{eq:sigma-equivalence-for-C} with the commutative square \eqref{eq:sigma-equivalence-for-CC}. Note that with our conventions, the classical limit of $T'_{i,-1}$ is equal to $T_i$. \hfill $\square$

\subsection{Tautological classes and the braid operators}\label{subsec:tautological-classes-braid}

In this section, we establish a compatibility between the truncated action of the braid group on the Yangian $\eY_\qv$ and the operators of multiplication by tautological classes. Our proof relies on the action of $\eY_\qv$ on the homology of Nakajima quiver varieties. 

\subsubsection{Braid group action on tautological bundles}\label{subsubsec:braid-group-tautological}

Let us begin by describing the action of the braid group $B_\qv$ on the collection of tautological bundles $\calV_j$ on $\Lambda_\qv$, for $j \in I$. We use the same notation for tautological bundles restriction to the open substacks $\iLambda_\qv$ and $\Lambda_\qv^{(i)}$.

Fix $i \in I$. The equivalence $\R S_i$ in \eqref{eq:isomorphism-S-w} induces an isomorphism of $\Ttilde$-equivariant Grothen\-die\-ck groups 
\begin{align}
	S_{i, \ast}\colon 
	\begin{tikzcd}[ampersand replacement=\&]
		\sfK_0^{\Ttilde}\big( \catCoh\big(\iLambda_\qv\big) \big) \ar{r}{\sim} \& \sfK_0^{\Ttilde}\big( \catCoh\big(\Lambdai_\qv\big) \big) \ . 
	\end{tikzcd}
\end{align}
For future use, we record the following fact.  
\begin{lemma}\label{lem:action-tautological-classes} 
	We have
	\begin{align}
		S_{i, \ast}([\calV_k])=\begin{cases}
			[\calV_k]  & \text{for } k \neq i\ , \\[4pt]
			\displaystyle \sum_{\genfrac{}{}{0pt}{}{e\in \doubleOmega}{e\colon i \to j}} \gamma_e [\calV_j] - q^2[\calV_i]  &\text{for } k=i\ .
		\end{cases}
	\end{align}
\end{lemma}
\begin{proof}
	Over objects of $\iLambda_\qv$, we have $\R S_{i, \ast}=S_{i, \ast}$, i.e., $\R S_{i, \ast}$ is in fact underived. From the explicit description of $S_i$ in \S\ref{subsec:equivalent-description-S-i}, we have the following short exact sequence of coherent sheaves on $\Lambda_\qv^{(i)}$:
	\begin{align}
		0\longrightarrow S_i(\calV_i) \longrightarrow \widetilde{\calV}_i \longrightarrow q^2 \calV_i \longrightarrow 0 \ ,
	\end{align}
	where
	\begin{align}
		\widetilde{\calV}_i\coloneqq \bigoplus_{\genfrac{}{}{0pt}{}{e\in \doubleOmega}{e\colon i \to j}} \gamma_e \calV_j\ .
	\end{align}
	The assertion follows.
\end{proof}

Let $\widehat\bS$ be the completion of $\bS$ relatively to the ideal of elements without degree zero term. We define the following generating series in $\widehat\bS$
\begin{align}
	\gamma_e\coloneqq \exp(\varepsilon_e)\ , \quad q^2\coloneqq \exp(\hbar) \ , \quad \ch(z_i) \coloneqq r_i+\sum_{\ell\geq 1} \frac{p_\ell(z_i)}{\ell!}
\end{align}
for $i\in I$ and $e\in \doubleOmega$. There is a unique $B_\qv$-action on the graded algebra $\bS$ such that for any $i\in I$ the element $T_i$ acts through the homogeneous automorphism of degree zero such that
\begin{align}\label{eq:braid-action-tautological-classes}
	T_i (\ch(z_k))\coloneqq 
	\begin{cases} 
		\ch(z_k)  &\text{for } k \neq i\ , \\[4pt]
		\displaystyle \sum_{\genfrac{}{}{0pt}{}{e\in \doubleOmega}{e\colon i \to j}} \gamma_e \ch(z_j) - q^2\ch(z_i)  &\text{for }  k = i\ .
	\end{cases}
\end{align}

Equivalently, if we introduce the following formal series in $\bS\llbracket s \rrbracket$
\begin{align}
	\gamma_e^s\coloneqq \exp(s\varepsilon_e)\ , \quad q^{2s}\coloneqq \exp(s\hbar)\ , \quad\text{and}\quad \ch(z_i,s) \coloneqq r_i+\sum_{\ell\geq 1} p_\ell(z_i)\frac{s^\ell}{\ell!}
\end{align}
then, the $B_\qv$-action on $\bS$ is characterized by
\begin{align}\label{eq:braid-action-tautological-classes-2}
	T_i (\ch(z_k,s))=\begin{cases} 
		\ch(z_k,s) & \text{for } k \neq i\ , \\[4pt] 
		\displaystyle \sum_{\genfrac{}{}{0pt}{}{e\in \doubleOmega}{e\colon i \to j}} \gamma_e^s \ch(z_j,s) - q^{2s}\ch(z_i,s)  &\text{for } k = i \ .
	\end{cases} 
\end{align}

\subsubsection{Yangians' actions on the Borel-Moore homology of Nakajima's quiver varieties}

We next recall some material on Yangians' actions on the Borel-Moore homology of Nakajima's quiver varieties. For any pair of dimension vectors $\bfv, \bfw$ in $\N I$, let $\frakM_\qv(\bfv,\bfw)$ be the Nakajima quiver variety associated to $\qv$ of dimension vectors $\bfv$ and $\bfw$. It carries two families of tautological vector bundles $\calV_i$ and $\calW_i$ over $\frakM_\qv(\bfv, \bfw)$, by varying of $i\in I$, yielding an action of the ring
\begin{align}
	\Hbullet_{\bfv,\bfw}\coloneqq H^\ast_{\sfG_\bfv \times \sfG_\bfw\times \Ttilde}(\mathsf{pt})\ ,
\end{align}
where $\sfG_\bfd\coloneqq \prod_{i\in I} \mathsf{GL}(d_i)$ for $\bfd\in \N I$, by multiplication with the Chern character of both families of bundles on the vector spaces
\begin{align}
	F_{\bfv, \bfw}\coloneqq H_\ast^{\sfG_\bfw\times \Ttilde}(\frakM_\qv(\bfv, \bfw)) \quad \text{and} \quad F_\bfw\coloneqq \bigoplus_{\bfv\in \N I}\, F_{\bfv, \bfw}\ .
\end{align}
By \cite{Varagnolo_Yangian}\footnote{Note that the construction of the representation in \textit{loc.cit.} is for the one-parameter Yangian of $\qv$, but it generalizes straightforwardly to our setting.}, there is a representation of the Yangian $\eY_\qv$ on $F_\bfw$ which is compatible with the $\cohaqv^{\Ttilde, \leqslant 0}$-action on $F_\bfw$ introduced in \cite[\S5.6]{SV_generators}, via the homomorphism $\Phi\colon \Y_\qv^{\frake,\leqslant 0}\to \cohaqv^{\Ttilde, \leqslant 0}$.

The $\zeroeY_\qv$-action on $F_\bfw$ is given as follows, see \cite[\S4]{Varagnolo_Yangian} and \cite[\S5.5]{YZ_Quantum_groups},
\begin{align}
	\begin{tikzcd}[ampersand replacement=\&]
		h_i(z)\ar[mapsto]{r} \& 1+ \displaystyle \sum_{\bfv\in \N I}\frac{\lambda_{-1/z}\Big(\calF_i(\bfv, \bfw) \Big)}{\lambda_{-1/z}\Big(q^2 \calF_i(\bfv, \bfw) \Big)}\cap -\\
		\kappa_{i, \ell} \ar[mapsto]{r} \&\ch_\ell(\calW_i) \cap -
	\end{tikzcd}\ ,
\end{align}
where
\begin{align}
	\calF_i(\bfv,\bfw)\coloneqq q^{-2} \calW_i -(1+q^{-2})\calV_i + \sum_{\genfrac{}{}{0pt}{}{e\in \Omega}{e\colon j\to i}}\gamma_e^{-1} \calV_j + \sum_{\genfrac{}{}{0pt}{}{e\in \Omega}{e\colon i\to j}}\gamma_{e^\ast}^{-1} \calV_j\ .
\end{align}

We denote by $\ev_z \colon \zeroeY_\qv \to \bS_z$ (resp.\ by $ \ev_y \colon \zeroeY_\qv \to \bS_y$) the evaluation at $p_\ell(y_i)=0$ (resp.\ at $z_{i,\ell}=0$) for any $i\in I$ and $\ell\in \N$. 

The elements $p_\ell(y_i)$ and $p_\ell(z_i)$ in $\zeroeY_\qv$ give rise to the formal series $\ch(y_i,s)$ and $\ch(z_i,s)$ of $\zeroeY_\qv\llbracket s \rrbracket$. From this point of view, the action of $\zeroeY_\qv$ on $F_{\bfv,\bfw}$ corresponds to the action induced by the algebra homomorphism (cf.\ \S~\ref{sec:extensionofPhi}):
\begin{align}
	\ev_{\bfv, \bfw}\colon \zeroeY_\qv &\longrightarrow \Hbullet_{\bfv, \bfw} \ ,\\[4pt]
	\ch(z_j,s)  &\longmapsto \ch(\calV_j,s)\ , \\[3pt]
	\ch(y_j,s) &\longmapsto \ch(\calW_j,s)\ .
\end{align}
Let us denote by $\kappa\colon \Y^0_\qv \to \zeroeY_\qv$ the canonical embedding (given by the first assignment in Theorem~\ref{thm:defkappa}). We shall need the following observations.

\begin{lemma}
	\hfill
	\begin{enumerate}\itemsep0.2cm
		\item The compositions
		\begin{align}
			\kappa_z\colon 
			\begin{tikzcd}[ampersand replacement=\&]
				\Y_\qv^0\ar{r}{\kappa} \& \zeroeY_\qv \ar{r}{\ev_z}\& \bS_z 
			\end{tikzcd}
			\quad\text{and}\quad	
			\kappa_y\colon 
			\begin{tikzcd}[ampersand replacement=\&]
				\Y_\qv^0\ar{r}{\kappa} \& \zeroeY_\qv \ar{r}{\ev_y} \& \bS_y 
			\end{tikzcd}
		\end{align}
		are isomorphisms.  
		
		\item The isomorphism $\sigma\coloneqq \kappa_z\circ \kappa_y^{-1}$  is given by
		\begin{align}\label{eq:tautological-z-y}
			\ch(y_j,s) \longmapsto \sum_{\genfrac{}{}{0pt}{}{e\in \doubleOmega}{e\colon j\to i}} \gamma_e^s \ch(z_i,s) - (1+q^{2s})\ch(z_j,s) = T_j(\ch(z_j, s))-\ch(z_j, s)\ . 
		\end{align}
	\end{enumerate}
\end{lemma}
Recall that we have an identification $\zeroeY_\qv \simeq \bS_z\otimes \bS_y$.

\subsubsection{Compatibility between truncated braid group action and tautological classes}

We are now ready to state the main result of this section. 

\begin{theorem}\label{thm:braid-Chern}
	Let $i\in I$, $u\in\iY_\qv$, and $c\in\bS$. We have $\overline{T}_i(c\cap u) = T_i(c) \cap \overline{T}_i(u)$.
\end{theorem}

\begin{proof}
	We shall prove that
	\begin{align}\label{eq:commutativity-tautological}
		\overline{T}_i(p_\ell(z_j)\cap u)=T_i(p_\ell(z_j)) \cap \overline{T}_i(u)
	\end{align}
	for any $u\in \iY$, $j\in I$, and $\ell\in \N$.
	
	First, note that if $j\neq i$ we have $T_i(p_\ell(z_j)\cap u) = p_\ell(z_j)\cap T_i(u)=T_i(p_\ell(z_j))\cap T_i(u)$ for any $u \in \eY_\qv$ since $p_\ell(z_j)\cap x_i^\pm=0$. In particular, Formula~\eqref{eq:commutativity-tautological} holds.
	
	\medskip
	
	From now on, we consider the case $j=i$. First,	under the isomorphism $\zeroeY_\qv\simeq \bS_z\otimes\bS_y$, we have $p_\ell(z_i)\cap u=[p_\ell(z_i),u]$ for any $u\in \eY_\qv$, hence
	\begin{align}
		T_i(p_\ell(z_i)\cap u)= [T_i(p_\ell(z_i)),T_i(u)]\ .
	\end{align}
	
	Let 
	\begin{align}
		\pr\colon \eY_\qv \longrightarrow \Y_\qv^-\quad\text{and}\quad \iii{\pi}\colon \Y_\qv^-\longrightarrow \iY_\qv
	\end{align}
	be the canonical projections. Set $t_{i,\ell}\coloneqq T_i(p_\ell(z_i))-p_\ell(z_i)$, which is an element of $\eY_i$, where we see $p_\ell(z_i)$ as an element in $\Y_i^{\frake,0}$. If $u\in\Y_\qv^-$, then
	\begin{align}\label{eq:proof-braid-1}
		\begin{split}
			\overline{T}_i(p_\ell(z_i)\cap \iii{\pi}(u))&=(\iii{\pi}\circ \pr)( [T_i(p_\ell(z_i)),T_i(u)])\\
			&= (\iii{\pi}\circ \pr) ( [p_\ell(z_i),T_i(u)])+ (\iii{\pi}\circ \pr) (\ad(t_{i,\ell})T_i(u))\\
			&=p_\ell(z_i)\cap \overline{T}_i(\iii{\pi}(u)) + (\iii{\pi}\circ \pr)(\ad(t_{i,\ell})T_i(u))\ .
		\end{split}
	\end{align}
	
	Note that $t_{i, \ell}\in \eY_i[0]$ and that
	\begin{align}
		\eY_i[0]=\Y_i^{\frake,0} \oplus (\I_i^-\zeroeY_i\I_i^+) \cap \eY_i[0]\ .
	\end{align}
	By construction, we have $(\iii{\pi}\circ \pr)(\I_i^-\zeroeY_i\I_i^+ v)=0 = (\iii{\pi} \circ \pr) (v \I_i^-\zeroeY_i\I_i^+ )$ for any $v\in \eY_\qv$. Hence 
	\begin{align}
		\iii{\pi} \circ\pr \circ \ad(t_{i,\ell})=\iii{\pi} \circ\pr \circ \ad(\rho_i(t_{i,\ell}))\ ,	
	\end{align}
	where $\rho_i\colon \eY_i[0]\to \zeroeY_i$ is the projection along $\I_i^-\zeroeY_i\I_i^+$. Moreover, by construction $\rho_i(t_{i,\ell})\in \Y_i^0$. 
	
	\medskip 
	
	Next, we compute the element $\rho_i(t_{i,\ell})=\rho_i((T_i-1)(p_{ \ell}(z_i)))$ in $\Y_i^0$. We first consider the  $A_1$ quiver consisting of the vertex $\{i\}$ instead of $\qv$. In order to distinguish them from $\Y^0,\kappa, \ldots$, we let $\Y^0_i, \kappa_i,\ldots$ refer to the same objects, but for the $A_1$ quiver. For any pair of non-negative integers $(v_i, w_i)$, we have $\frakM(v_i,w_i)=T^\ast \mathsf{Gr}(v_i,w_i)$. The vector space
	\begin{align}
		F_{w_i}= \bigoplus_{v_i=0}^{w_i} F_{v_i,w_i}
	\end{align}
	is an $\ell$-highest weight representation of $\Y_i^{\frake}$ with $\ell$-highest weight vector $u_+ \coloneqq [\frakM(0,w_i)]$ whose $\ell$-highest weight is given by the map $\ev_{0, w_i}\circ \kappa_i\colon \Y^0_i\to \Hbullet_{0, w_i}$ (cf.\ \cite{Varagnolo_Yangian}). The algebra $\eY_i$ acts on $F_{w_i}$ and the operator $T_i$ acts on $\eY_i$ and $F_{w_i}$ in a compatible way via Formula~\eqref{eq:Ti}. Set $u_-\coloneqq T_i^{-1}(u_+)$. Then $u_-\in F_{w_i,w_i}$ is a lowest weight vector.
	
	Consider the formal series $\ch(z_i,s)$ in $\Y^0_i\llbracket s\rrbracket$. Note that $\calV_i = \calW_i$ on $\frakM(w_i,w_i)$. Hence, we have
	\begin{align}\label{eq:braid-hw-vect}
		\begin{split}
			T_i(\ch(z_i,s))\cap u_+ &= T_i(\ch(z_i,s) \cap u_-) =T_i (\ch(\calV_i,s) \cap u_-)\\
			&= T_i (\ch(\calW_i,s) \cap u_-)=\ch(\calW_i,s)\cap u_+
		\end{split}\ .
	\end{align}
	Now consider the formal series $\displaystyle t_i(s)\coloneqq \sum_{\ell\in \N} t_{i,\ell}\,s^\ell$ in $\eY_i\llbracket s\rrbracket$. Noticing that $p_{\ell}(z_i)\cap u_+=0$ for any $\ell \geq 1$ and that the evaluation morphism 
	\begin{align}
		\eY_i[0]  \longrightarrow \Hbullet_{0, w_i} \ , \quad
		x  \longmapsto x\cap u_+ 
	\end{align}
	factors through $\rho_i$, we deduce from the equality \eqref{eq:braid-hw-vect} that 
	\begin{align}
		(\ev_{0,w_i}^i\circ \kappa_i)(\rho_i(t_i(s)))=\ch(\calW_i, s)\ ,
	\end{align}
	where $\ev_{v_i,w_i}^i\colon \Y_i^{\frake,0} \to \Hbullet_{v_i, w_i}$. Since this holds for any $w_i$, we have 
	\begin{align}
		\rho_i(t_i(s))= \kappa_{y_i}^{-1}((\ch(y_i,s))\in\Y_i^0 \llbracket s  \rrbracket \ ,
	\end{align}
	which completely determines $\rho_i(t_i(s))$.
	
	Going back to the setting of the whole quiver $\qv$, the previous analysis shows that
	\begin{align}
		\rho_i(t_i(s))= \kappa_{y_i}^{-1}((\ch(y_i,s))\in\bS \llbracket s  \rrbracket \ .
	\end{align}
	Since $\bS_y$ is central in $\eY_\qv$, and $\kappa_y^{-1}(\ch(y_i,s)) \in \mathsf{Span}\{ \ch(y_j,s), \ch(z_j,s)\,\vert\; j \in I\}$ we have
	\begin{align}
		\iii{\pi}\circ\pr \circ \ad(y)=	\iii{\pi}\circ\pr \circ \ad( \kappa_z(y))
	\end{align}
	for all $y\in \bS$. Thus, we obtain
	\begin{align}
		(\iii{\pi}\circ\pr)(\ad(\rho_i(t_i(s))))=(\iii{\pi}\circ\pr)(\ad(\kappa_z(\rho_i(t_i(s)))))=\iii{\pi}\circ\pr\circ \ad(\sigma(\ch(y_i,s)))\ .
	\end{align}
	
	Now, we go back to Formula~\eqref{eq:proof-braid-1}: for $u\in \Y_\qv^-$ we have
	\begin{align}
		\overline{T}_i(p_\ell(z_i)\cap \iii{\pi}(u))&=p_\ell(z_i)\cap \overline{T}_i(\iii{\pi}(u)) +\iii{\pi} \circ\pr \circ \ad(\rho_i(t_{i,\ell}))\\
		&=p_\ell(z_i)\cap \overline{T}_i(\iii{\pi}(u)) +(\iii{\pi} \circ\pr)([T_i(p_{i,\ell})- p_{i, \ell},T_i(u)]) \ ,
	\end{align}
	where the second equality follows from Formula~\ref{eq:tautological-z-y}. Therefore, the claim \eqref{eq:commutativity-tautological} follows.
\end{proof}

\subsection{Support and shuffle algebras}\label{subsec:support-shuffle-algebras} 

As explained in \S\ref{subsubsec:shuffle-algebra}, $\cohaqv^{\Ttilde}$ admits a canonical embedding into a shuffle algebra $\shuffle_\qv$. In particular, for any $\bfd \in \N I$ we have an inclusion $\cohaqvd^{\Ttilde} \subset \shuffle_\bfd$, and $\cohaqvd^{\Ttilde}$ is an ideal for the standard algebra structure on $\shuffle_\bfd$. The same holds for $\coha_{\bfd-\alpha_i}^{\Ttilde} \star \coha_{\alpha_i}^{\Ttilde}$ and $\coha_{\alpha_i}^{\Ttilde} \star \coha_{\bfd-\alpha_i}^{\Ttilde}$, and hence also for $\icohatilded$ and $\cohatildeid$. The aim of this section is to give an estimate on the supports of the latter $\shuffle_\bfd$-modules, viewed as coherent sheaves on $\Spec(\shuffle_\bfd)$, see Propositions~\ref{prop:support} and \ref{prop:support-} below.

\subsubsection{Partitions and strings}\label{subsubsec:partitions-strings} 

We begin by introducing some combinatorially defined subschemes of $\Spec(\shuffle_\bfd)$. Recall that for any $\bfd$, we have
\begin{align}
	\shuffle_\bfd\coloneqq \Hbullet_{\Ttilde}\Big[z_{i,\ell}\,\Big\vert\, i \in I\ , \ 1 \leq \ell \leq d_i\Big]^{\frakS_\bfd}\ .
\end{align}
Note that we have taken the torus to be $A_{\mathsf{max}}$. We also set
\begin{align}
	\shuffle_\qv\coloneqq \bigoplus_{\bfd\in \N I} \shuffle_\bfd\ .
\end{align}
Each $\shuffle_\bfd$ is a polynomial algebra for $\bfd\in \N I$, hence we can define
\begin{align}
	\bfC^{(\bfd)}\coloneqq \Spec(\shuffle_\bfd)\ .
\end{align}
We slightly generalize this definition. For a finite set $D$, define 
\begin{align}
	\shuffle_D\coloneqq \Hbullet_{\Ttilde}\big[y_a \,\big\vert\, a \in D\big]^{\frakS_D}\ ,
\end{align}
where $\frakS_{D}$ is the group of permutations of $D$. Set
\begin{align}
	\bfC^{(D)}\coloneqq \Spec(\shuffle_D)\ .
\end{align}
By abuse of notation, for any $n\in \N$, we set $\bfC^{(n)}\coloneqq \bfC^{(I_n)}$, where $I_n\coloneqq \{1, \ldots, n\}$. Furthermore, we write $\bfC\coloneqq \Spec(\Hbullet_{\Ttilde})$. We will view points of $\bfC^{(\bfd)}$ as parametrizing collections of $I$-colored points of $\bfC$ (possibly with multiplicities). We shall also sometimes use the notation
\begin{align}
	\bfC^{(D_1, D_2, \ldots, D_s)}\coloneqq \bfC^{(D_1)} \times \cdots \times \bfC^{(D_s)}\ .
\end{align}

\begin{definition}
	A point $\bfy\in\bfC^{(D)}$ is \textit{generic} if $y_a-y_b \not\in \Z \hbar$ for any $a,b\in D$, with $a \neq b$. Let $\bfC^{(D),\circ}\subset \bfC^{(D)}$ be the open subset of generic points. 
\end{definition}

Define
\begin{align}
	\mathsf{Comp}^D_n\coloneqq \Big\{\phi\colon D \to \N\,\Big\vert\, \sum_a \phi(a)=n\Big\}\quad \text{and} \quad \mathsf{Par}^{D}_n\coloneqq \mathsf{Comp}_n^D/\frakS_D 
\end{align}
for any $n\in \N$, which we shall call \textit{$D$-compositions of $n$} and \textit{$D$-partitions of $n$}, respectively. We will denote a $D$-partition by $\lambda=(\lambda_a)_{a\in D}$, remembering that these are defined up to permutation. 

\begin{definition}
	For $y \in \bfC$ and a positive integer $n$, the \textit{string of length $n$ starting at $y$} is the element
	\begin{align}
		\mathsf{str}(y,n)\coloneqq (y, y+\hbar, \ldots, y +n\hbar-\hbar) \in  \bfC^{(n)}\ .
	\end{align}
	A string of length zero is empty. 
\end{definition}

Given a partition $\lambda=(\lambda_1, \ldots, \lambda_s)$ of $n$, we have the tuple
\begin{align}
	\text{str}(y,\lambda)=(\text{str}(y,\lambda_1), \text{str}(y,\lambda_2), \ldots, \text{str}(y,\lambda_s)) \in 
	\bfC^{(\lambda)}\coloneqq \bfC^{(\lambda_1)} \times \cdots \times \bfC^{(\lambda_s)}\ .
\end{align}
We may view it as an element in $\bfC^{(n)}$ via the obvious map $\bfC^{(\lambda)}\to \bfC^{(n)}$. 
\begin{definition}
	Let $\bfy\in\bfC^{(D)}, \bfz \in \bfC^{(n)}$ with $\bfy$ generic, and let $\lambda$ be a $D$-partition of $n$. The pair $(\bfy,\bfz)$ is \textit{in left relative position $\lambda$} if there is a permutation $\sigma \in \frakS_D$ such that
	\begin{align}
		\bfz=\bigcup_{a\in D}\mathsf{str}(y_{\sigma(a)},\lambda_a)\ . \tag*{\qedhere} 
	\end{align}
\end{definition}
\begin{remark}
	Pictorially, $(\bfy,\bfz)$ is \textit{in left relative position $\lambda$} if the subset of $\bfC$ defining $\bfz$ is obtained by attaching strings of length $\lambda$ at the points of $\bfy$, in some order. 
\end{remark}

We let $X^{\circ}_{D,n,\lambda} \subset \bfC^{(D),\circ} \times \bfC^{(n)}$ be the (reduced) closed subvariety of pairs $(\bfy,\bfz)$ which are in left relative position $\lambda$. The first projection $X^{\circ}_{D,n,\lambda} \to \bfC^{(D),\circ}$ is a finite morphism. We consider the following subvarieties of $ \bfC^{(D)} \times \bfC^{(n)}$:
\begin{align}\label{eq:subvarieties}
	X_{D,n,\lambda}\coloneqq\overline{X^\circ_{D,n,\lambda}} \ ,\quad Z^\circ_{D,n}\coloneqq \bigcup_\lambda X^\circ_{D,n,\lambda}\ ,\quad Z_{D,n}\coloneqq \bigcup_\lambda X_{D,n,\lambda}\ .
\end{align}
Note that $X_{D,n,\lambda}$ and $Z_{D,n}$ are closed, reduced, of dimension $\vert D\vert\dim(\Ttilde)$, and that $X_{D,n,\lambda}$ is irreducible. 

Now, we introduce a finite set $D$ which will be associated to fixed $\bfd\in \N I$ and $i\in I$. First, set $\bfh\coloneqq \bfd-d_i\alpha_i$. For $j\in I$, with $j\neq i$, define
\begin{align}\label{eq:D}
	\Omega_{j,i}\coloneqq \{e \in \doubleOmega\,\vert\, e\colon j \to i\} \ ,\quad D_{\bfd,j}\coloneqq \{1, 2, \ldots, d_j\}\times \Omega_{j,i}\ ,\quad
	D\coloneqq \bigsqcup_{\genfrac{}{}{0pt}{}{j\in I}{j\neq i}} D_{\bfd,j}\ .
\end{align} 
The group $\frakS_{d_j}$ acts on $D_{\bfd,j}$. Thus, the group 
\begin{align}\label{eq:frakS}
	\frakS\coloneqq\prod_{\genfrac{}{}{0pt}{}{j\in I}{j\neq i}}\frakS_{d_j}
\end{align}
acts on the set $D$. An element in $D$ is a triple $(j,k,e)$ with $j\in I$, $k\in \{1,\ldots, d_j\}$, and $e\in\Omega_{j,i}$ (hence $i \neq j$ as $\qv$ is assumed to have no edge loops). 

Consider the surjective map
\begin{align}\label{eq:theta-map}
	\theta\colon \shuffle_D\longrightarrow \shuffle_\bfh
\end{align}
sending $y_a\in \shuffle_D$ with $a=(j, k, e)\in D$ to $z_{j, k}- \varepsilon_e\in \shuffle_\bfh$. Let $\bfC^{(\bfh)} \to \bfC^{(D)}$ be the corresponding injective map. Set
\begin{align}
	Z_{\bfd, i}\coloneqq \bfC^{(\bfh)} \times_{\bfC^{(D)}} Z_{D, d_i} \subset \bfC^{(\bfh)} \times \bfC^{(d_i)}=\bfC^{(\bfd)}\ .
\end{align}
Thus, $Z_{\bfd,i}$ is a (in general not reduced) closed subvariety of $\bfC^{(\bfd)}$ of codimension $d_i$. The $\frakS$-action on $D$ yields a $\frakS$-action on the set of $D$-compositions. For any $\phi\in\mathsf{Comp}^D_{d_i}$, set $O\coloneqq \frakS\cdot \phi$, and 
\begin{align}\label{eq:calC-calO}
	\calC_{O}\coloneqq \bigg\{ \Big(\bfz, \bigcup_{j,k,e}\mathsf{str}(z_{j,k}-\varepsilon_e,\phi(j,k,e))\Big)\,\bigg\vert\, \bfz \in \bfC^{(\bfh)}\bigg\} \ .
\end{align}
\begin{remark}\label{rem:C-O}
	Note that $\calC_{O}$ depends only on $O$ and not on $\phi$.
	
	Moreover, each $\calC_{O}$ is irreducible, closed, and of codimension $d_i$ in $\bfC^{(\bfd)}$. Hence it is an irreducible component of $Z_{\bfd,i}$. In addition, as $O$ runs among the $\frakS$-orbits in $\mathsf{Comp}^D_{d_i}$, the subvarieties $\calC_{O}$ cover $Z_{\bfd,i}$ and are pairwise distinct (though not disjoint).
\end{remark}

\subsubsection{Supports}

In this section, we shall prove the following estimate for the support of $\cohatildeid$ in $\bfC^{(\bfd)}$.
\begin{proposition}\label{prop:support}
	As a coherent sheaf on $\bfC^{(\bfd)}$, $\cohatildeid$ is pure of dimension $\sum_{j \neq i} d_j$ and $\mathsf{Supp}(\cohatildeid) \subseteq Z_{\bfd,i}$.  \end{proposition}
\begin{proof} 
	The subspace $\shuffle_{\bfd-\alpha_i}\star \shuffle_{\alpha_i}$ is an ideal of $\shuffle_\bfd$. As we shall show in Lemma~\ref{lem:equality-coha}, which we shall prove in \S\ref{subsubsec:injectivity}, the canonical map
	\begin{align}
		\cohatildeid= \coha_\bfd^{\Ttilde}/(\coha_{\bfd-\alpha_i}^{\Ttilde}\star \coha_{\alpha_i}^{\Ttilde}) \longrightarrow \shuffle_\bfd / (\shuffle_{\bfd-\alpha_i} \star \shuffle_{\alpha_i})\coloneqq \calE
	\end{align}	
	is injective. To finish the proof, it thus only remains to check that $\calE$ is isomorphic to $\scrO_{Z_{D,d_i}}$. We begin by introducing an auxiliary coherent sheaf $\calE_D$ on $\bfC^{(D)} \times \bfC^{(d_i)}$, from which $\calE$ will be obtained by base change. For this, consider the unique $\shuffle_D$-linear map
	\begin{align}\label{eq:omega}
		\omega\colon \shuffle_D \otimes \shuffle_{(d_i-1)\alpha_i} \otimes \shuffle_{\alpha_i} \longrightarrow \shuffle_D \otimes \shuffle_{d_i\alpha_i}
	\end{align}
	that fits into the commutative diagram
	\begin{align}\label{eq:commutative-diagram-shuffle}
		\begin{tikzcd}[ampersand replacement=\&]
			\shuffle_D \otimes \shuffle_{(d_i-1)\alpha_i} \otimes \shuffle_{\alpha_i} \ar[twoheadrightarrow]{d}{\Psi} \ar{r}{\omega} \& \shuffle_D \otimes \shuffle_{d_i\alpha_i} \ar[twoheadrightarrow]{d}{\Psi'}\\
			\shuffle_{\bfd-\alpha_i}\otimes \shuffle_{\alpha_i}\ar{r}{\ast} \& \shuffle_\bfd
		\end{tikzcd}\ ,
	\end{align}
	where 
	\begin{align}
		\Psi&\colon 	
		\begin{tikzcd}[ampersand replacement=\&, column sep=huge]
			\shuffle_D \otimes \shuffle_{(d_i-1)\alpha_i} \otimes \shuffle_{\alpha_i} \ar{r}{\theta\otimes \id \otimes \id} \& \shuffle_{\bfh}\otimes  \shuffle_{(d_i-1)\alpha_i} \otimes \shuffle_{\alpha_i} \simeq  \shuffle_{\bfd-\alpha_i}\otimes \shuffle_{\alpha_i}
		\end{tikzcd}\ ,\\[2pt]
		\Psi'&\colon 	
		\begin{tikzcd}[ampersand replacement=\&, column sep=huge]
			\shuffle_D\otimes \shuffle_{d_i\alpha_i} \ar{r}{\theta\otimes \id} \& \shuffle_{\bfh}\otimes \shuffle_{d_i\alpha_i} \simeq  \shuffle_\bfd
		\end{tikzcd}\ .
	\end{align}	
	Here, the map $\theta$ is introduced in Formula~\eqref{eq:theta-map}. Explicitly, we have
	\begin{multline}\label{eq:deftheta}
		\omega(1\otimes g\otimes h)(\bfy,z_1, \ldots, z_{d_i}) =\\ (-1)^{\langle \bfd-\alpha_i, \alpha_i\rangle} \sum_{\ell=1}^{d_i}h(z_\ell)g(z_1, \ldots, \widehat{z_\ell}, \ldots, z_{d_i})
		\prod_{a\in D} (z_\ell-y_a) \prod_{k \neq \ell} \frac{z_k-z_\ell-\hbar}{z_k-z_\ell}\ ,
	\end{multline}
	where $g\in\shuffle_{(d_i-1)\alpha_i}$, and $h\in\shuffle_{\alpha_i}$. Here, the notation $g(z_1, \ldots, \widehat{z_\ell}, \ldots, z_{d_i})$ means that the variable $\widehat{z_\ell}$ is removed from the tuple $(z_1, \ldots, z_{d_i})$ for any $\ell=1, \ldots, d_i$. 
	
	Note that the right-hand-side of the above formula is regular, although the summands have poles. Furthermore, the map \eqref{eq:omega} is $\shuffle_D \otimes \shuffle_{d_i\alpha_i}$-linear, where $\shuffle_{d_i\alpha_i}$ acts on the left-hand-side via the obvious embedding into $\shuffle_{(d_i-1)\alpha_i} \otimes \shuffle_{\alpha_i}$. We let $\calE_D$ be the cokernel of $\omega$. Note that 
	\begin{align}
		\scrO_{Z_{\bfd,i}} = (\theta^\ast  \otimes \id)(\scrO_{Z,D,d_i}) \in \catCoh(\bfC^{(\bfd)})\ ,
	\end{align}
	while, by virtue of the diagram~\eqref{eq:commutative-diagram-shuffle}, 
	\begin{align}
		\calE=(\theta^\ast  \otimes \id)(\calE_D)\ .
	\end{align}
	We are thus reduced to showing that $\calE_D \simeq \scrO_{Z,D,d_i}$. 
	
	From Formula~\eqref{eq:deftheta}, it is easy to see that for any $f \in \shuffle_D, g \in \shuffle_{(d_i-1)\alpha_i}$ and $h \in \shuffle_{\alpha_i}$ we have $\omega(f \otimes g \otimes h)(\bfy,\bfz)=0$ for any $(\bfy,\bfz) \in Z^{\heartsuit}_{D,d_i}$, where
	\begin{align}
		Z^{\heartsuit}_{D,d_i}\coloneqq \Big\{(\bfy,\bfz)\in Z_{D,d_i}^\circ\,\Big\vert\,z_\ell \neq z_k\ ,\, \forall \ell \neq k\Big\}\ .
	\end{align}
	Indeed, every term in the defining sum of $\omega(f\otimes  g \otimes h)$ is regular and vanishes. As $Z^{\heartsuit}_{D,d_i}$ is open and dense in $Z_{D,d_i}$, we deduce that there is a canonical surjective morphism 
	\begin{align}
		\kappa\colon \calE_D \to \scrO_{Z_{D,i}}\ .
	\end{align}
	In order to prove that $\kappa$ is injective, we use the (affine) projection $\pi\colon \bfC^{(D)} \times \bfC^{(d_i)} \to \bfC^{(D)}$.
	
	\begin{lemma}\label{lem:support}
		Set $N\coloneqq \vert D\vert$ and $c_{N,d_i}\coloneqq \big(\begin{smallmatrix}N+d_i-1\\ d_i\end{smallmatrix}\big)$. Then
		\begin{enumerate}\itemsep0.2cm
			\item \label{item:b-support} $\pi_\ast (\scrO_{Z_{D,d_i}})$ is of (generic) rank $c_{N,d_i}$.
			\item \label{item:a-support} $\pi_\ast (\calE_D)$ is a locally free coherent sheaf on $\bfC^{(D)}$ of rank $c_{N,d_i}$.
		\end{enumerate}
	\end{lemma}
	\begin{proof} 
		Let $(\bfy,\hbar)$ be a generic point of $\bfC^{(D)}$. The description of $Z_{D,d_i}$ in Formula~\eqref{eq:subvarieties} shows that the dimension of the fiber of $\scrO_{Z_{D,d_i}}$ at $(\bfy,\hbar)$ is equal to the number of $D$-compositions of $d_i$, which is $c_{N,d_i}$. This shows \eqref{item:b-support}. 
		
		Now, we prove \eqref{item:a-support}. Let us endow the algebras $\shuffle_D$ and $\shuffle_{d_i\alpha_i}$ with a grading by putting each of the variables $y_a, z_k$, and $\hbar$ in degree one. The map $\omega$ in Formula~\eqref{eq:omega} is homogeneous. We deduce that the $\shuffle_D\otimes \shuffle_{d_i\alpha_i}$-module $\calE_D$ is graded. 
		
		The algebra $\shuffle_D$ is a noetherian integral domain and $\calE_D$ is a finitely presented module over $\shuffle_D$. Thus, it is enough to show that the dimension of $\calE_D\otimes_{\shuffle_D} k(\frakp)$ is the same for all prime ideals $\frakp$ of $\shuffle_D$. Here $k(\frakp)$ denotes the residue field at $\frakp$. Thanks to the grading, it is enough to check this at a generic point and at the point $(\bfy, \hbar)=(0,0)$ in $\bfC^{(D)}$.
		
		Now, the fiber $\calE_D \otimes_{\shuffle_D} k(0,0)$ of $\calE_D$ over $(0,0)$ is isomorphic the quotient of $\C[z_1, \ldots, z_{d_i}]^{\frakS_{d_i}}$ by the ideal $I_N$ generated by the set
		\begin{align}
			\Big\{\sum_{\ell=1}^{d_i}(z_\ell)^N h(z_\ell)\,\Big\vert\, h\in\C[z]\Big\} \ .
		\end{align}
		The ideal $I_N$ is generated by the power-sum functions $p_j$ with $j\geq N$. We have
		\begin{align}
			\C[z_1, \ldots, z_{d_i}]^{\frakS_{d_i}}=I_N \oplus \bigoplus_{\lambda} \C m_\lambda\ ,
		\end{align}
		where the direct sum runs over the set of partitions of length $\leqslant d_i$ whose parts are strictly less than $N$, and $m_\lambda$ is the 
		monomial symmetric function. The number of such partitions is $c_{N,d_i}$. Thus, we have 
		\begin{align}
			\dim(\calE_D\otimes_{\shuffle_D} k(0,0)) =  c_{N,d_i} \ .
		\end{align}
		We now consider $\calE_D\otimes_{\shuffle_D}k(\bfy, \hbar)$ at a generic point $(\bfy, \hbar)$ in $\bfC^{(D),\circ}$. By upper semicontinuity, we have
		\begin{align}
			c_{N,d_i}=\dim(\calE_D\otimes_{\shuffle_D}k(0,0))\geq \dim(\calE_D\otimes_{\shuffle_D}k(\bfy, \hbar))\geq \dim(\C[Z_{D,d_i}]\otimes_{\shuffle_D}k(\bfy,\hbar)) = c_{N,d_i}\ .
		\end{align}
		It follows that $\dim(\calE_D\otimes_{\shuffle_D}k(\bfy, \hbar))=c_{N,d_i}$. This proves that $\calE_D$ is a locally free sheaf of rank $c_{N,d_i}$.
	\end{proof}
	
	We may now finish the proof of Proposition~\ref{prop:support}. The morphism $\pi_\ast (\kappa)$ is onto and generically an isomorphism since $\pi_\ast (\calE_D)$ and $\pi_\ast (\scrO_{Z_{D,i}})$ have the same generic rank. As $\pi_\ast (\calE_D)$ is locally free, in particular torsion-free, this implies that $\pi_\ast (\kappa)$ is injective. But then $\kappa$ is itself injective, and thus $\kappa$ is an isomorphism as desired. 
\end{proof}

A version of Proposition~\ref{prop:support} can be formulated for $\icohatilded$. We set
\begin{align}
	\mathsf{str}(n, y)\coloneqq (y, y-\hbar, \ldots, y -n\hbar+ \hbar) \in \bfC^{(n)}
\end{align}
for $y\in \bfC$ and $n\geq 1$. Moreover, if $\lambda$ is a partition of $n$, we set
\begin{align}
	\mathsf{str}(\lambda, y)\coloneqq (\mathsf{str}(\lambda_1, y), \ldots, \mathsf{str}(\lambda_s, y))\in \bfC^{(n)}\ .
\end{align}

\begin{definition}
	Let $\bfy\in\bfC^{(D)}$ be generic, $\lambda$ a $D$-partition of $n$, and $\bfz \in \bfC^{(n)}$. A pair $(\bfy,\bfz)$ is \textit{in right relative position $\lambda$} if there is a permutation $\sigma \in \frakS_D$ such that
	\begin{align}
		\bfz=\bigcup_{a\in D}\text{str}(\lambda_a,y_{\sigma(a)})\ . \tag*{\qedhere} 
	\end{align}	
\end{definition}

Let $X^\circ_{\lambda,D,n} \subset \bfC^{(D),\circ} \times \bfC^{(n)}$ be the (reduced) closed subvariety of all pairs which are in right relative position $\lambda$. We set
\begin{align}
	X_{\lambda,D,n}\coloneqq \overline{X^\circ_{\lambda,D,n}} \subset \bfC^{(D)} \times \bfC^{(n)}\ ,\quad W^\circ_{D,n}\coloneqq \bigcup_\lambda X^\circ_{\lambda,D,n}\ ,\quad W_{D,n}\coloneqq \bigcup_\lambda X_{\lambda,D,n}\ .
\end{align}

Fix $i\in I$ and $\bfd\in \Z I$. For $j\in I$, with $j\neq i$ and, we define 
\begin{align}
	D_{j,a}\coloneqq [1, a_j]\times \Omega_{j,i}\quad\text{and}\quad D=\bigsqcup_j D_{j,a}\ .
\end{align}
The group $\frakS$, introduced in \eqref{eq:frakS}, acts on the set $D$. Let $\bfh\coloneqq \bfd-d_i\alpha_i$. Consider the surjective map
\begin{align}\label{eq:theta-map-2}
	\theta\colon \shuffle_D\longrightarrow \shuffle_\bfh
\end{align}
sending $y_a\in \shuffle_D$ with $a=(j, k, e)\in D$ to $z_{j, k}+ \varepsilon_e\in \shuffle_\bfh$. Let $\bfC^{(\bfh)} \to \bfC^{(D)}$ be the corresponding injective map. Set
\begin{align}
	Z_{i, \bfd}\coloneqq \bfC^{(\bfh)} \times_{\bfC^{(D)}} W_{D, d_i} \subset \bfC^{(\bfd)}\ .
\end{align}
The irreducible components of $Z_{i, \bfd}$ are described as follows. For any $\frakS$-orbit $\calO\coloneqq\frakS \cdot \phi$ in $\mathsf{Comp}^D_n$ we put
\begin{align}\label{eq:CO}
	\calC_{\calO}=\bigg\{ \Big(\bfz, \bigcup_{j,k,e}\text{str}(\phi(j,k,e),z_{j,k}+\varepsilon_e)\Big)\,\bigg\vert\,\bfz \in \bfC^{(\bfh)}\bigg\}\ .
\end{align}
Each $\calC_{\calO}$ is irreducible, closed, and of codimension $d_i$ in $\bfC^{(\bfd)}$ as $\calO$ runs among the $\frakS$-orbits in $\mathsf{Comp}^D_{d_i}$. Hence it is an irreducible component of $Z_{i, \bfd}$. In addition, the sets $\calC_{\calO}$ cover $Z_{i, \bfd}$ and are pairwise distinct (though not disjoint).

The following is the analog of Proposition~\ref{prop:support} for $\icohatilded$.
\begin{proposition}\label{prop:support-}
	As a coherent sheaf on $\bfC^{(\bfd)}$, $\icohatilded$ is pure of dimension $\sum_{j \neq i} d_j$ and $\mathsf{Supp}(\icohatilded) \subseteq Z_{i,\bfd}$.  .
\end{proposition}

\subsubsection{Injectivity}\label{subsubsec:injectivity}

This section is devoted to the proof of the following result, which was needed in the derivation of Proposition~\ref{prop:support}.
\begin{lemma}\label{lem:equality-coha} 
	For $i\in I$ and $\bfd\in \N I$, we have 
	\begin{align}\label{eq:equality-coha}
		\cohaqvd^{\Ttilde} \cap (\shuffle_{\bfd-\alpha_i} \star \shuffle_{\alpha_i}) = \coha_{\bfd-\alpha_i}^{\Ttilde}\star \coha_{\alpha_i}^{\Ttilde}\ .
	\end{align}
\end{lemma}

Taking the direct sum with respect to all dimension vectors of \eqref{eq:equality-coha}, we obtain
\begin{align}
	\cohaqv^{\Ttilde} \cap (\shuffle_\qv \star \shuffle_{\alpha_i})= \cohaqv^{\Ttilde} \star \shuffle_{\alpha_i}\ .
\end{align}

To prove the lemma, we shall use the localized coproduct on the shuffle algebra $\shuffle_\qv$, constructed by Yang-Zhao in \cite[\S2]{YZ_Quantum_groups}, to which we refer for details. We consider the extended shuffle algebra $\shuffle_\qv^{\frake}\coloneqq \shuffle_\qv^0\ltimes \shuffle_\qv$, where
\begin{align}
	\shuffle_\qv^0\coloneqq \bigotimes_{j\in I} \Hbullet_{\Ttilde}\Big[h_{j,\ell}\,\vert\,\ell>0\Big]^{\frakS_\infty}\ .
\end{align}
The algebra $\shuffle_\qv^{\frake}$ has a natural grading $\displaystyle 	\shuffle_\qv^{\frake}=\bigoplus_{\bfd\in \N I}\shuffle_\bfd^{\frake}$, where
\begin{align}
	\shuffle_\bfd^{\frake}\coloneqq \shuffle_\qv^0\ltimes \shuffle_\bfd=\shuffle_\qv^0\Big[z_{j, k}\,\Big\vert\, j\in I\ , \ 0 <k\leq d_i\Big]^{\frakS_\bfd} \ .
\end{align}
It is equipped with a localized coproduct
\begin{align}\label{eq:coprod}
	\Delta \colon \shuffle_\bfc^{\frake}  \longrightarrow \bigoplus_{\bfd+\bfe=\bfc} \Big(\shuffle_\bfd^{\frake}\otimes \shuffle_\bfe^{\frake}\Big)_{\mathsf{loc}}\ ,
\end{align}
which takes values in a space of Laurent series in the variables $z_{j,k}, h_{j,k}$. The precise form is not important to us. It will suffice to know that it enjoys the following properties:
\begin{enumerate}\itemsep0.2cm
	\item \label{item:a} The algebra structure on $\shuffle_\qv^{\frake}$ extends to the right-hand-side of \eqref{eq:coprod} and $\Delta$ is a morphism of algebras.
	
	\item \label{item:b} We have
	\begin{align}
		\Delta\big( \shuffle_\qv \big) \subseteq \bigoplus_{\bfc, \bfd\in \N I} \shuffle_{\bfc, \bfd} \Big\llbracket z^{(2)}_{j,k}/z^{(1)}_{i,\ell}, z^{(2)}_{j,k}/h^{(1)}_{i,\ell}\,\Big\vert\, i,j \in I,\, k,\ell \in \N\Big\rrbracket\ ,
	\end{align}
	where $\shuffle_{\bfc, \bfd}\coloneqq \shuffle_{\bfc}^{\frake}\otimes \shuffle_\bfd$. Here, $x^{(1)}$ and $x^{(2)}$ stand for $x \otimes 1$ and $1 \otimes x$ respectively.
	
	\item \label{item:c} We have
	\begin{align}
		\Delta\big( \cohaqv^{\Ttilde} \big) \subseteq \Big( \big(\shuffle_\qv^0 \ltimes \cohaqv^{\Ttilde} \big) \otimes \cohaqv^{\Ttilde} \Big)_{\mathsf{loc}}
	\end{align}
	in the sense that 
	\begin{align}
		(\pi_m \otimes 1) \big(\Delta\big( \cohaqv^{\Ttilde} \big)\big) \subset \cohaqv^{\Ttilde} 
	\end{align}
	\begin{align}
		(1 \otimes \pi_m ) \big( \Delta\big( \cohaqv^{\Ttilde} \big)\big) \subset(\shuffle_\qv^0\ltimes \cohaqv^{\Ttilde})\Big[h_{j,\ell}^{-1}, z_{j,k}^{-1}\,\vert\, i,j \in I,\, k,\ell \in \N \Big] 
	\end{align}
	for every
	\begin{align}
		m =\prod_{j,k} z_{j,k}^{a_{j,k}} \cdot \prod_{j,k} h_{j,k}^{b_{j,k}}\ ,
	\end{align}
	where $\pi_m$ is the linear form picking the coefficient of the monomial $m$.
	
	\item \label{item:d} For any $x\in\shuffle_\bfd$, we have 
	\begin{align}
		\Delta_{\bfd,0}(x)=x \otimes 1\quad \text{and} \quad \Delta_{0,\bfd}(x) = 1 \otimes x+ P_x \ ,
	\end{align}
	where
	\begin{align}
		P_x \in \sum_{j,\ell} \big(h_{j,\ell}^{-1} \otimes \shuffle_\bfd\big)\Big\llbracket z^{(2)}_{j,k}/h^{(1)}_{i,\ell}\,\Big\vert\, i,j \in I, \,k,\ell \in \N\Big\rrbracket\ .
	\end{align} 
\end{enumerate}

\begin{proof}[Proof of Lemma~\ref{lem:equality-coha}] 
	Note that $\coha_{\alpha_i}^{\Ttilde} = \shuffle_{\alpha_i}$, hence it is enough to check that 
	\begin{align}
		\cohaqv^{\Ttilde} \cap (\shuffle_\qv \star \shuffle_{\alpha_i})\subset \cohaqv^{\Ttilde} \star \shuffle_{\alpha_i}\ .
	\end{align}
	Let us write $\coha_{>n\alpha_i}^{\Ttilde}\coloneqq \bigoplus_{m > n} \coha_{m\alpha_i}^{\Ttilde}$. We claim that 
	\begin{align}\label{eq:proof-coproduct-1}
		\cohaqv^{\Ttilde} \cap ( \shuffle_\qv \star \shuffle_{n\alpha_i}) \subseteq \cohaqv^{\Ttilde} \star \shuffle_{n\alpha_i} + \shuffle_\qv \star \shuffle_{>n\alpha_i}
	\end{align}
	for $n\geq 0$. Indeed, let $\bfc \in \N I$ and $\bfc=\bfd+n\alpha_i$. We choose a basis $\{h_\ell\}$ of $\shuffle_{n\alpha_i}$ and elements $g_\ell\in\shuffle_\bfd$. Set 
	\begin{align}
		f\coloneqq \sum_\ell g_\ell\star h_\ell\in\shuffle_\bfd \star \shuffle_{n\alpha_i}\ .
	\end{align}
	From properties~\eqref{item:a} and \eqref{item:d} above, we deduce
	\begin{align}\label{eq:proofcoprodloc}
		\Delta_{n\alpha_i, \bfd}(f) =&\sum_\ell \Delta_{0,\bfd}(g_\ell) \star \Delta_{n\alpha_i,0}(h_\ell) + 
		\sum_\ell\sum_{k = 1}^n  \Delta_{k\alpha_i,\bfd-k\alpha_i}(g_\ell) \star \Delta_{n\alpha_i-k\alpha_i,k\alpha_i}(h_\ell)\\
		&\in \sum_\ell \Big( h_\ell \otimes g_\ell + P_{g_\ell}\star h_\ell \Big)+ 
		\sum_\ell\sum_{k = 1}^n  \big(\shuffle^{\mathsf{e}}_{n\alpha_i} {\otimes} ( \shuffle_{\bfd-k\alpha_i} \star \shuffle_{k\alpha_i})\big)_{\mathsf{loc}}\ .
	\end{align}
	If $f\in \cohaqv^{\Ttilde}$ then by applying the linear form $(h_\ell)^\ast  \otimes 1$ to Formula~\eqref{eq:proofcoprodloc} and using Property~\eqref{item:c} we get
	\begin{align}
		g_\ell \in \cohaqvd^{\Ttilde} + \sum_{k=1}^n \shuffle_{\bfd-k\alpha_i} \star \shuffle_{k\alpha_i}\ .
	\end{align}
	Thus, Formula~\eqref{eq:proof-coproduct-1} holds.
	
	We can now finish the proof of Lemma~\ref{lem:equality-coha}. Note that if $n$ is positive then
	\begin{align}
		\coha_{\bfc-n\alpha_i}^{\Ttilde} \star \shuffle_{n\alpha_i}=&\coha_{\bfc-n\alpha_i}^{\Ttilde} \star (\shuffle_{\alpha_i})^{\star n}\\
		=&\coha_{\bfc-n\alpha_i}^{\Ttilde} \star (\coha_{\alpha_i})^{\star (n-1)}\star\shuffle_{\alpha_i} \subset\coha_{\bfc-\alpha_i}^{\Ttilde} \star \shuffle_{\alpha_i}\ .
	\end{align}
	Thus, by repeated use of \eqref{eq:proof-coproduct-1} we get
	\begin{align}
		\coha_\bfc^{\Ttilde} \cap \big(\shuffle_\qv \star \shuffle_{\alpha_i}\big) &\subset \big(\coha_{\bfc-\alpha_i}^{\Ttilde} \star \shuffle_{\alpha_i}\big) + \coha_\bfc^{\Ttilde} \cap \big(\shuffle_\qv \star \shuffle_{>\alpha_i}\big)\\ 
		&\subset \big(\coha_{\bfc-\alpha_i}^{\Ttilde} \star \shuffle_{\alpha_i}\big) + \coha_\bfc^{\Ttilde} \cap \big(\shuffle_\qv \star \shuffle_{>2\alpha_i}\big)\subset\cdots
	\end{align}
	On the other hand, for $n$ big enough we have $\shuffle_{\bfc-n\alpha_i}=\{0\}$ and hence $\coha_\bfc^{\Ttilde} \cap \big(\shuffle_\qv \star \shuffle_{\alpha_i}\big) = \coha_{\bfc-\alpha_i}^{\Ttilde} \star \shuffle_{\alpha_i}$ as wanted.
\end{proof}  

\subsection{Partial doubles at a fixed vertex}\label{subsec:partial-doubles}

Recall that we have fixed a vertex $i\in I$. In this section, we define a truncated braid group operator $\overline{T}_i \colon \icohatildeqv \to \coha^{\Ttilde,(i)}_\qv$. For this, we first define a `partial double' $\Dicohaqv$ of $\cohaqv^{\Ttilde}$, by completing the nilpotent subalgebra $\Y^-_i \subset \cohaqv^{\Ttilde}$ into the full Yangian $\Y_i$. For simplicity, we will denote by the same symbol $x^-_{j,\ell}$ the generators of the Yangian $\Y^-_\qv$ and their image under the morphism $\Phi\colon \Y^-_\qv \to \cohaqv^{\Ttilde}$.

\begin{definition}
	The \textit{partial double of $\cohaqv^{\Ttilde}$ at the vertex $i$} is the algebra $\Dicohaqv$ generated by $\cohaqv^{\Ttilde}$ together with elements $x_{i, \ell}^+, h_{i, \ell}$, with $\ell\in \N$, subject to the following relations:
	\begin{itemize}\itemsep0.2cm
		\item the elements $x_{i, \ell}^+$, $h_{i, \ell}$, and $x_{i, \ell}^-$ for $\ell\in \N$ satisfy the defining relations of the Yangian $\Y_i$ of type $\fraksl_2$.
		
		\item $[x_{i,k}^+, x_{j, \ell}^-]=0$ for all $k, \ell\in \N$ and $j,i\in I$ and $j\neq i$.
		
		\item the elements $h_{i, k}$ and $x_{j, \ell}^-$ satisfy the relations \eqref{eq:Yangian-7} for $k, \ell\in \N$.
	\end{itemize} 
	In a similar way, we define the \textit{partial double $\DiYqv$ of $\Y_\qv^-$ at the vertex $i$}. It is naturally a subalgebra of $\Y_\qv$.
\end{definition}
By construction, there is a surjective algebra homomorphism $\Phi\colon \DiYqv\to \Dicohaqv$ extending the algebra homomorphism introduced in Formula~\eqref{eq:Phi}.

Note that $(x_{i, 0}^+, h_{i, 0}, x_{i, 0}^-)$ forms a $\fraksl(2)$-triplet which acts locally nilpotently on $\Dicohaqv$, thus giving rise to an operator $T_i\in \Aut(\Dicohaqv)$. It fits into the commutative diagram
\begin{align}
	\begin{tikzcd}[ampersand replacement=\&]
		\DiYqv  \ar{r}{T_i}  \ar[swap, twoheadrightarrow]{d}{\Phi}\&\DiYqv \ar[twoheadrightarrow]{d}{\Phi}\\
		\Dicohaqv \ar{r}{T_i} \&   \Dicohaqv
	\end{tikzcd}\ ,
\end{align}
where $T_i\in \DiYqv$ is induced by Formula~\eqref{eq:Ti}.

Let us denote by $\coha_i^{\Ttilde,+}$ and $\coha_i^{\Ttilde,0}$ the subalgebras of $\Dicohaqv$ generated by $x_{i, \ell}^+$ and $h_{i, \ell}$, with $\ell\in \N$, respectively.
\begin{lemma}
	The multiplication map 
	\begin{align}
		m\colon \coha_i^{\Ttilde,+}\otimes \coha_i^{\Ttilde,0}\otimes \cohaqv^{\Ttilde} \longrightarrow\Dicohaqv
	\end{align}
	is an isomorphism.
\end{lemma}

\begin{proof}
	Let us denote by $\YMO$ the \textit{Maulik-Okounkov Yangian} of $\qv$. By \cite[Theorem~1.1--(a) and --(b)]{SV_Yangians=COHA}, it carries a triangular decomposition $\YMO^+ \otimes \YMO^0 \otimes \YMO^-$, and there is an algebra isomorphism $\cohaqv^{\Ttilde, \leq 0} \stackrel{\sim}{\longrightarrow} \YMO^{\leq 0}$. The defining relations of $\Dicohaqv$ imply the existence of a factorization $\cohaqv^{\Ttilde} \to \Dicohaqv\to \YMO$, which shows the injectivity of the map $\cohaqv^{\Ttilde} \to \Dicohaqv$. In addition, the elements of $\YMO^0$ correspond to operators of multiplication by Chern classes of the tautological bundles, as discussed in e.g. \cite[\S3.4.3]{SV_Yangians}. Finally, by \cite[Proposition~11.3.2]{MO_Yangian}, there is an inclusion $\eY_{i} \to \YMO$ corresponding to the vertex $i$. 
	
	Thus, there exists a commutative diagram
	\begin{align}
		\begin{tikzcd}[ampersand replacement=\&]
			\coha_i^{\Ttilde,+}\otimes \coha_i^{\Ttilde,0}\otimes \cohaqv^{\Ttilde}  \ar{r}{m}  \ar{d}\&\Dicohaqv \ar{d}\\
			\YMO^+\otimes\YMO^0\otimes \YMO^- \ar{r}{m} \&   \YMO
		\end{tikzcd}\ ,
	\end{align}
	with injective vertical homomorphisms. Since the bottom horizontal homomorphism is an isomorphism the top horizontal homomorphism is injective. The surjectivity can be proved by standard straightening arguments, because by \cite[Propositions~5.8 and 5.12]{SV_generators}, $\cohaqv^{\Ttilde}$ is generated by $x_{j, \ell}^-$ for $j\in I$ and $\ell\in \N$. 
\end{proof}

We may now define the truncated braid group operator $\overline T_i$ on $\cohaqv^{\Ttilde}$ as before. Consider the composition 
\begin{align}
	\begin{tikzcd}[ampersand replacement=\&]
		\cohaqv^{\Ttilde}\ar{r} \& \Dicohaqv \ar{r}{T_i} \& \Dicohaqv \ar{r}\& \cohaqv^{\Ttilde}\ ,
	\end{tikzcd}
\end{align}
where the last map is the projection induced by the triangular decomposition.

\begin{proposition}\label{prop:commutativity-Ti}
	The map $\overline T_i \colon \cohaqv^{\Ttilde}\to \cohaqv^{\Ttilde}$ descends to an isomorphism $\overline T_i\colon \icohatildeqv\to \cohatildeiqv$ such that the following diagram commutes:
	\begin{align}\label{eq:commutative-truncated}
		\begin{tikzcd}[ampersand replacement=\&]
			\iY_\qv  \ar{r}{\overline T_i}  \ar[swap,twoheadrightarrow]{d}{\Phi}\&\Yi_\qv \ar[twoheadrightarrow]{d}{\Phi}\\
			\icohatildeqv \ar{r}{\overline T_i} \&   \cohatildeiqv
		\end{tikzcd}\ ,
	\end{align}
	where $\overline T_i\colon \iY_\qv \to \Yi_\qv$ is defined in Proposition~\ref{prop:overlineTi}.
\end{proposition}
\begin{proof} 
	The fact that the map $\overline T_i \colon \cohaqv^{\Ttilde}\to \cohaqv^{\Ttilde}$ descends to a homomorphism $\overline T_i\colon \icohatildeqv\to \cohatildeiqv$ is proved just in the case of the Yangian, see Proposition~\ref{prop:overlineTi}. That $\overline T_i\colon \icohatildeqv\to \cohatildeiqv$ is an isomorphism follows from the fact that $\overline T_i\colon \iY_\qv \to \Yi_\qv$ is an isomorphism, that the map $\Phi$ is onto, and that $\icohatildeqv$ and $\cohatildeiqv$ have the same graded dimensions, see Proposition~\ref{prop:reflection-stacks}--\eqref{item:reflection-stacks-4}).
\end{proof}

\begin{remark}
	The same construction applied to $\DiYqv$ gives back the operator $\overline T_i$ defined in Formula~\eqref{eq:overlineT}.
\end{remark}

There is a canonical sign-twist isomorphism $\cohaqv^{\Ttilde} \simeq \coha_{s_i\qv}^{\Ttilde}$, compatible with the maps $\Y_\qv \to \cohaqv^{\Ttilde}$ and $\Y_\qv \to \coha_{s_i\qv}^{\Ttilde}$. It descends to an isomorphism $\cohatildeiqv \simeq \coha^{\Ttilde,(i)}_{s_i\qv}$. We will denote by the same letter the induced truncated braid group morphism $\overline{T}_i\colon\icohatildeqv \to \coha^{\Ttilde,(i)}_{s_i\qv}$.

\subsection{Proof of Theorem~\ref{thm:compatibility-braid-reflection-functors}}\label{subsec:proof-Theorem-compatibility}

In this section, we put together all the results of \S\ref{subsec:degree-zero-case}--\ref{subsec:partial-doubles} in order to finally prove the commutativity of the diagram~\eqref{eq:diag:T_i=S_i}. 

In view of Proposition~\ref{prop:commutativity-Ti}, Theorem~\ref{thm:compatibility-braid-reflection-functors} is a consequence of the following result:
\begin{proposition}\label{prop:S_i=T_i}
	The maps $S_{i,\, \ast}, \overline T_i\colon \icohatildeqv \to \coha^{\Ttilde,(i)}_{s_i\qv}$ are equal.
\end{proposition}

\begin{proof}
	Set 
	\begin{align}
		\Psi\coloneqq S_{i, \,\ast} - \overline T_i\colon \icohad \longrightarrow \cohai_{s_i(\bfd)}\ .
	\end{align}
	Our aim is to prove that $\Psi=0$. Since $\Ttilde=\Tmax$, all the torus weights $\epsilon_e$ satisfy the following condition
	\begin{align}\label{eq:gen}
		\epsilon_e+\epsilon_{e'} \neq 0\ ,
	\end{align}
	for $j, i\in I$, with $j\neq i$, $e\in \Omega_{i, j}$, and $e'\in \Omega_{j,i}$. Recall the automorphism $T_i \in \Aut(\bS)$ from \S\ref{subsec:tautological-classes-braid}. There is a commutative diagram
	\begin{align}
		\begin{tikzcd}[ampersand replacement=\&]
			\bS \ar{r}{T_i} \ar[d] \& \bS\ar[d]\\
			\End(\icohad) \ar{r}{S_{i,\, \ast}} \& \End(\cohai_{s_i(\bfd)})
		\end{tikzcd}\ ,
	\end{align}
	in which the vertical arrows are induced by the canonical actions of $\bS$ on $\icohad$ and $\cohai_{s_i(\bfd)}$. We use $T_i$ to equip the space $\cohai_{s_i(\bfd)}$ with the structure of a $\shuffle_\bfd$-module. From Theorem~\ref{thm:braid-Chern} we deduce that $\Psi$ is a morphism of $\shuffle_\bfd$-modules. Next, by Theorem~\ref{thm:compatibility-braid-reflection-functors-degree-zero}, the map $\Psi$ factors to a map of $\shuffle_\bfd$-modules 
	\begin{align}\label{eq:overline-Psi}
		\overline \Psi\colon \icohad/(\shuffle_\bfd\cdot  \tensor*[^{(i)}]{\coha}{_{0, \, \bfd}})\longrightarrow \cohai_{s_i(\bfd)}\ .
	\end{align}
	By comparing the supports of both sides in $\bfC^{(\bfd)}$, we will show that $\overline\Psi=0$, hence $\Psi=0$. 
	
	Let us first explicitly describe the image, under the embedding $\iota\colon \cohaqv^{\Ttilde} \to \shuffle_\qv$ introduced in Theorem~\ref{thm:iota}, of certain elements in $\cohazero^{\Ttilde}$.
	
	For any $\bfd\in\N I$, let $[\stackRep_\bfd(\qv)]\in \cohazero^{\Ttilde}$ be the fundamental class of the zero section of the map $\stackRep_\bfd(\doubleqv)\to \stackRep_\bfd(\qv)$. As $\qv$ carries no edge loops, we have $\iota[\stackRep_{d_j\alpha_j}(\qv)]= 1 \in \shuffle_{d_j\alpha_j}$ for $j\in I$ and $d_j\in \N$. Fix a total order $i_1 \prec i_2 \prec \cdots \prec i_r$ on the vertex set $I$. Then, for any $\bfd\in \N I$ we have
	\begin{align}\label{eq:ordered-product}
		\iota( [\stackRep_{d_{i_1}\alpha_{i_1}}(\qv)] \cdots  [\stackRep_{d_{i_r}\alpha_{i_r}}(\qv)])= \eta
		\cdot \prod_{s<t}\prod_{\genfrac{}{}{0pt}{}{e\in\doubleOmega}{e\colon i_s\to i_j}} \prod_{\genfrac{}{}{0pt}{}{1 \leq k \leq d_{i_s}}{1 \leq \ell \leq d_{i_t}}}(z_{i_s,k}-z_{j_t,\ell}+\varepsilon_e)
	\end{align}
	for some sign $\eta$, given by
	\begin{align}
		\eta\coloneqq \prod_{s=1}^{r-1}(-1)^{-\sum_{s<t}d_{i_s}d_{i_t}\#\{e\colon i_s\to i_{t}\in \Omega\}}\ .
	\end{align}
	We will simply write the above product as
	\begin{align}
		u_\prec\coloneqq \eta\cdot \prod_{\genfrac{}{}{0pt}{}{i\prec j}{e\colon i\to j}} (z_{i,\bullet}-z_{j,\bullet} + \varepsilon_e) \ .
	\end{align}
	
	By construction, the annihilator of the source of the map \eqref{eq:overline-Psi} in $\shuffle_\bfd$ contains $\tensor*[^{(i)}]{\coha}{_{0, \, \bfd}}$. Thus it contains all the elements $u_\prec$ in Formula~\eqref{eq:ordered-product}, for any total order $\prec$ on $I$. Therefore, the support of the source of the map \eqref{eq:overline-Psi} is contained in the intersection of the zero loci of all possible elements $u_\prec$. Hence, by Proposition~\ref{prop:support}, it suffices to check that for any irreducible component $\calC_O$ of $Z_{s_i(\bfd), i}$, there exists an element of the form $u_{\prec}$ which does not identically vanish on $\calC_O$, hence $Z_{s_i(\bfd), i}$ is not contained in the intersection of the zero loci of all possible elements $u_\prec$. We claim that this holds for any order $\prec$ in which $i$ is maximal. Indeed, for such an order $\prec$, we have
	\begin{align}\label{eq:Ti-ordered-product}
		\begin{split}
			\overline T_i(u_\prec)=&\eta
			\prod_{\genfrac{}{}{0pt}{}{j\prec k}{e\colon j \to k}}(z_{k,\bullet}-z_{j,\bullet} + \epsilon_e) \cdot 
			\prod_{\genfrac{}{}{0pt}{}{j,k\in I}{e\colon j \to i\, ,e'\colon k \to i}}(z_{k,\bullet}-z_{j,\bullet} + \epsilon_e+\epsilon_{e'})\cdot\\
			&\cdot \prod_{\genfrac{}{}{0pt}{}{j\in I}{e\colon i \to j}}(z_{i,\bullet}-z_{j,\bullet}+\hbar + \epsilon_e)^{-1} \ ,
		\end{split}
	\end{align}
	where $j,k$ run over elements in $I \smallsetminus\{i\}$. Using the description \eqref{eq:calC-calO}, it is easy to check that such an element does not identically vanish on any of the irreducible components $\calC_O$. This proves the claim.
\end{proof}

\section{Appendix: Sign twists of algebras}\label{sec:sign-twisting}

The following operation of twisting the multiplication in graded associative algebras is used in several instances in this work. Let $\Gamma\coloneqq\Z I$ be a free abelian group of finite rank and let $(A,\star)$ be a $\Gamma$-graded associative algebra, $A=\bigoplus_\gamma A_\gamma$. Let $\omega \colon \Gamma\times \Gamma \to \Z$ be a bilinear form. We define a new algebra $(A^\omega,\underset{\omega}{\star})$ by setting $A\coloneqq A^\omega$ and for homogeneous elements $x,y \in A^\omega$
\begin{align}
	x \underset{\omega}{\star} y \coloneqq (-1)^{\omega(\vert x\vert ,\vert y\vert )} x \star y\ .
\end{align}
This new product is associative. Moreover, the twisting operation is associative in the sense that $(A^\omega)^{\sigma}=A^{\omega+\sigma}$.

\begin{proposition}\label{prop:twist-sign-antisym}
	Let $A,\Gamma,\omega$ be as above and assume that $\omega$ is antisymmetric. Fix a basis $e_1, \ldots, e_n$ of $\Gamma$ and for each $\gamma=\sum \gamma_i e_i \in \Gamma$ we set
		\begin{align}
			u_\gamma\coloneqq (-1)^{\sum_{i <j}\gamma_i\gamma_j \omega(e_i,e_j)}\ .
		\end{align}
	Then the map $f \colon A \to A^\omega$, given by $x \mapsto  u_{\vert x\vert}x$ is an algebra isomorphism.
\end{proposition}

\begin{proof} 
	It suffices to check that $f$ is an algebra morphism, which reduces to the relation 
	\begin{align}
		u_{\gamma+\sigma} =u_\gamma u_\sigma (-1)^{\omega(\gamma,\sigma)}\ .
	\end{align}
	Setting $\varepsilon_{ij}=\omega(e_i,e_j)$ we have 
	\begin{align}
		\frac{u_{\gamma+\sigma}}{u_\gamma u_\sigma}&=(-1)^{\sum_{i<j} ((\gamma_i+\sigma_i)(\gamma_i+ \sigma_j)-\gamma_i\gamma_j - \sigma_i\sigma_j)\varepsilon_{ij}}=(-1)^{\sum_{i<j}(\sigma_i\gamma_j +\gamma_i\sigma_j)\varepsilon_{ij}}=(-1)^{\sum_{i<j}(\sigma_i\gamma_j -\gamma_i\sigma_j)\varepsilon_{ij}}\\
		&=(-1)^{\omega(\gamma,\sigma)}\ .	\tag*{\qedhere} 
	\end{align}
\end{proof}
As a corollary, we see that two bilinear forms $\Gamma\times\Gamma \to \Z$ which have the same symmetrization give rise to isomorphic twisted algebras. The isomorphism given above is explicit but not canonical (it depends on the choice of a basis). 

\begin{remark}
	Note that $u_{e_i}=1$ for any $i$. If $A$ is generated by $\bigoplus_{i} A_{e_i}$, this condition fully determines the isomorphism.
\end{remark}

\newpage
\part{COHA of a minimal resolution of a Kleinian singularity and affine Yangians}\label{part:ADE-case}

\subsubsection*{Notation} 

Throughout this Part, we fix a finite subgroup $G$ of $\mathsf{SL}(2, \C)$ and we let $\qvfin=(I_\sff, \Omega_\sff)$ be the corresponding finite type ADE Dynkin diagram, with $I_\sff=\{1, \ldots, e\}$ and an arbitrary orientation. Let $\qv=(I, \Omega)$ be the affine extension of $\qvfin$, with $I=I_\sff\cup \{0\}$. Let $\doubleqv=(I, \doubleOmega=\Omega\sqcup\{a^\ast\,\vert\, a\in \Omega\})$ be the double quiver associated to $\qv$.

\section{Lie theory of affine ADE quivers}\label{sec:affine-quivers}

Let us denote by $\frakg$ and $\frakgfin$ the affine Kac-Moody algebra and the semisimple Lie algebra associated to $\qv$ and $\qvfin$, respectively. This short section serves mainly to fix some standard notations related to (co)root systems of $\frakg, \frakgfin$, and the associated braid groups. 

\subsection{Affine and finite roots and coweights}\label{subsec:affine-quivers-roots-coweights}

Let $\{\alpha_0, \ldots, \alpha_e\}$ denote the simple roots of $\frakg$, so that the root and coroot systems of $\frakg$ are, respectively,
\begin{align}
	\rootlattice=\bigoplus_{i\in I} \Z \alpha_i \quad \text{and} \quad \corootlattice=\bigoplus_{i \in I} \Z \check{\alpha}_i\ .	
\end{align}
We denote by $\Delta \subset \rootlattice$ the subset of roots. We denote by $\Delta^{\mathsf{re}}$ and $\Delta^{\mathsf{im}}$ the subsets of real and imaginary roots. We use the notations $\rootlatticefin, \corootlatticefin, \rootsetfin$ for the analogous notions for $\frakgfin$. Denote by $\varphi$ the \textit{highest} positive root in $\rootsetfin$ and by $\Lvarphi$ its associated coroot. 

Let $\delta\in \rootlattice$ be the \textit{indivisible positive imaginary root}, i.e., the minimal element in the intersection of $\N I$ with the kernel of the canonical pairing \eqref{eq:canonical-pairing}, so that $\Delta^{\mathsf{im}}=\Z \delta$. We have
\begin{align}\label{eq:delta}
	\delta=\varphi+\alpha_0=\sum_{i\in I}r_i\, \alpha_i \ ,
\end{align}
where $r_0\coloneqq 1$. Then
\begin{align}\label{eq:Coxeter}
	h\coloneqq \sum_{i\in I}r_i
\end{align}
is the Coxeter number of $\qv$. The explicit values of $h$ and the $r_i$'s are:
\begin{align}
	\begin{array}{|c|c|c|}
		\hline
		& h & r_1, r_2, r_3, \ldots, r_{e-1}, r_e\\		
		\hline
		A_N & N+1 & 1, \ldots, 1\\
		\hline
		D_N & 2N-2 & 1, 2, 2, \ldots, 2, 1, 1\\
		\hline
		E_6 & 12 & 1, 2, 3, 2, 1, 2\\
		\hline
		E_7 & 16 & 1, 2, 4, 3, 2, 1, 2\\
		\hline
		E_8 & 30 & 2, 4, 6, 5, 4, 3, 2, 3\\
		\hline
	\end{array}
\end{align}

Let $\{\Llambda_i\,\vert\, i \neq 0\}$ be the fundamental coweights of $\qvfin$ and let $\{\Lomega_i\,\vert\, i \in I\}$ be the fundamental coweights of $\qv$. For $i \neq 0$, the coweight $\Llambda_i$ can be expressed as 
\begin{align}
	\Llambda_i= \Lomega_i-r_i\Lomega_0\ .
\end{align}
We put
\begin{align}\label{eq:r_i}
	\Lrho_\sff\coloneqq\sum_{i \neq 0}\Llambda_i \quad \text{and}\quad \Lrho =\sum_i \Lomega_i=\Lrho_\sff+(1+h)\Lomega_0\ .
\end{align} 
Let 
\begin{align}
	\coweightlattice\coloneqq \bigoplus_{i}\Z \Lomega_i \quad \text{and}\quad \coweightlattice^+\coloneqq \bigoplus_{i}\Z_{\geq 0} \Lomega_i
\end{align}
be the sets of coweights and \textit{dominant} coweights of $\qv$, respectively. The corresponding lattices and monoid for $\qvfin$ are denoted by $\coweightlatticefin, \coweightlatticefin^+$. A coweight $\Ltheta=\sum_i a_i\Llambda_i$ will be called \textit{strictly dominant} if $a_i >0$ for all $i$. 
For any $\bfd\in \Z I$, we write $\bfd=\bfd_\sff+n\delta$, with $\bfd_\sff\in \Z\rootsetfin$ and $n\in\Z$, and we call $(\Lrho_\sff, \bfd)$ the \textit{height of $\bfd_\sff$}. 

\subsection{Extended, affine, and finite braid groups}\label{subsec:braid-group-affine-quiver}

Let $W_\sfaf\coloneqq W_\qv$ be the \textit{affine Weyl group}, i.e., the Weyl group of the affine ADE quiver $\qv$. The Weyl group $W\coloneqq W_{\qvfin}$ of $\qvfin$ coincides with the subgroup of $W_\sfaf$ generated by the simple reflections $s_i$ with $i\neq 0$. Let $w_0$ be the longest element in $W$. 

Let $W_\sfex$ be the \textit{extended affine Weyl group}. By definition, we have $W_\sfex=W\ltimes\coweightlatticefin=\Gamma\ltimes W_\sfaf$, where $\Gamma$ is the group of outer automorphisms of the underlying diagram of $\qv$. We abbreviate $t_{\Llambda}\coloneqq 1\ltimes\Llambda$ for each coweight $\Llambda\in \coweightlatticefin$. The $W_\sfaf$-action on $\C I$ extends to a $W_\sfex$-action such that 
\begin{align}\label{eq:t_lambda}
	t_{\Llambda}(\bfd)=\bfd-(\Llambda,\bfd)\,\delta
\end{align}
for any $\bfd\in \C I$.

Let $B\coloneqq B_{\qvfin}$ be the \textit{braid group}, $B_\sfaf\coloneqq B_\qv$ the \textit{affine braid group}, and $B_\sfex$ the \textit{extended affine braid group}. The group $B_\sfex$ is generated by elements $T_w$ with $w\in W_\sfex$ modulo the relations 
\begin{align}\label{eq:defBex}
	T_vT_w=T_{vw}\quad \text{if }\,\ell(vw)=\ell(v) + \ell(w)\ ,
\end{align}
where $\ell(u)$ is the standard length function. We abbreviate $T_{\Llambda}\coloneqq T_{t_{\Llambda}}$ for each coweight $\Llambda\in \coweightlatticefin$ and $T_i\coloneqq T_{s_i}$ for each affine simple reflection. The affine braid group $B_\sfaf$ coincides with the subgroup of $B_\sfex$ generated by the $T_i$'s for $i\in I$, while the braid group $B$ coincides with the subgroup of $B_\sfex$ generated by the $T_i$'s with $i\neq 0$. 

Define the elements $L_{\Llambda}$ with $\Llambda\in \coweightlatticefin$ as 
\begin{align}\label{def:Llambda}
	L_{\Llambda}\coloneqq \begin{cases}
		T_{\Llambda} & \text{ if $\Llambda$ is dominant}\ , \\[2pt]
		T_{\Llambda_1}T_{\Llambda_2}^{-1} &\text{if  $\Llambda=\Llambda_1-\Llambda_2$, with $\Llambda_1$, $\Llambda_2$ are both dominant} \ .
	\end{cases}
\end{align}
It is known that $L_{\Llambda}= T_{\Llambda_1}T_{\Llambda_2}^{-1}$ does not depend on the choice of the decomposition of $\Llambda$. Moreover $L_0=1$.

We recall the following result, which will be used later on -- see, e.g., \cite[Chapter~3, \S3.3]{Mc}.
\begin{proposition}\label{prop:characterization-extended-affine-braid-group}
	The extended affine braid group $B_\sfex$ is generated by $B\cup\{L_{\Llambda}\,\vert\,\Llambda\in\coweightlatticefin\}$ subject to the following relations:
	\begin{align}\label{eq:extended-B-(1)}
		L_{\Llambda_1}L_{\Llambda_2}&=L_{\Llambda_2}L_{\Llambda_1}\\[2pt] \label{eq:extended-B-(2)}
		T_iL_{\Llambda}&=L_{\Llambda} T_i
	\end{align}
	when $s_i(\Llambda)=\Llambda$, and 
	\begin{align}\label{eq:extended-B-(3)}
		L_{\Llambda}=T_iL_{\Llambda-\Lalpha_i} T_i
	\end{align}
	when $s_i(\Llambda)=\Llambda-\Lalpha_i$. 
\end{proposition}

\begin{remark}\label{rem:braid-group}
	We have $T_0= L_{\Lvarphi} T_{s_\varphi}^{-1}$ by \cite[Formula~(3.2.8)]{Mc}, and there is a group isomorphism $B_\sfex\simeq \Gamma\ltimes B_\sfaf$. 
\end{remark}

Recall from \S\ref{sec:braid-group-action-Yangian} that there is an action of the braid group $B_\sfaf$ on $\sfU(\frakg)$. Combined with the obvious action of $\Gamma$, it gives rise to an action of $B_\sfex$ on $\sfU(\frakg)$. We will need the following explicit result for the action of translation elements.
\begin{proposition}\label{prop:action_translation_braid_Lie}
	Under the isomorphism $\frakg \simeq \uce\big(\frakgfin[s^{\pm 1}]\big)$, for any root $\alpha \in \Delta_\sff$, any $\Llambda \in \coweightlatticefin$, and any $n \in \Z$ we have 
	\begin{align}
		L_{\Llambda}(xs^n)=(-1)^{\langle \Llambda,\alpha\rangle}xs^{n-\langle \Llambda,\alpha\rangle}
	\end{align}
	for any $x \in \frakg_\alpha$.	
\end{proposition}
\begin{proof} 
For $\qv$ of type $A_1^{(1)}$ this can be checked by an explicit computation (it is enough to consider the action of $L_{\Lomega}=\gamma T_0$, where $\gamma$ is the automorphism exchanging the vertices $0$ and $1$ of $A_1^{(1)}$). 

\medskip

To deal with the general case, it is sufficient to assume that $\alpha \in \{\pm \alpha_i\,|\, i \neq 0\}$. This may in turn be deduced from the classical limit of \cite[Corollary after Proposition~3.8]{Beck_Quantum}, which yields, for any $i,j \in I\smallsetminus\{0\}$ a commutative diagram
\begin{align}
	\begin{tikzcd}[ampersand replacement=\&]
		\coweightlatticefin \arrow{r}{} \& \Aut(\sfU_j)\\
		(\coweightlatticefin)_{\fraksl_2} \arrow{r}{} \arrow{u}{u_i} \& \Aut(\widehat{\fraksl_2}) \arrow[swap]{u}{\gamma_j}
	\end{tikzcd}
\end{align}
where
\begin{align}
	\begin{tikzcd}[ampersand replacement=\&]
		\gamma_j\colon\widehat{\fraksl_2} \ar{r}{\sim}\& \sfU_j \subset \frakg
	\end{tikzcd}
\end{align}
is the canonical embedding corresponding to $j$ and $u_i\colon (\coweightlatticefin)_{\fraksl_2} \to \coweightlatticefin$ is defined by $\Lomega \mapsto \Lomega_i$. Notice that only the case $i=j$ is dealt with in \textit{loc. cit.} but the general case is similar.
\end{proof}

\section{Yangians of affine ADE quivers} \label{sec:Yangians-affine-quivers}

Yangians $\Y_\qv$ associated to arbitrary quivers and maximal number of quantum parameters were introduced in Part~\ref{part:Yangians}, see Definition~\ref{def:Yangian}. In the present section, we focus on the situation of an affine quiver $\qv$. We prove that, in that case, a two-parameters version of $\Y_\qv$ is isomorphic to the \textit{affine Yangian} introduced and studied in e.g. \cite{Varagnolo_Yangian, Guay07, GRW19, Kodera21, BT-Yangians}, which is a quantization of the elliptic Lie algebra $\frakgell\coloneqq \uce(\frakgfin[s^{\pm 1},t])$, a universal central extension of the double loop algebra $\frakgfin[s^{\pm 1},t]$. Recall that in our notation, $\frakgfin$ is the semisimple (finite-dimensional) Lie algebra associated to the finite type quiver $\qvfin$. We also prove that the map $\Y^-_\qv \to \cohaqv^T$ from the Yangian to the twisted equivariant nilpotent cohomological Hall algebra of $\qv$ given in Theorem~\ref{thm:Phi} is an isomorphism in that case. The definition of the affine Yangian in the $A_1^{(1)}$ case is a little bit different from that of other affine cases, we take extra care treating it.

\subsection{Definition of the affine Yangian}\label{subsec:Yangians-affine}

Fix the torus $T\coloneqq \C^\ast\times \C^\ast$. Let $\varepsilon_i$ be the first Chern classes of the tautological characters $\gamma_i\in \C^\ast$ for $i=1, 2$. Then $\Hbullet_T\coloneqq \Hbullet_T(\mathsf{pt})=\Q[\varepsilon_1, \varepsilon_2]$. We set $\hbar\coloneqq \varepsilon_1+\varepsilon_2$.

\begin{definition}\label{def:affine-Yangian}
	Let $\qv$ be an affine ADE quiver, $\qv\neq A_1^{(1)}$. The \textit{extended (two-parameters) Yangian $\eY_{\qv;\, \varepsilon_1, \varepsilon_2}$ of $\qv$} is the unital associative $\Hbullet_T$-algebra generated by $x_{i, \ell}^\pm, h_{i, \ell}, \kappa_{i, \ell}$, with $i \in I$ and $\ell \in \N$, subject to the relations
	\begin{itemize}\itemsep0.2cm
		\item for any $i\in I$ and $\ell\in \N$
		\begin{align}\label{eq:affine-Yangian-Lie-algebra-1} 
			\kappa_{i, \ell}&\in Z\Big(\eY_\qv\Big) \ ,
		\end{align}
		
		\item for any $i, j\in I$ and $r,s\in \N$
		\begin{align}\label{eq:affine-Yangian-Lie-algebra-2}
			\Big[h_{i,r}, h_{j,s}\Big]  & = 0\ ,\\[4pt]  \label{eq:affine-Yangian-Lie-algebra-3}
			\Big[x_{i,r}^{+}, x_{j,s}^{-}\Big] &= \delta_{i,j} h_{i, r+s} \ , \\[4pt] \label{eq:affine-Yangian-4}
			\Big[h_{i,0}, x_{j,r}^\pm\Big]&=\pm a_{i,j}x_{j,r}^\pm\ , \\[4pt]\label{eq:affine-Yangian-5}
			\Big[h_{i, r+1}, x_{j, s}^{\pm}\Big] - \Big[h_{i, r}, x_{j, s+1}^{\pm}\Big]& = \pm \frac{\hbar}{2} a_{i,j} \Big\{h_{i, r}, x_{j, s}^{\pm}\Big\}-m_{i,j}\frac{\varepsilon_1-\varepsilon_2}{2}\Big[h_{i,r},x_{j,s}^{\pm}\Big]\ ,\\[4pt] \label{eq:affine-Yangian-6}
			\Big[x_{i, r+1}^{\pm}, x_{j, s}^{\pm}\Big] - \Big[x_{i, r}^{\pm}, x_{j, s+1}^{\pm}\Big] &= \pm \frac{\hbar}{2} a_{i,j} \Big\{x_{i, r}^{\pm}, x_{j, s}^{\pm}\Big\}-m_{i,j}\frac{\varepsilon_1-\varepsilon_2}{2}\Big[x_{i,r}^{\pm},x_{j,s}^{\pm}\Big] \ ,
		\end{align}
		where
		\begin{align}
			m_{i,j}\coloneqq \begin{cases}
				1 & \text{if } i\to j\in \Omega\ , \\
				-1 & \text{if } j\to i \in \Omega \ , \\
				0 & \text{otherwise} \ .
			\end{cases}
		\end{align}
		
		\item Serre relations: 
		\begin{align}\label{eq:affine-Yangian-Serre}
			\sum_{\sigma \in \frakS_m}\Big[x_{i, r_{\sigma(1)}}^{\pm}, \Big[x_{i, r_{\sigma(2)}}^{\pm}, \Big[\cdots, \Big[x_{i, r_{\sigma(m)}}^{\pm}, x_{j,s}^{\pm}\Big]\cdots\Big]\Big]\Big] = 0 
		\end{align}
		for $i, j\in I$, with $i \neq j$, where $m\coloneqq 1 - a_{i,j}$ and $\frakS_m$ denotes the $m$-th symmetric group.
	\end{itemize}	
	
	The \textit{(two-parameters) Yangian $\Y_{\qv;\, \varepsilon_1, \varepsilon_2}$ of $\qv$} is the unital associative $\Hbullet_T$-algebra, which is the quotient of $\eY_{\qv; \varepsilon_1, \varepsilon_2}$ by the two-sided ideal generated by the $\kappa_{i, \ell}$'s, with $i \in I$ and $\ell \in \N$.
\end{definition}

\begin{remark}
	Assume that $\qv\neq A_1^{(1)}$. When $\varepsilon_1=\varepsilon_2$, i.e., when we restrict ourselves to the diagonal torus $\C^\ast \to \C^\ast \times \C^\ast$, the Yangian $\Y_{\qv;\, \varepsilon_1, \varepsilon_2}\vert_{\varepsilon_1=\varepsilon_2}$ is the usual affine Yangian which has appeared in the literature (cf.\ e.g. \cite{Varagnolo_Yangian,Guay07}). When $\qv=A_N^{(1)}$, with $N\geq 2$, the above definition recovers the definition of the corresponding two-parameters Yangian introduced by Kodera, see, e.g., \cite[Definition~2.1]{Kodera21}.
\end{remark}

We now turn to the special case $\qv=A_1^{(1)}$. We consider the affine Yangian introduced by Kodera in \cite[Definition~5.1]{Kodera_Fock19} and by Bershtein-Tsymbaliuk in \cite[\S1.2]{BT-Yangians}. Note that the relation \eqref{eq:affine-Yangian-sl(2)-3} was missing in Kodera's definition, while present in the one of Bershtein-Tsym\-ba\-liuk.

\begin{definition}\label{def:affine-Yangian-sl(2)}
	The \textit{extended (two-parameters) Yangian $\eY_{\qv;\, \varepsilon_1, \varepsilon_2}$ of $\qv=A_1^{(1)}$} is the unital associative $\Hbullet_T$-algebra generated by $x_{i, \ell}^\pm, h_{i, \ell}, \kappa_{i, \ell}$, with $i \in I$ and $\ell \in \N$, subject to the relations \eqref{eq:affine-Yangian-Lie-algebra-1}, \eqref{eq:affine-Yangian-Lie-algebra-2}, \eqref{eq:affine-Yangian-Lie-algebra-3}, and \eqref{eq:affine-Yangian-4} for any $i,j$, the relation \eqref{eq:affine-Yangian-5} for $i=j$, the Serre relations \eqref{eq:affine-Yangian-Serre}, and 
	\begin{align}\label{eq:affine-Yangian-sl(2)-1}
		\Big[h_{i,r+2}, x_{j,s}^{\pm}\Big] - 2\Big[h_{i,r+1}, x_{j,s+1}^{\pm}\Big] + \Big[h_{i,r}, x_{j,s+2}^{\pm}\Big] =&
		\mp(\varepsilon_1+\varepsilon_2)\Big\{h_{i, r+1}, x_{j,s}^\pm\Big\}\pm (\varepsilon_1+\varepsilon_2)\Big\{h_{i,r},x_{j, s+1}^\pm\Big\}\\[4pt]	 
		&\mp \varepsilon_1\varepsilon_2 \Big[h_{i, r}, x_{j,s}^\pm\Big] \ ,\\[6pt] \label{eq:affine-Yangian-sl(2)-2}
		\Big[x_{i,r+2}^\pm, x_{j,s}^{\pm}\Big] - 2\Big[x_{i,r+1}^\pm, x_{j,s+1}^{\pm}\Big] + \Big[x_{i,r}^\pm, x_{j,s+2}^{\pm}\Big] =& \mp(\varepsilon_1+\varepsilon_2)\Big\{x_{i, r+1}^\pm, x_{j,s}^\pm\Big\} \pm (\varepsilon_1+\varepsilon_2)\Big\{x_{i,r}^\pm, x_{j, s+1}^\pm\Big\}\\[4pt] 
		&\mp \varepsilon_1\varepsilon_2 \Big[x_{i, r}^\pm, x_{j,s}^\pm\Big]  \ , \\[6pt]
		\label{eq:affine-Yangian-sl(2)-3}
		\Big[h_{i,1},x_{j,r}^\pm\Big]=&\mp 2x_{j,r+1}^\pm\mp\hbar\Big\{h_{i,0},x_{j,r}^\pm\Big\}\ ,
	\end{align}
	for $i \neq j\in I$ and $r,s\in \N$.
	
	The \textit{(two-parameters) Yangian $\Y_{\qv;\, \varepsilon_1, \varepsilon_2}$ of $\qv=A_1^{(1)}$} is the unital associative $\Hbullet_T$-algebra, which is the quotient of $\eY_{\qv; \varepsilon_1, \varepsilon_2}$ by the two-sided ideal generated by the $\kappa_{i, \ell}$'s, with $i \in I$ and $\ell \in \N$.
\end{definition}

The subalgebras $\Y^{\pm}_{\qv;\varepsilon_1,\varepsilon_2}, \zeroeY_{\qv;\varepsilon_1,\varepsilon_2},\Y^{\leqslant 0}_{\qv;\varepsilon_1,\varepsilon_2}, \ldots$ are defined in the usual fashion. The affine Yangian $\Y_{\qv;\varepsilon_1,\varepsilon_2}$ carries two gradings: the \textit{horizontal grading} by the group $\Z I$ and the \textit{vertical}, or \textit{cohomological grading} by the monoid $\N$. The respective degrees of generators are as follows:
\begin{align}
	\text{deg}(x^{\pm}_{i,\ell})=(\pm \alpha_i, \ell)\ , \quad \text{deg}(h_{i,\ell})=\text{deg}(\kappa_{i,\ell})=(0, \ell)\ , \quad \text{deg}(\varepsilon_i)=(0,1)\ .
\end{align}
We shall sometimes denote by $(\Y_{\qv;\varepsilon_1,\varepsilon_2})_{\bfd}$ the graded piece of horizontal degree $\bfd \in \Z I$.

\begin{theorem}\label{thm:triangular-decomposition-affine} 
	The multiplication map yields the following isomorphisms
	\begin{align}\label{eq:triangular-affine}
		\begin{tikzcd}[ampersand replacement=\&]
			\Y_{\qv;\, \varepsilon_1, \varepsilon_2}^+ \otimes \zeroeY_{\qv;\, \varepsilon_1, \varepsilon_2} \otimes \Y_{\qv;\, \varepsilon_1, \varepsilon_2}^- \ar{r}{\sim}\& \eY_{\qv;\, \varepsilon_1, \varepsilon_2}
		\end{tikzcd}\ ,
	\end{align}	
	and
	\begin{align}\label{eq:half-triangular-affine} 
		\begin{tikzcd}[ampersand replacement=\&]
			\Y_{\qv;\, \varepsilon_1, \varepsilon_2}^+\otimes\zeroeY_{\qv;\, \varepsilon_1, \varepsilon_2}  \ar{r}{\sim}\& \Y_{\qv;\, \varepsilon_1, \varepsilon_2}^{\frake,\geqslant 0} 
		\end{tikzcd}
		 \quad \text{and}\quad
		 \begin{tikzcd}[ampersand replacement=\&]
		 		\zeroeY_{\qv;\, \varepsilon_1, \varepsilon_2}\otimes \Y_{\qv;\, \varepsilon_1, \varepsilon_2}^-  \ar{r}{\sim}\& \Y_{\qv;\, \varepsilon_1, \varepsilon_2}^{\frake,\leqslant 0}
		 \end{tikzcd}\ .
	\end{align}
\end{theorem}
\begin{proof}
	The arguments appearing in the proof of Theorem~\ref{thm:triangular-decomposition} apply here as well. We leave the details to the reader.
\end{proof}

\subsection{Elliptic Lie algebras}\label{subsec:elliptic-Lie-algebras}

We next describe the Lie algebras which arise as classical limits of the above Yangians. 

\begin{definition}\label{def:classical-Lie-algebra}
	Let $\qv$ be an affine ADE quiver. Let $\fraks_\qv$ be the Lie algebra generated by $X_{i, \ell}^\pm$ and $H_{i, \ell}$, with $i\in I$ and $\ell\in \N$, subject to the following defining relations:
	\begin{itemize}\itemsep0.2cm
		\item if $\qv\neq A_1^{(1)}$:
		\begin{align}\label{eq:classical-1}
			\Big[H_{i,r}, H_{j,s}\Big] &= 0\ ,\\[4pt] \label{eq:classical-2}
			\Big[X_{i,r}^{+}, X_{j,s}^{-}\Big] &= \delta_{i,j} H_{i, r+s}\ ,\\[4pt]\label{eq:classical-3}
			\Big[H_{i,r}, X_{j,s}^{\pm}\Big] &= \pm a_{i,j} X_{j,r+s}^{\pm}\ ,\\[4pt] \label{eq:classical-4}
			\Big[X_{i, r+1}^{\pm}, X_{j, s}^{\pm}\Big] &= \Big[X_{i, r}^{\pm}, X_{j, s+1}^{\pm}\Big]  \ ,
		\end{align}
		and
		\begin{align}\label{eq:classical-5}
			\ad\Big(X_{i, 0}^\pm\Big)^{1-a_{i,j}} \Big(X_{j, s}^\pm\Big)=0  \quad \text{for } i \neq j\ .
		\end{align}
		
		\item if $\qv=A_1^{(1)}$: the relation \eqref{eq:classical-4} for $i=j$, the relations \eqref{eq:classical-1}, \eqref{eq:classical-2}, \eqref{eq:classical-3}, \eqref{eq:classical-5}, as well as
		\begin{align} \label{eq:classical-A1-2}
			\Big[X_{i,r+2}^{\pm}, X_{j,s}^{\pm}\Big] - 2\Big[X_{i,r+1}^{\pm}, X_{j,s+1}^{\pm}\Big] + \Big[X_{i,r}^{\pm}, X_{j,s+2}^{\pm}\Big] & =0
		\end{align}
		for $i \neq j$.
	\end{itemize}
	Let $\fraks_\qv^0$ be the Lie subalgebra of $\fraks_\qv$ generated by the elements $H_{i, \ell}$ with $i\in I$ and $\ell\in \N$, and let $\fraks_\qv^\pm$ be the Lie subalgebra of $\fraks_\qv$ generated by the elements $X_{i, \ell}^\pm$ with $i\in I$ and $\ell\in \N$. 
\end{definition}

Note that $\fraks_\qv$ carries a $\Z I \times \Z$-grading with $\deg(X_{i, \ell}^\pm) = (\pm \alpha_i, -2\ell)$ and $\deg(H_{i, \ell})=(0,-2\ell)$ for any $i\in I$ and $\ell\in \N$. 

\begin{remark}
	The above definition coincides with \cite[Definition~2.5]{GRW19} when $\qv\neq A_1^{(1)}$. The case of affine type A quivers (including $A_1^{(1)}$) is also considered in \cite{Tsymbaliuk_Yangians}.
\end{remark}

Let $\Y\vert_0$ be the specialization $\Y_{\qv;\, \varepsilon_1, \varepsilon_2}/(\varepsilon_1, \varepsilon_2)\Y_{\qv;\, \varepsilon_1, \varepsilon_2}$. Define $\Y^-\vert_0$ similarly. Our main motivation to consider elliptic Lie algebras is the following result.

\begin{proposition}\label{prop:classical-limit}
	The assignment, for $i\in I$ and $\ell\in \N$, 
	\begin{align}\label{eq:assignment-classical}
		X_{i, \ell}^\pm \longmapsto x_{i, \ell}^\pm \quad \text{and}\quad H_{i, \ell}\longmapsto h_{i, \ell}\ ,
	\end{align}
	extends to an isomorphism of $\Z I \times \N$-graded $\Q$-algebras 
	\begin{align}\label{eq:homomorphism-classical}
		\begin{tikzcd}[ampersand replacement=\&]
			\sfU(\fraks_\qv)\ar{r}{\sim} \& \Y\vert_0
		\end{tikzcd}\ ,
	\end{align}
	and restricts to
	\begin{align}
		\begin{tikzcd}[ampersand replacement=\&]
			\sfU(\fraks_\qv^-)\ar{r}{\sim}\& \Y^-\vert_0
		\end{tikzcd}\ .
	\end{align}
\end{proposition}	

\begin{proof}
	The case $\qv\neq A_1^{(1)}$ is proved in \cite[Proposition~2.6]{GRW19}. We shall follow the same arguments as in \textit{loc.cit.} to prove the remaining case $\qv=A_1^{(1)}$. 
	
	First, note that the defining relations of $\Y\vert_0$ are of Lie type, hence $\Y\vert_0$ is isomorphic to the universal enveloping algebra of a Lie algebra $\fraks_\qv'$, where $\fraks_\qv'$ is the Lie algebra generated by $x_{i, \ell}^\pm$ and $h_{i, \ell}$, with $i\in I$ and $\ell\in \N$, subject to the relation \eqref{eq:classical-4} for $i=j$, the relations \eqref{eq:classical-1}, \eqref{eq:classical-2}, \eqref{eq:classical-A1-2}, 
	\begin{align}\label{eq:Yangian-classical-limit-1}
		\Big[h_{i,0}, x_{j,r}^{\pm}\Big] &= \pm a_{i,j} x_{j,r}^{\pm}\ , \\[2pt] \label{eq:Yangian-classical-limit-2}
		\Big[h_{i,1}, x_{j,r}^{\pm}\Big] &=\mp 2 x_{j, r+1}^\pm\ , \qquad\text{for } i \neq j \ ,\\[2pt] \label{eq:Yangian-classical-limit-3}
		\Big[h_{i, r+1}, x_{i, s}^{\pm}\Big] &= \Big[h_{i, r}, x_{i, s+1}^{\pm}\Big]\ , \\[2pt]\label{eq:Yangian-classical-limit-4}
		\Big[h_{i,r+2}, x_{j,s}^{\pm}\Big] &- 2\Big[h_{i,r+1}, x_{j,s+1}^{\pm}\Big] + \Big[h_{i,r}, x_{j,s+2}^{\pm}\Big]= 0\ , \qquad\text{for } i \neq j\ ,
	\end{align}
	and
	\begin{align}\label{eq:Serre-classical}
		\sum_{w \in \frakS_3}\Big[x_{i,r_{w(1)}}^{\pm}, \Big[x_{i,r_{w(2)}}^{\pm}, \Big[x_{i,r_{w(3)}}^{\pm}, x_{j,s}^{\pm}\Big]\Big]\Big] = 0\ , \qquad\text{for } i \neq j\ .
	\end{align}
	Note that if $i=j$, the relations \eqref{eq:Yangian-classical-limit-1} and  \eqref{eq:Yangian-classical-limit-3} are equivalent to \eqref{eq:classical-3}. On the other hand, when $i\neq j$, the relations \eqref{eq:Yangian-classical-limit-1}, \eqref{eq:Yangian-classical-limit-2}, and \eqref{eq:Yangian-classical-limit-4} are equivalent to \eqref{eq:classical-3}. Thus, the assignment \eqref{eq:assignment-classical} extends to an epimorphism of algebras $\sfU(\fraks_\qv) \longrightarrow \Y\vert_0$.
	
	To conclude, we need to show that the inverse assignment $x_{i, \ell}^\pm\mapsto X_{i, \ell}^\pm$ and $h_{i, \ell}\mapsto H_{i, \ell}$ extends to a homomorphism $\Y\vert_0 \longrightarrow \sfU(\fraks_\qv)$. We have to check that the relations of Definition~\ref{def:classical-Lie-algebra} imply the relations \eqref{eq:Yangian-classical-limit-1} -- \eqref{eq:Yangian-classical-limit-4}  and \eqref{eq:Serre-classical}. The former has been already proved in the above paragraph. We are left to deduce the Serre relation \eqref{eq:Serre-classical} from Definition~\ref{def:classical-Lie-algebra}. This follows by successively applying $\mathsf{ad}(H_{i,n})$ to \eqref{eq:classical-5}.
		\end{proof}

It is also interesting to consider presentations of the positive halves $\fraks_\qv^\pm$. The proof of the following result will be given in Appendix~\ref{sec:characterization-classical-limit}.
\begin{proposition}\label{prop:characterization-classical-limit}
	\hfill
	\begin{itemize}\itemsep0.2cm
		\item $\fraks_\qv^0$ is isomorphic to the free commutative Lie algebra generated by $H_{i, \ell}$, with $i\in I$ and $\ell\in \N$.
		\item $\fraks_\qv^\pm$ is isomorphic to the Lie algebra generated by the elements $X_{i, \ell}^\pm$ with $i\in I$ and $\ell\in \N$ satisfying the relations:
		\begin{enumerate}\itemsep0.2cm
			\item \label{item:tilde-relations-1} \eqref{eq:classical-4} and \eqref{eq:classical-5}, if $\qv\neq A_1^{(1)}$;
			\item \label{item:tilde-relations-2} \eqref{eq:classical-4} for $i=j$, \eqref{eq:classical-5}, \eqref{eq:classical-A1-2}, and
			\begin{align}\label{eq:classical-A1-3}
				\Big[\Big[X_{i, r+1}^{\pm}, X_{j, s}^{\pm}\Big], X_{j, \ell}^\pm\Big] - \Big[\Big[X_{i, r}^{\pm}, X_{j, s+1}^{\pm}\Big], X_{j, \ell}^\pm\Big] =0\ , \qquad\text{for} i \neq j 
			\end{align}
			if $\qv=A_1^{(1)}$.
		\end{enumerate}
	\end{itemize}
\end{proposition}

\subsection{Elliptic Lie algebras as double loop algebras}

We now give another description of $\fraks_\qv$ as the universal central extension of a double loop algebra. Consider the Lie algebra $\frakgfin[s^{\pm 1}, t]\coloneqq \frakgfin\otimes_\Q \Q[s^{\pm 1}, t]$. By \cite[Theorem~3.3]{Kassel84}, it admits a universal central extension defined by
\begin{align}
	\frakgell \coloneqq \uce\big(\frakgfin[s^{\pm 1}, t]\big)\ ,
\end{align}
where
\begin{align}
	\uce\big(\frakgfin[s^{\pm 1}, t]\big)\coloneqq \frakgfin[s^{\pm 1}, t]\oplus \big( \Omega^1\big(\Q[s^{\pm 1}, t]\big)/d\big(\Q[s^{\pm 1}, t]\big) \big)
\end{align}
as a vector space, with Lie bracket such that $\Omega^1\big(\Q[s^{\pm 1}, t]\big)/d\big(\Q[s^{\pm 1}, t]\big)$ is central and
\begin{align}
	[x\otimes a, y\otimes b]=[x, y]\otimes ab+(x, y)\cdot b(da)
\end{align}
for $x, y\in \frakgfin$ and $a,b \in \Q[s^{\pm 1}, t]$. Here, $(-,-)$ is an invariant bilinear form on $\frakgfin$. For a commutative $\Q$-algebra $A$ we denote by $\Omega^1(A)$ is the module of Kähler differentials of $A$ and by $d(A)$ the subspace of exact forms.

Note that $\Omega^1\big(\Q[s^{\pm 1}, t]\big)/d\big(\Q[s^{\pm 1}, t]\big)$ has a basis given by
\begin{align}
	\big\{s^{-1}ds, t^\ell s^k ds\, \vert\, \ell\in \Z_{\geqslant 1}, \ k\in \Z\big\} \ .
\end{align}
Set $c_\ell\coloneqq t^\ell s^{-1} ds$ for $\ell\in \N$ and $c_{k, \ell}\coloneqq t^\ell s^{k-1} ds$ for $\ell\in \N$, with $\ell\neq 0$, and $k\in \Z$, $k\neq 0$. Then $\frakgell$ has the following equivalent description:
\begin{align}
	\frakgell = \frakgfin[s^{\pm 1}, t]\oplus K\quad \text{with } K\coloneqq \bigoplus_{\ell\in \N}\Q c_\ell \oplus \bigoplus_{\genfrac{}{}{0pt}{}{\ell\in \N, \; \ell\geqslant 1}{k\in \Z, \; k\neq 0}} \Q c_{k, \ell} \ ,
\end{align} 
with $c_{\ell}, c_{k, \ell}$ central elements, and Lie bracket
\begin{align}
	[x\otimes s^k t^\ell, y\otimes s^h t^n]=
	\begin{cases}
		[x, y]\otimes t^{\ell+n}+ k(x, y)\cdot c_{\ell+n}  & \text{if } k+h=0\ ,\\[4pt]
		[x, y]\otimes s^{k+h} t^{\ell+n}+ (kh-\ell n) \cdot (x, y)\cdot c_{m+n, g+k} &\text{if } k+h\neq 0\ .
	\end{cases}
\end{align}
We equip the Lie algebra $\frakgell$ with the $\Z\times\Z I$-grading such that 
\begin{align}
	\deg\big(x\otimes s^k t^\ell\big)\coloneqq (-2\ell, \bfd+k\delta) \ , \quad \deg\big(c_{k, \ell}\big)\coloneqq (-2\ell, k\delta) \ , \quad \deg\big(c_{\ell}\big)\coloneqq (-2\ell, 0)\ ,
\end{align}
where $x$ belongs to the root space $(\frakgfin)_{\bfd}, k\in \Z$, and $\ell\in \N$. We'll call the first term of the grading the \textit{horizontal grading} and the second term the \textit{vertical grading}. 

The \textit{negative half} $\fraknell$ of $\frakgell$ is the Lie subalgebra spanned by the homogeneous elements whose horizontal grading belongs to $-\N I\smallsetminus\{0\}$. We have
\begin{align}\label{eq:negative-half}
	\fraknell= \frakn[t] \oplus K_- \quad\text{where } K_-\coloneqq \bigoplus _{k<0} \Q c_{k, \ell}\ ,
\end{align}
where $\frakn=s^{-1}\frakgfin[s^{-1}] \oplus \fraknfin$ is the (standard) \textit{negative nilpotent half} of $\frakg$.

Let $(\fraknell)_{\ell,\bfd}$ and $\sfU(\fraknell)_{\ell,\bfd}$ be the degree $(\ell,\bfd)$ pieces in the Lie algebra $\fraknell$ and its enveloping algebra $\sfU(\fraknell)$. 
\begin{lemma}\label{lem:dimension-Ln}
	One has
	\begin{align}
		\sum_{\ell\in \Z, \bfd\in \N I}\dim \sfU(\fraknell)_{\ell,-\bfd}\,t^\ell z^\bfd=\Exp\Bigg(\sum_{\bfd\in\N I} \frac{A_\bfd(t^{-2})}{1-t^{-2}}\,z^\bfd \Bigg)\ ,
	\end{align}
	where $A_{\bfd}(z)$ is the Kac's polynomial of the quiver $\qv$ with dimension vector $\bfd$ and $z^\bfd\coloneqq \prod_{i\in I} z_i^{d_i}$.
\end{lemma}
\begin{proof}
	This follows from the relation
	\begin{align}\label{eq:dim-Lie-Kac-pols}
		\sum_{\ell\in \Z, \ \bfd\in \N I} \dim (\fraknell)_{\ell,-\bfd}\,t^\ell z^\bfd =\sum_{\bfd\in\N I} \frac{A_\bfd(t^{-2})}{1-t^{-2}}\,z^\bfd\ ,
	\end{align}
	which is itself a consequence of Formula~\eqref{eq:negative-half} together with the equalities
	\begin{align}
		A_{\bfd}(z)=
		\begin{cases}
			1 & \text{if}\; \bfd \in \Delta^{\mathsf{re}}\ ,\\
			e+z & \text{if}\; \bfd \in \Delta^{\mathsf{im}}\ ,\\
			0 & \text{otherwise}\ .
		\end{cases}
	\end{align} 
	Recall that $e$ is the rank of $\frakgfin$.
\end{proof}

Let us denote by $X_i^\pm$ and $H_i$, with $i=1, \ldots, e$, the Chevalley generators for $\frakgfin$ normalized so that $(X_i^+,X_i^-) = 1$ and $H_i= [X_i^+,X_i^-]$. Let $X_{\pm \varphi}$ be root vectors of $\frakgfin$ for the roots $\pm\varphi$ normalized so that $(X_\varphi,X_{\textrm{-}\varphi}) = 1$, where $\varphi$ is the highest root of $\frakgfin$. Set $H_\varphi\coloneqq [X_\varphi^+,X_\varphi^-]$.

\begin{proposition}\label{prop:iso-s-uce-1}
	The assignment
	\begin{align}
		X_{i, \ell}^\pm &\longmapsto X_i^\pm \otimes t^\ell \quad \text{for } i=1, \ldots, e\;\text{and}\; \ell\in \N\ , \\
		X_{0, \ell}^\pm &\longmapsto X_{\mp\varphi}\otimes t^\ell s^{\pm 1}\quad \text{for } \ell\in \N\ , \\
		H_{i, \ell} &\longmapsto H_i\otimes t^\ell \quad \text{for } i=1, \ldots, e\;\text{and}\; \ell\in \N\ , \\
		H_{0, \ell} &\longmapsto H_\varphi\otimes t^\ell+t^\ell s^{-1} ds \quad \ell\in \N\ .
	\end{align}
	defines a Lie algebras isomorphism $\phi\colon \fraks_\qv \to \frakgell$. In particular, it restricts to an isomorphism of Lie algebras $\phi\colon \fraks_\qv^-\to \fraknell$.
\end{proposition}
\begin{proof}
	This may be shown using the same arguments as in \cite[\S3]{MRY90}, (see in particular Proposition 3.5 of \textit{loc.cit.}). We leave the details to the reader.
\end{proof}

\subsection{Standard filtration and relation with COHAs}\label{subsec:PBW-theorem}

The $\Hbullet_T$-algebras $\Y_{\qv;\, \varepsilon_1, \varepsilon_2}$ and $\Y_{\qv;\, \varepsilon_1, \varepsilon_2}^\pm$ admit increasing $\N$-filtrations, called the \textit{standard} filtrations, for which the generators $x^\pm_{i,\ell}$ and $h_{i,\ell}$ have degree $\ell$ and $\deg(\varepsilon_i)=0$ for $i=1,2$. Let $\gr\Y_{\qv;\, \varepsilon_1, \varepsilon_2}$ and $\gr\Y_{\qv;\, \varepsilon_1, \varepsilon_2}^\pm$ be the associated graded. They are $\Z I\times\Z$-graded in the obvious way. 
Let $\overline{x}_{i, \ell}^\pm$ and $\overline{h}_{i, \ell}$ denote the images of $x_{i, \ell}^\pm$ and $h_{i, \ell}^\pm$ in $\gr\Y_{\qv;\, \varepsilon_1, \varepsilon_2}$.

We may see $T$ as a subtorus of $\Tmax$ (introduced in Formula~\eqref{eq:Tmax}) via the embedding 
\begin{align}
	(\gamma_1, \gamma_2)\in T \longmapsto (\gamma_e=\gamma_1, \gamma_{e^\ast}=\gamma_2)_{e\in \Omega}\in \Tmax\ .
\end{align}
Then, $	\varepsilon_1=\varepsilon_e$ and $\varepsilon_2=\varepsilon_{e^\ast}$ for $e\in \Omega$. This yields a specialization morphism $\Hbullet_{\Tmax} \to \Hbullet_T$.

Let $\Y_\qv$ be the Yangian of $\qv$ introduced in Definition~\ref{def:Yangian} (i.e., the one with the extra cubic relations). Recall that $\Y_\qv$ and $\Y_\qv^\pm$ are equipped with standard filtrations defined in the same way (see \S\ref{sec:yangiangradingsandtruncation}). Let $\cohaqv^T$ be the nilpotent COHA as defined in \S\ref{sec:defCOHAnilp}, with respect to the torus $T$.

\begin{theorem}\label{thm:PBW}
	There is a $\Z I \times\Z $-graded $\Hbullet_T$-algebra isomorphism 
	\begin{align}
		\begin{tikzcd}[ampersand replacement=\&]
			\Psi\colon \sfU(\fraks_\qv^-)\otimes_\Q \Hbullet_T\simeq \sfU(\fraknell)\otimes_\Q \Hbullet_T \ar{r}{\sim}\& \gr \Y_{\qv;\, \varepsilon_1, \varepsilon_2}^-
		\end{tikzcd}
	\end{align}	
	given by the assignment
	\begin{align}
		X_{i, \ell}^- \longmapsto \overline{x}_{i, \ell}^- 
	\end{align}
	for $i\in I$ and $\ell\in \N$.
	
	Moreover, we have a chain of isomorphisms:
	\begin{align}
		\begin{tikzcd}[ampersand replacement=\&]
			\Y_{\qv;\, \varepsilon_1, \varepsilon_2}^- \ar{r}{\sim}\& \Y_\qv^-\otimes_{\Hbullet_{\Tmax}} \Hbullet_T \ar{r}{\sim} \& \cohaqv^T
		\end{tikzcd}\ ,
	\end{align}
	where the second map is the homomorphism $\Phi$ introduced in Theorem~\ref{thm:Phi}.
\end{theorem}

\begin{proof}
	Thanks to the defining relations of $\Y_{\qv;\, \varepsilon_1, \varepsilon_2}$ and Proposition~\ref{prop:characterization-classical-limit}, there exists a surjective homomorphism
	\begin{align}
		\Psi\colon \sfU(\fraks_\qv^-)\otimes_\Q \Hbullet_T\simeq \sfU(\fraknell)\otimes_\Q \Hbullet_T \longrightarrow \gr \Y_{\qv;\, \varepsilon_1, \varepsilon_2}^-\ .
	\end{align}
	To prove that $\sfU(\fraks_\qv^-)\otimes_\Q \Hbullet_T$ and $\gr \Y_{\qv;\, \varepsilon_1, \varepsilon_2}^-$ are isomorphic, it is enough to compare their graded dimensions. 
	
	Now, Remark~\ref{rem:specialized-relations} yields the existence of a surjective homomorphism
	\begin{align}
		\Y_{\qv;\, \varepsilon_1, \varepsilon_2}^- \longrightarrow \Y_\qv^-\ .
	\end{align}
	Because of Theorem~\ref{thm:Phi}, we have also a surjective homomorphism 
	\begin{align}
		\Phi\colon \Y_\qv^- \longrightarrow\cohaqv^T\ .
	\end{align}
	Therefore, we have the following chain of inequalities, in any fixed degree $(\bfd, \ell)$:
	\begin{align}
		\mathsf{grdim}\big( \sfU(\fraks_\qv^-)\otimes_\Q \Hbullet_T \big) \geq \mathsf{grdim}\big( \gr \Y_{\qv;\, \varepsilon_1, \varepsilon_2}^- \big) \geq \mathsf{grdim} \big( \Y_\qv^- \big) \geq \mathsf{grdim} \big( \cohaqv^T \big) \ .
	\end{align}
	Now, the graded dimension of $\sfU(\fraks_\qv^-)\otimes_\Q \Hbullet_T$ is computed in Lemma~\ref{lem:dimension-Ln} and it coincides with the graded dimension of $\cohaqv^T$ by \cite[Theorem~A-(d)]{SV_generators}. Hence
	\begin{align}
		\mathsf{grdim}\big( \sfU(\fraks_\qv^-)\otimes_\Q \Hbullet_T \big) = \mathsf{grdim}\big( \gr \Y_{\qv;\, \varepsilon_1, \varepsilon_2}^- \big) = \mathsf{grdim} \big( \Y_\qv^- \big) = \mathsf{grdim} \big( \cohaqv^T \big) \ .
	\end{align}
	From this we obtain all the claims in the Theorem.
\end{proof}	

\begin{corollary}\label{cor:PBW}
	There is a $\Z I\times\Z $-graded $\Hbullet_T$-algebra isomorphism 
	\begin{align}
		\begin{tikzcd}[ampersand replacement=\&]
			\Psi\colon \sfU(\fraks_\qv)\otimes_\Q \Hbullet_T\simeq \sfU(\frakgell)\otimes_\Q \Hbullet_T \ar{r}{\sim}\& \gr \Y_{\qv;\, \varepsilon_1, \varepsilon_2}
		\end{tikzcd}
	\end{align}	
	given by the assignment
	\begin{align}
		X_{i, \ell}^\pm \longmapsto \overline{x}_{i, \ell}^\pm \quad \text{and} \quad H_{i, \ell} \longmapsto \overline{h}_{i, \ell},
	\end{align}
	for $i\in I$ and $\ell\in \N$. 
	
	Likewise, there is an isomorphism of  $\Z I\times\Z $-graded $\Hbullet_T$-algebras
	\begin{align}
		\Y_{\qv;\, \varepsilon_1, \varepsilon_2} \simeq \Y_\qv \otimes \Hbullet_T\ .
	\end{align}
\end{corollary}

\begin{proof}
	The first statement follows from the triangular decomposition of the Yangian (cf.\ Theorem~\ref{thm:triangular-decomposition-affine}) and Theorem~\ref{thm:PBW}. The proof of the second claim is similar, based on Theorems~\ref{thm:triangular-decomposition}.
\end{proof}

\begin{remark}
	When $\qv \neq A_1^{(1)}$, the first statement of Corollary~\ref{cor:PBW} appears as \cite[Theorem~6.9]{GRW19}.
\end{remark}

From now on and until \S\ref{sec:limit-affine-Yangian}, we will only work with equivariant parameters belonging to $\Hbullet_T$ and hence simply denote $\Y_{\qv; \varepsilon_1, \varepsilon_2}$ by $\Y_\qv$.

We conclude this Section with the following refinement of Proposition~\ref{prop:overlineTi}. The standard filtration of $\Y_\qv$ induces filtrations on the quotients $\iY_\qv$ and $\Yi_\qv$. The map $\overline T_i$ preserve these filtrations and descends to a map $\gr(\overline{T}_i)$ between the associated graded.

\begin{corollary}\label{cor:PBW-consequences}
	\hfill
	\begin{enumerate}\itemsep0.2cm
		\item\label{item:PBW-consequences-1} The map 
		\begin{align}\label{eq:mapYi}
			\Yi_\qv \longrightarrow  \Y_\qv /(\I^{\geqslant 0} \Y_\qv+ \Y_\qv\I_i^- )\ ,
		\end{align}
		induced by the canonical embedding $\Y_\qv^-\subset \Y_\qv$, is an isomorphism.  
		
		\item\label{item:PBW-consequences-2} The map $\overline T_i\colon \iY_\qv \to \Yi_\qv$, introduced in Proposition~\ref{prop:overlineTi}, is an isomorphism.
		
		\item\label{item:PBW-consequences-3} There are isomorphisms
		\begin{align}
			\gr(\iY_\qv)&\simeq  \big( \sfU(\fraknell)/\sum_{\ell \geq 0} x^-_it^\ell \sfU(\fraknell) \big)\otimes_\Q \Hbullet_T\ ,\\
			\gr(\Yi_\qv)&\simeq  \big( \sfU(\fraknell)/\sum_{\ell \geq 0} \sfU(\fraknell) x^-_it^\ell \big)\otimes_\Q \Hbullet_T\ ,
		\end{align}
		and the map $\gr(\overline T_i)$ coincides with the (truncation of the) standard action of the braid group on the enveloping algebra $\sfU(\frakgell)$.
	\end{enumerate}	
\end{corollary}

\begin{proof} 
	We start by proving \eqref{item:PBW-consequences-1}. By Theorem~\ref{thm:triangular-decomposition-affine}, the composition $\Y_\qv^- \to \Y_\qv \to \Y_\qv/ \I^{\geqslant 0} \Y_\qv$ is an isomorphism. In addition $\I^{\geqslant 0} \Y_\qv + \Y_\qv\I_i^{-}=\I^{\geqslant 0} \Y_\qv + \Y^-_\qv\I_i^{-}$. The result follows by definition of $\Yi_\qv$. We turn to \eqref{item:PBW-consequences-2}. In light of Theorem~\ref{thm:compatibility-braid-reflection-functors}, it is enough to note that the map $\Phi\colon \Y^-_\qv \simeq \cohaqv^T$ descends to isomorphisms  $\iY_\qv \simeq \icohaqv$ and $\Yi_\qv \simeq \coha^{T,(i)}_\qv$.
	
	Finally, we show \eqref{item:PBW-consequences-3}. We have
	\begin{align}
		\gr(\iY_\qv)=\gr(\Y^-\qv)/\gr(\I^-_i \cdot \Y_\qv^-)\ .
	\end{align} 
	Note that, \textit{a priori}, we only have an inclusion
	\begin{align}
		\gr_j(\I^-_i \cdot \Y_\qv^-) \supseteq \bigoplus_{k+h=j} \gr_k(\I^-_i) \cdot \gr_h(\Y^-_\qv)\ ,
	\end{align}
	hence only a surjection  
	\begin{align}
		\big( \sfU(\fraknell)/\sum_{\ell \geq 0} x^-_it^\ell \sfU(\fraknell) \big)\otimes_\Q \Hbullet_T \longrightarrow \gr(\iY)\ .
	\end{align}
	To see that this map is an isomorphism, we compare the graded characters. By Proposition~\ref{prop:reflection-stacks}--\eqref{item:reflection-stacks-4}, we have
	\begin{align}
		\sum_{\bfd,\ell}\dim\,\gr(\iY)_{\textrm{-}\bfd,\ell}\, q^{\ell}z^{\bfd}=\frac{1}{(1-q^{-1})^2}\Exp\Big( \frac{1}{1-q^{-1}}\sum_{\bfd\in \N I\smallsetminus \N \alpha_i} A_{\bfd}(q^{-1}) z^\bfd \Big)\ .
	\end{align}
	The graded character of $\big( \sfU(\fraknell)/\sum_{\ell \geq 0} x^-_it^\ell \sfU(\fraknell) \big)\otimes_\Q \Hbullet_T$ is given by the same expression thanks to the PBW theorem. The same argument works for $\Yi_\qv$. The identification of $\gr(\overline{T}_i)$ follows by construction.
\end{proof}

\section{Quotients of the Yangian and of the COHA of an affine ADE quiver}\label{sec:quotients-affine-coha}

In this section, we  consider the structure of quotients of the COHA $\cohaqv^T$ associated to open substacks of $\dLambda_\qv$ determined by Harder-Narasimhan conditions, and describe the corresponding quotients of the affine Yangians $\Y_\qv^-$.

\subsection{Quotients of the nilpotent quiver COHA}\label{sec:quotients-nilpotent-quiver-COHA}

In this section we shall introduce quotients of the equivariant nilpotent quiver COHA arising from Harder-Narasimhan strata and study their algebraic structures. Our quiver $\qv$ is fixed to be of affine type, but the results are valid for any quiver (even one with edge-loops) and any torus $T \subseteq \Tmax$.
\subsubsection{Homology of Harder-Narasimhan strata}

\begin{definition}
	Let $\kappa\subset\Q$ be an interval. We denote by $\modPi^{\kappa}$ the full subcategory of $\modPi$ consisting of those finite-dimensional representations of $\Pi_\qv$ for which the $\Ltheta$-slopes of all the HN factors\footnote{See \S\ref{subsec:stabquiver} for the definition of \textit{HN factors}.} belong to $\kappa$.
\end{definition}
For any $\ell\in \Q$, we set $\modPi^{\ell}\coloneqq \modPi^{\{\ell\}}$, $\modPi^{>\ell}\coloneqq \modPi^{(\ell,+\infty)}$, and $\modPi^{\leqslant\ell}\coloneqq \modPi^{(-\infty,\ell]}$. The following is clear. 

\begin{lemma}\label{lem:torsion-pair-quiver}
	The pair $\big(\modPi^{>\ell},\modPi^{\leqslant\ell}\big)$ is a torsion pair on the category $\modPi$, which is open in the sense of Lieblich \cite[Definition~A.2]{AB_Moduli_Bridgeland}. 
\end{lemma}

Given an interval $\kappa\subset\Q$, let $\nilpPi^\kappa$ be the category of nilpotent finite-dimensional representations $M$ of $\Pi_\qv$ whose Harder-Narasimhan-factors all have $\Ltheta$-slope in $\kappa$ and let $\Lambda_\qv^\kappa$ be the corresponding moduli stack. 
\begin{lemma}
	The stack $\Lambda_\qv^\kappa$ is an open substack of $\Lambda_\qv$.
\end{lemma}
\begin{proof}
	This follows from Lemma~\ref{lem:torsion-pair-quiver} in the case of intervals of the form $[\ell,\infty)$, $(\ell,\infty)$, $(-\infty,\ell)$ or $(-\infty,\ell)$; the general case now follows from the formula 
	\begin{align}
		\Lambda_\qv^{\kappa_1 \cap \kappa_2}=\Lambda_\qv^{\kappa_1} \times_{\Lambda_\qv}\Lambda_\qv^{\kappa_2}\ . \tag*{\qedhere} 
	\end{align}
\end{proof}

\begin{notation}
	We denote by $\dLambda_\qv^\kappa$ the canonical derived enhancement of $\Lambda_\qv^\kappa$ inside $\dLambda_\qv$. For $\bfd \in \N I$, we let $\Lambda_\bfd^\kappa$ (resp.\ $\dLambda_\bfd^\kappa$) be the connected component of $\Lambda_\qv^\kappa$ (resp.\ $\dLambda_\qv^\kappa$) of representations of dimension $\bfd$. We abbreviate $\Lambda_\qv^\ell\coloneqq \Lambda_\qv^{\{\ell\}}$ and $\dLambda_\qv^\ell\coloneqq \dLambda_\qv^{\{\ell\}}$.
\end{notation}

\begin{definition}
	Given two disjoint intervals $\kappa, \kappa'$ of $\R$, we write $\kappa'<\kappa$ if for any $a\in \kappa'$ and for any $b\in \kappa$ we have $a<b$.
\end{definition}

Given disjoint intervals $\kappa_1,\ldots,\kappa_s$ with $\kappa_1<\kappa_2<\cdots<\kappa_s$, define the derived stack $\dLambda_\qv^{\kappa_1, \ldots, \kappa_s}$ as
\begin{align}
	\begin{tikzcd}[ampersand replacement = \&]
		\dLambda_\qv^{\kappa_1,\kappa_2} \arrow[swap]{d}{q^{\kappa_1, \kappa_2}} \arrow{r} \& \dLambdaext_\qv \arrow{d}{q} \\
		\dLambda_\qv^{\kappa_1} \times \dLambda_\qv^{\kappa_2}\arrow{r} \& \dLambda_\qv\times  \dLambda_\qv
	\end{tikzcd} 
\end{align}
for $s=2$, and for arbitrary $s\geq 2$ as
\begin{align}\label{eq:HN-strata-Lambda}
	\dLambda_\qv^{\kappa_1, \ldots, \kappa_s} \coloneqq \dLambda_\qv^{\kappa_1, \kappa_2 \sqcup \cdots \sqcup \kappa_s} \times_{\dLambda_\qv^{\kappa_2 \sqcup \cdots \sqcup \kappa_s}}  \dLambda_\qv^{\kappa_2, \ldots, \kappa_s}\ .
\end{align}
Denote by $\Lambda_\qv^{\kappa_1, \ldots, \kappa_s}$ the classical truncation of $\dLambda_\qv^{\kappa_1, \ldots, \kappa_s}$. Moreover, for $\bfd_1,\ldots,\bfd_s\in\N I$, let $\Lambda^{\kappa_1,\ldots,\kappa_s}_{\bfd_1,\ldots,\bfd_s}$ (resp.\ $\dLambda^{\kappa_1,\ldots,\kappa_s}_{\bfd_1,\ldots,\bfd_s}$) be the connected component of $\Lambda_\qv^{\kappa_1, \ldots, \kappa_s}$ (resp.\ $\dLambda_\qv^{\kappa_1, \ldots, \kappa_s}$) parametrizing representations $M$ having a (necessarily unique) filtration $0\eqqcolon M_{s+1} \neq M_s \subset M_{s-1} \subset \cdots \subset M_1=M$ for which $M_i/M_{i+1}$ has dimension vector $\bfd_i$ and the $\Ltheta$-slopes of all its HN-factors belong to $\kappa_i$ for each $i=1,\dots, s$. Note that $\Lambda^{\kappa_1,\ldots,\kappa_s}_{\bfd_1,\ldots,\bfd_s}$ is a locally closed substack of $\Lambda_{\sum_i\bfd_i}$.

If each $\kappa_i$ is a singleton, then $\Lambda^{\kappa_1,\ldots,\kappa_s}_{\bfd_1,\ldots,\bfd_s}$ is a Harder-Narasimhan strata. In general, $\Lambda^{\kappa_1,\ldots,\kappa_s}_{\bfd_1,\ldots,\bfd_s}$ is a finite union of Harder-Narasimhan strata. For $M \in \Lambda^{\kappa_1,\ldots,\kappa_s}_{\bfd_1,\ldots,\bfd_s}$ as above, we denote by $M^i\coloneqq M_i/M_{i+1}$ the $i$th factor. Finally, $\Lambda_{\bfd_1, \ldots, \bfd_s}$ denotes the locally closed substack of $\Lambda_{\sum \bfd_i}$ parametrizing objects of HN-type $(\bfd_1, \ldots, \bfd_s)$. We denote by $\dLambda_{\bfd_1, \ldots, \bfd_s}$ the corresponding derived enhancement.

\begin{lemma}\label{lem:HNstackbundle}
	Fix $\kappa_1, \ldots, \kappa_s$, resp. $\bfd_1, \ldots, \bfd_s$  as above.
	\begin{enumerate}\itemsep0.2cm
		\item \label{item:HNstackbundle-1} $\dLambda_\qv^{\kappa_1, \ldots, \kappa_s}$ is an iterated vector bundle stack over $\dLambda_\qv^{\kappa_1} \times \cdots \times \dLambda_\qv^{\kappa_s}$.
		
		\item \label{item:HNstackbundle-2} $\dLambda_{\bfd_1, \ldots, \bfd_s}$ is an iterated vector bundle stack over $\dLambda_{\bfd_1}^{\Ltheta\textrm{-}\mathsf{ss}} \times \cdots \times \dLambda_{\bfd_s}^{\Ltheta\textrm{-}\mathsf{ss}}$.
	\end{enumerate}
\end{lemma}

\begin{proof}
	To prove \eqref{item:HNstackbundle-1}, it suffices to deal with the case $s=2$ because of Formula~\eqref{eq:HN-strata-Lambda}. For $s=2$, the projection map $M \mapsto (M_1/M_2, M_2)$ induces an equivalence
	\begin{align}
		\dLambda_\qv^{\kappa_1,\kappa_2} \simeq \Spec\Sym\big( (\R \calHom(\frakM^1, \frakM_2)[1])^\vee\big)\ ,
	\end{align}
	where $\R \calHom(\frakM^1, \frakM_2)$ is viewed as a complex over $\dLambda_\qv^{\kappa_1} \times \dLambda_\qv^{\kappa_2}$, where $\frakM^1$ is the universal object of $\dLambda_\qv^{\kappa_1}$ while $\frakM_2$ is the universal object of $\dLambda_\qv^{\kappa_2}$ (see e.g. \cite[Proposition~3.6]{Porta_Sala_Hall}). It is enough to prove that $\R \calHom(\frakM^1, \frakM_2)[1]$ has tor-amplitude $[0, 1]$\footnote{In homological notation.}. Furthermore, to perform the computation of the tor-amplitude it is enough to consider the two families of representations $\frakM^1$ and $\frakM_2$ over a point, i.e., it is enough to show that $\Ext^2(M^1, M_2)=0$ for $M^1\in \nilpPi^{\kappa_1}$ and $M_2\in \nilpPi^{\kappa_2}$. The properties of the Harder-Narasimhan filtration, the conditions on $\kappa_1, \kappa_2$ ensure that $\Hom(M_2, M^1)$ vanishes. By the 2-Calabi-Yau property (cf.\ Theorem~\ref{thm:PiQis2CY}), this implies that $\Ext^2(M^1, M_2)=0$ as wanted.
	
	Statement~\eqref{item:HNstackbundle-2} is obtained by specializing \eqref{item:HNstackbundle-1} to the case of $\kappa_i$ being singletons and fixing the connected component of $\dLambda_\qv^{\kappa_i}$.
\end{proof}

\begin{proposition}\label{prop:pure}
	\hfill
	\begin{enumerate}\itemsep0.2cm
		\item \label{item:pure-1} For any interval $\kappa$, the mixed Hodge structure of the stack $\dLambda_\qv^\kappa$ is pure and the restriction morphism $\HBMbulletT(\dLambda_\qv) \to \HBMbulletT(\dLambda^\kappa_\qv)$ is surjective.
		
		\item \label{item:pure-2} The map 
		\begin{align}
			\dLambda^{\kappa_1,\ldots,\kappa_s}_{\bfd_1,\ldots,\bfd_s} \longrightarrow \prod_{i=1}^s \dLambda^{\kappa_i}_{\bfd_i}\ , \quad M\longmapsto (M^1,\ldots,M^s)\ 
		\end{align}
		is an iterated vector bundle stack. In particular, the mixed Hodge structure of the stack $\dLambda^{\kappa_1,\ldots,\kappa_s}_{\bfd_1,\ldots,\bfd_s}$ is pure.
	\end{enumerate}
\end{proposition}

\begin{proof}
	For a general pair $(\kappa,\bfd)$, the stack $\Lambda_\bfd^\kappa$ is a finite union of Harder-Narasimhan strata. By Lemma~\ref{lem:HNstackbundle}, each such Harder-Narasimhan stratum is an iterated vector bundle stack over a product of stacks of the form $\dLambda^{\ell_i}_{\bfd_i}$.  The mixed Hodge structure of $\dLambda_\bfe^\ell$ is pure for any $\bfe, \ell$ by \cite[Theorem~3.2--(b)]{SV_generators} or \cite[Theorem~6.4]{Davison_Integrality}. It follows that the same holds for any Harder-Narasimhan stratum, and hence also for $\dLambda^{\kappa}_\bfd$. Then, Statement~\eqref{item:pure-1} is a consequence of the purity of $\dLambda^\kappa_\qv$ and its closed complement $\dLambda_\qv \smallsetminus \dLambda_\qv^{\kappa}$, which is itself, for any fixed $\bfd$, a finite union of Harder-Narasimhan strata.
	
	Statement~\eqref{item:pure-2} is a consequence of Lemma~\ref{lem:HNstackbundle}.
\end{proof}

\begin{notation}
	For any geometric classical stack $\frakX$, which is a quotient stack, its Borel-Moore homology is finite-dimensional in each degree as $\Q$-vector space. Thus, we can define
	\begin{align}
		P(\frakX,t)\coloneqq \sum_i \dim(\sfH^T_i(\frakX))t^i\ .
	\end{align}
	For any interval $\kappa$, we also set
	\begin{align}
		P_{\Lambda^\kappa_\qv}(z,t)\coloneqq \sum_{\bfd\in \N I} P(\Lambda^\kappa_\bfd,t)\,t^{(\bfd,\bfd)}z^\bfd\ .
	\end{align}
	The factor $t^{(\bfd,\bfd)}$ is added to take into account the cohomological degree shift in the grading of $\cohaqvd^T$.
\end{notation}

Recall that $A_\bfd(t) \in \N[t]$ is the Kac polynomial counting absolutely indecomposable representations of the quiver $\qv$ of dimension vector $\bfd$ over finite fields, see, e.g., \cite{BSV_Nilpotent}. The following is a refinement of Theorem~\ref{thm:purityandKac}.
\begin{proposition}\label{prop:coha-kac}  
	The following hold:
	\begin{enumerate}\itemsep0.2cm
		\item \label{item:coha-kac-2} For any partition $\Q =\kappa_1 \sqcup \cdots \sqcup \kappa_s$ with $\kappa_i < \kappa_j$ for $i <j$, we have
		\begin{align}
			P_{\Lambda_\qv}(z,t)=\prod_{k=1}^s P_{\Lambda_\qv^{\kappa_k}}(z,t)\ .
		\end{align}
		\item \label{item:coha-kac-3} For any $\ell \in \Q$, we have
		\begin{align}\label{eq:Poincqrepolsemistablelambda}
			P_{\Lambda^\ell_\qv}(z,t)=(1-t^{-2})^{-2}\Exp\Bigg(\frac{1}{1-t^{-2}}\sum_{\genfrac{}{}{0pt}{}{\bfd\in\N I}{\mu_{\Ltheta}(\bfd)=\ell}} A_\bfd(t^{-2})\,z^\bfd \Bigg)\ .
		\end{align}
		\item \label{item:coha-kac-4} For any subset $\kappa \subset \Q$ (not necessarily an interval) we have
		\begin{align}\label{eq:Poincqrepolsemistablelambda-kappa}
			P_{\Lambda^\kappa_\qv}(z,t)=(1-t^{-2})^{-2}\Exp\Bigg(\frac{1}{1-t^{-2}}\sum_{\genfrac{}{}{0pt}{}{\bfd\in\N I}{\mu_{\Ltheta}(\bfd)\in \kappa}} A_\bfd(t^{-2})\,z^\bfd \Bigg)\ .
		\end{align}
	\end{enumerate}
\end{proposition}

\begin{proof}
	
	Let us begin with \eqref{item:coha-kac-2}. By Proposition~\ref{prop:pure} we have
	\begin{align}
		P_{\Lambda_\qv}(z,t)&=\sum_{\bfd_1, \ldots, \bfd_s} P(\Lambda^{\kappa_1, \ldots, \kappa_s}_{\bfd_1,\ldots, \bfd_s},t)\,t^{(\bfd,\bfd)}z^{\sum_i \bfd_i}=\sum_{\bfd_1, \ldots, \bfd_s}\sum_i \dim\,\sfH^T_i(\Lambda^{\kappa_1, \ldots, \kappa_s}_{\bfd_1,\ldots, \bfd_s})\,t^{i+(\bfd,\bfd)}z^{\sum_i \bfd_i}\\[4pt]
		&=\sum_{\bfd_1, \ldots, \bfd_s}\sum_i \dim\,\sfH^T_i\Big(\prod_j\Lambda^{\kappa_j}_{\bfd_j}\Big)\,t^{i+\sum_j (\bfd_j,\bfd_j)}z^{\sum_i \bfd_i}=\prod_\ell P_{\Lambda^{\kappa_\ell}_\qv}(z,t)\ .
	\end{align}
	In the above calculation, we have made use of the relations
	\begin{align}
		P(\Lambda_\bfd,t)=\sum_{\bfd_1 + \cdots + \bfd_s=\bfd}P(\Lambda^{\kappa_1, \ldots, \kappa_s}_{\bfd_1, \ldots, \bfd_s},t)\ ,
	\end{align}
	which comes from the purity of each strata $\Lambda^{\kappa_1, \ldots, \kappa_s}_{\bfd_1, \ldots, \bfd_s}$ and
	\begin{align}
		\dim \sfH^T_{i-(\bfd,\bfd)}(\Lambda^{\kappa_1, \ldots, \kappa_s}_{\bfd_1, \ldots, \bfd_s})=\dim \sfH^T_{i-\sum_j (\bfd_j,\bfd_j)}\Big(\prod_j\Lambda^{\kappa_j}_{\bfd_j}\Big)\ ,
	\end{align}
	which comes from the fact that the morphism $\Lambda^{\kappa_1, \ldots, \kappa_s}_{\bfd_1, \ldots, \bfd_s} \to \prod_j\Lambda^{\kappa_j}_{\bfd_j}$ is an iterated vector bundle stack of rank $-\sum_{u>v} (\bfd_u,\bfd_v)$ (see Lemma~\ref{lem:HNstackbundle} and \cite[Remark~\ref*{torsion-pairs-rem:Gysin_vector_bundle}]{DPS_Torsion-pairs}).
	
	\medskip
	
	Let's turn to \eqref{item:coha-kac-3}. By \eqref{item:coha-kac-2}, we have a factorization into an infinite product
	\begin{align}\label{eq:infiniteproductHN}
		P_{\Lambda_\qv}(z,t)=\prod_{\ell \in \Q} P_{\Lambda_\qv^{\ell}}(z,t)\ .
	\end{align}
	Note that for any fixed $\bfd$ only finitely many slopes $\mu_{\Ltheta}(\bfe)$ may arise  for $\bfe \leq \bfd$, hence the above infinite product is well-defined and convergent. We may use the Harder-Narasimhan recursion formula to uniquely determine $P_{\Lambda^{\ell}_\qv}(z,t)$ from Equation~\eqref{eq:infiniteproductHN}; in other words, there is a unique family of power series 
	\begin{align}
		P_{\Lambda^\ell_\qv}(z,t) =\sum_{\substack{k,\bfd\in\N I\\ \mu_{\Ltheta}(\bfd)=\ell}} \alpha_{k,\bfd}t^k\,z^\bfd\ ,
	\end{align} 
	with $\ell \in \Q$, for which Equation~\eqref{eq:infiniteproductHN} holds. It only remains to observe that the right-hand-side of Formula~\eqref{eq:Poincqrepolsemistablelambda} obviously satisfies this condition. Statement \eqref{item:coha-kac-4} is proved in the same way.
\end{proof}

Thanks to Proposition~\ref{prop:coha-kac}, we may relate the existence of semistable representations of the preprojective algebra to the value of the Kac polynomial $A_\bfd(t)$, which is nonzero precisely when $\bfd$ belongs to the root system of $\qv$. In particular, we have the following.  
\begin{corollary} 
	A $\Ltheta$-semistable representation of $\Pi_\qv$ of slope zero is necessarily of dimension $\bfd \in \N \delta$.
\end{corollary}	

\begin{proof} 
	Let $\bfd \in \N I$ be a root of $\qv$, which we may write as $\bfd=n\delta + \bfe$ where $\bfe$ either belongs to $\N \Delta_\sff$ or $-\N \Delta_\sff$. But then unless $\bfe=0$, we will have $\Ltheta(\bfe) \neq 0$. This means that $\bfd \in \N \delta$.
\end{proof}	

\begin{remark}
	\hfill
	\begin{itemize}\itemsep0.2cm
		\item Proposition~\ref{prop:coha-kac} is valid for an arbitrary quiver $\qv$ (even one including edge-loops, as long as one replaces the Kac polynomials by its nilpotent version, see \cite{BSV_Nilpotent}). 
		\item For any collection ordered disjoint intervals (or even subsets)  $\kappa_1, \ldots, \kappa_s$ we have
		\begin{align}
			P_{\Lambda^{\kappa_1, \ldots, \kappa_s}_\qv}(z,t)=\prod_{k=1}^s P_{\Lambda_\qv^{\kappa_k}}(z,t)\ . \tag*{\qedhere} 
		\end{align}
	\end{itemize}
\end{remark}

\subsubsection{Algebraic structures on quotients of the COHA}

We examine in this section the structure on the quotients 
\begin{align}
	\coha_{\Ltheta,\kappa}^T\coloneqq \HBMbulletT(\dLambda_\qv^\kappa)
\end{align}
induced by the multiplication in $\cohaqv^T$. This heavily depends on the type of interval $\kappa$.

\begin{proposition}\label{prop:COHA}
	\hfill
	\begin{enumerate}\itemsep0.2cm
		\item \label{item:COHA-1} For any $\ell \in \Q$, $\coha^T_{\Ltheta, \{\ell\}}$ is a $\Z\times\N I$-graded $\Hbullet_T$-algebra.
		
		\item \label{item:COHA-2} The canonical restriction 
		\begin{align}
			\bigoplus_{\genfrac{}{}{0pt}{}{\bfd\in \N I}{\mu_{\Ltheta}(\bfd)=\ell}} \HBMbulletT(\dLambda_\bfd) \longrightarrow \coha^T_{\Ltheta, \{\ell\}}
		\end{align}
		is an algebra homomorphism, where the source has an algebra structure induced by the one of $\cohaqv^T$.
		
		\item \label{item:COHA-4} For any interval $\kappa=[a,b]$, $\coha^T_{\Ltheta, \kappa}$ is a right $\Z\times\N I$-graded  $\coha^T_{\Ltheta, \{b\}}$-module and a left $\Z\times\N I$-graded $\coha^T_{\Ltheta, \{a\}}$-module.
		
		\item \label{item:COHA-5} For any interval $\kappa=(a,b]$, the restriction map $\cohaqv^T \to \coha^T_{\Ltheta,\kappa}$ is surjective and its kernel is equal to 
		\begin{align}
			\K_{\Ltheta, \kappa}\coloneqq \sum_{\mu_{\Ltheta}(\bfd)>b} \cohaqv^T \cdot \cohaqvd^T\, + \sum_{\mu_{\Ltheta}(\bfd)\leqslant a} \cohaqvd^T \cdot \cohaqv^T\ .
		\end{align}
		Similar results hold if we replace we replace $(a,b]$ by $[a,b), [a,b]$ or $(a,b)$, allowing $a=b$.
		
		\item \label{item:COHA-3} For $\kappa=(-\infty, b]$, $\coha_{\Ltheta, \kappa}^T$ is a left $\Z\times\N I$-graded $\cohaqv^T$-module. For $\kappa=[a,\infty)$, $\coha_{\Ltheta, \kappa}^T$ is a right $\Z\times\N I$-graded $\cohaqv^T$-module. Similar results hold for $\kappa=(-\infty,b), (a,\infty)$.
	\end{enumerate}
\end{proposition}

\begin{proof}
	The category of all $\Ltheta$-semistable representations of a fixed $\Ltheta$-slope $\nu$ is an abelian subcategory of the category $\modPi$ which is stable under extension. Statement~\eqref{item:COHA-1} follows. In addition, if
	\begin{align}
		0\longrightarrow x\longrightarrow y \longrightarrow z \longrightarrow 0
	\end{align} 
	is an exact sequence in $\modPi$ in which $x,y,z$ are all of slope $\nu$ and $y$ is semistable then both $x$ and $y$ have to be semistable as well (indeed, $x$ belongs to $\modPi^{\leq \nu}$ and has slope $\nu$, which forces it to be semistable, and likewise for $z$). We deduce that there is, for any $\bfd,\bfe$ of slope $\nu$, a diagram
	\begin{align}
		\begin{tikzcd}[ampersand replacement=\&]
			\dLambda^\nu_\bfd \times \dLambda^\nu_\bfe \ar{d}{j} \& \dLambda^{\mathsf{ext},\nu}_{\bfd,\bfe} \arrow[swap]{l}{q^\circ} \ar{r}{p^\circ} \ar{d}{j'}\&
			\dLambda^\nu_{\bfd+\bfe} \ar{d}{j''}\\
			\dLambda_\bfd \times \dLambda_\bfe  \& \dLambdaext_{\bfd,\bfe} \arrow[swap]{l}{q} \ar{r}{p}\&
			\dLambda_{\bfd+\bfe}
		\end{tikzcd}\ ,
	\end{align}
	in which all the maps $j,j',j''$ are open immersions, and all squares are cartesian. The fact that the restriction map is a morphism of algebra follows by base change.
	
	The two parts of \eqref{item:COHA-4} are similar, we only prove the first. Fix a short exact sequence in $\modPi$ of the form
	\begin{align}
		0\longrightarrow x\longrightarrow y \longrightarrow z \longrightarrow 0\ ,
	\end{align}
	with $y\in \modPi^{\kappa}$ and $x$ of slope $b$. Arguing as above, we deduce that $x$ is semistable. The map $x\to y$ thus factors as $x\to y^b\to y$, where $y^b$ is the maximal subobject of $y$ which belongs to $\modPi^b$. We deduce that $z$ fits in an exact sequence 
	\begin{align}
		0\longrightarrow y^b/x\longrightarrow z \longrightarrow y' \longrightarrow 0\ ,
	\end{align}
	where $y' \in \Lambda^{(a,b)}_\qv$ and hence $z$ belongs to $\modPi^{\kappa}$. We therefore have the following diagram, for any $\bfd$ and any $\bfe$ of slope $b$
	\begin{align}
		\begin{tikzcd}[ampersand replacement=\&]
			\dLambda^{\kappa}_\bfd \times \dLambda^b_\bfe \ar{d}{j} \& \dLambda^{\mathsf{ext},\kappa}_{\bfd,\bfe} \arrow[swap]{l}{q^\circ} \ar{r}{p^\circ} \ar{d}{j'}\&
			\dLambda^{\kappa}_{\bfd+\bfe} \ar{d}{j''}\\
			\dLambda_\bfd \times \dLambda_\bfe  \& \dLambdaext_{\bfd,\bfe} \arrow[swap]{l}{q} \ar{r}{p}\&
			\dLambda_{\bfd+\bfe}
		\end{tikzcd}\ ,
	\end{align}
	where all the squares are cartesian and the morphisms $j,j',j''$ are open immersions. Again, \eqref{item:COHA-4} follows by base change.
	
	We now turn to \eqref{item:COHA-5}. Surjectivity is proved in Proposition~\ref{prop:pure}--\eqref{item:pure-1}. The kernel of the restriction morphism is the image of the map $\HBMbulletT(\Lambda_\qv \smallsetminus \Lambda_\qv^{\kappa}) \to \HBMbulletT(\Lambda_\qv)$, which is injective for purity reasons. We will show that $\HBMbulletT(\Lambda_\qv \smallsetminus \Lambda_\qv^{\kappa})=\K_{\Ltheta,\kappa}$.  Fix $\bfd, \bfe \in \N I$ and assume that $\mu_{\Ltheta}(\bfd) >b$. If 
	\begin{align}
		0\longrightarrow x\longrightarrow y \longrightarrow z \longrightarrow 0\ ,
	\end{align}
	is a short exact sequence of nilpotent $\Pi_\qv$-modules with $x$ of dimension $\bfd$ and $z$ of dimension $\bfe$ then the maximal slope submodule of $y$ has slope at least $\mu_{\Ltheta}(\bfd)>b$. Hence $y \in \Lambda_{\bfd + \bfe} \smallsetminus \Lambda_{\bfd + \bfe}^{\kappa}$. In other words, there is a factorization
	\begin{align}
		\trunc{p}\colon \trunc{\dLambdaext_{\bfe,\bfd}} \longrightarrow \Lambda_{\bfd + \bfe} \smallsetminus \Lambda_{\bfd + \bfe}^{\kappa} \longrightarrow \Lambda_{\bfd+\bfe} 
	\end{align}
	which proves that $\coha_\bfe^T \cdot \cohaqvd^T \subseteq \HBMbulletT(\Lambda_\qv \smallsetminus \Lambda^\kappa_\qv)$. The inclusion $\cohaqvd^T \cdot \coha_\bfe^T \subseteq \HBMbulletT(\Lambda_\qv \smallsetminus \Lambda^\kappa_\qv)$ for $\mu_{\Ltheta}(\bfd) \leqslant a$ is proved in a similar fashion, using minimal slope quotients. Thus $\K_{\Ltheta,\kappa} \subseteq  \HBMbulletT(\Lambda_\qv \smallsetminus \Lambda^\kappa_\qv)$. To finish the proof, we need to show the reverse inclusion. 
	
	We shall argue by induction on Harder-Narasimhan strata, so we first introduce some preliminary material. Let $HN(\bfd)$ stand for the (finite) set of Harder-Narasimhan types of weight $\bfd$, i.e., tuples $\underline{\bfd}=(\bfd_1, \ldots, \bfd_s)$ of elements of $\N I$ satisfying
	\begin{align}
		\sum_k \bfd_k=\bfd \quad \text{and}\quad \mu_{\Ltheta}(\bfd_1) < \cdots < \mu_{\Ltheta}(\bfd_s)\ .
	\end{align}
	Denote by $S_{\underline{\bfd}}$ the locally closed Harder-Narasimhan strata of $\Lambda_\bfd$ corresponding to $\underline{\bfd}$. We have
	\begin{align}
		\Lambda^\kappa_\bfd=\bigsqcup_{\underline{\bfd}} S_{\underline{\bfd}}\ ,
	\end{align}
	where the union ranges over all HN types $\underline{\bfd}=(\bfd_1, \ldots, \bfd_s)$ for which
	\begin{align}\label{eq:HNcondition}
		a < \mu_{\Ltheta}(\bfd_1) < \cdots < \mu_{\Ltheta}(\bfd_s) \leqslant b\ .
	\end{align}
	Let $HN'(\bfd)$ be the subset of HN types not satisfying \eqref{eq:HNcondition}. We define a total order on $HN(\bfd)$ as follows: for $\underline{\bfd}=(\bfd_1, \ldots, \bfd_s)$ and $\nu \in \Q$ we put 
	\begin{align}
		\vert\underline{\bfd}\vert_{\geqslant \nu}=\sum_{\mu_{\Ltheta}(\bfd_i) \geq \nu} \bfd_i\ .
	\end{align}
	We put $\underline{\bfd} \prec \underline{\bfd'}$ if there exists $\nu$ such that $\vert\underline{\bfd}\vert_{\geqslant \eta}=\vert\underline{\bfd}'\vert_{\geqslant \eta}$ for all $\eta > \nu$ while $\vert\underline{\bfd}\vert_{\geqslant \nu} > \vert\underline{\bfd}'\vert_{\geqslant \nu}$.
	
	To conclude the proof, we need the following lemma.
	\begin{lemma}\label{lem:HNorderclosed} 
		Fix $\underline{\bfd} \in HN(\bfd)$.
		\begin{enumerate}\itemsep0.2cm
			\item\label{item:HNClosure-1} The substack $S_{\preceq \underline{\bfd}}\coloneqq \bigsqcup_{\underline{\bfe} \preceq \underline{\bfd}} S_{\underline{\bfe}}$ is closed,
			\item \label{item:HNClosure-2} The multiplication map factors as $\coha_{\bfd_1}^T \otimes \cdots \otimes \coha_{\bfd_s}^T \to \HBMbulletT(S_{\preceq \underline{\bfd}}) \to \cohaqvd^T$, and the composition $\coha_{\bfd_1}^T \otimes \cdots \otimes \coha_{\bfd_s}^T \to \HBMbulletT(S_{\preceq \underline{\bfd}}) \to \HBMbulletT(S_{ \underline{\bfd}})$ is surjective.
		\end{enumerate}
	\end{lemma}
	
	\begin{proof}[Proof of Lemma~\ref{lem:HNorderclosed}] 
		To prove \eqref{item:HNClosure-1} it is enough to show that for any $\underline{\bfd}$, $\overline{S_{\underline{\bfd}}} \subset S_{\preceq \underline{\bfd}}$. Consider the induction diagram
		\begin{align}\label{eq:proofLemmaHNclosure}
			\begin{tikzcd}[ampersand replacement=\&]
				\dLambda^{\mathsf{ss}}_{\bfd_1} \times \cdots \times \dLambda^{\mathsf{ss}}_{\bfd_s} \ar{d}{j} \& 	\dLambda^{\mathsf{ext},\mathsf{ss}}_{\underline{\bfd}} \arrow[swap]{l}{q^\circ} \ar{r}{p^\circ} \ar{d}{j'}\&
				\dLambda^{\ell_1, \ldots, \ell_s}_{\bfd_1, \ldots, \bfd_s} \ar{d}{h}\\
				\dLambda_{\bfd_1} \times \cdots \times \dLambda_{\bfd_s}   \& \dLambdaext_{\underline{\bfd}} \arrow[swap]{l}{q} \ar{r}{p}\&
				\dLambda_{\bfd}
			\end{tikzcd}\ ,
		\end{align}
		where $\ell_i=\mu_{\Ltheta}(\bfd_i)$ and where the leftmost square is a pull back. Note that $p^\circ$ is an equivalence and $q^\circ$ is a  vector bundle stack morphism (see Proposition~\ref{prop:pure}). Moreover, $\Lambda^{\ell_1, \ldots, \ell_s}_{\bfd_1, \ldots, \bfd_s}=S_{\underline{\bfd}}$. As $p$ is proper, $\overline{S_{\underline{\bfd}}}$ is contained in the image of $p$, hence it suffices to prove that there exists a factorization
		\begin{align}\label{eq:pfactorization}
			\trunc{p}\colon  \trunc{\dLambdaext_{\underline{\bfd}}} \stackrel{p'}{\longrightarrow}  S_{\preceq \underline{\bfd}} \longrightarrow \Lambda_{\bfd}\ . 
		\end{align}
		Let $x \in \Lambda_\bfd$ and assume that $x$ has a filtration $x_{s+1}=0 \neq x_s \subset x_{s-1} \subset \cdots \subset x_1=x$ with $x_i/x_{i+1} \in \Lambda_{\bfd_i}$. Unless $x \in S_{\underline{\bfd}}$, there exists $j$ such that $x_k/x_{k+1} \in \Lambda^{\mathsf{ss}}_{\bfd_k}$ for $k >j$ while $x_j/x_{j+1} \in S_{\underline{\bfn}}$ with $\underline{\bfn} =(\bfn_1, \ldots, \bfn_r) \neq (\bfd_j)$. But then $\mu_{\Ltheta}(\bfn_r) > \mu_{\Ltheta}(\bfd_j)=\ell_j$. We deduce the existence of a subrepresentation $x_{j} \subset x' \subset x_{j+1}$ such that $\mu_{\Ltheta -\mathsf{min}}(x') > \ell_j$ while  $\dim(x') > \sum_{k >j}\bfd_k$, which implies that $x\in S_{\underline{\bfe}}$ for some $\underline{\bfe} \prec \underline{\bfd}$. This proves the factorization \eqref{eq:pfactorization}, from which we deduce both \eqref{item:HNClosure-1} and the first part of \eqref{item:HNClosure-2}. 
		
		The last part of \eqref{item:HNClosure-2} follows from the base change formula 
		\begin{align}
			p^\circ_\ast \circ (q^{\circ})^!\circ j^\ast= (j')^\ast\circ p'_\ast\circ q^!\ ,
		\end{align}
		where $j'\colon S_{\underline{\bfd}} \to S_{\preceq \underline{\bfd}}$ is the open embedding. Observe that, by Proposition~\ref{prop:pure}, the map $j^\ast \colon \HBMbulletT(\dLambda_{\bfd_1} \times \cdots \times \dLambda_{\bfd_s}) \to \HBMbulletT(\dLambda_{\bfd_1}^{\mathsf{ss}} \times \cdots \times \dLambda^{\mathsf{ss}}_{\bfd_s})$ is surjective.
	\end{proof}
	
	We return to \eqref{item:COHA-5}. Assume that $\HBMbulletT(\Lambda_\qv \smallsetminus \Lambda^\kappa_\qv) \not\subseteq \K_{\Ltheta,\kappa}$ and let $\underline{\bfd}=(\bfd_1,\ldots,\bfd_s) \in HN'(\bfd)$ be minimal such that $\HBMbulletT(S_{\preceq \underline{\bfd}}) \not\subseteq \K_{\Ltheta, \kappa}$. By hypothesis, we have $\mu_{\Ltheta}(\bfd_1) \leqslant a$ or $\mu_{\Ltheta}(\bfd_s) >b$. In either case, $\coha_{\bfd_1}^T \cdots \coha_{\bfd_s}^T \subset \K_{\Ltheta,\kappa}$. By the minimality of $\underline{\bfd}$, $\HBMbulletT(S_{\underline{\bfe}}) \subset \K_{\Ltheta,\kappa}$ for any $\underline{\bfe}\prec \underline{\bfd}$, and thus $\HBMbulletT(S_{\prec \underline{\bfd}}) \subset \K_{\Ltheta,\kappa}$. But then, by Lemma~\ref{lem:HNorderclosed}--\eqref{item:HNClosure-2}, we have 
	\begin{align}
		\HBMbulletT(S_{\preceq \underline{\bfd}})= \HBMbulletT(S_{\prec \underline{\bfd}}) + \coha_{\bfd_1}^T \cdots \coha_{\bfd_s}^T \subseteq \K_{\Ltheta,\kappa}\ ,
	\end{align} 
	a contradiction. This yields the reverse inclusion $\HBMbulletT(\Lambda_\qv \smallsetminus \Lambda^\kappa_\qv) \subseteq \K_{\Ltheta,\kappa}$ and finishes the proof of Statement~\eqref{item:COHA-5}.
	
	Finally, Statement~\eqref{item:COHA-3} follows in a similar fashion from the fact that the category $\modPi^{(-\infty,b]}$ is closed under taking subobjects (since it is the torsion-free part of a torsion pair -- see Lemma~\ref{lem:torsion-pair-quiver}), and likewise $\modPi^{[a,\infty)}$ is closed under taking quotients.
\end{proof}

\subsection{Quotients of the Yangian}\label{sec:quotients-Yangians}

We now examine the counterparts, for the Yangian, of the quotients considered in the previous section. For any interval $\kappa=(a,b]$ (including the case $a=-\infty$ and/or $b=\infty$), we define
\begin{align}\label{eq:quotient-Yangian-kappa}
	\J_{\Ltheta,\kappa}\coloneqq \sum_{\mu_{\Ltheta}(\bfd)>b}\Y_\qv^-\Y^-_{\textrm{-}\bfd}+\sum_{\mu_{\Ltheta}(\bfd)\leqslant a}\Y^-_{\textrm{-}\bfd}\Y_\qv^- \qquad \text{and}\qquad \Y_{\Ltheta, \kappa}\coloneqq \Y_\qv^-/\J_{\Ltheta, \kappa}\ .
\end{align}
We extend this definition to other types of intervals in an obvious manner. Note that $\Y_{\Ltheta, \kappa}$ is $\N \times \Z I$-graded. Moreover, $\J_{\Ltheta,\kappa} \subset \J_{\Ltheta,\kappa'}$ whenever $\kappa' \subset \kappa$. 

\begin{proposition}\label{prop:various_COHA_actions_truncated}
	The following holds.
	\begin{enumerate}\itemsep0.2cm
		\item \label{item:yangian-kappa-action1} For any $b \in\Q$, $\Y_{\Ltheta,\{b\}}$ inherits from $\Y_\qv$ of an associative algebra structure.
		
		\item \label{item:yangian-kappa-action3} For $\kappa=[a,b]$, there is a canonical right action of $\Y_{\Ltheta,\{b\}}$ and a canonical left action of $\Y_{\Ltheta,\{a\}}$ on $\Y_{\Ltheta,\kappa}$.
		
		\item \label{item:yangian-kappa-action2} For $\kappa=(-\infty,b]$, there is a canonical left action of $\Y^-_\qv$ on $\Y_{\Ltheta,\kappa}$. For $\kappa=[a,+\infty)$ there is a canonical right action of $\Y^-_\qv$ on $\Y_{\Ltheta, \kappa}$. Similar results hold for intervals of the form $(-\infty, b)$ or $(a,+\infty)$.
	\end{enumerate}	
\end{proposition}

\begin{proof}
	By definition, $(\Y_{\Ltheta, \{b\}})_\bfd$ is nonzero only when $\mu_{\Ltheta}(\bfd)=b$, hence we may view $\Y_{\Ltheta,\{b\}}$ as the quotient
	of 
	\begin{align}
		U\coloneqq\bigoplus_{\mu_{\Ltheta}(\bfd)=b} \Y^-_\bfd
	\end{align}
	by $U \cap \J_{\Ltheta,\{b\}}$. Observe that if $\bfe$ satisfies $\mu_{\Ltheta}(\bfe)<b$ (resp. $\mu_{\Ltheta}(\bfe)>b$) and $\bfd$ satisfies $\mu_{\Ltheta}(\bfd)=b$ then $\mu_{\Ltheta}(\bfe+\bfd)<b$ (resp. $\mu_{\Ltheta}(\bfe+\bfd)>b$). It follows that $U \cap\J_{\Ltheta,\{b\}}$ is an ideal of $U$, hence \eqref{item:yangian-kappa-action1}. 
	
	Statement~\eqref{item:yangian-kappa-action3} is proved in a similar way, while Statement~\eqref{item:yangian-kappa-action2} is obvious from the definition of $\J_{\Ltheta, \kappa}$.
\end{proof}

Combining the above with Proposition~\ref{prop:COHA} yields the following.
\begin{corollary}\label{cor:isom-PHI-for-quotients} 
	\hfill
	\begin{enumerate}\itemsep0.2cm
		\item \label{item:Phi-for-quotients1}
		For any interval $\kappa \subset \Q$, the isomorphism $\Phi\colon \Y^-_\qv \to \cohaqv^T$ yields an isomorphism of vector spaces $\Phi_\kappa\colon \Y_{\Ltheta,\kappa} \to \coha_{\Ltheta,\kappa}^T$.
		
		\item \label{item:Phi-for-quotients2} For any $b\in \R$, the map $\Phi_b\colon \Y_{\Ltheta,\{b\}} \to \coha_{\Ltheta,\{b\}}^T$ is an isomorphism of algebras.
		
		\item \label{item:Phi-for-quotients3} For any interval $\kappa$, the isomorphism $\Phi_\kappa$ intertwines the various right or left actions of $\Y_{\Ltheta,\{b\}}$ and $\coha_{\Ltheta,\{b\}}^T$ or $\Y^-_\qv$ and $\cohaqv^T$.
	\end{enumerate}
\end{corollary}

To gain a better understanding of the quotients $\Y_{\Ltheta,\kappa}$, we now describe their classical limit in terms of double-loop root systems. For an interval $\kappa \subset \Q$, we put
\begin{align}
	\Delta_{\Ltheta, \kappa}\coloneqq \Big\{\bfd\in \N I\, \Big\vert\, \mu_{\Ltheta}(\bfd) \in \kappa \Big\}\ .
\end{align}	
Note that $\Delta_{\Ltheta,\kappa}$ is stable under addition inside $\Delta$  if $\bfd, \bfe \in \Delta_{\Ltheta,\kappa}$ are such that $\bfd + \bfe \in \Delta$ then $\bfd + \bfe \in \Delta_{\Ltheta, \kappa}$. In particular, 
\begin{align}
	(\fraknell)_{\Ltheta,\kappa}\coloneqq \bigoplus_{\bfd \in \Delta_{\Ltheta,\kappa}} (\fraknell)_{\textrm{-}\bfd}
\end{align} 
is a Lie subalgebra of $\fraknell$. For any $\kappa$, the standard filtration on $\Y_\qv^-$ induces one on $\J_{\Ltheta,\kappa}$, and descends to a filtration on $\Y_{\Ltheta,\kappa}$. We denote by $\gr \Y_{\Ltheta,\kappa}$ the associated graded. Recall the isomorphism $\Psi\colon \sfU(\fraknell)\otimes \Hbullet_T \to \gr\Y^-_\qv$.

\begin{proposition} 
	The following hold:
	\begin{enumerate}\itemsep0.2cm
		\item \label{item:Phi-for-gr-quotients1} For any $\kappa$ there is a canonical isomorphism $\Psi_\kappa\colon \sfU((\fraknell)_{\Ltheta,\kappa})\otimes_\Q \Hbullet_T \to \gr\Y_{\Ltheta,\kappa}$.
		\item \label{item:Phi-for-gr-quotients2} For any $\nu \in\Q$, there is a canonical graded algebra isomorphism $\Psi_\nu\colon \sfU((\fraknell)_{\Ltheta,\{\nu\}})\otimes_\Q \Hbullet_T \to \gr\Y_{\Ltheta,\{\nu\}}$.
		\item \label{item:Phi-for-gr-quotients3} The above isomorphisms intertwine the right or left actions of $\sfU((\fraknell)_{\Ltheta,\kappa})\otimes_\Q \Hbullet_T$ and $\gr\Y_{\Ltheta,\kappa}$ or of $\sfU((\fraknell)_{\Ltheta,\{\nu\}})\otimes_\Q \Hbullet_T$ and $\gr\Y_{\Ltheta,\{\nu\}}$.
	\end{enumerate}
\end{proposition}

\begin{proof}
	The compatibility with the algebra structure or actions being clear from the construction, we focus on \eqref{item:Phi-for-gr-quotients1}. For $\kappa=(a,b]$ there is a canonical embedding
	\begin{align}
		\gr\J_{\Ltheta,\kappa} \subseteq & \sum_{\mu_{\Ltheta}(\bfd)>b}\gr\Y_\qv^-\cdot \gr\Y^-_{\textrm{-}\bfd}+\sum_{\mu_{\Ltheta}(\bfd)\leqslant a}\gr\Y^-_{\textrm{-}\bfd}\cdot \gr\Y_\qv^-\\
		&=\sum_{\mu_{\Ltheta}(\bfd)>b}\sfU(\fraknell)\cdot \sfU(\fraknell)_{\textrm{-}\bfd}+\sum_{\mu_{\Ltheta}(\bfd)\leqslant a}\sfU(\fraknell)_{\textrm{-}\bfd}\cdot \sfU(\fraknell)\ .
	\end{align}
	It follows that there is a canonical surjection
	\begin{align}
		\begin{tikzcd}[ampersand replacement=\&]
			\gr\Y_{\Ltheta,\kappa} =\gr\Y^-/ \gr\J_{\Ltheta,\kappa} \ar{r}{\pi}\& \sfU(\fraknell)\Big/\Big( 
			\sum_{\mu_{\Ltheta}(\bfd)>b}\sfU(\fraknell)\cdot \sfU(\fraknell)_{\textrm{-}\bfd}+\sum_{\mu_{\Ltheta}(\bfd)\leqslant a}\sfU(\fraknell)_{\textrm{-}\bfd}\cdot \sfU(\fraknell)
			\Big)
		\end{tikzcd} \ ,
	\end{align}
	and the right-hand-side coincides with $\sfU((\fraknell)_{\Ltheta,\kappa})$ thanks to the PBW theorem
	\begin{align}
		\sfU(\fraknell)=\sfU((\fraknell)_{\Ltheta, (-\infty,a]})\otimes \sfU((\fraknell)_{\Ltheta,\kappa})\otimes \sfU((\fraknell)_{\Ltheta, (b,\infty)}) 
	\end{align}
	together with the observations that
	\begin{align}
		\sum_{\mu_{\Ltheta}(\bfd)\leqslant a} \sfU(\fraknell)_{\textrm{-}\bfd}\cdot \sfU(\fraknell)&= \sfU((\fraknell)_{\Ltheta,(-\infty,a]})_+\cdot \sfU(\fraknell)\ ,\\ \sum_{\mu_{\Ltheta}(\bfd)>b} \sfU(\fraknell)\cdot \sfU(\fraknell)_{\textrm{-}\bfd}&= \sfU(\fraknell) \cdot \sfU((\fraknell)_{\Ltheta,(b,\infty)})_+\ ,
	\end{align}
	where as before $(-)_+$ denotes the augmentation ideal. It it thus enough to prove that $\gr\Y_{\Ltheta,\kappa}$ and $\sfU((\fraknell)_{\Ltheta,\kappa})$ have the same graded dimension. For this, we may use the isomorphism of graded vector spaces $\Y_{\Ltheta,\kappa} \simeq \coha_{\Ltheta,\kappa}^T$ and Proposition~\ref{prop:coha-kac}--\eqref{item:coha-kac-4}, which gives us
	\begin{align}
		\sum_{\ell,\bfd} \dim((\coha_{\Ltheta,\kappa}^T)_{\bfd,\ell})z^\bfd t^\ell=(1-t^{-2})^{-2}\Exp\Bigg(\frac{1}{1-t^{-2}}\sum_{\substack{\bfd\in\N I\\ \mu_{\Ltheta}(\bfd)\in\kappa}} A_\bfd(t^{-2})\,z^\bfd \Bigg)\ ,
	\end{align}
	while
	\begin{align}
		\sum_{\ell,\bfd} \dim((\sfU((\fraknell)_{\Ltheta,\kappa}) \otimes \Hbullet_T)_{\textrm{-}\bfd,\ell})z^\bfd t^\ell=(1-t^{-2})^{-2}\Exp\Bigg(\sum_{\substack{\bfd\in\N I, \ell\\ \mu_{\Ltheta}(\bfd)\in\kappa}} \dim(\fraknell)_{\textrm{-}\bfd,\ell}\,z^\bfd \Bigg)\ .
	\end{align}
	The proof now follows from the relation~\eqref{eq:dim-Lie-Kac-pols}, see Lemma~\ref{lem:dimension-Ln}.
\end{proof}

\section{Perverse coherent sheaves on resolutions of ADE singularities}\label{sec:perverse-coherent-sheaves}

In this section we relate, via the McKay correspondence, the derived category of coherent sheaves on a Kleinian resolution of singularities $\rsv$ with the derived category of representations of the preprojective algebra $\Pi_\qv$ of the associated affine quiver. We introduce for this the perverse $t$-structure on $\catDb(\catCoh(\rsv))$, following \cite{VdB_Flops}. We also relate (complexes of) coherent sheaves on $\rsv$ supported on the exceptional locus with nilpotent representations of $\Pi_\qv$.

\subsection{Geometry of the resolutions}\label{subsec:geometry-resolutions}

The quotient $\ssv\coloneqq\C^2/G$ has an isolated singularity at the origin. Let $\pi\colon \rsv\to \ssv$ be the minimal resolution of singularities of $\ssv$.
Fix a diagonal torus $A\subset \mathsf{GL}(2,\C)$ centralizing $G$. For $G$ of type $A$, the diagonal torus $A$ could be $\{1\}$, $\G_m$, or $\G_m\times\G_m$, while for $G$ of types $D$ or $E$, it could be $\{1\}$ or $\G_m$. The group $A$ acts on $\ssv$ in the obvious way and this action lifts to a $A$-action on $\rsv$ such that the map $\pi$ is $A$-equivariant (one may see this using e.g. the Nakajima quiver varieties interpretation of the map $\pi$).

Let $C\coloneqq \pi^{-1}(0)$ and denote by $C_{\red}$ its reduced variety. The irreducible components $C_1, \ldots, C_e$ of $C_{\red}$ are isomorphic to $\PP^1$. The following is well-known.
\begin{proposition}
	\hfill
	\begin{enumerate}\itemsep0.2cm
		\item \label{Mckay-(finite)} The vertices of $\qvfin$ are in bijection with the irreducible components $C_i$ of $C_{\red}$. Two vertices are joined by an edge if and only if the corresponding components intersect. The intersection is transverse and consists of one point. The intersection matrix of the $C_i$'s is equal to the opposite of the Cartan matrix of $\qvfin$.
		\item There is a bijection $i\mapsto \eta_i$ from the set $I$ of vertices of $\qv$ to the set $\widehat{G}$ of irreducible representations of $G$, and for any $i,j\in I$ the number of arrows in $\overline{\Omega}$ from $i$ to $j$ is the dimension of the vector space $\Hom_G(\eta_i\otimes\C^2,\eta_j)$.
	\end{enumerate}
\end{proposition}

We have the following equality in $\Pic(\rsv)$:
\begin{align}
	C= \sum_{i=1}^e\, r_i C_i\ ,
\end{align}
where the numbers $r_i$ are introduced in Formula~\eqref{eq:r_i}.
%
\begin{lemma}\label{lem:PicMcKay}
	There is a canonical isomorphism 
	\begin{align}\label{eq:isoPic}
		\begin{tikzcd}[ampersand replacement=\&]
			\Pic(\rsv)\ar{r}{\sim}\& \coweightlatticefin\ ,
		\end{tikzcd}
	\end{align}
	mapping $\scrO_Y(-C_i)$ to $\Lalpha_i$ for $i=1, \ldots, e$. 
\end{lemma}

\begin{notation}
	We denote by $\calL_{\Llambda}$ the line bundle associated to the coweight $\Llambda\in\coweightlatticefin$.
\end{notation}

\subsection{Perverse coherent sheaves on the resolution}\label{subsec:psheavesres} 

Let $\scrC$ be the so-called \textit{null category} of the pair $(\ssv,\rsv)$, i.e., the abelian subcategory of $\catCoh(\rsv)$ consisting of objects $E$ such that $\R\pi_\ast E=0$. Following \cite[\S3.1]{VdB_Flops}, we introduce the following torsion pair on $\catCoh(\rsv)$ (denoted in \textit{loc.\ cit.} by $(\scrT_{\textrm{-}1}, \scrF_{\textrm{-}1})$):
\begin{align}
	\scrT&\coloneqq \{\calF\in \catCoh(\rsv)\, \vert\, \R^1\pi_\ast \calF=0\,\text{ and }\, \Hom(\calF, \scrC)=0\}\ ,\\[2pt]
	\scrF&\coloneqq \{\calF\in \catCoh(\rsv)\, \vert\, \R^0\pi_\ast \calF=0\}\ .
\end{align}
We denote by $\catP(\rsv/\ssv) \subset \catDb(\catCoh(\rsv))$  the heart of the tilted $t$-structure induced by $(\scrT, \scrF)$. We call it the \textit{Bridgeland's perverse $t$-structure}. For $E\in \catP(\rsv/\ssv)$, all cohomology sheaves $\calH^i(E)$ vanish except for $i=-1, 0$, $\calH^{-1}(E)\in \scrF$, and $\calH^0(E)\in \scrT$. We call $E$ a \textit{perverse coherent sheaf}.

Consider now the abelian subcategory $\catCohps(\rsv)$ of $\catCoh(\rsv)$ consisting of properly supported coherent sheaves. Then the pair $(\scrT\cap \catCohps(\rsv), \scrF\cap \catCohps(\rsv))$ is a torsion pair of $\catCohps(\rsv)$. We denote by $\catPps(\rsv/\ssv)$ the corresponding tilted heart in $\catDb(\catCohps(\rsv))\simeq \catDbps(\catCoh(\rsv))$, where the latter is the subcategory of $\catDb(\catCoh(\rsv))$ consisting of complexes with properly supported cohomology. Thus, $\catPps(\rsv/\ssv)$ is the abelian subcategory of $\catP(\rsv/\ssv)$ consisting of objects with proper support. 
\begin{remark}\label{rem:torsion-pair-ps}
	By \cite[\S2.1]{YoshiokaPerverse}, see also \cite[Appendix~A]{DPS_McKay}, one has the following characterization of these intersections:
	\begin{align}
		\scrT\cap \catCohps(\rsv) = &\{\calF\in \catCohps(\rsv)\, \vert\, \Hom(\calF, \scrO_{C_i}(-1))=0 \text{ for } i=1, \ldots, e\} \ ,\\[2pt]
		\scrF\cap \catCohps(\rsv)  = &\{\calF\in \catCohps(\rsv)\, \vert\, \calF \text{ is a successive extension of subsheaves of the }\\[2pt] 
		&\scrO_{C_i}(-1)\text{'s for }i=1, \ldots, e\} \ . 	
	\end{align}
	Hence, $\calH^{-1}(E)$ is a pure one-dimensional coherent sheaf for $E\in \catPps(\rsv/\ssv)$.
\end{remark}

We introduce the following subcategories of $\catP(\rsv/\ssv)$ (resp.\ $\catPps(\rsv/\ssv)$):
\begin{itemize}\itemsep0.2cm
	\item $\catP^0(\rsv/\ssv)$ (resp.\ $\catPps^0(\rsv/\ssv)$) is the full subcategory of $\catP(\rsv/\ssv)$ (resp.\ $\catPps(\rsv/\ssv)$) consisting of perverse coherent sheaves $E$ with $\calH^{-1}(E)=0$;
	\item $\catP^{-1}(\rsv/\ssv)$ (resp.\ $\catPps^{-1}(\rsv/\ssv)$) is the full subcategory of $\catP(\rsv/\ssv)$ (resp.\ $\catPps(\rsv/\ssv)$) consisting of perverse coherent sheaves $E$ with $\calH^{0}(E)=0$.
\end{itemize}

By construction, the pairs $(\catP^0(\rsv/\ssv), \catP^{-1}(\rsv/\ssv)[-1])$ and $(\catPps^0(\rsv/\ssv), \catPps^{-1}(\rsv/\ssv)[-1])$ are torsion pairs of $\catCoh(\rsv)$ and $\catCohps(\rsv)$, respectively. From the defining properties of torsion pairs (see also \cite[Lemma~2.22]{DPS_McKay}), we have:
\begin{lemma}\label{lem:perverse-torsion-pair}
	\hfill
	\begin{itemize}\itemsep0.2cm
		\item Let $E\in \catPps^0(\rsv/\ssv)$ and $\calF\coloneqq \calH^0(E)$. Then, any quotient $\calF\twoheadrightarrow \calF''$ in $\catCohps(\rsv)$ corresponds to an object in $\catPps^0(\rsv/\ssv)$.
		\item Let $E\in \catPps^{-1}(\rsv/\ssv)$ and $\calF\coloneqq \calH^{-1}(E)$. Then any subsheaf $\calF'\hookrightarrow \calF$ in $\catCohps(\rsv)$ corresponds to an object in $\catPps^{-1}(\rsv/\ssv)$.
	\end{itemize}
\end{lemma} 

\begin{proposition}\label{prop:zero-dimensional-perverse}
	If the support\footnote{Recall that the support of an object in $\catDb(\catCoh(\rsv))$ is the union of the supports of its cohomology sheaves.} of $E\in \catPps(\rsv/\ssv)$ is zero-dimensional, then $E\in \catPps^0(\rsv/\ssv)$. Vice versa, if $\calF\in \catCohps(\rsv)$ is zero-dimensional, then $\calF\in \catPps^0(\rsv/\ssv)$.
\end{proposition}

\begin{proof}
	We prove the first statement. Let $E\in \catPps(\rsv/\ssv)$ be a perverse coherent sheaf whose support is zero-dimensional. Thus both $\calH^{-1}(E)$ and $\calH^0(E)$ are zero-dimensional. On the other hand, since $\R^0\pi_\ast\calH^{-1}(E)=0$, we get $\calH^{-1}(E)=0$, therefore $E\simeq \calH^0(E)$. 
	
	Conversely, let $\calF\in \catCohps(\rsv)$ be zero-dimensional. We have $\R^1\pi_\ast \calF=0$ since $H^1(\rsv, \calF)=0$. Now, consider $\calE\in \scrC$. By definition, $\R\pi_\ast \calE=0$. Hence, $\calE$ does not admit any zero-dimensional subsheaf. Thus $\Hom(\calF, \calE)=0$. Therefore, $\calF\in \scrT\cap \catCohps(\rsv)=\catPps^0(\rsv/\ssv)$.
\end{proof}

\subsection{Tilting equivalence}

Following \cite[Lemma~3.4.4]{VdB_Flops}, let $D_i$ be the unique Cartier divisor of $\rsv$ such that $D_i \cdot C_j=\delta_{i, j}$ for $i,j=1,\ldots, e$. Set $D\coloneqq D_1+\cdots+D_e$. Under the isomorphism \eqref{eq:isoPic}, the line bundle $\scrO_\rsv(D_i)$ corresponds to the fundamental coweights $\Llambda_i$ of $\qvfin$ for $i=1, \ldots, e$. Moreover, $\dim\,H^1(\rsv, \scrO_\rsv(-D_i))=r_i-1$ where $r_i$ is as in Formula~\eqref{eq:delta}.

For each $i=1, \ldots, e$, let $\calE_i$ be the locally free sheaf obtained as the universal extension
\begin{align}\label{eq:projgenZN_components}
	0\longrightarrow \scrO_\rsv^{r_i-1}\longrightarrow \calE_i\longrightarrow \scrO_\rsv(D_i) \longrightarrow 0
\end{align}
associated to a basis of $H^1(\rsv, \scrO_\rsv(-D_i))$. We set $\calE_0 \coloneqq \scrO_\rsv$. 

\begin{remark}
	By construction, the $\calE_i$'s are the locally free sheaves defined by Gonzalez-Sprinberg and Verdier \cite{GSV_McKay}. As shown in \textit{loc.cit.}, they generate, over $\Z$, the Grothendieck group of locally free sheaves on $\rsv$. 
\end{remark}
By \cite[Theorem~3.5.5]{VdB_Flops}, each $\calE_i$ is an indecomposable projective object in $\catP(\rsv/\ssv)$ and the direct sum 
\begin{align}\label{eq:projgenZN}
	\calP \coloneqq \bigoplus_{i=0}^{e}\, \calE_i
\end{align}
is a projective generator in $\catP(\rsv/\ssv)$. 

Set $R\coloneqq \C[\A^2]$ and $R_\eta\coloneqq (\eta^\vee \otimes R)^G$ for any irreducible representation $\eta$ of $G$. The following hold (see e.g. \cite{CH_DeformedPreprojectiveAlgebras, Wemyss_McKay}):
\begin{align}
	R\rtimes G=\End_{R^G}(R)\ ,\quad R^G=\C[\ssv]\ ,\quad \Pi_\qv=\End(\calP)=\End_{R^G}(\pi_\ast \calP)\ ,\quad\text{and} \quad\pi_\ast \calP=\bigoplus_{\eta\in\widehat{G}}R_\eta \ .
\end{align}
We equip the ring $R$ with the obvious $A$-action. This yields an $A$-action on the algebra $\Pi_\qv$ by algebra automorphisms. Hence, the group $A$ acts on the categories $\catCoh(\rsv)$ and $\ModPi$ and therefore on the triangulated categories $\catDb(\catCoh(\rsv))$ and $\catDb(\ModPi)$. By combining this with \cite[Corollary~3.2.8]{VdB_Flops}, we have the following result.
\begin{theorem}\label{thm:tau}
	$\calP$ is a tilting object of $\catDb(\catCoh(\rsv))$. The functors 
	\begin{align}\label{eq:tau}
		\tau\coloneqq \R\Hom(\calP,-)\colon \catDb(\catCoh(\rsv)) &\longrightarrow \catDb(\ModPi)\ , \\[4pt] 
		(-) \otimes^{\LL}_{\Pi_\qv} \calP \colon \catDb(\ModPi) &\longrightarrow \catDb(\catCoh(\rsv)) 
	\end{align}
	determine an equivalence of derived categories, which is equivariant with respect to the action of $A$. The same holds if we consider the categories $\catDbps(\catCoh(\rsv))\simeq \catDb(\catCohps(\rsv))$ and $\catDbps(\ModPi)\simeq \catDb(\modPi)$. These functors restrict to an $A$-equivariant equivalence of abelian categories
	\begin{align}
		\begin{tikzcd}[ampersand replacement=\&]
			\catP(\rsv/\ssv) \ar{r}{\sim} \&  \ModPi
		\end{tikzcd}\ ,
	\end{align}
	as well as  
	\begin{align}\label{eq:derived-equivalence-ps}
		\begin{tikzcd}[ampersand replacement=\&]
			\catPps(\rsv/\ssv) \ar{r}{\sim} \& \modPi
		\end{tikzcd}\ .
	\end{align}
\end{theorem}

\begin{remark}\label{rem:simple}
	By \cite[Proposition~3.5.7]{VdB_Flops} the simple modules $\sigma_0, \ldots, \sigma_e$ associated to the nodes of the affine quiver correspond to the \textit{spherical objects}, in the sense of Definition~\ref{def:spherical},
	\begin{align}\label{eq:sphericalobj}
		\sigma_0 \simeq \tau(\scrO_{C})\quad\text{and}\quad	\sigma_i \simeq \tau(\scrO_{C_i}(-1)[1]) \ ,
	\end{align}
	for $1\leq i \leq e$. Set
	\begin{align}
		\calI_j\coloneqq \begin{cases}
			\scrO_{C} & \text{for } j=0\ ,\\
			\scrO_{C_j}(-1)[1] & \text{for } j=1, \ldots, e\ .
		\end{cases}
	\end{align}
	These complexes satisfy the orthogonality relations 
	\begin{align}\label{eq:orthrel} 
		\R\Hom(\calE_i, \calI_j) = \delta_{i,j}\, {\underline \C}\ , 
	\end{align}
	for $0\leq i, j \leq e$.
\end{remark}

Given an object $E\in \catPps(\rsv/\ssv)$, the categorical equivalence \eqref{eq:derived-equivalence-ps} implies that $\R\Hom(\calE_k, E)$ is a one-term complex of amplitude $[0, \, 0]$, which we denote by $V_k$, for any $k=0, \ldots, e$. Thus, an object $E$ of $\catPps(\rsv/\ssv)$ is mapped to a representation $M$ of $\Pi_\qv$ with underlying $\Z_{e+1}$-graded vector space 
\begin{align}\label{eq:vectorspace}
	V\coloneqq \bigoplus_{k=0}^{e}\ V_k\ , 
\end{align}
where the $k$-th summand has degree $k$. We denote by $d_k(E)$ the dimension of $V_k$ for $0\leq k \leq e$. We set $\bfd(E)\coloneqq (d_k(E)) \in \N I$ and call it the \textit{dimension vector} of $E$. Since $E$ is properly supported one has
\begin{align}\label{eq:ch1}
	\ch_1(E)\coloneqq -\ch_1(\calH^{-1}(E))+ \ch_1(\calH^{0}(E))=\sum_{i=1}^{e}\, n_i C_i\quad\text{and}\quad \chi(\scrO_\rsv, E)=n\ ,
\end{align}
for some $n, n_i\in \Z$, $i=1, \ldots, e$.
\begin{lemma}[{\cite[Lemma~2.12]{DPS_McKay}}]\label{lem:dimension-vector} 
	We have 
	\begin{align}\label{eq:dimvect} 
		r_k d_0(E)-d_k(E)=n_k\ ,
	\end{align}
	for any $1\leq k \leq e$ and 
	\begin{align}\label{eq:dNchi}
		(\Lomega_0, \bfd(E)) = d_0(E) = n  \ .
	\end{align}
\end{lemma}

\subsection{Nilpotent perverse coherent sheaves and nilpotent representations}

The goal of this section is to establish a version of the tilting equivalence \eqref{eq:derived-equivalence-ps} for \textit{nilpotent} objects. 

Denote by $\nilpPi$ the abelian category of finite-dimensional nilpotent representations of $\Pi_\qv$ (cf.\ Definition~\ref{def:nilpotent-representation} and Remark~\ref{rem:nilpotency}), and by $\catPnil(\rsv/\ssv)$ the abelian category of perverse coherent sheaves on $\rsv$ set-theoretically supported on the exceptional curve $C$. They are Serre subcategories of $\modPi$ and $\catPps(\rsv/\ssv)$. 

Let $\catDb_C(\catCoh(\rsv))$ denote the bounded derived category of complexes of coherent sheaves on $\rsv$ with cohomology sheaves set-theoretically supported on $C$ and let $\catDb_{\ps,\,\mathsf{nil}}(\ModPi)$ denote the bounded derived category of complexes of $\Pi_\qv$-representations with nilpotent cohomology objects.
\begin{theorem}\label{thm:nilpotent}
	The functors 
	\begin{align}
		\tau\coloneqq \R\Hom(\calP,-)\colon \catDb_C(\catCoh(\rsv)) &\longrightarrow \catDb_{\ps,\,\mathsf{nil}}(\ModPi)\ , \\[4pt] 
		(-) \otimes^{\LL}_{\Pi_\qv} \calP \colon \catDb_{\ps,\,\mathsf{nil}}(\ModPi) &\longrightarrow \catDb_C(\catCoh(\rsv)) 
	\end{align}
	determine an equivalence of derived categories, which is equivariant with respect to the action of $A$. This restricts to an $A$-equivariant equivalence of abelian categories
	\begin{align}
		\begin{tikzcd}[ampersand replacement=\&]
			\catPnil(\rsv/\ssv) \ar{r}{\sim} \&  \nilpPi
		\end{tikzcd}\ .
	\end{align}
\end{theorem}
In the remaining part of the section, we will show that the categories $\nilpPi$ and $\catPnil(\rsv/\ssv)$ are mapped to each other under the equivalence $\tau$. This will imply Theorem~\ref{thm:nilpotent}.

\medskip

Let us recall two results from \cite{SY13}, which will be useful later on. Denote by $\scrS$ the set of isomorphism classes of simple representations of $\Pi_\qv$, which are not isomorphic to one of the $\{\sigma_0, \ldots, \sigma_e\}$. 
\begin{lemma}[{\cite[Lemma~2.31]{SY13}}]\label{lem:nilpotency} 
	\hfill
	\begin{enumerate}\itemsep0.2cm
		
		\item \label{item:nilpotency-1} The dimension vector of any representation in $\modPi$ in $\scrS$ is $\delta$. 
		
		\item \label{item:nilpotency-2}  If the dimension vector of a representation $M$ of $\Pi_\qv$ is $\delta$ then $M$ is either simple or nilpotent. 
		
		\item \label{item:nilpotency-3}  Suppose $M,N$ are simple representations of $\Pi_\qv$ so that the isomorphism class $[M]$ of $M$ belongs to $\scrS$ and $M,N$ are not 
		isomorphic. Then 	
		\begin{align}
			\Ext^1_{\Pi_\qv}(M,N) =0\quad \text{and} \quad \Ext^1_{\Pi_\qv}(N,M) =0\ .
		\end{align}
	\end{enumerate}
\end{lemma} 

\begin{remark}\label{rem:nilpotency-extensions} 
	Under the assumptions of Lemma~\ref{lem:nilpotency}--\eqref{item:nilpotency-3}, we have
	\begin{align}
		\Ext^k_{\Pi_\qv}(M,N) =0\quad \text{and} \quad \Ext^k_{\Pi_\qv}(N,M) =0
	\end{align}
	for all $k\in \Z$. Indeed, since $M,N$ are non-isomorphic simple modules, we get
	\begin{align}
		\Ext^0_{\Pi_\qv}(M,N) =0\quad \text{and} \quad \Ext^0_{\Pi_\qv}(N,M) =0\ .
	\end{align}
	Since the abelian category $\modPi$ is 2-Calabi-Yau (cf.\ \cite[Lemma~2.18]{SY13}), by Serre duality, the above equation implies the vanishing result for $k=2$ as well. 
\end{remark} 

For any $S\in \scrS$, we define the full subcategory $\catmod_S(\Pi_\qv)$ of $\modPi$ which consists of finite-dimensional representations of $\Pi_\qv$ whose composition factors consist only of $S$. The previous lemma leads to the following. 
\begin{proposition}[{\cite[Proposition~2.32]{SY13}}]\label{prop:nilpotency-decomposition}
	There is an orthogonal decomposition of abelian categories 
	\begin{align}
		\modPi \simeq \catmod_{\scrS}(\Pi_\qv)\oplus \nilpPi \ ,
	\end{align}
	where
	\begin{align}
		\catmod_{\scrS}(\Pi_\qv)\coloneqq \bigoplus_{S\in \scrS} \catmod_S(\Pi_\qv)\ .
	\end{align}
\end{proposition} 

\begin{remark}\label{rem:dimension-vector-skyscraper}
	Let $p\in \rsv$ be a closed point. Note that $\scrO_p$ belongs to the abelian category of properly supported perverse coherent sheaves $\catPps(\rsv/\ssv)$. Therefore $\tau(\scrO_p)$ is a finite-dimensional representation of $\Pi_\qv$. Since each direct summand of the local projective generator $\calP$, defined in \eqref{eq:projgenZN}, has rank $\mathsf{rk}(\calE_i)= r_i$, for $0\leq i \leq e$, the dimension vector of $\tau(\scrO_p)$ is $\delta$.
\end{remark}

\begin{lemma}\label{lem:nilpotency-mckay-I} 
	Let $p\in \rsv$ be a closed point. 
	\begin{enumerate}\itemsep0.2cm
		\item \label{item:nilpotency-mckay-I-1} If $p \in \rsv \smallsetminus C_{\red}$, then $\tau(\scrO_p)$ belongs to $\scrS$.
		
		\item  \label{item:nilpotency-mckay-I-2} If $p \in C_{\red}$, then $\tau(\scrO_p)$ is nilpotent.
	\end{enumerate}
\end{lemma} 

\begin{proof}
	Suppose that $p \in C_i$ for some $1\leq i \leq e$. Then there is a an exact sequence 
	\begin{align}
		\scrO_{C_i} \longrightarrow \scrO_p \longrightarrow \scrO_{C_i}(-1)[1] 
	\end{align}
	in $\catPps(\rsv/\ssv)$. By applying $\tau$, one obtains a surjective morphism $\tau(\scrO_p) \to \sigma_i$ with nontrivial kernel. This shows that $\tau(\scrO_p)$ is not simple. Then Lemma~\ref{lem:nilpotency}--\eqref{item:nilpotency-2} shows that $\tau(\scrO_p)$ is nilpotent. 
	
	Now, suppose that $p \in \rsv \smallsetminus C_\red$. We have to show that $\tau(\scrO_p)$ is simple. Let us assume that it is not. Then it must be nilpotent by Lemma~\ref{lem:nilpotency}--\eqref{item:nilpotency-2}, therefore its composition factors belong to $\{\sigma_0, \cdots, \sigma_e\}$, by Remark~\ref{rem:nilpotency}. Since all objects $\tau^{-1}(\sigma_i)$ are set-theoretically supported in $C_{\red}$, this is a contradiction.
\end{proof}

\begin{lemma}\label{lem:nilpotency-mckay-II} 
	Suppose that $M$ is a simple representation of $\Pi_\qv$ of dimension vector $\delta$. Then $\tau^{-1}(M) \simeq \scrO_p$ for 
	some closed point $p\in \rsv$ where 
	\begin{enumerate}\itemsep0.2cm
		\item $p\in C_{\red}$ if $M$ is nilpotent, or
		\item $p \in \rsv \smallsetminus C_{\red}$ if $M$ is simple.
	\end{enumerate}
\end{lemma} 

\begin{proof}
	Since $M$ is simple, it either belongs to $\tau(\catPps^0(\rsv/\ssv))$ or to $\tau(\catPps^{-1}(\rsv/\ssv))$. Since there are no coherent sheaves on $\rsv$ of dimension vector $-\delta$, we get that $M$ belongs to $\tau(\catPps^0(\rsv/\ssv))$. But then $\calF\coloneqq \tau^{-1}(M) \in \catCoh(\rsv)$ is of dimension vector $\delta$, hence by Lemma~\ref{lem:dimension-vector} satisfies $\ch_1(\calF)=0$ and $\chi(\calF)=1$, from which we conclude that $\calF \simeq \scrO_p$ for some closed point $p \in\rsv$.  Then the claim follows from Lemmas~\ref{lem:nilpotency}--\eqref{item:nilpotency-2} and \ref{lem:nilpotency-mckay-I}.  
	\end{proof}

\begin{lemma}\label{lem:nilpotency-mckay-III} 
	Let $1 \leq i \leq e$ and let $\calL\subset \scrO_{C_i}(-1)$ be an arbitrary subsheaf. Then $\calL[1]$ belongs to $\catPps(\rsv/\ssv)$ and $\tau(\calL[1])$ is nilpotent. 
\end{lemma}

\begin{proof}
	The fact that $\calL[1]$ belongs to $\catPps(\rsv/\ssv)$ follows from Remark~\ref{rem:torsion-pair-ps}. Now, we prove that $\tau(\calL[1])$ is nilpotent. First, one has an exact sequence in $\catCohps(\rsv)$ 
	\begin{align}
		0\longrightarrow \calL \to \scrO_{C_i}(-1) \longrightarrow \calQ \longrightarrow 0
	\end{align}
	in $\catCohps(\rsv)$, where $\calQ$ is a zero-dimensional sheaf with $\mathsf{Supp}(\calQ) \subset C_i$. This yields the exact sequence 
	\begin{align}
		0\longrightarrow \calQ\longrightarrow \calL[1] \longrightarrow \scrO_{C_i}(-1)[1] \longrightarrow 0 
	\end{align}
	in $\catPps(\rsv/\ssv)$. Therefore, one has an exact sequence 
	\begin{align}
		0\longrightarrow \tau(\calQ) \longrightarrow \tau(\calL[1]) \longrightarrow \sigma_i \longrightarrow 0
	\end{align}
	in $\modPi$. Since $\nilpPi$ is a Serre category, thus closed by extensions, and $\tau(\calQ)$ can be realized as an extension of skyscraper sheaves $\scrO_p$ with $p\in C_{\red}$, $\tau(\calQ)$ is a nilpotent representation by Lemma~\ref{lem:nilpotency-mckay-I}. Therefore, $\tau(\calL[1])$ is nilpotent as well.
\end{proof}

\begin{proposition}\label{prop:nilpotency-mckay-IV} 
	Let $E$ be a perverse coherent sheaf, which is set-theoretically supported on $C$. Then $\tau(E)$ is a nilpotent finite-dimensional representation of $\Pi_\qv$. 
\end{proposition} 

\begin{proof}
	One has the canonical exact sequence 
	\begin{align}
		0\longrightarrow \calH^{-1}(E)[1] \longrightarrow E \longrightarrow \calH^0(E) \longrightarrow 0 
	\end{align}
	in $\catPps(\rsv/\ssv)$. Remark~\ref{rem:torsion-pair-ps} shows that $\calH^{-1}(E)$ belongs to smallest Ext closed subcategory of $\catCohps(\rsv)$ generated by subsheaves of $\scrO_{C_i}(-1)$, for $1\leq i \leq e$. Then Lemma~\ref{lem:nilpotency-mckay-III} implies that $\tau(\calH^{-1}(E)[1])$ is nilpotent. Therefore it suffices to show that $\tau(\calH^0(E))$ is also nilpotent, since the category $\nilpPi$ is Ext closed.
	
	By Proposition~\ref{prop:nilpotency-decomposition} one has a direct sum decomposition 
	\begin{align}
		\tau(\calH^0(E))\simeq M_1 \oplus M_2
	\end{align}
	where $M_1$ is belongs to $\modPii{\scrS}$ and $M_2$ is nilpotent. This yields a direct sum decomposition 
	\begin{align}
		\calH^0(E) \simeq \tau^{-1}(M_1) \oplus \tau^{-1}(M_2) 
	\end{align}
	in $\catPps(\rsv/\ssv)$. Since all the simple composition factors of $M_1$ belong to $\scrS$, Lemma~\ref{lem:nilpotency-mckay-II} implies that $\tau^{-1}(M_1)$ has a filtration in $\catPps(\rsv/\ssv)$ so that all successive quotients are zero-dimensional sheaves supported in the complement $\rsv\smallsetminus C_{\red}$. Hence $\tau^{-1}(M_1)$ is also a zero-dimensional sheaf supported in the complement $\rsv\smallsetminus C_{\red}$. Since $E$ is set-theoretically supported on $C$, also $\calH^0(E)$ is set-theoretically supported on $C$, hence $M_1$ must be identically zero. Thus $\tau(\calH^0(E))\simeq M_2$ is nilpotent. 
\end{proof}

The next goal is to prove the converse of Proposition~\ref{prop:nilpotency-mckay-IV}. First note the following:
\begin{lemma}\label{lem:nilpotency-mckay-V} 
	Let $E\in \catPps^0(\rsv/\ssv)$. Suppose that $E$ has a filtration 
	\begin{align}
		0\subset E_1 \subset \cdots \subset E_\ell = E
	\end{align}
	in $\catPps(\rsv/\ssv)$ so that all successive quotients belong to $\{ \scrO_{C_1}(-1)[1], \ldots,  \scrO_{C_e}(-1)[1], \scrO_{C}\}$ up to isomorphism. Then $\mathsf{Supp}(E) \subset C$. 
\end{lemma} 

\begin{proof}
	We will first prove by descending induction on $1\leq i \leq \ell$ that all subobjects $E_i\subset E$, for $1\leq i \leq \ell$, belong to $\catPps^0(\rsv/\ssv)$. Set $F_\ell \coloneqq E_\ell/E_{\ell-1}$. Then one has two cases:
	\begin{enumerate}[label=(\roman*)]\itemsep0.2cm
		\item Suppose that $F_\ell \simeq \scrO_C$. Since $E_\ell=E$ has amplitude $[0, 0]$ by assumption, the long exact sequence associated to 
		\begin{align}
			0\longrightarrow E_{\ell-1}\longrightarrow E_{\ell} \longrightarrow F_\ell\longrightarrow 0
		\end{align}
		yields $\calH^{-1}(E_{\ell-1})=0$. 
		
		\item Suppose that $F_\ell\simeq \scrO_{C_i}(-1)[1]$ for some $1\leq i \leq e$. Then one has an exact triangle 
		\begin{align}
			\scrO_{C_i}(-1) \longrightarrow E_{\ell-1} \longrightarrow E_\ell 
		\end{align}
		in $\catDb(\catCoh(\rsv))$. This implies again that $\calH^{-1}(E_{\ell-1})=0$. 
	\end{enumerate}
	The inductive step is completely analogous.
	
	Next we prove by increasing induction on $1\leq i \leq \ell$ that $\mathsf{Supp}(E_i)\subset C$. Since $E_1$ is a coherent sheaf, by assumption one must have $E_1 \simeq \scrO_C$. Set $F_2=E_2/E_1$. One has two cases:
	\begin{enumerate}[label=(\roman*)]\itemsep0.2cm
		\item $F_2 \simeq \scrO_{C}$. Then, clearly, $\mathsf{Supp}(E_2)\subset C$. 
		
		\item $F_2 \simeq \scrO_{C_i}(-1)[1]$ for some $1\leq i \leq e$. Then one has the following exact triangle 
		\begin{align}
			\scrO_{C_i}(-1)\to E_1\to E_2 
		\end{align}
		in $\catDb(\rsv)$. Since all objects have amplitude $[0, 0]$ this is an exact sequence in $\catCohps(\rsv)$. Using the induction hypothesis, this implies that $\mathsf{Supp}(E_2)\subset C$. 
	\end{enumerate}
	The inductive step is completely analogous. 
\end{proof}

\begin{proposition}\label{prop:nilpotency-mckay-VI} 
	Let $M$ be a finite-dimensional nilpotent representation of $\Pi_\qv$. Then $\tau^{-1}(M)$ is set-theoretically supported on $C$. 
\end{proposition} 

\begin{proof}
	Let $E=\tau^{-1}(M)$. By Remark~\ref{rem:torsion-pair-ps}, $\calH^{-1}(E)$ is set-theoretically supported on $C$. Then Proposition~\ref{prop:nilpotency-mckay-IV} shows that $M_1\coloneqq \tau(\calH^{-1}(E)[1])$ is a nilpotent representation. Moreover, one has an exact sequence 
	\begin{align}
		0\longrightarrow M_1 \longrightarrow M \longrightarrow M_2 \longrightarrow 0 \ ,
	\end{align}
	where $M_2\simeq \tau(\calH^0(E))$. Since the subcategory of nilpotent representations is abelian, $M_2$ is also nilpotent. In particular, $M_2$ admits a filtration so that all successive quotients belong to  $\{\sigma_0, \ldots, \sigma_e\}$ up to isomorphism. This implies that $\calH^0(E)$ admits a filtration in $\catPps(\rsv/\ssv)$ so that all successive quotients belong to $\{ \scrO_{C_1}(-1)[1], \ldots,  \scrO_{C_e}(-1)[1], \scrO_{C}\}$ up to isomorphism. Then Lemma~\ref{lem:nilpotency-mckay-V} shows that also $\calH^0(E)$ is set-theoretically supported on $C$. 
\end{proof}

The combination of Propositions~\ref{prop:nilpotency-mckay-IV} and \ref{prop:nilpotency-mckay-VI} yields a proof of Theorem~\ref{thm:nilpotent}.

\section{Braid group actions on $\catDbps(\catCoh(\rsv))$}\label{sec:braid-group-derived-categories}

In this section, we use the McKay correspondence to transfer the well-known action of the extended affine braid group $B_\sfex$ on $\catDbps(\catMod(\qv))$ to an action on $\catDbps(\catCoh(\rsv))$ by means of spherical twist functors associated to the objects $\scrO_C$ and $\scrO_{C_i}(-1)$ for $i=1, \ldots, e$. Our main result is the identification of the autoequivalence $\R S_{\Llambda}$ corresponding to a translation in $B_\sfex$ with the autoequivalence given by tensoring with line bundles $\calL_{\Llambda}$. 

\subsection{Recollection on braid group actions on $\catDbps(\ModPi)$}

The theory of derived reflection functors has been introduced in \S\ref{subsec:reflection-functors}. Now, we recall how to use them to define a braid group actions on $\catDbps(\ModPi)$. 
\begin{theorem}[{\cite[Theorem~6.6-(2)]{IR08}}]\label{thm:twist}
	The assignment $T_i\mapsto \LL S_i$ for $i\in I$ defines an action of the affine braid group $B_\sfaf$ on $\catDbps(\ModPi)$. 
\end{theorem}
Note that this statement has been generalized to any quiver without edge-loops, see \cite[Propositions~III.1.4 and III.1.5, and Theorem~III.1.9]{BIRS09}.

Now, we shall introduce an equivalent way to formulate the affine braid group action via \textit{twist functors}. Let $\scrC$ be a $\C$-linear 2-Calabi-Yau triangulated category. 
\begin{definition}
	Let $x$ be an object in $\scrC$. The \textit{twist functor} $\bfT_x$ is the endofunctor of $\scrC$ acting as 
	\begin{align}\label{eq:twist-functor} 
		\begin{tikzcd}[ampersand replacement=\&]
			\bfT_x(y) \coloneqq \mathsf{Cone}\Big( \R\Hom(x, y)\otimes^\LL_\C x) \ar{r}{\ev} \& y\Big)
		\end{tikzcd}\ ,
	\end{align}
	while the \textit{dual twist functor} $\bfT'_x$ is the endofunctor of $\scrC$ acting as
	\begin{align}\label{eq:dual-twist-functor}
		\begin{tikzcd}[ampersand replacement=\&]
			\bfT'_x(y) \coloneqq \mathsf{Cone}\Big(y  \ar{r}{\ev} \& \R\Hom(y, x)^\vee \otimes^\LL_\C x\Big)[-1]
		\end{tikzcd}\ . 
	\end{align}
\end{definition}

\begin{definition}\label{def:spherical}
	Let $x$ be an object in $\scrC$. We say that $x$ is \textit{(2-)spherical} if we have 
	\begin{align}
		\R \Hom(x,x)=\C\oplus\C[-2] \ . \tag*{\qedhere} 
	\end{align}
\end{definition}
According to \cite[Proposition~2.10]{ST01}, if $x$ is 2-spherical, the functor $\bfT_x$ is an equivalence, with quasi-inverse $\bfT'_x$.

Let $\scrC=\catDbps(\ModPi)$. Then, for any $i\in I$, the simple $\Pi_\qv$-module  $\sigma_i$ with dimension $\alpha_i$ is spherical in $\catDbps(\ModPi)$.
\begin{lemma}\label{lem:twistPi} 
	One has $\bfT'_{\sigma_i}(\Pi_\qv) \simeq I_i$.
\end{lemma} 

\begin{proof}
	Note that  
	\begin{align}
		\Ext_{\Pi_\qv}^k(\Pi_\qv, \sigma_i) \simeq \begin{cases}
			\C & \text{for } k=0 \ , \\
			0 & \text{otherwise} \ .
		\end{cases}
	\end{align}
	Then the evaluation morphism in Formula~\eqref{eq:dual-twist-functor} reduces to 
	\begin{align}
		\begin{tikzcd}[ampersand replacement=\&]
			\Pi_\qv \ar{r}{\ev} \& \R\Hom(\Pi_\qv,\sigma_i)^\vee \otimes_\C \sigma_i \simeq \sigma_i 
		\end{tikzcd}
	\end{align}
	which coincides with the canonical projection. Hence $\bfT'_{\sigma_i}(\Pi_\qv) \simeq I_i$. 
\end{proof}
The above lemma together with \cite[Theorem~6.14]{IR08} yield the following.
\begin{proposition}\label{prop:tilt_twist_thm} 
	The functors $\bfT'_{\sigma_i}$ and $\LL S_i$ are isomorphic for any $i\in I$. Hence, also the functors $\bfT_{\sigma_i}$ and $\R S_i$ are isomorphic for any $i\in I$.
\end{proposition} 

We now extend the action of $B_\sfaf$ on $\catDbps(\ModPi)$ to an action of $B_\sfex$. Recall that $\Gamma$ denotes the group of outer automorphisms of $\qv$. Any automorphism $\pi\in\Gamma$ gives rise to an algebra automorphism of the algebra $\Pi_\qv$, hence to an autoequivalence $\pi$ of the category $\ModPi$ such that $\pi(\sigma_i)=\sigma_{\pi(i)}$ for each vertex $i\in I$. This yields an action of $\Gamma$ on $\catDbps(\ModPi)$. 
\begin{corollary}\label{cor:twist}
	The above actions of $B_\sfaf$ and $\Gamma$ give rise to an action of the extended affine braid group $B_\sfex$ on $\catDbps(\ModPi)$. 
\end{corollary}

\begin{remark}
	Using the canonical lift $W_\sfex=\Gamma \ltimes W_\sfaf \to B_\sfex$ we may thus define functors $\R S_w, \LL S_w$ for any element $w\in W_\sfex$. For each $\Llambda\in \coweightlatticefin$ and $i\in I$ we abbreviate 
	\begin{align}
		\R S_{\Llambda}& \coloneqq \R S_{t_{\Llambda}}\quad \text{and}\quad \R S_i \coloneqq \R S_{s_i}\ ,\\[2pt]
		\LL S_{\Llambda}& \coloneqq \LL S_{t_{\Llambda}}\quad \text{and}\quad \LL S_i \coloneqq \LL S_{s_i}\ . \tag*{\qedhere} 
	\end{align}
\end{remark}

\subsection{Action of the extended affine Braid group on $\catDbps(\catCoh(\rsv))$}\label{sec:action aff braid dbcoh}

Let $\Aut(\catDbps(\catCoh(\rsv)))$ be the group of exact self-equivalences of $\catDbps(\catCoh(\rsv))$. By Theorems~\ref{thm:tau} and \ref{thm:twist}, Proposition~\ref{prop:tilt_twist_thm}, and Corollary~\ref{cor:twist} there is a group homomorphism 
\begin{align}
	\rho_{\mathsf{IR}}\colon B_\sfex\longrightarrow \Aut(\catDbps(\catCoh(\rsv)))\ ,
\end{align}
such that
\begin{align}
	\rho_{\mathsf{IR}}(T_i)&= \bfT'_ {\scrO_{C_i}(-1)}\quad\text{and}\quad \rho_{\mathsf{IR}}(T_i^{-1})=\bfT_{\scrO_{C_i}(-1)}\quad\text{for } i\neq 0\ , \\ 
	\rho_{\mathsf{IR}}(T_0)&=\bfT'_{\scrO_C}\quad\text{and}\quad \rho_{\mathsf{IR}}(T_0^{-1})=\bfT_{\scrO_C}\ .
\end{align}
Note that the sheaves $\scrO_C$ and $\scrO_{C_i}(-1)$, for $i=1, \ldots, e$, admit an $A$-equivariant sheaf structure. Thus, the subgroup $\rho(B_\sfex)$ of $\Aut(\catDbps(\catCoh(\rsv)))$ centralizes $A$.

\begin{proposition}\label{prop:braid}	
	\hfill
	\begin{enumerate}\itemsep0.2cm
		\item \label{item:prop:braid-(1)} There is a group homomorphism $\rho_\sfex\colon B_\sfex\longrightarrow \Aut(\catDbps(\catCoh(\rsv)))$ such that
		\begin{align}\label{eq:rho-ex}
			\rho_\sfex(T_i)&= \bfT'_ {\scrO_{C_i}(-1)}\quad\text{and}\quad \rho_\sfex(T_i^{-1})=\bfT_{\scrO_{C_i}(-1)}\quad\text{for } i\neq 0\ , \\ 	\rho_\sfex(L_{\Llambda})&=\calL_{\Llambda}\otimes - \quad\text{for } \lambda\in\coweightlatticefin \ ,
		\end{align}
		with respect to the description of $B_\sfex$ given in Proposition~\ref{prop:characterization-extended-affine-braid-group}.
				
		\item \label{item:prop:braid-(2)}  We have $\rho_\sfex=\rho_{\mathsf{IR}}$.
	\end{enumerate}
\end{proposition}

\begin{proof}
	Let us first prove Claim \eqref{item:prop:braid-(1)} .
	
	Since the defining relations of the braid group $B_\sfaf$ follows from Theorem~\ref{thm:twist}, we must check the relations \eqref{eq:extended-B-(1)}, \eqref{eq:extended-B-(2)}, and \eqref{eq:extended-B-(3)} in Proposition~\ref{prop:characterization-extended-affine-braid-group}.
	
	Relation \eqref{eq:extended-B-(1)} is obvious. To prove \eqref{eq:extended-B-(2)}, note that the condition 
	$s_i(\Llambda)=\Llambda$ implies that $\calL_{\Llambda}\otimes\scrO_{C_i}(-1)=\scrO_{C_i}(-1)$. Hence, we have
	\begin{align}
		\bfT'_{\scrO_{C_i}(-1)}\circ(\calL_{\Llambda}\otimes-)=(\calL_{\Llambda} \otimes-)\circ\bfT'_{\scrO_{C_i}(-1)} \ .
	\end{align}
	
	To prove \eqref{eq:extended-B-(3)}, note that the condition $s_i(\Llambda)=\Llambda-\Lalpha_i$ implies that $\calL_{\Llambda}\otimes\scrO_{C_i}(-1)=\scrO_{C_i}$ and $\calL_{\Llambda-\Lalpha_i}\otimes\scrO_{C_i}(-1)=\scrO_{C_i}(-2)$. Hence, we have
	\begin{align}
		\bfT'_{\scrO_{C_i}(-1)}\circ(\calL_{\Llambda-\Lalpha_i}\otimes-)\circ\bfT'_{\scrO_{C_i}(-1)}
		&=\bfT'_{\scrO_{C_i}(-1)}\circ\bfT'_{\calL_{\Llambda-\Lalpha_i}\otimes\scrO_{C_i}(-1)}\circ(\calL_{\Llambda-\Lalpha_i}\otimes-)\\[2pt]
		&=\bfT'_{\scrO_{C_i}(-1)}\circ\bfT'_{\scrO_{C_i}(-2)}\circ(\calL_{\Llambda-\Lalpha_i}\otimes-)\\[2pt]
		&=(\calL_{\Lalpha_i}\otimes-)\circ(\calL_{\Llambda-\Lalpha_i}\otimes-)\\[2pt]
		&=\calL_{\Llambda}\otimes-\ ,
	\end{align}
	where the third equality follows e.g. from \cite[Lemma~4.15-(i)-(2)]{IU05}.
	
	\medskip
	
	Now, we concentrate on Claim \eqref{item:prop:braid-(2)}. It is enough to check that $\rho_\sfex(T_0)=\bfT'_{\scrO_C}=\rho_{\mathsf{IR}}(T_0)$. 
	
	We have $T_0=L_{\Lvarphi} T_{s_\varphi}^{-1}$, where $\varphi$ is the highest positive root in $\rootsetfin$ (cf.\ Remark~\ref{rem:braid-group}). Fix a reduced expression of $s_\varphi$ in $W$ of the form 	\begin{align}
		s_\varphi=s_{i_1}\cdots s_{i_\ell}s_js_{i_\ell}\cdots s_{i_1}\ .
	\end{align}
	We have $s_{i_1}\cdots s_{i_\ell}(\Lalpha_j)=\Lvarphi$ and
	\begin{align}
		T_0=T_{i_1}T_{i_2}\cdots T_{i_\ell}L_{\Lalpha_j}T_j^{-1}T_{i_\ell}^{-1}\cdots T_{i_2}^{-1}T_{i_1}^{-1}\ .
	\end{align}
	By applying \eqref{eq:rho-ex}, we have
	\begin{align}
		\rho_\sfex(L_{\Lalpha_j}T_j^{-1})=(\scrO_\rsv(-C_j)\otimes-)\circ \bfT_{\scrO_{C_j}(-1)}=(\bfT'_{\scrO_{C_j}}\circ\bfT'_{\scrO_{C_j}(-1)})\big) \circ \bfT_{\scrO_{C_j}(-1)} =\bfT'_{\scrO_{C_j}}\ ,
	\end{align}
	where the second equality follows from \cite[Lemma~4.15--(i)--(2)]{IU05}. Thus, we are reduced to prove that
	\begin{align}\label{eq:equality}
		\bfT'_{\scrO_C}=\bfT'_{\scrO_{C_{i_1}}(-1)}\circ\bfT'_{\scrO_{C_{i_2}}(-1)}\circ
		\cdots\circ \bfT'_{\scrO_{C_{i_\ell}}(-1)}\circ
		\bfT'_{\scrO_{C_j}}\circ
		\bfT_{\scrO_{C_{i_\ell}}(-1)}\circ\cdots\circ\bfT_{\scrO_{C_{i_2}}(-1)}\circ\bfT_{\scrO_{C_{i_1}}(-1)}\ .
	\end{align}
	If we set 
	\begin{align}
		\Phi\coloneqq \bfT'_{\scrO_{C_{i_1}}(-1)}\circ\bfT'_{\scrO_{C_{i_2}}(-1)}\circ
		\cdots\circ \bfT'_{\scrO_{C_{i_\ell}}(-1)}\ ,
	\end{align}
	the formula above becomes
	\begin{align}
		\bfT'_{\scrO_C}=\Phi\circ \bfT'_{\scrO_{C_j}}\circ\Phi^{-1} = \bfT'_{\Phi(\scrO_{C_j})}\ ,
	\end{align}
	where the second equality follows from \cite[Lemma~4.14--(i)]{IU05}. Thus, we are reduced to prove that 
	\begin{align}
		\scrO_C = \Phi(\scrO_{C_j})\ ,
	\end{align}
	or, equivalently, that
	\begin{align}\label{eq:form5}
		E\coloneqq\bfT_{\scrO_{C_{i_\ell}}(-1)}\circ\cdots\circ\bfT_{\scrO_{C_{i_2}}(-1)}\circ\bfT_{\scrO_{C_{i_1}}(-1)}(\scrO_C)= \scrO_{C_j}\ .
	\end{align}
	
	To prove the equality above, first, we shall show that $E$ is a perverse coherent sheaf concentrated in amplitude $[0,\, 0]$. Define
	\begin{align}
		\beta_r\coloneqq s_{i_\ell}s_{i_{\ell-1}}\cdots s_{i_{r+1}}(\alpha_{i_r})\in\Delta\quad\text{for}\quad r=0,\ldots,\ell\ .
	\end{align}
	Here, we set $i_0\coloneqq 0$. Since $s_{i_\ell}\dots s_{i_1}s_0$ is a reduced expression in $W_\sfaf$ with $i_\ell,\ldots,i_1\neq 0$, we deduce that
	$\beta_\ell,\ldots,\beta_0$ are positive roots and that 
	\begin{align}\label{eq:form6}
		\beta_0=s_{i_\ell}s_{i_{\ell-1}}\cdots s_{i_1}(\alpha_0)=\delta-\alpha_j \quad\text{and}\quad
		\beta_\ell=\alpha_{i_\ell}\ .
	\end{align}
	
	We fix a coweight $\Ltheta=\sum_{i\in I}\theta_i\,\Lomega_i$ in $\coweightlattice\otimes_\Z \R$ such that for all $i\neq 0$
	\begin{align}\label{eq:condition-1}
		\theta_i=(\Ltheta,\alpha_i)>0=(\Ltheta,\beta_0)\ .
	\end{align}
	We consider the Bridgeland's stability condition whose heart is $\calA\coloneqq \modPi$ and whose central charge is the group homomorphism
	\begin{align}
		Z\colon K_0(\modPi)=\Z I&\longrightarrow\C\ ,\\[4pt] 
		\bfd&\longmapsto -D(\bfd)+\imath R(\bfd)=-(\Ltheta,\bfd)+\imath(\Lrho,\bfd)\ .
	\end{align}
	Recall that the simple $\Pi_\qv$-module $\sigma_0$ with dimension $\alpha_0$ is a spherical object in $\catDbps(\ModPi)$, and that the $B_\sfaf$-action on the category $\catDbps(\ModPi)$ is given by twist functors as in Theorem~\ref{thm:twist}. Thus, the object
	\begin{align}
		M\coloneqq \bfT_{\sigma_{i_\ell}}\cdots\bfT_{\sigma_{i_1}}(\sigma_0)\in \catDbps(\ModPi)
	\end{align}
	is spherical. Note that the argument of $Z(\beta_0)$ is $\pi/2$ because $D(\beta_0)=0$, while the argument of $Z(\beta_r)$ is strictly bigger than $\pi/2$ for each $r>0$ because $D(\beta_r)>0$. Here, the argument is taken in $[0, \pi)$. Hence, we apply \cite[Proposition~4.2]{BDL21} and we deduce that $M$ belongs to the abelian subcategory $\modPi$, is of class $\beta_0$ and that $M$ is Bridgeland's $\Ltheta$-semistable, (hence slope $\Ltheta$-semistable of slope $(\Ltheta, \beta_0)=0$).
	
	The tilting functor $\tau$ in Theorem~\ref{thm:tau} is compatible with twist functors, i.e., $\bfT_{\tau(x)}\circ \tau \simeq \tau \circ \bfT_x$ for any spherical object $x\in \catDbps(\catMod(\rsv))$. Thus, $\tau(E)=M$. Hence, the object $E$ is perverse and $\tau(E)$ is $\Ltheta$-semistable with $(\Ltheta,\tau(\ch(E)))=(\Ltheta, \beta_0)=0$. Now, Conditions~\eqref{eq:condition-1} imply
	\begin{align}
		(\Ltheta, \delta)=(\Ltheta,\beta_0)+(\Ltheta, \alpha_j)= (\Ltheta, \alpha_j)>0\ .
	\end{align} 
	Thus, by \cite[Proposition~2.23]{DPS_McKay}, we obtain that $E\simeq \calH^0(E)\in \catPps^0(\rsv/\ssv)$ and moreover $E$ is a pure one-dimensional sheaf on $\rsv$. It remains to prove that $E$ is isomorphic to $\scrO_{C_j}$.
	
	We have $\tau(\ch(\scrO_{C_j}(-1)[1]))=\alpha_j$ and $\tau(\ch(\scrO_p))=\delta$ for each closed point $p\in C_j$ (cf.\ Remark~\ref{rem:dimension-vector-skyscraper}). Hence, we get
	\begin{align}
		\tau(\ch(\scrO_{C_j}))=\tau(\ch(\scrO_p))-\tau(\ch(\scrO_{C_j}(-1)[1]))=\delta-\alpha_j\ .
	\end{align}
	Thus, $[E]=[\scrO_{C_j}]$ as K-theory classes in $\sfK_0(\rsv)$. Therefore, $E$ is a pure sheaf scheme-theoretically supported on $C_j$, hence it is the pushforward of a degree zero line bundle on $C_j$. Since $C_j$ is rational, we must have $E=\scrO_{C_j}$.
\end{proof}

\begin{remark}
	A three-dimensional version of Proposition~\ref{prop:braid} is proved in \cite[Theorem~7.4]{HW_Stability}. As pointed out to us by M. Wemyss, one could ``slice down'' the 3-fold to get the surface we are considering and one should get an equivalent proof of Proposition~\ref{prop:braid}.
\end{remark}

\begin{remark}
	When $\qv=A_e^{(1)}$, the equality $\rho_\sfex(T_0)=\bfT'_{\scrO_C}$ may also be proved using \cite[Lemma~4.16]{IU05} and the following formula for $T_0$ in $B_\sfaf$
	\begin{align}
		T_0=T_1T_2\cdots T_{e-1}L_{\Lalpha_e}T_e^{-1}\cdots T_2^{-1}T_1^{-1}\ ,
	\end{align}
	which follows from the relations in \S\ref{subsec:braid-group-affine-quiver}.
\end{remark}

By Theorems~\ref{thm:twist} and \ref{thm:tau}, and Proposition~\ref{prop:braid}	the tilting equivalence $\tau\colon \catDb(\catCoh(\rsv)) \to \catDb(\ModPi)$ intertwines the $B_\sfex$-action on both sides. This yields the following.
\begin{proposition}\label{prop:braid-tilting-intertwiner} 
	Let $\lambda\in \coweightlatticefin$. The tilting equivalence $\tau$ intertwines the autoequivalence $\bfS_{\lambda}\coloneqq (\calL_{\lambda}\otimes -) $ of $\catDbps(\catCoh(\rsv))$ with the autoequivalence $\R S_{\lambda}$ of $\catDbps(\ModPi)$, i.e., 
	there is an equivalence of functors 
	\begin{align}
		\tau \circ \bfS_{\lambda} \simeq \R S_{\lambda} \circ \tau\ .
	\end{align}
\end{proposition}

\section{Semistability and braid group action}\label{sec:stab-braiding}

In this section, we provide a refinement of the tilting equivalences of Theorems~\ref{thm:tau} and \ref{thm:nilpotent}, which matches subcategories of semistable objects in $\catDbps(\catCoh(\rsv))$ and $\catDbps(\ModPi)$ for appropriate stability conditions. We also describe the behavior of the categories of semistable objects under the action of the extended affine braid group $B_\sfex$. We present this last result in the setting of Bridgeland stability conditions on a triangulated category, as in Part~\ref{part:COHA-stability-condition}.

\subsection{Torsion pairs with respect to $\omega$-stability}

Let $\omega$ be an integral polarization on $\rsv$. For any properly supported coherent sheaf $\calF$ on $\rsv$ with $\ch_1(\calF)\neq 0$, the \textit{$\omega$-slope of $\calF$} is given by
\begin{align}
	\mu_\omega(\calF) \coloneqq \frac{\chi(\calF)}{\omega \cdot \ch_1(\calF)}\ ,
\end{align}
(if $\calF$ is zero-dimensional, its $\omega$-slope is $+\infty$ by definition).
\begin{definition}
	Let $\calF$ be a pure properly supported coherent sheaf on $\rsv$ of dimension one. We say that $\calF$ is \textit{($\omega$-)(semi)stable}\footnote{Here we use the notation $(\leq)$ following \cite[Notation~1.2.5]{HL_Moduli}: if in a statement the word ``(semi)stable'' appears together with relation signs ``$(\leq)$'', the statement encodes in fact two assertions: one about semistable sheaves and relation signs ``$\leq$'' and one about stable sheaves and relation signs ``$<$'', respectively.}  if for any nonzero subsheaf $\calF'$ of $\calF$ we have $\mu_\omega(\calF')(\leq)\mu_\omega(\calF)$. 
\end{definition}

Recall that any coherent sheaf $\calF$ in $\catCohps(\rsv)$ has a canonical filtration 
\begin{align}\label{eq:HNfiltration-1}
	0 \subseteq \calF_0 \subset \calF_1 \subset \cdots \subset \calF_k=\calF
\end{align}
with $k \geq 0$, where 
\begin{itemize}\itemsep0.2cm
	\item $\calF_0\subset \calF$ is the maximal zero-dimensional subsheaf of $\calF$, and 
	\item for any $i >0$, $\calF_i/\calF_{i-1}$ is $\omega$-semistable and 
	\begin{align}
		\mu_\omega(\calF_1/\calF_0) > \mu_\omega(\calF_2/\calF_1) > \cdots > \mu_\omega(\calF_k/\calF_{k-1}) \ ,
	\end{align}
\end{itemize} 
We call the filtration \eqref{eq:HNfiltration-1} the \textit{Harder-Narasimhan filtration} of $\calF$. The successive quotients $\calF_s/\calF_{s-1}$, for $0\leq s\leq k$ are the \textit{Harder-Narasimhan factors} of $\calF$. For any properly supported coherent sheaf $\calF$ we set 
\begin{align}
	\mu_{\omega\textrm{-}\mathsf{min}}(\calF) \coloneqq \mu_{\omega} (\calF_k/\calF_{k-1}) \quad\text{and}\quad 
	\mu_{\omega\textrm{-}\mathsf{max}}(\calF) \coloneqq \begin{cases}\mu_{\omega}(\calF_1) & \text{if}\,\calF_0=0 \\ +\infty & \text{if}\, \calF_0 \neq 0\end{cases} \ .
\end{align}

\begin{definition}
	Let $\kappa\subset\Q$ be an interval. We denote by $\catCohps(\rsv)^{\kappa}$ the full subcategory of $\catCohps(\rsv)$ consisting of those properly supported sheaves on $\rsv$ for which the $\omega$-slopes of all the HN-factors belong to $\kappa$.
\end{definition}

\begin{notation}
	For any $\ell\in \Q$, we set $\catCohps(\rsv)^{\ell}\coloneqq \catCohps(\rsv)^{\{\ell\}}$, $\catCohps(\rsv)^{>\ell}\coloneqq \catCohps(\rsv)^{(\ell,+\infty)}$, and $\catCohps(\rsv)^{\leqslant\ell}\coloneqq \catCohps(\rsv)^{(-\infty,\ell]}$.
\end{notation}

The following is clear (cf.\ \cite[Example~A.4--(2)]{AB_Moduli_Bridgeland}). 
\begin{lemma} 
	For any $\ell\in \Q$, $\big(\catCohps(\rsv)^{>\ell}, \catCohps(\rsv)^{\leqslant\ell}\big)$ is a torsion pair of $\catCohps(\rsv)$, which is open in the sense of Lieblich \cite[Definition~A.2]{AB_Moduli_Bridgeland}. 
\end{lemma}

\subsection{Semistability and tilting}\label{subsec:semistability-tilting}

The goal of this section is the proof of the following result. Let $\Ltheta = \sum_{i \in I} \Ltheta_i \Lomega_i\in \coweightlattice$ be a coweight such that
\begin{align}\label{eq:Lthetaconditions}
	\Ltheta_i >0\ \text{ for }\ 1\leq i \leq e \quad \text{and} \quad \Ltheta_0=-\sum_{i=1}^{e} r_i \Ltheta_i\ .
\end{align}
\begin{remark}
	For $\Ltheta\in \coweightlattice$ satisfying \eqref{eq:Lthetaconditions} we get $\langle \Ltheta, \delta\rangle=0$. We denote by $\Lthetafin\coloneqq \sum_{i\neq 0} \Ltheta_i \Lomega_i \in \coweightlatticefin$ the finite coweight canonically associated to $\Ltheta$.
\end{remark}
Let $\omega \in \Pic_\Q(\rsv)$ be a rational polarization on $\rsv$ such that
\begin{align}
	\Ltheta_i= \langle \omega, C_i\rangle 
\end{align}	
for $i=1, \ldots, e$.

\begin{theorem}\label{thm:tilting-stability}
	Under the above assumptions, the tilting functor \eqref{eq:tau} yields  equivalences of categories 
	\begin{align}
		\catCohps^{>0}(\rsv)&\simeq \catPps^0(\rsv/\ssv)\simeq \catmod^{\leqslant 0}(\Pi_\qv) \ ,\\[2pt]
		\catCohps^{\leqslant 0}(\rsv)[1]&\simeq \catPps^{-1}(\rsv/\ssv)\simeq \catmod^{> 0}(\Pi_\qv) \ ,
	\end{align}
	for which the first equivalence maps the subcategory of zero-dimensional sheaves onto $\catmod^0(\Pi_\qv)$.
\end{theorem}

Before proving the theorem, let us start by recalling \cite[Corollary~2.17]{DPS_McKay} applied to the above $\Ltheta$:
\begin{lemma}\label{thetadeglemm} 
	Let $E$ be an object of $\catPps(\rsv,\ssv)$. Then one has 
	\begin{align}\label{eq:theta-deg}
		\deg_{\Ltheta}(E) &= \omega \cdot \big( \ch_1(\calH^{-1}(E)) -\ch_1(\calH^0(E))\big)\ ,\\[2pt]
		d_0(E) &=  \dim\, H^0(\calH^0(E)) + \dim\, H^1(\calH^{-1}(E))\ .
	\end{align} 
\end{lemma}

\begin{proposition}\label{prop:zeta-stability} 
	Let $E$ be a $\Ltheta$-semistable object of $\catPps(\rsv,\ssv)$ with $\deg_{\Ltheta}(E) \leq 0$. Then $E$ belongs to $\catPps^0(\rsv/\ssv)$.
\end{proposition} 

\begin{proof}
	Suppose $E$ is $\Ltheta$-semistable and $\calH^{-1}(E) \neq 0$. Using Formula~\eqref{eq:theta-deg},  
	\begin{align}
		\deg_{\Ltheta}(\calH^{-1}(E)[1]) = \omega \cdot \ch_1(\calH^{-1}(E))\ .
	\end{align}
	Since $E$ is $\Ltheta$-semistable with $\deg_{\Ltheta}(E)\leq 0$, this implies that $\omega\cdot \ch_1(\calH^{-1}(E))\leq 0$. This leads to a contradiction since $\calH^{-1}(E)$ must be pure of dimension one as shown in Remark~\ref{rem:torsion-pair-ps}. Therefore, $E$ belongs to $\catPps^0(\rsv/\ssv)$. 
\end{proof}

\begin{corollary}\label{cor:zeta-stability}
	There is an equivalence of categories
	\begin{align}
		\catPps^0(\rsv/\ssv)\simeq \catmod^{\leqslant 0}(\Pi_\qv)\ .
	\end{align}
\end{corollary} 

\begin{proof}
	Let $E$ be an object of $\catPps(\rsv,\ssv)\simeq \modPi$. We shall prove that $E$ belongs to $\catPps^0(\rsv/\ssv)$ if and only if $\mu_{\Ltheta\textrm{-}\mathsf{max}}(E) \leq 0$. 
	
	($\Rightarrow$) Since $\calH^{-1}(E)=0$, it follows immediately that $\calH^{-1}(F)=0$ for any sub-object $F\subset E$ in $\catPps(\rsv/\ssv)$. Hence Equation~\eqref{eq:theta-deg} yields 
	\begin{align}
		\deg_{\Ltheta}(F) = - \omega \cdot \ch_1(F) \leq 0\ .
	\end{align}
	
	($\Leftarrow$) Conversely, let 
	\begin{align}
		0 = E_0 \subset E_1 \subset \cdots \subset E_k=E
	\end{align}
	be the Harder-Narasimhan filtration of $E$ relative to $\Ltheta$-stability. Thus, $\mu_{\Ltheta}(E_i/E_{i-1}) \leq \mu_{\Ltheta\textrm{-}\mathsf{max}}(E)\leq 0$ for all $i$. Hence all $E_i/E_{i-1}$ belong to $\catPps^0(\rsv/\ssv)$ by Proposition~\ref{prop:zeta-stability}. This implies that $E$ belongs to $\catPps^0(\rsv/\ssv)$, since this category is closed under extensions.
\end{proof}

\begin{lemma}\label{lem:zeta-stability}
	Let $E\in \catPps^0(\rsv/\ssv)$ and set $\calF\coloneqq \calH^0(E)$.
	\begin{itemize}\itemsep0.2cm
		\item $\deg_{\Ltheta}(E) =0$ if and only if $\calF$ is a zero-dimensional sheaf. 
		\item If $\mu_{\Ltheta}(E)<0$, then $\mu_\omega(\calF)>0$. 
	\end{itemize}
\end{lemma}

\begin{proof}
	By using Equation~\eqref{eq:theta-deg}, we get $\deg_{\Ltheta}(E)= - \omega \cdot \ch_1(\calF)$, which implies the first statement. 
	
	Now, let us prove the second statement. Assume that $\ch_1(\calF)=\sum_{i=1}^e n_i C_i$. Note that
	\begin{align}\label{eq:slope}
		\mu_{\Ltheta}(E)&=\frac{\deg_{\Ltheta}(E)}{(\Lrho, \bfd(E))}=-\frac{\omega\cdot \ch_1(\calF)}{\chi(\calF)+\sum_{i=1}^e(r_i\, \chi(\calF)-n_i)}\\[4pt]
		&=\left(\frac{\sum_{i=1}^e n_i}{\omega\cdot \ch_1(\calF)}-\mu_\omega(\calF)(\Lrho, \delta)\right)^{-1} \ ,
	\end{align}
	where the second equality follows from Lemma~\ref{lem:dimension-vector} and Formula~\eqref{eq:theta-deg}. Now, the quantity $\sum_{i=1}^r n_i$ coincides with $D\cdot \ch_1(\calF)$, hence it is positive. If $\mu_{\Ltheta}(E)<0$, then 
	\begin{align}
		\mu_\omega(\calF)>\frac{1}{(\Lrho, \delta)}\frac{\sum_{i=1}^e n_i}{\omega\cdot \ch_1(\calF)}>0\ . \tag*{\qedhere} 
	\end{align}
\end{proof}

\begin{proposition}\label{prop:sheaf-stability} 
	There is an equivalence of categories $\catCohps^{> 0}(\rsv)\simeq \catPps^0(\rsv/\ssv)$.
\end{proposition} 

\begin{proof}
	Let $\calF$ be a coherent sheaf on $\rsv$ with proper support. We shall prove that $\calF$, seen as a one-term complex of amplitude $[0,0]$, belongs to $\catPps^0(\rsv/\ssv)$ if and only if $\mu_{\omega\textrm{-}\mathsf{min}}(\calF)>0$.
	
	If $\calF$ is zero-dimensional, the assertion follows from Proposition~\ref{prop:zero-dimensional-perverse}. Assume now that $\calF$ is not zero-dimensional. Consider the Harder-Narasimhan filtration of $\calF$:
	\begin{align}
		0 \subseteq \calF_0 \subset \calF_1 \subset \cdots \subset \calF_k=\calF\ .
	\end{align}
	
	Since $\calF\twoheadrightarrow \calF/\calF_{k-1}$ is a quotient in $\catCohps(\rsv)$, by Lemma~\ref{lem:perverse-torsion-pair} we have $\calF/\calF_{k-1}\in \catPps^0(\rsv/\ssv)$. Set $E_{k-1}\coloneqq \calF_{k-1}\in \catPps^0(\rsv/\ssv)$. Thus, $\mu_{\Ltheta\textrm{-}\mathsf{max}}(\calF/\calF_{k-1}) \leq 0$, hence $\mu_{\Ltheta}(\calF/\calF_{k-1}) \leq 0$. If $\deg_{\Ltheta}(\calF/\calF_{k-1})=0$, by Lemma~\ref{lem:zeta-stability}, the complex $E/E_{k-1}$ is zero-dimensional, hence we obtain $\mu_\omega( \calF/\calF_{k-1})=\infty$. If $\mu_{\Ltheta}(\calF/\calF_{k-1}) < 0$, by Lemma~\ref{lem:zeta-stability}, we get $\mu_\omega( \calF/\calF_{k-1})>0$. But then $\mu_{\omega\textrm{-}\mathsf{min}}(\calF)=\mu_\omega( \calF/\calF_{k-1})>0$.
	
	\medskip
	
	Let us prove the reverse implication. Assume that each $\calF_{i+1}/\calF_i$ is $\omega$-semistable of positive slope for $i=0, \ldots, k-1$. By \cite[Proposition~2.26]{DPS_McKay}, each quotient corresponds to an object in $\catPps^0(\rsv/\ssv)$. Note that the zero-dimensional sheaf $\calF_0$ corresponds canonically to an object in $\catPps^0(\rsv/\ssv)$. Since $\catPps^0(\rsv/\ssv)$ is closed under extensions, we have $\calF\in\catPps^0(\rsv/\ssv)$ as well. 
\end{proof}

\begin{proof}[Proof of Theorem~\ref{thm:tilting-stability}]
	Corollary~\ref{cor:zeta-stability} and Proposition~\ref{prop:sheaf-stability} yield that the tilting functor \eqref{eq:tau} yields an equivalence of categories 
	\begin{align}
		\catCohps^{>0}(\rsv)\simeq \catPps^0(\rsv/\ssv)\simeq \catmod^{\leqslant 0}(\Pi_\qv) \ ,
	\end{align}
	which maps the subcategory of zero-dimensional sheaves onto $\catmod^0(\Pi_\qv)$. 
	
	Now, recall  that $(\catPps^{-1}(\rsv, \ssv), \catPps^0(\rsv/\ssv))$ is a torsion pair of $\catPps(\rsv/\ssv)\simeq\modPi$. In addition, the pairs $(\catmod^{>0}(\Pi_\qv), \catmod^{\leqslant 0}(\Pi_\qv))$ and $(\catCohps^{>0}(\rsv), \catCohps^{\leqslant 0}(\rsv))$ are torsion pairs of $\modPi$ and $\catCohps(\rsv)$, respectively. By the orthogonality relation of torsion pairs, it follows that the tilting functor \eqref{eq:tau} yields an equivalence of categories
	\begin{align}
		\catCohps^{\leqslant 0}(\rsv)[1]\simeq \catPps^{-1}(\rsv/\ssv)\simeq \catmod^{> 0}(\Pi_\qv) \ . \tag*{\qedhere} 
	\end{align}
\end{proof}

\subsection{Bridgeland's stability conditions and braid group action}\label{sec:Bridgeland-stab-braid-group}

\subsubsection{Exact autoequivalences preserving phase order}

Let $\widetilde{\GL}^+(2, \R)$ be the universal cover of the group $\GL^+(2,\R)$ of real $2 \times 2$-matrices with positive determinant. It can be described explicitly in the following way. Let $\R\to S^1$ be the universal cover given by $\phi \mapsto e^{\imath\pi\phi}$; then $\widetilde{\GL}^+(2, \R)$ is the set of pairs $(T, f)$, where $T\in \GL^+(2,\R)$ and $f\colon \R\to \R$ is a choice of a lift of the induced action of $T$ on $S^1$. In particular, the map $f$ is increasing and such that $f(\phi+1)=f(\phi)+1$. 

By \cite[Lemma~8.2]{B07}, given a triangulated category $\scrD$, there is a right action of $\widetilde{\GL}^+(2, \R)$ on the space $\mathsf{Stab}(\scrD)$ of \textit{Bridgeland's pre-stability conditions}, in the sense of Definition~\ref{def:pre-stability-condition}, on $\scrD$ given as:
\begin{align}
	(Z, \calP)\cdot (T, f) \mapsto (Z_{(T,f)}, \calP_{(T,f)})\ ,
\end{align}
where $Z_{(T,f)}\coloneqq T^{-1}\circ Z$ and $\calP_{(T,f)}(\phi)\coloneqq  \calP(f(\phi))$ for any $\phi\in \R$. In particular, the action of $(T, f)$ ``permutes'' the categories of semistable objects with respect to the pre-stability condition $(Z, \calP)$.

Following again \cite[Lemma~8.2]{B07}, recall that there is also a left action of the group $\Aut(\scrD)$ of exact autoequivalences of $\scrD$ on $\mathsf{Stab}(\scrD)$ given by:
\begin{align}
	\Phi\cdot (Z, \calP)\coloneqq (Z_\Phi, \calP_\Phi)\ ,
\end{align}
where $Z_\Phi \coloneqq Z\circ \Phi^{-1}$ and $\calP_\Phi(\phi)\coloneqq \Phi(\calP(\phi))$ for any $\phi\in \R$. Thus, in general, the action of $\Phi$ does not preserve the family of categories of semistable objects. We introduce a criterium to ensure that an exact autoequivalence of $\scrD$ permutes the subcategories $\calP(\phi)$ of semistable objects. 

The semi-closed upper half plane, which is the union of the strict upper-half plane and the strictly negative real line:
\begin{align}
	H\coloneqq \HH\cup \R_{<0}= \{r e^{\imath \pi \nu}\,\vert\, r>0\, \text{ and }\, 0<\nu\leq 1\}\subset \C \ .
\end{align}
\begin{definition}\label{def:order-preserving}
	Let $\Phi\in \Aut(\scrD)$. We say that \textit{$\Phi$ preserves the phase order in $e^{\imath\pi\lambda}H$} for some $\lambda\in \R$ if $\Phi$ restricts to an autoequivalence of the abelian category $\calP((\lambda,\lambda +1])$ and there exists $T\in \GL^+(2, \R)$ such that 
	\begin{itemize}\itemsep0.2cm
		\item $T^{-1}(e^{\imath\pi\lambda}H)\subset e^{\imath\pi\lambda}H$.
		\item $Z_\Phi= T^{-1} \circ Z$ as morphisms from $\sfK_0(\scrD)$ to $\C$. \qedhere
	\end{itemize}
\end{definition}
Let $\Phi$, $\lambda$, and $T$ be as in the above definition. We denote by $f_T\colon \R \to \R$ the unique continuous bijection such that $f_T((\lambda, \lambda+1])\subset (\lambda, \lambda+1]$ and 
\begin{align}
	f_T(\phi)= \frac{1}{\pi}\arg\left(T\left(e^{\imath\pi \phi}\right)\right)
\end{align}
as maps on $S^1$. In particular, $(T, f_T)\in \widetilde{\GL}^+(2, \R)$. Hence, $f_T$ is increasing and $f_T(\phi+1)=f_T(\phi)+1$.

\begin{proposition}\label{prop:order-preserving}
	Let $\Phi$ be an exact autoequivalence of $\scrD$ that preserves the phase order in $e^{\imath\pi\lambda}H$ for some $\lambda\in \R$ and let $T$ be corresponding element in $\GL^+(2,\R)$. Then, for any $\phi \in \R$ and for any interval $I\subset \R$, the functor $\Phi$ restricts to equivalences 
	\begin{align}
		\begin{tikzcd}[ampersand replacement=\&]
			\Phi\colon \calP(\phi) \ar{r}{\sim}\&  \calP(f_T(\phi))
		\end{tikzcd}
		\quad\text{and}\quad 
		\begin{tikzcd}[ampersand replacement=\&]
			\Phi\colon \calP(I) \ar{r}{\sim}\&  \calP(f_T(I))
		\end{tikzcd}\ .
	\end{align}
	In particular, $(Z_\Phi, \calP_\Phi)=(Z_{(T, f_T)}, \calP_{(T, f_T)})$.
\end{proposition}

\begin{proof}
	Since $\Phi$ preserves the phase order in $e^{\imath \pi \lambda}H$, for any $E, F \in \calP((\lambda,\lambda +1])$, one has that $\arg(Z(F)) \leq \arg(Z(E))$ if and only if $\arg(Z(T(F))) \leq \arg(Z(T(E)))$. Since $\Phi$ is exact in $\calP((\lambda,\lambda +1])$, we deduce that the equivalence $\Phi$ maps semistable objects of phase $\phi \in (\lambda,\lambda+1]$ to semistable objects of phase $f_T(\phi)$. Since $\Phi$ is exact in $\scrD$, we deduce that $\Phi$ maps the category $\calP(I)$ to $\calP(f_T(I))$ for any interval $I\subset (\lambda, \lambda+1]$.
	
	Now, note that for any $I\subset (\lambda+k, \lambda+1+k]$ for some $k\in \Z$, we have $\Phi(\calP(I))\simeq \calP(f_T(I))$ by using the relation between slicing and the shift functor (see e.g. \cite[Definition~3.3]{B07}). Finally, for any $I\subset \R$, we have
	\begin{align}
		\Phi(\calP(I))&\simeq \langle \Phi(\calP(\phi))\, \vert\, \phi \in I\rangle\simeq \langle \Phi(\calP(I\cap (\lambda+k, \lambda+1+k]))\, \vert\, k\in \Z\rangle\\[2pt]
		&\simeq \langle \calP(f_T(I\cap (\lambda+k, \lambda+1+k]))\, \vert\, k\in \Z\rangle\simeq \calP(f_T(I)) \ .
	\end{align}
	Thus, the first part of the statement follows. The second part is then evident.
\end{proof}

For any $\phi\in\R$, let us by $\tau_\phi$ the $t$-structure on $\scrD$ such that $\calP( (\phi, +\infty))$ is the connective part of $\tau_\phi$. 
\begin{corollary}
	Under the same hypotheses of Proposition~\ref{prop:order-preserving}, $\Phi$ is a $t$-exact autoequivalence of $\scrD$ with respect to the $t$-structures $\tau_\phi$ and $\tau_{f_T(\phi)}$ for any $\phi\in \R$. 
\end{corollary}

\subsubsection{Stability condition for quivers and braid group action}\label{subsubsec:stability-condition-quiver}

Now, let $\Ltheta\in \coweightlattice$ be a coweight satisfying conditions~\eqref{eq:Lthetaconditions}, and let $\Ltheta_\sff \in \coweightlatticefin$ be the corresponding finite coweight. Consider the stability condition $Z_{\Ltheta}$ on $\modPi$ given by
\begin{align}
	Z_{\Ltheta}\colon \sfK_0(\modPi)\simeq\Z I &\longrightarrow \C\ ,\\[2pt]
	\bfd & \longmapsto -D_{\Ltheta}(\bfd)+\imath R(\bfd)=-(\Ltheta, \bfd)+\imath (\Lrho, \bfd) \ .
\end{align}
It satisfies the Harder-Narasimhan property (cf.\ \cite[Example~5.5]{B07}). Hence, by \cite[Proposition~5.3]{B07}, it corresponds to a pre-stability condition $(Z_{\Ltheta}, \calP_{\Ltheta}$) on $\catDbps(\ModPi)$. By construction, $\modPi = \calP_{\Ltheta}((0,1])$.

\begin{proposition}\label{prop:preserve-order}
	The exact autoequivalence $\R S_{\textrm{-}2\Lthetafin}$ of $\catDbps(\ModPi)$ preserves the phase order in $e^{-\imath\pi/2}H$.
\end{proposition}

\begin{proof}
	Tensoring by the line bundle $\calL_{\textrm{-}2\Ltheta}$ is an autoequivalence of $\catDbps(\catCoh(\rsv))\simeq \catDb(\catCohps(\rsv))$ which descends to an autoequivalence of the abelian category $\catCohps(\rsv)$. Via the tilting equivalence $\tau$ we obtain that the autoequivalence $\R S_{\textrm{-}2\Lthetafin}$ restricts to an autoequivalence of the abelian category $\calP_{\Ltheta}((-1/2,1/2])$ (cf.\ Proposition~\ref{prop:braid-tilting-intertwiner}). 
	
	Now, consider 
	\begin{align}
		T\coloneqq \begin{pmatrix}
			1 & 0\\
			2h & 1
		\end{pmatrix}\in \SL(2, \Z)\ ,
	\end{align}
	where $h$ is the Coxeter number introduced in Formula~\eqref{eq:Coxeter}. By Formula~\eqref{eq:t_lambda}, for any $\bfd\in \Z I\simeq \sfK_0(\modPi)$ we see that 
	\begin{align}
		Z_{\Ltheta}(\R S_{\textrm{-}2\Lthetafin}^{-1}(\bfd))&=Z_{\Ltheta}(\R S_{2\Lthetafin}(\bfd))=-(\Ltheta, \bfd)+2(\Ltheta, \bfd)(\Ltheta, \delta)+\imath (\Ltheta, \bfd) -2(\Ltheta, \bfd)(\Lrho, \delta)\imath\\
		&= T^{-1}(Z_{\Ltheta}(\bfd)) \ ,
	\end{align}
	where the last equality follows from the fact that $(\Ltheta, \delta)$ vanishes and $h=(\Lrho,\delta)$. Hence, 
	\begin{align}
		Z_{\Ltheta}\circ \R S_{\textrm{-}2\Lthetafin}^{-1}= T^{-1} \circ Z_{\Ltheta}\ .
	\end{align}
	Moreover, $T$ preserves the half plane $\R_{\geqslant 0}+\imath\R$. In particular, $T^{-1}(e^{-\imath\pi/2}\HH)\subset e^{-\imath\pi/2}\HH$. Thus, the claim follows.
\end{proof}

We abbreviate $f\coloneqq f_T \colon \R \to \R$. Note that the map $f$ preserves the interval $(-\frac{1}{2}, \frac{1}{2}]$. Propositions~\ref{prop:order-preserving} and \ref{prop:preserve-order} yield the following.
\begin{lemma}\label{lem:derived-reflection-equivalence}
	For any interval $I\subset \R$, the derived equivalence $\R S_{\textrm{-}2\Lthetafin}$ of $\catDbps(\ModPi)$ descends to an equivalence
	\begin{align}
		\begin{tikzcd}[ampersand replacement=\&]
			\R S_{\textrm{-}2\Lthetafin}\colon \calP_{\Ltheta}(I)\ar{r}{\sim}\& \calP_{\Ltheta}(f(I))
		\end{tikzcd}\ .
	\end{align}
\end{lemma}

Set 
\begin{align}\label{eq:kappa}
	\kappa_0\coloneqq (-\infty, 0]\quad\text{and} \quad \kappa_\ell\coloneqq (-1/(2h\ell), 0] \quad \text{for } \ell \in \N\ , \ell \geq 1\ .
\end{align}
Note that $\kappa_{\ell+1} \subset \kappa_\ell \subset \kappa_0$ for any $\ell\in \N, \ell\geq 1$. 

We also set
\begin{align}
	\nu_k\coloneqq f^k(0)\in\left(-\frac{1}{2},\frac{1}{2}\right]
\end{align}
for $k\in \Z$. 
\begin{remark}
	For any $\phi\in [0, 1/2)$, we have
	\begin{align}
		f^k(\phi)=\frac{1}{\pi}\arg\left(T^k\left(e^{\imath\pi \phi}\right)\right)\ .
	\end{align}
	Thus, we have
	\begin{align}\label{eq:nu}
		\nu_k=\frac{1}{\pi}\arg\big(T^k(1)\big)=\frac{1}{\pi}\arctan(2hk)
	\end{align}
	for any $k\in \Z$. Hence one has $\nu_k<\nu_{k+1}$ for any $k\in \Z$ while
	\begin{align}
		\inf\big\{\nu_k\,\vert\, k\in \Z\big\}=-\frac{1}{2} \quad \text{and} \quad \sup\big\{\nu_k\,\vert\, k\in \Z\big\}=\frac{1}{2}\ . \tag*{\qedhere} 
	\end{align}
\end{remark}

\begin{theorem}\label{thm:tau-slicing}
	Let $\Ltheta\in \coweightlattice$ be a coweight satisfying the conditions~\eqref{eq:Lthetaconditions} and let $\omega \in \Pic_\Q(\rsv)$ be a rational polarization on $\rsv$ such that $\theta_i= \langle \omega, C_i\rangle$ for $i=1, \ldots, e$.
	\begin{enumerate}\itemsep0.2cm
		\item \label{thm:tau-slicing-1} 
		The functor $\tau$ restricts to  equivalences of additive categories
		\begin{align}
			\begin{tikzcd}[ampersand replacement=\&, row sep=tiny]
				\tau \colon \catCohps(\rsv)^{>0}\ar{r}{\sim} \&\modPi^{\leqslant 0}= \calP_{\Ltheta}\left( \left(0, \frac{1}{2}\right] \right)  \\
				\tau\colon \catCohps(\rsv)^{\leqslant 0}[1]\ar{r}{\sim} \& \modPi^{>0}= \calP_{\Ltheta}\left( \left(\frac{1}{2}, 1\right]\right)
			\end{tikzcd}\ .
		\end{align}
		In particular, the functor $\tau$ restricts to an  equivalence of abelian categories
		\begin{align}
			\begin{tikzcd}[ampersand replacement=\&]
				\tau\colon \catCohps(\rsv) \ar{r}{\sim} \& \calP_{\Ltheta}\left( \left(-\frac{1}{2}, \frac{1}{2}\right] \right)
			\end{tikzcd} \ .
		\end{align}
		
		\item \label{thm:tau-slicing-2} For any $k\in \Z$, we have  equivalences of abelian categories
		\begin{align}
			\begin{tikzcd}[ampersand replacement=\&, row sep=tiny]
				(\R S_{\textrm{-}2\Lthetafin})^k \colon \modPi^{\leqslant 0} \ar{r}{\sim}\& \calP_{\Ltheta}\left(\left(\nu_{\textrm{-}k}, \frac{1}{2}\right]\right)\\
				(\R S_{\textrm{-}2\Lthetafin})^k \colon \modPi^{> 0} \ar{r}{\sim}\& \calP_{\Ltheta}\left(\left(\frac{1}{2},\nu_{\textrm{-}k} + 1\right]\right) \\
				(\R S_{\textrm{-}2\Lthetafin})^k \colon \modPi \ar{r}{\sim} \& \calP_{\Ltheta}\left(\left(\nu_{\textrm{-}k}, 1 + \nu_{\textrm{-}k}\right]\right)
			\end{tikzcd}\ .
		\end{align}
		
		\item \label{thm:tau-slicing-3} For $\ell\in \N$, with $\ell\geq 1$, we have $\modPi^{\kappa_\ell}=\calP_{\Ltheta}((\nu_\ell,1/2])$.
		
		\item \label{thm:tau-slicing-4} For any $k\in \Z$, the square of  equivalences
		\begin{align}
			\begin{tikzcd}[ampersand replacement=\&]
				\catCohps(\rsv)^{>2k} \ar{d}{\bfS_{\textrm{-}2\Lthetafin}} \ar{r}{\tau} \& \calP_{\Ltheta}\left(\left(\nu_k, \frac{1}{2}\right]\right) \ar{d}{\R S_{\textrm{-}2\Lthetafin}} \\
				\catCohps(\rsv)^{>2k-2}  \ar{r}{\tau} \& \calP_{\Ltheta}\left(\left(\nu_{k-1}, \frac{1}{2}\right]\right) 
			\end{tikzcd}
		\end{align}
		is commutative. 
	\end{enumerate}
\end{theorem}

\begin{proof}
	We start by proving Claim \eqref{thm:tau-slicing-1}. Recall that the slope associated to $Z_{\Ltheta}$ is 
	\begin{align}
		\mu_{\Ltheta}(\bfd)=\frac{D_{\Ltheta}(\bfd)}{R(\bfd)}=\frac{(\Ltheta, \bfd)}{(\Lrho, \bfd)}\ .
	\end{align}
	Then, an object $N \in \modPi$ is slope-semistable if and only if it is Bridgeland's semistable with respect to $Z_{\Ltheta}$, i.e., if and only if 
	\begin{align}\label{eq:slope-slicing}
		N \in \calP_{\Ltheta}(\lambda)\quad\text{with}\quad \lambda =\frac{1}{\pi} \arg\left(-\mu_{\Ltheta}(N)R(N)+\imath R(N)\right)\in (0,1]\ ,
	\end{align}
	for $\mu_{\Ltheta}(N)\in \R$, while $\lambda=1/2$ if $\mu_{\Ltheta}(N)=+\infty$. Here, $R(N)=R(\bfd)$ if $\bfd$ is the dimension vector of $N$.
	
	Now the first part of the claim follows from this observation together with Theorem~\ref{thm:tilting-stability}. Since 
	\begin{align}
		\modPi^{>0}[-1] = \calP_{\Ltheta}\left( \left(\frac{1}{2}, 1\right]\right)[-1]= \calP_{\Ltheta}\left( \left(-\frac{1}{2}, 0\right]\right)\ ,
	\end{align}
	and $\catCohps(\rsv)$ is the extension-closed subcategory of $\catDb(\catCohps(\rsv))$ generated by $\catCohps(\rsv)^{\leqslant 0}$ and $\catCohps(\rsv)^{> 0}$, we get that
	the functor $\tau$ restricts to an $A$-equivariant equivalence of abelian categories
	\begin{align}
		\begin{tikzcd}[ampersand replacement=\&]
			\tau\colon \catCohps(\rsv) \ar{r}{\sim}\& \calP_{\Ltheta}\left( \left(-\frac{1}{2}, \frac{1}{2}\right] \right) 
		\end{tikzcd}\ .
	\end{align}
	
	\medskip
	
	To prove Claim \eqref{thm:tau-slicing-2}, we apply Lemma~\ref{lem:derived-reflection-equivalence} successively $k$ times to the subintervals $(0, 1/2]$ and $(-1/2,0]$ of $(-1/2, 1/2]$. This gives the first two equivalences. Note that $\calP_{\Ltheta}( (1/2, 1]) = \calP_{\Ltheta}((-1/2, 0])[1]$ by definition of slicing. As $\modPi$ is the extension-closed subcategory of $\catDb(\modPi)$ generated by $\modPi^{\leqslant 0}$ and $\modPi^{> 0}$ we get the last equivalence.
	
	\medskip
	
	Claim \eqref{thm:tau-slicing-3} follows from Formulas~\eqref{eq:slope-slicing} and \eqref{eq:nu}.
	
	\medskip 
	
	Finally, note that 
	\begin{align}
		\mu_\omega(\bfS_{\textrm{-}2\Lthetafin}(\calF))=\mu_\omega(\calF\otimes \calL_{\textrm{-}2\Lthetafin})=\frac{\chi(\calF)-2\ch_1(\calF)\cdot \omega}{\ch_1(\calF)\cdot \omega}=\mu_\omega(\calF)-2\ .
	\end{align}
	Claim~\eqref{thm:tau-slicing-4} follows from this together with Proposition~\ref{prop:braid-tilting-intertwiner}.
\end{proof}

\subsection{Braid group action and the quotients of the Yangian and of the COHA}\label{subsec:quotients}

We shall now consider the action of the extended affine braid group $B_\sfex$ on $\catDbps(\ModPi)$, and the various quotients of $\Y^-_\qv$ and $\cohaqv^T$ associated to the image under these of the standard heart $\modPi$ or of the exact subcategory $\modPi^{\leqslant 0}$.

Let $\Ltheta \in\coweightlattice$ be again a coweight satisfying the conditions~\eqref{eq:Lthetaconditions}. For $\ell \in \N$, let $\kappa_\ell$ be as in Formula~\eqref{eq:kappa}. We define the following quotients of $\cohaqv^T$:
\begin{align}\label{eq:quotient-COHA-ell}
	\coha_{\Ltheta, (\ell)}^T\coloneqq \HBMbulletT(\dLambda_\qv^{\kappa_\ell})\ .
\end{align}
There is a sequence of projections
\begin{align}\label{eq:chain-coha}
	\coha_{\Ltheta,(0)}^{T} \twoheadrightarrow \cdots \twoheadrightarrow \coha_{\Ltheta,(\ell)}^{T} \twoheadrightarrow \coha_{\Ltheta,(\ell+1)}^{T} \twoheadrightarrow \cdots \ .
\end{align}	
Note that $\coha_{\Ltheta, (\ell)}^T$ maps surjectively to $\HBMbulletT(\dLambda_\qv^0)$  for any $\ell\in \N$.

Thanks to the isomorphisms of Corollary~\ref{cor:isom-PHI-for-quotients}, there is an isomorphic chain of projections
\begin{align}\label{eq:chain-yangian}
	\Y_{\Ltheta,(0)} \twoheadrightarrow \cdots \twoheadrightarrow \Y_{\Ltheta, (\ell)}\twoheadrightarrow \Y_{\Ltheta, (\ell+1)}\twoheadrightarrow \cdots\ .
\end{align}
where for $\ell\in \N$, we set
\begin{align}\label{eq:quotient-Yangian-ell}
	\Y_{\Ltheta, (\ell)}\coloneqq \Y_{\Ltheta, \kappa_\ell}\ .
\end{align}
Note that $\Y_{\Ltheta, (\ell)}$ maps surjectively to $\Y_{\Ltheta, \{0\}}$ for $\ell\in \N$.

All the maps featuring in the chains of projections \eqref{eq:chain-coha} and \eqref{eq:chain-yangian} are $\N \times \Z I$-graded. Moreover, there are canonical isomorphisms 
\begin{align}
	\begin{tikzcd}[ampersand replacement=\&]
		\Phi_{(\ell)}\colon \Y_{\Ltheta,(\ell)} \ar{r}{\sim} \& \coha_{\Ltheta,(\ell)}^T
	\end{tikzcd}\ .
\end{align}
By Theorem~\ref{thm:tau-slicing}--\eqref{thm:tau-slicing-2} and --\eqref{thm:tau-slicing-3} there is an equivalence of stacks
\begin{align}
	\begin{tikzcd}[ampersand replacement=\&]
		(\R S_{\textrm{-}2\Lthetafin})^{\ell} \colon \dLambda_\qv^{\kappa_\ell}  \ar{r}{\sim}\& \dLambda_\qv^{\kappa_0}
	\end{tikzcd} \ .
\end{align}
This equivalence induces, for each $\bfd \in \N I$, an isomorphism of vector spaces which we denote by the same symbol
\begin{align}
	\begin{tikzcd}[ampersand replacement=\&]
		(\R S_{\textrm{-}2\Lthetafin})^{\ell} \colon\big(\coha_{\Ltheta,(\ell)}^{T}\big)_{\bfd} \ar{r}{\sim}\& \big(\coha_{\Ltheta,(0)}^{T}\big)_{t_{\textrm{-}2\ell\Lthetafin }(\bfd)}
	\end{tikzcd}\ .
\end{align} 

Our next task is to get a better understanding of the restriction map $\coha_{\Ltheta, (\ell)}^{T} \to \coha_{\Ltheta, (\ell+1)}^{T}$ using the isomorphisms $\Phi_{(\ell)}$.
Let $\Gamma$ be the group of outer automorphisms of the underlying diagram of $\qv$. It acts by algebra automorphisms on $\Y_\qv$ in the obvious way. This yields an action of the extended affine braid group $B_\sfex$ on $\Y_\qv$ by algebra automorphisms.

For any $w\in W_\sfex$, recall from \S\ref{sec:braid-on-Yangian} that we define a linear operator $\overline T_w$ given by the composition
\begin{align}
	\overline T_w\colon
	\begin{tikzcd}[ampersand replacement=\&]
		\Y_\qv^-\arrow{r} \& \Y_\qv\arrow{r}{T_w} \&\Y_\qv\arrow{r}{\pr} \&\Y_\qv^-
	\end{tikzcd}\ ,
\end{align} 
where the last map is the projection \eqref{eq:projection}. Let $B^+_\sfex\subset B_\sfex$ be the \textit{positive extended affine braid monoid}, i.e., the submonoid generated by the elements $T_w$ with $w\in W_\sfex$. The following can be proved as Proposition~\ref{prop:truncated-action}.
\begin{proposition}\label{prop:truncated-action-affine}
	There exists a homomorphism $B^+_\sfex \to \End(\Y_\qv^-)$ of monoids, which sends $T_w$ to $\overline T_w$. It defines an action of the positive braid monoid $B^+_\sfex$ on $\Y_\qv^-$.
\end{proposition}	

We now come to the main result of this section. For future use, it will be important to allow arbitrary twistings $\Theta$ in the definition of the COHA multiplication.

\begin{theorem}\label{thm:Stheta=Ttheta}
	Let $\Theta\colon \Z I \times \Z I \to \Z/2\Z$ be a twist and assume that $t_{\Llambda}\cdot \Theta=\Theta$ for any $\Llambda \in \coweightlatticefin$.
	Then, the diagram
	\begin{align}\label{eq:Rtheta=Stheta_twisted}
		\begin{tikzcd}[ampersand replacement=\&, column sep=large]
			\coha_{\Ltheta,(2)}^{T,\Theta} \ar{r}{\R S_{\textrm{-}4\Lthetafin}} \& \coha_{\Ltheta,(0)}^{T,\Theta}\\
			\Y_{\Ltheta,(2)} \ar{u}{\Phi_{(2)}} \& \Y_{\Ltheta,(0)} \arrow[swap]{u}{\Phi_{(0)}}\\
			\Y_\qv^-  \ar{r}{\overline T_{\textrm{-}4\Lthetafin}} \ar{u} \& \Y_\qv^-  \ar{u}
		\end{tikzcd}
	\end{align}
	commutes. Here, we set $\coha_{\Ltheta,(2)}^{T,\Theta}\coloneqq \coha_{\Ltheta,(2)}^T$ and $\coha_{\Ltheta,(0)}^{T,\Theta}\coloneqq \coha_{\Ltheta,(0)}^T$.
\end{theorem}

\begin{proof} 
	Recall from Formula~\eqref{eq:def-of-Lambda-uv} that we defined, for every pair of Weyl group elements $u,v \in W_\qv$ an open substack $\dLambda^u_w \subset \dLambda_\qv$. When we have only $u$ (resp.\ $w$), we simply write $\dLambda^u$ (resp. $\dLambda_w$).
	
	Because $\R S_{\textrm{-}2\Lthetafin}$ restricts to an equivalence of categories 
	\begin{align}
		\R S_{\textrm{-}2\Ltheta}\colon \nilpPi^{\kappa_1}\longrightarrow \nilpPi^{\kappa_0}
	\end{align}
	there are open embeddings
	\begin{align}\label{eq:open-embeddings-stability-kappa}
		\dLambda_\qv^{\kappa_1}\subset \dLambda^{t_{\textrm{-}2\Lthetafin}} \quad\text{and}\quad \dLambda_\qv^{\kappa_0}\subset\dLambda_{t_{\textrm{-}2\Lthetafin}}\ .
	\end{align}
	Indeed, by Proposition~\ref{prop:torsion-pairs}, the stack $\dLambda_{t_{\textrm{-}2\Lthetafin}}$ parametrizes nilpotent representations $M$ satisfying $H^1(\R S_{\textrm{-}2\Ltheta}(M))=0$. Likewise, $\dLambda^{t_{\textrm{-}2\Lthetafin}}$ parametrizes nilpotent representations $M$ satisfying $H^{-1}(\R S_{2\Ltheta}(M))=0$. Taking the pushforward, we get the following commutative diagram
	\begin{align}\label{eq:diagram-C}
		\begin{tikzcd}[ampersand replacement=\&, column sep=large]
			\coha_{\Ltheta, (1)}^{T,\Theta}\ar{r}{\R S_{\textrm{-}2\Lthetafin}} \& \coha_{\Ltheta, (0)}^{T,\Theta}\\
			\tensor*[^{(t_{\textrm{-}2\Lthetafin})}]{\coha}{^T_\qv} \ar{r}{\R S_{\textrm{-}2\Lthetafin}} \ar{u} \& \coha_\qv^{T, (t_{\textrm{-}2\Lthetafin})} \ar{u}
		\end{tikzcd}\ ,
	\end{align}
	where the vertical maps are induced by the open restrictions. 
	
	Consider the reduced expression $t_{\textrm{-}2\Lthetafin}=s_{i_1}s_{i_2}\cdots s_{i_r}$. By Propositions~\ref{prop:derived-reflection-functors-factorization} and \ref{prop:derived-underived-reflection-functors}, the lower map in \eqref{eq:diagram-C}	factors as
	\begin{align}
		\begin{tikzcd}[ampersand replacement=\&]
			\tensor*[^{(t_{\textrm{-}2\Lthetafin})}]{\coha}{^T_\qv} \ar{r}{\R S_{i_1}} \& \tensor*[^{(s_{i_2}\cdots s_{i_r})}]{\coha}{_{s_{i_1}\qv}^{T, (s_{i_1})}} \ar{r}{\R S_{i_2}} \& \tensor*[^{(s_{i_3}\cdots s_{i_r})}]{\coha}{_{s_{i_2}s_{i_1}\qv}^{T, (s_{i_1}s_{i_2})}} \ar{r}{\R S_{i_3}} \& \cdots \ar{r}{\R S_{i_r}} \& \coha_{t_{2\Lthetafin}\qv}^{T, (t_{\textrm{-}2\Lthetafin})} 
		\end{tikzcd}\ .
	\end{align}
	Observe that $\coha_\qv^{T, (t_{\textrm{-}2\Lthetafin})}\simeq \coha_{t_{2\Lthetafin}\qv}^{T, (t_{\textrm{-}2\Lthetafin})}$ because $t_{2\Lthetafin}\cdot \Theta=\Theta$. We may complete this chain of maps into a large commutative diagram by adding commutative squares such as below (for clarity, we only write it for the second step):
	\begin{align}
		\begin{tikzcd}[ampersand replacement=\&]
			\coha_{\Ltheta,(1)}^T \ar{rrrr}{\R S_{\textrm{-}2\Lthetafin}} \& \& \& \& \coha_{\Ltheta,(0)}^T\\
			\tensor*[^{(t_{\textrm{-}2\Lthetafin})}]{\coha}{^T_\qv} \ar{u} \ar{r}{\R S_{i_1}} \& \tensor*[^{(s_{i_2}\cdots s_{i_r})}]{\coha}{_{s_{i_1}\qv}^{T, (s_{i_1})}} \ar{r}{\R S_{i_2}} \& \tensor*[^{(s_{i_3}\cdots s_{i_r})}]{\coha}{_{s_{i_2}s_{i_1}\qv}^{T, (s_{i_1}s_{i_2})}} \ar{r}{\R S_{i_3}} \& \cdots \ar{r}{\R S_{i_r}} \& \coha_{t_{2\Lthetafin}\qv}^{T, (t_{\textrm{-}2\Lthetafin})} \ar{u}\\
			{} \& \tensor*[^{(s_{i_2})}]{\coha}{^T_{s_{i_1}\qv}} \ar{r}{\R S_{i_2}} \ar[u] \& \coha_{s_{i_2}s_{i_1}\qv}^{T, (s_{i_2})}\ar[u]\& {}\& {}\\
			{} \&\Y_\qv^-  \ar{r}{\overline T_{i_2}} \ar[u] \& \Y_\qv^-\ar[u] \& {} \& {}
		\end{tikzcd}\ .
	\end{align}
	
	The bottom square commutes (up to a sign) by Proposition~\ref{prop:overlineTi} and Theorem~\ref{thm:compatibility-braid-reflection-functors}, the middle square commutes via the isomorphisms \eqref{eq:isomorphism-S-uvw} and the top square is \eqref{eq:diagram-C}. Gathering all theses diagrams together yields a commutative square
	\begin{align}\label{eq:1-0}
		\begin{tikzcd}[ampersand replacement=\&, column sep=large]
			\coha_{\Ltheta,(1)}^{T,\Theta} \ar{r}{\R S_{\textrm{-}2\Lthetafin}} \& \coha_{\Ltheta,(0)}^{T,\Theta}\\
			\Y_\qv^- \ar{u}{\Phi^{\Llambda}_\Theta} \ar{r}{\overline T_{\textrm{-}\Lthetafin}} \& \Y_\qv^-  \arrow[swap]{u}{\Phi_\Theta}
		\end{tikzcd}
	\end{align}
	for some coweight $\Llambda$, which may be taken modulo $2 \coweightlatticefin$. Repeating this process composing the resulting diagrams yields a diagram 
	\begin{align}
		\begin{tikzcd}[ampersand replacement=\&, column sep=large]
			\coha_{\Ltheta,(2)}^{T,\Theta} \ar{r}{\R S_{\textrm{-}4\Lthetafin}} \& \coha_{\Ltheta,(0)}^{T,\Theta}\\
			\Y_\qv^- \ar{u}{\Phi^{\Llambda'}_\Theta} \ar{r}{\overline T_{\textrm{-}4\Lthetafin}} \& \Y_\qv^-  \arrow[swap]{u}{\Phi_\Theta}
		\end{tikzcd}
	\end{align}
	where $\Llambda'=\Llambda + t_{2\Lthetafin}(\Llambda)\equiv 2 \Llambda \equiv 0$ in $\coweightlatticefin/2\coweightlatticefin$, as wanted.
\end{proof}

We shall denote by the same symbol $\overline{T}_{\textrm{-}4\Lthetafin}$ the isomorphism $ \Y_{\Ltheta, (2)}\longrightarrow \Y_{\Ltheta, (0)}$ fitting in the diagram
\begin{align}
	\begin{tikzcd}[ampersand replacement=\&]
		\Y_{\Ltheta, (2)}\ar{r}{\overline T_{\textrm{-}4\Lthetafin}} \& \Y_{\Ltheta, (0)}\\
		\Y_\qv^-  \ar{r}{\overline T_{\textrm{-}4\Lthetafin}} \ar{u} \& \Y_\qv^-  \ar{u}
	\end{tikzcd}\ .
\end{align}

\begin{remark}\label{rem:arbitray-coweight-Stheta=Ttheta} 
	Theorem~\ref{thm:Stheta=Ttheta} and it proof remain valid for an arbitrary $\Llambda$ in place of $2\Lthetafin \in \coweightlatticefin$, provided that the definitions of the intervals $\kappa_0, \kappa_1$ are suitably modified in order for the inclusions in Formula~\eqref{eq:open-embeddings-stability-kappa} to remain true, and one considers a multiple $n\Llambda$ in place of $4\Lthetafin$. 
\end{remark}

\section{Limit of affine Yangians and COHA of an ADE surface}\label{sec:limit-affine-Yangian}

In this section, we put the pieces together and finally obtain an algebraic description of the nilpotent COHA of resolutions of Kleinian singularities in terms of affine Yangians, see our main Theorems~\ref{thm:coha-surface-as-limit1}, \ref{thm:coha-surface-as-limit2}, and \ref{thm:coha-surface-as-limit3}.

\subsection{Recollection of nilpotent cohomological Hall algebras}

Let $X$ be a quasi-projective sche\-me, $Z$ a closed subscheme of $X$. Let $\dstackCohpsnil(\widehat{X}_Z)$ be the derived stack  of properly supported coherent sheaves on $X$ whose set-theoretic support lies on $Z$. As shown in \cite{DPSSV-1}, this stack exhibits an intricate geometric structure. Indeed, it arises naturally as the colimit as $n$ tends to infinity of the derived stacks of coherent sheaves scheme-theoretically supported on thickenings $nZ$ of $Z$, and consequently is not a (derived) geometric stack. 

The following result characterizes the geometry of $\dstackCohpsnil(\widehat{X}_Z)$:
\begin{theorem}[{\cite[Theorem~\ref*{foundation-thmintro:B}]{DPSSV-1}}]\label{thm:nilpotent_coherent}
	Let $X$ be a smooth quasi-projective scheme, $Z \subset X$ a closed subscheme and assume that $X$ admits a projective compactification that contains $Z$ as a closed subscheme. Then:
	\begin{enumerate}[label=(\roman*)]\itemsep=0.2cm
		\item The underlying reduced stack $\Red{\dstackCohpsnil(\widehat{X}_Z)}$ is geometric, and the canonical morphism
		\begin{align}
			\Red{\dstackCohpsnil(\widehat{X}_Z)} \longrightarrow \dstackCohps(X) 
		\end{align}
		is a closed immersion.
		
		\item The derived stack $\dstackCohpsnil(\widehat{X}_Z)$ coincides with the formal completion of $\dstackCohps(X)$ along $\Red{\dstackCohpsnil(\widehat{X}_Z)}$.
	\end{enumerate}
\end{theorem}
This theorem shows that $\Red{\dstackCohpsnil(\widehat{X}_Z)}$ generalizes the \textit{global nilpotent cone}, i.e., the moduli stack of nilpotent Higgs sheaves on a smooth projective complex curve $C$, when $X=T^\ast C$ and $Z=C\subset X$ is the zero section. 

Now, the construction of the COHA $\coha_{X,Z}$ of coherent sheaves on $X$ set-theoretically supported on $Z$ in \cite{DPSSV-1} proceeds in two stages. First, we introduce the abstract class of \textit{admissible indgeometric} derived stacks, which is essentially characterized by the property that their underlying reduced stack is geometric. Theorem~\ref{thm:nilpotent_coherent} shows that both $\dstackCohpsnil(\widehat{X}_Z)$ and the derived stack $\calS_n\dstackCohpsnil(\widehat{X}_Z)$, classifying $n$-flags of objects in $\dstackCohpsnil(\widehat{X}_Z)$ for any $n \geq 2$, are admissible indgeometric. Second, we adapt the motivic formalism developed by A. Khan for derived stacks (see \cite{Khan_VFC, Khan_Voevodsky_criterion} and references therein) and show that it gives rise to a well-behaved Borel-Moore homology theory of admissible indgeometric stacks. In particular, the Borel-Moore homology space $\HBMbullet(\dstackCohpsnil(\widehat{X}_Z))$ is defined as a topological abelian group and satisfies the standard functorialities. We prove the following:

\begin{theorem}[{\cite[Theorem~\ref*{foundation-thmintro:existscoha}]{DPSSV-1}}]\label{thm:existscoha}
	Let $(X,Z)$ be as in Theorem~\ref{thm:nilpotent_coherent}, with $\dim(X)=2$. Then
	\begin{enumerate}[label=(\roman*)]\itemsep0.2cm
		\item In the diagram
		\begin{align}
			\begin{tikzcd}[ampersand replacement=\&]			
				\dstackCohpsnil(\widehat{X}_Z) \times \dstackCohpsnil(\widehat{X}_Z) \& \calS_2\dstackCohpsnil(\widehat{X}_Z) \ar[swap]{l}{q} \ar{r}{p} \& \dstackCohpsnil(\widehat{X}_Z)
			\end{tikzcd}
		\end{align}
		the morphism $q$ is representable by finitely connected (in the sense of \cite[Definition~\ref*{foundation-def:admissible-rpas-connected}]{DPSSV-1}), quasi-compact, lci derived geometric stacks, and the morphism $p$ is locally rpas (in the sense of \textit{loc. cit.}). The topological vector space $\coha_{\widehat{X}_Z}\coloneqq \HBMbullet(\dstackCohpsnil(\widehat{X}_Z))$, equipped with the multiplication 
		\begin{align}
			\mu\coloneqq (p_\ast\circ q^!)\circ \boxtimes\ ,
		\end{align}
		is an associative algebra, where 
		\begin{align}
			\boxtimes\colon \HBMbullet( \dstackCohpsnil(\widehat{X}_Z) \widehat{\otimes} \HBMbullet(\dstackCohpsnil(\widehat{X}_Z))\to \HBMbullet( \dstackCohpsnil(\widehat{X}_Z) \times \dstackCohpsnil(\widehat{X}_Z) )
		\end{align}
		is the exterior product.
		
		\item \label{item:COHA-ii} Let $(X',Z')$ be another pair satisfying the same properties as $(X,Z)$. Any isomorphism of formal completions $\widehat{X}_Z \simeq \widehat{X'}_{Z'}$ induces an isomorphism of algebras $\coha_{X,Z} \simeq \coha_{X',Z'}$.
		
		\item \label{item:COHA-iii} Let $j\colon Z' \subset Z$ be a closed subset. The direct image $j_\ast$ gives rise to a continuous algebra morphism $\coha_{X,Z'} \to \coha_{X,Z}$. 
	\end{enumerate}
\end{theorem}
If there is an algebraic torus $T$ acting on $X$ such that $Z$ is $T$-invariant, the above theorem extends verbatim to the equivariant setting. Moreover, the above construction recovers the cohomological Hall algebra of nilpotent Higgs sheaves on a smooth projective complex $C$ studied in \cite{Sala_Schiffmann}, when $X=T^\ast C$ and $Z=C$.

\subsection{Nilpotent cohomological Hall algebra of $\widehat{\rsv}_C$}\label{subsec:coha-rsv} 

Recall the resolution of Kleinian singularity $\rsv \to \ssv\coloneqq\C^2/G$ and the diagonal torus $A\subset \mathsf{GL}(2,\C)$ centralizing $G$. Note that the torus $A$ could be $\{1\}$, $\G_m$, or $\G_m\times\G_m$,  for $G$ of type $A$, while it could be $\{1\}$ or $\G_m$ for $G$ of types $D$ or $E$. 

Let $\dstackCohpsnil(\widehat{\rsv}_C)$ be the derived stack  of properly supported coherent sheaves on $\rsv$ whose set-theoretic support lies on $C$ (cf.\ Theorem~\ref{thm:nilpotent_coherent}). First note that the derived McKay equivalence, Theorem~\ref{thm:nilpotent}, induces an isomorphism $G_0( \catCoh_C(\rsv) )\simeq G_0( \nilpPi )\simeq  \Z I$. Thus, the derived stack $\dstackCohpsnil(\widehat{\rsv}_C)$ decomposes into a disjoint union of open and closed substacks of the form $\dstackCohpsnil(\widehat{\rsv}_C;\bfd)$, where the latter is the substack parametrizing sheaves of class $\bfd\in \Z I$. 

By Theorem~\ref{thm:existscoha}, we consider the equivariant nilpotent cohomological Hall algebra $\coha_{\rsv,C}^A$ with $\Q$-coefficients. The topological algebra $\coha_{\rsv,C}^A$ is a $\Z \times \Z I$-graded $\Hbullet_A$-algebra. The first grading (\textit{vertical}) corresponds to the homology grading; the second one (\textit{horizontal}) is the class in the Grothendieck group $G_0(\catCoh_C(\rsv))$ of the category $\catCoh_C(\rsv)$ of coherent sheaves set-theoretically supported on $C$. As a topological vector space, $\coha^A_{\rsv,C}$ may be realized as 
\begin{align}
	\coha_{\rsv,C}^A\simeq \lim_{\calU_i} \HBMbulletA(\calU_i)\ ,
\end{align} 
where $\{\calU_i\}_i$  runs over an admissible open exhaustion of $\dstackCohpsnil(\widehat{\rsv}_C)$, and the limit is taken in the category of graded vector spaces (cf.\ \cite[Remark~\ref*{foundation-rem:quasi-compact_topology}]{DPSSV-1}). The first grading is shifted in the following way:
\begin{align}
	(\coha_{\rsv,C}^A)_{\bfd,k} \coloneqq \sfH^A_{k-(\bfd,\bfd)}(\dstackCohpsnil(\widehat{\rsv}_C;\bfd))\ .
\end{align}

 In order to relate it to quantum groups, we will twist the multiplication by the sign associated to the bilinear form 
 \begin{align}
 	\Theta\colon K_0(\catCoh_C(\rsv))\times  K_0(\catCoh_C(\rsv)) \simeq \Z I\times \Z I \longrightarrow \Z/2\Z
 \end{align}
defined by
\begin{align}\label{eq:twisting_form_ADE_resolution}
	\Theta(\alpha_i,\alpha_j)=\langle \alpha_i,\alpha_j\rangle\ , \qquad \Theta(\delta,\delta)=0\ , \qquad \Theta(\alpha_i,\delta)=\Theta(\delta,\alpha_i)=0 
\end{align}
for all $i,j \in I_\sff$. 

The symmetrization of the bilinear map $\Theta$ is equal to the Euler form $(-,-)$ and moreover $t_{\Llambda}\cdot \Theta=\Theta$ for any $\Llambda \in \coweightlatticefin$, as $\Z\delta$ belongs to the kernel of $\Theta$. As a consequence, Theorems~\ref{thm:PBW} and~\ref{thm:Stheta=Ttheta} are valid for this twist, see Proposition~\ref{prop:sign-twists-I} and Corollary~\ref{cor:Phi-sign-twists}. 

\begin{remark}
	Note that because $\Theta(\delta,\delta)=0$, the subalgebra corresponding to the derived stack of coherent sheaves on $\rsv$, set-theoretically supported on $C$, which are zero-dimensional sheaves is untwisted. 
\end{remark}

\begin{notation}
	In order to unburden the notation, we will refrain from using the notation $\coha_{\rsv,C}^{A,\Theta}$ and simply write $\coha_{\rsv,C}^A$.
\end{notation}

\subsection{Limit of COHA of affine quivers}\label{subsec:limit-coha-affine-quiver}

We will now use the tilting equivalence $\tau$ together with the braid group action to relate $\coha_{\rsv,C}^A$ with the twisted equivariant nilpotent cohomological Hall algebra of the quiver $\qv$. 

We begin by fixing notations concerning the torus. Recall that $T\coloneqq \C^\ast \times \C^\ast$ acts on all the moduli stacks or moduli spaces of representations of $\Pi_\qv$. In particular, it acts on $\rsv$, which may be realized as the moduli space of semistable $\delta$-dimensional representations for a suitable stability parameter. There is a subtorus of $T$ which maps isomorphically onto $A$. For simplicity, we shall denote it by $A$ as well. Note that because $\Lambda_\qv$ is equivariantly formal and $\Hbullet(\dLambda_\qv)$ is cohomologically pure (see \cite[Theorem~A]{SV_generators}), the natural map $\HBMbulletT(\dLambda_\qv) \otimes_{\Hbullet_T} \Hbullet_A \to \HBMbulletA(\dLambda_\qv)$ is an isomorphism (cf.\ Proposition~\ref{prop:surjectivity} applied to $\Lambda_\qv$). Base changing from $\Hbullet_T$ to $\Hbullet_A$ we obtain analogues of all the results of \S\ref{sec:Yangians-affine-quivers} and \S\ref{sec:quotients-affine-coha} relative to the torus $A$ instead of $T$.

We shall use the framework developed in Part~\ref{part:COHA-stability-condition} to approximate the derived moduli stack $\dstackCoh(\widehat{\rsv}_C)$ of coherent sheaves on $\rsv$ set-theoretically supported on $C$ with the derived moduli stacks associated to the images of the category $\nilpPi^{\leqslant 0}$ under the functor $(\R S_{\textrm{-}2\Lthetafin})^k$ when $k$ tends to $+\infty$. Note that for $k >0$ this functor maps objects of $\nilpPi^{\leqslant 0}$ to actual complexes of $\Pi_\qv$-modules. 

We set 
\begin{align}
	\scrD\coloneqq \Pi_\qv\Mod\quad \text{and} \quad \scrD_0\coloneqq \catDb_{\ps,\,\mathsf{nil}}(\ModPi)\ ,
\end{align}
where the former is the $\infty$-category of modules over $\Pi_\qv$. Recall that for $\ell \in \Z$ we have (cf.\ Formula~\eqref{eq:nu})
\begin{align}
	\nu_\ell\coloneqq \frac{1}{\pi}\arctan(2h\ell)\ .
\end{align}
Following \S\ref{sec:limiting-COHA}, for any $\ell\in \Z$, define the derived stack 
\begin{align}
	\elldLambdaqv\coloneqq \dstackCohps(\scrD_0,\tau_{\nu_\ell}) \cap \dstackCohps(\scrD_0,\tau_{1/2}) \ ,
\end{align}
where $\tau_\phi$ is the $t$-structure on $\scrD_0$ associated to the slicing $\calP_{\Ltheta}$ introduced in \S\ref{subsubsec:stability-condition-quiver} for any $\phi\in \R$ (cf.\ Remark~\ref{rem:slicing-t-structures}) and $\dstackCohps(\scrD_0,\tau_\phi)$ is given in Construction~\ref{construction:derived-stack-D_0}. Then, $\elldLambdaqv$ parametrizes the finite-dimensional nilpotent objects belonging to the category
\begin{align}
	\calP^{\mathsf{nil}}_{\Ltheta}((\nu_\ell, 1/2])\coloneqq \catDb_{\ps, \, \mathsf{nil}}(\ModPi)\cap \calP_{\Ltheta}((\nu_\ell, 1/2])\ .
\end{align}
When $\ell \in \N$, Theorem~\ref{thm:tau-slicing}--\eqref{thm:tau-slicing-3} gives a $A$-equivariant equivalence of derived stacks
\begin{align}
	\dLambda^{\kappa_\ell}_\qv \simeq \elldLambdaqv\ ,
\end{align}
while for $\ell<0$, the stack $\elldLambdaqv$ parametrizes complexes of $\Pi_\qv$-modules rather than actual modules. 

For any $\ell \in \Z$, by Theorem~\ref{thm:tau-slicing}--\eqref{thm:tau-slicing-4} the tilting equivalence $\tau$ induces $A$-equivariant isomorphisms
\begin{align}
	\dstackCohpsnil(\widehat{\rsv}_C)^{>2\ell} \simeq \elldLambdaqv\ .
\end{align}
We now extend the definition of $\coha_{\Ltheta, (\ell)}^A$ given in Formula~\eqref{eq:quotient-COHA-ell} by setting for any $\ell\in \Z$
\begin{align}
	\coha_{\Ltheta, (\ell)}^A\coloneqq \HBMbulletA(\elldLambdaqv )\ .
\end{align}
Thus, we have a chain of projections
\begin{align}\label{eq:chain-cohas-biinfinite}
	\cdots \twoheadrightarrow \coha_{\Ltheta,(-1)}^A \twoheadrightarrow \coha_{\Ltheta,(0)}^A \twoheadrightarrow \coha_{\Ltheta,(1)}^A \twoheadrightarrow \cdots \ .
\end{align}	
Note that all the maps $\rho_{\ell-1,\ell}\colon  \coha_{\Ltheta,(\ell-1)}^A \to  \coha_{\Ltheta,(\ell)}^A$ occuring in \eqref{eq:chain-cohas-biinfinite} are induced by the open embeddings $\elldLambdaqv\to \tensor*[^{(\ell-1)}]{\dLambda}{_\qv}$ and hence $\N \times \Z I$-graded. Define
\begin{align}
	\coha_\infty^A\coloneqq \lim_\ell \coha_{\Ltheta,(\textrm{-}\ell)}^A
\end{align}
with respect to the maps $\rho_{\ell-1,\ell}$, where the limit is equipped with the quasi-compact topology again. This is a topological $\N \times \Z I$-graded vector space.
\begin{theorem}\label{thm:coha-surface-as-limit1}
	\hfill
	\begin{enumerate}\itemsep0.2cm
		\item \label{item:coha-surface-as-limit1.1} There is a canonical $\N \times \Z I$-graded algebra structure on $\coha_\infty^A$ induced by the COHA multiplication on $\coha_\qv^A$.
		
		\item \label{item:coha-surface-as-limit1.2} There is a canonical graded algebra isomorphism $\Theta\colon \coha_{\rsv, C}^A \simeq \coha_\infty^A$.
	\end{enumerate}
	In particular, $\coha_\infty^A$ does not depend on the choice of a (strictly dominant) $\Lthetafin$.
\end{theorem}

\begin{proof}
	We apply the framework developed in \S\ref{sec:limiting-COHA}. For $k\in \N$, set $a_k\coloneqq \nu_{\textrm{-}k}+1$. Then, we have
	\begin{align}
		\lim_{k\to +\infty} a_k= -\frac{1}{2}+1=\frac{1}{2}\eqqcolon a_\infty \ .
	\end{align}
	Let $\Lambda\coloneqq \Z I$ and let $v\colon K_0(\scrD_0) \to \Lambda$ be the map that associates to the K-theory class of a nilpotent finite-dimensional representation of $\Pi_\qv$ its dimension vector.
	
	Now, since 
	\begin{align}\label{eq:equivalence}
		\R S_{\textrm{-}2k\Lthetafin} (\modPi) = \calP_{\Ltheta}((a_k-1,a_k])
	\end{align}
	for any $k\in \N$ by Theorem~\ref{thm:tau-slicing}--\eqref{thm:tau-slicing-2} and the standard $t$-structure on $\Pi_\qv\Mod$ is open, Assumption~\ref{assumption:limiting_2_Segal_stack} holds.
	
	Again, thanks to the equivalence~\eqref{eq:equivalence}, Assumption~\ref{assumption:quasi-compact_interval} holds for any $k\in \N$ since it is evidently true for $k=0$: in this case, Harder-Narasimhan strata of the moduli stack $\Lambda_\qv$ are known to be quasi-compact and locally closed. Thus, Corollary~\ref{cor:assumption-holds} holds as well. Since Assumption~\ref{assumption:limiting_CoHA_I}--\eqref{assumption:limiting_CoHA_I-1}) holds for $k=0$, by the equivalence~\eqref{eq:equivalence} it holds for any $k\in \N$. Moreover, Assumption~\ref{assumption:limiting_CoHA_I}--\eqref{assumption:limiting_CoHA_I-2} holds for $k=k'=0$, by using again  the equivalence~\eqref{eq:equivalence} and Corollary~\ref{cor:assumption-holds}, we obtain that Assumption~\ref{assumption:limiting_CoHA_I} holds for any $k\in \N$.
	
	Thus, we can apply Proposition~\ref{prop:limiting_CoHA} and we obtain an $A$-equivariant limiting cohomological Hall algebra
	\begin{align}
		\coha_{\scrD_0, \tau_\infty^+}^A\coloneqq \bigoplus_{\bfd \in \Z I} \lim_{\ell} \colim_{k \geqslant \ell} \sfH_\bullet^A\big( \dstackCohps\big( \scrD_0, (a_\ell-1, a_k]; \bfv \big) \big) 
	\end{align}
	as a $\Lambda$-graded-vector space, endowed with the quasi-compact topology.
	
	Now, Theorem~\ref{thm:limiting_vs_limit} yields 
	\begin{align}
		\colim_{k \geqslant \ell} \sfH_\bullet^A\big( \dstackCohps\big( \scrD_0, (a_\ell-1, a_k]; \bfd \big) \big) \simeq \sfH_\bullet^A\big( \dstackCohps\big( \scrD_0, (a_\ell-1, 1/2]; \bfd \big) \big) = \sfH_\bullet^A\big( \tensor*[^\ell]{\dLambda}{_\bfd} \big) \ .
	\end{align}
	Therefore, as $\Lambda$-graded vector spaces, $\coha_{\scrD_0, \tau_\infty^+}^A$ is isomorphic to $\coha_\infty^A$, endowed with the quasi-compact topology. This proves Statement~\eqref{item:coha-surface-as-limit1.1}.
	
	\medskip
	
	Now, notice that $\calP_{\Ltheta}((-1/2, 1/2])$ is equivalent to the abelian category $\catCoh_C(\rsv)$ of coherent sheaves on $\rsv$ set-theoretically supported on $\rsv$. Statement~\eqref{item:coha-surface-as-limit1.2} follows from Corollary~\ref{cor:limiting_COHA_vs_COHA_limit}, since its hypotheses are all verified by $\tau_\infty$, which corresponds to the standard $t$-structure on $\catDb_C(\catCoh(\rsv))$. 
\end{proof}	

\begin{remark}
	Let us explicitly describe the multiplication in $\coha_\infty^A$. For any $\ell \in\N$, let us denote by $\coha_{\Ltheta, (a_\ell-1, a_\ell]}^A$ the cohomological Hall algebra of the category $\calP^{\mathsf{nil}}_{\Ltheta}((a_\ell-1,a_\ell])$. There is a surjective morphism 
	\begin{align}
		\coha_{\Ltheta, (a_\ell-1, a_\ell]}^A \longrightarrow \coha_{\Ltheta,(\ell)}^A\ .
	\end{align}
	Let $x=(x_\ell)_\ell$ and $y=(y_\ell)_\ell$ be elements of $\coha_\infty^A$. Let us fix some arbitrary lifts $\tilde{x}_\ell, \tilde{y}_\ell$ of $x_\ell,y_\ell$ in $\coha_{\Ltheta, (a_\ell-1, a_\ell]}^A$ for $\ell \in \N$. Let us also denote by $\pi_{n,\ell}\colon \coha_{\Ltheta,(n)}^A \to \coha_{\Ltheta,(\ell)}^A$ the canonical restriction morphism when $n <\ell$. Then for any fixed $\ell$ we have
	\begin{align}
		(x \cdot y)_\ell= \pi_{n,\ell}(\tilde{x}_n \cdot \tilde{y}_n)
	\end{align}
	for any $n \ll \ell$, with the multiplication on the right-hand-side taking place in $\coha_{\Ltheta, (a_n-1, a_n]}^A$.
\end{remark}

\subsection{Limit of affine Yangians}\label{subsec:def-limit-affine-Yangian}

We now recast the limit $\coha_\infty^A$ as a limit of quotients of the affine Yangian. To do that, we shall use the results described in \S\ref{sec:quotients-affine-coha} and \S\ref{subsec:quotients}. For simplicity we still denote by $\Y_\qv^-,\ldots$ the specializations $\Y_\qv^- \otimes_{\sfH^\bullet_T} \Hbullet_A,\ldots$. Theorem~\ref{thm:Stheta=Ttheta} implies that we have a commutative diagram 
\begin{align}\label{diagram_restrict_coha-yangian}
	\begin{tikzcd}[ampersand replacement=\&, column sep=large]
		\coha_{\Ltheta,(0)}^A \ar{r}{} \& \coha_{\Ltheta,(2)}^A \ar{r}{\R S_{\textrm{-}4\Lthetafin}} \& \coha_{\Ltheta,(0)}^A\\
		\Y_{\Ltheta,(0)} \ar{u}{\Phi_{(0)}} \ar{r}{} \& \Y_{\Ltheta,(2)} \ar{u}{\Phi_{(2)}} \ar{r}{\overline{T}_{\textrm{-}4\Lthetafin}} \& \Y_{\Ltheta,(0)}  \arrow[swap]{u}{\Phi_{(0)}}\\
		\Y^-_\qv \ar{u}{} \ar{r}{\id} \& \Y_\qv^-  \ar{r}{\overline T_{\textrm{-}4\Lthetafin}} \ar{u} \& \Y_\qv^-  \ar{u}
	\end{tikzcd}\ ,
\end{align}
in which the maps $\coha_{\Ltheta,(0)}^A \to \coha_{\Ltheta,(2)}^A$ and $\Y_{\Ltheta,(0)} \to \Y_{\Ltheta,(2)}$ are the canonical restriction maps. The middle row of diagram~\eqref{diagram_restrict_coha-yangian} is explicitly given as the composition
\begin{align}
	\begin{tikzcd}[ampersand replacement=\&]
		\Y^-_\qv / \J_{\Ltheta,\kappa_0} \ar[two heads]{r} \& \Y^-_\qv / \J_{\Ltheta,\kappa_2} \ar{r}{T_{\textrm{-}4\Lthetafin}} \& T_{\textrm{-}4\Lthetafin}(\Y^-_\qv) / T_{\textrm{-}4\Lthetafin}(\J_{\Ltheta,\kappa_2}) \ar{r}{\mathsf{pr}} \& \Y^-_\qv / \J_{\Ltheta,\kappa_0}
	\end{tikzcd}\ ,
\end{align}
where the two rightmost maps are isomorphisms by Corollary~\ref{cor:PBW-consequences}. 

For an arbitrary $\ell \in \Z$, we define
\begin{align}
	\Y_{\Ltheta,(2\ell)}\coloneqq T_{4\ell\Lthetafin}(\Y^-_\qv)/T_{4\ell\Lthetafin}(\J_{\Ltheta,\kappa_0})\ .
\end{align}
\begin{remark}
	The above definition agrees with the definition of $\Y_{\Ltheta,(2\ell)}$ given in Formula~\eqref{eq:quotient-Yangian-ell} for $\ell\in \N$.
\end{remark}

We have, by \textit{transport de structure}, an isomorphism $T_{\textrm{-}4\Lthetafin}\colon \Y_{\Ltheta,(2\ell)} \to \Y_{\Ltheta,(2\ell-2)}$ which is $\N$-graded but acts as $t_{\textrm{-}4\Lthetafin}$ on the weight as well as a restriction map $\pi_{2\ell-2,2\ell}\colon \Y_{\Ltheta,(2\ell-2)} \to \Y_{\Ltheta,(2\ell)}$ which is a map of $\N \times \Z I$-graded vector spaces. We may extend in a unique way the isomorphisms $\Phi_{(2\ell)}\colon \Y_{\Ltheta,(2\ell)} \to\coha_{\Ltheta,(2\ell)}^A$ to all values of $\ell \in \Z$, in a way compatible with the restriction maps and braiding isomorphisms, namely so that the following diagrams commute: 
\begin{align}\label{diagram_restrict_coha-yangian2}
	\begin{tikzcd}[ampersand replacement=\&, column sep=large]
		\coha_{\Ltheta,(2\ell)}^A \ar{r}{\R S_{\textrm{-}4k\Lthetafin}} \& \coha_{\Ltheta,(2\ell-2k)}^A \\
		\Y_{\Ltheta,(2\ell)} \ar{u}{\Phi_{(2\ell)}} \ar{r}{\overline{T}_{\textrm{-}4k\Lthetafin}} \& 
		\Y_{\Ltheta,(2\ell-2k)} \ar{u}{\Phi_{(2\ell-2k)}}
	\end{tikzcd}
	\quad\text{and}\quad
	\begin{tikzcd}[ampersand replacement=\&, column sep=large]
		\coha_{\Ltheta,(2\ell-2)}^A \ar{r}{\rho_{2\ell-2,2\ell}} \& \coha_{\Ltheta,(2\ell)}^A \\
		\Y_{\Ltheta,(2\ell-2)} \ar{r}{\pi_{2\ell-2,2\ell}} \ar{u}{\Phi_{(2\ell-2)}} \& 
		\Y_{\Ltheta,(2\ell)} \ar{u}{\Phi_{(2\ell)}}
	\end{tikzcd}
\end{align}
for all $\ell \in \Z$ and $k \in \N$. 

Define
\begin{align}
	\Y_\infty^+ \coloneqq \lim_\ell \Y_{\Ltheta,(\textrm{-}2\ell)}\ .
\end{align}
The limit is equipped with the quasi-compact topology again. The following result holds.
\begin{theorem}\label{thm:coha-surface-as-limit2}
	The following holds:
	\begin{enumerate}\itemsep0.2cm
		\item \label{item:coha-surface-as-limit2.1} There is a canonical $\N \times \Z I$-graded algebra structure on $\Y_\infty^+$ induced by the multiplication on $\Y_\qv$.
		
		\item \label{item:coha-surface-as-limit2.2} There is a canonical $\N \times \Z I$-graded algebra isomorphism 
		\begin{align}
			\begin{tikzcd}[ampersand replacement=\&]
				\Phi_\infty\colon \coha_\infty^A\ar{r}{\sim} \& \Y_\infty^+
			\end{tikzcd}\ ,
		\end{align}
		induced by the maps $\Phi_{(2\ell)}$.
	\end{enumerate}
\end{theorem} 

\begin{notation}
	We denote by $\Theta_{\rsv,C}\colon\coha_{\rsv, C}^A \simeq \Y_\infty^+$ the composition of $\Theta$ and $\Phi_\infty$:
	\begin{align}
		\begin{tikzcd}[ampersand replacement=\&, baseline=(current bounding box.south)]
			\coha_{\rsv, C}^A \ar{r}{\Theta} \ar{rd}[swap]{\Theta_{\rsv,C}} \& \coha_\infty^A \ar{d}{\Phi_\infty}\\
			\& \Y_\infty^+
		\end{tikzcd} \tag*{\qedhere} 
	\end{align}	
\end{notation}

\begin{remark}
	It follows from the above theorem that $\Y_\infty^+$ is independent of the choice of $\Lthetafin$ in the strictly dominant chamber. This can also be seen directly.
\end{remark}

\begin{remark}
	For the reader's comfort, let us explicitly describe the multiplication in $\Y_\infty^+$. Let $x=(x_{2\ell})_\ell$ and $y=(y_{2\ell})_\ell$ be elements of $\Y_\infty^+$. Let us fix some arbitrary lifts $\tilde{x}_{2\ell}, \tilde{y}_{2\ell}$ of $x_{2\ell},y_{2\ell}$ in $T_{4\ell\Lthetafin}(\Y^-_\qv)$ for $\ell \in \Z$. For $n<\ell$, let us also denote by $\pi_{2n,2\ell}\coloneqq \pi_{2\ell-2,2\ell} \circ \cdots \circ \pi_{2n,2n+2}\colon \Y_{\Ltheta,(2n)} \to \Y_{\Ltheta,(2\ell)}$. Then for any fixed $\ell$ we have
	\begin{align}
		(x \cdot y)_\ell= \pi_{2n,2\ell}(\tilde{x}_{2n} \cdot \tilde{y}_{2n})
	\end{align}
	for any $n \ll \ell$, with the multiplication on the right-hand-side taking place in $T_{4n\Lthetafin}(\Y^-_\qv)$.
\end{remark}

\begin{remark}\label{rem:third-limit-realization-coha}
	It is also possible to give presentation of $\Y_\infty^+$ purely in terms of $\Y_{\Ltheta,(0)}$ using the isomorphisms $T_{\textrm{-}4\ell\Lthetafin}\colon \Y_{\Ltheta,(2\ell)} \to \Y_{\Ltheta,(0)}$. Namely, for any $\bfd \in \Z I$ we have
	\begin{align}
		(\Y_\infty^+)_{\bfd}\simeq \lim_\ell (\Y_{\Ltheta,(0)})_{t_{\textrm{-}4\ell\Lthetafin(\bfd)}}
	\end{align}	
	with transition maps
	\begin{align}
		\begin{split}
			\tilde{\pi}_{2\ell,2\ell+2}\colon (\Y_{\Ltheta,(0)})_{t_{\textrm{-}4\ell\Lthetafin(\bfd)}} &\longrightarrow (\Y_{\Ltheta,(0)})_{t_{\textrm{-}4(\ell+1)\Lthetafin(\bfd)}}\\
			x &\longmapsto \overline{T}_{\textrm{-}4\Lthetafin}(x)
		\end{split}\ .
	\end{align}
	Again, the multiplication is induced by that on $\Y^-_\qv$. In other words, we may represent (and make computations with) elements in $\Y_\infty^+$ as collections $(x_{2\ell})_\ell$ of elements in $\Y^-_\qv$ (rather than in a translate $T_{4\ell\Lthetafin}(\Y^-_\qv)$) but with increasing weights. 
\end{remark}

\begin{remark}\label{remark :semi-inifinite-PBW-projection}
	The collection of twisted $PBW$ projection maps $\mathsf{pr}_{2\ell} \coloneqq T_{4\ell\Lthetafin} \circ \mathsf{pr} \circ T_{\textrm{-}4\ell\Lthetafin} \colon \Y_\qv \to T_{4\ell\Lthetafin}(\Y^-_\qv)$ gives rise, after further projection to $T_{4\ell\Lthetafin}(\Y^-_\qv/\J^-_{\Ltheta,\kappa_0})$, to a morphism $\Y_\qv \to \Y_\infty^+$, which may be interpreted as a `semi-infinite' $PBW$-projection map. 
\end{remark}

\subsection{Action of $\coweightlatticefin$ }\label{subsec:limit-tensor-line-bundles}

The abelian group $\coweightlatticefin$ acts by tensor product on $\catCoh(\widehat{\rsv}_C)$ and hence by algebra automorphisms on $\coha^A_{\rsv,C}$. Since the twist $\Theta$ is invariant under $t_{\Llambda}$, the tensor product is compatible with the twist in the multiplication. In this paragraph, we make explicit the induced action on $\Y_\infty^+$. Fix $\Llambda \in \coweightlatticefin$ and consider the composition
\begin{align}\label{eq:composition-T-lambda}
	\begin{tikzcd}[ampersand replacement=\&]
		\Y^-_\qv \ar{r}{L_{\Llambda}} \& \Y_\qv \ar{r}{\mathsf{pr}_\ell}\& T_{4\ell\Lthetafin}(\Y^-_\qv)
	\end{tikzcd}
\end{align}
where $\mathsf{pr}_\ell\coloneqq T_{4\ell\Lthetafin}\circ \mathsf{pr} \circ T_{\textrm{-}4\ell\Lthetafin}$ is the twisted $PBW$ projection map to $T_{2\ell\Lthetafin}(\Y^-_\qv)$. Because $L_{\Llambda}$ acts on the weights as the translation element $t_{\Llambda}$ of the Weyl group, and because $t_{\Llambda}$ preserves $\Delta_\sff^- \times \Z\delta$, the composition \eqref{eq:composition-T-lambda} gives rise to a map $L_{\Llambda,2\ell}\colon \Y^-_\qv/\J^-_{\Ltheta,\kappa_0} \to T_{4\ell\Lthetafin}(\Y^-_\qv/\J^-_{\Ltheta,\kappa_0})=\Y_{\Ltheta,(2\ell)}$. Moreover, by construction the diagram
\begin{align}
	\begin{tikzcd}[ampersand replacement=\&, column sep=large]
		\Y^-_\qv/\J^-_{\Ltheta,\kappa_0} \ar{r}{L_{\Llambda,2\ell-2}} \arrow[swap]{dr}{L_{\Llambda,2\ell}} \& \Y_{\Ltheta,(2\ell-2)}
		\ar{d}{\pi_{2\ell-2,2\ell}}	\\
		\& \Y_{\Ltheta,(2\ell)}
	\end{tikzcd}
\end{align}
commutes. This shows that the collection of maps $L_{\Llambda,2\ell}$ gives rise to a well-defined map $\Y^-_\qv/\J_{\Ltheta,\kappa_0} \to \Y_\infty^+$. The same process allows us to define maps $T_{4\ell\Lthetafin}(\Y^-_\qv/\J_{\Ltheta,\kappa_0}) \to \Y_\infty^+$ which fit into a commutative diagram
\begin{align}
	\begin{tikzcd}[ampersand replacement=\&, column sep=large]
		T_{4(\ell-1)\Lthetafin}(\Y^-_\qv/\J^-_{\Ltheta,\kappa_0}) \ar{r} \arrow[swap]{d}{\pi_{2\ell-2,2\ell}} \& \Y_\infty^+\\
		T_{4\ell\Lthetafin}(\Y^-_\qv/\J^-_{\Ltheta,\kappa_0})  \arrow{ur}
		\& 
	\end{tikzcd}
\end{align}
which finally yields a map $L_{\Llambda}\colon \Y_\infty^+\to \Y_\infty^+$.

\begin{proposition}\label{prop:coha-tensor-product-intertwine} 
	\hfill
	\begin{enumerate}\itemsep0.2cm
		\item \label{item:action-coweight-affine-Yangian-limit} For any $\Llambda,\Lmu \in \coweightlatticefin$ we have $L_{\Llambda} L_{\Lmu}=L_{\Llambda+\Lmu} \in \End(\Y_\infty^+)$, i.e., $\coweightlatticefin$ acts on $\Y_\infty^+$. 
		\item \label{item:Theta'-intertwines-coweight-action} The map $\Theta_{\rsv,C}$ intertwines the actions of $\coweightlatticefin$ on $\coha_{\rsv, C}^A$ and  $\Y_\infty^+$.
	\end{enumerate}
\end{proposition}   

\begin{proof}
	Statement~\eqref{item:action-coweight-affine-Yangian-limit} follows from the fact that braid operators $T_{\Llambda},T_{\Lmu},T_{4\ell\Lthetafin} \in \Aut(\Y_\qv)$ all commute and preserve weights in $-\Delta \times \Z\delta$. Statement~\eqref{item:Theta'-intertwines-coweight-action} for dominant $\Llambda$ is a consequence of Theorem~\ref{thm:Stheta=Ttheta} and Proposition~\ref{prop:braid}--\ref{item:prop:braid-(1)}, see Remark~\ref{rem:arbitray-coweight-Stheta=Ttheta}. The case of an arbitrary $\Llambda$ then follows from the definition of $L_{\Llambda}$ (cf.\ Formula~\eqref{def:Llambda}).
\end{proof}

\begin{remark}
	From Proposition~\ref{prop:coha-tensor-product-intertwine}--\eqref{item:Theta'-intertwines-coweight-action} we deduce that $\coweightlatticefin$ acts on $\Y_\infty^+$ by algebra automorphisms. This may also be checked directly from the definition, using the fact that $\coweightlatticefin$ acts on $\Y_\qv$ by algebra automorphisms and preserves weights in $-\Delta \times \Z\delta$.
\end{remark}

\subsection{Action of tautological classes on $\coha_{\rsv,C}^A$} \label{subsec:coha-rsv-tautological} 

In this section we will define certain tautological cohomology classes on $\dstackCohpsnil(\widehat{\rsv}_C)$, forming a commutative family of operators on $\coha_{\rsv,C}^A$. Recall that the cohomology ring of $\dstackCohpsnil(\widehat{\rsv}_C)$ can be described as follows, see \cite[Notation~\ref*{torsion-pairs-notation:homology_and_cohomology} and Remark~\ref*{torsion-pairs-rem:quasi-compact_topology}]{DPS_Torsion-pairs},
\begin{align}
	\Hbullet_A(\dstackCohpsnil(\widehat{\rsv}_C))\simeq \lim_{\calU_i} \Hbullet_A(\calU_i)\ ,
\end{align}
where $\{\calU_i\}$ is an admissible open exhaustion of $\dstackCohpsnil(\widehat{\rsv}_C)$, and the limit is taken in the category of $\Z \times \Z I$-graded rings.  There is a similar description of the cohomology ring $\Hbullet_A(\dstackCohpsnil(\widehat{\rsv}_C) \times \rsv)$.

There is a natural action 
\begin{align}
	\Hbullet_A(\dstackCohpsnil(\widehat{\rsv}_C)) \otimes  \HBMbulletA(\dstackCohpsnil(\widehat{\rsv}_C)) \to  \HBMbulletA(\dstackCohpsnil(\widehat{\rsv}_C))
\end{align}
given by the cap product (cf.\ \cite[Construction~\ref*{torsion-pairs-constr:operations}--(\ref*{torsion-pairs-constr:operations-cap-product})]{DPS_Torsion-pairs}). 

Next, we consider tautological sheaves. The stack $\dstackCohpsnil(\widehat{\rsv}_C)$ carries a universal sheaf $\scrU \in \catCoh(\widehat{\rsv}_C\times \dstackCohpsnil(\widehat{\rsv}_C))$. It will be more convenient to consider a slightly different sheaf. Let $\scrE' \in \catCoh(\rsv\times \dstackCoh_{C,\ps}(\rsv))$ be the universal sheaf of $\dstackCoh_{C,\ps}(\rsv)$. Let $\jmath\colon \widehat{\rsv}_C \to \rsv$ be the canonical morphism, and let 
\begin{align}
	\bfjmathhat_\ast \colon \dstackCohpsnil(\widehat{\rsv}_C) \longrightarrow \dstackCoh_{C,\ps}(\rsv)
\end{align}
be the morphism constructed in \cite[\S\ref*{foundation-subsec:nilpotent-vs-set-theoretic}]{DPSSV-1}. Thus we have a morphism
\begin{align}
	\begin{tikzcd}[ampersand replacement=\&, column sep=large]
		\widehat{\rsv}_C\times \dstackCohpsnil(\widehat{\rsv}_C)  \ar{r}{\jmath\times \id}\& \rsv\times \dstackCohpsnil(\widehat{\rsv}_C)   \ar{r}{\id\times \bfjmathhat_\ast }\& \rsv\times \dstackCoh_{C,\mathsf{ps}}(\rsv) 
	\end{tikzcd} \ .
\end{align}
Set $\scrE\coloneqq(\jmath\times \id)_\ast(\scrU) \in \catCoh(\rsv\times\dstackCohpsnil(\widehat{\rsv}_C))$. Note that $\scrE$ is coherent because $\scrU$ is coherent. \cite[Theorem~\ref*{foundation-thm:nilpotent_into_set_theoretic_moduli}]{DPSSV-1} yields the following.
\begin{lemma}
	We have 
	\begin{align}
		\scrE \simeq \big(\id \times \bfjmathhat_\ast \big)^\ast (\scrE')\ .
	\end{align}
\end{lemma}

We are now ready to define some tautological cohomology classes on $\dstackCohpsnil(\widehat{\rsv}_C)$. We fix a smooth $A$-equivariant compactification $h\colon \rsv \to \overline{\rsv}$ and consider the natural extension 
\begin{align}
	\overline{\scrE} \in \catCoh( \overline{\rsv}\times \dstackCohpsnil(\widehat{\rsv}_C))
\end{align}
of $\scrE$. Note that $\widehat{\rsv}_C \simeq \widehat{\overline{\rsv}}_C$. Locally on $\catCoh(\overline{\rsv}\times \dstackCohpsnil(\widehat{\rsv}_C))$, the sheaf $\overline{\scrE}$ is of finite perfect amplitude, hence we may define the Chern character 
\begin{align}
	\ch(\overline{\scrE}) \in  \sfH^\bullet_A(\overline{\rsv}\times \dstackCohpsnil(\widehat{\rsv}_C))\ .
\end{align}

We consider the canonical projections
\begin{align}
	\begin{tikzcd}[ampersand replacement=\&]
		\& \overline{\rsv}\times \dstackCohpsnil(\widehat{\rsv}_C) 
		\ar{dl}[swap]{\ev_{\overline{\rsv}}}\ar{dr}{\ev}\& \\
		\overline{\rsv} \& \& \dstackCohpsnil(\widehat{\rsv}_C) 
	\end{tikzcd}\ .
\end{align}
Given any class $\alpha \in \Hbullet_A(\overline{\rsv})$, let 
\begin{align}
	\overline{\nu}_\alpha\coloneqq \ev_!(\ev_{\overline{\rsv}}^\ast(\alpha) \cup \ch(\overline{\scrE})) \ . 
\end{align}
Then the assignment $\alpha \mapsto \overline{\nu}_\alpha$ defines a homomorphism of $\HBMbulletA$-modules 
\begin{align}\label{eq:tautclassesA}
	\overline{\nu}\colon \Hbullet_A(\overline{\rsv})\longrightarrow  \Hbullet_A(\dstackCohpsnil(\widehat{\rsv}_C))\ .
\end{align}
Now, we will show that this homomorphism factors through the restriction map
\begin{align}
	h^\ast\colon \Hbullet_A(\overline{\rsv})\longrightarrow \sfH^\bullet_A(\rsv)\ .
\end{align}

Let $f\colon C\to \overline{\rsv}$ be the canonical closed immersion. Let $U \coloneqq \overline{\rsv}\smallsetminus C$. Note that the restriction 
\begin{align}
	\ch(\overline{\scrE}) \big\vert_{(\overline{\rsv}\smallsetminus X) \times \dstackCohpsnil(\widehat{\rsv}_C)} =0
\end{align}
since 
$\overline{\scrE}\big\vert_{(\overline{\rsv}\smallsetminus X) \times \dstackCohpsnil(\widehat{\rsv}_C)} =0$. Therefore, the long exact sequence associated to the closed immersion $f\times \id$ shows that 
\begin{align}
	\ch(\overline{\scrE})  =  (f\times \id)_\ast(c) 
\end{align}
for some (not necessarily unique) relative cohomology class $c\in \sfH^\bullet_A(\overline{\rsv} \times\dstackCohpsnil(\widehat{\rsv}_C), U\times\dstackCohpsnil(\widehat{\rsv}_C))$. 

Now, let $f'\colon C \to \rsv$ be the natural closed immersion and let $U' \coloneqq \rsv \smallsetminus C$. Recall that excision yields an isomorphism of relative cohomology groups
\begin{align}
	\varepsilon\colon \sfH^\bullet_A(\overline{\rsv} \times\dstackCohpsnil(\widehat{\rsv}_C), U\times\dstackCohpsnil(\widehat{\rsv}_C)  ) \longrightarrow \sfH^\bullet_A(X \times\dstackCohpsnil(\widehat{\rsv}_C), U'\times\dstackCohpsnil(\widehat{\rsv}_C)  )\ .
\end{align}
Moreover, one has cup-products 
\begin{multline}
	\sqcup\colon \sfH^\bullet_A(\overline{\rsv} \times\dstackCohpsnil(\widehat{\rsv}_C)) 
	\otimes \sfH^\bullet_A(\overline{\rsv} \times\dstackCohpsnil(\widehat{\rsv}_C), U\times\dstackCohpsnil(\widehat{\rsv}_C))  \longrightarrow\\
	\sfH^\bullet_A(\overline{\rsv} \times\dstackCohpsnil(\widehat{\rsv}_C), U\times\dstackCohpsnil(\widehat{\rsv}_C) )
\end{multline}
and 
\begin{multline}
	\sqcup'\colon \sfH^\bullet_A(X \times\dstackCohpsnil(\widehat{\rsv}_C)) \otimes \sfH^\bullet_A(X \times\dstackCohpsnil(\widehat{\rsv}_C), U'\times\dstackCohpsnil(\widehat{\rsv}_C))  \longrightarrow\\
	\sfH^\bullet_A(\overline{\rsv} \times\dstackCohpsnil(\widehat{\rsv}_C), U'\times\dstackCohpsnil(\widehat{\rsv}_C) ) 
\end{multline}
which are naturally compatible with $\varepsilon$. Finally, note that the restriction map $h^\ast\colon \sfH^\bullet_A(\overline{\rsv}) \to \sfH^\bullet_A(\rsv)$ is surjective since the equivariant cohomology of $\rsv$ is pure. 

\begin{lemma}\label{lem:restoSA} 
	Let $\ev_{\rsv}\colon \rsv\times \dstackCohpsnil(\widehat{\rsv}_C)\to X$ be the canonical projection. Given $\beta \in \sfH^\bullet_A(\rsv)$, one has 
	\begin{align}\label{eq:relformulaA}
		(\ev_{\overline{\rsv}}^\ast(\alpha)\cup \ch(\overline{\scrE})) = (f\times \id)_\ast (\varepsilon^{-1}(\ev_{\rsv}^\ast (\beta) \sqcup' \varepsilon(c))) 
	\end{align}
	for any $\alpha \in \sfH^\bullet_A(\overline{\rsv})$ and any $c\in \HBMbulletA(\overline{\rsv} \times\dstackCohpsnil(\widehat{\rsv}_C), 
	U \times\dstackCohpsnil(\widehat{\rsv}_C))$ so that $h^\ast(\alpha)=\beta$ and $(f\times \id)_\ast(c) = \ch(\overline{\scrE})$. 
\end{lemma} 

\begin{proof} 
	Given  $\alpha \in \sfH^\bullet_A(\overline{\rsv})$ and any $c\in \HBMbulletA(\overline{\rsv} \times\dstackCohpsnil(\widehat{\rsv}_C), 
	U \times\dstackCohpsnil(\widehat{\rsv}_C))$ so that $h^\ast(\alpha)=\beta$ and $(f\times \id)_\ast(c) = \ch(\overline{\scrE})$, one has 
	\begin{align} 
		(f\times \id)_\ast (\ev_{\overline{\rsv}}^\ast(\alpha) \sqcup c) = \ev_{\overline{\rsv}}^\ast(\alpha) \cup \ch(\overline{\scrE})\ , 
	\end{align} 
	as well as 
	\begin{align}
		\varepsilon(\ev_{\overline{\rsv}}^\ast(\alpha) \sqcup \ch(\overline{\scrE})) = \ev_{\rsv}^\ast(\beta) \sqcup' \varepsilon(c)\ .
	\end{align}
	This implies identity~\eqref{eq:relformulaA}.
	
	Next, suppose $c_1,c_2\in \HBMbulletT(C \times\dstackCohpsnil(\widehat{\rsv}_C) )$ satisfy the conditions 
	\begin{align}
		(f\times\id)_\ast(c_k)= \ch(\overline{\scrE})
	\end{align}
	for $1\leq k \leq 2$. Then one has 
	\begin{align} 
		(f\times\id)_\ast \, \varepsilon^{-1}(\ev_{\rsv}^\ast( \beta) \sqcup' \varepsilon(c_1-c_2)) & = (f\times \id)_\ast(\ev_{\overline{\rsv}}^\ast(\alpha) \sqcup (c_1-c_2))
		\\ &= \ev_{\overline{\rsv}}^\ast(\alpha) \cup (f\times \id)_\ast (c_1-c_2) =0
	\end{align} 
	for any $\alpha \in \sfH_T^\bullet(\overline{\rsv})$ so that $h^\ast(\alpha) = \beta$. 
\end{proof} 

As a consequence of Lemma~\ref{lem:restoSA},  we obtain the following.
\begin{corollary}\label{cor:restoSB} 
	The morphism \eqref{eq:tautclassesA} factors as $\overline{\nu} = \nu \circ h^\ast$ where 
	\begin{align}
		\nu\colon \sfH^\bullet(\rsv) \longrightarrow \sfH^\bullet_A( \dstackCohpsnil(\widehat{\rsv}_C) )\ , 
	\end{align}
	is the unique $\sfH^\bullet_A$-module homomorphism mapping $\beta \in \sfH^\bullet(\rsv)$ to 
	\begin{align}
		\nu_\beta\coloneqq \ev_!(\ev_{\overline{\rsv}}^\ast(\alpha) \cup \ch(\overline{\scrE})) 
	\end{align}
	for any $\alpha \in  \sfH^\bullet(\overline{\rsv})$ so that $\beta = h^\ast(\alpha)$. 
\end{corollary} 
For any $\ell\in \N$, define 
\begin{align}
	\nu_{\ell, \beta} \coloneqq \ev_!(\ev_{\overline{\rsv}}^\ast(\alpha) \cup \ch_\ell(\overline{\scrE}))  \in \sfH^\bullet_A( \dstackCohpsnil(\widehat{\rsv}_C) )\ ,
\end{align}
where $\alpha \in  \sfH^\bullet(\overline{\rsv})$ is an arbitrary lift of $\beta$. 

We define the \textit{ring of universal tautological classes} as a commutative polynomial algebra in infinitely many variables
\begin{align}\label{eq:universal-tautological}
	\bS_{\rsv,C}\coloneqq \sfH^\bullet_A\left[\underline{\ch}_{\ell,\beta_i}, \underline{\ch}_{k,1}\;\vert\; i=1, \ldots, e, \ell \geq 1, k \geq 2\right]\ ,
\end{align}
with $\beta_i \coloneqq \ch_1(\calE_i)$, where $\calE_i$ is the $i$-th indecomposable projective object in $\catP(\rsv/\ssv)$, introduced in Formula~\ref{eq:projgenZN_components}, for $1\leq i\leq e$. We equip it with a grading such that
\begin{align}
	\deg(\underline{\ch}_{\ell,\gamma})\coloneqq \deg(\gamma) + 2\ell -4\ .
\end{align} 
There is a graded ring homomorphism 
\begin{align}\label{eq:grringhom}
	\begin{split}
		\iota\colon \bS_{\rsv,C} &\longrightarrow \sfH^\bullet_A(\dstackCohpsnil(\widehat{\rsv}_C))\\
		\underline{\ch}_{\ell,\gamma} &\longmapsto \nu_{\ell,\gamma}
	\end{split}\ .
\end{align}
This in turn yields an action by cap product 
\begin{align}
	\begin{split}
		\bS_{\rsv,C} \otimes \coha_{\rsv,C}^A &\longrightarrow \coha_{\rsv,C}^A\\ 
		x \otimes c& \longmapsto x \cap c
	\end{split}\ .
\end{align}
Define a coproduct on $\bS_{\rsv,C}$ by the formula
\begin{align}
	\Delta(\underline{\ch}_{\ell,\gamma})=\underline{\ch}_{\ell,\gamma} \otimes 1 +1 \otimes \underline{\ch}_{\ell,\gamma} 
\end{align}
for all $\ell, \gamma$. 
\begin{proposition}\label{prop:coha-surface-module-algebra}
	The action of $\bS_{\rsv,C}$ on $\coha_{\rsv,C}^A$ is Hopf, i.e., we have
	\begin{align}
		c \cap (x_1 \cdot x_2)=\sum_i (c'_i \cap x_1) \cdot (c''_i \cap x_2)
	\end{align} 
	where we have used Sweedler's notation $\Delta(x)=\sum_i c'_i \otimes c''_i$.
\end{proposition}
The proof of Proposition~\ref{prop:coha-surface-module-algebra} is given in Appendix~\ref{sec:hopfact}.

\begin{remark} 
	Using the quasi-compact open exhaustion of $\dstackCohpsnil(\widehat{\rsv}_C)$ constructed in \cite[\S\ref*{foundation-subsec:admissibility}]{DPSSV-1}, one can prove the Künneth formula 
	\begin{align}
		\Hbullet_A(\overline{\rsv} \times \dstackCohpsnil(\widehat{\rsv}_C)) \simeq \sfH^{\bullet}_A(\overline{\rsv}) \widehat{\otimes} \sfH^\bullet_A(\dstackCohpsnil(\widehat{\rsv}_C))\ .
	\end{align}
	Then the class $\overline{\nu}_\alpha$ can be equivalently constructed as 
	\begin{align}
		\overline{\nu}_\alpha = \langle \alpha, \ch(\overline{\scrE})\rangle\in \sfH^\bullet_A(\dstackCohpsnil(\widehat{\rsv}_C)) \ ,
	\end{align}
	where $\langle -, - \rangle\colon \Hbullet_A (\overline{\rsv})\otimes\Hbullet_A(\overline{\rsv})  \to  \Hbullet_A(\overline{\rsv})$ is the intersection pairing. This leads to an alternative of proof of Proposition~\ref{prop:coha-surface-module-algebra}, proceeding by analogy to \cite[Proposition~1.18]{MMSV}. The above approach uses only the six operations.
\end{remark}

\subsection{Compatibility of action of tautological classes}\label{subsec:limit-tautological-classes}

Recall the ring $\bS=\Hbullet_A[p_\ell(z_i)\;\vert \; i \in I, \ell \geq 0]$ from \S\ref{sec:the-algebra-S}. It may be understood as the ring of universal tautological classes for the stack $\Lambda_\qv$. There is a representation of $B_\qv$ on $\bS$ and an action $\cap\colon \bS \otimes \coha_\qv^A \to \coha_\qv^A$ which descends to an action of $\bS$ on $\coha^A_{\Ltheta,\kappa}$ for any interval $\kappa$, in particular on $\coha^A_{\Ltheta,(\ell)}$ for any $\ell \geq 0$. More generally, if $\ell <0$, the stack $\elldLambdaqv$ parametrizing objects in $\calP^{\mathsf{nil}}_{\Ltheta}((\nu_\ell,1+\nu_\ell])$ carries tautological complexes $\scrV_i$ for $i \in I$; we may define an action of $\bS$ on $\coha^A_{\Ltheta,(\ell)}$ via
\begin{align}
	\ch(z_i) =p_0(z_i) + \sum_{\ell \geq 1} \frac{p_\ell(z_i)}{\ell!} \longmapsto \ch(\calV_i)\ . 
\end{align}
By construction, these actions of $\bS$ on $\coha^A_{\Ltheta,(\ell)}$ are compatible with each other: the restrictions of the actions of $\bS$ on $\coha^A_{\Ltheta,(\ell)}$ and $\coha^A_{\Ltheta,(k)}$ coincide on the stack parametrizing objects in $\calP^{\mathsf{nil}}_{\Ltheta}(J)$, where $J = (\nu_\ell, \nu_\ell +1] \cap  (\nu_k, \nu_k +1]$. It entails that there is an action 
\begin{align}
	\cap\colon\bS \otimes \coha_\infty^A \longrightarrow \coha_\infty^A\ .
\end{align}

From Lemma~\ref{lem:action-tautological-classes} we deduce the following:
\begin{lemma} For any $\ell \in \Z, k \in \N$, the diagram
	\begin{align}\label{eq:coha-chern-braid}
		\begin{tikzcd}[ampersand replacement=\&, column sep=large]
			\bS \otimes \coha_{\Ltheta,(\ell)}^A \ar{rr}{T_{\textrm{-}2k\Lthetafin} \otimes \R S_{\textrm{-}2k\Lthetafin}} \ar{d}{\cap} \& \&\bS \otimes \coha_{\Ltheta,(\ell-k)}^A \ar{d}{\cap} \\
			\coha_{\Ltheta,(\ell)}^A \ar{rr}{\R S_{\textrm{-}2k\Lthetafin}} \& \& \coha_{\Ltheta,(\ell-k)}^A
		\end{tikzcd}
	\end{align}
	is commutative, where the braid group action on $\bS$ is introduced \S\ref{subsubsec:braid-group-tautological}.
\end{lemma}
Note that $\rsv$ is equivariantly formal, i.e., one has an isomorphism 
\begin{align}\label{eq:equivformal}
	\sfH^\bullet_A(\rsv,\Q) \simeq \sfH^\bullet(\rsv,\Q)\otimes \sfH^\bullet_A\ .
\end{align}
Furthermore, the cohomology $\sfH^2(\rsv)$ is freely generated by the classes $\beta_i \coloneqq \ch_1(\calE_i)$, where $\calE_i$ is the 
$i$-th indecomposable projective object in $\catP(\rsv/\ssv)$, introduced in Formula~\ref{eq:projgenZN_components}, for $1\leq i\leq e$. Using the isomorphism \eqref{eq:equivformal}, let us write 
\begin{align}
	\mathsf{Td}_\rsv= u \cdot 1 +  \sum_i u_i\cdot \beta_i\ , 
\end{align}
where $u, u_1, \ldots, u_e \in \sfH^\bullet_A$. Because $\rsv$ is symplectic and non-proper, we have $\mathsf{Td}_X=1$ if $A=\{\id\}$, hence $u \in 1+\sfH^{>0}_A, u_i \in \sfH^{>0}_A$. Setting $\underline{\ch}_\gamma=\sum_\ell \underline{\ch}_{\ell,\gamma}$, there is a unique algebra isomorphism
\begin{align}
	\sigma \colon \bS \longrightarrow \bS_{\rsv,C} 
\end{align}
satisfying
\begin{align}\label{eq:identification_taut-classes}
	\sigma({\ch}_{z_i})&= r_iu\underline{\ch}_1 + \sum_{j=1}^e r_i u_j \underline{\ch}_{\beta_j}-u\underline{\ch}_{\beta_i}, \qquad (i =1, \ldots ,e)\\
	\sigma({\ch}_{z_0})&=u\underline{\ch}_1 +\sum_{j=1}^e u_j \underline{\ch}_{\beta_j}\ ,
\end{align}
where $r_i=\mathsf{rk}(\calE_i)$. We use it to identify $\bS_{\rsv,C}$ and $\bS$.

\begin{proposition}\label{prop:coha-tautological class theta compat}
	The map $\Theta$ is $\bS$-equivariant, i.e., the following diagram is commutative
	\begin{align}\label{eq:isom-chern-braid}
		\begin{tikzcd}[ampersand replacement=\&, column sep=large]
			\bS_{\rsv,C} \otimes \coha_{\rsv,C}^A \ar{rr}{\sigma \otimes \Theta} \ar{d}{\cap} \& \&\bS \otimes \coha_\infty^A \ar{d}{\cap} \\
			\coha_{\rsv,C}^A \ar{rr}{\Theta} \& \& \coha_\infty^A
		\end{tikzcd}\ .
	\end{align}
\end{proposition}
\begin{proof}
	By construction, the isomorphism $\Theta$ is constructed using the quasi-compact open exhaustion of $\dstackCohpsnil(\widehat{\rsv}_C)$ by the substacks $\dstackCohpsnil(\widehat{\rsv}_C)^{<\ell}$ for $\ell \in \Z$, together with the isomorphisms $\tau_\ell\colon  \dstackCohpsnil(\widehat{\rsv}_C)^{<\ell} \simeq \elldLambdaqv$. The algebras $\bS_{\rsv,C}$ and $\bS$ act on $\sfH_\bullet^A(\dstackCohpsnil(\widehat{\rsv}_C)^{<\ell})$ and $\sfH_\bullet^A(\elldLambdaqv)$, respectively, by cap product with some explicit tautological cohomology classes. To prove the Proposition, it thus suffices to show that the morphism $\sigma$ identifies these two families of tautological classes.
	
	By Theorem~\ref{thm:tau}, 
	\begin{align}
		\scrV_i=\tau_\ast(\ev_\ast (\ev_{\rsv}^\ast \calE_i^\vee \otimes \scrE)) \in \sfG_0^A(\elldLambdaqv)
	\end{align}
	for any $i \in I$, where $\ev$ and $\ev_{\rsv}$ are the projections from $\dstackCohpsnil(\widehat{\rsv}_C)\times \rsv$ to the two factors, respectively. . For $0\leq i \leq e$, let $\overline{\calE}_i\in \catDb_A(\overline{\rsv})$ be a $A$-equivariant perfect complex extending $\calE_i$ to $\overline{\rsv}$. Using the Grothendieck-Riemann-Roch formula \cite[Corollary~3.25]{Khan_VFC}, we get
	\begin{align}
		\ch(\scrV_i)=\tau_\ast(\ch( \ev_\ast (\ev_{\rsv}^\ast \calE_i^\vee \otimes \scrE) ))=\tau_\ast \ev_!( \ev_{\overline{\rsv}}^\ast(\ch(\overline{\calE}_i^\vee) \cup \mathsf{Td}_{\overline{\rsv}})\cup \ch(\overline{\scrE}))\ .
	\end{align}
	Using the isomorphism \eqref{eq:equivformal}, a straightforward computation yields
	\begin{align}
		h^\ast(\ch(\overline{\calE}_i^\vee) \cup \mathsf{Td}_{\overline{\rsv}}) = r_i u\cdot 1 + \sum_{j=1}^e r_iu_j\beta_j - u \beta_i\ .
	\end{align}
	Therefore, applying Corollary~\ref{cor:restoSB}, one obtains 
	\begin{align}
		\ev_!( \ev_{\overline{\rsv}}^\ast(\ch(\overline{\calE}_i^\vee) \cup \mathsf{Td}_{\overline{\rsv}}) = r_iu \nu_{1} + \sum_{j=1}^e r_iu_j\nu_{\beta_j} - u \nu_{\beta_i}\ .
	\end{align}
	Given Formula~\eqref{eq:grringhom}, this proves the claim. 
\end{proof}

Recall from \S\ref{sec:the-algebra-S} that there is also an action of $\bS$ on $\Y_\qv$, provided by the Formulas~\eqref{eq:S-action-Yangian}. By Theorem~\ref{thm:Phi}, the isomorphisms $\Phi_{(\ell)}$ are $\bS$-equivariant. We thus obtain the following direct corollary.
\begin{corollary}
	The isomorphism $\Theta_{\rsv,C}\colon\coha_{\rsv, C}^A \to \Y_\infty^+$ is $\bS$-equivariant.
\end{corollary}

\subsection{Explicit description of $\coha_{\rsv, C}$}\label{subsec:classical-limit}

At the moment we do not know a full presentation of $\Y_\infty^+$ (or of $\coha_{\rsv, C}^A$) by generators and relations. However, as we will see in this Section, it is possible to describe explicitly $\coha_{\rsv, C}^A$ when $A=\{\id\}$. As the forgetful map $\HBMbulletA(\dLambda_\qv) \otimes_{\Hbullet_A} \Q \simeq \HBMbullet(\dLambda_\qv)$ is functorial, it induces an isomorphism of algebras $\HBMbulletA(\elldLambdaqv) \otimes_{\Hbullet_A} \Q \simeq \HBMbullet(\elldLambdaqv)$ for any $\ell\in\Z$, hence by Theorem~\ref{thm:coha-surface-as-limit1} we get
\begin{align}
	\coha_{\rsv, C}^A\otimes_{\Hbullet_A} \Q \simeq \coha_{\rsv, C}\ .
\end{align}
In view of Propositions~\ref{prop:classical-limit} and \ref{prop:iso-s-uce-1}, this gives an isomorphism of algebras
\begin{align}\label{eq:cohaquiver-classical-limit}
	\cohaqv \simeq \sfU(\fraknell)\ .
\end{align}
Since Theorem~\ref{thm:coha-surface-as-limit1} works regardless of torus actions, the non-equivariant versions of Theorems~\ref{thm:coha-surface-as-limit1} and \ref{thm:coha-surface-as-limit2} yield an isomorphism of topological algebras
\begin{align}\label{eq:cohasurface-classical-limit}
	\coha_{\rsv, C} \simeq \lim_\ell \Y_{\Ltheta,(\textrm{-}2\ell)}\otimes_{\Hbullet_A} \Q\ ,
\end{align}
where the left-hand-side is now the non-equivariant cohomological Hall algebra $\coha_{\rsv, C}$. 

We now use the isomorphisms~\eqref{eq:cohaquiver-classical-limit} and \eqref{eq:cohasurface-classical-limit} to obtain an explicit description of $\coha_{\rsv, C}$. From \S\ref{sec:quotients-Yangians}, we have 
\begin{align}
	\Y_{\Ltheta,(0)}\otimes_{\sfH^\bullet_A} \Q \simeq \sfU((\fraknell)_{\Ltheta,(0)})\quad \text{where} \quad (\fraknell)_{\Ltheta,(0)}\coloneqq\bigoplus_{\beta \in \Delta_{(0)}} (\fraknell)_{\beta}\ ,
\end{align}
with
\begin{align}
	\Delta_{(0)}\coloneqq \Big\{\alpha +n\delta \;\Big\vert\; \alpha \in \Delta^+_\sff \cup \{0\},\; n < 0\Big\}\ .
\end{align}
For $\ell\in \N$, the subquotient $\Y_{\Ltheta,(\textrm{-}2\ell)}\otimes_{\sfH^\bullet_A} \Q$ of $\sfU(\frakgell)$ is obtained from $\Y_{\Ltheta,(0)}\otimes_{\sfH^\bullet_A} \C$ by applying the automorphism $T_{\textrm{-}4\ell\Lthetafin}$. It follows that
\begin{align}
	\Y_{\Ltheta,(\textrm{-}2\ell)}\otimes_{\sfH^\bullet_A} \Q \simeq \sfU((\fraknell)_{\Ltheta,(\textrm{-}2\ell)}) \quad \text{where} \quad (\fraknell)_{\Ltheta,(\textrm{-}2\ell)}\coloneqq \bigoplus_{\beta \in \Delta_{\Ltheta, (\textrm{-}2\ell)}} (\fraknell)_{\beta}\ ,
\end{align}
with
\begin{align}
	\Delta_{\Ltheta, (\textrm{-}2\ell)}\coloneqq\big\{\alpha +n\delta \;\big\vert\; \alpha \in \Delta^+_\sff \cup \{0\},\; n -4\ell (\Lthetafin,\alpha) <0\big\}\ .
\end{align}

Define
\begin{align}
	\frakr_{(\textrm{-}2\ell)}\coloneqq\bigoplus_{\beta \in \Delta_{\Ltheta, (\textrm{-}2\ell)} \smallsetminus \Delta_{\Ltheta, (\textrm{-}2(\ell-1))}} (\fraknell)_{\beta}
\end{align}
so that $(\fraknell)_{\Ltheta,(\textrm{-}2\ell)}=(\fraknell)_{\Ltheta,(\textrm{-}2(\ell-1))} \oplus \frakr_{(\textrm{-}2\ell)}$. The PBW isomorphism 
\begin{align}
	\sfU((\fraknell)_{\Ltheta,(\textrm{-}2\ell)}) \simeq  \sfU(\frakr_{(\textrm{-}2\ell)})\otimes \sfU((\fraknell)_{\Ltheta,(\textrm{-}2(\ell-1))})
\end{align}
gives rise to a projection morphism $\sfU((\fraknell)_{\Ltheta,(\textrm{-}2\ell)}) \to \sfU(\frakn_{\Ltheta,(\textrm{-}2(\ell-1))})$, and we have 
\begin{align}
	\coha_{\rsv, C} \simeq \Y_\infty^+\otimes_{\sfH^\bullet_A} \Q\simeq \lim_\ell \sfU((\fraknell)_{\Ltheta,(\textrm{-}2\ell)})\ .
\end{align}

Note that, as subsets of $\Delta$, we have $\Delta_{\Ltheta,(0)} \subset \Delta_{\Ltheta,(\textrm{-}2)} \subset \cdots $ and
\begin{align}
	\Delta_{\Ltheta}\coloneqq \bigcup_\ell \Delta_{\Ltheta,(\textrm{-}2\ell)}= (\Delta^+_\sff \times \Z \delta) \cup (-\N \delta)\ .
\end{align}
In terms of the isomorphism between $\frakgell$ and the universal central extension of $\frakgfin[s^{\pm 1}, t])$, we have
\begin{align}\label{eq:Ln'}
	\fraknell^+\coloneqq \bigoplus_{\beta \in \Delta_{\Ltheta}} (\fraknell)_\beta= \fraknfin [s^{\pm 1},t] \oplus s^{-1}\frakhfin [s^{-1},t] \oplus \bigoplus_{k <0} \Q c_{k,\ell}\ .
\end{align} 

As we next explain, the limit $\lim_\ell \sfU((\fraknell)_{\Ltheta,(\ell)})$ may also be described as a completion of the subalgebra $\sfU(\fraknell^+) \subset \sfU(\frakgell)$. For this we shall use the following general lemma. 
\begin{lemma}\label{lem:completion-graded-Lie-algebra}
	Let $\frakm =\bigoplus_{\alpha \in \Gamma}\frakm_\alpha$ be a Lie algebra graded by a subset $\Gamma$ of a $\Z$-lattice $L$. Let $r,d\colon L \to \Z$ be two linear forms, and assume that $r(\Gamma) \subset \N$. Set $\mu(\alpha)\coloneqq d(\alpha)/r(\alpha)$ and set $\frakm_{<k}\coloneqq\bigoplus_\alpha\frakm_\alpha$, where $\alpha$ runs over the set of all weights such that $\mu(\alpha)<k$. Then, the multiplication in $\sfU(\frakm)$ yields a complete topological algebra structure on the topological vector space
	\begin{align}
		\widehat{\sfU}(\frakm)\coloneqq \lim_k \sfU(\frakm)/\left( \frakm_{<k}\cdot \sfU(\frakm) \right)\ ,
	\end{align}
	the limit being taken in the category of topological $\Z I$-graded vector spaces.
\end{lemma}

\begin{proof}
	Choose a basis $\{x_r\,\vert\,r\in\N\}$ of $\frakm$ consisting of weight vectors such that $x_r\in\frakm_{\alpha_r}$, with
	\begin{align}
		\mu(\alpha_0)\geq\mu(\alpha_1)\geq\mu(\alpha_2)\geq\cdots
	\end{align}
	We call \textit{normal monomial} an element of $\sfU(\frakm)$ of the form
	\begin{align}
		x_{r_1}x_{r_2}\cdots x_{r_i}
	\end{align}
	with $r_1\geq r_2\geq\cdots\geq r_i\geq 0$. The normal monomials form a basis of $\sfU(\frakm)$.
	
	Let $u=(u_\ell)$ and $v=(v_\ell)$ be elements of $\widehat{\sfU}(\frakm)$ of weight $\alpha,\beta\in\Gamma$ respectively. Let $\tilde u_\ell$ and $\tilde v_\ell$ be lifts of $u_\ell$ and $v_\ell$ in $\sfU(\frakm)$ of weight $\alpha$ and $\beta$, respectively, for each integer $\ell\leq 0$. We consider the obvious projection 
	\begin{align}
		\pi_\ell\colon \sfU(\frakm)\longrightarrow \sfU(\frakm)\,/\,(\,\frakm_{<\ell}\cdot \sfU(\frakm))\ .
	\end{align}
	We must check that for each large enough integer $m$ (depending on $\ell$, $\alpha$ and $\beta$), the element 
	\begin{align}
		\pi_\ell(\tilde u_{\ell-m} \cdot \tilde v_{\ell-m})
	\end{align}
	does not depend on $m$. To do that, we fix for every $n \leqslant 0$ elements $x_n$ and $y_n$ in $\frakm_{<n}\cdot \sfU(\frakm)$ of degree $\alpha$ and $\beta$ respectively, and we prove that
	\begin{align}
		\pi_\ell\big((\tilde u_{\ell-m}+x_{\ell-m})\cdot (\tilde v_{\ell-m}+y_{\ell-m})\big)=\pi_\ell(\tilde u_{\ell-m} \cdot \tilde v_{\ell-m})
	\end{align}
	for $m\gg 0$. Thus, we must check that
	\begin{align}
		\tilde u_{\ell-m}\cdot y_{\ell-m}\in\frakm_{<\ell}\cdot \sfU(\frakm)
	\end{align}
	for $m\gg 0$. We may assume that $\tilde u_{\ell-m}=x_r$ and $y_{\ell-m}=x_s$, where $x_r=x_{r_1}x_{r_2}\cdots x_{r_i}$ and $x_s=x_{s_1}x_{s_2}\cdots x_{s_j}$ are normal monomials of weight $\alpha$ and $\beta$ with $\mu(\alpha_{s_1})<\ell-m$. Write
	\begin{align}
		x_rx_s=x_{s_1}x_{r_1}x_{r_2}\cdots x_{r_i}x_{s_2}\cdots x_{s_j}+
		[x_{r},x_{s_1}]x_{s_2}\cdots x_{s_j}\ .
	\end{align}
	The first term belongs to 
	\begin{align}
		\frakm_{<\ell-m}\cdot \sfU(\frakm)\subset \frakm_{<\ell}\cdot \sfU(\frakm)\ .
	\end{align}
	We claim that $[x_r,x_{s_1}]$ belongs also to $\frakm_{<\ell}\cdot \sfU(\frakm)$ if $m$ is large enough. Indeed, the element $[x_r,x_{s_1}]$ in $\sfU(\frakm)$ is a linear combination of normal monomials. Let $x_{t_1}x_{t_2}\cdots x_{t_k}$ be any normal monomial entering the decomposition of $[x_r,x_{s_1}]$. We must check that $\mu(\alpha_{t_1})<\ell$ if $m$ is large enough. The slope function is convex, i.e., we have
	\begin{align}
		\mu(\sigma_1)\leq \mu(\sigma_2)\Rightarrow \mu(\sigma_1)\leq \mu(\sigma_1+\sigma_2)\leq \mu(\sigma_2)\ ,
	\end{align}
	for all $\sigma_1,\sigma_2\in\Gamma$. Hence, the claim, and the proposition, follows from the computation
	\begin{align}
		\mu(\alpha_{s_1})\leq\mu(\alpha_{s_1}+\alpha)=\frac{r(\alpha_{s_1})\mu(\alpha_{s_1})}{r(\alpha_{s_1}+\alpha)}+
		\frac{r(\alpha)\mu(\alpha)}{r(\alpha_{s_1}+\alpha)}\leq\frac{\ell-m}{r(\alpha+\beta)}+\mu(\alpha)\ .
	\end{align}
	In the above, we used the fact that $r(\Gamma) \in \N$ so that $r(\alpha_{s_1}) \leq r(\alpha_{s_1}+\alpha)$ and $ r(\alpha) \leq r(\alpha_{s_1}+\alpha)$.
\end{proof}

We summarize the main results of this section in the following Theorem. Let $\widehat{\sfU}(\fraknell^+)$ be the completion, in the sense of Lemma~\ref{lem:completion-graded-Lie-algebra}, of $\sfU(\fraknell^+)$, with respect to the slope function $\mu_{\Ltheta}$. Recall that $\fraknell^+$ is introduced in Formula~\eqref{eq:Ln'}.
\begin{theorem}\label{thm:coha-surface-as-limit3} 
	There is a canonical isomorphism of complete topological algebras
	\begin{align}
		\begin{tikzcd}[ampersand replacement=\&]
			\Theta_{\rsv,C}\colon\coha_{\rsv, C} \ar{r}{\sim}\& \widehat{\sfU}(\fraknell^+)
		\end{tikzcd}\ .
	\end{align}
	Moreover, $\Theta_{\rsv,C}$ intertwines the action of tautological classes in the same sense as in \S\ref{subsec:limit-tautological-classes}. Likewise, it intertwines the actions of the lattice $\coweightlatticefin$ by tensor product and by braid operators in the sense of Proposition~\ref{prop:coha-tensor-product-intertwine}.
\end{theorem}

In the next Section, we will describe explicitly the image under $\Theta_{\rsv,C}$ of a set of generators of $\coha_{\rsv, C}$.   

\begin{remark}
	The above computation also allows one to understand the classical limit, in the sense of the standard filtration, of $\coha_{\rsv, C}^A$. More precisely, the standard filtration of $\Y_\qv$ induces one on both $\Y^-_\qv$ and the quotients $\Y_{\Ltheta,(2\ell)}$ for $\ell \geq 0$. The transition map $\Y_{\Ltheta,(0)} \to \Y_{\Ltheta,(2)}$ is compatible with this filtration because the braid operators $T_i$ and and their truncated versions $\overline{T}_i$ are, by Formula~\eqref{eq:Ti}, compatible with this filtration. This induces by Remark~\ref{rem:third-limit-realization-coha} a filtration on $\coha_{\rsv, C}$. We then have an isomorphism
	\begin{align}
		\begin{tikzcd}[ampersand replacement=\&]
			\gr\, \Theta_{\rsv,C}\colon\gr\coha_{\rsv, C}^A  \ar{r}{\sim}\& \gr \Y_\infty^+\simeq \lim_\ell\gr \Y_{\Ltheta,(2\ell)} \simeq \lim_\ell \sfU((\fraknell)_{\Ltheta,(2\ell)}\big)\otimes \Hbullet_A \simeq \widehat{\sfU}(\fraknell^+)\otimes \Hbullet_A 
		\end{tikzcd}\ . \tag*{\qedhere} 
	\end{align}
\end{remark}

\section{Some explicit computations of fundamental classes}\label{sec:explicit-computations}

In this Section we compute explicitely the image under $\Theta,\Theta_{\rsv,C}$ of certain geometrically defined elements of $\coha_{\rsv, C}^A$, namely the subalgebras corresponding to \textit{zero-dimensional} Hecke modifications and to certain fundamental classes of irreducible components of the stack $\dstackCohpsnil(\widehat{\rsv}_C)$. As we expect that these generate $\coha_{\rsv, C}^A$ under the operations of Hall multiplication and multiplication by tautological classes, this in principle fully characterizes the isomorphism $\Theta, \Theta_{\rsv,C}$. 

\medskip

Under the tilting equivalence $\tau$ we identify $G_0(\catCoh_C(\rsv))$ with $\Z I$, see Theorem~\ref{thm:nilpotent}. We denote by $\dstackCohpsnil(\widehat{\rsv}_C;\bfd)$ the substack of $\dstackCohpsnil(\widehat{\rsv}_C)$ parametrizing those sheaves with K-theory class $\bfd\in \Z I$. 

\subsection{The subalgebra corresponding to zero-dimensional sheaves}

\begin{notation}
	We denote by 
	\begin{align}
		(\coha^A_{\rsv,C})_{\N\delta}\coloneqq\bigoplus_{\ell\in \N} \sfH_\bullet^A(\dstackCohpsnil(\widehat{\rsv}_C;\ell\delta))\ ,
	\end{align}
	the cohomological Hall algebra of zero-dimensional sheaves on $\rsv$ with set-theoretic support in $C$. 
\end{notation}
As being zero-dimensional and semistable of slope $+\infty$ are equivalent conditions for a nilpotent coherent sheaf on $\rsv$, the algebra isomorphism $\Theta$ restricts to an algebra isomorphism
\begin{align}
	\begin{tikzcd}[ampersand replacement=\&]
		\Theta_{\N \delta}\colon (\coha^A_{\rsv,C})_{0,\bullet} \ar{r}{\sim} \& \coha^A_{\Ltheta,\{0\}}=\HBMbulletA(\dLambda^0_\qv)
	\end{tikzcd}\ ,
\end{align}
where the right-hand-side is the equivariant nilpotent cohomological Hall algebra of semistable nilpotent representations of slope $0$. In terms of Yangians, we can write things rather explicitly. Since $\Theta_{\rsv,C}$ is compatible with the grading, we have 
\begin{align}
	(\coha^A_{\rsv,C})_{\N\delta}\simeq \bigoplus_{\ell \in \N} \Big(\Y^+_{\infty}\Big)_{\ell\delta}\ .
\end{align}
To unburden notations, set 
\begin{align}
	\Y^-_{\N\delta}\coloneqq\bigoplus_{\ell \geq 0} \Y^-_{\ell\delta}\ .
\end{align}
First note that for any $\kappa=(a,0]$ we have
\begin{align}
	\J_{\Ltheta,\kappa} \cap \Y^-_{\N\delta}=\sum_{\substack{\bfe+\bfd\in\N\delta\\ \mu_{\Ltheta}(\bfe)<0,\\ \mu_{\Ltheta}(\bfd)>0}} \Y^-_\bfe \cdot \Y^-_\bfd\ ,
\end{align}
which is in fact independent of $a<0$. Note that $\J_{\Ltheta,\kappa} \bigcap \Y^-_{\N\delta}$ is a two-sided ideal in $\Y^-_{\N\delta}$, hence 
\begin{align}
	\Big(\Y^-_{\Ltheta,(2\ell)}\Big)_{\N\delta}\coloneqq\bigoplus_k \Big(\Y^-_{\Ltheta,(2\ell)}\Big)_{k\delta}
\end{align}
is naturally a graded algebra. Furthermore, 
\begin{align}
	T_{4\Lthetafin}\colon \Big(\Y_{\Ltheta,(2\ell)}\Big)_{\N\delta} \longrightarrow \Big(\Y_{\Ltheta,(2\ell+2)}\Big)_{\N\delta}
\end{align}
is a graded algebra isomorphism for any $\ell$. This proves the following.
\begin{proposition}
	$\Theta_{\rsv,C}$ restricts to an algebra isomorphism
	\begin{align}
		\begin{tikzcd}[ampersand replacement=\&]
			\Big(\Theta_{\rsv,C}\Big)_{\N\delta}\colon (\coha^A_{\rsv,C})_{\N\delta} \ar{r}{\sim}\& \Y^-_{\N\delta}\big/\big((\Y^-_{<0}\cdot \Y^-_{>0}) \cap \Y^-_{\N\delta}\big)
		\end{tikzcd} \ ,
	\end{align}
	where 
	\begin{align}
		\Y^-_{<0}\coloneqq\sum_{\genfrac{}{}{0pt}{}{\bfe \in \N I}{\mu_{\Ltheta}(\bfe)<0}} \Y^-_\bfe, \qquad \text{and}\quad \Y^-_{>0}\coloneqq\sum_{\genfrac{}{}{0pt}{}{\bfd \in \N I}{\mu_{\Ltheta}(\bfd)>0}} \Y^-_\bfd\ .
	\end{align}
\end{proposition}
When $A=\{\id\}$, the situation is even simpler: Theorem~\ref{thm:coha-surface-as-limit3} yields the following.
\begin{corollary}\label{cor:non-equivariant-zero-dimensional}
	 We have
	\begin{align}
		\begin{tikzcd}[ampersand replacement=\&]
			\Big(\Theta_{\rsv,C}\Big)_{\N\delta}\colon (\coha_{\rsv,C})_{\N\delta}\ar{r}{\sim} \& \sfU(\frakhfin \, s^{-1}[s^{-1},t] \oplus K_- )
		\end{tikzcd}\ ,
	\end{align}
	where $K_-$ is as in Formula~\eqref{eq:negative-half}. 
\end{corollary}

\subsection{Fundamental classes of irreducible components} 

Recall that by \cite[Corollary~\ref*{foundation-cor:Cohnil_as_formal_completion-1}]{DPSSV-1}, the canonical morphism $\dstackCohpsnil(\widehat{\rsv}_C) \to \dstackCohps(\rsv)$ is formally étale. Thus, we have 
\begin{align}
	\dim(\stackCohpsnil(\widehat{\rsv}_C;\bfd))=\dim(\stackCoh_{\mathsf{ps}}(\rsv;\bfd))=-(\bfd,\bfd)
\end{align}
for any $\bfd\in \Z$. On the other hand, we have a chain of closed immersions
\begin{align}\label{eq:inclusions-red-classical-der}
	\Red{\dstackCohpsnil(\widehat{\rsv}_C;\bfd)} =\Red{\stackCohpsnil(\widehat{\rsv}_C;\bfd)} \subset \stackCohpsnil(\widehat{\rsv}_C;\bfd) \subset \dstackCohpsnil(\widehat{\rsv}_C;\bfd)
\end{align}
which induce isomorphisms in Borel-Moore homology (cf.\ \cite[Remark~\ref*{torsion-pairs-rem:BM-reduced-stack}]{DPS_Torsion-pairs}). As described in \cite[\S\ref*{foundation-subsec:admissibility}]{DPSSV-1}, an admissible open exhaustion of $\dstackCohpsnil(\widehat{\rsv}_C, \bfd)$ is formed by the open substacks $\dstackCohpsnil(\widehat{\rsv}_C, \bfd)^{>\ell}$ parametrizing those sheaves with minimal slope greater than $\ell$. These correspond to the substacks $\frakU_\ell(\widehat{\rsv}_C;\bfd)$ introduced in \cite[Formula~(\ref*{foundation-eq:U_alpha})]{DPSSV-1}. 

Tensoring by line bundles, we see that
\begin{align}
	\dstackCohpsnil(\widehat{\rsv}_C, \bfd)^{>\ell}\simeq \dstackCohpsnil(\widehat{\rsv}_C, \bfe)^{>0} \simeq \dLambda_\bfe^{\leq 0}\ ,
\end{align}
for some $\bfe \in \N I$, where the last equivalence follows from Theorem~\ref{thm:tilting-stability}. From \cite{Lusztig-Affine-Quivers-IHES}, see also \cite[Theorem~4.11]{Schiffmann-lectures-canonical-bases} it follows that $\Red{\dstackCohpsnil(\widehat{\rsv}_C;\bfd)}$ is pure dimensional of dimension $-\frac{1}{2}(\bfd,\bfd)$ (with infinitely many irreducible components). Denote by $\sfIrr_{\rsv,C}(\bfd)$ the collection of irreducible components of $\Red{\dstackCohpsnil(\widehat{\rsv}_C;\bfd)}$ (or, equivalently, of $\Red{\stackCohpsnil(\widehat{\rsv}_C;\bfd)}$). Unless $\bfd \in \N\delta$, this is an infinite set. Given the degree shift in the definition of $\coha_{\rsv,C}^A$ and its construction as a limit, we deduce that the subalgebra of $\coha_{\rsv,C}^A$ consisting of elements of vertical degree $0$ is 
\begin{align}
	\big(\coha_{\rsv,C}^A\big)_{\bullet,0}\coloneqq \bigoplus_\bfd \sfH^A_{\textrm{-}(\bfd,\bfd)}(\stackCohpsnil(\widehat{\rsv}_C;\bfd),\Q)=\bigoplus_\bfd\Big\{\sum_{\frakZ \in \sfIrr_{\rsv,C}(\bfd)} u_\frakZ [\frakZ]\;\vert \; u_\frakZ \in \Q\Big\}\ .
\end{align}
Unless otherwise specified, we will always consider the fundamental class $[\frakZ]\coloneqq [\Red{\frakZ}]$ of a \textit{reduced} irreducible component $\frakZ\in \sfIrr_{\rsv,C}(\bfd)$. 

Since $\Theta_{\rsv,C}$ is $\Z \times \Z I$-graded, it restricts to an algebra isomorphism
\begin{align}
	\begin{tikzcd}[ampersand replacement=\&]
		\big(\Theta_{\rsv,C}\big)_{\bullet,0} \colon \big(\coha_{\rsv,C}^A\big)_{\bullet,0} \ar{r}{\sim}\& \widehat{\sfU}(\fraknellzero^+)
	\end{tikzcd}\ ,
\end{align}
where 
\begin{align}\label{eq:defun'}
	\fraknellzero^+\coloneqq \fraknfin[s^{\pm 1}] \oplus s^{-1}\frakhfin [s^{-1}]\ .
\end{align}

\begin{remark}
	In passing from a derived stack to its reduction, it is useful to keep in mind the following property which is valid for any geometric derived stack: the closed immersions \eqref{eq:inclusions-red-classical-der} induce equivalences of small étale sites (see \cite[Corollary~2.2.2.9]{TV_HAG-II}, or \cite[Remark~1.6.6(a)]{Hennion-Holstein-Robalo}.
\end{remark}

\begin{remark}
	It is well known that the smooth locus of $\Red{\dLambda}_\qv$ is a Lagrangian substack of $\stackRep(\Pi_\qv)$ (this follows for instance from  \cite{Lusztig-Affine-Quivers-IHES}). However, $\Red{\dLambda}_\qv$ is not expected to be derived lci and one may wonder about the existence of a canonical intermediate (derived) stack 
	\begin{align}
		\Red{\dLambda}_\qv \subset \tensor*[^{\mathsf{lag}}]{\dLambda}{_\qv} \subset \dstackRep(\Pi_\qv)\ ,
	\end{align}
	Lagrangian in the sense of \cite{Shifted_symplectic}. The same question arises for $\Red{\dstackCohpsnil(\widehat{\rsv}_C)}$ in $\dstackCoh_{\ps}(\rsv)$. We will not need this.
\end{remark}

It is an interesting question to describe precisely the image by $\Theta_{\rsv,C}$ of the fundamental class $[\frakZ]$ of a given irreducible component $\frakZ \in \sfIrr_{\rsv,C}(\bfd)$. By analogy with the case of quivers, where this basis of $\big(\coha_\qv^A\big)_{\bullet,0}$ is known as the \textit{semi-canonical basis}, one cannot expect to have explicit expressions for all elements. We will however obtain precise results for zero-dimensional sheaves and sheaves set-theoretically supported on a single (arbitrary) component of $C$ and of rank one on that component. As we expect that these generate $(\coha^A_{\rsv,C})_{\bullet,0}$ as a topological algebra, this fully characterizes the isomorphism $\big(\Theta_{\rsv,C}\big)_{\bullet,0}$. Note also that because of the compatibility with the action of tautological classes, we expect that $\Theta_{\rsv,C}$ is uniquely determined by its restriction $\big(\Theta_{\rsv,C}\big)_{\bullet,0}$.

Our main theorem, which will be proved in \S\ref{sec:proof-thm-class-irred}, reads as follows:
\begin{theorem}\label{thm:class_irred_component} 
	\hfill
	\begin{enumerate}\itemsep0.2cm
		\item \label{item:class_irred_component-1} For $i \in I_\sff$ and $n \in \N$, let $\frakY_{i,(n)}\in \sfIrr_{\rsv,C}(n\delta)$ be the irreducible component containing the zero-dimensional coherent sheaves on $\rsv$, of length $n$, whose support is a single point $x \in C_i$. Then 
		\begin{align}
			\Theta_{\rsv,C}([\frakY_{i,(n)}])=(-1)^{n-1}h_is^{-n}\ .
		\end{align}
		
		\item \label{item:class_irred_component-2} For $i \in I_\sff$ and $n \in \Z$, denote by $\frakZ_{i,n} \in \sfIrr_{\rsv,C}(-\alpha_i+n\delta)$ the unique irreducible component containing the coherent sheaf $\scrO_{C_i}(n-1)$ (equivalently, the irreducible component parametrizing sheaves of class $-\alpha_i+n\delta$ scheme-theoretically supported on $C_i$). Then 
		\begin{align}\label{eq:Theta'(Z)}
			\sum_{n \in \Z} (-1)^n\Theta_{\rsv,C}([\frakZ_{i,n}])u^{-n}= \left( \sum_{n \in\Z} x_i^+s^{-n}u^{-n}\right) \cdot \exp\left(\sum_{k \geq 1}h_is^{-k}\frac{u^{-k}}{k}\right)\ .
		\end{align} 
	\end{enumerate}
\end{theorem}

\begin{corollary}\label{cor:generation_HAXC}
	$\coha_{\rsv,C}^A$ is topologically generated, as an $\bS_{\rsv,C}$-algebra, by the fundamental classes $[\frakY_{i,d}]$, $[\frakZ_{i,n}]$ for $i \in I_\sff, d \in \N$ and $n\in \Z$. 
\end{corollary}

\begin{proof}
	In view of Theorem~\ref{thm:coha-surface-as-limit3}, it is enough to prove this after passing to the associated graded, i.e., in the non-equivariant situation, which we now assume is the case. Thanks to Theorem~\ref{thm:class_irred_component}, it suffices to show that the $\bS_{\rsv,C}$-subalgebra $U$ of $\widehat{\sfU}(\fraknell^+)$ generated by the Fourier modes of the right-hand-side of \eqref{eq:Theta'(Z)} is dense. By definition, $U$ contains the elements
	\begin{align}
		x_i^+s^{n} + \sum_{k \geq 1} \frac{1}{k}x_i^+s^{n+k}h_is^{-k}
	\end{align}
	for any $i,n$ as well as the subalgebra $\Q[h_is^{-\ell}]_{i,\ell}$. From this it follows that $x_i^+s^n \in \overline{U}$ for any $i$ and any $n$. Indeed, by definition of the topology on $\widehat{\sfU}(\fraknell^+)$ any sequence $(u_k)_k$ of the form 
	\begin{align}
		u_k=x_i^+s^{n+k}P_k
	\end{align}
	with $P_k=P_k(h_i, h_is^{-1}, \ldots)$ a polynomial in generators $h_i, h_is^{-1}, \ldots$ of $s$-degree $-k$ tends to zero as $k \to \infty$ and we may apply a triangular basis argument. We deduce that $\widehat{\sfU}(\fraknellzero^+) \subset \overline{U}$ and finally that $\widehat{\sfU}(\fraknell^+) \subset \overline{U}$ since $\widehat{\sfU}(\fraknell^+)$ is generated by $\bS_{\rsv,C}$ and $\widehat{\sfU}(\fraknellzero^+)$.
\end{proof}

Let $i\in I_\sff$ and let
\begin{align}
	\gamma_i \colon \fraksl_2 \longrightarrow\frakgfin \quad \text{and} \quad \gamma_{i,\mathsf{ell}}\colon \frakgl_{2,\mathsf{ell}} \longrightarrow \frakgell 
\end{align}
be the canonical embeddings of the corresponding root subalgebras. 

\begin{corollary}\label{cor:functoriality-Theta}
	For any irreducible component $C_i$ of $C_\red$, we have a commutative square
	\begin{align}\label{diag:functorialityTheta}
		\begin{tikzcd}[ampersand replacement=\&]
			\coha_{T^\ast\PP^1,\PP^1} \arrow[swap]{d}{\Theta_{T^\ast\PP^1,\PP^1}} \arrow{r}{\sim}\& \coha_{\rsv,C_i} \arrow{r}{} \& \coha_{\rsv,C} \arrow{d}{\Theta_{\rsv,C}}\\
			\widehat{\sfU}(\frakgl_{2,\mathsf{ell}}^+) \arrow{rr}{\gamma_{i,\mathsf{ell}}} \& \&\widehat{\sfU}(\frakgell^+)
		\end{tikzcd}
	\end{align}
	of algebra morphisms, where the top left (resp.\ right) arrow is induced by the isomorphism of formal neighborhoods $\widehat{T^\ast\PP^1}_{\PP^1} \simeq \widehat{\rsv}_{C_i}$ (resp.\ by the functoriality of COHAs with respect to closed embeddings).
\end{corollary}

\begin{proof}
	There are obvious algebra embeddings and identifications
	\begin{align}
		\begin{tikzcd}[ampersand replacement=\&]
		\bS_{T^\ast\PP^1,\PP^1} \arrow[swap]{d}{\simeq} \arrow{r}{\sim}\& \bS_{\rsv,C_i} \arrow{r}{} \& \bS_{\rsv,C} \arrow{d}{\simeq}\\
		\bS_{A_1^{(1)}} \arrow{rr}{} \& \& \bS_{\qv}
	\end{tikzcd}
	\end{align}
	with respect to which the morphism $\gamma_{i,\mathsf{ell}}$ is equivariant. The corollary is now a consequence of the fact that, by Theorem~\ref{thm:class_irred_component} and Corollary~\ref{cor:generation_HAXC}, the diagram \eqref{diag:functorialityTheta} commutes for a set of topological generators.
\end{proof}

\medskip

The proof of Theorem~\ref{thm:class_irred_component} hinges on a few geometric computations, which are valid in greater generality than for the pair $(\rsv,C)$ and which may be of independent interest. We begin with these.

\subsection{Geometric computations in dimension zero}\label{subsec:geometric-computations-dim-0}

In this section we work under the following assumption:
\begin{assumption}\label{assumption:symplectic-surface}
	$S$ is a smooth quasi-projective symplectic complex surface and $C \subset S$ is a reduced curve whose irreducible components $C_j$ for $j=1, \ldots, f$ are assumed to be smooth and proper. We will also assume for simplicity that the intersections $C_i\cap C_j$ are transverse. 
\end{assumption}
We denote by $(-,-)$ the \textit{Mukai pairing} on $K_0(S)_\Q$ and by $[\calF] \in K_0(S)_{\Q}$ the class of a coherent sheaf $\calF$. In accordance with the case of Kleinian resolutions, we set $\alpha_j\coloneqq -[\scrO_{C_j}(-1)]$ for $j=1, \ldots, f$ and denote by $\delta$ the class of a zero-dimensional sheaf of length one. 

\begin{remark}
	As $S$ is symplectic, the Mukai pairing $(-,-)$ is symmetric, and may be written as the symmetrization of a bilinear form $\langle-,-\rangle$. By Proposition~\ref{prop:twist-sign-antisym}, the twisted algebra $\coha_{S,C}^{\langle\textrm{-},\textrm{-}\rangle}$ is independent (up to canonical isomorphism) on the choice of $\langle-,-\rangle$. In other words, there is a `natural' twisting of the algebra structure on $\coha_{S,C}$, similar to Formula~\eqref{eq:twisting_form_ADE_resolution}.
	For simplicity, we do not consider such twisting in this paragraph. This choice is validated by the fact that the commutation relations which we consider here all involve zero-dimensional sheaves, and that the class $\delta$ of a point lies in the kernel of $(-,-)$. 
\end{remark}

\begin{remark}\label{rem:completion}
	We will often use the following easy fact: for a smooth affine variety $U$ and a smooth closed subvariety $X \subset U$, the formal completions $\widehat{U}_X$ and $\widehat{(N_{U/X})}_X$ are isomorphic.  
\end{remark}

\begin{definition}
	Let $\lambda=(\lambda_1 \geq \cdots \geq \lambda_s)$ be a partition of an integer $n$. We say that a zero-dimensional sheaf $\calT$ on $S$ is \textit{of type $\lambda=(\lambda_1 \geq \cdots \geq \lambda_s)$} if $\mathsf{Supp}(\calT)=\sum_{i=1}^s \lambda_i x_i\in \Sym^n(S)$ for distinct elements $x_1, \ldots, x_s$.
\end{definition}

For $j\in\{1, \ldots, f\}$, define
\begin{align}
	\stackCohpsnil(\widehat{S}_{C_j})_{\N\delta}\coloneqq \bigsqcup_{n\in \N}\stackCohpsnil(\widehat{S}_{C_j}; n\delta)\quad \text{and}\quad \dstackCohpsnil(\widehat{S}_{C_j})_{\N\delta}\coloneqq \bigsqcup_{n\in \N}\dstackCohpsnil(\widehat{S}_{C_j}; n\delta)\ .
\end{align}
For a partition $\lambda$, let $\Red{\stackCohpsnil(\widehat{S}_{C_j})}_{\N\delta}^\lambda$ be the substack of $\Red{\stackCohpsnil(\widehat{S}_{C_j})}_{\N\delta}$ parametrizing those zero-dimensional sheaves of type $\lambda$ and let $\frakY_{j,\lambda}$ be its closure in $\Red{\stackCohpsnil(\widehat{S}_{C_j})}_{\N\delta}$.

\begin{proposition}\label{prop:zero-dim_Higgs_sheaves_classification} 
	Under Assumption~\ref{assumption:symplectic-surface}, the following holds:
	\begin{enumerate}\itemsep0.2cm
		\item \label{item:zero-dim_Higgs_sheaves_classification-1} For $j\in \{1, \ldots, f\}$, the reduced substacks $\frakY_{j,\lambda}$, by varying of the partition $\lambda$, form a full collection of irreducible components of $\Red{\stackCohpsnil(\widehat{S}_{C_j})}_{\N\delta}$.They are all of dimension zero.
		
		\item \label{item:zero-dim_Higgs_sheaves_classification-2} There is an algebra isomorphism $\Psi_j\colon \Q[z_1,z_2,\ldots]^{\mathfrak{S}_{\infty}} \to \big(\coha_{S,C_j}\big)_{\N\delta,0}$ such that for any $\lambda$, $\Psi_j(m_\lambda)=[\frakY_{j,\lambda}]$. Here $(m_\lambda)_\lambda$ stands for the basis of monomial symmetric functions.
		
		\item \label{item:zero-dim_Higgs_sheaves_classification-3} The canonical map $\big(\coha_{S,C_j}\big)_{\N\delta,0} \to \big(\coha_{S,C}\big)_{\N\delta,0}$ is injective, and there is an isomorphism of (commutative) graded algebras 
		\begin{align}
			\begin{tikzcd}[ampersand replacement=\&]
				\displaystyle\bigotimes_{j=1}^s \big(\coha_{S,C_j}\big)_{\N\delta,0} \ar{r}{\sim}\& \big(\coha_{S,C}\big)_{\N\delta,0}
			\end{tikzcd}\ .
		\end{align}
	\end{enumerate}
\end{proposition}

\begin{proof} 
	Statement~\eqref{item:zero-dim_Higgs_sheaves_classification-1} is local on $C$ and is easily reduced to the case of an affine surface $S$. As $S$ is symplectic and $C_j$ smooth, we have that $\widehat{S}_{C_j}\simeq \widehat{T^\ast C_j}_{C_j}$ (cf.\ Remark~\ref{rem:completion}). Thus, $\stackCohpsnil(\widehat{S}_{C_j})_{\N\delta}$ is isomorphic to $\stackCohpsnil(\widehat{T^\ast C_j}_{C_j})_{\N\delta}$, where the latter is equivalent to the stack of zero-dimensional nilpotent Higgs sheaves on $C_j$ (cf.\ \cite[\S2.4.3]{Porta_Sala_Hall}). Statement~\eqref{item:zero-dim_Higgs_sheaves_classification-1} thus follows from the explicit description of irreducible components of the global nilpotent cone of smooth curves given in \cite[Corollary~2.5]{Bozec-nilp-cone}, together with the fact that a generic zero-dimensional sheaf on $C_j$ is isomorphic to a sum $\scrO_{x_1} \oplus \cdots \oplus \scrO_{x_n}$ of simple skyscraper sheaves with distinct supports $x_1, \ldots, x_n$.
	
	\medskip
	
	We turn to Statement~\eqref{item:zero-dim_Higgs_sheaves_classification-2}. We will compute the structure constants of the multiplication in $\big(\coha_{S,C_j}\big)_{\N\delta,0}$. To unburden the notation we drop the index $j$ throughout. Let us call \textit{uniserial} a zero-dimensional sheaf $\calT$ on $S$ which admits a unique filtration $0=\calT_0  \subset\calT_1 \subset \cdots \subset \calT_\ell=\calT$ for which $\calT_i / \calT_{i-1}$ is simple (i.e., of length one) for all $i$. This condition is satisfied by a generic coherent sheaf on $S$ with punctual support. For instance, if $S=\Spec \,\C[x,y]$ then a sheaf $\calT$ supported at $0$ is uniserial as soon as $x$ or $y$ acts on $H^0(S, \calT)$ as a regular nilpotent operator.
	
	Let us fix partitions $\lambda, \mu, \nu$ with $\vert \nu\vert =\vert \mu\vert +\vert \lambda\vert $. Let $\bfY^{\nu, \circ}$ be an open substack of $\dstackCohpsnil(\widehat{S}_C)$ parametrizing sheaves of the form $\calG_1 \oplus \cdots \oplus \calG_{\ell(\nu)}$ where $\calG_i$ is a uniserial sheaf of length $\nu_i$ for $i=1, \ldots, \ell(\nu)$. The existence of $\bfY^{\nu, \circ}$ follows e.g. from the description of the irreducible components of the classical moduli stack of nilpotent Higgs sheaves on $C$ -- see \cite[\S3]{Schiffmann_Kac-polynomial}. Consider the induction diagram	
	\begin{align}\label{eq:proof_zero-dim-irred}
		\begin{tikzcd}[ampersand replacement=\&, column sep=large]
			\dstackCohpsnil(\widehat{S}_C)_{\N\delta} \times \dstackCohpsnil(\widehat{S}_C)_{\N\delta} \& \calS_2\dstackCohpsnil(\widehat{S}_C)_{\N\delta}
			\ar[swap]{l}{q} \ar{r}{p} \& \dstackCohpsnil(\widehat{S}_C)_{\N\delta}\\
			\Red{\frakY}_\lambda \times \Red{\frakY}_\mu\ar{u}{i_\lambda \times i_\mu} \&
			\bfX \ar{u}{i_{\lambda, \mu}} \ar[swap]{l}{q_{\lambda, \nu}} \ar{r}{p_{\lambda, \mu}} \& \dstackCohpsnil(\widehat{S}_C)_{\N\delta} \ar[equal]{u}{}\\
			\& \bfX^\circ \ar{u}{j_\nu} \ar{r}{p_\nu} \& \bfY^{\nu, \circ} \ar{u}{j_\nu^\circ}
		\end{tikzcd}
	\end{align}
	where $\bfX$ (resp.\ $\bfX^\circ$) is a geometric derived stack defined so as to make the left (resp.\ bottom) square a pullback, $i_\lambda,i_\mu, i_{\lambda, \mu}$ are closed immersions and $j_\nu,j_\nu^\circ$ are open embeddings. Here, $p\coloneqq \partial_1$ and $q\coloneqq \partial_0\times \partial_2$, following the notation of Formula~\eqref{eq:partial}.
	
	Note the the irreducible components of the derived stacks appearing in the diagram~\eqref{eq:proof_zero-dim-irred} are of dimension zero. The derived lci morphism $q$ is of dimension zero, and the map $p$ is representable by proper derived schemes, of dimension zero as well. By base change we have
	\begin{align}
		(j_\nu^\circ)^\ast ([\frakY_\lambda] \star [\frakY_\mu])=	(j_\nu^\circ)^\ast p_\ast q^!(i_\lambda \times i_\mu)_\ast ([\frakY_\lambda] \otimes [\frakY_\mu])=(p_\nu)_\ast j_\nu^\ast (q_{\lambda, \nu})^!([\frakY_\lambda] \otimes [\frakY_\mu])\ .
	\end{align}
	
	From the defining property of uniserial sheaves one deduces that $\Red{p}_\nu$ is a finite morphism. Comparing dimensions, it follows that $\Red{\bfX}^\circ$ is of pure dimension zero as well. It follows that $j_\nu^\ast (q_{\lambda, \nu})^!([\frakY_\lambda] \otimes [\frakY_\mu])=[\bfX^\circ]$. Therefore we have
	\begin{align}
		(p_\nu)_\ast j_\nu^\ast (q_{\lambda, \nu})^!([\frakY_\lambda] \otimes [\frakY_\mu])=(p_\nu)_\ast([\bfX^\circ])=n_{\lambda,\mu}^\nu [\bfY^{\nu, \circ}]\ ,
	\end{align}
	where $n_{\lambda,\mu}^\nu$ is the cardinality of the generic fiber of $p_\nu$. This cardinality is the number $n_{\lambda,\mu}^\nu$ of subsheaves of type $\mu$ and cotype $\lambda$ of a sheaf $\calG_1 \oplus \cdots \oplus \calG_{\ell(\nu)}$ in $\bfY^{\nu, \circ}$, which is equal to
	\begin{align}
		\Big\vert \Big\{ (d_1, \ldots, d_{\ell(\nu)})\in \N^{\ell(\nu)}\;\Big\vert\; (d_1, \ldots, d_{\ell(\nu)}), (\nu_1-d_1, \ldots, \nu_{\ell(\nu)}-d_{\ell(\nu)})\text{ are permutations of } \mu,\lambda\Big\}\Big\vert\ .
	\end{align}
	This matches the structure constants of multiplication in $\Q[z_1,z_2,\ldots]^{\mathfrak{S}_{\infty}}$ in the basis of monomial functions, i.e.,
	\begin{align}
		m_\lambda m_\mu=\sum_\nu n_{\lambda,\mu}^\nu m_\nu\ ,
	\end{align}
	which proves \eqref{item:zero-dim_Higgs_sheaves_classification-2}. 
	
	\medskip
	
	Finally we deal with Statement~\eqref{item:zero-dim_Higgs_sheaves_classification-3}. It is not hard to see that the irreducible components of $\stackCohpsnil(\widehat{S}_C)_{\N\delta}$ are parametrized by tuples of partitions $(\lambda^j)_{j=1, \ldots, f}$: an irreducible component $\frakY_{\lambda^1, \ldots, \lambda^f}$ asssociated to a tuple of partition  and given by $(\lambda^j)_{j=1, \ldots, f}$ is the closure of the substack of $\stackCohpsnil(\widehat{S}_C)$ parametrizing zero-dimensional sheaves of the form $\calG = \calG^1 \oplus \cdots \oplus \calG^f$ such that $\calG^i$ is of type $\lambda^i$ and $\mathsf{Supp}(\calG^i)\subset C_i$, for $i=1, \ldots, f$.

	A computation similar to (and simpler than) that proving Statement~\eqref{item:zero-dim_Higgs_sheaves_classification-2} above shows that the elements $[\frakY_{j,\lambda}]$ all commute with one another (for $j=1, \ldots ,f$ and all partitions $\lambda$), and that 
	\begin{align}
		[\frakY_{\lambda^1,\ldots,\lambda^f}]= [\frakY_{1,\lambda^1}] \star \cdots \star [\frakY_{f,\lambda^f}]\ .
	\end{align}
	Statement~\eqref{item:zero-dim_Higgs_sheaves_classification-3} follows.
\end{proof}

\begin{remark}
	The above Proposition may be viewed as an extension of \cite[Chapter~7]{Nakajima_Lectures} from the Hilbert scheme of points on a surface to the full stack of zero-dimensional sheaves (both supported on a fixed curve).
\end{remark}

\begin{remark}
	The stacks $\frakY_{j,\lambda}$ are in general not smooth but we can construct smooth open substacks $\frakY_{j,\lambda}^\circ \subset \frakY_{j,\lambda}$ as follows. We treat the case $\lambda=(d)$ only, the general case being similar. 
	
	Any open affine subset $S' \subset S$ and any isomorphism $\iota\colon \widehat{S'}_{C'_j} \simeq (\widehat{T^\ast C'_j})_{C'_j}$ gives rise to a projection $\pi\colon \widehat{S'}_{C_j'} \to C_j'$. Here $C'_j\coloneqq C_j \cap S'$. Let $\frakY_{j,(d)}^{S',\iota} \subset \frakY_{j,(d)}$ be the open substack parametrizing those zero-dimensional sheaves $\calT$ such that $\iota_\ast (\calT)$ is \textit{regular}, i.e., isomorphic to the indecomposable $d$-fold self-extension $\scrO_x^{(d)}$ of a simple sheaf $\scrO_x$ for some $x \in C'_i$. In a suitable étale local chart, we may identify the pair $(S',C_i)$ with $(\A^2,\A^1)$, in which case we have an equivalence of reduced classical stacks
	\begin{align}
		\frakY_{j,(d)}^{\A^2,\iota}\simeq \Big\{(u,v) \in \frakgl_d^2\;\Big\vert\; [u,v]=0, \; v \in \C \id+\scrO^{\mathsf{reg}}\ , \ u \in \calN_d \Big\}\Big/\GL_d \simeq T^\ast \scrO^{\mathsf{reg}}/\GL_d\ ,
	\end{align}
	where $\calN_d$ is the nilpotent cone and $\scrO^{\mathsf{reg}} \subset \calN_d$ is the regular nilpotent orbit.
	
	In particular, $\frakY_{j,(d)}^{S',\iota}$ is smooth. We now set 
	\begin{align}\label{eq:def_open_Yilambda}
		\frakY_{j,(d)}^\circ\coloneqq \bigcup_{S',\iota} \frakY_{j,(d)}^{S',\iota}\ ,
	\end{align}
	where the union ranges over all pairs $(S',\iota)$ as above.
\end{remark}

\subsection{Irreducible components in dimension one}

We work under Assumption~\ref{assumption:symplectic-surface}. Fix $j\in \{1,\ldots,f\}$ and $\ell \in \Z$, and put $\alpha\coloneqq-\alpha_j+(\ell+1)\delta$ (thus $[\scrO_{C_j}(\ell)]=\alpha$ in K-theory). Let $g_j$ be the genus of $C_j$. 

For $\lambda$ a (possibly empty) partition, let $\frakZ^\lambda_{j,\ell}$ be the closure of the substack of $\Red{\stackCohpsnil(\widehat{S}_{C_j})}$ parametrizing coherent sheaves $\calF$ for which there exists $\calL \in \stackPic_{\ell-\vert \lambda\vert}(C_j)$ and $\calT\in \frakY_{j,\lambda}$ and a short exact sequence 
\begin{align}
	0 \longrightarrow \calT \longrightarrow \calF \longrightarrow  i_{C_j,\ast}(\calL)\longrightarrow 0\ .
\end{align}
When $\lambda$ is empty, we simply write $\frakZ_{j,\ell}$, which is the smooth reduced classical stack classifying coherent sheaves of class $\alpha$ scheme-theoretically supported on $C_j$.

Let us now assume given a family $\boldsymbol{\lambda}$ of partitions $\lambda^{(i)}$ for $i=1,\ldots,f$, and denote by $\frakZ_{j,\ell}^{\boldsymbol{\lambda}}\subset \Red{\stackCohpsnil(\widehat{S}_{C_j})}$ the closure of the reduced substack parametrizing sheaves of the form
\begin{align}
	\calF \oplus \bigoplus_{i \neq j} \calT_i\ ,
\end{align}
where $\calF \in \frakZ_{j,\ell-\vert \boldsymbol{\lambda}\vert}^{\lambda^{(j)}}$ and $\calT_i \in \frakY_{i,\lambda^{(i)}}$ for $i=1, \ldots,f$. Here, $\vert \boldsymbol{\lambda}\vert\coloneqq \sum_i \vert \lambda^i\vert$.

\begin{proposition}\label{prop:irreducible-components-in-rank-1} 	
	Fix $j\in \{1, \ldots, f\}$ and $\ell\in Z$. Under Assumption~\ref{assumption:symplectic-surface}, the following holds:
	\begin{enumerate}\itemsep0.2cm
		\item \label{item:irreducible-components-in-rank-1-1} The substacks $\{\frakZ^\lambda_{j,\ell}\}_\lambda$ as $\lambda$ runs over all partitions form a complete family of irreducible components of $\Red{\stackCohpsnil(\widehat{S}_{C_j})}$ of class $\alpha$. They are all of dimension $g_j-1$. 
		
		\item \label{item:irreducible-components-in-rank-1-2} The substacks $\{\frakZ^{\boldsymbol{\lambda}}_{j,\ell}\}_{\boldsymbol{\lambda}}$ as $\boldsymbol{\lambda}$ runs over all $f$-tuples of partitions form a complete family of irreducible components of $\Red{\stackCohpsnil(\widehat{S}_{C})}$ of class $\alpha$. They are all of dimension $g_j-1$. 
	\end{enumerate}
\end{proposition}

To prove the above proposition, we need a preliminary result. Let $\dstackCohps(S)^{\dim\geqslant 1}$ be the derived moduli stack of coherent sheaves of dimension $\geq 1$. It is an open substack of $\dstackCohps(S)$. Then, we define the derived stack $\dstackCohpsnil(\widehat{S}_{C_j})^{\dim\geqslant 1}$ by the pullback
\begin{align}
	\begin{tikzcd}[ampersand replacement=\&]
		\dstackCohpsnil(\widehat{S}_{C_j})^{\dim\geqslant 1}\ar{rr} \ar{d}\& \& \dstackCohps(S)^{\dim\geqslant 1}\ar{d}\\
		\dstackCohpsnil(\widehat{S}_{C_j}) \ar{r}{\bfjmathhat_\ast} \&\dstackCoh_{C,\ps}(S) \ar{r}\& \dstackCohps(S)
	\end{tikzcd}\ .
\end{align}
Then, $\dstackCohpsnil(\widehat{S}_{C_j})^{\dim\geqslant 1}$ is an open substack of $\dstackCohpsnil(\widehat{S}_{C_j})$. Similarly, we define $\dstackCohpsnil(\widehat{S}_C)^{\dim\geqslant 1}$.

\begin{lemma} 
	The substack $\stackPic_\ell(C_j)$ is smooth of dimension $g_j-1$, open in $\Red{\stackCohpsnil(\widehat{S}_{C_j})}$ for $j=1, \ldots, f$.
\end{lemma}

\begin{proof} 
	Denoting by $-\alpha_j + \delta \in K_0(S)_\Q$ the class of $\scrO_{C_j}$, it is thus enough to show that
	\begin{align}
		\Red{\stackCohpsnil(\widehat{S}_{C_j}; -\alpha_k+(\ell+1)\delta)}^{\dim\geqslant 1}
	\end{align}
	is equal to $\stackPic_\ell(C_j)$. This boils down to checking that a pure one-dimensional coherent sheaf on $S$ with set-theoretical support on $C_j$ and of class $-\alpha_j+(l+1)\delta$ is the pushforward to $S$ of a line bundle on $C_j$. This is a local statement for which we may assume that $S$ and $C_j$ are affine and $S=T^\ast C_j$. In this case, this reduces to the well-known fact that a nilpotent Higgs bundle of rank $1$ has trivial Higgs field. The smoothness and dimension computation also follow.
\end{proof}
Because of the lemma, we shall set 
\begin{align}
	\dstackPic_\ell(C_j)\coloneqq \dstackCohpsnil(\widehat{S}_{C_j}; -\alpha_k+(\ell+1)\delta))^{\dim\geqslant 1}\ .
\end{align}
 Its reduction is $\stackPic_\ell(C_j)$.

\begin{proof}[Proof of Proposition~\ref{prop:irreducible-components-in-rank-1}]
	We start with Statement~\eqref{item:irreducible-components-in-rank-1-1}, in the case of the empty partition.
	
	Let us consider the partial Harder-Narasimhan stratification
	\begin{align}
		\Red{\stackCohpsnil(\widehat{S}_{C_j})}=\bigsqcup_{d \geq 0} \frakX_d\ ,
	\end{align}
	where $\frakX_d$ is the locally closed substack of $\Red{\stackCohpsnil(\widehat{S}_{C_j};\alpha)}$ parametrizing sheaves $\calF$ of class $\alpha$ whose canonical zero-dimensional subsheaf $\calF^0$ is of length $d$. The semisimplification morphism 
	\begin{align}
		\pi_d\colon \frakX_d \longrightarrow \Red{\stackCohpsnil(\widehat{S}_{C_j};\alpha-d\delta)} \times \Red{\stackCohpsnil(\widehat{S}_{C_j};d\delta)}
	\end{align} 
	sending a coherent sheaf $\calF$ to the pair $(\calF^1\coloneqq\calF/\calF^0, \calF^0)$, is a vector bundle stack of rank $(\alpha-d\delta,d\delta)=0$. Indeed, for any pair $(\calF^1,\calF^0)$ consisting of a pure one-dimensional sheaf $\calF^1$ and a zero-dimensional sheaf $\calF^0$ we have $\Ext^2(\calF^1,\calF^0)=\Hom(\calF^0,\calF^1\otimes \omega_S)^\vee=0$ hence the complex $\R\Hom(\calF^1,\calF^0)[1]$ has homological tor-amplitude $[0,1]$. From this we deduce that for any partition $\lambda$ of $d$, the substack
	\begin{align}
		\pi_d^{-1}(\stackPic_{\ell-d}(C_j) \times \frakY_{j,\lambda})\coloneqq \big(\stackPic_{\ell-d}(C_j) \times \frakY_{j,\lambda}\big)\times_{\Red{\stackCohpsnil(\widehat{S}_{C_j};\alpha-d\delta)} \times \Red{\stackCohpsnil(\widehat{S}_{C_j};d\delta)}} \frakX_d
	\end{align}
	is irreducible and of dimension $g_j-1$. Statement \eqref{item:irreducible-components-in-rank-1-1} easily follows. 
	
	To prove Statement~\eqref{item:irreducible-components-in-rank-1-2} it suffices to observe that any coherent sheaf $\calF$ of class $\alpha$ splits as a direct sum of a coherent sheaf (set-theoretically) supported on $C_j$ of dimension one and zero-dimensional sheaves $\calT_i$ supported on $C_i$ for $i \neq j$.
\end{proof}


\subsection{Some dimension zero -- dimension one commutation relations}\label{sec:some-dim-0-dim-1-commutations-rel}

Still keeping the setup of the previous paragraphs, we now turn to commutation relations involving stacks of one-dimensional and zero-dimensional sheaves.

\begin{proposition}\label{prop:commutations-in-rank-0-and-rank-1}
	Under Assumption~\ref{assumption:symplectic-surface}, for $i, j\in \{1, \ldots, f\}$, $\ell\in \Z$, and $d\in \N$, we have
	\begin{itemize}
		\item if $i\neq j$:
		\begin{align}\label{eq:commutations-relations-rank-0-rank-1-i-j}
			\left[ [\frakY_{i,(d)}], [\frakZ_{j,\ell}] \right]= (C_i \cdot C_j) [\frakZ_{j,\ell+d}]\ .
		\end{align}
		
		\item if $i=j$ and $\widehat{S}_{C_i} \simeq (\widehat{T^\ast C_i})_{C_i}$:
		\begin{align}\label{eq:commutations-relations-rank-0-rank-1-i}
			\left[ [\frakY_{i,(d)}], [\frakZ_{i,\ell}] \right]= (C_i \cdot C_i) [\frakZ_{i,\ell+d}]\ .
		\end{align}
	\end{itemize}
\end{proposition}

\begin{conjecture}
	Equality~\eqref{eq:commutations-relations-rank-0-rank-1-i} holds also without the assumption on $\widehat{S}_{C_i}$.
\end{conjecture}

\subsubsection{Proof of Equality~\eqref{eq:commutations-relations-rank-0-rank-1-i-j}}
	
The proof follows from the following lemma.
\begin{lemma}\label{lem:commutations-relations}
	For $i \neq j$ we have
	\begin{align}\label{eq:commutations-relations-1}
		[\frakY_{i,(d)}]\cdot [\frakZ_{j,\ell}] = &[\frakZ_{j,\ell+d}^{\boldsymbol{\lambda}}] + \vert C_i\cap C_j\vert  [\frakZ_{j,\ell+d}]\ , \\[2pt] \label{eq:commutations-relations-2} 
		[\frakZ_{j,\ell}]\cdot[\frakY_{i,(d)}] = &[\frakZ_{j,\ell+d}^{\boldsymbol{\lambda}}]\ ,
	\end{align}
	where $\boldsymbol{\lambda}=(\lambda^{(h)})$ with $\lambda^{(i)}=(d)$ and $\lambda^{(h)}=\emptyset$ for $h \neq i$. 
\end{lemma}

\begin{proof}
	Consider extension $\calF$ of a one-dimensional sheaf $\calG$ scheme-theoretically supported on $C_j$ by a zero-dimensional sheaf $\calT$ set-theoretically supported at one point of $C_i$ of the form 
	\begin{align}
		0 \longrightarrow \calT \longrightarrow \calF \longrightarrow  \calG \longrightarrow 0 \ .
	\end{align}
	If $\calT$ is set-theoretically supported on the intersection $C_i\cap C_j$, then $\calF$ is set-theoretically supported on $C_j$ and scheme-theoretically supported on $C_j\smallsetminus C_i\cap C_j$; otherwise, $\calF\simeq \calG \oplus \calT$. It follows from the description of irreducible components given in Proposition~\ref{prop:irreducible-components-in-rank-1} that $[\frakY_{i,(d)}] \cdot [\frakZ_{j,\ell}] \in \Q [\frakZ_{j,\ell+d}^{\boldsymbol{\lambda}}]$. Base changing to the open locus of $\frakZ_{j,\ell+d}^{\boldsymbol{\lambda}}$ consisting of sheaves $\calG \oplus \calT$ with $\mathsf{Supp}(\calT)\not\subset C_j$ one easily obtains Equation~\eqref{eq:commutations-relations-2}.
	
	A similar argument involving short exact sequences
	\begin{align}
		0 \longrightarrow \calG \longrightarrow \calF \longrightarrow  \calT \longrightarrow 0
	\end{align}
	shows that $[\frakY_{i,(d)}]\cdot [\frakZ_{j,\ell}] \in \Q [\frakZ_{j,\ell+d}^{\boldsymbol{\lambda}}] \oplus \Q [\frakZ_{j,\ell+d}]$, and that the coefficient of $[\frakZ_{j,\ell+d}^{\boldsymbol{\lambda}}]$ is equal to $1$. It remains to compute the coefficient of $[\frakZ_{j,\ell+d}]$. Let us set $\alpha\coloneqq -\alpha_j+(\ell+1)\delta$ and $\gamma\coloneqq\alpha+d\delta$. We consider the following diagram
	\begin{align}\label{eq:proof_one-dim-commutation-1}
		\begin{tikzcd}[ampersand replacement=\&, column sep=large]
			\dstackCohpsnil(\widehat{S}_C;d\delta) \times \dstackCohpsnil(\widehat{S}_C;\alpha) \& \calS_2\dstackCohpsnil(\widehat{S}_C;d\delta,\alpha) \ar[swap]{l}{q} \ar{r}{p} \& \dstackCohpsnil(\widehat{S}_C;\gamma)\\
			\frakY_{i,(d)} \times \frakZ_{j,\ell}\ar{u}{\imath} \& \bfW \ar{u}{h} \ar[swap]{l}{q_{i, j}} \ar{r}{p_{i, j}}  \& \dstackCohpsnil(\widehat{S}_C;\gamma) \ar[equal]{u}{}\\
			\frakY_{i,(d)}^\circ \times \frakZ_{j,\ell} \ar{u}{\jmath_{i,j}}\& \bfX \ar[swap]{l}{q_{i, j}^\circ} \ar{u}{\jmath_{i,j}^\circ} \ar{r}{p_{i, j}^\circ}  \& \dstackPic_{\ell+d}(C_j) \ar{u}{\jmath}
		\end{tikzcd}
	\end{align}
	where the upper left and bottom squares are cartesian and where $\jmath_{i,j},\jmath_{i,j}^\circ,\jmath$ are open immersions. Here, $p\coloneqq \partial_1$ and $q\coloneqq \partial_0\times \partial_2$, following the notation of Formula~\eqref{eq:partial}.
	
	Note that a zero-dimensional quotient $\calT$ of a line bundle on $C_j$ with support consisting of a point is necessarily \textit{regular}, i.e., $\calT \simeq \scrO_x^{(d)}$ for some $x\in C_j$, hence the factorization of the map $\bfX \to \frakY_{i,(d)} \times \frakZ_{j,\ell}$ through the open embedding $\frakY_{i,(d)}^\circ \times \frakZ_{j,\ell} \to \frakY_{i,(d)} \times \frakZ_{j,\ell}$. 
	
	By base change, we have
	\begin{align}
		\jmath^\ast \left([\frakY_{i,(d)}]\star [\frakZ_{j,\ell}]\right)&=\jmath^\ast p_\ast q^!\imath_\ast ([\frakY_{i,(d)}]\boxtimes [\frakZ_{j,\ell}])=\jmath^\ast (p_{i, j})_\ast (q_{i, j})^!([\frakY_{i,(d)}]\boxtimes [\frakZ_{j,\ell}])\\
		&=(p_{i, j}^\circ)_\ast (q_{i, j}^\circ)^!\jmath_{i,j}^\ast ([\frakY_{i,(d)}]\boxtimes [\frakZ_{j,\ell}])=(p_{i, j}^\circ)_\ast ([\bfX])
	\end{align}
	Here we used the fact that $\frakY_{i,(d)}^\circ \times \frakZ_{j,\ell}$ is smooth and classical, hence the virtual fundamental class coincides with the usual fundamental class.
	
	The map $q$ being of dimension zero, so is the map $q_{i, j}$, and thus
	\begin{align}
		\dim(\bfX)=\dim(\bfW)=\dim(\frakY_{i,(d)})+\dim(\frakZ_{j,\ell})=g_j-1\ .
	\end{align}
	On the other hand, since any short exact sequence
	\begin{align}
		0 \longrightarrow \calF \longrightarrow \calL \longrightarrow  \calT \longrightarrow 0
	\end{align}
	with $\calL$ a line bundle on $C_j$ and $\calT$ a zero-dimensional sheaf on $C_i$ of length $d$ supported at a point is of the form
	\begin{align}
		0 \longrightarrow \calL(-dx) \longrightarrow \calL \longrightarrow  \scrO_x^{(d)} \longrightarrow 0
	\end{align}
	for some $x \in C_i \cap C_j$, the classical truncation of $p_{i, j}^\circ$ factors as
	\begin{align}
		\begin{tikzcd}[ampersand replacement=\&, column sep=large]
			\trunc{\bfX}\ar[equal]{d} \ar{r}{\trunc{(p_{i, j}}^\circ)}\& \stackPic_{\ell+d}(C_j)\\
			\displaystyle\bigsqcup_{x \in C_i \cap C_j} \trunc{\bfX}_x \ar{r}{\sim} \& \stackPic_{\ell+d}(C_j) \times (C_i \cap C_j) \ar{u}{\mathsf{pr}_1}
		\end{tikzcd}\ .
	\end{align} 
	Here, $\trunc{\bfX}_x$ is defined via the pullback of $\partial_2\colon \bfX\to \frakY_{i,(d)}^\circ$ with respect to the inclusion $x\to C_i\cap C_j$. It follows that $\dim(\trunc{\bfX})=\dim(\stackPic_{\ell+d}(C_j))=g_j-1=\dim(\bfX)$. This in turn implies that $\bfX$ is classical by e.g. \cite[Proposition~2.14]{Porta_Yu_Non-archimedian_quantum}\footnote{Note that the result is written for derived analytic stacks, but it holds in our setting verbatim.}. Denoting by $\iota\colon \stackPic_{\ell+d}(C_j) \to \dstackPic_{\ell+d}(C_j)$ the closed embedding, we have
	\begin{align}
		(p_{i, j}^\circ)_\ast ([\bfX])=\iota_\ast \trunc{(p_{i, j}}^\circ)_\ast ([\trunc{\bfX}])=\vert C_i \cap C_j\vert  \iota_\ast ([\stackPic_{\ell+d}(C_j)])=\vert C_i \cap C_j\vert  [\frakZ_{j,\ell+d}]\ ,
	\end{align}
	as wanted.
\end{proof}

\subsubsection{Proof of Equality~\eqref{eq:commutations-relations-rank-0-rank-1-i}}

Before proving Equality~\eqref{eq:commutations-relations-rank-0-rank-1-i}, we need several preliminary results.

Set $C=C_i$. Let $\stackCoh(C)$ be the classical moduli stack of coherent sheaves on $C$. This is a smooth geometric classical stack. The connected component classifying sheaves of rank $r$ and degre $d$ will be denoted $\stackCoh(C;r,d)$. For any $\alpha,\beta$ in the \textit{numerical Grothendieck group} $K_0^{\mathsf{num}}(C)$ there is an induction diagram
\begin{align}\label{eq:conv_diagram_dim_one}
	\begin{tikzcd}[ampersand replacement=\&]
		\stackCoh(C;\alpha) \times \stackCoh(C;\beta) \ar[leftarrow]{r}{q_{\alpha,\beta}} \& \stackCoh^{\mathsf{ext}}(C;\alpha,\beta) \ar{r}{p_{\alpha,\beta}} \& \stackCoh(C;\alpha+\beta)
	\end{tikzcd}
\end{align}
where $\stackCoh^{\mathsf{ext}}(C;\alpha,\beta)$ is the stack classifying short exact sequences
\begin{align}
	0 \longrightarrow \calF \longrightarrow \calG \longrightarrow  \calH \longrightarrow 0
\end{align}
of coherent sheaves on $C$, where $\calF,\calH$ are of classes $\beta$ ad $\alpha$ respectively. Here, $q_{\alpha, \beta}=\partial_0\times \partial_2$ and $p_{\alpha, \beta}=\partial_1$ in the notation of Formula~\eqref{eq:partial}.

Let $\catDb_c(\stackCoh(C))$ be the derived category of constructible complexes of $\Q$-vector spaces on $\stackCoh(C)$. Laumon considered an \textit{Eisenstein series} monoidal functor (cf.\ \cite[Lecture~5]{Schiffmann-lectures-canonical-bases} and references therein)
\begin{align}
	\mathsf{Eis}=\bigoplus_{\alpha,\beta} \mathsf{Eis}_{\alpha,\beta}\colon \catDb_c(\stackCoh(C)) \boxtimes \catDb_c(\stackCoh(C)) \longrightarrow \catDb_c(\stackCoh(C))\ ,  
\end{align} 
where
\begin{align}
	 \mathsf{Eis}_{\alpha,\beta}\coloneqq (p_{\alpha,\beta})_!q_{\alpha,\beta}^\ast [-\langle \alpha,\beta\rangle] \ .
\end{align}
He also defined a category of semisimple complexes $\frakQ_C \subset \catDb_c(\stackCoh(C))$ containing the constant sheaves $\Q_{\stackCoh(C;\alpha)}$ for any $\alpha \in K_0^{\mathsf{num}}(C)$, stable under shifts and $\mathsf{Eis}$. The graded Grothendieck group $K_0(\frakQ_C)$ is the quotient of the split Grothendieck group of $\frakQ_C$ by the relation $\left[\PP[1]\right]=\upsilon[\PP]$; it is naturally a $\Q[\upsilon,\upsilon^{-1}]$-algebra. 

As explained in \cite[Remark~2.41]{Porta_Sala_Hall}, there is a canonical equivalence 
\begin{align}
T^\ast[0]\dstackCoh(C) \simeq \dstackCoh_{\mathsf{Dol}}(C)\simeq \dstackCohps(T^\ast C)\ ,
\end{align}
where $\dstackCoh_{\mathsf{Dol}}(C)$ is the\textit{ Dolbeault moduli stack} of $C$ (see \cite[\S2.4.3]{Porta_Sala_Hall}) and the second equivalence follows from the spectral correspondence for Higgs sheaves \cite[Lemma~6.8]{Simpson_Moduli_II}. It is known that the singular supports of all objects in $\frakQ_C$ lie in the substack $\stackCoh_{C,\mathsf{ps}}(T^\ast C) \simeq \stackCohpsnil(\widehat{T^\ast C}_C)$. We thus get a map
\begin{align}
	\mathsf{CC} \colon K_0(\frakQ_C) \longrightarrow (\coha_{T^\ast C,C})_0\ .
\end{align}
As explained in \cite[\S1.3]{Hennecart_Geometric} this is an algebra morphism. Note that, since we consider the untwisted COHA $\coha_{T^\ast C,C}$, there is no need to twist the multiplication in $K_0(\frakQ_C)$. The category $\frakQ_C$ is well-studied (see e.g. \cite[Lecture~5]{Schiffmann-lectures-canonical-bases} and \cite{Schiffmann_Duke_ellipticII} for genus one curves). It carries a $\N$-grading 
\begin{align}
	\frakQ_C=\bigsqcup_r \frakQ_C^r
\end{align}
corresponding to the rank of the coherent sheaves on $C$. One has the following result (cf.\ \cite[\S2.4]{Schiffmann-lectures-canonical-bases}).
\begin{proposition}\label{prop:algebra_Hecke_operators}
	There is an algebra isomorphism 
	\begin{align}\label{eq:iso-Psi}
		\begin{tikzcd}[ampersand replacement=\&]
			\Psi\colon K_0(\frakQ_C^0) \ar{r}{\sim} \& \Q[x_1,x_2, \ldots]^{\mathfrak{S}_\infty} \otimes \Q[\upsilon,\upsilon^{-1}]
		\end{tikzcd}
	\end{align}
	which sends $\Q_{\stackCoh(C;0,d)}$ to the elementary symmetric function $e_d$ for any $d \geq 1$.
\end{proposition} 
The monoidal subcategory $\frakQ_C^0$ corresponding to rank $0$ sheaves plays a special role as it yields Hecke operators. More precisely, let 
\begin{align}
	\stackBun(C)=\bigsqcup_\alpha \stackBun(C;\alpha)
\end{align}
be the open substack of $\stackCoh(C)$ parametrizing vector bundles on $C$. Let $\frakQ_C^{\mathsf{bun}}$ be the image of $\frakQ_C$ under the restriction functor $\catDb_c(\stackCoh(C)) \to \catDb_c(\stackBun(C))$. One has the following (cf.\ \cite[\S5.5]{Schiffmann-lectures-canonical-bases}).

\begin{proposition}
	The Eisenstein functor $\mathsf{Eis}$ induces a monoidal action of $\frakQ_C^0$ on $\frakQ^{\mathsf{bun}}_C$.
\end{proposition}
For $u \in K_0(\frakQ^0_C)$ and $c \in K_0(\frakQ^{\mathsf{bun}}_C)$ we denote by $u \otimes c \mapsto u \bullet c$ the corresponding Hecke action. 

We will deduce the desired relation \eqref{eq:commutations-relations-rank-0-rank-1-i} from a Hecke relation by applying the morphism $\mathsf{CC}$. The precise Hecke relation which we will use is the following. 

\begin{proposition}\label{prop:Hecke-relation}
	For any $d\in \N$, $d\geq 1$, set $T_d\coloneqq \Psi^{-1}(p_d)$. Then, we have 
	\begin{align}\label{eq:commutator-0-1}
		\left[ [{T}_d], [\Q_{\stackCoh(C;1,\ell)}] \right]=-((-\upsilon)^d-2g + (-\upsilon)^{-d})[\Q_{\stackCoh(C;1,\ell+d)}]\ .
	\end{align}
\end{proposition}
This result is without doubt well-known to experts but we provide a proof as we could not locate one in the literature. Before giving it, we need some preliminary results.

Consider the \textit{constant term functor}
\begin{align}
	\mathsf{CT}=\prod_{\alpha,\beta} \mathsf{CT}_{\alpha,\beta}\colon \frakQ_C \longrightarrow \frakQ_C \widehat{\boxtimes} \frakQ_C  
\end{align}
given by
\begin{align}
	\mathsf{CT}_{\alpha,\beta}\coloneqq (q_{\alpha,\beta})_!p_{\alpha,\beta}^\ast [-\langle \alpha,\beta\rangle] \ .
\end{align}
It equips $K_0(\frakQ_C)$ and $K_0(\frakQ^0_C)$ with the structure of bialgebras (understood in a twisted and topological sense in the former case). Moreover, the map $\Psi$ introduced in Formula~\eqref{eq:iso-Psi} is a morphism of bialgebras. 

We introduce a truncation $\mathsf{CT}^{\mathsf{bun}, 0}\colon \frakQ_C \to \frakQ^{\mathsf{bun}}_C \widehat{\boxtimes} \frakQ^0_C$ by projecting to $\frakQ^{\mathsf{bun}}_C$ in the first component and to $\frakQ^{0}_C$ in the second.
\begin{lemma}\label{lem:injectivity} 
	The functor $\mathsf{CT}^{\mathsf{bun}, 0}$ is essentially injective on objects, i.e., $K_0(\mathsf{CT}^{\mathsf{bun}, 0})$ is injective.
\end{lemma}

\begin{proof}
	Let 
	\begin{align}
		\stackCoh(C;\alpha)=\bigsqcup_{d \geq 0} \stackCoh(C;\alpha)^{\mathsf{tor}=d}
	\end{align}
	be the stratification of $\stackCoh(C;\alpha)$ according to the length of the canonical zero-dimensional subsheaf. Then 
	\begin{align}
		\stackCoh(C;\alpha)^{\mathsf{tor} \geq d} \coloneqq \bigsqcup_{\ell \geq d} \stackCoh(C;\alpha)^{\mathsf{tor}=\ell}\quad\text{and}\quad \stackCoh(C;\alpha)^{\mathsf{tor} \leq d} \coloneqq \bigsqcup_{\ell \leq d} \stackCoh(C;\alpha)^{\mathsf{tor}=\ell}
	\end{align}
	are respectively closed and open in $\stackCoh(C;\alpha)$. 
	
	Now let us fix a simple complex $\PP \in \catDb_c(\stackCoh(C;\alpha))$ and let $d$ be the greatest integer such that $\mathsf{Supp}(\PP) \subset \stackCoh(C;\alpha)^{\mathsf{tor} \geq d}$. Put $\delta\coloneqq (0,1)$ and $\beta\coloneqq \alpha-d\delta$. There exists a diagram pulled back from the diagram~\eqref{eq:conv_diagram_dim_one}:
	\begin{align}
		\begin{tikzcd}[ampersand replacement=\&, column sep=large]
			\stackCoh(C;\beta) \times \stackCoh(C;d\delta) \& \stackCoh(C;\beta,d\delta)^{\mathsf{ext},\mathsf{tor}\geq d} \ar[swap]{l}{q} \ar{r}{p} \& \stackCoh(C;\alpha)^{\mathsf{tor} \geq d}\\
			\stackBun(C;\beta) \times \stackCoh(C;d\delta) \ar{u}{j} \& \stackCoh(C;\beta,d\delta)^{\mathsf{ext},\mathsf{tor}= d} \ar[swap]{l}{q^\circ} \ar{r}{p^\circ} \ar{u}{j'} \& \stackCoh(C;\alpha)^{\mathsf{tor}= d} \ar{u}{j''}
		\end{tikzcd}\ .
	\end{align} 
	Since a coherent sheaf in $\stackCoh(C;\alpha)^{\mathsf{tor}=d}$ has a unique subsheaf of length $d$, the map $p^\circ$ is an isomorphism. Since $\Ext^{1}_C(\calF,\calT)=0$ for any vector bundle $\calF$ and any torsion sheaf $\calT$, the map $q^\circ$ is a unipotent gerbe. It follows that 
	\begin{align}
		q^\circ_!(p^\circ)^\ast  \colon  \catDb_c(\stackCoh(C;\alpha))^{\mathsf{tor}=d} \longrightarrow \catDb_c(\stackBun(C;\beta)) \boxtimes \catDb_c(\stackCoh(C;d\delta))
	\end{align}
	is an equivalence. By construction $(j'')^\ast (\PP)$ is a simple complex hence $j^\ast (\mathsf{CT}^{\mathsf{bun}, 0}(\PP))=q^\circ_!(p^\circ j'')^\ast (\PP)$ is simple as well, and fully determines $\PP$. The Lemma easily follows.
\end{proof}

To unburden the notation, let us set 
\begin{align}
	\mathbf{1}_\alpha=[\Q_{\stackCoh(C;\alpha)}]\quad \text{and}\quad \mathbf{1}^{\mathsf{bun}}_\alpha=[\Q_{\stackBun(C;\alpha)}]
\end{align}
for $\alpha\in K_0^{\mathsf{num}}(C)$, and
\begin{align}
	\Delta=K_0(\mathsf{CT})\quad \text{and}\quad  \Delta'=K_0(\mathsf{CT}^{\mathsf{bun}, 0})\ .
\end{align}
Then from the relation
\begin{align}
	\mathsf{CT}_{\alpha,\beta}(\Q_{\stackCoh(C;\alpha+\beta)})=\Q_{\stackCoh(C;\alpha)} \boxtimes \Q_{\stackCoh(C;\beta)}[-\langle \beta,\alpha \rangle]
\end{align} 
(see \cite[Lemma 1.8]{Schiffmann-lectures-canonical-bases}) we get
\begin{align}
	\Delta'(\mathbf{1}_{\alpha})&=\sum_{n \geq 0} \upsilon^{-\langle \alpha,\alpha-n\delta\rangle} \mathbf{1}^{\mathsf{bun}}_{\alpha-n\delta} \otimes \mathbf{1}_{n\delta}\ ,\\[2pt]
	\Delta(T_d)&=T_d \otimes 1 + 1 \otimes T_d\ ,
\end{align}
and therefore
\begin{align}
	\Delta'([T_d, \mathbf{1}_\alpha])=\sum_{n \geq 0}\upsilon^{-\langle \alpha,\alpha-n\delta\rangle} T_d \bullet \mathbf{1}^{\mathsf{bun}}_{\alpha-n\delta} \otimes \mathbf{1}_{n\delta}\ .
\end{align}

\begin{proof}[Proof of Proposition~\ref{prop:Hecke-relation}]
	Thanks to Lemma~\ref{lem:injectivity}, it is enough to show that 
	\begin{align}
		T_d \bullet \mathbf{1}^{\mathsf{bun}}_{(1,a)}=-((-\upsilon)^d-2g + (-\upsilon)^{-d})\mathbf{1}^{\mathsf{bun}}_{(1,a+d)}
	\end{align}
	for any $a\in \Z$. To prove this, we first compute $\mathbf{1}_{d\delta} \bullet \mathbf{1}^{\mathsf{bun}}_{(1,a)}$. Set $\alpha\coloneqq (1,a+d)$ and $\beta\coloneqq (1,a)$. There is an isomorphism
	\begin{align}
		\begin{tikzcd}[ampersand replacement=\&, column sep=large]
			\stackBun(C;\alpha) \times \Sym^d(C) \ar{r}{\sim} \ar{dr}[swap]{\mathsf{pr}_1}\& \stackCoh^{\mathsf{ext}}(C;\beta,d\delta) \underset{\stackCoh(C;\alpha)}{\times}\stackBun(C;\alpha) \ar{d}{p}\\ 
			\& \stackBun(C;\alpha)
		\end{tikzcd}
	\end{align}
	provided by the map $(\calF, D) \mapsto (\calF(-D) \subset \calF)$. It follows that
	\begin{align}
		\mathbf{1}_{d\delta} \bullet \mathbf{1}^{\mathsf{bun}}_{(1,a)}=\upsilon^d P_c(\Sym^d(C),-\upsilon^{-1})\mathbf{1}^{\mathsf{bun}}_{(1,a+d)}\ ,
	\end{align}
	where $P_c(-,t)$ stands for the Poincaré polynomial with compact support. The relation between the elementary and power sum functions may be written in terms of generating series as
	\begin{align}
		P(-z)=-z\frac{E'(z)}{E(z)}\ , 
	\end{align}
	where
	\begin{align}
		E(z)\coloneqq 1+\sum_{n\geq 1} e_n z^n\quad \text{and} \quad P(z)\coloneqq \sum_{n\geq 1} p_nz^n\ .
	\end{align}
	From this and the well-known relation
	\begin{align}
		\sum_{d \geq 0} P_c(\Sym^dC,t)z^d=\frac{(1-zt)^{2g}}{(1-z)(1-t^2z)}\ ,
	\end{align}
	where $g$ is the genus of $C$ we deduce by a direct computation that
	\begin{align}
		T_d \bullet \mathbf{1}^{\mathsf{bun}}_{(1,a)}=-((-\upsilon)^d-2g + (-\upsilon)^{-d})\mathbf{1}^{\mathsf{bun}}_{(1,a+d)}
	\end{align}
	as wanted. This finishes the proof of the Proposition.
\end{proof}
	
\begin{proof}[Proof of Equality~\eqref{eq:commutations-relations-rank-0-rank-1-i-j}]
	We now apply the $\mathsf{CC}$ map to the Equality~\ref{eq:commutator-0-1}. By Propositions~\ref{prop:algebra_Hecke_operators} and \ref{prop:zero-dim_Higgs_sheaves_classification}--\eqref{item:zero-dim_Higgs_sheaves_classification-2}, we have
	\begin{align}
		\mathsf{CC}([\Q_{\stackCoh(C;d\delta)}])=[\frakY_{(1^d)}]\quad\text{and} \quad \mathsf{CC}([T_d])=[\frakY_{(d)}]\ .
	\end{align}
	Together with $\mathsf{CC}(\Q_{\stackCoh(C;(1,\ell))})=[\frakZ_\ell]$ and the fact that $\mathsf{CC}(\upsilon^a c)=(-1)^a\mathsf{CC}(c)$ we finally obtain 
	\begin{align}
		\left[ [\frakY_{(d)}],[\frakZ_{\ell}]\right]=(2g-2)[\frakZ_{\ell+d}]
	\end{align}
	as wanted.
\end{proof}	
	
\subsection{Proof of Theorem~\ref{thm:class_irred_component}}\label{sec:proof-thm-class-irred}

We now come back to the setting of the pair $(\rsv,C)$, with $\rsv$ being the resolution of the Kleinian singularity of type $\qvfin$. We start by proving both statements of Theorem~\ref{thm:class_irred_component} in the case $\qvfin$ is of type $A_1$. We will then use the functoriality to extend it to the other types.

\subsubsection{Type $A_1$}

We assume that $\rsv=T^\ast \PP^1$ with $C \subset \rsv$ being the zero section, and $\qv$ is the Kronecker quiver with vertex set $I=\{0,1\}$ and two arrows $x,y\colon 1 \mapsto 0$. For the stability condition we choose $\Lthetafin=\Lrho_{\sf{f}}=\frac{1}{2}\Lalpha_1$, so that $\Ltheta=\frac{1}{2}(\Lalpha_1-\Lalpha_0)$.

We begin with Statement~\eqref{item:class_irred_component-1}. We will drop the index $j$ from the notation. Instead of $[\frakY_{(n)}]$, we will compute $[\frakY_{(1^n)}]$, and use Proposition~\ref{prop:zero-dim_Higgs_sheaves_classification}--\eqref{item:zero-dim_Higgs_sheaves_classification-2}. For this, we shall use the tilting equivalence $\tau$ of Theorem~\ref{thm:nilpotent}. 

Observe that the equivalence $\tau$ sends a zero-dimensional sheaf schematically supported on $C=\PP^1$ of length $\ell$ to a regular representation of $\qv$ (viewed as a $\Pi_\qv$-representation) of dimension $\ell\delta$. The regular representations of $\qv$ of dimension $n\delta$ form a dense open substack of $\stackRep(\qv; d\delta)$ and the latter is irreducible. In particular, the stack equivalence $\stackCohpsnil(\widehat{\rsv}_C)^{>0} \simeq \Lambda^{\leq 0}_\qv$ identifies $\frakY_{(1^n)}$ with $\stackRep(\qv,n\delta)^{\leq 0}$. It follows that
\begin{align}
	\Theta([\frakY_{(1^n)}])=[\stackRep(\qv,n\delta)^{\leq 0}]\ .
\end{align}
Let us now translate this in terms of enveloping algebras. We have 
\begin{align}
	\frakgfin=\fraksl_2\quad \text{and} \quad \fraknellzero^+= \Q e[s^{\pm 1}] \oplus \Q hs^{-1}[s^{\pm 1}]\ .
\end{align}
Here we use the traditional $e,f$ in place of $x_1^{\pm 1}$; note that $x_0^-=es^{-1}$. A simple computation using the twisted COHA multiplication shows that
\begin{align}	
	\Phi([\stackRep(\qv,n\delta)^{\leq 0}])=(-1)^n(x_1^-)^{(n)}(x_0^-)^{(n)}=(-1)^nf^{(n)}(es^{-1})^{(n)} \in \sfU(L\fraksl_2) \subset \Y^-_\qv\ ,
\end{align}
where we have used the notation $z^{(a)}\coloneqq z^a/a!$. We deduce that
\begin{align}
	\Theta_{\rsv,C}([\frakY_{(1^n)}])=(-1)^nf^{(n)}(es^{-1})^{(n)} \in \widehat{\sfU}(\fraknell^+)\ .
\end{align}
We will use the following Lemma.
\begin{lemma}We have
	\begin{align}\label{eq:Kr-computation-h}
		\sum_{\ell \geq 0}  (-1)^{\ell} f^{(\ell)}(es^{-1})^{(\ell)}u^{-\ell}=\exp\left( \sum_{k\geq 1} hs^{-k} \frac{u^{-k}}{k}\right)\ . 
	\end{align}
\end{lemma}

\begin{proof}
	Set $Y(u) \coloneqq \sum_{\ell \geq 0}  (-1)^{\ell}f^{(\ell)}(es^{-1})^{(\ell)}u^{-\ell}$.
	
	We shall use the coproduct map. Set $\J\coloneqq \J_{\Ltheta, \kappa_0}$, where $\kappa_0\coloneqq (-\infty, 0]$ and $\J_{\Ltheta, \kappa_0}$ is defined in Formula~\eqref{eq:quotient-Yangian-kappa}. Since $\J$ is a coideal of $\sfU^-(L\fraksl_2)$, the comultiplication $\Delta$ descends to $\sfU^-(L\fraksl_2)/\J$. For weight reasons, one sees that $\Delta(Y(u))=Y(u) \otimes Y(u)$, from which we deduce that 
	\begin{align}
		Y(u)=\exp\left(\sum_{k \geq 1} c_k hs^{-k}u^{-k}\right)
	\end{align}
	for some coefficients $c_k \in \Q$. Finally, to see that $c_k=1/k$, we argue by induction using the relation
	\begin{align}
		\left(f^\ell(es^{-1})^{\ell}\right)_{[hs^{-\ell}]}&=-\left( f^{\ell-1} \Big(\ell(es^{-1})^{\ell-1} hs^{-1} + \ell(\ell-1)(es^{-1})^{\ell-2} es^{-2}\Big)\right)_{[hs^{-\ell}]}\\
		&=-\ell(\ell-1)\left(f^{\ell-1}\cdot(es^{-1})^{\ell-2}\cdot es^{-2}\right)_{[hs^{-\ell}]}\ .
	\end{align} 
	We leave the details of the calculation to the reader.
\end{proof}

Note that under the isomorphism of Proposition~\ref{prop:zero-dim_Higgs_sheaves_classification}--\eqref{item:zero-dim_Higgs_sheaves_classification-2}, $[\frakY_{(1^n)}]$ and $[\frakY_{(n)}]$ are respectively mapped to the symmetric functions $e_n$ and $p_n$. Comparing Equality~\eqref{eq:Kr-computation-h} with the well-known relation in the ring of symmetric functions
\begin{align}
	\sum_{\ell\geq 0}e_\ell u^\ell=\exp\left(\sum_{k \geq 1}(-1)^{k-1}p_k \frac{u^k}{k}\right)
\end{align}
we obtain the desired equality $\Theta_{\rsv,C}([\frakY_{(n)}])=(-1)^{n-1}hs^{-n}$. This proves Statement~\eqref{item:class_irred_component-1} for type $A_1$.

\medskip

The proof of Statement~\eqref{item:class_irred_component-2} proceeds in an analogous fashion. We have $\tau(\scrO_C)=S_0$. We claim that for any $m \geq 0$, 
\begin{align}
	\R S_{m\Lthetafin}(S_0)=\tau(\scrO_C(m))=P^{(m)}\ ,
\end{align}
where $P^{(m)}$ is the unique indecomposable preprojective representation of $\qv$ of dimension $\alpha_0+m\delta$, viewed as a nilpotent $\Pi_\qv$-module. Indeed, $\tau(\scrO_C(m))$ is the unique indecomposable and semistable $\Pi_\qv$-module of dimension $\alpha_0+m\delta$ and one easily checks that $P^{(m)}$ is semistable. Moreover, the Zariski closure of the point of $\Lambda_{\alpha_0+m\delta}$ corresponding to $P^{(m)}$ is $\stackRep(\qv,\alpha_0+m\delta)$, which is an irreducible component of $\Lambda_{\alpha_0+m\delta}$. We deduce that
\begin{align}
	\Theta([\frakZ_{4k+1}])=\left( \R S_{\textrm{-}4k\Lthetafin}([\stackRep(\qv,\alpha_0+4k\delta)])\right)_k \in \lim \coha_{\Ltheta,(\textrm{-}2k)}^A\ .
\end{align}

A direct computation shows that
\begin{align}
	\Phi([\stackRep(\qv,\alpha_0+4k\delta)]) =(x_1^-)^{(4k)} (x_0^-)^{(4k+1)}=f^{(4k)} (es^{-1})^{(4k+1)}\in \sfU^-(L\fraksl_2) \subset \Y_\qv^-\ .
\end{align}
Putting it all together it follows that 
\begin{align}\label{eq:theta'_irreduc_comp_type_A1}
	\Theta_{\rsv,C}([\frakZ_{4k+1}])&=\left( \overline{T}_{\textrm{-}4k\Lthetafin}((f)^{(4k)} (es^{-1})^{(4k+1)})\right)_k \\
	&=\left((fs^{-4k})^{(4k)} (es^{4k-1})^{(4k+1)}\right)_{2k}\in \widehat{\sfU}(\fraknellzero^+)\ ,
\end{align}
see Remark~\ref{rem:third-limit-realization-coha}. In order to make the right-hand-side of \eqref{eq:theta'_irreduc_comp_type_A1} explicit, we prove the following identity.
\begin{lemma}\label{prop:funky_identity_sl_2}
	In $\Y_{\Ltheta,(0)}$ we have
	\begin{align}\label{eq:funky_identity_in_Lsl_2}
		\sum_{\ell \geq 0} (-1)^\ell f^{(\ell)} (es^{-1})^{(\ell+1)} u^{-\ell-1}=\left(\sum_{k \geq 1} es^{-k}u^{-k}\right) \exp\left(\sum_{k \geq 1} hs^{-k} \frac{u^{-k}}{k}\right)\ .
	\end{align}
\end{lemma}

\begin{proof}
	First observe that the graded piece of $\Y_{\Ltheta,(0)}$ of vertical degree $0$ is equal to $\sfU^-(L\fraksl_2)/\J$, where $\J$ is the left ideal generated by $fs^a$ for $a \leq 0$. By the PBW theorem, every element of $\sfU^-(L\fraksl_2)/\J$ of degree $\alpha_1 + k\delta$ (resp.\ $k\delta$) may uniquely be written as $es^{-n} P(hs^{-1}, hs^{-2}, \ldots)$ for some integer $n \geq 1$ and some polynomial $P$ (resp.\ as $P((hs^{-1}, hs^{-2}, \ldots)$). Let us denote by $X(u)$ the left-hand-side of Equality~\eqref{eq:funky_identity_in_Lsl_2}. By the above, it is enough to show that 
	\begin{align}
		\Delta'(X(u))=\left(\sum_{k\geq 1} es^{-k}u^{-k}\right) \otimes \exp\left(\sum_{k \geq 1} hs^{-k} \frac{u^{-k}}{k}\right) + \ldots
	\end{align}
	where $\Delta'$ is the graded component of $\Delta$ of bidegree in $(\alpha_1 + \Z \delta,\Z\delta)$ and where $\ldots$ stands for terms whose first component is of the form $es^{-l}P(hs^{-1}, hs^{-2},\ldots)$ with $P \neq 1$. From the Leibniz formula we obtain
	\begin{align}
		\Delta'(f^{(\ell)} (es^{-1})^{(\ell+1)})=\sum_{k=0}^\ell  f^{(k)} (es^{-1})^{(k+1)} \otimes f^{(\ell-k)} (es^{-1})^{(\ell-k)}\ ,
	\end{align}
	from which we deduce
	\begin{align}
		\Delta'(X(u))=X(u) \otimes \left(\sum_{k \geq 0} (-1)^kf^{(k)} (es^{-1})^{(k)}\right)\ .
	\end{align}
	Denote as usual by $x_{[a]}$ the coefficient of $a$ in $x$ (with respect to a fixed natural basis). It is easy to see from the commutation relations that the coefficient
	\begin{align}
		\left(f^\ell \cdot es^{-a_1}\cdots es^{-a_{\ell+1}}\right)_{[es^{-\sum a_i}]}
	\end{align}
	is independent of $(a_1, \ldots, a_{\ell+1})$. A simple induction based on the identity
	\begin{align}
		\left(f^\ell(es^{-1})^{\ell+1}\right)_{[es^{-\ell-1}]}&=-\left( f^{\ell-1} \Big((\ell+1)(es^{-1})^{\ell} hs^{-1} + \ell(\ell+1)(es^{-1})^{\ell} es^{-2}\Big)\right)_{[es^{-\ell-1}]}\\
		&=-\ell(\ell+1)\left(f^{\ell-1}\cdot(es^{-1})^{\ell}\cdot es^{-2}\right)_{[es^{-\ell-1}]}
	\end{align} 
	thus shows that
	\begin{align}
		X(u)=\sum_{k\geq 1} es^{-k}u^{-k} + \ldots
	\end{align}
	where $\ldots$ denotes terms of the form $es^{-l}P(hs^{-1}, hs^{-2},\ldots)$ with $P \neq 1$. We conclude using Formula~\eqref{eq:Kr-computation-h}.
\end{proof}

As a direct corollary of Lemma~\ref{prop:funky_identity_sl_2}, applying the correct power of $\overline{T}_{4\Lthetafin}$ to the terms occuring in \eqref{eq:funky_identity_in_Lsl_2} we finally obtain
\begin{align}
	\Theta_{\rsv,C}([\frakZ_1])=\sum_{n \geq -1}es^{n}\cdot P_{n-1}(hs^{-1}, hs^{-2},\ldots)\ ,
\end{align}
where the polynomials $P_k(hs^{-1},hs^{-2},\ldots)$ are defined through the relation 
\begin{align}
	\exp\left(\sum_{k \geq 1} hs^{-k}\frac{u^{-k}}{k}\right)=\sum_{k \geq 1} P_k(hs^{-1},hs^{-2},\ldots) u^{-k}\ .
\end{align}
Observe that the restriction of $\overline{T}_{4\Lthetafin}$ to $\sfU(s^{-1}\frakhfin[s^{-1}])$ is the identity. By tensoring by line bundles and applying the braid group action, we can deduce the explicit formulas for all $\Theta_{\rsv,C}([\frakZ_m])$'s for any $m\in \Z$. This provides the desired formula~\eqref{eq:Theta'(Z)} for $\rsv=T^\ast \PP^1$ and proves the Theorem for type $A_1$.

\subsubsection{General case}

Let us next consider the case of an arbitrary affine type $\qv$. Again, we first establish Statement~\eqref{item:class_irred_component-1}. We will use Proposition~\ref{prop:commutations-in-rank-0-and-rank-1}, and characterize $[\frakY_{i,(n)}]$ through its commutations relations with $[\frakZ_{i,\ell}]$. 

Set $\sfU_0\coloneqq\Q[h_is^{-n}\;\vert\;i \in I_\sff, n \geq 1]$. For weight reasons, $\Theta_{\rsv,C}$ restricts to an isomorphism
\begin{align}
	\begin{tikzcd}[ampersand replacement=\&]
		\Q\left[[\frakY_{i,(n)}]\;\vert\; i \in I_\sff, n \geq 1\right]\ar{r}{\sim} \& \sfU_0
	\end{tikzcd}\ , 
\end{align}
and we have
\begin{align}
	[\frakZ_{i,\ell}]\in \bigoplus_{k\geq 0} \left(e_is^{-\ell+k} \cdot \sfU_0 \right)\ .
\end{align}
In particular, 
\begin{align}
	(\Theta_{\rsv,C})^{-1}(e_is^{-\ell}) =\sum_{k \geq 0}[\frakZ_{i,\ell+k}] P_k
\end{align}
for suitable elements $P_k \in \Q[[\frakY_{i,(n)}]\;\vert\; i \in I_\sff, n \geq 1]$. Using Proposition~\ref{prop:action_translation_braid_Lie} we see that $P_k$ is independent of $\ell$. By Proposition~\ref{prop:commutations-in-rank-0-and-rank-1} (taking into account the twist in the multiplication), we have for any $i,j \in I_\sff$
\begin{align}\label{eq:relation theta(e_i) and Yj}
	\left[[\frakY_{j,(n)}], (\Theta_{\rsv,C})^{-1}(e_{i}s^{-\ell})\right]&=\left[[\frakY_{j,(n)}], \sum_k [\frakZ_{i,\ell+k}]P_k\right]= \sum_k \left[[\frakY_{j,(n)}],[\frakZ_{i,\ell+k}]\right]P_k\\
	&= (C_i \cdot C_j)\sum_k [\frakZ_{i,\ell+k+n}] P_k=(-1)^{n-1} (-C_i \cdot C_j) (\Theta_{\rsv,C})^{-1}(e_is^{-\ell+n}) \ .
\end{align}

We now make use of the following elementary observation.
\begin{lemma}\label{lem:injectivity-ad}
	The map 
	\begin{align}
		\bigoplus_i \mathsf{ad}(e_{i})\colon \sfU_0 \longrightarrow \bigoplus_i \sfU(L\frakgfin)
	\end{align}
	is injective.
\end{lemma}

\begin{proof} 
	Let $Q \in \sfU_0$ and let us write $Q\coloneqq\sum_{k=0}^N Q_k$ with $Q_k \in \sfU_0$ being a sum of monomials in the $\Lomega_i s^{-n}$ of length $k$. Here $\Lomega_i \in \frakh_\sff$ is the $i$th fundamental coweight. Assume that $[Q,e_i]=0$ for all $i \in I_\sff$. We have
	\begin{align}\label{eq:commutation-lemma}
		\begin{split}
			[e_i, \Lomega_{j_1} s^{-n_1} \cdots \Lomega_{j_k} s^{-n_k}]&=-\sum_{\ell=1}^k\delta_{i,j_\ell} \Lomega_{j_1} s^{-n_1} \cdots e_{j_\ell}s^{-n_\ell} \cdots \Lomega_{j_k} s^{-n_k}\\
			&=-\sum_{\ell=1}^k\delta_{i,j_\ell} e_{j_\ell}s^{-n_\ell} \Lomega_{j_1} s^{-n_1} \cdots \widehat{\Lomega_{j_\ell}s^{-n_\ell}} \cdots \Lomega_{j_k} s^{-n_k}  + \text{l.o.t.}\ ,
		\end{split}
	\end{align}
	where $\text{l.o.t.}$ stands for a sum of monomials of length $<k$ and $\widehat{\Lomega_{j_\ell}s^{-n_\ell}}$ means that $\Lomega_{j_\ell}s^{-n_\ell}$ does not appear in the product. From Formula~\eqref{eq:commutation-lemma} we see that the only terms contributing to monomials of length $N$ come from $Q_N$ and that these terms vanish only if $Q_N=0$. The proof follows by induction.
\end{proof}
Lemma~\ref{lem:injectivity-ad} together with Formula~\eqref{eq:relation theta(e_i) and Yj} imply Statement~\eqref{item:class_irred_component-1} in the general case.

\medskip

Now, we address Statement~\eqref{item:class_irred_component-2} for an arbitrary affine quiver $\qv$. We will be dealing at the same time with objects corresponding to type $A_1^{(1)}$ (resp.\ $T^\ast \PP^1$) and $\qv$ (resp.\ $\rsv$): we will distinguish these by adding the subscript $T^\ast\PP^1$ or $\rsv$ in the context of surfaces, and $\Kr$ or $\qv$ in the context of quivers. 

Let us denote by $\pr_i\colon \coweightlattice_{\sff,\qv} \to \coweightlattice_{\sff,\Kr}$ the map defined by $\Lomega_j \mapsto \delta_{i,j}\Lomega_1$. We will fix stability parameters $\Ltheta_\qv \in \coweightlattice_\qv$ and $\Ltheta_\Kr \in \coweightlattice_\Kr$ satisfying Formula~\eqref{eq:Lthetaconditions} such that $\pr_i(\Ltheta_{\sff,\qv})=\Ltheta_{\sff,\Kr}$.

By Theorem~\ref{thm:existscoha}--\ref{item:COHA-ii} and --\ref{item:COHA-iii} (cf.\ \cite[Theorem~\ref*{foundation-thm:nilpotent-COHA}--(\ref*{foundation-item:nilpotent-COHA-1}) and --(\ref*{foundation-item:nilpotent-COHA-2})]{DPSSV-1}), the isomorphism $\widehat{T^\ast\PP^1}_{\PP^1} \simeq \widehat{\rsv}_{C_i}$ and the embedding $C_i \to C$ give rise to an algebra morphism 
\begin{align}
	\iota_\ast\colon \coha_{T^\ast\PP^1,\PP^1}^A\simeq \coha_{\rsv,C_i}^A \longrightarrow\coha_{\rsv,C}^A
\end{align}
which sends the fundamental class of an irreducible component $\frakZ$ of $\stackCohpsnil(\widehat{T^\ast\PP^1}_{\PP^1})\simeq\stackCohpsnil(\widehat{\rsv}_{C_i})$ to the fundamental class of $\frakZ$ viewed as an irreducible component of $\stackCohpsnil(\widehat{\rsv}_C)$. Restricting to the part of vertical degree zero and composing with the algebra isomorphism $\Theta_{\rsv,C}$, we thus obtain a morphism 
\begin{align}
	\iota_{\ast, 0}\colon \widehat{\sfU}(\tensor*[_\Kr]{\frakn}{_{\mathsf{ell}, 0}^+}) \longrightarrow \widehat{\sfU}(\tensor*[_\qv]{\frakn}{_{\mathsf{ell}, 0}^+})\ ,
\end{align}
where $\tensor*[_\qv]{\frakn}{_{\mathsf{ell}, 0}^+} \subset \widehat{\frakg}_\sff$ is defined in Formula~\eqref{eq:defun'} and $\tensor*[_\Kr]{\frakn}{_{\mathsf{ell}, 0}^+}\coloneqq \Q e[s^{\pm 1}] \oplus s^{-1}\Q h[s^{-1}] \subset \widehat{\fraksl}_2$. 

Note the following properties:
\begin{enumerate}[label=(\roman*)]\itemsep0.2cm
	\item \label{item:iota-i} the functor $\iota_\ast\colon \catCoh_{C_i}(\rsv)\to \catCoh_C(\rsv)$ is equivariant with respect to tensoring by line bundles in the sense that $\calL_{\Llambda} \otimes \iota_\ast(\calF) \simeq \iota_\ast(\calL_{\pr_i(\Llambda)} \otimes \calF)$ for any $\calF\in  \catCoh_{C_i}(\rsv)$ and $\Llambda\in \coweightlattice_{\sff,\qv}$,
	\item \label{item:iota-ii} this functor restricts to fully faithful functors
	\begin{align}
		\catCoh_{\PP^1}(T^\ast \PP^1)^{> 0} \to  \catCoh_C(\rsv)^{> 0} \quad \text{and} \quad \catCoh_{\PP^1}(T^\ast \PP^1)^{\leqslant 0} \to \catCoh_C(\rsv)^{\leqslant 0} \ .
	\end{align}
\end{enumerate}

Statement~\eqref{item:class_irred_component-2} for the surface $\rsv$ now amounts to the equalities
\begin{align}\label{eq:image_by_iota_Lie_algebra1}
	\iota_{\ast, 0}(es^{n})&=x_i^+s^n \ ,\\[2pt]
\label{eq:image_by_iota_Lie_algebra2}
	\iota_{\ast, 0}(hs^n)&=h_is^n
\end{align}
for all $n$. Equality~\eqref{eq:image_by_iota_Lie_algebra2} being a consequence of Statement~\eqref{item:class_irred_component-1}, we focus on Equality~\eqref{eq:image_by_iota_Lie_algebra1}. For weight reasons, we have
\begin{align}
	\iota_{\ast, 0}(es^{n})=\sum_{j \geq 0} x_i^+s^{n+j}P_{j,n}(\underline{h})
\end{align} 
where $P_{j,n}(\underline{h}) \in \sfU(s^{-1}\frakh_0[s^{-1}])$ is of weight $-j\delta$. Using the property~\ref{item:iota-i} and Proposition~\ref{prop:coha-tensor-product-intertwine} we deduce that $\iota_{\ast, 0}$ is equivariant with respect to the action of $\coweightlattice_{\sff,\qv}$. From this and from the relations (see Proposition~\ref{prop:action_translation_braid_Lie})
\begin{align}
	L_{\Llambda}(es^{n})&=(-1)^{\langle \alpha_1,\Llambda\rangle}es^{n-(\alpha_1,\Llambda)} \quad \text{and} \quad 	L_{\Llambda}(hs^{n})=hs^n \ ,\\
	L_{\Lmu}(x_i^+s^{-n})&=(-1)^{\langle \alpha_i,\Lmu\rangle}x_i^+s^{n-(\alpha_i,\Lmu)} \quad \text{and} \quad L_{\Lmu}(h_is^{n})=h_is^n
\end{align}
for any $\Llambda \in \coweightlattice_{\sff,\Kr}, \Lmu \in \coweightlattice_{\sff,\qv}$, we deduce that $P_{j,n}(\underline{h})$ is independent of $n$ (we henceforth drop the $n$ from the notation).

We will now use nilpotent cohomological Hall algebras of affine quivers. From property~\ref{item:iota-ii}, it follows that functor induced by $\iota_\ast$ at the level of bounded derived categories preserves the perverse coherent hearts, hence it induces a fully faithful functor 
\begin{align}
	\tensor*[_\tau]{\iota}{}\colon \catmod(\Pi_{\Kr}) \longrightarrow \modPi \ .
\end{align}
This gives rise to a map of derived stacks $\tensor*[_\tau]{\iota}{}\colon\dLambda_{\Kr} \to \dLambda_\qv$, which by abuse of notation we denote in the same way. The morphism $\tensor*[_\tau]{\iota}{}$ restricts to a closed embedding $\tensor*[_\tau]{\iota}{_{\leq 0}} \colon \dLambda_{\Kr}^{\leqslant 0} \to \dLambda_\qv^{\leqslant 0}$ which fits in the commutative diagram
\begin{align}
	\begin{tikzcd}[ampersand replacement=\&, column sep=large]
		\Lambda^{\leq 0}_{\Kr} \ar{r}{\tensor*[_\tau]{\iota}{_{\leq 0}}} \ar{d}{\simeq} \& \Lambda^{\leq 0}_\qv \ar{d}{\simeq}\\
		\dstackCohpsnil(\widehat{\rsv}_{C_i})^{>0} \ar{r}{\iota} \& \dstackCohpsnil(\widehat{\rsv}_{C})^{>0}
	\end{tikzcd}\ .
\end{align}
We have identifications
\begin{align}
	\sfH_{\mathsf{top}}^A(\dLambda^{\leqslant 0}_{\Kr})\simeq \sfU(\frakn_{\Kr}^-)/I_{\Kr}\quad \text{and} \quad \sfH_{\mathsf{top}}^A(\dLambda^{\leqslant 0}_\qv)\simeq \sfU(\frakn_\qv^-)/I_\qv \ ,
\end{align}
where $I_{\Kr}$ and $I_\qv$ are the left ideals respectively generated by
\begin{align}
	\bigoplus_{\bfd \in -\Delta_{\mathsf{f,Kr}} \times \N \delta} \frakn_{\Kr,\bfd} \quad \text{and} \bigoplus_{\bfd \in -\Delta_{\sff,\qv} \times \N \delta} \frakn_{\qv,\bfd} \ .
\end{align}
In particular, the projection map $\pi_0 \colon \widehat{\sfU}(\tensor*[_\qv]{\frakn}{_{\mathsf{ell}, 0}^+}) \to \sfU(\tensor*[_\qv]{\frakn}{_{\mathsf{ell}, 0}^+})/I_\qv$ satisfies
\begin{align}
	\pi_0(\iota_{\ast, 0}(es^{n}))=\sum_{\textrm{-}n > j \geq 0} x_i^+s^{n+j}P_{j}(\underline{h})\ .
\end{align}

To conclude, we will use the adjoint action of the  $\fraksl_2$ subalgebra corresponding to the simple root $\alpha_i$. By construction, 
\begin{align}\label{eq:imageirrcomp}
	\tensor*[_\tau]{\iota}{}(S_1)=\tensor*[_\tau]{\iota}{}(S_1[-1])[1]=\tensor*[_\tau]{\iota}{}(\scrO_{\PP^1}(-1))[1]=\scrO_{C_i}(-1)[1]=S_i
\end{align}
hence at the stacky level the map $\tensor*[_\tau]{\iota}{}$ restricts to an isomorphism 
\begin{align}
	\begin{tikzcd}[ampersand replacement=\&]
		(\dLambda_{\Kr})_{\N\alpha_1} \ar{r}{\sim} \& (\dLambda_\qv)_{\N\alpha_i}
	\end{tikzcd}\ .
\end{align}
It is easy to see that the induced map $\tensor*[_\tau]{\iota}{_{\leq 0, \ast}}\colon \HBMbulletA(\dLambda^{\leq 0}_{\Kr}) \to \HBMbulletA(\dLambda^{\leq 0}_\qv)$ is equivariant with respect to the left action of the algebras $\HBMbulletA(\dLambda_{\Kr})_{\N\alpha_1} \simeq \HBMbulletA(\dLambda_\qv)_{\N\alpha_i}$ (see Proposition~\ref{prop:various_COHA_actions_truncated}--\eqref{item:yangian-kappa-action2}). Note that 
\begin{align}
	f u = [f,u] \in \sfU(\frakn_{\Kr}^-)/I_{\Kr}\quad \text{and} \quad f_i v = [f_i,v] \in \sfU(\frakn_\qv^-)/I_\qv
\end{align}
for any $u \in \sfU(\frakn_{\Kr}^-)$ and $v \in \sfU(\frakn_\qv^-)$. Hence we get
\begin{align}
	\tensor*[_\tau]{\iota}{_{\leq 0, \ast}}(fes^n)=\tensor*[_\tau]{\iota}{_{\leq 0, \ast}}([f,es^n])=-2\tensor*[_\tau]{\iota}{_{\leq 0, \ast}}(hs^n)=-2h_is^n\ ,
\end{align}
and also
\begin{align}
	\tensor*[_\tau]{\iota}{_{\leq 0, \ast}}(fes^n)=f_i\tensor*[_\tau]{\iota}{_{\leq 0, \ast}}(es^n)=\left[f_i,\sum_{\textrm{-}n > j \geq 0} x_i^+s^{n+j}P_{j}(\underline{h})\right]=-2\sum_{\textrm{-}n > j \geq 0} h_is^{n+j}P_{j}(\underline{h})\ .
\end{align}
As this is true for all $n<0$ we deduce that $P_0(\underline{h})=1$ and $P_j(\underline{h})=0$ for all $j>0$ as wanted. This concludes the proof of Statement~\eqref{item:class_irred_component-2} in the general case. Thus, Theorem~\ref{thm:class_irred_component} is proved.

\section{Appendix: Proof of Proposition~\ref{prop:characterization-classical-limit}}\label{sec:characterization-classical-limit}

Let $\fraks_\qv$ be the Lie algebra introduced in Definition~\ref{def:classical-Lie-algebra}. Let $\fraks_\qv^0$ be the Lie subalgebra of $\fraks_\qv$ generated by the elements $H_{i, \ell}$ with $i\in I$ and $\ell\in \N$, while let $\fraks_\qv^\pm$ be the Lie subalgebra of $\fraks_\qv$ generated by the elements $X_{i, \ell}^\pm$ with $i\in I$ and $\ell\in \N$. 

Let us first recall the statement of Proposition~\ref{prop:characterization-classical-limit}.
\begin{proposition}
	\hfill
	\begin{itemize}\itemsep0.2cm
		\item $\fraks_\qv^0$ is isomorphic to the free commutative Lie algebra generated by $H_{i, \ell}$, with $i\in I$ and $\ell\in \N$.
		\item $\fraks_\qv^\pm$ is isomorphic to the Lie algebra generated by the elements $X_{i, \ell}^\pm$ with $i\in I$ and $\ell\in \N$ satisfying the relations:
		\begin{enumerate}\itemsep0.2cm
			\item \label{item:tilde-relations-1-proof}
			\begin{align}\label{eq:1}
				\Big[X_{i, r+1}^{\pm}, X_{j, s}^{\pm}\Big] = \Big[X_{i, r}^{\pm}, X_{j, s+1}^{\pm}\Big]  \ ,
			\end{align}
			and 
			\begin{align}\label{eq:2}
				\ad\Big(X_{i, 0}^\pm\Big)^{1-a_{i,j}} \Big(X_{j, r}^\pm\Big)=0  \quad \text{for } i \neq j\ ,
			\end{align}
			if $\qv\neq A_1^{(1)}$;
			
			\item \label{item:tilde-relations-2-proof} \eqref{eq:1} for $i=j$, \eqref{eq:2}, 
			\begin{align}\label{eq:3}
				\Big[X_{i,r+2}^{\pm}, X_{i+1,s}^{\pm}\Big] - 2\Big[X_{i,r+1}^{\pm}, X_{i+1,s+1}^{\pm}\Big] + \Big[X_{i,r}^{\pm}, X_{i+1,s+2}^{\pm}\Big]  =0\ ,
			\end{align}
			and
			\begin{align}\label{eq:4}
				\Big[\Big[X_{i, r+1}^{\pm}, X_{j, s}^{\pm}\Big], X_{j, \ell}^\pm\Big] - \Big[\Big[X_{i, r}^{\pm}, X_{j, s+1}^{\pm}\Big], X_{j, \ell}^\pm\Big] =0\ ,
			\end{align}
			if $\qv=A_1^{(1)}$.
		\end{enumerate}
	\end{itemize}
\end{proposition}

We will provide a proof of the proposition in the remaining part of the section.

Denote by $\tilde \fraks_\qv^0$ the free commutative Lie algebra generated by $H_{i, \ell}$, with $i\in I$ and $\ell\in \N$, and by $\tilde \fraks_\qv^\pm$ the Lie algebra generated by the elements $X_{i, \ell}^\pm$ with $i\in I$ and $\ell\in \N$ satisfying the relations \eqref{item:tilde-relations-1-proof} and \eqref{item:tilde-relations-2-proof}. We shall prove that $\fraks_\qv^0 \simeq \tilde \fraks_\qv^0$ and $\fraks_\qv^\pm \simeq \tilde \fraks_\qv^\pm$.

Define the Lie algebra $\tensor*[_1]{\fraks}{_\qv}$ as the Lie algebra generated by $X_{i, \ell}^\pm$ and $H_{i, \ell}$, with $i\in I$ and $\ell \in \N$, subject to the relations \eqref{eq:classical-1}, \eqref{eq:classical-2}, and \eqref{eq:classical-3} in Definition~\ref{def:classical-Lie-algebra}.  Moreover, define $\tensor*[_1]{\tilde \fraks}{_\qv^\pm}$ as the free Lie algebra generated by $X_{i, \ell}^\pm$ and $\tensor*[_1]{\tilde \fraks}{_\qv^0}$ as the free commutative Lie algebra generated by $H_{i, \ell}$, with $i\in I$ and $\ell \in \N$. On the other hand, define the Lie algebras $\tensor*[_2]{\fraks}{_\qv}$ and $\tensor*[_2]{\tilde \fraks}{_\qv^\pm}$ as $\fraks_\qv$ and $\tilde \fraks_\qv^\pm$ with the relation \eqref{eq:2} removed. 

For $i=1, 2$, let $\tensor*[_i]{\fraks}{_\qv^\pm}$ and $\tensor*[_i]{\fraks}{_\qv^0}$ be the Lie subalgebras of $\tensor*[_i]{\fraks}{_\qv}$ generated by $X_{i, \ell}^\pm$ and $H_{i, \ell}$, with $i\in I$ and $\ell\in \N$, respectively. We get a triangular decomposition
\begin{align}
	\tensor*[_1]{\fraks}{_\qv}= \tensor*[_1]{\fraks}{_\qv^+} \oplus \tensor*[_1]{\fraks}{_\qv^0}\oplus \tensor*[_1]{\fraks}{_\qv^-} \ .
\end{align}
Note that $\tensor*[_1]{\tilde \fraks}{_\qv^\pm}\simeq \tensor*[_1]{\fraks}{_\qv^\pm}$ and $\tensor*[_1]{\tilde \fraks}{_\qv^0}\simeq \tensor*[_1]{\fraks}{_\qv^0}$. In particular, the above triangular decomposition becomes
\begin{align}
	\tensor*[_1]{\fraks}{_\qv}= \tensor*[_1]{\tilde \fraks}{_\qv^+} \oplus \tensor*[_1]{\tilde \fraks}{_\qv^0}\oplus \tensor*[_1]{\tilde \fraks}{_\qv^-} \ .
\end{align}	
Let $\tensor*[_1]{\tau}{^\pm}$ be the ideal of $\tensor*[_1]{\fraks}{_\qv^\pm}$ such that $\tensor*[_2]{\tilde \fraks}{_\qv^\pm}=\tensor*[_1]{\fraks}{_\qv^\pm}/\tensor*[_1]{\tau}{^\pm}$, i.e., $\tensor*[_1]{\tau}{^\pm}$ is generated by the relation \eqref{eq:1} if $\qv\neq A_1^{(1)}$, while is generated by \eqref{eq:1} for $i=j$, \eqref{eq:3} and \eqref{eq:4} if $\qv=A_1^{(1)}$.

\begin{lemma}
	$[\tensor*[_1]{\fraks}{_\qv} , \tensor*[_1]{\tau}{^\pm}]\subset \tensor*[_1]{\tau}{^\pm}$.
\end{lemma}

\begin{proof}
	First, assume that $\qv=A_1^{(1)}$. Recall that $\tensor*[_1]{\tau}{^\pm}$ is the Lie ideal of $\tensor*[_1]{\fraks}{_\qv^\pm}$ generated by
	\begin{align}
		A_{i, r, s}&\coloneqq \Big[X_{i, r+1}^{\pm}, X_{i, s}^{\pm}\Big] - \Big[X_{i, r}^{\pm}, X_{i, s+1}^{\pm}\Big]\ ,\\[4pt]
		A_{i, r, s}' &\coloneqq \Big[X_{i,r+2}^{\pm}, X_{i+1,s}^{\pm}\Big] - 2\Big[X_{i,r+1}^{\pm}, X_{i+1,s+1}^{\pm}\Big] + \Big[X_{i,r}^{\pm}, X_{i+1,s+2}^{\pm}\Big] \ ,\\[4pt]
		A_{i, j, r, s, \ell}&\coloneqq \Big[\Big[X_{i, r+1}^{\pm}, X_{j, s}^{\pm}\Big], X_{j, \ell}^\pm\Big] - \Big[\Big[X_{i, r}^{\pm}, X_{j, s+1}^{\pm}\Big], X_{j, \ell}^\pm\Big]  \ ,
	\end{align}
	with $i\neq j\in I$ and $r, s, \ell\in \N$. Choose $k\in I$ and $u\in \N$. The element $\Big[H_{k, u}, A_{i, r, s}\Big]$ is a linear combination of $A_{i, r+u,s}$ and $A_{i, r, s+u}$, the element $\Big[H_{k, u}, A_{i, r, s}'\Big]$ is a linear combination of $A_{i, r+u, s}$ and $A_{i, r, s+u}$, while $\Big[H_{k, u}, A_{i, j, r, s, \ell}\Big]$ is a linear combination of $A_{i, j, r+u, s, \ell}$, $A_{i, j, r,s+u, \ell}$, and $A_{i, j, r, s, \ell+u}$. Hence, 
	$\Big[H_{k, u}, A_{i, r, s}\Big], \Big[H_{k, u}, A_{i, r, s}'\Big], \Big[H_{k, u}, A_{i, j, r, s, \ell}\Big]$ belong to $\tensor*[_1]{\tau}{^\pm}$.
	
	Further, note that
	\begin{align}
		\Big[X_{k, u}^\mp, \Big[X_{i, r}^{\pm}, X_{j, s}^{\pm}\Big]\Big]&=-\Big[X_{j, s}^{\pm},\Big[X_{k, u}^\mp, X_{i, r}^{\pm}\Big]\Big]-\Big[X_{i, r}^{\pm},\Big[X_{j, s}^{\pm},X_{k, u}^\mp \Big]\Big]\\[4pt]
		&=-\delta_{k, i}\Big[X_{j, s}^{\pm}, H_{i, u+r}\Big]-\delta_{k, j}\Big[X_{i, r}^{\pm}, H_{j, u+s}\Big]\\[4pt]
		&= \delta_{k, i}a_{i, j}X_{j, r+s+u}^\pm+\delta_{k, j} a_{j, i}X_{i, r+u+s}^\pm\ .
	\end{align}
	Thus, $\Big[X_{k, u}^\mp, A_{i, r, s}\Big], \Big[X_{k, u}^\mp, A_{i, r, s}'\Big], \Big[X_{k, u}^\mp, A_{i, j, r, s, \ell}\Big]$ vanish, hence they belong to $\tensor*[_1]{\tau}{^\pm}$. Thus, an easy induction on the degree of any element $X^\pm$ in $\tensor*[_1]{\fraks}{_\qv^\pm}$ implies that all elements
	\begin{align}
		\Big[H_{k, u}, \Big[X^\pm, B\Big]\Big] \quad \text{and} \quad \Big[X_{k, u}^\mp, \Big[X^\pm, B\Big]\Big]
	\end{align}
	with $k\in I$, $u\in \N$, and $B= A_{i, r, s}, A_{i, r, s}', A_{i, j, r, s, \ell}$, belongs to $\tensor*[_1]{\tau}{^\pm}$. The proof when $\qv\neq A_1^{(1)}$ uses the same arguments and it is simpler.
\end{proof}	

Thanks to the above lemma, the image of $\tensor*[_1]{\tau}{^\pm}$ by the inclusion $\tensor*[_1]{\fraks}{_\qv^\pm}\subset \tensor*[_1]{\fraks}{_\qv}$ is an ideal. We deduce that
\begin{align}
	\tensor*[_2]{\fraks}{_\qv} &=\tensor*[_1]{\fraks}{_\qv} / \big(\tensor*[_1]{\fraks}{_\qv^+} + \tensor*[_1]{\fraks}{_\qv^-} \big)\\
	&=\big(\tensor*[_1]{\fraks}{_\qv^+} \oplus \tensor*[_1]{\fraks}{_\qv^0}\oplus \tensor*[_1]{\fraks}{_\qv^-}\big)/\tensor*[_1]{\fraks}{_\qv} \big(\tensor*[_1]{\fraks}{_\qv^+} + \tensor*[_1]{\fraks}{_\qv^-} \big)\\
	&=\big(\tensor*[_1]{\fraks}{^+}/\tensor*[_1]{\fraks}{^+} \big)\oplus \tensor*[_1]{\fraks}{^0}\oplus \big( \tensor*[_1]{\fraks}{^-}/\tensor*[_1]{\fraks}{^-} \big)\\
	&=\tensor*[_2]{\tilde \fraks}{_\qv^+} \oplus \tensor*[_2]{\tilde \fraks}{_\qv^0}\oplus \tensor*[_2]{\tilde\fraks}{_\qv^-}\ .
\end{align}
In particular, we also have $\tensor*[_2]{\tilde \fraks}{_\qv^\pm}=\tensor*[_2]{\fraks}{_\qv^\pm}$ and $\tensor*[_2]{\tilde \fraks}{_\qv^0}=\tensor*[_2]{\tilde \fraks}{_\qv^0}$.

Let $\tensor*[_2]{\fraks}{_\qv^\pm}$ be the ideal of $\tensor*[_2]{\fraks}{_\qv^\pm}$ such that $\tilde \fraks_\qv^\pm=\tensor*[_2]{\fraks}{_\qv^\pm}/\tensor*[_2]{\tau}{^\pm}$, i.e., the ideal generated by the relation \eqref{eq:2}.

\begin{lemma}
	$[\tensor*[_2]{\fraks}{_\qv} , \tensor*[_2]{\tau}{^\pm}]\subset \tensor*[_2]{\tau}{^\pm}$.
\end{lemma}

\begin{proof}
	Recall that $\tensor*[_2]{\tau}{^\pm}$ is the Lie ideal of $\tensor*[_2]{\fraks}{_\qv^\pm}$ generated by
	\begin{align}
		B_{i,j,s}\coloneqq \ad\Big(X_{i,0}^\pm\Big)^{1-{a_{i,j}}}\Big(X_{j,s}^\pm\Big)\ ,
	\end{align}
	with $i\neq j\in I$ and $s\in \N$. Choose any $k\in I$ and $u\in \N$. It is easy to see that for any element $X^\pm$ in $\tensor*[_2]{\fraks}{_\qv^\pm}$ the Lie bracket $\Big[H_{k,u},[X^\pm,B_{i,j,s}\Big]\Big]$ belongs to $\tensor*[_2]{\tau}{^\pm}$. To prove that $\Big[X_{k,u}^\mp,\Big[X^\pm,B_{i,j,s}\Big]\Big]$ belongs also to $\tensor*[_2]{\tau}{^\pm}$, by an induction  on the degree of $X^\pm$, it is enough to check that $\Big[X_{k,u}^\mp, B_{i,j,s}\Big]\in \tensor*[_2]{\tau}{^\pm}$. This is a case by case computation for $a_{i,j}=0,-1$ or $-2$.
\end{proof}	

\begin{proof}[Proof of Proposition~\ref{prop:characterization-classical-limit}]
	The previous lemma yields that the image of $\tensor*[_2]{\tau}{^\pm}$ by the inclusion $\tensor*[_2]{\fraks}{_\qv^\pm}\subset \tensor*[_2]{\fraks}{_\qv}$ is an ideal. We deduce as above that
	\begin{align}
		\fraks_\qv=\tilde\fraks_\qv^+\oplus\tilde\fraks_\qv^0\oplus\tilde\fraks_\qv^- \ , \quad
		\tilde\fraks_\qv^\pm=\fraks_\qv^\pm \ , \quad
		\tilde\fraks_\qv^0=\fraks_\qv^0\ .
	\end{align}
	This proves the proposition.
\end{proof}

\section{Appendix: Proof of Proposition~\ref{prop:coha-surface-module-algebra}}\label{sec:hopfact} 

Given a decomposition $(\gamma, n) = (\gamma_1, n_1) + (\gamma_2, n_2)$, for ease of exposition, we will employ the notation 
\begin{align}
	\frakX\coloneqq \dstackCoh({\widehat X}_C; \gamma, n)\qquad\text{and}\qquad  \frakX_i \coloneqq \dstackCoh({\widehat X}_C; \gamma_i, n_i) 
\end{align}
for $1\leq i\leq 2$. Let 
\begin{align}
	\begin{tikzcd}[ampersand replacement=\&]
		\& \overline{\rsv} \times \frakX\ar{dl}[swap]{\ev_{\overline{\rsv}}}\ar{dr}{\ev}\& \\
		\overline{\rsv}\& \& \frakX
	\end{tikzcd}
	\qquad 
	\begin{tikzcd}[ampersand replacement=\&]
		\& \overline{\rsv} \times \frakX_i \ar{dl}[swap]{\ev_{\overline{\rsv}, i}}\ar{dr}{\ev_i}\& \\
		\overline{\rsv}\& \& \frakX_i
	\end{tikzcd}
\end{align}
denote the canonical projections, with $1\leq i\leq 2$.

Given $\alpha \in \sfH^\bullet_A(\overline{\rsv})$, recall that $\overline{\nu}_\alpha \in \sfH^\bullet_A(\frakX)$ 
is defined by 
\begin{align}
	\overline{\nu}_\alpha = \ev_!(\ev_{\overline{\rsv}}^\ast(\alpha) \cup \ch(\overline{\scrE})) 
\end{align}
where $\overline{\scrE}$ is the universal sheaf constructed in \S\ref{subsec:coha-rsv-tautological}. Define
\begin{align}
	\overline{\nu}^{(i)}_\alpha  = (\ev_i)_!(\ev_{\overline{\rsv}, i}^\ast(\alpha)\cup \ch(\overline{\scrE}_i))
\end{align}
where $\overline{\scrE}_i$ is the universal sheaf over $\overline{\rsv} \times \frakX_i$, for $1\leq i \leq 2$.

The proof of Proposition \ref{prop:coha-surface-module-algebra} will use some preliminary results derived from \cite{Khan_VFC, Porta_Yu_Non-archimedian_quantum}. First, as a consequence of  \cite[Proposition~4.10-(3)]{Porta_Yu_Non-archimedian_quantum}, one has the following.
\begin{lemma}\label{lem:hopfactA} 
	For any $x\in \HBMbulletA(\frakX)$, for any invariants $(\gamma, n)$, we have
	\begin{align}
		\overline{\nu}_\alpha \cap x = \ev_\ast((\ev_{\overline{\rsv}}^\ast(\alpha) \cup \ch(\overline{\scrE}))\cap \ev^!(x))\ .
	\end{align}
\end{lemma} 
We next record the following compatibility results for the external product defined in \cite{Porta_Yu_Non-archimedian_quantum}. 
\begin{lemma}\label{lem:qvirtpullback} 
	Let $\frakY$ be an admissible indgeometric derived stack and let $r\colon \overline{\rsv}\times \frakX \to \frakY$ be the canonical projection. Assume that there is a torus action $A$ on $\frakY$ so that $r$ is $A$-equivariant. Let 
	\begin{align}
		\kappa\colon \HBMbulletA(\overline{\rsv}) \otimes \HBMbulletA(\frakY) \longrightarrow \HBMbulletA(\overline{\rsv}\times \frakY) 
	\end{align}
	be the external product defined in \cite[Definition~3.10-(2)]{Porta_Yu_Non-archimedian_quantum}. Then 
	\begin{align}\label{eq:qvirtpullback} 
		r^!(x) = \kappa([\overline{\rsv}]\otimes y)
	\end{align}
	for all $y \in \HBMbulletA(\frakY)$. 
\end{lemma} 

\begin{proof} 
	By \cite[Remark 3.7]{Khan_VFC}, one has $r^!(x)= r^!(1) \circ y$, where 
	\begin{align}
		\circ\colon \sfH^A_k(\overline{\rsv}\times \frakY/\frakY) \otimes \sfH^A_{k'}(\frakY) \longrightarrow \sfH^A_{k+k'}(\overline{\rsv}\times \frakY) 
	\end{align}
	is the composition product defined in \cite[\S2.2.5]{Khan_VFC}. Moreover, applying \cite[Proposition~4.7]{Porta_Yu_Non-archimedian_quantum} to the pull-back diagram 
	\begin{align}
		\begin{tikzcd}[ampersand replacement=\&]
			\overline{\rsv} \times \frakY \ar{d}{} \ar{r}{r} \& \frakY\ar{d}{b}\\
			\overline{\rsv} \ar{r}{a} \& \Spec(\C)
		\end{tikzcd}\ ,
	\end{align}
	one obtains the identity $b^\ast([\overline{\rsv}]) = r^!(1)$. Here $b^\ast \colon \HBMbulletA(\frakY) \to \HBMbulletA(\overline{\rsv}\times \frakY)$ is the change of base map defined in \cite[\S2.2.3]{Khan_VFC}. 
	
	Finally, applying  \cite[Proposition~4.10-(4)]{Porta_Yu_Non-archimedian_quantum} to the same  pull-back diagram 
	yields 
	\begin{align}
		b^\ast([\overline{\rsv}]) \circ y = \kappa([\overline{\rsv}] \otimes y)\ .
	\end{align}
	This proves the claim.
\end{proof} 

\begin{lemma}\label{lem:functkunneth} 
	Let 
	\begin{align}
		\kappa &\colon \HBMbulletA(\frakX_1)\otimes \HBMbulletA(\frakX_2) \longrightarrow
		\HBMbulletA(\frakX_1\times \frakX_2)
	\end{align}
	and
	\begin{align}	
		\kappa_1&\colon \HBMbulletA(\overline{\rsv}\times \frakX_1) \otimes \HBMbulletA(\frakX_2) \longrightarrow 
		\HBMbulletA(\overline{\rsv} \times \frakX_1\times \frakX_2) \ ,\\
		\kappa_2 &\colon \HBMbulletA(\frakX_1) \otimes \HBMbulletA(\overline{\rsv}\times \frakX_2) \longrightarrow
		\HBMbulletA(\overline{\rsv} \times \frakX_1\times \frakX_2) 
	\end{align}
	be the external products defined in \cite[Definition~3.10-(2)]{Porta_Yu_Non-archimedian_quantum}. Let 
	\begin{align}
		r_i\colon \overline{\rsv} \times \frakX_1\times \frakX_2 \longrightarrow \overline{\rsv}  \times \frakX_i 
	\end{align}
	for $1\leq i \leq 2$, and 
	\begin{align}
		\ev_{12}\colon \overline{\rsv} \times \frakX_1 \times \frakX_2\longrightarrow \frakX_1 \times \frakX_2 
	\end{align}
	be the canonical projections. Then, the following identities hold:
	\begin{align}\label{eq:kunnethcapA}
		r_1^\ast(c_1) \cap \kappa_1(\ev_1^!(x_1)\otimes x_2) &= \kappa_1((c_1 \cap \ev_1^!(x_1))\otimes x_2)\ ,\\
		\label{eq:kunnethcapB}
		r_2^\ast(c_2) \cap \kappa_2(x_1 \otimes \ev_2^!(x_2)) &= \kappa_2(x_1 \otimes (c_2 \cap \ev_2^!(x_2))\ ,
	\end{align}
	for any $x_i \in \HBMbulletA(\frakX_i)$ and any $c_i \in \Hbullet_A(\overline{\rsv} \times \frakX_i)$, with $1\leq i\leq 2$.
	
	Furthermore, for any $x_i \in \HBMbulletA(\frakX_i)$ and any $y_i \in \HBMbulletA(\overline{\rsv} \times \frakX_i)$, with $1\leq i\leq 2$, one also has 
	\begin{align}\label{eq:kunnethprojA} 
		(\ev_{12})_\ast(\kappa_1(y_1\otimes x_2)) &= \kappa((\ev_1)_\ast(y_1)\otimes x_2) \ , \\ \label{eq:kunnethprojB}
		(\ev_{12})_\ast(\kappa_2(x_1\otimes y_2)) &= \kappa(x_1 \otimes (\ev_2)_\ast(y_2))\ . 
	\end{align} 
\end{lemma} 

\begin{proof} 
	Identity~\eqref{eq:kunnethcapA} is a consequence of \cite[Proposition~4.10-(5)]{Porta_Yu_Non-archimedian_quantum}, using the identity $(r_1)_\ast(c_1) =\kappa_1( c_1 \otimes 1_{\frakX_2})$, where $1_{\frakX_2}$ is the unit element in the cohomology ring of ${\frakX_2}$. The proof of identity~\eqref{eq:kunnethcapB} is completely analogous. 
	
	\medskip
	
	Identity~\eqref{eq:kunnethprojA} follows from \cite[Proposition 4.10-(6)]{Porta_Yu_Non-archimedian_quantum} applied with 
	$X = \frakX_1$, $Y = \frakX_2$,  $X' = \overline{\rsv} \times \frakX_1$, $Y' = \frakX_2$, $f = \ev_1$, and $ g =\id_{\frakX_2}$. Then we get $h = \ev_{12}$, and the cited statement gives exactly identity~\eqref{eq:kunnethprojA}. The proof of identity~\eqref{eq:kunnethprojB} is completely analogous. 
\end{proof} 

Consider now the commutative diagram 
\begin{align}\label{eq:hopfdiagA}
	\begin{tikzcd}[ampersand replacement=\&, column sep=large]
		\&\& \overline{\rsv}\&\& \\
		\& \& \overline{\rsv} \times \frakX_1\times \frakX_2 \ar{dl}[swap]{r_1}\ar{dr}{r_2}
		\ar{dd}{\ev_{12}}\ar{u}{\ev_{\overline{\rsv}, 12}} \& \& \\
		\overline{\rsv} \& \overline{\rsv}\times \frakX_1\ar{dd}{\ev_1}\ar{l}[swap]{\ev_{\overline{\rsv}, 1}} \& \& \overline{\rsv} \times 
		\frakX_2\ar{dd}{\ev_2}\ar{r}{\ev_{\overline{\rsv}, 2}} \& \overline{\rsv} \\
		\&\&\frakX_1\times \frakX_2 \ar{dl}[swap]{\pi_1} 
		\ar{dr}{\pi_2}\&\& \\
		\&\frakX_1 \& \& \frakX_2\&
	\end{tikzcd}
\end{align}
where both squares are pull-back. As in Lemma~\ref{lem:qvirtpullback}, let 
\begin{align}
	\kappa\colon \HBMbulletA(\frakX_1)\otimes \HBMbulletA(\frakX_2) \longrightarrow\HBMbulletA(\frakX_1\times \frakX_2) 
\end{align}
denote the external product defined in \cite[Definition~3.10-(2)]{Porta_Yu_Non-archimedian_quantum}.  Abusing notation, set also $\overline{\nu}_\alpha(x) \coloneqq \overline{\nu}_\alpha\cap x$ for any $\alpha \in \sfH^\bullet_A(\overline{\rsv})$ and $x\in \sfH^\bullet_A(\frakX)$. Then, the following holds.
\begin{lemma}\label{lem:hopfB} 
	One has the identity 
	\begin{multline}\label{eq:hopfB} 
		\kappa(\overline{\nu}^{(1)}_\alpha(x_1) \otimes \overline{\nu}^{(2)}_\alpha(x_2))= \\
		(\ev_{12})_\ast\big(((r_1^\ast\ch(\overline{\scrE}_1) + r_2^\ast\ch(\overline{\scrE}_2))\cup \ev_{\overline{\rsv}, 12}^\ast(\alpha) )
		\cap \ev_{12}^!(\kappa(x_1\otimes x_2)))
	\end{multline}
	for any $x_i \in \HBMbulletA(\frakX_i)$, with $1\leq i\leq 2$. 
\end{lemma} 

\begin{proof} 
	By substitution, the left hand side of identity \eqref{eq:hopfB} reads 
	\begin{align} 
		\kappa(\ch^{(1)}_\alpha(x_1) \otimes \ch^{(2)}_\alpha(x_2))=& 
		\kappa\big( (\ev_1)_\ast(\ch(\overline{\scrE}_1)\cup \ev_{\overline{\rsv}, 1}^\ast(\alpha)) \cap \ev_1^!(x_1)) \otimes x_2\\ &+ 
		x_1\otimes (\ev_2)_\ast(\ch(\overline{\scrE}_2)\cup \ev_{\overline{\rsv}, 2}^\ast(\alpha)) \cap \ev_2^!(x_2)))\big)\ . 
	\end{align} 
	In order to rewrite the right hand side of  \eqref{eq:hopfB}, as in Lemma~\ref{lem:functkunneth}, let 
	\begin{align}
		\kappa_1\colon \HBMbulletA(\overline{\rsv}\times \frakX_1) \otimes \HBMbulletA(\frakX_2) \longrightarrow 
		\HBMbulletA(\overline{\rsv} \times \frakX_1\times \frakX_2)\\
		\kappa_2\colon  
		\HBMbulletA(\frakX_1) \otimes \HBMbulletA(\overline{\rsv}\times\frakX_2) \longrightarrow
		\HBMbulletA(\overline{\rsv} \times \frakX_1\times \frakX_2) 
	\end{align}
	be the external products defined in \cite[Definition~3.10-(2)]{Porta_Yu_Non-archimedian_quantum}. 
	Since the external product is associative, Lemma~\ref{lem:qvirtpullback} yields the identities 
	\begin{align}\label{eq:kunnethidA} 
		\ev_{12}^!(\kappa(x_1\otimes x_2)) = \kappa_1(\ev_1^!(x_1)\otimes x_2)\quad \text{and}\quad 
		\ev_{12}^!(\kappa(x_1\otimes x_2)) = \kappa_2(x_1 \otimes \ev_2^!(x_2))\ . 
	\end{align}
	Then the right hand side of identity~\eqref{eq:hopfB} reads 
	\begin{multline} 
		(\ev_{12})_\ast\big(((r_1^\ast\ch(\overline{\scrE}_1) + r_2^\ast\ch(\overline{\scrE}_2))\cup \ev_{\overline{\rsv}, 12}^\ast(\alpha) )
		\cap \ev_{12}^!(\kappa(x_1\otimes x_2)))
		=\\
		(\ev_{12})_\ast\big(r_1^\ast(\ch(\overline{\scrE}_1)\cup \ev_{\overline{\rsv}, 1}^\ast(\alpha)) \cap \kappa_1(\ev_1^!(x_1)\otimes x_2)
		+ r_2^\ast(\ch(\overline{\scrE}_2)\cup \ev_{\overline{\rsv}, 2}^\ast(\alpha)) \cap \kappa_2(x_1 \otimes 
		\ev_2^!(x_2))\big)\ .
	\end{multline} 
	Next, note that the first part of Lemma~\ref{lem:functkunneth} yields the identities 
	\begin{align}
		r_1^\ast(\ch(\overline{\scrE}_1)\cup \ev_{\overline{\rsv}, 1}^\ast(\alpha)) \cap \kappa_1(\ev_1^!(x_1)\otimes x_2)= 
		\kappa_1(((\ch(\overline{\scrE}_1)\cup \ev_{\overline{\rsv}, 1}^\ast(\alpha))\cap \ev_1^!(x_1))\otimes x_2)\ ,\\
		r_2^\ast(\ch(\overline{\scrE}_2)\cup \ev_{\overline{\rsv}, 2}^\ast(\alpha)) \cap \kappa_2(x_1 \otimes 
		\ev_2^!(x_2)) = \kappa_2(x_1 \otimes ((\ch(\overline{\scrE}_2)\cup \ev_{\overline{\rsv}, 2}^\ast(\alpha)) \cap \ev_2^!(x_2)))\ .
	\end{align}
	Hence, setting 
	\begin{align}
		y_i = (\ch(\overline{\scrE}_i)\cup \ev_{\overline{\rsv}, i}^\ast(\alpha))\cap \ev_i^!(x_i)\ ,
	\end{align}
	for $1\leq i \leq 2$, one obtains
	\begin{multline} 
		(\ev_{12})_\ast\big(((r_1^\ast\ch(\overline{\scrE}_1) + r_2^\ast\ch(\overline{\scrE}_2))\cup \ev_{\overline{\rsv}, 12}^\ast(\alpha) )\cap \ev_{12}^!(\kappa(x_1\otimes x_2)))	= \\
		(\ev_{12})_\ast \big(\kappa_1(y_1\otimes x_2) + \kappa_2(x_1 \otimes y_2)\big)\ .
	\end{multline}
	Using the second part of Lemma~\ref{lem:functkunneth} we further obtain 
	\begin{multline} 
		(\ev_{12})_\ast\big(((r_1^\ast\ch(\overline{\scrE}_1) + r_2^\ast\ch(\overline{\scrE}_2))\cup \ev_{\overline{\rsv}, 12}^\ast(\alpha) )
		\cap \ev_{12}^!(\kappa(x_1\otimes x_2)))= \\ 
		\kappa((\ev_1)_\ast(y_1)\otimes x_2) + \kappa(x_1 \otimes (\ev_2)_\ast(y_2))\ .
	\end{multline}
	This proves identity~\eqref{eq:hopfB}. 
\end{proof}

The next lemma concludes the proof of Proposition~\ref{prop:coha-surface-module-algebra}. 
\begin{lemma}\label{lem:hopfaction}
	Let 
	\begin{align}
		\mu\colon \HBMbulletA(\frakX_1) \otimes \HBMbulletA(\frakX_2) \longrightarrow \HBMbulletA(\frakX)
	\end{align}
	be the Hall product. Then one has 
	\begin{align}\label{eq:hopfX} 
		\overline{\nu}_\alpha(\mu(x_1\otimes x_2)) = \mu(\overline{\nu}^{(1)}_\alpha(x_1) \otimes x_2 + 
		x_1 \otimes \overline{\nu}^{(2)}_\alpha(x_2))
	\end{align}
	for any $x_i \in \HBMbulletA(\frakX_i)$, with $1\leq i\leq 2$. 
\end{lemma} 
\begin{proof} 
	Let 
	\begin{align}
		\begin{tikzcd}[ampersand replacement=\&]
			\dstackCohext \ar{d}[swap]{q} \ar{r}{p} \& 
			\frakX\\
			\frakX_1 \times \frakX_2 \& 
		\end{tikzcd} 
	\end{align}
	be the Hall convolution diagram, where $\dstackCohext$ denotes the derived moduli stack of extensions.
	
	This yields the pull-back diagrams 
	\begin{align}\label{eq:hopfdiagB}
		\begin{tikzcd}[ampersand replacement=\&]
			\overline{\rsv} \times \dstackCohext \ar{d}[swap]{\id_{\overline{\rsv}} \times q} \ar{r}{\ev_\sfex} \& 
			\dstackCohext \ar{d}{q} \\
			\overline{\rsv} \times \frakX_1 \times \frakX_2 \ar{r}{\ev_{12}}\& \frakX_1\times \frakX_2 
		\end{tikzcd} 
		\quad \text{and}\quad
		\begin{tikzcd}[ampersand replacement=\&]
			\overline{\rsv} \times \dstackCohext \ar{d}[swap]{\id_{\overline{\rsv}}\times p} \ar{r}{\ev_\sfex} \& \dstackCohext
			\ar{d}{p} \\
			\overline{\rsv} \times \frakX\ar{r}{\ev}\& \frakX 
		\end{tikzcd} \ .
	\end{align}
	
	Recall that $r_i\colon \overline{\rsv} \times \frakX_1\times \frakX_2 \to \overline{\rsv} \times \frakX_i$ denotes the 
	canonical projection for $1\leq i\leq 2$.
	Then diagrams \eqref{eq:hopfdiagB} yield the identity
	\begin{align}\label{eq:hopfC} 
		(\id_{\overline{\rsv}} \times q)^\ast((r_1^\ast\ch(\overline{\scrE}_1) + r_2^\ast\ch(\overline{\scrE}_2))
		\cup \ev_{\overline{\rsv}, 12}^\ast(\alpha)) = 
		(\id_{\overline{\rsv}} \times p)^\ast\ch(\overline{\scrE})\cup s^\ast(\alpha) \ ,
	\end{align}
	where $s\colon \overline{\rsv} \times \dstackCohext \to \overline{\rsv}$ is the canonical projection. 
	
	Using \cite[Theorem~A.5-(iv)]{Khan_VFC}, one has 
	\begin{align}
		q^! \,(\ev_{12})_\ast = (\ev_\sfex)_\ast\, (\id_{\overline{\rsv}} \times q)^! \ .
	\end{align}
	Hence identity~\eqref{eq:hopfB} yields 
	\begin{multline}\label{eq:hopfCB} 
		q^! \, \kappa(\overline{\nu}^{(1)}_\alpha(x_1) \otimes \overline{\nu}^{(2)}_\alpha(x_2))= \\
		q^! \, (\ev_{12})_\ast\big( ((r_1^\ast\ch(\overline{\scrE}_1) + r_2^\ast\ch(\overline{\scrE}_2))\cup \ev_{\overline{\rsv}, 12}^\ast(\alpha) )
		\cap \ev_{12}^!(\kappa(x_1\otimes x_2)))=\\
		(\ev_\sfex)_\ast\, (\id_{\overline{\rsv}} \times q)^!  \big( ((r_1^\ast\ch(\overline{\scrE}_1) + r_2^\ast\ch(\overline{\scrE}_2))\cup \ev_{\overline{\rsv}, 12}^\ast(\alpha) ) \cap \ev_{12}^!(\kappa(x_1\otimes x_2))\big)\ .
	\end{multline}
	By \cite[Proposition~4.10-(1)]{Porta_Yu_Non-archimedian_quantum}, one  has 
	\begin{multline}
		(\ev_\sfex)_\ast\, (\id_{\overline{\rsv}} \times q)^!  \big( ((r_1^\ast\ch(\overline{\scrE}_1) + r_2^\ast\ch(\overline{\scrE}_2))\cup \ev_{\overline{\rsv}, 12}^\ast(\alpha) ) \cap \ev_{12}^!(\kappa(x_1\otimes x_2))\big) =\\
		(\ev_\sfex)_\ast\, \big( (\id_{\overline{\rsv}} \times q)^\ast ((r_1^\ast\ch(\overline{\scrE}_1) + r_2^\ast\ch(\overline{\scrE}_2))\cup \ev_{\overline{\rsv}, 12}^\ast(\alpha) ) \cap (\id_{\overline{\rsv}} \times q)^! \, \ev_{12}^!(\kappa(x_1\otimes x_2))\big) \ .
	\end{multline} 
	Hence identities~\eqref{eq:hopfC} and \eqref{eq:hopfCB} yield
	\begin{multline} 
		q^! \, \kappa(\overline{\nu}^{(1)}_\alpha(x_1) \otimes \overline{\nu}^{(2)}_\alpha(x_2))= \\
		(\ev_\sfex)_\ast \big((\id_{\overline{\rsv}} \times p)^\ast\ch(\overline{\scrE}) \cup r^\ast(\alpha) )\cap 
		(\id_{\overline{\rsv}} \times q)^! \, \ev_{12}^!(\kappa(x_1\otimes x_2))\big) =\\
		(\ev_\sfex)_\ast \big((\id_{\overline{\rsv}} \times p)^\ast\ch(\overline{\scrE}) \cup s^\ast(\alpha) )\cap 
		\ev_\sfex^!\, q^! (\kappa(x_1\otimes x_2))\big) \ .
	\end{multline} 
	Furthermore, $s^\ast(\alpha) = (\id_{\overline{\rsv}}\times q)^\ast\, \ev^\ast(\alpha)$. Therefore we are left with
	\begin{align}\label{eq:hopfD}
		q^! \, \kappa(\overline{\nu}^{(1)}_\alpha(x_1) \otimes \overline{\nu}^{(2)}_\alpha(x_2))= 
		(\ev_\sfex)_\ast \big((\id_{\overline{\rsv}} \times p)^\ast(\ch(\overline{\scrE}) \cup \ev_\rsv^\ast(\alpha)) \cap \ev_\sfex^!\, q^! (\kappa(x_1\otimes x_2))\big) \ .
	\end{align}
	Applying $p_\ast$, we further obtain 
	\begin{align} 
		p_\ast\, q^! \, \kappa(\overline{\nu}^{(1)}_\alpha(x_1) \otimes \overline{\nu}^{(2)}_\alpha(x_2))= 
		p_\ast\, (\ev_\sfex)_\ast \big((\id_{\overline{\rsv}} \times p)^\ast(\ch(\overline{\scrE}) \cup \ev_\rsv^\ast(\alpha)) \cap 
		\ev_\sfex^!\, q^! (\kappa(x_1\otimes x_2))\big) \ .
	\end{align}
	Using the second diagram in Equation~\eqref{eq:hopfdiagB}, this yields 
	\begin{align} 
		p_\ast\, q^! \, \kappa(\overline{\nu}^{(1)}_\alpha(x_1) \otimes \overline{\nu}^{(2)}_\alpha(x_2))= 
		\ev_\ast\, (\id_{\overline{\rsv}} \times p)_\ast \big((\id_{\overline{\rsv}} \times p)^\ast(\ch(\overline{\scrE}) \cup \ev_\rsv^\ast(\alpha)) \cap 
		\ev_\sfex^!\, q^! (\kappa(x_1\otimes x_2))\big)\ .
	\end{align} 
	Applying \cite[Proposition~4.10-(2)]{Porta_Yu_Non-archimedian_quantum} to the right hand side of the above equation 
	yields 
	\begin{align} 
		p_\ast\, q^! \, \kappa(\overline{\nu}^{(1)}_\alpha(x_1) \otimes \overline{\nu}^{(2)}_\alpha(x_2))=
		\ev_\ast\, \big((\ch(\overline{\scrE}) \cup \ev_\rsv^\ast(\alpha)) \cap (\id_{\overline{\rsv}} \times p)_\ast\, \ev_\sfex^!\, q^! (\kappa(x_1\otimes x_2))\big)\ .
	\end{align} 
	Using again the second diagram in \eqref{eq:hopfdiagB},  we finally obtain 
	\begin{align} 
		p_\ast\, q^! \, \kappa(\overline{\nu}^{(1)}_\alpha(x_1) \otimes \overline{\nu}^{(2)}_\alpha(x_2))=
		\ev_\ast\, \big((\ch(\overline{\scrE}) \cup \ev_\rsv^\ast(\alpha)) \cap \ev^!\, p_\ast\, q^! (\kappa(x_1\otimes x_2))\big)\ ,
	\end{align}
	which proves the claim.
\end{proof}

\section{Appendix: Purity and formality}\label{sec:purity-formality}

Let $G$ be a complex linear group, $T$ a complex torus, and $X$ an algebraic $G\times T$-variety over $\C$. We consider the following quotient stacks
\begin{align}
	\frakX\coloneqq X/G\quad\text{and}\quad \frakX_T\coloneqq \frakX\times T\ . 
\end{align}
Let $A_T$ be a mixed Hodge module over $\frakX_T$. Applying the forgetful functor to $A_T$ we get a mixed Hodge module $A$ over the stack $\frakX$. The following is well-known. Let us recall a proof.
\begin{proposition} \label{prop:surjectivity}
	Assume that the cohomology group $\sfH^q(\frakX,A)$ is pure of weight $q$ for each integer $q$. Then, the following hold:
	\begin{enumerate}\itemsep0.2cm
		\item \label{item:surjectivity-1} $\sfH^q(\frakX_T,A_T)$ is pure of weight $q$ for each integer $q$.
		\item \label{item:surjectivity-2} If $\Hbullet(\frakX_T,A_T)$ is free as an $\Hbullet_T$-module, then, for each integer $q$, the restriction yields an isomorphism
		\begin{align}
			\sfH^q(\frakX_T,A_T)\otimes_{\Hbullet_T} \Q\longrightarrow \sfH^q(\frakX, A)\ .
		\end{align} 
	\end{enumerate}
\end{proposition}

\begin{proof}
	The Leray Spectral sequence for the fibration $\frakX_T\to B T$ is
	\begin{align}
		E_2^{p,q}\coloneqq \sfH^p(BT,\Q)\otimes \sfH^q(\frakX,A)\longrightarrow \sfH^{p+q}(\frakX_T,A_T)\ .
	\end{align}
	Since $\sfH^p_T$ and $\sfH^q(\frakX,A)$ are pure of weights $p$ and $q$, the term $E_2^{p,q}$ is pure of weight $p+q$. The differentials of the spectral sequence are strictly compatible with weights. Thus, they vanish. Hence $\sfH^{p+q}(\frakX_T,A_T)$ is pure of weight $p+q$. This proves \eqref{item:surjectivity-1}.
	
	Now, consider the Cartesian square
	\begin{align}
		\begin{tikzcd}[ampersand replacement=\&]
			\frakX \ar{r}\ar{d}\& \frakX_T\ar{d}\\
			\Spec(\C) \ar{r} \& BT
		\end{tikzcd}\ .
	\end{align}
	The Eilenberg-Moore spectral sequence attached to this diagram is
	\begin{align}
		E_2^{q,\bullet}\coloneqq\Tor_{\textrm{-}q}^{\Hbullet_T}(\Hbullet(\frakX_T,A_T),\Q)\longrightarrow \sfH^{\bullet+q}(\frakX,A)\ .
	\end{align}
	It degenerates at $E_2$, see, e.g., \cite[Theorem~1.3]{FW_weights}. Hence, we have
	\begin{align}
		\Hbullet(\frakX,A) = \Tor_{\textrm{-}\bullet}^{\Hbullet_T}(\Hbullet(\frakX_T,A_T),\Q)\ ,
	\end{align}
	proving \eqref{item:surjectivity-2}.
\end{proof}

\newpage

\providecommand{\bysame}{\leavevmode\hbox to3em{\hrulefill}\thinspace}

\end{document}